\documentclass[12pt]{article}
\usepackage{amsmath,amssymb,amsfonts}     
\usepackage{graphics}                     
\usepackage{graphicx}                     
\date{}
\newtheorem{lemma}{Lemma}[section]
\newtheorem{theorem}{Theorem}[section]

\newcommand{\adfQED}{\hfill $\square$}

\newcommand{\adfmod}[1]{~(\mathrm{mod}~#1)}

\newcommand{\adfhide}[1]{}

\newcommand{\adfTgap}{\vskip 2mm}

\newcommand{\adfPgap}{\vskip 2mm}

\newcommand{\adfind}{\makebox[8mm]{}}

\newenvironment{adfenumerate}{
\begin{enumerate}
\setlength{\itemsep}{0.5mm}
\setlength{\parskip}{0mm}
\setlength{\parsep}{0mm}
}{
\end{enumerate}
}

\begin{document}
\title{\bf Group divisible designs with block size four and type $g^u m^1$ -- III\\}
\author{A.~D.~Forbes\\
        School of Mathematics and Statistics,\\
        The Open University,\\
        Walton Hall, Milton Keynes MK7 6AA, UK\\
        }
\maketitle
\begin{abstract}
\noindent We deal with group divisible designs that have block size 4 and group type $g^u m^1$, where $g \equiv 2$ or 4 (mod 6).
We show that the necessary conditions for the existence of a 4-GDD of type $g^u m^1$ are sufficient when
$g$ = 14, 20, 22, 26, 28, 32, 34, 38, 40, 44, 46, 50, 52, 58, 62, 68, 76, 88, 92, 100,
104, 116, 124, 136, 152, 160, 176, 184, 200, 208, 224, 232, 248, 272, 304, 320, 368, 400, 448, 464 and 496.
Using these results we go on to show that the necessary conditions are sufficient for
$g = 2^t q^s$, $q$ = 19, 23, 25, 29, 31, $s$, $t$ = 1, 2, \dots, as well as for
$g = 2^t q$, $q$ = 2, 5, 7, 11, 13, 17, $t$ = 1, 2, \dots, with possible exceptions
$56^9 m^1$, $80^9 m^1$ and $112^9 m^1$ for a few large values of $m$.
\end{abstract}




\section{Introduction}\label{sec:Introduction}

A {\em group divisible design}, $K$-GDD, of type $g_1^{u_1} \dots g_r^{u_r}$ is an ordered triple ($V,\mathcal{G},\mathcal{B}$) such that:
\begin{adfenumerate}
\item[(i)]{$V$
is a base set of cardinality $u_1 g_1 + \dots + u_r g_r$;}
\item[(ii)]{$\mathcal{G}$
is a partition of $V$ into $u_i$ subsets of cardinality $g_i$, $i = 1, \dots, r$, called \textit{groups};}
\item[(iii)]{$\mathcal{B}$
is a non-empty collection of subsets of $V$ with cardinalities $k \in K$, called \textit{blocks}; and}
\item[(iv)]{each pair of elements from distinct groups occurs in precisely one block but no pair of
elements from the same group occurs in any block.}
\end{adfenumerate}
A {\em parallel class} in a group divisible design is a subset of the block set in which each element of the base set appears exactly once.
A $k$-GDD is called {\em resolvable}, and is denoted by $k$-RGDD, if the entire set of blocks can be partitioned into parallel classes.
If there exist $k$ mutually orthogonal Latin squares (MOLS) of side $q$, then there exists
a $(k+2)$-GDD of type $q^{k+2}$ and
a $(k+1)$-RGDD of type $q^{k+1}$, \cite[Theorem III.3.18]{AbelColbournDinitz2007}.
As is well known, there exist $q - 1$ MOLS of side $q$ whenever $q$ is a prime power.

The existence spectrum for $k$-GDDs has been completely determined for
3-GDDs of type $g^u$, \cite{Hanini1975}, \cite[Theorem IV.4.1]{Ge2007}, and
3-GDDs of type $g^u m^1$, \cite{ColbournHoffmanRees1992}, \cite[Theorem IV.4.4]{Ge2007}.
Furthermore, the necessary conditions for the existence of 4-GDDs of type $g^u$, namely
\begin{equation}
\label{eqn:4-GDD g^u necessary}
u \ge 4,~~ g(u - 1) \equiv 0 \adfmod{3} \textrm{~~and~~} g^2 u (u - 1) \equiv 0 \adfmod{12},
\end{equation}
are known to be sufficient except that there are no 4-GDDs of types $2^4$ and $6^4$, \cite{BrouwerSchrijverHanani1977}, \cite[Theorem IV.4.6]{Ge2007}.
The next case one would naturally consider is that of 4-GDDs of type $g^u m^1$, and here a partial solution has been achieved
in the sense that for each $g$ there are at most a small number of $u$ where existence remains undecided,
\cite{GeRees2002}, \cite{GeReesZhu2002}, \cite{GeLing2004}, \cite{GeRees2004}, \cite{GeLing2005}, \cite{Schuster2010},
\cite{WeiGe2013}, \cite{WeiGe2014}, \cite{Schuster2014}, \cite{WeiGe2015}, \cite{ForbesForbes2018}, \cite{Forbes2019}.

Assuming $m > 0$ and $m \neq g$,
the necessary conditions for the existence of a 4-GDD of type $g^u m^1$ are
\begin{equation}
\label{eqn:4-GDD g^u m^1 necessary}
\left\{\begin{array}{@{}ll}
&u \ge 4,~~ m \le g(u - 1)/2,\\
&gu \equiv g (u - 1) + m \equiv 0 \adfmod{3}, \textrm{~and~}\\
&g^2 u (u-1) + 2 g u m \equiv 0 \adfmod{12},
\end{array} \right.
\end{equation}
which simplify to
\begin{equation*}
\label{eqn:4-GDD g^u m^1 g = 2 or 4 mod 6 necessary}
u \ge 4,~~ m \le g(u - 1)/2,~~ u \equiv m - g \equiv0 \adfmod{3}
\end{equation*}
when $g \equiv 2$ or $4 \adfmod{6}$.
As for sufficiency, it is convenient to quote \cite[Theorem 1.2]{Forbes2019}.
%
\begin{theorem}
\label{thm:4-GDD g^u m^1 existence}
Let
\begin{eqnarray*}
\label{eqn:Q}
\mathcal{Q} &=& \{3, 5, 7, 11, 13, 17, 19, 23, 29, 31, 37, 41, 43\}, \\
\label{eqn:P}
\mathcal{P}   &=& \{n: n \equiv 1 \adfmod{4},\; n \ge 5\} \cup \mathcal{Q}.
\end{eqnarray*}
The necessary conditions (\ref{eqn:4-GDD g^u m^1 necessary}) for the existence of a 4-GDD of type $g^u m^1$ with $m > 0$ and $m \neq g$ are sufficient
except that there is no 4-GDD of type $2^6 5^1$ and except possibly for the following:
\begin{enumerate}
\item[]{$g \equiv 9 \adfmod{12}$, $g \ge 141$, $g/3$ is not divisible by any integer in $\mathcal{P}$, \\
      $u = 8$ and $m' < m < (7g - 21)/2$, where $m' = (17g - 6)/5$ if $g \equiv 3 \adfmod{30}$, $m' = 3g$ otherwise;}
\item[]{$g \equiv 1, 5 \adfmod{6}$, $g \ge 37$: \\
      $u = 9$, \\
      $u \in \{12, 24, 15, 27, 39, 51, 21, 33\}$ and $0 < m < g$, \\
      $u = 24$ when $g$ is not divisible by any integer in $\mathcal{P}$ and
      $m'' < m < (23g - 3)/2$, where $m'' = (56g - 6)/5$ if $g \equiv 1$ or $11 \adfmod{30}$, $m'' = 10g$ otherwise;}
\item[]{$g \equiv 2, 4 \adfmod{6}$, $g \ge 14$, $g \neq 16$:\\
      $u \in \{6, 9\}$, \\
      $u \in \{12, 15, 18, 21, 27\}$ and $0 < m < g$.}
\end{enumerate}
\end{theorem}
\noindent{\bf Proof}~ See \cite{Forbes2019} and the references therein.
~\adfQED

\adfPgap
\noindent In this paper we address residue classes $g \equiv 2,4 \adfmod{6}$. The main results are gathered together in the following.
\begin{theorem}
\label{thm:4-GDD g^u m^1 new existence}
The necessary conditions (\ref{eqn:4-GDD g^u m^1 necessary}) for the existence of a 4-GDD of type $g^u m^1$ with
$m > 0$ and $m \neq g$ are sufficient for
$$g = 2^t q,~~ q \in \{1, 5, 7, 11, 13, 17, 19, 23, 25, 29, 31\},~~ t = 0, 1, \dots,$$
and
$$g = 2^t q^s,~~ q \in \{19, 23, 25, 29, 31\},~~ s, t = 1, 2, \dots,$$
except that there is no 4-GDD of type $2^6 5^1$, and except possibly for the following:
\begin{adfenumerate}
\item[]{$g = 56$, $u = 9$, $m \in \{206, 209, 215, 218, 221\}$;}
\item[]{$g = 80$, $u = 9$, $m \in \{299, 311, 317\}$;}
\item[]{$g = 112$, $u = 9$, $m \in \{433, 439, 445\}$.}
\end{adfenumerate}
\end{theorem}
\noindent{\bf Proof}~ This follows from Theorems~\ref{thm:4-GDD g^u m^1 existence},
\ref{thm:4-GDD g^u m^1 g = 14 20 22 ...},
\ref{thm:4-GDD (2^t q)^u m^1 q = 2 5 7 ...},
\ref{thm:4-GDD (2^t q^s)^u m^1 q = 19 25 31} and
\ref{thm:4-GDD (2^t q^s)^u m^1 q = 23 29}.
~\adfQED


\section{New 4-GDDs}\label{sec:new 4-GDDs}

The 4-GDDs that we use to prove our theorems are given here, grouped into lemmas.
The blocks of the designs are generated by appropriate mappings from sets of base blocks,
which are collected in the appendix together with the instructions for expanding them.
If the point set of the 4-GDD has $v$ elements, it is represented by $Z_v = \{0, 1, \dots, v-1\}$ partitioned into groups as indicated.
The expression $a\adfmod{b}$ denotes the integer $n$ such that $0 \le n < b$ and $b \,|\, n - a$.
Where applicable, addition in $\mathbb{Z}_{gu/3} \times \mathbb{Z}_3$ is denoted by $\oplus$ and for brevity we represent element $(a,b)$ of this group by the number $3a+b$.
As an aid to checking the correctness of the designs we also provide the block generation instructions in a compact format: $(v,M,T)$,
where $v$ is the number of points, $M$ encodes the mapping and $T$ encodes the group type.
We leave it for the interested reader to determine the structure of $M$ and $T$ from the corresponding English descriptions.
It is possible that some of the designs might be constructible by other means.

\adfhide{To save space most of the lengthy appendix has been omitted.
However, the unabridged preprint of this paper, available on ArXiv at
\begin{center}
{\sf http://arxiv.org/abs/[[TBA]]},
\end{center}
includes the appendix in its entirety, and hence the proofs of the lemmas in this section are complete.}

\adfhide{\begin{center} 
\begin{tabular}{lr}
General designs:              &  19 \\
Designs for type $14^u m^1$:  &  21 \\
Designs for type $20^u m^1$:  &   3 \\
Designs for type $22^u m^1$:  &  31 \\
Designs for type $26^u m^1$:  &  32 \\
Designs for type $28^u m^1$:  &   3 \\
Designs for type $34^u m^1$:  &  15 \\
Designs for type $38^u m^1$:  &  32 \\
Designs for type $40^u m^1$:  &  11 \\
Designs for type $46^u m^1$:  &  24 \\
Designs for type $50^u m^1$:  &  16 \\
Designs for type $52^u m^1$:  &   1 \\
Designs for type $58^u m^1$:  &   5 \\
Designs for type $62^u m^1$:  &  13 \\
Designs for type $76^u m^1$:  &   2 \\
TOTAL:                        & 228
\end{tabular}
\end{center}}

\begin{lemma}
\label{lem:4-GDD general}
There exist 4-GDDs of types
$4^2 10^5$,
$4^3 10^4$,
$4^4 10^3$,
$4^5 10^2$,
$4^2 10^4 1^1$,
$8^2 2^8$,
$8^3 2^7$,
$8^4 2^6$,
$8^5 2^5$,
$8^6 2^4$,
$8^7 2^3$,
$8^2 2^{11}$,
$8^3 2^6 5^1$,
$8^5 2^4 5^1$,
$8^7 2^2 5^1$,
$8^5 14^1 20^1$,
$20^4 8^2 2^1$,
$20^4 8^2 5^1$,
$20^5 8^1 5^1$.
\end{lemma}
{\bf Proof}~ The designs are presented in Appendix~\ref{app:4-GDD general}.
~\adfQED

\begin{lemma}
\label{lem:4-GDD 14^u m^1}
There exist 4-GDDs of types
$14^{6} 8^1$,
$14^{6} 11^1$,
$14^{6} 17^1$,
$14^{6} 20^1$,
$14^{6} 23^1$,
$14^{6} 26^1$,
$14^{6} 29^1$,
$14^{6} 32^1$,
$14^{9} 11^1$,
$14^{9} 17^1$,
$14^{9} 20^1$,
$14^{9} 23^1$,
$14^{9} 26^1$,
$14^{9} 29^1$,
$14^{9} 32^1$,
$14^{9} 38^1$,
$14^{9} 41^1$,
$14^{9} 44^1$,
$14^{9} 47^1$,
$14^{9} 50^1$,
$14^{9} 53^1$.
\end{lemma}
{\bf Proof}~ The designs are presented in Appendix~\ref{app:4-GDD 14^u m^1}.
~\adfQED

\begin{lemma}
\label{lem:4-GDD 20^u m^1}
There exist 4-GDDs of types
$20^{9} 11^1$,
$20^{9} 17^1$,
$20^{9} 23^1$.
\end{lemma}
{\bf Proof}~ The designs are presented in Appendix~\ref{app:4-GDD 20^u m^1}.
~\adfQED

\begin{lemma}
\label{lem:4-GDD 22^u m^1}
There exist 4-GDDs of types
$22^{6} 25^1$,
$22^{6} 28^1$,
$22^{6} 31^1$,
$22^{6} 34^1$,
$22^{6} 37^1$,
$22^{6} 40^1$,
$22^{6} 43^1$,
$22^{6} 46^1$,
$22^{6} 49^1$,
$22^{6} 52^1$,
$22^{9} 7^1$,
$22^{9} 10^1$,
$22^{9} 13^1$,
$22^{9} 16^1$,
$22^{9} 19^1$,
$22^{12} 7^1$,
$22^{12} 10^1$,
$22^{12} 13^1$,
$22^{12} 16^1$,
$22^{12} 19^1$,
$22^{15} 4^1$,
$22^{15} 10^1$,
$22^{15} 13^1$,
$22^{15} 16^1$,
$22^{15} 19^1$,
$22^{18} 4^1$,
$22^{18} 7^1$,
$22^{18} 10^1$,
$22^{18} 13^1$,
$22^{18} 16^1$,
$22^{18} 19^1$.
\end{lemma}
{\bf Proof}~ The designs are presented in Appendix~\ref{app:4-GDD 22^u m^1}.
~\adfQED

\begin{lemma}
\label{lem:4-GDD 26^u m^1}
There exist 4-GDDs of types
$26^{6} 8^1$,
$26^{6} 11^1$,
$26^{6} 14^1$,
$26^{6} 17^1$,
$26^{6} 20^1$,
$26^{6} 23^1$,
$26^{6} 29^1$,
$26^{6} 32^1$,
$26^{6} 35^1$,
$26^{6} 38^1$,
$26^{6} 41^1$,
$26^{6} 44^1$,
$26^{6} 47^1$,
$26^{6} 50^1$,
$26^{6} 53^1$,
$26^{6} 56^1$,
$26^{6} 59^1$,
$26^{6} 62^1$,
$26^{9} 11^1$,
$26^{9} 14^1$,
$26^{9} 17^1$,
$26^{9} 20^1$,
$26^{9} 23^1$,
$26^{12} 14^1$,
$26^{12} 17^1$,
$26^{12} 20^1$,
$26^{12} 23^1$,
$26^{15} 17^1$,
$26^{15} 20^1$,
$26^{15} 23^1$,
$26^{18} 20^1$,
$26^{18} 23^1$.
\end{lemma}
{\bf Proof}~ The designs are presented in Appendix~\ref{app:4-GDD 26^u m^1}.
~\adfQED

\begin{lemma}
\label{lem:4-GDD 28^u m^1}
There exist 4-GDDs of types
$28^{9} 19^1$,
$28^{9} 25^1$,
$28^{9} 31^1$.
\end{lemma}
{\bf Proof}~ The designs are presented in Appendix~\ref{app:4-GDD 28^u m^1}.
~\adfQED

\begin{lemma}
\label{lem:4-GDD 34^u m^1}
There exist 4-GDDs of types
$34^{6}  7^1$,
$34^{6} 10^1$,
$34^{6} 13^1$,
$34^{6} 16^1$,
$34^{6} 19^1$,
$34^{6} 22^1$,
$34^{6} 73^1$,
$34^{6} 76^1$,
$34^{6} 79^1$,
$34^{6} 82^1$,
$34^{9} 19^1$,
$34^{12} 25^1$,
$34^{12} 28^1$,
$34^{12} 31^1$,
$34^{15} 31^1$.
\end{lemma}
{\bf Proof}~ The designs are presented in Appendix~\ref{app:4-GDD 34^u m^1}.
~\adfQED

\begin{lemma}
\label{lem:4-GDD 38^u m^1}
There exist 4-GDDs of types
$38^{6} 5^1$,
$38^{6} 8^1$,
$38^{6} 11^1$,
$38^{6} 14^1$,
$38^{6} 17^1$,
$38^{6} 20^1$,
$38^{6} 23^1$,
$38^{6} 26^1$,
$38^{6} 83^1$,
$38^{6} 86^1$,
$38^{6} 89^1$,
$38^{6} 92^1$,
$38^{9} 11^1$,
$38^{9} 14^1$,
$38^{9} 17^1$,
$38^{9} 20^1$,
$38^{9} 23^1$,
$38^{12} 14^1$,
$38^{12} 17^1$,
$38^{12} 20^1$,
$38^{12} 23^1$,
$38^{15} 17^1$,
$38^{15} 20^1$,
$38^{15} 23^1$,
$38^{15} 26^1$,
$38^{15} 29^1$,
$38^{18} 20^1$,
$38^{18} 23^1$,
$38^{18} 26^1$,
$38^{18} 29^1$,
$38^{18} 32^1$,
$38^{18} 35^1$.
\end{lemma}
{\bf Proof}~ The designs are presented in Appendix~\ref{app:4-GDD 38^u m^1}.
~\adfQED

\begin{lemma}
\label{lem:4-GDD 40^u m^1}
There exist 4-GDDs of types
$40^{6} 73^1$,
$40^{6} 79^1$,
$40^{6} 82^1$,
$40^{6} 91^1$,
$40^{6} 94^1$,
$40^{6} 97^1$,
$40^{9} 139^1$,
$40^{9} 142^1$,
$40^{9} 151^1$,
$40^{9} 154^1$,
$40^{9} 157^1$.
\end{lemma}
{\bf Proof}~ The designs are presented in Appendix~\ref{app:4-GDD 40^u m^1}.
~\adfQED

\begin{lemma}
\label{lem:4-GDD 46^u m^1}
There exist 4-GDDs of types
$46^{6} 13^1$,
$46^{6} 16^1$,
$46^{6} 19^1$,
$46^{6} 22^1$,
$46^{6} 73^1$,
$46^{6} 76^1$,
$46^{6} 79^1$,
$46^{6} 82^1$,
$46^{6} 85^1$,
$46^{6} 88^1$,
$46^{6} 91^1$,
$46^{6} 94^1$,
$46^{6} 97^1$,
$46^{6} 100^1$,
$46^{6} 103^1$,
$46^{6} 106^1$,
$46^{6} 109^1$,
$46^{6} 112^1$,
$46^{9} 19^1$,
$46^{12} 25^1$,
$46^{12} 28^1$,
$46^{12} 31^1$,
$46^{12} 34^1$,
$46^{15} 31^1$.
\end{lemma}
{\bf Proof}~ The designs are presented in Appendix~\ref{app:4-GDD 46^u m^1}.
~\adfQED

\begin{lemma}
\label{lem:4-GDD 50^u m^1}
There exist 4-GDDs of types
$50^{6} 8^1$,
$50^{6} 11^1$,
$50^{6} 14^1$,
$50^{6} 17^1$,
$50^{6} 23^1$,
$50^{6} 26^1$,
$50^{6} 29^1$,
$50^{6} 32^1$,
$50^{6} 38^1$,
$50^{6} 113^1$,
$50^{6} 116^1$,
$50^{6} 119^1$,
$50^{6} 122^1$,
$50^{9} 11^1$,
$50^{9} 17^1$,
$50^{9} 23^1$.
\end{lemma}
{\bf Proof}~ The designs are presented in Appendix~\ref{app:4-GDD 50^u m^1}.
~\adfQED

\begin{lemma}
\label{lem:4-GDD 52^u m^1}
There exists a 4-GDD of type
$52^{15} 34^1$.
\end{lemma}
{\bf Proof}~ The design is presented in Appendix~\ref{app:4-GDD 52^u m^1}.
~\adfQED

\begin{lemma}
\label{lem:4-GDD 58^u m^1}
There exist 4-GDDs of types
$58^{6} 133^1$,
$58^{6} 136^1$,
$58^{6} 139^1$,
$58^{6} 142^1$,
$58^{9} 19^1$.
\end{lemma}
{\bf Proof}~ The designs are presented in Appendix~\ref{app:4-GDD 58^u m^1}.
~\adfQED

\begin{lemma}
\label{lem:4-GDD 62^u m^1}
There exist 4-GDDs of types
$62^{6} 23^1$,
$62^{6} 26^1$,
$62^{6} 29^1$,
$62^{6} 32^1$,
$62^{6} 35^1$,
$62^{6} 38^1$,
$62^{6} 41^1$,
$62^{6} 44^1$,
$62^{6} 47^1$,
$62^{6} 143^1$,
$62^{6} 146^1$,
$62^{6} 149^1$,
$62^{6} 152^1$.
\end{lemma}
{\bf Proof}~ The designs are presented in Appendix~\ref{app:4-GDD 62^u m^1}.
~\adfQED

\begin{lemma}
\label{lem:4-GDD 76^u m^1}
There exist 4-GDDs of types
$76^{6} 13^1$,
$76^{6} 16^1$.
\end{lemma}
{\bf Proof}~ The designs are presented in Appendix~\ref{app:4-GDD 76^u m^1}.
~\adfQED


\section{4-GDDs of type $g^u m^1$ with $g \equiv 2, 4 \adfmod{6}$}\label{sec:Theorems g = 2, 4 (mod 6)}

A {\em double group divisible design}, $k$-DGDD, of type
$(g_1, h_1^w)^{u_1} (g_2, h_2^w)^{u_2} \dots (g_r, h_r^w)^{u_r}$, $g_i = w h_i$, $i = 1, 2, \dots, r$, is an ordered quadruple
($V,\mathcal{G},\mathcal{H},\mathcal{B}$) such that:
\begin{adfenumerate}
\item[(i)]{$V$
is a base set of cardinality $u_1 g_1 + u_2 g_2 + \dots + u_r g_r$;}
\item[(ii)]{$\mathcal{G}$
is a partition of $V$ into $u_i$ subsets of cardinality $g_i$, $i = 1, 2, \dots, r$, the \textit{groups};}
\item[(iii)]{$\mathcal{H}$
is a partition of $V$ into $w$ subsets of cardinality $h_1 u_1 + h_2 u_2 + \dots + h_r u_r$, the \textit{holes},
where each group of size $g_i$ intersects each hole in $h_i$ points, $i = 1, 2, \dots, r$;}
\item[(iv)]{$\mathcal{B}$
is a non-empty collection of subsets of $V$ of cardinality $k$, the \textit{blocks};}
\item[(v)]{%
each pair of elements of $V$ not in the same group and not in the same hole occurs in precisely one block; and}
\item[(vi)]{%
a pair of elements from the same group or from the same hole does not occur in any block.}
\end{adfenumerate}
We quote a theorem concerning the existence of 4-DGDDs of type $(h v, h^v)^t$ and
we describe several constructions.

\begin{theorem}
\label{thm:4-DGDD existence}
There exists a 4-DGDD of type $(h v, h^v)^t$ if and only if $t, v \ge 4$ and
$(t - 1)(v - 1) h \equiv 0 \adfmod{3}$ except for $(h,v,t) \in \{(1,4,6), (1,6,4)\}$.
\end{theorem}
\noindent {\bf Proof~} See \cite{GeWei2004} and \cite{CaoWangWei2009}.
~\adfQED

\begin{theorem}
\label{thm:(u+1)-GDD q^(u+1) + h^u d^1 d in D}
Let $h$, $u$ and $q$ be positive integers and let $D$ be a set of non-negative integers.
Suppose there exists a $(u + 1)$-GDD of type $q^{u + 1}$.
Suppose also that there exists a 4-GDD of type $h^u d^1$ for each $d \in D$.
Then there exists a 4-GDD of type $(hq)^u (d_1 + d_2 + \dots + d_q)^1$ for $d_1, d_2, \dots, d_q \in D$.
\end{theorem}
\noindent {\bf Proof~}
Take the $(u+1)$-GDD of type $q^{u+1}$,
inflate the $q$ points in one group by $d_1, d_2, \dots, d_q \in D$ and inflate all other points by $h$.
Overlay the inflated blocks with 4-GDDs of types $h^u d_i^1$, $i = 1, 2, \dots, q$, as appropriate.
~\adfQED

\begin{theorem}
\label{thm:(u+1)-GDD q^(u+1) + h^u d^1}
Let $h$, $u$ and $q$ be positive integers with $h \equiv 2 \textrm{~or~} 4 \adfmod{6}$, $u \equiv 0 \adfmod 3$, $u \ge 6$ and $q$ not divisible by $3$.
Suppose there exists a $(u + 1)$-GDD of type $q^{u + 1}$.
Suppose also there exists a 4-GDD of type $h^u d^1$ for each admissible $d$, and denote the smallest such $d$ by $d_{\min}$.
Then there exists a 4-GDD of type $(hq)^u m^1$ for every $m \equiv hq \adfmod{3}$ in the range $q d_{\min} \le m \le q h (u - 1)/2$.
\end{theorem}
\noindent {\bf Proof~}
Use Theorem~\ref{thm:(u+1)-GDD q^(u+1) + h^u d^1 d in D} with $D$ as the set of admissible $d$ for a 4-GDD of type $h^u d^1$.
Then $D = \{d_{\min}, d_{\min} + 3, \dots, h(u - 1)/2\}$,
and $m = d_1 + d_2 + \dots + d_q$, $d_1$, $d_2$, \dots, $d_q$ $\in D$, can take any value in $\{qd_{\min}$, $qd_{\min} + 3$, \dots, $qh(u - 1)/2\}$.
~\adfQED

\adfPgap
In Theorem~\ref{thm:(u+1)-GDD q^(u+1) + h^u d^1} $d_{\min}$ is either 1 or 2, depending on $h \adfmod{3}$.
Hence the lower bound on the size of $m$ is $q$ if $h \equiv 4 \adfmod{6}$, or $2q$  if $h \equiv 2 \adfmod{6}$.
Furthermore, the maximum admissible value $m = q h (u - 1)/2$ is actually attained.
Note that this construction cannot be used if $(h, u) = (2, 6)$.

\begin{theorem}
\label{thm:4-GDD (gi)^(ti) + (hgi)^u m^1}
Suppose there exists a 4-GDD of type $g_1^{t_1} g_2^{t_2} \dots g_r^{t_r}$.
Let $h$, $m$ and $u$ be integers such that $h \ge 1$, $m \ge 0$, $u \ge 4$ and $(u,h) \neq (6,1)$.
Suppose there exists a 4-GDD of type $(hg_i)^u m^1$ for $i = 1$, $2$, \dots, $r$.
Then there exists a 4-GDD of type $(hg_1 t_1 + hg_2 t_2 + \dots + hg_r t_r)^u m^1$.
\end{theorem}
\noindent {\bf Proof~}
This is \cite[Constructions 1.8--1.10]{GeReesZhu2002}, or see \cite[Corollary 2.7 \& Construction 2.8]{WeiGe2015}.
By Theorem~\ref{thm:4-DGDD existence}, since $u \ge 4$ and $(u,h) \neq (6,1)$, there exists a 4-DGDD of type $(hu, h^u)^4$.
Take the 4-GDD of type $g_1^{t_1} g_2^{t_2} \dots g_r^{t_r}$,
inflate each point by a factor of $hu$ and overlay each inflated block with the 4-DGDD
to create a 4-DGDD of type $(u hg_1, (hg_1)^u)^{t_1}$ $(u hg_2, (hg_2)^u)^{t_2}$ \dots\ $(u hg_r, (hg_r)^u)^{t_r}$.
Adjoin $m$ points and
overlay the groups plus the $m$ points with 4-GDDs of types $(hg_1)^u m^1$, $(hg_2)^u m^1$, \dots, $(hg_r)^u m^1$, as appropriate.
~\adfQED

\begin{theorem}
\label{thm:4-GDD (rg)^u (rm)^1}
Suppose there exists a 4-GDD of type $g^u m^1$ and let $r \ge 3$, $r \neq 6$ be an integer.
Then there exists a 4-GDD of type $(rg)^u (rm)^1$.
\end{theorem}
\noindent {\bf Proof~}
Take the 4-GDD of type $g^u m^1$, inflate its points by $r$ and overlay the inflated blocks with 4-GDDs of type $r^4$.
~\adfQED

\adfPgap
The next three theorems describe specific constructions that we use for proving the first of our main results, Theorem~\ref{thm:4-GDD g^u m^1 g = 14 20 22 ...}.
Theorem~\ref{thm:4-GDD g^u m^1 u = 21 27} deals with 4-GDDs of type $g^u m^1$ where $g \equiv 2$ or $4 \adfmod{6}$ and $u = 21$ or $27$.
Theorem~\ref{thm:4-GDD h^u d^1 + (g/4)^u m^1 -> g^u m^1} helps with $g = 8q$ or $g = 16q$ whenever there exists a 10-GDD of type $q^{10}$, and
Theorem~\ref{thm:4-GDD construction 2-5-8} exploits the availability of 4-GDDs of type $2^r 8^s$ with $r, s \ge 2$.

\begin{theorem}
\label{thm:4-GDD g^u m^1 u = 21 27}
Suppose $g \equiv 2$ or $4 \adfmod{6}$ and
suppose there exists a 4-GDD of type $g^{6} d^1$ for all admissible $d < g$.
Then there exist 4-GDDs of types $g^{21} m_1^1$ and $g^{27} m_2^1$ for all admissible $m_1$ and $m_2$ respectively.
\end{theorem}
\noindent {\bf Proof~}
By Theorem~\ref{thm:4-GDD g^u m^1 existence}, we may assume $m_1, m_2 < g$ and $m_1 \equiv m_2 \equiv g \adfmod{3}$.
Take a 4-GDD of type $g^{15} (6g + m_1)^1$, which exists by Theorem~\ref{thm:4-GDD g^u m^1 existence}
since $g \le 6g + m_1 \le 7g$, and overlay the big group with a 4-GDD of type $g^6 m_1^1$.
The result is a 4-GDD of type $g^{21} m_1^1$.
Repeat the construction with a 4-GDD of type $g^{21}(6g + m_2)^1$ to obtain a 4-GDD of type $g^{27} m_2^1$.
~\adfQED

\begin{theorem}
\label{thm:4-GDD h^u d^1 + (g/4)^u m^1 -> g^u m^1}
Suppose $g = hq$, where $h \in \{8,16\}$ and $q$ is not divisible by $3$.
Suppose also there exists a 10-GDD of type $q^{10}$ as well as
a 4-GDD of type $(g/4)^v d^1$ for $v \in \{6, 9, 12, 15, 18\}$ and every admissible $d \le g - 3$.
Then there exists a 4-GDD of type $g^u m^1$ for all admissible $u$ and $m$.
\end{theorem}
\noindent {\bf Proof~}
By Theorems~\ref{thm:4-GDD g^u m^1 existence} and \ref{thm:4-GDD g^u m^1 u = 21 27}, it suffices to consider
$u \in \{6,9\}$ and all admissible $m$ as well as $u \in \{12,15,18\}$ and all admissible $m \le g - 3$.

Take a 4-DGDD of type $(g, (g/4)^4)^u$, which exists by Theorem~\ref{thm:4-DGDD existence}.
Adjoin $m \le \min\{g(u - 1)/8, g - 3\}$ points, where $m \equiv g \adfmod{3}$, and
overlay the holes plus the new points with 4-GDDs of type $(g/4)^u m^1$.
This creates a 4-GDD of type $g^u m^1$.
The upper bound for $m$ is $5g/8$ when $u = 6$ and $g - 3$ when $u \ge 9$.
Hence the theorem is proved for $u \in \{12, 15, 18\}$.

For $u \in \{6, 9\}$, we have $u < q$.
Moreover, by Theorem~\ref{thm:4-GDD g^u m^1 existence}, there exists a 4-GDD of type $h^u d^1$ for $d \equiv h \adfmod{3}$ and $0 < d \le h(u - 1)/2$.
Hence, by Theorem~\ref{thm:(u+1)-GDD q^(u+1) + h^u d^1},
there exists a 4-GDD of type $g^u m^1$ for every admissible $m \ge 2q = 2g/h$.
This range overlaps the range of $m$ in the 4-DGDD construction, and hence the proof is complete.
~\adfQED

\begin{theorem}
\label{thm:4-GDD construction 2-5-8}
Let $q$, $u$ and $c$ be integers.
Suppose $g \equiv 2$ or $4 \adfmod{6}$, $g = 8c + 2(q - c)$, $0 \le c \le q$, $9 \le u \le q$, $ u \equiv 0 \adfmod{3}$, and
suppose there exists a $(u + 1)$-GDD of type $q^{u + 1}$.

{\normalfont (i)} If there exists a 4-GDD of type $8^{\gamma} 2^{u + 1 - \gamma}$ for each integer $\gamma$ satisfying
\begin{equation}
\label{eqn:4-GDD construction 2-5-8 (i)}
   \max\{2, u + c - q\} \le \gamma \le \min\{c, u - 1\},
\end{equation}
then there exists a 4-GDD of type $g^u m^1$ for $2q \le m \le 2q + u - 3$, $m \equiv g \adfmod{3}$.

{\normalfont (ii)} If there exists a 4-GDD of type $8^{\gamma} 2^{u + 1 - \gamma}$ for each integer $\gamma$ satisfying
\begin{equation}
\label{eqn:4-GDD construction 2-5-8 (ii)}
   \max\{2, u + c - q\} \le \gamma \le \min\{c + 1, u - 1\},
\end{equation}
then there exists a 4-GDD of type $g^u m^1$ for $2q \le m \le 8q + u - 9$, $m \equiv g \adfmod{3}$.
\end{theorem}
\noindent {\bf Proof~}
First we deal with (ii).
Take the $(u + 1)$-GDD of type $q^{u + 1}$.
Choose a point, $A$, and inflate it by a factor chosen from $\{2, 5, \dots, u - 1\}$.
If $c > 0$, choose $c$ blocks containing $A$ and inflate their points other than $A$ by 8.
Inflate all points other than $A$ in the same group as $A$ by factors chosen from $\{2, 8\}$.
Inflate the remaining points, if any, by 2.
Observe that blocks which do not contain $A$ and have
inflation pattern $8^{r} 2^{s}$ with both $r \ge 2$ and $s \ge 2$
actually have inflation pattern $8^{\gamma} 2^{u + 1 - \gamma}$ with $\gamma$
satisfying both $2 \le \gamma \le \min\{c + 1, u + 1\}$ and $2 \le u + 1 - \gamma \le \min\{q - c + 1, u + 1\}$,
equivalent to (\ref{eqn:4-GDD construction 2-5-8 (ii)}).
The diagram shows the $u + 1$ groups as columns with the
$q$ blocks of the parallel class defined by $A$ on the left of the vertical line.
\begin{center}
\begin{tabular}{ccccc|l}
 2 & 2 & \dots & 2 & 2 & $\{2, 8\}$ \\
\dots & \dots & \dots & \dots & \dots & \dots  \\
 2 & 2 & \dots & 2 & 2 & $\{2, 8\}$ \\
 8 & 8 & \dots & 8 & 8 & $\{2, 8\}$ \\
\dots & \dots & \dots & \dots & \dots & \dots  \\
 8 & 8 & \dots & 8 & 8 & $\{2, 8\}$ \\
 8 & 8 & \dots & 8 & 8 & $\{2, 5, \dots, u - 1\}$
\end{tabular}
\end{center}
Overlay the inflated blocks with 4-GDDs of types
$2^{u} a^1$, $8^u a^1$, $a \in \{2, 5, \dots, u - 1\}$ from Theorem~\ref{thm:4-GDD g^u m^1 existence}, as well as
$8^{\gamma} 2^{u + 1 - \gamma}$ for $\gamma$ satisfying (\ref{eqn:4-GDD construction 2-5-8 (ii)}).

The construction for (i) is the same except that the points other than $A$ in the group containing $A$ are all inflated by 2,
and consequently the range of $m$ is smaller.
The blocks which do not contain $A$ and have
inflation pattern $8^{r} 2^{s}$ with $r, s \ge 2$
have inflation pattern $8^{\gamma} 2^{u + 1 - \gamma}$ with $\gamma$
satisfying $2 \le \gamma \le \min\{c, u\}$ and $2 \le u + 1 - \gamma \le \min\{q - c + 1, u + 1\}$,
or, equivalently, (\ref{eqn:4-GDD construction 2-5-8 (i)}).
~\adfQED

\adfPgap
We are now ready to prove the main theorems of this section.

\begin{theorem}
\label{thm:4-GDD g^u m^1 g = 14 20 22 ...}
Suppose
\begin{eqnarray*}
g &\in& \{14, 20, 22, 26, 28, 32, 34, 38, 40, 44, 46, 50, 52, 56, 58, 62, 68, 76, 80, \\
  &  &~~  88, 92, 100, 104, 112, 116, 124, 136, 152, 160, 176, 184, 200, 208, 224, \\
  &  &~~ 232, 248, 272, 304, 320, 368, 400, 448, 464, 496\},
\end{eqnarray*}
$u \equiv 0 \adfmod{3}$, $m \equiv g \adfmod{3}$ and $m \le g(u - 1)/2$.
Then there exists a 4-GDD of type $g^u m^1$ except possibly for the following:
\begin{adfenumerate}
\item[]{$g = 56$, $u = 9$, $m \in \{206, 209, 215, 218, 221\}$;}
\item[]{$g = 80$, $u = 9$, $m \in \{299, 311, 317\}$;}
\item[]{$g = 112$, $u = 9$, $m \in \{433, 439, 445\}$.}
\end{adfenumerate}
\end{theorem}
\noindent {\bf Proof~}
By Theorem~\ref{thm:4-GDD g^u m^1 existence}, \cite{BrouwerSchrijverHanani1977}, \cite{GeReesZhu2002} and Theorem~\ref{thm:4-GDD g^u m^1 u = 21 27},
we need to consider only:
\begin{adfenumerate}
\item[]{$u = 6$, $m \equiv g \adfmod{3}$, $3 < m < 5g/2$, $m \neq g$;}
\item[]{$u = 9$, $m \equiv g \adfmod{3}$, $3 < m < 4g$, $m \neq g$;}
\item[]{$u \in \{12, 15, 18\}$, $m \equiv g \adfmod{3}$, $3 < m < g$.}
\end{adfenumerate}
We deal with each value of $g$ in turn. For the existence of 4-GDD types used in the constructions,
see either Theorem~\ref{thm:4-GDD g^u m^1 existence}, or \cite[Theorem IV.4.6]{Ge2007} for types $h^v$, or earlier results of this theorem.
For whenever Theorem~\ref{thm:(u+1)-GDD q^(u+1) + h^u d^1 d in D} or Theorem~\ref{thm:(u+1)-GDD q^(u+1) + h^u d^1} is invoked, and elsewhere,
recall that if $q$ is a prime power, then there exists $q - 1$ MOLS of side $q$ and hence a $(q + 1)$-GDD of type $q^{q + 1}$ and a $d$-RGDD of type $q^{d}$ for $2 \le d \le q$.
For the existence of 4-DGDDs of type $(hv, h^v)^u$, see Theorem~\ref{thm:4-DGDD existence}.

\adfTgap
\noindent {\bf\boldmath $g = 14$~}
For type $14^6 5^1$, see \cite[Lemma 4.13]{WeiGe2015}.

For type $14^u m^1$, $u = 9$ and $m \in \{5, 8\}$, or $u \in \{12, 15, 18\}$ and $m \in \{5, 8, 11\}$,
take a 4-DGDD of type $(14, 2^7)^u$, adjoin $m$ points and overlay the holes plus the new points with 4-GDDs of type $2^u m^1$.

For type $14^9 35^1$, use Theorem~\ref{thm:4-GDD (rg)^u (rm)^1} with a 4-GDD of type $2^9 5^1$.

For type $14^u m^1$, where \\
\adfind $u = 6$ and $m \in \{8, 11, \dots, 32\} \setminus \{14\}$, or \\
\adfind $u = 9$ and $m \in \{11, 14, \dots, 53\} \setminus \{14, 35\}$,
see Lemma~\ref{lem:4-GDD 14^u m^1}.

\adfTgap
\noindent {\bf\boldmath $g = 20$~}
For all except types $20^9 11^1$, $20^9 17^1$, $20^9 23^1$ and $20^{21} 17^1$, see \cite[Theorem 5.11]{Schuster2014}.

For $20^{21} 17^1$, use Theorem~\ref{thm:4-GDD g^u m^1 u = 21 27}.

The other three exceptions are given by Lemma~\ref{lem:4-GDD 20^u m^1}.

\adfTgap
\noindent {\bf\boldmath $g = 22$~}
For type $22^6 m^1$, $m \le 19$, see \cite[Lemma 4.10]{WeiGe2015}.

For types $22^9 4^1$ and $22^{12} 4^1$, see \cite[Lemma 4.12]{WeiGe2015}.

For type $22^9 m^1$, $m \ge 22$, use Theorem~\ref{thm:(u+1)-GDD q^(u+1) + h^u d^1} with $h = 2$, $u = 9$ and $q = 11$.

For type $22^{15} 7^1$, take a 4-DGDD of type $(22, 1^{22})^{15}$, adjoin 7 points and overlay holes plus the new points with 4-GDDs of type $1^{15} 7^1$.

For type $22^u m^1$, where \\
\adfind $u = 6$ and $m \in \{25, 28, \dots, 52\}$, or \\
\adfind $u = 9$ and $m \in \{7, 10, 13, 16, 19\}$, or \\
\adfind $u = 12$ and $m \in \{7, 10, 13, 16, 19\}$, or \\
\adfind $u = 15$ and $m \in \{4, 10, 13, 16, 19\}$, or \\
\adfind $u = 18$ and $m \in \{4, 7, 10, 13, 16, 19\}$,
see Lemma~\ref{lem:4-GDD 22^u m^1}.

\adfTgap
\noindent {\bf\boldmath $g = 26$~}
For type $26^6 5^1$, see \cite[Lemma 4.13]{WeiGe2015}.

For type $26^9 m^1$, $m \ge 26$, use Theorem~\ref{thm:(u+1)-GDD q^(u+1) + h^u d^1} with $h = 2$ and $q = 13$.

For type $26^u m^1$, $u \in \{9, 12, 15, 18\}$, $m \le u - 1$,
take a 4-DGDD of type $(26, 2^{13})^u$, adjoin $m$ points and overlay the holes plus the new points with 4-GDDs of type $2^u m^1$.

For type $26^u m^1$, where \\
\adfind $u = 6$ and $m \in \{8, 11, \dots, 62\} \setminus \{26\}$, or \\
\adfind $u = 9$ and $m \in \{11, 14, 17, 20, 23\}$, or \\
\adfind $u = 12$ and $m \in \{14, 17, 20, 23\}$, or \\
\adfind $u = 15$ and $m \in \{17, 20, 23\}$, or \\
\adfind $u = 18$ and $m \in \{20, 23\}$,
see Lemma~\ref{lem:4-GDD 26^u m^1}.

\adfTgap
\noindent {\bf\boldmath $g = 28$~}
For all except types $28^9 19^1$, $28^9 25^1$ and $28^9 31^1$, see \cite[Theorem 5.17]{Schuster2014}.

The three exceptions are given by Lemma~\ref{lem:4-GDD 28^u m^1}.

\adfTgap
\noindent {\bf\boldmath $g = 32$~}
For all except type $32^9 m^1$, see \cite[Theorem 3.12]{Schuster2010}.

For type $32^9 m^1$, $m \le 32$,
take a 4-DGDD of type $(32, 8^4)^9$, adjoin $m$ points and overlay the holes plus the new points with 4-GDDs of type $8^9 m^1$.

For type $32^9 m^1$, $m \ge 32$, use Theorem~\ref{thm:(u+1)-GDD q^(u+1) + h^u d^1} with $h = 2$, $u = 9$ and $q = 16$.

\adfTgap
\noindent {\bf\boldmath $g = 34$~}
For type $34^u m^1$, $u \in \{9, 12, 15, 18\}$, $m \le 2(u - 1)$,
use Theorem~\ref{thm:4-GDD (gi)^(ti) + (hgi)^u m^1} with 4-GDDs of types $4^6 10^1$, $4^u m^1$ and $10^u m^1$.

For type $34^6 4^1$, see \cite[Lemma 4.12]{WeiGe2015}.

For type $34^6 m^1$, $25 \le m \le 70$:
Take a 7-GDD of type $7^7$,
select a point, $A$, and inflate it by a factor of $a_0$ chosen from $\{1, 4, 7, 10\}$.
Select a block containing $A$ and inflate all other points in this block by 10.
Inflate all other points in the same group as $A$ by $a_j$, where each $a_j$ is chosen from $\{4,10\}$, $j$ = 1, 2, \dots, 6.
Inflate the remaining points by 4.
\begin{center}
\begin{tabular}{ccccccl}
4 & 4 & 4 & 4 & 4 & 4 & $a_6 \in \{4,10\}$ \\
\dots & \dots & \dots & \dots & \dots & \dots & \dots \\
4 & 4 & 4 & 4 & 4 & 4 & $a_1 \in \{4,10\}$ \\
10 & 10 & 10 & 10 & 10 & 10 & $a_0 \in \{1, 4, 7, 10\}$
\end{tabular}
\end{center}
Overlay the inflated blocks with 4-GDDs of types $4^6 a_0^1$, $10^6 a_0^1$, $4^7$, $4^6 10^1$ and $4^5 10^2$, as appropriate;
see Lemma~\ref{lem:4-GDD general} for type $4^5 10^2$.
The result is a 4-GDD of type $34^6 m^1$ where,
by suitable choice of $a_0$, $a_1$, \dots, $a_6$, $m = \sum_{i=0}^6 a_i$ can take any value in $\{25, 28, \dots, 70\}$.

For type $34^9 m^1$, $m \ge 34$, use Theorem~\ref{thm:(u+1)-GDD q^(u+1) + h^u d^1} with $h = 2$, $u = 9$ and $q = 17$.

For type $34^9 m^1$, $22 \le m \le 31$,
use Theorem~\ref{thm:4-GDD construction 2-5-8}(ii) with $q = 11$, $c = 2$, $u = 9$, and
4-GDDs of types $8^2 2^8$ and $8^3 2^7$ from Lemma~\ref{lem:4-GDD general}.

For type $34^u m^1$, where \\
\adfind $u = 6$ and $m \in \{7, 10, 13, 16, 19, 22\} \cup \{73, 76, 79, 82\}$, or \\
\adfind $u = 9$ and $m = 19$, or \\
\adfind $u = 12$ and $m \in \{25, 28, 31\}$, or \\
\adfind $u = 15$ and $m = 31$,
see Lemma~\ref{lem:4-GDD 34^u m^1}.

\adfTgap
\noindent {\bf\boldmath $g = 38$~}
For type $38^u m^1$, $u \in \{9, 12, 15, 18\}$, $m \le u - 1$,
take a 4-DGDD of type $(38, 2^{19})^u$, adjoin $m$ points and overlay the holes plus the new points with 4-GDDs of type $2^u m^1$.

For type $38^6 m^1$, $29 \le m \le 80$, we use the construction of $g = 34$.
Take a 7-GDD of type $8^7$,
select a point, $A$, and inflate it by a factor of $a_0 \in \{1, 4, 7, 10\}$.
Inflate all other points in one block containing $A$ by 10,
inflate all other points in the same group as $A$ by $a_j \in \{4,10\}$, $j$ = 1, 2, \dots, 7, and
inflate the remaining points by 4.
Overlay the inflated blocks with 4-GDDs of types $4^6 a_0^1$, $10^6 a_0^1$, $4^7$, $4^6 10^1$ and $4^5 10^2$,
the last one from Lemma~\ref{lem:4-GDD general}.

For type $38^9 m^1$, $m \ge 38$,
use Theorem~\ref{thm:(u+1)-GDD q^(u+1) + h^u d^1} with $h = 2$, $u = 9$ and $q = 19$.

For type $38^u m^1$, $u \in \{9, 12, 15\}$, $m \in \{32, 35\}$,
use Theorem~\ref{thm:4-GDD construction 2-5-8}(i) with $q = 16$ and $c = 1$.

For type $38^u m^1$, $u \in \{9, 12\}$, $m \in \{26, 29\}$,
use Theorem~\ref{thm:4-GDD construction 2-5-8}(i) with $q = 13$, $c = 2$, and
4-GDDs of type $8^2 2^8$ and $8^2 2^{11}$, from Lemma~\ref{lem:4-GDD general}.

For type $38^u m^1$, where \\
\adfind $u = 6$  and $m \in \{5, 8, \dots, 26\} \cup \{83, 86, 89, 92\}$, or \\
\adfind $u = 9$  and $m \in \{11, 14, 17, 20, 23\}$, or \\
\adfind $u = 12$ and $m \in \{14, 17, 20, 23\}$, or \\
\adfind $u = 15$ and $m \in \{17, 20, 23, 26, 29\}$, or \\
\adfind $u = 18$ and $m \in \{20, 23, 26, 29, 32, 35\}$,
see Lemma~\ref{lem:4-GDD 38^u m^1}.

\adfTgap
\noindent {\bf\boldmath $g = 40$~}
For all except types $40^6 m^1$ and $40^9 m^1$, see \cite[Theorem 7.2]{Schuster2010}.

For type $40^u m^1$, $u \in \{6, 9\}$, $m \le 5(u - 1)$,
take a 4-DGDD of type $(40, 10^4)^u$, adjoin $m$ points and overlay the holes plus the new points with 4-GDDs of type $10^u m^1$.

For type $40^6 m^1$, $25\le m \le 70$:
Take a 7-GDD of type $7^{7}$,
choose a point, $A$, and inflate it by a factor of $a \in \{1, 4, 7, 10\}$.
Inflate the other points in the same group as $A$ by factors chosen from $\{4, 10\}$.
Choose two blocks containing $A$ and inflate their points other than $A$ by 10.
Inflate the remaining points by 4.
\begin{center}
\begin{tabular}{cccccc|c}
 4 &  4 &  4 &  4 &  4 &  4 & $\{4, 10\}$ \\
 \dots & \dots & \dots & \dots & \dots & \dots & \dots \\
 4 &  4 &  4 &  4 &  4 &  4 & $\{4, 10\}$ \\
10 & 10 & 10 & 10 & 10 & 10 & $\{4, 10\}$ \\
10 & 10 & 10 & 10 & 10 & 10 & $a \in \{1, 4, 7, 10\}$
\end{tabular}
\end{center}
Overlay the inflated blocks with 4-GDDs of types
$4^6 d^1$, $10^6 d^1$, $d \in \{1, 4, 7, 10\}$, $4^4 10^3$ and $4^5 10^2$, the last two from Lemma~\ref{lem:4-GDD general}.

For type $40^9 m^1$, $34 \le m \le 136$,
use Theorem~\ref{thm:4-GDD construction 2-5-8}(ii) with $q = 17$, $c = 1$, and
a 4-GDD of type $8^2 2^8$ from Lemma~\ref{lem:4-GDD general}.

For type $40^u m^1$, $(u, m) \in \{(6, 76), (6, 88), (9, 148)\}$,
use Theorem~\ref{thm:4-GDD (rg)^u (rm)^1} with a 4-GDD of type $10^u (m/4)^1$.

For type $40^u m^1$, $(u, m) \in \{(6, 85), (9, 145)\}$,
use Theorem~\ref{thm:4-GDD (rg)^u (rm)^1} with a 4-GDD of type $8^u (m/5)^1$.

For type $40^u m^1$, where \\
\adfind $u = 6$ and $m \in \{73, 79, 82, 91, 94, 97\}$, or \\
\adfind $u = 9$ and $m \in \{139, 142, 151, 154, 157\}$, see Lemma~\ref{lem:4-GDD 40^u m^1}.

\adfTgap
\noindent {\bf\boldmath $g = 44$~}
For type $44^6 m^1$, $m \in \{5, 8\}$, see \cite[Lemma 2.8]{ForbesForbes2018}.

For type $44^u m^1$, $u \in \{9, 15,18\}$, $m \le u - 1$,
take a 4-DGDD of type $(44, 2^{22})^u$, adjoin $m$ points and overlay the holes plus the new points with 4-GDDs of type $2^u m^1$.

For type $44^u m^1$, $u \in \{6, 9\}$, $m \ge 11$,
use Theorem~\ref{thm:(u+1)-GDD q^(u+1) + h^u d^1} with $h = 4$ and $q = 11$.

For type $44^{12} m^1$, $m < 44$,
take a 4-DGDD of type $(44, 11^{4})^{12}$, adjoin $m$ points and overlay the holes plus the new points with 4-GDDs of type $11^{12} m^1$.

For type $44^{15} 17^1$,
take a 4-DGDD of type $(44, 11^{4})^{15}$, adjoin $17$ points and overlay the holes plus the new points with 4-GDDs of type $11^{15} 17^1$.

For type $44^u m^1$, $u \in \{15, 18\}$, $20 \le m \le 44$:
Take a 10-RGDD of type $(u + 1)^{10}$ and remove a parallel class to create a 10-DGDD of type $(u + 1, 1^{u + 1})^{10}$.
Select four points, $A_1, A_2, A_3, A_4$, all in the same hole, and
inflate $A_i$ by a factor of $a_i \in \{2, 5, 8\}$, $i = 1$, 2, 3, 4.
Inflate all other points in the four groups containing the points $A_i$ by 8.
Inflate the remaining points by 2.
\begin{center}
\begin{tabular}{cccccccccc}
2 & 2 & 2 & 2 & 2 & 2 & 8 & 8 & 8 & 8 \\
\dots & \dots & \dots & \dots & \dots & \dots & \dots & \dots & \dots & \dots \\
2 & 2 & 2 & 2 & 2 & 2 & 8 & 8 & 8 & 8 \\
2 & 2 & 2 & 2 & 2 & 2 & $a_1$ & $a_2$ & $a_3$ & $a_4$
\end{tabular}
\end{center}
The hole containing the inflated points $A_i$ has size $m = 12 + a_1 + a_2 + a_3 + a_4$, and the other $u$ holes have size 44 each.
Overlay the inflated blocks with 4-GDDs of types $2^6 8^4$, $2^7 8^3$ and $2^6 8^3 5^1$ from Lemma~\ref{lem:4-GDD general}.
Overlay the inflated groups with 4-GDDs of types $2^{u + 1}$, $8^u 2^1$, $8^u 5^1$ and $8^{u + 1}$.

\adfTgap
\noindent {\bf\boldmath $g = 46$~}
For type $46^u m^1$, $u \in \{6, 9, 12, 15, 18\}$, $m \le 2(u - 1)$,
use Theorem~\ref{thm:4-GDD (gi)^(ti) + (hgi)^u m^1} with $h = 2$ and 4-GDDs of types $2^9 5^1$, $4^u m^1$ and $10^u m^1$.

For type $46^6 m^1$, $25 \le m \le 70$:
Take a 7-GDD of type $7^{7}$,
choose a point, $A$, and inflate it by a factor of $a \in \{1, 4, 7, 10\}$.
Inflate the other points in the same group as $A$ by factors chosen from $\{4, 10\}$.
Choose three blocks containing $A$ and inflate their points other than $A$ by 10.
Inflate the remaining points by 4.
Observe that blocks containing $A$ have their other points either all inflated by 10 or all inflated by 4.
Observe also that blocks not containing $A$ have $r$ points inflated by 10 and $7-r$ inflated by 4, $r = 2$, 3, 4.
The diagram shows the groups as columns and
the parallel class defined by $A$ as rows on the left of the vertical line.
\begin{center}
\begin{tabular}{cccccc|c}
  4 &  4 &  4 &  4 &  4 &  4 & $\{4, 10\}$ \\
 \dots & \dots & \dots & \dots & \dots & \dots & \dots \\
  4 &  4 &  4 &  4 &  4 &  4 & $\{4, 10\}$ \\
 10 & 10 & 10 & 10 & 10 & 10 & $\{4, 10\}$ \\
 10 & 10 & 10 & 10 & 10 & 10 & $\{4, 10\}$ \\
 10 & 10 & 10 & 10 & 10 & 10 & $a \in \{1, 4, 7, 10\}$
\end{tabular}
\end{center}
Overlay the inflated blocks with 4-GDDs of types
$4^6 a^1$, $10^6 a^1$, $a \in \{1, 4, 7, 10\}$, $4^3 10^4$, $4^4 10^3$ and $4^5 10^2$,
the last three from Lemma~\ref{lem:4-GDD general}.

For type $46^9 m^1$, $m \ge 46$, use Theorem~\ref{thm:(u+1)-GDD q^(u+1) + h^u d^1} with $h = 2$, $u = 9$ and $q = 23$.

For type $46^9 m^1$, $22 \le m \le 43$,
use Theorem~\ref{thm:4-GDD construction 2-5-8}(ii) with $q = 11$, $c = 4$, and
4-GDDs of types $8^2 2^8$, $8^3 2^7$, $8^4 2^6$ and $8^5 2^5$ from Lemma~\ref{lem:4-GDD general}.

For type $46^u m^1$, $u \in \{12, 15, 18\}$, $m \in \{37, 40, 43\}$:
Remove a parallel class from a 10-RGDD of type $(u + 1)^{10}$ to create a 10-DGDD of type $(u + 1, 1^{u + 1})^{10}$.
Select a point $A$ and inflate it by a factor of $m - 36 \in \{1, 4, 7\}$.
Inflate all other points in the same group as $A$ by 10, and inflate the remaining points by 4.
\begin{center}
\begin{tabular}{cccccccccc}
4 & 4 & 4 & 4 & 4 & 4 & 4 & 4 & 4 & 10 \\
\dots & \dots & \dots & \dots & \dots & \dots & \dots & \dots & \dots & \dots \\
4 & 4 & 4 & 4 & 4 & 4 & 4 & 4 & 4 & 10 \\
4 & 4 & 4 & 4 & 4 & 4 & 4 & 4 & 4 & $m - 36$
\end{tabular}
\end{center}
Overlay the inflated blocks with 4-GDDs of types $4^9 (m - 36)^1$ and $4^9 10^1$.
Overlay the inflated groups with 4-GDDs of types $4^{u + 1}$ and $10^u (m - 36)^1$.

For type $46^{15} 34^1$,
Take a 16-GDD of type $17^{16}$.
Inflate $15$ points in one block by 14 and inflate all other points by 2.
\begin{center}
\begin{tabular}{c@{~}c@{~}c@{~}c@{~}c@{~}c@{~}c@{~}c@{~}c@{~}c@{~}c@{~}c@{~}c@{~}c@{~}c@{~}c}
2 & 2 & 2 & 2 & 2 & 2 & 2 & 2 & 2 & 2 & 2 & 2 & 2 & 2 & 2 & 2  \\
\dots & \dots & \dots & \dots & \dots & \dots & \dots & \dots & \dots & \dots & \dots & \dots & \dots & \dots & \dots & \dots \\
2 & 2 & 2 & 2 & 2 & 2 & 2 & 2 & 2 & 2 & 2 & 2 & 2 & 2 & 2 & 2  \\
14 & 14 & 14 & 14 & 14 & 14 & 14 & 14 & 14 & 14 & 14 & 14 & 14 & 14 & 14 & 2
\end{tabular}
\end{center}
Overlay the inflated blocks with 4-GDDs of types $2^{16}$, $2^{15} 14^1$ and $14^{15} 2^1$.

For type $46^u m^1$, where \\
\adfind $u = 6$ and $m \in \{13, 16, 19, 22\} \cup \{73, 76, \dots, 112\}$, or \\
\adfind $u = 9$ and $m = 19$, or \\
\adfind $u = 12$ and $m \in \{25, 28, 31, 34\}$, or \\
\adfind $u = 15$ and $m = 31$, see Lemma~\ref{lem:4-GDD 46^u m^1}.

\adfTgap
\noindent {\bf\boldmath $g = 50$~}
For type $50^u m^1$, $u \in \{9, 12, 15, 18\}$, $m \le u - 1$,
take a 4-DGDD of type $(50, 2^{25})^u$, adjoin $m$ points and overlay the holes plus the new points with 4-GDDs of type $2^u m^1$.

For type $50^6 m^1$, $41 \le m \le 110$:
Take a 7-GDD of type $11^7$,
select a point, $A$, and, as for $g = 34$, inflate it by a factor of $a_0 \in \{1, 4, 7, 10\}$.
Inflate all other points in one block containing $A$ by 10,
inflate all other points in the same group as $A$ by $a_j \in \{4,10\}$, $j$ = 1, 2, \dots, 10, and
inflate the remaining points by 4.
Overlay the inflated blocks with 4-GDDs of types $4^6 a_0^1$, $10^6 a_0^1$, $4^7$, $4^6 10^1$ and $4^5 10^2$,
the last one from Lemma~\ref{lem:4-GDD general}.

For type $50^6 m^1$, $m \in \{5, 20, 35\}$,
use Theorem~\ref{thm:4-GDD (rg)^u (rm)^1} with a 4-GDD of type $10^{6} (m/5)^1$.

For type $50^9 m^1$, $m \ge 50$,
use Theorem~\ref{thm:(u+1)-GDD q^(u+1) + h^u d^1} with $h = 2$, $u = 9$ and $q = 25$.

For type $50^9 m^1$, $26 \le m \le 47$,
use Theorem~\ref{thm:4-GDD construction 2-5-8}(ii) with $q = 13$, $c = 4$, and
4-GDDs of types $8^2 2^8$, $8^3 2^7$, $8^4 2^6$ and $8^5 2^5$ from Lemma~\ref{lem:4-GDD general}.

For type $50^{u} m^1$, $(u, m) \in \{(9, 14), (9, 20), (12, 14), (12, 17), (15, 17)\}$,
take a 4-DGDD of type $(50, 5^{10})^{u}$, adjoin $m$ points and overlay the holes plus the new points with 4-GDDs of type $5^{u} m^1$.

For type $50^u m^1$, $u \in \{12, 15, 18\}$, $ 20 \le m \le 53$:
Remove a parallel class from a 10-RGDD of type $(u + 1)^{10}$ to create a 10-DGDD of type $(u + 1, 1^{u + 1})^{10}$.
Select six points, $A_0$, $A_1$, \dots, $A_5$, all in the same hole.
Inflate $A_0$ by a factor of $a_0 \in \{2, 5\}$ and
inflate $A_j$ by $a_j \in \{2, 8\}$, $j = 1$, 2, \dots, 5.
Inflate all other points in the five groups containing $A_j$, $j = 1$, 2, \dots, 5, by 8 and
inflate the remaining points by 2.
\begin{center}
\begin{tabular}{cccccccccc}
2 & 2 & 2 & 2 & 2 & 8 & 8 & 8 & 8 & 8 \\
\dots & \dots & \dots & \dots & \dots & \dots & \dots & \dots & \dots & \dots \\
2 & 2 & 2 & 2 & 2 & 8 & 8 & 8 & 8 & 8 \\
2 & 2 & 2 & 2 & $a_0$ & $a_1$ & $a_2$ & $a_3$ & $a_4$ & $a_5$
\end{tabular}
\end{center}
The hole containing the inflated points $A_i$ has size $m = 8 + a_0 + a_1 + \dots + a_5$, and the other $u$ holes have size 50 each.
Overlay the inflated blocks with 4-GDDs of types
$2^6 8^4$, $2^5 8^5$ and $2^4 8^5 5^1$ from Lemma~\ref{lem:4-GDD general}.
Overlay the inflated groups with 4-GDDs of types
$2^{u + 1}$, $2^u 5^1$, $8^u 2^1$ and $8^{u + 1}$.

For type $50^u m^1$, where \\
\adfind $u = 6$ and $m \in \{8, 11, 14, 17, 23, 26, 29, 32, 38\} \cup \{113, 116, 119, 122\}$, or \\
\adfind $u = 9$ and $m \in \{11, 17, 23\}$,
see Lemma~\ref{lem:4-GDD 50^u m^1}.

\adfTgap
\noindent {\bf\boldmath $g = 52$~}
For type $52^u m^1$, $u \in \{6, 9, 12, 15, 18\}$, $m \le 2(u - 1)$,
take a 4-DGDD of type $(52, 4^{13})^u$, adjoin $m$ points and overlay the holes plus the new points with 4-GDDs of type $4^u m^1$.

For type $52^u m^1$, $u \in \{6, 9, 12\}$, $m \ge 13$,
use Theorem~\ref{thm:(u+1)-GDD q^(u+1) + h^u d^1} with $h = 4$ and $q = 13$.

For type $52^u m^1$, $u \in \{15, 18\}$, $m \in \{37, 40, 43, 46, 49\}$, we proceed as for $g = 46$.
Remove a parallel class from a 10-RGDD of type $(u + 1)^{10}$ to create a 10-DGDD, of type $(u + 1, 1^{u + 1})^{10}$.
Select a point, $A$, and inflate it by a factor of $m - 36 \in \{1, 4, 7, 10, 13\}$.
Inflate all other points in the same group as $A$ by 16, and inflate the remaining points by 4.
\begin{center}
\begin{tabular}{cccccccccc}
4 & 4 & 4 & 4 & 4 & 4 & 4 & 4 & 4 & 16 \\
\dots & \dots & \dots & \dots & \dots & \dots & \dots & \dots & \dots & \dots \\
4 & 4 & 4 & 4 & 4 & 4 & 4 & 4 & 4 & 16 \\
4 & 4 & 4 & 4 & 4 & 4 & 4 & 4 & 4 & $m - 36$
\end{tabular}
\end{center}
Overlay the inflated blocks and inflated groups with 4-GDDs of types
$4^9 (m - 36)^1$, $4^9 16^1$, $4^{u + 1}$ and $16^u (m - 36)^1$.

For type $52^{15} 31^1$,
take a 4-DGDD of type $(52, 13^{4})^{15}$, adjoin $31$ points and overlay the holes plus the new points with 4-GDDs of type $13^{15} 31^1$.

For type $52^{15} 34^1$,
see Lemma~\ref{lem:4-GDD 52^u m^1}.

\adfTgap
\noindent {\bf\boldmath $g = 56$~}
For type $56^u m^1$, $u \in \{6, 9, 12, 15, 18\}$, $m \le 7(u - 1)$,
take a 4-DGDD of type $(56, 14^{4})^u$, adjoin $m$ points and overlay the holes plus the new points with 4-GDDs of type $14^u m^1$.

For type $56^6 m^1$, $m \ge 14$,
use Theorem~\ref{thm:(u+1)-GDD q^(u+1) + h^u d^1} with $h = 8$ and $q = 7$.

For type $56^9 m^1$, $50 \le m \le 200$,
use Theorem~\ref{thm:4-GDD construction 2-5-8}(ii) with $q = 25$, $c = 1$, and
a 4-GDD of type $8^2 2^8$ from Lemma~\ref{lem:4-GDD general}.

For type $56^9 212^1$,
use Theorem~\ref{thm:4-GDD (rg)^u (rm)^1} with a 4-GDD of type $14^{9} 53^1$.

For type $56^9 203^1$,
use Theorem~\ref{thm:4-GDD (rg)^u (rm)^1} with a 4-GDD of type $8^{9} 29^1$.

The remaining types $56^9 m^1$, $m \in \{206, 209, 215, 218, 221\}$, are possible exceptions.

\adfTgap
\noindent {\bf\boldmath $g = 58$~}
For type $58^u m^1$, $u \in \{6, 9, 12, 15, 18\}$, $m \le 2(u - 1)$,
use Theorem~\ref{thm:4-GDD (gi)^(ti) + (hgi)^u m^1} with $h = 2$ and 4-GDDs of types $2^{12} 5^1$, $4^u m^1$ and $10^u m^1$.

For type $58^u m^1$, $u \in \{6, 12, 15, 18\}$, $13 \le m \le 64$:
Remove a parallel class from a 7-RGDD of type $(u + 1)^{7}$ to create a 7-DGDD of type $(u + 1, 1^{u + 1})^{7}$.
Select a hole and select two of its points, $B_1$ and $B_2$, say.
Inflate $B_1$ and $B_2$ by factors of $b_1$ and $b_2$ respectively, chosen from $\{4, 10\}$.
Inflate the other five points in the same hole by $a_i$, chosen from $\{1, 4, 10\}$, $i = 1$, 2, \dots, 5.
Inflate all other points in the same groups as $B_1$ and $B_2$ by 4.
Inflate the remaining points by 10.
\begin{center}
\begin{tabular}{ccccccc}
10 & 10 & 10 & 10 & 10 & 4 & 4 \\
\dots & \dots & \dots & \dots & \dots & \dots & \dots \\
10 & 10 & 10 & 10 & 10 & 4 & 4 \\
$a_1$ & $a_2$ & $a_3$ & $a_4$ & $a_5$ & $b_1$ & $b_2$
\end{tabular}
\end{center}
Overlay the inflated blocks with 4-GDDs of types
$10^6 4^1$, $10^5 4^2$, $10^4 4^2 1^1$ and $10^4 4^3$, the last three from Lemma~\ref{lem:4-GDD general}.
Overlay the inflated groups with 4-GDDs of types
$4^{u + 1}$, $4^u 10^1$, $10^u 1^1$, $10^u 4^1$ and $10^{u + 1}$.

For type $58^6 m^1$, $49 \le m \le 130$:
Take a 7-GDD of type $13^{7}$.
Select a point, $A$, and inflate it by $a_0 \in \{1, 4, 7, 10\}$.
Select a block containing $A$ and inflate its points other than $A$ by 10.
Inflate all other points in the same group as $A$ by $a_j \in \{4, 10\}$, $j = 1$, 2, \dots, 12.
Inflate the remaining points by 4.
\begin{center}
\begin{tabular}{ccccccl}
4 & 4 & 4 & 4 & 4 & 4 & $a_{12} \in \{4, 10\}$ \\
\dots & \dots & \dots & \dots & \dots & \dots & \dots \\
4 & 4 & 4 & 4 & 4 & 4 & $a_{1} \in \{4, 10\}$ \\
10 & 10 & 10 & 10 & 10 & 10 & $a_0 \in \{1, 4, 7, 10\}$
\end{tabular}
\end{center}
Overlay the inflated blocks and groups with 4-GDDs of types
$f^6 d^1$, $f \in \{4, 10\}$, $d \in \{1, 4, 7, 10\}$, and $4^5 10^2$ from Lemma~\ref{lem:4-GDD general}.

For type $58^9 m^1$, $m \ge 58$, use Theorem~\ref{thm:(u+1)-GDD q^(u+1) + h^u d^1} with $h = 2$, $u = 9$ and $q = 29$.

For type $58^9 m^1$, $22 \le m \le 88$,
use Theorem~\ref{thm:4-GDD construction 2-5-8}(ii) with $q = 11$, $c = 6$, and
4-GDDs of types $8^4 2^6$, $8^5 2^5$, $8^6 2^4$ and $8^7 2^3$ from Lemma~\ref{lem:4-GDD general}.

For type $58^u m^1$, where \\
\adfind $u = 6$ and $m \in \{133, 136, 139, 142\}$, or \\
\adfind $u = 9$ and $m = 19$, see Lemma~\ref{lem:4-GDD 58^u m^1}.

\adfTgap
\noindent {\bf\boldmath $g = 62$~}
For type $62^u m^1$, $u \in \{6, 9, 12, 15, 18\}$, $m \le 4(u - 1)$,
use Theorem~\ref{thm:4-GDD (gi)^(ti) + (hgi)^u m^1} with $h = 2$ and 4-GDDs of types $4^{6} 7^1$, $8^u m^1$ and $14^u m^1$.

For type $62^6 m^1$, $50 \le m \le 140$:
Take a 7-GDD of type $7^{7}$,
select a point, $A$, and inflate it by a factor of $a \in \{2, 5, \dots, 20\}$.
Choose a block containing $A$ and inflate the other six points by 14.
Inflate all other points in the same group as $A$ by values chosen from $\{8, 20\}$.
Inflate the remaining points by 8.
\begin{center}
\begin{tabular}{ccccccc}
8 & 8 & 8 & 8 & 8 & 8 & $\{8, 20\}$ \\
\dots & \dots & \dots & \dots & \dots & \dots & \dots \\
8 & 8 & 8 & 8 & 8 & 8 & $\{8, 20\}$ \\
14 & 14 & 14 & 14 & 14 & 14 & $\{2, 5, \dots, 20\}$
\end{tabular}
\end{center}
Overlay the inflated blocks with 4-GDDs of types
$8^6 a^1$, $14^6 a^1$, $a \in \{2, 5, \dots, 20\}$, and $8^5 14^1 20^1$ from Lemma~\ref{lem:4-GDD general}.

For type $62^9 m^1$, $m \ge 62$, use Theorem~\ref{thm:(u+1)-GDD q^(u+1) + h^u d^1} with $h = 2$, $u = 9$ and $q = 31$.

For type $62^9 m^1$, $26 \le m \le 59$,
use Theorem~\ref{thm:4-GDD construction 2-5-8}(ii) with $q = 13$, $c = 6$, and
4-GDDs of types $8^2 2^8$, $8^3 2^7$, $8^4 2^6$, $8^5 2^5$, $8^6 2^4$ and $8^7 2^3$ from Lemma~\ref{lem:4-GDD general}.

For type $62^{u} m^1$, $u \in \{12, 15\}$, $20 \le m \le 65$:
Remove a parallel class from a 10-RGDD of type $(u + 1)^{10}$ to create a 10-DGDD of type $(u + 1, 1^{u + 1})^{10}$.
Choose a hole, \{$A_1$, $A_2$, \dots, $A_{10}$\},
inflate $A_i$ by $a_i \in \{2, 8\}$, $i = 1$, 2, \dots, 7, and
inflate $A_8$ by $a_8 \in \{2, 5\}$.
Inflate all other points in the groups containing $A_8$, $A_9$ and $A_{10}$ by 2.
Inflate the remaining points by 8.
\begin{center}
\begin{tabular}{cccccccccc}
8 & 8 & 8 & 8 & 8 & 8 & 8 & 2 & 2 & 2 \\
\dots & \dots & \dots & \dots & \dots & \dots & \dots & \dots & \dots & \dots \\
8 & 8 & 8 & 8 & 8 & 8 & 8 & 2 & 2 & 2 \\
$\{2,8\}$ & $\{2,8\}$ & $\{2,8\}$ & $\{2,8\}$ & $\{2,8\}$ & $\{2,8\}$ & $\{2,8\}$ & $\{2,5\}$ & 2 & 2
\end{tabular}
\end{center}
Overlay the inflated groups with 4-GDDs of types
$8^{u} 2^1$, $8^{u + 1}$, $2^{u} 5^1$ and $2^{u+1}$.
Overlay the inflated blocks with 4-GDDs of types
$8^6 2^4$, $8^7 2^3$ and $8^7 2^2 5^1$ from Lemma~\ref{lem:4-GDD general}.

For type $62^u m^1$, where
$m \in \{23, 26, \dots, 47\} \cup \{143, 146, 149, 152\}$, see Lemma~\ref{lem:4-GDD 62^u m^1}.

\adfTgap
\noindent {\bf\boldmath $g = 68$~}
For type $68^u m^1$, $u \in \{6, 9, 12, 15, 18\}$, $m \le 4(u - 1)$,
use Theorem~\ref{thm:4-GDD (gi)^(ti) + (hgi)^u m^1} with $h = 2$ and 4-GDDs of types $4^6 10^1$, $8^u m^1$ and $20^u m^1$.

For type $68^u m^1$, $u \in \{6, 9, 12, 15\}$, $m \ge 17$,
use Theorem~\ref{thm:(u+1)-GDD q^(u+1) + h^u d^1} with $h = 4$ and $q = 17$.

\adfTgap
\noindent {\bf\boldmath $g = 76$~}
For type $76^6 m^1$, $u \in \{6, 9, 12, 15, 18\}$, $m \le 2(u - 1)$,
take a 4-DGDD of type $(76, 4^{19})^u$, adjoin $m$ points and overlay the holes plus the new points with 4-GDDs of type $4^u m^1$.

For type $76^u m^1$, $u \in \{6, 9, 12, 15, 18\}$, $m \ge 19$,
use Theorem~\ref{thm:(u+1)-GDD q^(u+1) + h^u d^1} with $h = 4$ and $q = 19$.

For types $76^6 13^1$ and $76^6 16^1$, see Lemma~\ref{lem:4-GDD 76^u m^1}.

\adfTgap
\noindent {\bf\boldmath $g = 80$~}
For type $80^u m^1$, $u \in \{6, 9, 12, 15, 18\}$, $m \le 10(u - 1)$,
take a 4-DGDD of type $(80, 20^{4})^u$, adjoin $m$ points and overlay the holes plus the new points with 4-GDDs of type $20^u m^1$.

For type $80^6 m^1$, $m \ge 8$,
use Theorem~\ref{thm:(u+1)-GDD q^(u+1) + h^u d^1} with $h = 10$ and $q = 8$.

For type $80^9 m^1$, $74 \le m \le 296$,
use Theorem~\ref{thm:4-GDD construction 2-5-8}(ii) with $q = 37$, $c = 1$, and
a 4-GDD of type $8^2 2^8$ from Lemma~\ref{lem:4-GDD general}.

For type $80^9 m^1$, $m \in \{302, 308, 314\}$,
use Theorem~\ref{thm:(u+1)-GDD q^(u+1) + h^u d^1 d in D} with $h = 5$, $u = 9$, $q = 16$ and $D = \{14, 20\}$.

For type $80^9 305^1$,
use Theorem~\ref{thm:4-GDD (rg)^u (rm)^1} with a 4-GDD of type $16^{9} 61^1$.

The remaining types $80^9 m^1$, $m \in \{299, 311, 317\}$, are possible exceptions.

\adfTgap
\noindent {\bf\boldmath $g = 88$~}
Use Theorem~\ref{thm:4-GDD h^u d^1 + (g/4)^u m^1 -> g^u m^1} with $h = 8$ and $q = 11$.

\adfTgap
\noindent {\bf\boldmath $g = 92$~}
For type $92^u m^1$, $u \in \{6, 9, 12, 15, 18\}$, $m \le 20$,
use Theorem~\ref{thm:4-GDD (gi)^(ti) + (hgi)^u m^1} with $h = 2$ and 4-GDDs of types $4^9 10^1$, $8^u m^1$ and $20^u m^1$.

For type $92^u m^1$, $u \in \{6, 9, 12, 15, 18\}$, $m \ge 23$, use Theorem~\ref{thm:(u+1)-GDD q^(u+1) + h^u d^1} with $h = 4$ and $q = 23$.

\adfTgap
\noindent {\bf\boldmath $g = 100$~}
For type $100^u m^1$, $u \in \{6, 9, 12, 15, 18\}$, $m \le 22$,
take a 4-DGDD of type $(100, 10^{10})^u$, adjoin $m$ points and overlay the holes plus the new points with 4-GDDs of type $10^u m^1$.

For type $100^u m^1$, $u \in \{6, 9, 12, 15, 18\}$, $m \ge 25$,
use Theorem~\ref{thm:(u+1)-GDD q^(u+1) + h^u d^1} with $h = 4$ and $q = 25$.

\adfTgap
\noindent {\bf\boldmath $g = 104$~}
Use Theorem~\ref{thm:4-GDD h^u d^1 + (g/4)^u m^1 -> g^u m^1} with $h = 8$ and $q = 13$.

\adfTgap
\noindent {\bf\boldmath $g = 112$~}
For type $112^u m^1$, $u \in \{6, 9, 12, 15, 18\}$, $m \le 14(u - 1)$,
take a 4-DGDD of type $(112, 28^{4})^u$, adjoin $m$ points and overlay the holes plus the new points with 4-GDDs of type $28^u m^1$.

For type $112^6 m^1$, $m \ge 7$,
use Theorem~\ref{thm:(u+1)-GDD q^(u+1) + h^u d^1} with $h = 16$ and $q = 7$.

For type $112^9 m^1$, $106 \le m \le 424$,
use Theorem~\ref{thm:4-GDD construction 2-5-8}(ii) with $q = 53$, $c = 1$, and
a 4-GDD of type $8^2 2^8$ from Lemma~\ref{lem:4-GDD general}.

For type $112^9 m^1$, $m \in \{430, 436, 442\}$,
use Theorem~\ref{thm:(u+1)-GDD q^(u+1) + h^u d^1 d in D} with $h = 7$, $ u = 9$, $q = 16$ and $D = \{22, 28\}$.

For type $112^9 427$,
use Theorem~\ref{thm:4-GDD (rg)^u (rm)^1} with a 4-GDD of type $16^{9} 61^1$.

The remaining types $112^9 m^1$, $m \in \{433, 439, 445\}$, are possible exceptions.

\adfTgap
\noindent {\bf\boldmath $g = 116$~}
For type $116^u m^1$, $u \in \{6, 9, 12, 15, 18\}$, $m \le 4(u - 1)$,
use Theorem~\ref{thm:4-GDD (gi)^(ti) + (hgi)^u m^1} with $h = 2$ and 4-GDDs of types $4^{12} 10^1$, $8^u m^1$ and $20^u m^1$.

For type $116^u m^1$, $u \in \{6, 9, 12, 15, 18\}$, $m \ge 29$, use Theorem~\ref{thm:(u+1)-GDD q^(u+1) + h^u d^1} with $h = 4$ and $q = 29$.

For type $116^6 m^1$, $m \in \{23, 26\}$, the two cases not covered,
the following construction turned out to be easier than trying to build $116^6 23^1$ directly---and
it also gives $116^6 26^1$.
Take a 7-RGDD of type $7^{7}$ and remove a parallel class to create a 7-DGDD of type $(7, 1^7)^7$.
Select a hole and denote its points by $A_1$, $A_2$, \dots, $A_7$.
Inflate $A_1$, $A_2$, $A_3$ by a factor of 5, inflate $A_4$ by 2 or 5, and inflate $A_5$, $A_6$, $A_7$ by 2.
Inflate all other points in the same groups as $A_1$ and $A_2$ by 8, and
inflate the remaining points by 20.
\begin{center}
\begin{tabular}{ccccccc}
 8 &  8 & 20 & 20 & 20 & 20 & 20 \\
\dots & \dots & \dots & \dots & \dots & \dots & \dots \\
 8 &  8 & 20 & 20 & 20 & 20 & 20 \\
 5 &  5 &  5 & \{2,5\} & 2 & 2 & 2
\end{tabular}
\end{center}
Overlay the inflated groups and blocks with 4-GDDs of types
$8^6 5^1$, $20^6 2^1$, $20^6 5^1$,
$20^4 8^2 2^1$, $20^4 8^2 5^1$ and $20^5 8^1 5^1$. 
See Lemma~\ref{lem:4-GDD general} for the last three.

\adfTgap
\noindent {\bf\boldmath $g = 124$~}
For type $124^u m^1$, $u \in \{6, 9, 12, 15, 18\}$, $m \le 28$,
use Theorem~\ref{thm:4-GDD (gi)^(ti) + (hgi)^u m^1} with $h = 4$ and 4-GDDs of types $4^{6} 7^1$, $16^u m^1$ and $28^u m^1$.

For type $124^u m^1$, $u \in \{6, 9, 12, 15, 18\}$, $m \ge 31$,
use Theorem~\ref{thm:(u+1)-GDD q^(u+1) + h^u d^1} with $h = 4$ and $q = 31$.

\adfTgap
\noindent {\bf\boldmath $g = 136$~}
Use Theorem~\ref{thm:4-GDD h^u d^1 + (g/4)^u m^1 -> g^u m^1} with $h = 8$ and $q = 17$.

\adfTgap
\noindent {\bf\boldmath $g = 152$~}
Use Theorem~\ref{thm:4-GDD h^u d^1 + (g/4)^u m^1 -> g^u m^1} with $h = 8$ and $q = 19$.

\adfTgap
\noindent {\bf\boldmath $g = 160$~}
For type $160^u m^1$, $u \in \{6, 9, 12, 15, 18\}$, $m \le 20(u - 1)$,
take a 4-DGDD of type $(160, 40^{4})^u$, adjoin $m$ points and overlay the holes plus the new points with 4-GDDs of type $40^u m^1$.

For type $160^u m^1$, $u \in \{6, 9\}$, $m \ge 16$,
use Theorem~\ref{thm:(u+1)-GDD q^(u+1) + h^u d^1} with $h = 10$ and $q = 16$.

\adfTgap
\noindent {\bf\boldmath $g = 176$~}
Use Theorem~\ref{thm:4-GDD h^u d^1 + (g/4)^u m^1 -> g^u m^1} with $h = 16$ and $q = 11$.

\adfTgap
\noindent {\bf\boldmath $g = 184$~}
Use Theorem~\ref{thm:4-GDD h^u d^1 + (g/4)^u m^1 -> g^u m^1} with $h = 8$ and $q = 23$.

\adfTgap
\noindent {\bf\boldmath $g = 200$~}
Use Theorem~\ref{thm:4-GDD h^u d^1 + (g/4)^u m^1 -> g^u m^1} with $h = 8$ and $q = 25$.

\adfTgap
\noindent {\bf\boldmath $g = 208$~}
Use Theorem~\ref{thm:4-GDD h^u d^1 + (g/4)^u m^1 -> g^u m^1} with $h = 16$ and $q = 13$.

\adfTgap
\noindent {\bf\boldmath $g = 224$~}
For type $224^u m^1$, $u \in \{6, 9, 12, 15, 18\}$, $m \le 16(u - 1)$,
take a 4-DGDD of type $(224, 32^{7})^u$, adjoin $m$ points and overlay the holes plus the new points with 4-GDDs of type $32^u m^1$.

For type $224^u m^1$, $u \in \{6, 9, 12\}$, $m \ge 32$,
use Theorem~\ref{thm:(u+1)-GDD q^(u+1) + h^u d^1} with $h = 14$ and $q = 16$.

\adfTgap
\noindent {\bf\boldmath $g = 232$~}
Use Theorem~\ref{thm:4-GDD h^u d^1 + (g/4)^u m^1 -> g^u m^1} with $h = 8$ and $q = 29$.

\adfTgap
\noindent {\bf\boldmath $g = 248$~}
Use Theorem~\ref{thm:4-GDD h^u d^1 + (g/4)^u m^1 -> g^u m^1} with $h = 8$ and $q = 31$.

\adfTgap
\noindent {\bf\boldmath $g = 272$~}
Use Theorem~\ref{thm:4-GDD h^u d^1 + (g/4)^u m^1 -> g^u m^1} with $h = 16$ and $q = 17$.

\adfTgap
\noindent {\bf\boldmath $g = 304$~}
Use Theorem~\ref{thm:4-GDD h^u d^1 + (g/4)^u m^1 -> g^u m^1} with $h = 16$ and $q = 19$.

\adfTgap
\noindent {\bf\boldmath $g = 320$~}
For type $320^u m^1$, $u \in \{6, 9, 12, 15, 18\}$, $m \le 29$,
take a 4-DGDD of type $(320, 32^{10})^u$, adjoin $m$ points and overlay the holes plus the new points with 4-GDDs of type $32^u m^1$.

For type $320^u m^1$, $u \in \{6, 9, 12, 15, 18\}$, $m \ge 32$,
use Theorem~\ref{thm:(u+1)-GDD q^(u+1) + h^u d^1} with $h = 10$ and $q = 32$.

\adfTgap
\noindent {\bf\boldmath $g = 368$~}
Use Theorem~\ref{thm:4-GDD h^u d^1 + (g/4)^u m^1 -> g^u m^1} with $h = 16$ and $q = 23$.

\adfTgap
\noindent {\bf\boldmath $g = 400$~}
Use Theorem~\ref{thm:4-GDD h^u d^1 + (g/4)^u m^1 -> g^u m^1} with $h = 16$ and $q = 25$.

\adfTgap
\noindent {\bf\boldmath $g = 448$~}
For type $448^u m^1$, $u \in \{6, 9, 12, 15, 18\}$, $m \le 61$,
take a 4-DGDD of type $(448, 64^{7})^u$, adjoin $m$ points and overlay the holes plus the new points with 4-GDDs of type $64^u m^1$.

For type $448^u m^1$, $u \in \{6, 9, 12, 15, 18\}$, $m \ge 64$,
use Theorem~\ref{thm:(u+1)-GDD q^(u+1) + h^u d^1} with $h = 14$ and $q = 32$.

\adfTgap
\noindent {\bf\boldmath $g = 464$~}
Use Theorem~\ref{thm:4-GDD h^u d^1 + (g/4)^u m^1 -> g^u m^1} with $h = 16$ and $q = 29$.

\adfTgap
\noindent {\bf\boldmath $g = 496$~}
Use Theorem~\ref{thm:4-GDD h^u d^1 + (g/4)^u m^1 -> g^u m^1} with $h = 16$ and $q = 31$.
~\adfQED

\begin{theorem}
\label{thm:4-GDD (2^t q)^u m^1 q = 2 5 7 ...}
Let $q \in \{2, 5, 7, 11, 13, 17, 19, 23, 25, 29, 31\}$, let $t$ be a positive integer and let $g = 2^t q$.
Suppose $u \equiv 0 \adfmod{3}$, $u \ge 6$, $m \equiv g \adfmod{3}$ and $0 < m \le g(u - 1)/2$.
Then there exists a 4-GDD of type $g^u m^1$
with the possible exceptions stated in Theorem~\ref{thm:4-GDD g^u m^1 g = 14 20 22 ...}, namely
$56^9 m^1$, $m \in \{206, 209, 215, 218, 221\}$,
$80^9 m^1$, $m \in \{299, 311, 317\}$ and
$112^9 m^1$, $m \in \{433, 439, 445\}$.
\end{theorem}
\noindent {\bf Proof~}
The result follows from Theorem~\ref{thm:4-GDD g^u m^1 g = 14 20 22 ...} for $t \le 6$ if $q \in \{5, 7\}$ and for $t \le 4$ otherwise.
By Theorems~\ref{thm:4-GDD g^u m^1 existence} and \ref{thm:4-GDD g^u m^1 u = 21 27},
we may assume $u \in \{6, 9, 12, 15, 18\}$ and $m < g$ when $u \in \{12, 15, 18\}$.

Let $t \ge 7$ if $q \in \{5, 7\}$, $t \ge 5$ if $q \in \{2, 11, 13, 17, 19, 25, 29, 31\}$, and assume this theorem holds for smaller $t$.

Take a 4-DGDD of type $(2^t q, (2^{t - 2}q)^4)^u$, adjoin $m$ points and
overlay the holes plus the new points with 4-GDDs of type $(2^{t - 2}q)^u m^1$ to
create a 4-GDD of type $(2^t q)^u m^1$ for admissible $m \le 2^{t - 3} q (u - 1)$.
Observe that
$2^{t - 3}q(u - 1) = 5 g/8$ when $u = 6$ and
$2^{t - 3}q(u - 1) \ge g$ when $u \ge 9$.
Hence the theorem is proved for $u \in \{12, 15, 18\}$.
(Note that this construction does not use 4-GDDs of type $g^u m^1$ with $g \in \{56, 80, 112\}$.)

For $u \in \{6, 9\}$, we have $2^{t - 1} > u$. Hence, by Theorem~\ref{thm:(u+1)-GDD q^(u+1) + h^u d^1} with $(h, q) = (2q, 2^{t - 1})$,
there exists a 4-GDD of type $g^u m^1$ for all admissible $m \ge 2^{t}$.
Since $2^{t} = g/q \le g/2$ this range overlaps the interval for $m$ in the 4-DGDD construction and hence the theorem is proved.
~\adfQED

\adfPgap
Thus Theorems~\ref{thm:4-GDD g^u m^1 existence}, \ref{thm:4-GDD g^u m^1 g = 14 20 22 ...} and \ref{thm:4-GDD (2^t q)^u m^1 q = 2 5 7 ...}
assert the existence of 4-GDDs of type $g^u m^1$, $g = 2^t q$, $t \ge 0$ and prime power $q \le 32$,
for all admissible $u$ and $m$, except for
$$(g,u,m) \in \{(2,4,0), (2,3,2), (6,4,0), (6,3,6), (2,6,5)\}$$
and with a few possible exceptions when $g \in \{56, 80, 112\}$, $u = 9$ and $m \ge 4g - 21$.
To illustrate one way of extending Theorem~\ref{thm:4-GDD (2^t q)^u m^1 q = 2 5 7 ...}, we show that
there exists a 4-GDD of type $(2^t q^s)^u m^1$ for $s$, $t \ge 1$, $q \in \{19, 23, 25, 29, 31\}$ and all admissible $u$ and $m$.

\begin{theorem}
\label{thm:4-GDD (2^t q^s)^u m^1 q = 19 25 31}
Let $s$ and $t$ be positive integers, and let $g = 2^t q^s$, $q \in \{19, 25, 31\}$.
Suppose $u \equiv 0 \adfmod{3}$, $u \ge 6$, $m \equiv g \adfmod{3}$ and $0 < m \le g(u - 1)/2$.
Then there exists a 4-GDD of type $g^u m^1$.
\end{theorem}
\noindent {\bf Proof~}
For $s = 1$ and $t \ge 1$, the result follows from Theorem~\ref{thm:4-GDD (2^t q)^u m^1 q = 2 5 7 ...}.
Assume $s \ge 2$ and that this theorem holds for smaller $s$.
By Theorems~\ref{thm:4-GDD g^u m^1 existence} and \ref{thm:4-GDD g^u m^1 u = 21 27}, we may also assume
$u \in \{6, 9, 12, 15, 18\}$. Let $t \ge 1$.

Take a 4-DGDD of type $(2^t q^s, (2^t q^{s - 1})^{q})^u$, adjoin $m$ points and
overlay the holes plus the new points with 4-GDDs of type $(2^t q^{s - 1})^u m^1$ to
create a 4-GDD of type $(2^t q^s)^u m^1$ for admissible $m \le 2^{t - 1} q^{s - 1} (u - 1)$.

By Theorem~\ref{thm:(u+1)-GDD q^(u+1) + h^u d^1} with $h = 2^t q^{s - 1}$,
there exists a 4-GDD of type $(2^t q^s)^u m^1$ for admissible $m \ge 2q$.

Since $2q \le 2^{t - 1} q^{s - 1} (u - 1)$ the ranges of $m$ in the two constructions overlap and therefore the theorem is proved.
~\adfQED

\begin{theorem}
\label{thm:4-GDD (2^t q^s)^u m^1 q = 23 29}
Let $s$ and $t$ be positive integers, and let $g = 2^t q^s$, $q \in \{23, 29\}$.
Suppose $u \equiv 0 \adfmod{3}$, $u \ge 6$, $m \equiv g \adfmod{3}$ and $0 < m \le g(u - 1)/2$.
Then there exists a 4-GDD of type $g^u m^1$.
\end{theorem}
\noindent {\bf Proof~}
By Theorems~\ref{thm:4-GDD g^u m^1 existence}, \ref{thm:4-GDD g^u m^1 u = 21 27} and \ref{thm:4-GDD (2^t q)^u m^1 q = 2 5 7 ...},
we may assume $s \ge 2$ and $u \in$ \{6, 9, 12, 15, 18\}.
First, assume $s = 2$ and $t = 1$.

By Theorem~\ref{thm:(u+1)-GDD q^(u+1) + h^u d^1} with $h = 2 q$, since $u < q$ and $2q \equiv 1 \adfmod{3}$,
there exists a 4-GDD of type $(2q^2)^u m^1$ for admissible $m \ge q$.

For $(2 \cdot 23^2)^u m^1$ and $m < 23$, use Theorem~\ref{thm:4-GDD (gi)^(ti) + (hgi)^u m^1} with $h = 2$, $q = 23$ and
4-GDDs of types $10^{51} 19^1$, $20^u m^1$ and $38^u m^1$
from Theorems~\ref{thm:4-GDD g^u m^1 existence} and \ref{thm:4-GDD g^u m^1 g = 14 20 22 ...}.

For $(2 \cdot 29^2)^u m^1$ and $m < 29$, use Theorem~\ref{thm:4-GDD (gi)^(ti) + (hgi)^u m^1} with $h = 2$, $q = 29$ and
4-GDDs of types $10^{81} 31^1$, $20^u m^1$ and $62^u m^1$
from Theorems~\ref{thm:4-GDD g^u m^1 existence} and \ref{thm:4-GDD g^u m^1 g = 14 20 22 ...}.

Next, assume $s = 2$, $t \ge 2$ and that this theorem holds for $s = 2$ and smaller $t$.

By Theorem~\ref{thm:(u+1)-GDD q^(u+1) + h^u d^1} with $h = 2^{t} q$,
there exists a 4-GDD of type $(2^t q^2)^u m^1$ for admissible $m \ge 2 q$.

Take a 4-DGDD of type $(2^t q^2, (2^{t - 1} q)^{2q})^u$,
adjoin $m$ points and overlay the holes plus the new points with 4-GDDs of type $(2^{t - 1} q)^u m^1$.
Noting that $u - 1 \ge 5$, this gives a 4-GDD of type $(2^t q^2)^u m^1$ for admissible $m \le 5 \cdot 2^{t - 2} q$.
Since $t \ge 2$ this interval overlaps the previous one.

Finally, assume $s \ge 3$ and that this theorem holds for smaller $s$. Let $t \ge 1$.

By Theorem~\ref{thm:(u+1)-GDD q^(u+1) + h^u d^1} with $h = 2^t q^{s - 1}$,
there exists a 4-GDD of type $(2^t q^s)^u m^1$ for admissible $m \ge 2 q$.

Take a 4-DGDD of type $(2^t q^s, (2^t q^{s - 2})^{q^2})^u$,
adjoin $m$ points and overlay the holes plus the new points with 4-GDDs of type $(2^t q^{s - 2})^u m^1$.
This gives a 4-GDD of type $(2^t q^s)^u m^1$ for admissible $m \le 5 q^{s - 2}$. Since $s \ge 3$ this interval overlaps the previous one.
~\adfQED


\newpage
\appendix
%
%
%
%
%
%
%
%
%
\newcommand{\ADFvfyParStart}[1]{{\par\noindent#1}}
%
\newcommand{\adfDgap}{\vskip 1.75mm}      
\newcommand{\adfLgap}{\vskip 0.75mm}      
\newcommand{\adfsplit}{\par}              




\section{4-GDDs for the proof of Lemma \ref{lem:4-GDD general}}
\label{app:4-GDD general}
\adfhide{
$ 4^2 10^5 $,
$ 4^3 10^4 $,
$ 4^4 10^3 $,
$ 4^5 10^2 $,
$ 4^2 10^4 1^1 $,
$ 8^2 2^8 $,
$ 8^3 2^7 $,
$ 8^4 2^6 $,
$ 8^5 2^5 $,
$ 8^6 2^4 $,
$ 8^7 2^3 $,
$ 8^2 2^{11} $,
$ 8^3 2^6 5^1 $,
$ 8^5 2^4 5^1 $,
$ 8^7 2^2 5^1 $,
$ 8^5 14^1 20^1 $,
$ 20^4 8^2 2^1 $,
$ 20^4 8^2 5^1 $ and
$ 20^5 8^1 5^1 $.
}

\adfDgap
\noindent{\boldmath $ 4^{2} 10^{5} $}~
With the point set $Z_{58}$ partitioned into
 residue classes modulo $4$ for $\{0, 1, \dots, 39\}$,
 $\{40, 41, \dots, 49\}$, and
 residue classes modulo $2$ for $\{50, 51, \dots, 57\}$,
 the design is generated from

\adfLgap 
$(53, 50, 40, 34)$,
$(23, 4, 6, 25)$,
$(21, 4, 3, 34)$,
$(30, 11, 17, 0)$,\adfsplit
$(6, 8, 7, 17)$,
$(26, 23, 9, 20)$,
$(17, 36, 3, 2)$,
$(45, 20, 3, 18)$,\adfsplit
$(46, 38, 35, 21)$,
$(46, 24, 30, 39)$,
$(40, 3, 22, 32)$,
$(48, 2, 0, 5)$,\adfsplit
$(46, 8, 3, 14)$,
$(47, 2, 1, 19)$,
$(45, 21, 26, 31)$,
$(43, 21, 28, 7)$,\adfsplit
$(40, 18, 19, 5)$,
$(44, 21, 18, 32)$,
$(46, 9, 2, 27)$,
$(42, 37, 4, 2)$,\adfsplit
$(43, 16, 2, 25)$,
$(0, 7, 40, 51)$,
$(4, 9, 15, 40)$,
$(6, 37, 40, 52)$,\adfsplit
$(2, 7, 41, 50)$,
$(4, 7, 48, 56)$,
$(3, 37, 49, 51)$,
$(0, 35, 37, 53)$,\adfsplit
$(4, 27, 29, 53)$,
$(5, 27, 38, 42)$,
$(4, 17, 44, 54)$,
$(0, 14, 15, 42)$,\adfsplit
$(1, 16, 34, 48)$,
$(2, 24, 54, 55)$,
$(4, 18, 35, 49)$,
$(1, 4, 39, 57)$,\adfsplit
$(2, 35, 44, 51)$,
$(5, 6, 35, 54)$,
$(2, 15, 33, 52)$,
$(1, 35, 50, 55)$,\adfsplit
$(0, 25, 46, 55)$,
$(0, 22, 31, 57)$,
$(1, 22, 28, 53)$,
$(2, 29, 39, 56)$,\adfsplit
$(1, 3, 41, 56)$,
$(1, 26, 42, 52)$,
$(0, 1, 49, 54)$,
$(3, 9, 16, 44)$,\adfsplit
$(0, 26, 39, 44)$,
$(0, 18, 33, 50)$,
$(0, 19, 43, 56)$,
$(0, 13, 47, 52)$,\adfsplit
$(6, 28, 50, 57)$,
$(8, 18, 29, 52)$,
$(6, 16, 46, 53)$,
$(6, 9, 47, 55)$,\adfsplit
$(8, 39, 48, 51)$,
$(3, 38, 47, 53)$,
$(3, 13, 48, 57)$

\adfLgap \noindent by the mapping:
$x \mapsto x + 10 j \adfmod{40}$ for $x < 40$,
$x \mapsto (x + 5 j \adfmod{10}) + 40$ for $40 \le x < 50$,
$x \mapsto (x - 50 + 2 j \adfmod{8}) + 50$ for $x \ge 50$,
$0 \le j < 4$.
\ADFvfyParStart{(58, ((59, 4, ((40, 10), (10, 5), (8, 2)))), ((10, 4), (10, 1), (4, 2)))} 

\adfDgap
\noindent{\boldmath $ 4^{3} 10^{4} $}~
With the point set $Z_{52}$ partitioned into
 residue classes modulo $4$ for $\{0, 1, \dots, 39\}$, and
 residue classes modulo $3$ for $\{40, 41, \dots, 51\}$,
 the design is generated from

\adfLgap 
$(44, 40, 14, 12)$,
$(41, 51, 26, 7)$,
$(40, 41, 1, 28)$,
$(45, 49, 38, 17)$,\adfsplit
$(40, 45, 3, 2)$,
$(45, 50, 11, 36)$,
$(47, 48, 36, 3)$,
$(46, 45, 9, 19)$,\adfsplit
$(40, 47, 23, 37)$,
$(42, 50, 38, 5)$,
$(42, 40, 30, 24)$,
$(50, 40, 26, 35)$,\adfsplit
$(1, 3, 6, 16)$,
$(1, 7, 18, 36)$,
$(0, 6, 25, 40)$,
$(0, 26, 39, 49)$,\adfsplit
$(4, 6, 37, 39)$,
$(2, 19, 36, 51)$,
$(2, 15, 16, 33)$,
$(4, 26, 27, 43)$,\adfsplit
$(5, 16, 18, 19)$,
$(2, 25, 27, 46)$,
$(1, 23, 34, 49)$,
$(2, 8, 29, 43)$,\adfsplit
$(0, 1, 11, 46)$,
$(3, 5, 8, 49)$,
$(1, 8, 38, 51)$,
$(2, 5, 28, 39)$,\adfsplit
$(4, 9, 18, 47)$,
$(3, 13, 28, 34)$,
$(0, 18, 27, 41)$,
$(0, 3, 37, 38)$,\adfsplit
$(0, 2, 9, 31)$,
$(2, 7, 17, 48)$,
$(0, 10, 33, 51)$,
$(2, 12, 23, 50)$,\adfsplit
$(1, 14, 27, 32)$,
$(3, 4, 29, 50)$,
$(0, 13, 19, 42)$,
$(0, 7, 29, 44)$,\adfsplit
$(1, 24, 39, 50)$,
$(1, 15, 22, 48)$,
$(0, 21, 22, 50)$,
$(0, 14, 17, 35)$,\adfsplit
$(0, 5, 15, 47)$,
$(4, 5, 34, 51)$,
$(1, 4, 35, 45)$

\adfLgap \noindent by the mapping:
$x \mapsto x + 10 j \adfmod{40}$ for $x < 40$,
$x \mapsto (x - 40 + 3 j \adfmod{12}) + 40$ for $x \ge 40$,
$0 \le j < 4$.
\ADFvfyParStart{(52, ((47, 4, ((40, 10), (12, 3)))), ((10, 4), (4, 3)))} 

\adfDgap
\noindent{\boldmath $ 4^{4} 10^{3} $}~
With the point set $Z_{46}$ partitioned into
 residue classes modulo $3$ for $\{0, 1, \dots, 29\}$,
 residue classes modulo $3$ for $\{30, 31, \dots, 41\}$, and
 $\{42, 43, 44, 45\}$,
 the design is generated from

\adfLgap 
$(1, 2, 3, 45)$,
$(4, 30, 31, 45)$,
$(0, 4, 29, 32)$,
$(1, 9, 29, 31)$,\adfsplit
$(0, 19, 26, 30)$,
$(0, 1, 14, 39)$,
$(2, 4, 18, 37)$,
$(2, 9, 28, 42)$,\adfsplit
$(2, 7, 24, 44)$,
$(3, 4, 39, 40)$,
$(3, 29, 36, 44)$,
$(1, 18, 23, 40)$,\adfsplit
$(0, 16, 23, 41)$,
$(0, 8, 28, 40)$,
$(0, 2, 13, 33)$,
$(3, 35, 37, 42)$,\adfsplit
$(0, 22, 31, 35)$,
$(0, 11, 37, 44)$,
$(2, 34, 36, 41)$,
$(1, 11, 27, 41)$,\adfsplit
$(0, 7, 17, 34)$,
$(1, 12, 30, 34)$,
$(0, 5, 38, 43)$,
$(1, 6, 38, 44)$,\adfsplit
$(45, 0, 10, 20)$

\adfLgap \noindent by the mapping:
$x \mapsto x + 5 j \adfmod{30}$ for $x < 30$,
$x \mapsto (x - 30 + 2 j \adfmod{12}) + 30$ for $30 \le x < 42$,
$x \mapsto (x +  j \adfmod{3}) + 42$ for $42 \le x < 45$,
$45 \mapsto 45$,
$0 \le j < 6$
 for the first 24 blocks,
$0 \le j < 2$
 for the last block.
\ADFvfyParStart{(46, ((24, 6, ((30, 5), (12, 2), (3, 1), (1, 1))), (1, 2, ((30, 5), (12, 2), (3, 1), (1, 1)))), ((10, 3), (4, 3), (4, 1)))} 

\adfDgap
\noindent{\boldmath $ 4^{5} 10^{2} $}~
With the point set $Z_{40}$ partitioned into
 residue classes modulo $5$ for $\{0, 1, \dots, 19\}$, and
 residue classes modulo $2$ for $\{20, 21, \dots, 39\}$,
 the design is generated from

\adfLgap 
$(0, 1, 20, 39)$,
$(0, 2, 23, 34)$,
$(0, 8, 30, 35)$,
$(0, 6, 17, 19)$,\adfsplit
$(0, 4, 28, 29)$,
$(0, 9, 26, 37)$,
$(0, 7, 31, 38)$,
$(0, 3, 33, 36)$,\adfsplit
$(1, 5, 21, 26)$,
$(1, 9, 23, 36)$,
$(1, 7, 30, 33)$

\adfLgap \noindent by the mapping:
$x \mapsto x + 2 j \adfmod{20}$ for $x < 20$,
$x \mapsto (x + 2 j \adfmod{20}) + 20$ for $x \ge 20$,
$0 \le j < 10$.
\ADFvfyParStart{(40, ((11, 10, ((20, 2), (20, 2)))), ((4, 5), (10, 2)))} 

\adfDgap
\noindent{\boldmath $ 4^{2} 10^{4} 1^{1} $}~
With the point set $Z_{49}$ partitioned into
 residue classes modulo $4$ for $\{0, 1, \dots, 39\}$,
 residue classes modulo $2$ for $\{40, 41, \dots, 47\}$, and
 $\{48\}$,
 the design is generated from

\adfLgap 
$(41, 44, 4, 39)$,
$(42, 45, 18, 8)$,
$(46, 47, 26, 11)$,
$(46, 45, 32, 7)$,\adfsplit
$(48, 47, 9, 35)$,
$(48, 46, 22, 3)$,
$(48, 16, 34, 27)$,
$(48, 0, 1, 38)$,\adfsplit
$(1, 34, 28, 39)$,
$(9, 2, 27, 28)$,
$(0, 14, 29, 35)$,
$(0, 11, 18, 33)$,\adfsplit
$(6, 8, 9, 39)$,
$(1, 3, 18, 41)$,
$(3, 8, 34, 37)$,
$(5, 18, 27, 40)$,\adfsplit
$(3, 16, 17, 38)$,
$(2, 8, 25, 41)$,
$(5, 8, 15, 44)$,
$(2, 13, 24, 46)$,\adfsplit
$(1, 15, 24, 26)$,
$(4, 14, 27, 47)$,
$(0, 5, 34, 46)$,
$(1, 4, 31, 46)$,\adfsplit
$(2, 3, 4, 29)$,
$(3, 24, 25, 47)$,
$(1, 22, 32, 35)$,
$(0, 2, 19, 44)$,\adfsplit
$(1, 6, 12, 47)$,
$(0, 21, 22, 27)$,
$(2, 16, 35, 37)$,
$(2, 33, 39, 43)$,\adfsplit
$(0, 6, 15, 41)$,
$(6, 16, 29, 40)$,
$(0, 9, 31, 40)$,
$(3, 13, 36, 42)$,\adfsplit
$(0, 7, 17, 42)$,
$(0, 23, 25, 26)$,
$(0, 37, 39, 47)$,
$(0, 3, 30, 43)$,\adfsplit
$(1, 27, 36, 43)$

\adfLgap \noindent by the mapping:
$x \mapsto x + 10 j \adfmod{40}$ for $x < 40$,
$x \mapsto (x + 2 j \adfmod{8}) + 40$ for $40 \le x < 48$,
$48 \mapsto 48$,
$0 \le j < 4$.
\ADFvfyParStart{(49, ((41, 4, ((40, 10), (8, 2), (1, 1)))), ((10, 4), (4, 2), (1, 1)))} 

\adfDgap
\noindent{\boldmath $ 8^{2} 2^{8} $}~
With the point set $Z_{32}$ partitioned into
 residue classes modulo $2$ for $\{0, 1, \dots, 15\}$, and
 residue classes modulo $8$ for $\{16, 17, \dots, 31\}$,
 the design is generated from

\adfLgap 
$(26, 25, 27, 20)$,
$(28, 27, 29, 22)$,
$(30, 29, 31, 24)$,
$(16, 31, 17, 26)$,\adfsplit
$(27, 21, 9, 2)$,
$(1, 2, 16, 18)$,
$(0, 1, 17, 27)$,
$(0, 9, 16, 20)$,\adfsplit
$(1, 4, 23, 26)$,
$(1, 6, 25, 30)$,
$(1, 12, 20, 22)$,
$(1, 14, 21, 31)$,\adfsplit
$(0, 3, 23, 30)$,
$(2, 5, 26, 30)$,
$(0, 11, 18, 29)$,
$(3, 14, 18, 27)$,\adfsplit
$(0, 7, 22, 26)$,
$(5, 6, 22, 24)$,
$(4, 5, 18, 21)$,
$(6, 7, 18, 20)$,\adfsplit
$(3, 6, 16, 28)$,
$(6, 15, 27, 31)$,
$(5, 14, 25, 29)$,
$(3, 12, 17, 24)$,\adfsplit
$(4, 7, 24, 27)$,
$(2, 13, 19, 24)$,
$(3, 4, 19, 22)$,
$(0, 15, 21, 24)$,\adfsplit
$(2, 3, 20, 29)$,
$(4, 13, 20, 31)$,
$(0, 5, 19, 31)$,
$(2, 11, 17, 23)$,\adfsplit
$(4, 15, 17, 29)$,
$(2, 15, 22, 25)$,
$(2, 7, 28, 31)$,
$(0, 13, 25, 28)$

\adfLgap \noindent by the mapping:
$x \mapsto x + 8 j \adfmod{16}$ for $x < 16$,
$x \mapsto (x + 8 j \adfmod{16}) + 16$ for $x \ge 16$,
$0 \le j < 2$.
\ADFvfyParStart{(32, ((36, 2, ((16, 8), (16, 8)))), ((8, 2), (2, 8)))} 

\adfDgap
\noindent{\boldmath $ 8^{3} 2^{7} $}~
With the point set $Z_{38}$ partitioned into
 residue classes modulo $6$ for $\{0, 1, \dots, 11\}$,
 $\{12, 13\}$, and
 residue classes modulo $3$ for $\{14, 15, \dots, 37\}$,
 the design is generated from

\adfLgap 
$(0, 2, 17, 33)$,
$(0, 21, 23, 25)$,
$(1, 20, 21, 31)$,
$(1, 3, 17, 36)$,\adfsplit
$(12, 21, 22, 35)$,
$(1, 5, 14, 33)$,
$(12, 16, 20, 33)$,
$(1, 18, 22, 29)$,\adfsplit
$(0, 20, 22, 30)$,
$(0, 14, 16, 24)$,
$(0, 3, 15, 32)$,
$(0, 4, 11, 36)$,\adfsplit
$(1, 19, 23, 24)$,
$(0, 1, 13, 27)$,
$(0, 18, 19, 35)$,
$(0, 5, 12, 34)$,\adfsplit
$(1, 4, 15, 34)$

\adfLgap \noindent by the mapping:
$x \mapsto x + 2 j \adfmod{12}$ for $x < 12$,
$x \mapsto (x +  j \adfmod{2}) + 12$ for $12 \le x < 14$,
$x \mapsto (x - 14 + 4 j \adfmod{24}) + 14$ for $x \ge 14$,
$0 \le j < 6$.
\ADFvfyParStart{(38, ((17, 6, ((12, 2), (2, 1), (24, 4)))), ((2, 6), (2, 1), (8, 3)))} 

\adfDgap
\noindent{\boldmath $ 8^{4} 2^{6} $}~
With the point set $Z_{44}$ partitioned into
 residue classes modulo $3$ for $\{0, 1, \dots, 23\}$,
 $\{24, 25, \dots, 31\}$, and
 residue classes modulo $6$ for $\{32, 33, \dots, 43\}$,
 the design is generated from

\adfLgap 
$(0, 2, 10, 32)$,
$(0, 1, 8, 37)$,
$(0, 5, 13, 34)$,
$(0, 4, 11, 24)$,\adfsplit
$(0, 14, 19, 27)$,
$(0, 31, 35, 36)$,
$(0, 22, 30, 39)$,
$(0, 23, 26, 38)$,\adfsplit
$(0, 29, 40, 42)$,
$(0, 25, 41, 43)$,
$(1, 3, 23, 25)$,
$(1, 2, 6, 29)$,\adfsplit
$(1, 14, 31, 33)$,
$(1, 15, 27, 36)$,
$(1, 5, 30, 41)$,
$(2, 15, 30, 38)$,\adfsplit
$(1, 34, 38, 43)$,
$(1, 18, 24, 42)$,
$(1, 11, 28, 35)$,
$(2, 3, 36, 37)$,\adfsplit
$(2, 19, 35, 43)$,
$(3, 11, 33, 38)$,
$(2, 28, 33, 42)$

\adfLgap \noindent by the mapping:
$x \mapsto x + 4 j \adfmod{24}$ for $x < 24$,
$x \mapsto (x + 4 j \adfmod{8}) + 24$ for $24 \le x < 32$,
$x \mapsto (x - 32 + 2 j \adfmod{12}) + 32$ for $x \ge 32$,
$0 \le j < 6$.
\ADFvfyParStart{(44, ((23, 6, ((24, 4), (8, 4), (12, 2)))), ((8, 3), (8, 1), (2, 6)))} 

\adfDgap
\noindent{\boldmath $ 8^{5} 2^{5} $}~
With the point set $Z_{50}$ partitioned into
 residue classes modulo $5$ for $\{0, 1, \dots, 39\}$, and
 residue classes modulo $5$ for $\{40, 41, \dots, 49\}$,
 the design is generated from

\adfLgap 
$(0, 1, 2, 34)$,
$(0, 3, 22, 26)$,
$(0, 7, 14, 38)$,
$(0, 6, 9, 18)$,\adfsplit
$(1, 5, 18, 46)$,
$(2, 31, 39, 41)$,
$(1, 30, 40, 47)$,
$(2, 3, 45, 49)$,\adfsplit
$(1, 22, 42, 43)$,
$(2, 15, 44, 46)$,
$(1, 3, 27, 29)$,
$(0, 13, 21, 44)$,\adfsplit
$(1, 7, 19, 25)$,
$(0, 4, 37, 47)$,
$(0, 12, 29, 41)$,
$(0, 16, 39, 49)$,\adfsplit
$(0, 27, 31, 43)$,
$(0, 8, 19, 42)$

\adfLgap \noindent by the mapping:
$x \mapsto x + 4 j \adfmod{40}$ for $x < 40$,
$x \mapsto (x +  j \adfmod{10}) + 40$ for $x \ge 40$,
$0 \le j < 10$.
\ADFvfyParStart{(50, ((18, 10, ((40, 4), (10, 1)))), ((8, 5), (2, 5)))} 

\adfDgap
\noindent{\boldmath $ 8^{6} 2^{4} $}~
With the point set $Z_{56}$ partitioned into
 residue classes modulo $6$ for $\{0, 1, \dots, 47\}$,
 residue classes modulo $3$ for $\{48, 49, \dots, 53\}$, and
 $\{54, 55\}$,
 the design is generated from

\adfLgap 
$(0, 1, 14, 46)$,
$(0, 2, 9, 25)$,
$(0, 4, 33, 41)$,
$(0, 3, 5, 32)$,\adfsplit
$(0, 7, 20, 40)$,
$(0, 10, 13, 29)$,
$(0, 22, 26, 35)$,
$(0, 11, 38, 43)$,\adfsplit
$(0, 17, 31, 39)$,
$(0, 23, 27, 52)$,
$(0, 34, 44, 49)$,
$(0, 37, 48, 55)$,\adfsplit
$(0, 47, 50, 54)$,
$(0, 45, 51, 53)$,
$(1, 2, 4, 49)$,
$(1, 3, 6, 53)$,\adfsplit
$(1, 10, 18, 48)$,
$(1, 22, 36, 50)$,
$(1, 11, 21, 51)$,
$(1, 35, 39, 52)$,\adfsplit
$(4, 30, 48, 53)$,
$(5, 22, 39, 51)$,
$(2, 22, 47, 48)$,
$(3, 4, 43, 51)$,\adfsplit
$(2, 29, 39, 50)$,
$(2, 7, 27, 34)$,
$(2, 6, 35, 54)$,
$(2, 19, 30, 45)$,\adfsplit
$(1, 27, 44, 54)$,
$(3, 23, 28, 44)$,
$(2, 3, 36, 37)$,
$(6, 7, 39, 46)$,\adfsplit
$(2, 28, 31, 55)$,
$(4, 13, 14, 21)$,
$(4, 12, 23, 45)$,
$(1, 5, 28, 30)$,\adfsplit
$(1, 34, 45, 47)$,
$(1, 29, 38, 55)$

\adfLgap \noindent by the mapping:
$x \mapsto x + 8 j \adfmod{48}$ for $x < 48$,
$x \mapsto (x +  j \adfmod{6}) + 48$ for $48 \le x < 54$,
$x \mapsto (x +  j \adfmod{2}) + 54$ for $x \ge 54$,
$0 \le j < 6$.
\ADFvfyParStart{(56, ((38, 6, ((48, 8), (6, 1), (2, 1)))), ((8, 6), (2, 3), (2, 1)))} 

\adfDgap
\noindent{\boldmath $ 8^{7} 2^{3} $}~
With the point set $Z_{62}$ partitioned into
 residue classes modulo $6$ for $\{0, 1, \dots, 47\}$,
 $\{48, 49, \dots, 55\}$, and
 residue classes modulo $3$ for $\{56, 57, \dots, 61\}$,
 the design is generated from

\adfLgap 
$(0, 1, 14, 46)$,
$(0, 2, 9, 25)$,
$(0, 3, 4, 32)$,
$(0, 5, 27, 43)$,\adfsplit
$(0, 7, 8, 21)$,
$(0, 10, 11, 51)$,
$(0, 15, 17, 55)$,
$(0, 22, 26, 41)$,\adfsplit
$(0, 23, 28, 53)$,
$(0, 29, 33, 49)$,
$(0, 31, 34, 50)$,
$(0, 35, 37, 54)$,\adfsplit
$(0, 38, 48, 61)$,
$(0, 39, 52, 57)$,
$(0, 44, 58, 59)$,
$(0, 45, 56, 60)$,\adfsplit
$(1, 2, 4, 48)$,
$(1, 6, 29, 52)$,
$(2, 5, 15, 52)$,
$(3, 11, 20, 52)$,\adfsplit
$(1, 3, 10, 41)$,
$(1, 15, 22, 61)$,
$(1, 20, 27, 53)$,
$(1, 12, 23, 57)$,\adfsplit
$(1, 39, 50, 60)$,
$(1, 35, 38, 59)$,
$(1, 36, 54, 56)$,
$(1, 5, 21, 44)$,\adfsplit
$(1, 28, 51, 58)$,
$(3, 14, 28, 36)$,
$(2, 12, 28, 31)$,
$(4, 14, 31, 39)$,\adfsplit
$(4, 6, 21, 50)$,
$(2, 36, 37, 51)$,
$(4, 37, 45, 57)$,
$(2, 22, 27, 54)$,\adfsplit
$(3, 7, 22, 23)$,
$(5, 7, 46, 49)$,
$(2, 11, 39, 61)$,
$(3, 47, 51, 57)$,\adfsplit
$(3, 30, 53, 59)$,
$(2, 7, 29, 34)$,
$(2, 13, 30, 57)$,
$(2, 21, 35, 56)$,\adfsplit
$(5, 6, 14, 55)$,
$(2, 6, 19, 58)$,
$(2, 10, 53, 60)$

\adfLgap \noindent by the mapping:
$x \mapsto x + 8 j \adfmod{48}$ for $x < 48$,
$x \mapsto (x + 4 j \adfmod{8}) + 48$ for $48 \le x < 56$,
$x \mapsto (x - 56 +  j \adfmod{6}) + 56$ for $x \ge 56$,
$0 \le j < 6$.
\ADFvfyParStart{(62, ((47, 6, ((48, 8), (8, 4), (6, 1)))), ((8, 6), (8, 1), (2, 3)))} 

\adfDgap
\noindent{\boldmath $ 8^{2} 2^{11} $}~
With the point set $Z_{38}$ partitioned into
 residue classes modulo $10$ for $\{0, 1, \dots, 19\}$,
 $\{20, 21\}$, and
 residue classes modulo $2$ for $\{22, 23, \dots, 37\}$,
 the design is generated from

\adfLgap 
$(25, 26, 20, 4)$,
$(29, 30, 21, 9)$,
$(34, 23, 21, 11)$,
$(22, 27, 20, 16)$,\adfsplit
$(37, 32, 20, 0)$,
$(1, 2, 20, 23)$,
$(3, 5, 7, 20)$,
$(8, 9, 20, 28)$,\adfsplit
$(0, 2, 22, 31)$,
$(3, 4, 22, 23)$,
$(0, 5, 23, 28)$,
$(3, 6, 28, 31)$,\adfsplit
$(4, 7, 8, 31)$,
$(5, 9, 31, 32)$,
$(6, 8, 23, 32)$,
$(7, 9, 23, 26)$,\adfsplit
$(1, 3, 30, 37)$,
$(0, 4, 18, 30)$,
$(1, 5, 22, 35)$,
$(2, 7, 16, 30)$,\adfsplit
$(5, 17, 25, 30)$,
$(8, 19, 30, 33)$,
$(2, 17, 26, 29)$,
$(6, 7, 18, 29)$,\adfsplit
$(5, 8, 29, 36)$,
$(3, 15, 29, 34)$,
$(0, 12, 14, 29)$,
$(1, 16, 29, 32)$,\adfsplit
$(4, 19, 28, 29)$,
$(1, 18, 27, 34)$,
$(2, 8, 25, 34)$,
$(4, 5, 6, 34)$,\adfsplit
$(0, 1, 8, 15)$,
$(2, 6, 9, 15)$,
$(2, 5, 24, 33)$,
$(4, 15, 21, 35)$,\adfsplit
$(1, 6, 14, 19)$,
$(0, 6, 26, 27)$,
$(0, 9, 13, 34)$,
$(6, 21, 33, 36)$,\adfsplit
$(0, 3, 16, 33)$,
$(0, 7, 21, 24)$,
$(3, 14, 24, 25)$,
$(1, 7, 25, 36)$,\adfsplit
$(0, 11, 19, 25)$,
$(1, 4, 17, 24)$,
$(1, 12, 13, 28)$,
$(0, 17, 35, 36)$,\adfsplit
$(2, 19, 27, 32)$,
$(3, 9, 17, 27)$,
$(3, 18, 32, 35)$,
$(2, 13, 18, 21)$,\adfsplit
$(2, 14, 28, 35)$

\adfLgap \noindent by the mapping:
$x \mapsto x + 10 j \adfmod{20}$ for $x < 20$,
$x \mapsto x$ for $20 \le x < 22$,
$x \mapsto (x - 22 + 8 j \adfmod{16}) + 22$ for $x \ge 22$,
$0 \le j < 2$.
\ADFvfyParStart{(38, ((53, 2, ((20, 10), (2, 2), (16, 8)))), ((2, 10), (2, 1), (8, 2)))} 

\adfDgap
\noindent{\boldmath $ 8^{3} 2^{6} 5^{1} $}~
With the point set $Z_{41}$ partitioned into
 residue classes modulo $3$ for $\{0, 1, \dots, 23\}$,
 residue classes modulo $6$ for $\{24, 25, \dots, 35\}$, and
 $\{36, 37, 38, 39, 40\}$,
 the design is generated from

\adfLgap 
$(0, 2, 22, 24)$,
$(1, 2, 18, 25)$,
$(0, 1, 8, 29)$,
$(0, 4, 11, 28)$,\adfsplit
$(0, 14, 19, 27)$,
$(0, 5, 33, 40)$,
$(2, 3, 34, 40)$,
$(0, 10, 32, 37)$,\adfsplit
$(0, 13, 23, 39)$,
$(0, 30, 35, 38)$,
$(0, 31, 34, 36)$,
$(1, 3, 14, 38)$,\adfsplit
$(1, 6, 23, 36)$,
$(2, 30, 31, 39)$,
$(2, 27, 32, 35)$,
$(1, 15, 27, 37)$,\adfsplit
$(1, 5, 24, 35)$,
$(3, 11, 27, 29)$,
$(1, 9, 30, 32)$,
$(3, 7, 28, 32)$

\adfLgap \noindent by the mapping:
$x \mapsto x + 4 j \adfmod{24}$ for $x < 24$,
$x \mapsto (x + 2 j \adfmod{12}) + 24$ for $24 \le x < 36$,
$x \mapsto x$ for $x \ge 36$,
$0 \le j < 6$.
\ADFvfyParStart{(41, ((20, 6, ((24, 4), (12, 2), (5, 5)))), ((8, 3), (2, 6), (5, 1)))} 

\adfDgap
\noindent{\boldmath $ 8^{5} 2^{4} 5^{1} $}~
With the point set $Z_{53}$ partitioned into
 residue classes modulo $5$ for $\{0, 1, \dots, 39\}$,
 residue classes modulo $4$ for $\{40, 41, \dots, 47\}$, and
 $\{48, 49, 50, 51, 52\}$,
 the design is generated from

\adfLgap 
$(52, 44, 47, 4)$,
$(51, 40, 45, 23)$,
$(49, 46, 47, 15)$,
$(48, 45, 7, 8)$,\adfsplit
$(50, 47, 22, 13)$,
$(50, 42, 16, 29)$,
$(48, 42, 22, 21)$,
$(42, 41, 38, 20)$,\adfsplit
$(46, 44, 7, 21)$,
$(41, 47, 36, 5)$,
$(11, 25, 34, 8)$,
$(25, 6, 32, 38)$,\adfsplit
$(33, 25, 12, 1)$,
$(35, 33, 29, 36)$,
$(33, 9, 27, 10)$,
$(3, 4, 15, 48)$,\adfsplit
$(0, 6, 9, 48)$,
$(2, 4, 5, 42)$,
$(2, 9, 25, 44)$,
$(1, 17, 35, 44)$,\adfsplit
$(3, 25, 27, 47)$,
$(1, 5, 8, 50)$,
$(5, 19, 38, 52)$,
$(5, 9, 37, 51)$,\adfsplit
$(4, 25, 36, 37)$,
$(0, 14, 37, 50)$,
$(4, 8, 26, 40)$,
$(0, 16, 33, 44)$,\adfsplit
$(3, 9, 17, 42)$,
$(0, 19, 28, 42)$,
$(3, 14, 16, 40)$,
$(0, 2, 38, 40)$,\adfsplit
$(2, 26, 28, 51)$,
$(4, 17, 28, 49)$,
$(1, 28, 34, 37)$,
$(0, 8, 27, 36)$,\adfsplit
$(3, 7, 36, 52)$,
$(1, 7, 9, 47)$,
$(0, 7, 26, 41)$,
$(1, 18, 19, 43)$,\adfsplit
$(2, 18, 29, 41)$,
$(2, 19, 36, 49)$,
$(0, 29, 31, 47)$,
$(2, 6, 34, 47)$,\adfsplit
$(0, 4, 32, 43)$,
$(1, 3, 29, 32)$,
$(0, 1, 22, 52)$,
$(0, 12, 13, 34)$,\adfsplit
$(1, 4, 13, 45)$,
$(0, 3, 21, 49)$,
$(0, 11, 24, 51)$

\adfLgap \noindent by the mapping:
$x \mapsto x + 10 j \adfmod{40}$ for $x < 40$,
$x \mapsto (x + 2 j \adfmod{8}) + 40$ for $40 \le x < 48$,
$x \mapsto x$ for $x \ge 48$,
$0 \le j < 4$.
\ADFvfyParStart{(53, ((51, 4, ((40, 10), (8, 2), (5, 5)))), ((8, 5), (2, 4), (5, 1)))} 

\adfDgap
\noindent{\boldmath $ 8^{7} 2^{2} 5^{1} $}~
With the point set $Z_{65}$ partitioned into
 residue classes modulo $7$ for $\{0, 1, \dots, 55\}$,
 residue classes modulo $2$ for $\{56, 57, 58, 59\}$, and
 $\{60, 61, 62, 63, 64\}$,
 the design is generated from

\adfLgap 
$(60, 58, 3, 47)$,
$(64, 57, 17, 1)$,
$(61, 59, 34, 21)$,
$(63, 58, 54, 34)$,\adfsplit
$(62, 59, 41, 14)$,
$(57, 56, 41, 52)$,
$(0, 1, 2, 56)$,
$(0, 3, 30, 57)$,\adfsplit
$(1, 4, 16, 58)$,
$(0, 4, 6, 59)$,
$(1, 3, 40, 59)$,
$(2, 5, 18, 58)$,\adfsplit
$(4, 5, 7, 56)$,
$(5, 8, 9, 57)$,
$(6, 8, 53, 58)$,
$(7, 9, 53, 59)$,\adfsplit
$(7, 22, 24, 57)$,
$(9, 11, 52, 58)$,
$(8, 26, 39, 59)$,
$(9, 12, 13, 56)$,\adfsplit
$(0, 5, 20, 31)$,
$(0, 8, 17, 48)$,
$(1, 5, 10, 34)$,
$(0, 9, 34, 47)$,\adfsplit
$(0, 10, 11, 19)$,
$(0, 13, 15, 33)$,
$(1, 7, 11, 33)$,
$(1, 13, 18, 47)$,\adfsplit
$(3, 4, 9, 19)$,
$(2, 7, 47, 48)$,
$(2, 8, 21, 33)$,
$(2, 13, 19, 38)$,\adfsplit
$(3, 5, 11, 62)$,
$(5, 13, 22, 64)$,
$(5, 25, 38, 61)$,
$(5, 37, 50, 63)$,\adfsplit
$(1, 9, 20, 25)$,
$(2, 3, 6, 25)$,
$(0, 12, 22, 39)$,
$(0, 25, 40, 45)$,\adfsplit
$(1, 6, 30, 39)$,
$(0, 16, 53, 64)$,
$(2, 35, 53, 55)$,
$(3, 13, 21, 53)$,\adfsplit
$(1, 24, 53, 63)$,
$(8, 11, 41, 60)$,
$(9, 27, 39, 40)$,
$(0, 18, 44, 63)$,\adfsplit
$(3, 27, 49, 63)$,
$(0, 36, 41, 51)$,
$(0, 24, 29, 60)$,
$(0, 43, 54, 61)$,\adfsplit
$(0, 26, 37, 52)$,
$(0, 32, 38, 55)$,
$(0, 23, 46, 50)$,
$(2, 27, 45, 54)$,\adfsplit
$(3, 37, 46, 61)$,
$(2, 22, 41, 61)$,
$(1, 26, 32, 41)$,
$(2, 17, 36, 40)$,\adfsplit
$(4, 35, 38, 41)$,
$(1, 27, 31, 51)$,
$(3, 8, 20, 32)$,
$(3, 7, 18, 36)$,\adfsplit
$(1, 21, 37, 45)$,
$(1, 23, 35, 62)$,
$(1, 48, 49, 54)$,
$(1, 38, 44, 46)$,\adfsplit
$(6, 23, 52, 64)$,
$(4, 21, 40, 64)$,
$(7, 37, 38, 54)$,
$(2, 10, 50, 62)$,\adfsplit
$(6, 24, 35, 36)$,
$(4, 20, 23, 60)$,
$(2, 12, 20, 24)$,
$(4, 24, 26, 50)$,\adfsplit
$(2, 26, 49, 60)$,
$(4, 12, 34, 62)$

\adfLgap \noindent by the mapping:
$x \mapsto x + 14 j \adfmod{56}$ for $x < 56$,
$x \mapsto (x +  j \adfmod{4}) + 56$ for $56 \le x < 60$,
$x \mapsto x$ for $x \ge 60$,
$0 \le j < 4$.
\ADFvfyParStart{(65, ((78, 4, ((56, 14), (4, 1), (5, 5)))), ((8, 7), (2, 2), (5, 1)))} 

\adfDgap
\noindent{\boldmath $ 8^{5} 14^{1} 20^{1} $}~
With the point set $Z_{74}$ partitioned into
 residue classes modulo $5$ for $\{0, 1, \dots, 39\}$,
 $\{40, 41, \dots, 59\}$, and
 $\{60, 61, \dots, 73\}$,
 the design is generated from

\adfLgap 
$(0, 1, 40, 70)$,
$(0, 2, 43, 72)$,
$(0, 3, 46, 73)$,
$(1, 5, 43, 71)$,\adfsplit
$(0, 4, 11, 49)$,
$(0, 6, 32, 57)$,
$(0, 9, 17, 55)$,
$(0, 13, 19, 22)$,\adfsplit
$(0, 12, 59, 60)$,
$(0, 16, 58, 61)$,
$(0, 29, 56, 65)$,
$(0, 33, 51, 69)$,\adfsplit
$(0, 27, 45, 62)$,
$(0, 39, 52, 67)$,
$(0, 23, 48, 68)$,
$(0, 21, 44, 66)$,\adfsplit
$(1, 3, 29, 50)$,
$(1, 19, 48, 64)$,
$(1, 17, 46, 60)$

\adfLgap \noindent by the mapping:
$x \mapsto x + 2 j \adfmod{40}$ for $x < 40$,
$x \mapsto (x +  j \adfmod{20}) + 40$ for $40 \le x < 60$,
$x \mapsto (x +  j \adfmod{10}) + 60$ for $60 \le x < 70$,
$x \mapsto (x - 70 +  j \adfmod{4}) + 70$ for $x \ge 70$,
$0 \le j < 20$.
\ADFvfyParStart{(74, ((19, 20, ((40, 2), (20, 1), (10, 1), (4, 1)))), ((8, 5), (20, 1), (14, 1)))} 

\adfDgap
\noindent{\boldmath $ 20^{4} 8^{2} 2^{1} $}~
With the point set $Z_{98}$ partitioned into
 residue classes modulo $2$ for $\{0, 1, \dots, 15\}$,
 residue classes modulo $4$ for $\{16, 17, \dots, 95\}$, and
 $\{96, 97\}$,
 the design is generated from

\adfLgap 
$(96, 9, 2, 27)$,
$(85, 92, 4, 9)$,
$(55, 89, 2, 15)$,
$(94, 65, 10, 9)$,\adfsplit
$(97, 66, 85, 87)$,
$(97, 40, 91, 73)$,
$(97, 64, 79, 38)$,
$(28, 37, 78, 2)$,\adfsplit
$(50, 20, 87, 17)$,
$(75, 13, 41, 78)$,
$(25, 90, 64, 11)$,
$(41, 15, 80, 66)$,\adfsplit
$(0, 21, 70, 80)$,
$(0, 24, 30, 55)$,
$(16, 25, 43, 70)$,
$(16, 19, 85, 90)$,\adfsplit
$(17, 18, 75, 92)$,
$(0, 33, 38, 75)$,
$(17, 38, 48, 55)$,
$(0, 25, 52, 82)$,\adfsplit
$(18, 19, 40, 53)$,
$(16, 77, 94, 95)$,
$(18, 20, 39, 89)$,
$(0, 16, 58, 73)$,\adfsplit
$(0, 26, 43, 77)$,
$(0, 63, 66, 88)$,
$(0, 28, 34, 57)$,
$(0, 23, 36, 37)$,\adfsplit
$(17, 28, 39, 82)$,
$(0, 39, 53, 84)$,
$(0, 31, 76, 78)$,
$(0, 22, 64, 83)$,\adfsplit
$(0, 19, 41, 48)$,
$(0, 32, 59, 86)$,
$(0, 29, 54, 91)$,
$(16, 79, 82, 89)$,\adfsplit
$(0, 61, 67, 92)$,
$(0, 42, 51, 81)$,
$(16, 27, 29, 62)$,
$(0, 69, 74, 87)$,\adfsplit
$(16, 21, 31, 54)$

\adfLgap \noindent by the mapping:
$x \mapsto x +  j \adfmod{16}$ for $x < 16$,
$x \mapsto (x - 16 + 5 j \adfmod{80}) + 16$ for $16 \le x < 96$,
$x \mapsto (x +  j \adfmod{2}) + 96$ for $x \ge 96$,
$0 \le j < 16$.
\ADFvfyParStart{(98, ((41, 16, ((16, 1), (80, 5), (2, 1)))), ((8, 2), (20, 4), (2, 1)))} 

\adfDgap
\noindent{\boldmath $ 20^{4} 8^{2} 5^{1} $}~
With the point set $Z_{101}$ partitioned into
 residue classes modulo $2$ for $\{0, 1, \dots, 15\}$,
 residue classes modulo $4$ for $\{16, 17, \dots, 95\}$, and
 $\{96, 97, 98, 99, 100\}$,
 the design is generated from

\adfLgap 
$(0, 1, 25, 40)$,
$(0, 3, 60, 70)$,
$(0, 7, 30, 85)$,
$(16, 17, 30, 75)$,\adfsplit
$(0, 16, 90, 95)$,
$(16, 18, 25, 55)$,
$(0, 17, 31, 80)$,
$(16, 19, 21, 50)$,\adfsplit
$(0, 18, 65, 76)$,
$(16, 29, 34, 35)$,
$(16, 27, 37, 70)$,
$(17, 19, 40, 74)$,\adfsplit
$(17, 20, 23, 34)$,
$(17, 39, 68, 90)$,
$(18, 20, 83, 89)$,
$(17, 28, 70, 83)$,\adfsplit
$(17, 24, 95, 100)$,
$(17, 55, 78, 96)$,
$(17, 35, 54, 98)$,
$(18, 39, 80, 98)$,\adfsplit
$(18, 19, 45, 97)$,
$(0, 5, 44, 54)$,
$(0, 21, 64, 94)$,
$(0, 22, 49, 84)$,\adfsplit
$(0, 24, 39, 58)$,
$(0, 28, 33, 59)$,
$(0, 37, 38, 79)$,
$(0, 34, 47, 81)$,\adfsplit
$(0, 69, 72, 91)$,
$(0, 51, 74, 88)$,
$(0, 42, 89, 98)$,
$(16, 26, 41, 79)$,\adfsplit
$(0, 23, 61, 68)$,
$(0, 48, 78, 87)$,
$(0, 43, 77, 92)$,
$(0, 46, 52, 73)$,\adfsplit
$(0, 53, 71, 100)$,
$(0, 63, 66, 96)$,
$(16, 33, 46, 87)$,
$(0, 41, 86, 97)$,\adfsplit
$(0, 36, 83, 93)$,
$(0, 32, 62, 67)$,
$(0, 26, 57, 99)$,
$(0, 27, 56, 82)$

\adfLgap \noindent by the mapping:
$x \mapsto x +  j \adfmod{16}$ for $x < 16$,
$x \mapsto (x - 16 + 5 j \adfmod{80}) + 16$ for $16 \le x < 96$,
$x \mapsto (x +  j \adfmod{4}) + 96$ for $96 \le x < 100$,
$100 \mapsto 100$,
$0 \le j < 16$.
\ADFvfyParStart{(101, ((44, 16, ((16, 1), (80, 5), (4, 1), (1, 1)))), ((8, 2), (20, 4), (5, 1)))} 

\adfDgap
\noindent{\boldmath $ 20^{5} 8^{1} 5^{1} $}~
With the point set $Z_{113}$ partitioned into
 residue classes modulo $5$ for $\{0, 1, \dots, 99\}$,
 $\{100, 101, \dots, 107\}$, and
 $\{108, 109, 110, 111, 112\}$,
 the design is generated from

\adfLgap 
$(112, 100, 66, 9)$,
$(108, 102, 71, 14)$,
$(110, 103, 55, 89)$,
$(111, 105, 60, 94)$,\adfsplit
$(108, 88, 24, 30)$,
$(109, 93, 29, 35)$,
$(111, 38, 81, 62)$,
$(112, 43, 86, 67)$,\adfsplit
$(112, 22, 76, 39)$,
$(108, 27, 81, 44)$,
$(109, 83, 5, 66)$,
$(110, 88, 10, 71)$,\adfsplit
$(112, 88, 2, 29)$,
$(108, 93, 7, 34)$,
$(110, 50, 61, 2)$,
$(111, 55, 66, 7)$,\adfsplit
$(108, 83, 12, 10)$,
$(109, 88, 17, 15)$,
$(14, 82, 50, 93)$,
$(19, 87, 55, 98)$,\adfsplit
$(58, 95, 7, 54)$,
$(63, 0, 12, 59)$,
$(84, 6, 30, 2)$,
$(89, 11, 35, 7)$,\adfsplit
$(65, 39, 78, 66)$,
$(70, 44, 83, 71)$,
$(28, 36, 72, 94)$,
$(33, 41, 77, 99)$,\adfsplit
$(14, 90, 81, 37)$,
$(19, 95, 86, 42)$,
$(13, 81, 10, 52)$,
$(18, 86, 15, 57)$,\adfsplit
$(80, 48, 64, 76)$,
$(85, 53, 69, 81)$,
$(69, 82, 71, 13)$,
$(74, 87, 76, 18)$,\adfsplit
$(42, 8, 84, 80)$,
$(47, 13, 89, 85)$,
$(0, 8, 39, 82)$,
$(0, 6, 22, 89)$,\adfsplit
$(1, 15, 29, 32)$,
$(0, 16, 29, 92)$,
$(0, 31, 93, 99)$,
$(0, 7, 21, 69)$,\adfsplit
$(0, 17, 19, 38)$,
$(1, 7, 8, 19)$,
$(0, 46, 49, 98)$,
$(1, 38, 39, 100)$,\adfsplit
$(2, 9, 16, 38)$,
$(1, 9, 17, 107)$,
$(0, 9, 23, 56)$,
$(7, 79, 88, 106)$,\adfsplit
$(0, 53, 79, 102)$,
$(5, 49, 97, 101)$,
$(3, 26, 49, 105)$,
$(0, 18, 41, 101)$,\adfsplit
$(0, 36, 67, 107)$,
$(0, 88, 97, 100)$,
$(0, 48, 81, 103)$,
$(0, 27, 33, 51)$,\adfsplit
$(0, 26, 44, 87)$,
$(0, 57, 66, 104)$,
$(0, 14, 77, 106)$,
$(0, 37, 83, 86)$,\adfsplit
$(1, 3, 67, 85)$,
$(6, 27, 68, 101)$,
$(1, 23, 24, 27)$,
$(3, 54, 65, 87)$,\adfsplit
$(3, 76, 77, 84)$,
$(3, 16, 95, 102)$,
$(3, 66, 94, 107)$,
$(4, 56, 58, 75)$,\adfsplit
$(2, 64, 65, 96)$,
$(2, 23, 25, 76)$,
$(3, 34, 55, 96)$,
$(2, 18, 36, 106)$,\adfsplit
$(1, 22, 65, 102)$,
$(1, 28, 35, 54)$,
$(1, 2, 48, 95)$,
$(2, 8, 75, 107)$,\adfsplit
$(1, 88, 92, 94)$,
$(1, 75, 84, 98)$,
$(1, 14, 45, 104)$,
$(1, 53, 62, 105)$,\adfsplit
$(1, 4, 52, 55)$,
$(2, 3, 15, 104)$,
$(2, 35, 74, 101)$,
$(2, 14, 43, 103)$,\adfsplit
$(2, 83, 94, 100)$

\adfLgap \noindent by the mapping:
$x \mapsto x + 10 j \adfmod{100}$ for $x < 100$,
$x \mapsto (x - 100 + 4 j \adfmod{8}) + 100$ for $100 \le x < 108$,
$x \mapsto (x - 108 + 2 j \adfmod{5}) + 108$ for $x \ge 108$,
$0 \le j < 10$.
\ADFvfyParStart{(113, ((89, 10, ((100, 10), (8, 4), (5, 2)))), ((20, 5), (8, 1), (5, 1)))} 

\section{4-GDDs for the proof of Lemma \ref{lem:4-GDD 14^u m^1}}
\label{app:4-GDD 14^u m^1}
\adfhide{
$ 14^6 8^1 $,
$ 14^6 11^1 $,
$ 14^6 17^1 $,
$ 14^6 20^1 $,
$ 14^6 23^1 $,
$ 14^6 26^1 $,
$ 14^6 29^1 $,
$ 14^6 32^1 $,
$ 14^9 11^1 $,
$ 14^9 17^1 $,
$ 14^9 20^1 $,
$ 14^9 23^1 $,
$ 14^9 26^1 $,
$ 14^9 29^1 $,
$ 14^9 32^1 $,
$ 14^9 38^1 $,
$ 14^9 41^1 $,
$ 14^9 44^1 $,
$ 14^9 47^1 $,
$ 14^9 50^1 $ and
$ 14^9 53^1 $.
}

\adfDgap
\noindent{\boldmath $ 14^{6} 8^{1} $}~
With the point set $Z_{92}$ partitioned into
 residue classes modulo $6$ for $\{0, 1, \dots, 83\}$, and
 $\{84, 85, \dots, 91\}$,
 the design is generated from

\adfLgap 
$(84, 5, 4, 80)$,
$(84, 12, 39, 49)$,
$(85, 68, 72, 39)$,
$(85, 71, 25, 22)$,\adfsplit
$(86, 38, 27, 5)$,
$(86, 61, 46, 0)$,
$(87, 43, 44, 18)$,
$(87, 21, 34, 77)$,\adfsplit
$(88, 37, 36, 70)$,
$(88, 51, 83, 8)$,
$(89, 15, 35, 61)$,
$(89, 70, 80, 18)$,\adfsplit
$(90, 34, 29, 13)$,
$(90, 33, 48, 20)$,
$(91, 50, 11, 64)$,
$(91, 48, 51, 7)$,\adfsplit
$(67, 71, 0, 63)$,
$(1, 63, 28, 30)$,
$(40, 25, 53, 36)$,
$(0, 2, 9, 23)$,\adfsplit
$(0, 13, 15, 53)$,
$(0, 7, 20, 59)$,
$(0, 5, 8, 19)$,
$(0, 16, 35, 45)$,\adfsplit
$(0, 11, 28, 68)$,
$(0, 41, 64, 75)$,
$(1, 21, 23, 44)$,
$(0, 31, 50, 76)$,\adfsplit
$(0, 14, 39, 70)$,
$(0, 10, 57, 83)$,
$(0, 21, 22, 38)$,
$(1, 33, 52, 77)$,\adfsplit
$(0, 32, 47, 49)$,
$(0, 40, 74, 77)$,
$(0, 29, 44, 79)$,
$(0, 58, 65, 81)$,\adfsplit
$(1, 4, 8, 9)$,
$(1, 26, 57, 64)$,
$(1, 10, 27, 62)$,
$(1, 32, 65, 69)$,\adfsplit
$(1, 11, 38, 40)$,
$(1, 15, 58, 80)$,
$(1, 51, 56, 76)$

\adfLgap \noindent by the mapping:
$x \mapsto x + 6 j \adfmod{84}$ for $x < 84$,
$x \mapsto x$ for $x \ge 84$,
$0 \le j < 14$.
\ADFvfyParStart{(92, ((43, 14, ((84, 6), (8, 8)))), ((14, 6), (8, 1)))} 

\adfDgap
\noindent{\boldmath $ 14^{6} 11^{1} $}~
With the point set $Z_{95}$ partitioned into
 residue classes modulo $6$ for $\{0, 1, \dots, 83\}$, and
 $\{84, 85, \dots, 94\}$,
 the design is generated from

\adfLgap 
$(84, 51, 77, 22)$,
$(85, 66, 22, 65)$,
$(86, 31, 24, 35)$,
$(87, 65, 49, 15)$,\adfsplit
$(88, 78, 17, 73)$,
$(89, 28, 65, 60)$,
$(90, 22, 30, 20)$,
$(91, 37, 72, 23)$,\adfsplit
$(0, 1, 59, 62)$,
$(0, 2, 69, 77)$,
$(0, 3, 47, 68)$,
$(0, 4, 21, 37)$,\adfsplit
$(0, 27, 58, 73)$,
$(0, 14, 25, 45)$,
$(0, 10, 19, 33)$,
$(0, 9, 41, 80)$,\adfsplit
$(0, 13, 20, 92)$,
$(0, 22, 43, 56)$,
$(1, 4, 23, 50)$,
$(1, 26, 46, 77)$,\adfsplit
$(1, 11, 28, 80)$,
$(0, 28, 29, 94)$,
$(0, 38, 82, 93)$

\adfLgap \noindent by the mapping:
$x \mapsto x + 3 j \adfmod{84}$ for $x < 84$,
$x \mapsto x$ for $x \ge 84$,
$0 \le j < 28$.
\ADFvfyParStart{(95, ((23, 28, ((84, 3), (11, 11)))), ((14, 6), (11, 1)))} 

\adfDgap
\noindent{\boldmath $ 14^{6} 17^{1} $}~
With the point set $Z_{101}$ partitioned into
 residue classes modulo $6$ for $\{0, 1, \dots, 83\}$, and
 $\{84, 85, \dots, 100\}$,
 the design is generated from

\adfLgap 
$(98, 62, 28, 9)$,
$(99, 13, 72, 20)$,
$(100, 76, 65, 9)$,
$(84, 38, 40, 79)$,\adfsplit
$(84, 48, 35, 69)$,
$(84, 68, 11, 66)$,
$(84, 12, 34, 29)$,
$(84, 19, 10, 51)$,\adfsplit
$(84, 23, 46, 67)$,
$(84, 64, 8, 7)$,
$(84, 1, 56, 45)$,
$(0, 1, 3, 88)$,\adfsplit
$(0, 5, 13, 93)$,
$(0, 44, 47, 84)$,
$(0, 16, 62, 85)$,
$(0, 20, 39, 91)$,\adfsplit
$(0, 38, 64, 86)$,
$(0, 27, 61, 76)$,
$(0, 14, 31, 35)$,
$(0, 10, 41, 80)$,\adfsplit
$(1, 4, 17, 26)$,
$(0, 8, 28, 77)$,
$(0, 15, 73, 83)$,
$(0, 7, 26, 59)$,\adfsplit
$(0, 4, 9, 55)$,
$(0, 23, 33, 70)$

\adfLgap \noindent by the mapping:
$x \mapsto x + 3 j \adfmod{84}$ for $x < 84$,
$x \mapsto (x +  j \adfmod{14}) + 84$ for $84 \le x < 98$,
$x \mapsto x$ for $x \ge 98$,
$0 \le j < 28$.
\ADFvfyParStart{(101, ((26, 28, ((84, 3), (14, 1), (3, 3)))), ((14, 6), (17, 1)))} 

\adfDgap
\noindent{\boldmath $ 14^{6} 20^{1} $}~
With the point set $Z_{104}$ partitioned into
 residue classes modulo $5$ for $\{0, 1, \dots, 69\}$,
 $\{70, 71, \dots, 83\}$, and
 $\{84, 85, \dots, 103\}$,
 the design is generated from

\adfLgap 
$(80, 43, 102, 56)$,
$(4, 71, 58, 37)$,
$(8, 25, 32, 64)$,
$(78, 34, 0, 85)$,\adfsplit
$(90, 19, 28, 22)$,
$(59, 41, 0, 91)$,
$(2, 91, 72, 4)$,
$(28, 75, 2, 29)$,\adfsplit
$(10, 58, 73, 98)$,
$(88, 44, 52, 40)$,
$(85, 45, 3, 22)$

\adfLgap \noindent by the mapping:
$x \mapsto x +  j \adfmod{70}$ for $x < 70$,
$x \mapsto (x +  j \adfmod{14}) + 70$ for $70 \le x < 84$,
$x \mapsto (x - 84 + 2 j \adfmod{20}) + 84$ for $x \ge 84$,
$0 \le j < 70$.
\ADFvfyParStart{(104, ((11, 70, ((70, 1), (14, 1), (20, 2)))), ((14, 5), (14, 1), (20, 1)))} 

\adfDgap
\noindent{\boldmath $ 14^{6} 23^{1} $}~
With the point set $Z_{107}$ partitioned into
 residue classes modulo $6$ for $\{0, 1, \dots, 83\}$, and
 $\{84, 85, \dots, 106\}$,
 the design is generated from

\adfLgap 
$(98, 10, 74, 3)$,
$(99, 71, 66, 31)$,
$(100, 14, 4, 24)$,
$(101, 37, 3, 53)$,\adfsplit
$(102, 41, 12, 49)$,
$(103, 21, 77, 34)$,
$(104, 64, 21, 35)$,
$(105, 54, 47, 70)$,\adfsplit
$(106, 46, 51, 59)$,
$(84, 69, 50, 66)$,
$(84, 54, 53, 22)$,
$(84, 19, 14, 81)$,\adfsplit
$(84, 83, 60, 21)$,
$(84, 0, 33, 44)$,
$(0, 1, 38, 75)$,
$(1, 8, 58, 83)$,\adfsplit
$(1, 53, 68, 95)$,
$(1, 71, 74, 85)$,
$(0, 2, 35, 88)$,
$(0, 20, 59, 63)$,\adfsplit
$(1, 2, 40, 59)$,
$(0, 28, 32, 93)$,
$(1, 29, 50, 88)$,
$(0, 4, 26, 57)$,\adfsplit
$(0, 46, 55, 97)$,
$(0, 61, 76, 86)$,
$(0, 67, 70, 95)$,
$(0, 15, 40, 73)$,\adfsplit
$(0, 19, 82, 96)$

\adfLgap \noindent by the mapping:
$x \mapsto x + 3 j \adfmod{84}$ for $x < 84$,
$x \mapsto (x +  j \adfmod{14}) + 84$ for $84 \le x < 98$,
$x \mapsto x$ for $x \ge 98$,
$0 \le j < 28$.
\ADFvfyParStart{(107, ((29, 28, ((84, 3), (14, 1), (9, 9)))), ((14, 6), (23, 1)))} 

\adfDgap
\noindent{\boldmath $ 14^{6} 26^{1} $}~
With the point set $Z_{110}$ partitioned into
 residue classes modulo $6$ for $\{0, 1, \dots, 83\}$, and
 $\{84, 85, \dots, 109\}$,
 the design is generated from

\adfLgap 
$(105, 52, 2, 79)$,
$(105, 3, 12, 29)$,
$(106, 69, 38, 46)$,
$(106, 83, 31, 24)$,\adfsplit
$(107, 80, 19, 18)$,
$(107, 40, 41, 45)$,
$(108, 20, 23, 78)$,
$(108, 51, 49, 4)$,\adfsplit
$(109, 54, 81, 59)$,
$(109, 40, 37, 38)$,
$(84, 2, 34, 11)$,
$(84, 57, 19, 48)$,\adfsplit
$(84, 31, 65, 51)$,
$(84, 47, 45, 58)$,
$(84, 25, 0, 46)$,
$(84, 24, 82, 38)$,\adfsplit
$(84, 7, 26, 77)$,
$(84, 62, 83, 52)$,
$(84, 78, 27, 14)$,
$(84, 32, 59, 79)$,\adfsplit
$(84, 29, 64, 39)$,
$(84, 30, 13, 28)$,
$(84, 60, 21, 50)$,
$(84, 43, 54, 75)$,\adfsplit
$(85, 50, 45, 28)$,
$(85, 71, 21, 61)$,
$(85, 47, 27, 1)$,
$(85, 14, 53, 12)$,\adfsplit
$(85, 9, 52, 36)$,
$(85, 39, 74, 25)$,
$(85, 59, 22, 26)$,
$(85, 0, 35, 79)$,\adfsplit
$(85, 33, 82, 65)$,
$(85, 4, 24, 62)$,
$(85, 58, 57, 72)$,
$(85, 2, 31, 18)$,\adfsplit
$(85, 7, 80, 34)$,
$(85, 13, 41, 48)$,
$(1, 14, 70, 83)$,
$(1, 9, 32, 77)$,\adfsplit
$(0, 50, 51, 86)$,
$(0, 8, 28, 81)$,
$(1, 44, 81, 92)$,
$(2, 27, 83, 101)$,\adfsplit
$(2, 21, 59, 95)$,
$(0, 56, 63, 71)$,
$(0, 10, 39, 80)$,
$(1, 51, 68, 101)$,\adfsplit
$(0, 19, 44, 98)$,
$(0, 23, 32, 101)$,
$(1, 59, 64, 80)$,
$(1, 69, 76, 104)$,\adfsplit
$(1, 10, 75, 95)$,
$(1, 5, 34, 98)$,
$(0, 22, 61, 95)$,
$(0, 4, 15, 37)$,\adfsplit
$(0, 3, 31, 47)$,
$(0, 40, 43, 65)$,
$(0, 34, 53, 89)$,
$(0, 76, 83, 92)$,\adfsplit
$(0, 11, 52, 104)$

\adfLgap \noindent by the mapping:
$x \mapsto x + 6 j \adfmod{84}$ for $x < 84$,
$x \mapsto (x + 3 j \adfmod{21}) + 84$ for $84 \le x < 105$,
$x \mapsto x$ for $x \ge 105$,
$0 \le j < 14$.
\ADFvfyParStart{(110, ((61, 14, ((84, 6), (21, 3), (5, 5)))), ((14, 6), (26, 1)))} 

\adfDgap
\noindent{\boldmath $ 14^{6} 29^{1} $}~
With the point set $Z_{113}$ partitioned into
 residue classes modulo $6$ for $\{0, 1, \dots, 83\}$, and
 $\{84, 85, \dots, 112\}$,
 the design is generated from

\adfLgap 
$(112, 80, 18, 55)$,
$(84, 44, 31, 30)$,
$(84, 77, 18, 56)$,
$(84, 11, 3, 67)$,\adfsplit
$(84, 65, 4, 78)$,
$(84, 54, 40, 21)$,
$(84, 14, 22, 41)$,
$(84, 35, 64, 50)$,\adfsplit
$(84, 57, 8, 61)$,
$(84, 72, 19, 34)$,
$(84, 39, 36, 32)$,
$(84, 83, 25, 51)$,\adfsplit
$(84, 76, 24, 29)$,
$(84, 82, 55, 27)$,
$(84, 62, 16, 13)$,
$(0, 9, 43, 84)$,\adfsplit
$(0, 7, 16, 98)$,
$(0, 13, 15, 90)$,
$(0, 29, 49, 97)$,
$(0, 25, 26, 76)$,\adfsplit
$(1, 11, 22, 110)$,
$(0, 20, 22, 96)$,
$(0, 2, 53, 89)$,
$(0, 11, 45, 107)$,\adfsplit
$(1, 23, 68, 100)$,
$(0, 23, 67, 91)$,
$(0, 56, 61, 87)$,
$(1, 8, 17, 88)$,\adfsplit
$(0, 17, 27, 110)$,
$(0, 21, 65, 68)$,
$(0, 40, 79, 83)$,
$(0, 41, 73, 108)$

\adfLgap \noindent by the mapping:
$x \mapsto x + 3 j \adfmod{84}$ for $x < 84$,
$x \mapsto (x +  j \adfmod{28}) + 84$ for $84 \le x < 112$,
$112 \mapsto 112$,
$0 \le j < 28$.
\ADFvfyParStart{(113, ((32, 28, ((84, 3), (28, 1), (1, 1)))), ((14, 6), (29, 1)))} 

\adfDgap
\noindent{\boldmath $ 14^{6} 32^{1} $}~
With the point set $Z_{116}$ partitioned into
 residue classes modulo $6$ for $\{0, 1, \dots, 83\}$, and
 $\{84, 85, \dots, 115\}$,
 the design is generated from

\adfLgap 
$(112, 42, 41, 73)$,
$(112, 56, 45, 64)$,
$(113, 66, 74, 55)$,
$(113, 10, 69, 5)$,\adfsplit
$(114, 29, 0, 49)$,
$(114, 20, 40, 45)$,
$(115, 39, 12, 7)$,
$(115, 50, 83, 64)$,\adfsplit
$(84, 16, 77, 55)$,
$(84, 41, 10, 25)$,
$(84, 17, 12, 21)$,
$(84, 50, 60, 33)$,\adfsplit
$(84, 43, 68, 64)$,
$(84, 65, 38, 61)$,
$(84, 31, 9, 6)$,
$(84, 4, 66, 62)$,\adfsplit
$(84, 56, 27, 7)$,
$(84, 72, 3, 74)$,
$(84, 79, 29, 70)$,
$(84, 39, 78, 11)$,\adfsplit
$(84, 76, 42, 2)$,
$(84, 40, 5, 57)$,
$(85, 22, 63, 56)$,
$(85, 44, 17, 28)$,\adfsplit
$(85, 13, 57, 11)$,
$(85, 58, 1, 3)$,
$(85, 9, 10, 24)$,
$(85, 48, 31, 68)$,\adfsplit
$(85, 79, 76, 60)$,
$(85, 25, 35, 38)$,
$(85, 47, 61, 69)$,
$(85, 36, 29, 20)$,\adfsplit
$(85, 42, 4, 49)$,
$(85, 33, 74, 12)$,
$(85, 40, 30, 41)$,
$(85, 65, 50, 81)$,\adfsplit
$(86, 73, 30, 69)$,
$(86, 10, 50, 18)$,
$(86, 2, 71, 43)$,
$(86, 48, 20, 23)$,\adfsplit
$(86, 59, 51, 14)$,
$(86, 75, 37, 4)$,
$(86, 26, 19, 3)$,
$(0, 50, 71, 98)$,\adfsplit
$(3, 64, 77, 106)$,
$(0, 55, 64, 86)$,
$(0, 4, 41, 90)$,
$(1, 34, 56, 98)$,\adfsplit
$(0, 37, 51, 106)$,
$(3, 5, 34, 90)$,
$(0, 23, 81, 87)$,
$(1, 57, 64, 103)$,\adfsplit
$(1, 59, 68, 95)$,
$(0, 53, 61, 99)$,
$(1, 16, 27, 41)$,
$(2, 3, 53, 87)$,\adfsplit
$(0, 1, 47, 107)$,
$(1, 28, 74, 91)$,
$(1, 75, 80, 99)$,
$(0, 38, 40, 91)$,\adfsplit
$(2, 21, 65, 111)$,
$(0, 13, 14, 63)$,
$(1, 4, 32, 111)$,
$(0, 35, 52, 111)$,\adfsplit
$(0, 28, 75, 103)$,
$(0, 26, 58, 65)$,
$(0, 33, 82, 95)$

\adfLgap \noindent by the mapping:
$x \mapsto x + 6 j \adfmod{84}$ for $x < 84$,
$x \mapsto (x + 4 j \adfmod{28}) + 84$ for $84 \le x < 112$,
$x \mapsto x$ for $x \ge 112$,
$0 \le j < 14$.
\ADFvfyParStart{(116, ((67, 14, ((84, 6), (28, 4), (4, 4)))), ((14, 6), (32, 1)))} 

\adfDgap
\noindent{\boldmath $ 14^{9} 11^{1} $}~
With the point set $Z_{137}$ partitioned into
 residue classes modulo $9$ for $\{0, 1, \dots, 125\}$, and
 $\{126, 127, \dots, 136\}$,
 the design is generated from

\adfLgap 
$(126, 76, 99, 32)$,
$(126, 7, 36, 95)$,
$(127, 102, 119, 14)$,
$(127, 40, 105, 115)$,\adfsplit
$(128, 95, 43, 93)$,
$(128, 92, 88, 42)$,
$(129, 8, 42, 27)$,
$(129, 88, 11, 103)$,\adfsplit
$(130, 107, 24, 92)$,
$(130, 118, 99, 103)$,
$(131, 56, 96, 88)$,
$(131, 105, 125, 43)$,\adfsplit
$(132, 104, 90, 55)$,
$(132, 29, 52, 93)$,
$(133, 19, 58, 86)$,
$(133, 24, 57, 89)$,\adfsplit
$(134, 98, 97, 15)$,
$(134, 54, 107, 34)$,
$(135, 123, 43, 50)$,
$(135, 29, 58, 6)$,\adfsplit
$(136, 38, 41, 30)$,
$(136, 21, 115, 46)$,
$(19, 92, 30, 97)$,
$(125, 36, 67, 38)$,\adfsplit
$(46, 58, 32, 102)$,
$(106, 111, 72, 51)$,
$(72, 65, 112, 79)$,
$(70, 27, 92, 1)$,\adfsplit
$(31, 108, 10, 114)$,
$(93, 82, 79, 9)$,
$(1, 95, 90, 65)$,
$(94, 84, 34, 27)$,\adfsplit
$(108, 39, 70, 8)$,
$(88, 80, 10, 101)$,
$(3, 114, 20, 117)$,
$(56, 13, 122, 53)$,\adfsplit
$(19, 111, 117, 32)$,
$(0, 1, 13, 21)$,
$(0, 19, 24, 44)$,
$(0, 4, 66, 114)$,\adfsplit
$(0, 25, 42, 71)$,
$(0, 30, 77, 110)$,
$(0, 73, 79, 95)$,
$(0, 107, 112, 113)$,\adfsplit
$(0, 103, 122, 124)$,
$(0, 58, 61, 116)$,
$(0, 93, 98, 104)$,
$(0, 35, 55, 85)$,\adfsplit
$(0, 74, 94, 123)$,
$(0, 87, 100, 125)$,
$(0, 41, 51, 75)$,
$(1, 5, 61, 103)$,\adfsplit
$(1, 32, 80, 87)$,
$(1, 76, 86, 116)$,
$(1, 38, 62, 75)$,
$(1, 88, 104, 105)$,\adfsplit
$(1, 39, 69, 119)$,
$(1, 11, 111, 125)$,
$(1, 3, 94, 101)$,
$(3, 40, 59, 107)$,\adfsplit
$(3, 4, 47, 51)$,
$(3, 11, 71, 113)$,
$(2, 35, 52, 86)$,
$(2, 27, 76, 113)$,\adfsplit
$(2, 40, 71, 82)$,
$(2, 14, 53, 105)$,
$(4, 10, 65, 106)$

\adfLgap \noindent by the mapping:
$x \mapsto x + 6 j \adfmod{126}$ for $x < 126$,
$x \mapsto x$ for $x \ge 126$,
$0 \le j < 21$.
\ADFvfyParStart{(137, ((67, 21, ((126, 6), (11, 11)))), ((14, 9), (11, 1)))} 

\adfDgap
\noindent{\boldmath $ 14^{9} 17^{1} $}~
With the point set $Z_{143}$ partitioned into
 residue classes modulo $9$ for $\{0, 1, \dots, 125\}$, and
 $\{126, 127, \dots, 142\}$,
 the design is generated from

\adfLgap 
$(126, 65, 62, 33)$,
$(126, 112, 12, 79)$,
$(127, 95, 57, 26)$,
$(127, 16, 37, 0)$,\adfsplit
$(128, 85, 45, 78)$,
$(128, 95, 4, 20)$,
$(129, 31, 14, 16)$,
$(129, 96, 27, 47)$,\adfsplit
$(130, 54, 76, 14)$,
$(130, 71, 69, 31)$,
$(131, 115, 60, 113)$,
$(131, 100, 105, 110)$,\adfsplit
$(132, 68, 46, 54)$,
$(132, 117, 91, 71)$,
$(133, 95, 105, 56)$,
$(133, 4, 72, 91)$,\adfsplit
$(134, 88, 74, 95)$,
$(134, 81, 67, 120)$,
$(135, 107, 18, 93)$,
$(135, 13, 104, 34)$,\adfsplit
$(136, 53, 27, 42)$,
$(136, 68, 1, 94)$,
$(137, 3, 25, 76)$,
$(137, 20, 12, 35)$,\adfsplit
$(138, 12, 116, 95)$,
$(138, 91, 52, 99)$,
$(139, 74, 112, 0)$,
$(139, 25, 27, 83)$,\adfsplit
$(140, 75, 88, 35)$,
$(140, 102, 104, 79)$,
$(141, 33, 29, 43)$,
$(141, 42, 106, 14)$,\adfsplit
$(142, 115, 9, 64)$,
$(142, 98, 108, 23)$,
$(100, 40, 57, 9)$,
$(105, 36, 16, 66)$,\adfsplit
$(58, 46, 98, 39)$,
$(9, 67, 71, 113)$,
$(94, 62, 25, 95)$,
$(56, 0, 80, 60)$,\adfsplit
$(104, 44, 112, 92)$,
$(19, 116, 117, 83)$,
$(46, 43, 35, 40)$,
$(0, 4, 82, 105)$,\adfsplit
$(0, 3, 6, 52)$,
$(0, 10, 34, 38)$,
$(0, 29, 44, 94)$,
$(1, 16, 62, 85)$,\adfsplit
$(2, 46, 59, 111)$,
$(3, 11, 70, 100)$,
$(0, 40, 51, 119)$,
$(3, 4, 47, 88)$,\adfsplit
$(0, 28, 59, 84)$,
$(0, 17, 65, 88)$,
$(4, 23, 29, 53)$,
$(0, 95, 107, 124)$,\adfsplit
$(0, 24, 50, 92)$,
$(1, 13, 29, 32)$,
$(0, 47, 91, 113)$,
$(0, 5, 33, 79)$,\adfsplit
$(0, 1, 25, 35)$,
$(1, 2, 89, 121)$,
$(0, 15, 49, 125)$,
$(0, 21, 71, 97)$,\adfsplit
$(0, 32, 85, 115)$,
$(0, 12, 43, 121)$,
$(0, 13, 62, 78)$,
$(1, 14, 57, 61)$,\adfsplit
$(1, 33, 75, 122)$,
$(2, 8, 15, 27)$,
$(1, 45, 56, 111)$,
$(2, 63, 87, 93)$,\adfsplit
$(1, 63, 86, 116)$

\adfLgap \noindent by the mapping:
$x \mapsto x + 6 j \adfmod{126}$ for $x < 126$,
$x \mapsto x$ for $x \ge 126$,
$0 \le j < 21$.
\ADFvfyParStart{(143, ((73, 21, ((126, 6), (17, 17)))), ((14, 9), (17, 1)))} 

\adfDgap
\noindent{\boldmath $ 14^{9} 20^{1} $}~
With the point set $Z_{146}$ partitioned into
 residue classes modulo $9$ for $\{0, 1, \dots, 125\}$, and
 $\{126, 127, \dots, 145\}$,
 the design is generated from

\adfLgap 
$(140, 106, 92, 72)$,
$(141, 92, 75, 28)$,
$(142, 96, 119, 76)$,
$(143, 40, 36, 116)$,\adfsplit
$(144, 7, 102, 32)$,
$(145, 77, 79, 117)$,
$(126, 115, 5, 60)$,
$(126, 64, 79, 112)$,\adfsplit
$(126, 91, 110, 125)$,
$(126, 74, 85, 109)$,
$(126, 40, 96, 45)$,
$(126, 20, 111, 53)$,\adfsplit
$(126, 35, 4, 84)$,
$(126, 50, 72, 24)$,
$(126, 76, 97, 59)$,
$(126, 39, 99, 2)$,\adfsplit
$(126, 56, 9, 100)$,
$(126, 65, 78, 75)$,
$(0, 8, 73, 133)$,
$(0, 74, 97, 138)$,\adfsplit
$(0, 1, 5, 107)$,
$(0, 2, 28, 122)$,
$(0, 6, 59, 101)$,
$(0, 7, 62, 65)$,\adfsplit
$(0, 14, 44, 119)$,
$(0, 13, 41, 98)$,
$(1, 38, 50, 98)$,
$(1, 14, 53, 61)$,\adfsplit
$(1, 8, 13, 88)$,
$(0, 42, 103, 125)$,
$(0, 11, 64, 105)$,
$(0, 40, 43, 110)$,\adfsplit
$(0, 12, 49, 50)$,
$(0, 52, 82, 92)$,
$(0, 19, 76, 87)$,
$(0, 10, 94, 111)$,\adfsplit
$(0, 24, 57, 124)$,
$(0, 16, 22, 30)$

\adfLgap \noindent by the mapping:
$x \mapsto x + 3 j \adfmod{126}$ for $x < 126$,
$x \mapsto (x +  j \adfmod{14}) + 126$ for $126 \le x < 140$,
$x \mapsto x$ for $x \ge 140$,
$0 \le j < 42$.
\ADFvfyParStart{(146, ((38, 42, ((126, 3), (14, 1), (6, 6)))), ((14, 9), (20, 1)))} 

\adfDgap
\noindent{\boldmath $ 14^{9} 23^{1} $}~
With the point set $Z_{149}$ partitioned into
 residue classes modulo $9$ for $\{0, 1, \dots, 125\}$, and
 $\{126, 127, \dots, 148\}$,
 the design is generated from

\adfLgap 
$(147, 56, 93, 72)$,
$(147, 11, 118, 79)$,
$(148, 24, 70, 41)$,
$(148, 13, 104, 57)$,\adfsplit
$(126, 78, 90, 97)$,
$(126, 20, 108, 52)$,
$(127, 65, 25, 75)$,
$(127, 92, 16, 107)$,\adfsplit
$(128, 93, 121, 29)$,
$(128, 122, 6, 35)$,
$(129, 114, 30, 22)$,
$(129, 88, 112, 15)$,\adfsplit
$(130, 87, 80, 30)$,
$(130, 35, 51, 111)$,
$(131, 11, 25, 54)$,
$(131, 76, 114, 61)$,\adfsplit
$(132, 51, 108, 93)$,
$(132, 119, 116, 41)$,
$(133, 83, 104, 1)$,
$(133, 89, 29, 30)$,\adfsplit
$(134, 57, 40, 55)$,
$(134, 81, 77, 33)$,
$(135, 114, 62, 21)$,
$(135, 121, 27, 30)$,\adfsplit
$(136, 62, 13, 86)$,
$(136, 22, 99, 59)$,
$(137, 17, 66, 22)$,
$(137, 91, 52, 57)$,\adfsplit
$(138, 65, 1, 102)$,
$(138, 112, 106, 93)$,
$(139, 113, 58, 20)$,
$(139, 115, 5, 74)$,\adfsplit
$(140, 2, 23, 73)$,
$(140, 32, 58, 100)$,
$(141, 13, 88, 37)$,
$(141, 9, 62, 122)$,\adfsplit
$(142, 109, 24, 99)$,
$(142, 75, 36, 41)$,
$(143, 103, 14, 124)$,
$(143, 12, 65, 116)$,\adfsplit
$(144, 40, 33, 30)$,
$(144, 24, 25, 86)$,
$(145, 52, 63, 80)$,
$(145, 83, 91, 79)$,\adfsplit
$(146, 110, 12, 70)$,
$(146, 55, 22, 44)$,
$(106, 59, 51, 90)$,
$(72, 70, 120, 14)$,\adfsplit
$(0, 2, 64, 94)$,
$(0, 4, 30, 52)$,
$(0, 8, 13, 102)$,
$(0, 11, 14, 28)$,\adfsplit
$(0, 15, 31, 41)$,
$(0, 6, 107, 109)$,
$(0, 49, 55, 60)$,
$(0, 26, 51, 56)$,\adfsplit
$(0, 79, 80, 122)$,
$(0, 43, 112, 113)$,
$(0, 44, 61, 92)$,
$(0, 40, 65, 71)$,\adfsplit
$(0, 23, 106, 119)$,
$(0, 67, 86, 105)$,
$(0, 35, 47, 93)$,
$(1, 31, 87, 98)$,\adfsplit
$(1, 21, 43, 95)$,
$(1, 29, 49, 106)$,
$(1, 26, 27, 123)$,
$(1, 63, 101, 124)$,\adfsplit
$(1, 4, 8, 14)$,
$(2, 33, 35, 76)$,
$(2, 57, 89, 113)$,
$(2, 10, 51, 116)$,\adfsplit
$(1, 23, 34, 80)$,
$(2, 45, 46, 69)$,
$(3, 17, 70, 82)$,
$(1, 47, 69, 89)$,\adfsplit
$(2, 4, 93, 105)$,
$(1, 9, 15, 61)$,
$(3, 34, 94, 101)$

\adfLgap \noindent by the mapping:
$x \mapsto x + 6 j \adfmod{126}$ for $x < 126$,
$x \mapsto (x +  j \adfmod{21}) + 126$ for $126 \le x < 147$,
$x \mapsto x$ for $x \ge 147$,
$0 \le j < 21$.
\ADFvfyParStart{(149, ((79, 21, ((126, 6), (21, 1), (2, 2)))), ((14, 9), (23, 1)))} 

\adfDgap
\noindent{\boldmath $ 14^{9} 26^{1} $}~
With the point set $Z_{152}$ partitioned into
 residue classes modulo $9$ for $\{0, 1, \dots, 125\}$, and
 $\{126, 127, \dots, 151\}$,
 the design is generated from

\adfLgap 
$(140, 111, 19, 74)$,
$(141, 91, 32, 105)$,
$(142, 37, 68, 45)$,
$(143, 120, 16, 44)$,\adfsplit
$(144, 46, 75, 14)$,
$(145, 69, 11, 67)$,
$(146, 37, 111, 71)$,
$(147, 60, 5, 10)$,\adfsplit
$(148, 76, 77, 120)$,
$(149, 7, 86, 12)$,
$(150, 18, 61, 59)$,
$(151, 62, 78, 43)$,\adfsplit
$(126, 112, 42, 75)$,
$(126, 104, 26, 100)$,
$(126, 47, 57, 59)$,
$(126, 98, 88, 74)$,\adfsplit
$(126, 37, 27, 52)$,
$(126, 11, 44, 123)$,
$(126, 25, 72, 93)$,
$(126, 125, 45, 85)$,\adfsplit
$(126, 40, 34, 105)$,
$(0, 1, 38, 136)$,
$(0, 4, 102, 134)$,
$(0, 5, 11, 132)$,\adfsplit
$(0, 13, 59, 138)$,
$(1, 8, 61, 130)$,
$(0, 7, 66, 85)$,
$(0, 16, 39, 96)$,\adfsplit
$(0, 49, 100, 125)$,
$(1, 13, 34, 98)$,
$(1, 14, 17, 25)$,
$(1, 23, 58, 107)$,\adfsplit
$(0, 26, 64, 94)$,
$(1, 4, 43, 104)$,
$(0, 15, 44, 88)$,
$(0, 32, 98, 115)$,\adfsplit
$(0, 3, 78, 109)$,
$(0, 6, 62, 101)$,
$(0, 17, 92, 113)$,
$(0, 8, 12, 119)$,\adfsplit
$(0, 20, 42, 77)$

\adfLgap \noindent by the mapping:
$x \mapsto x + 3 j \adfmod{126}$ for $x < 126$,
$x \mapsto (x +  j \adfmod{14}) + 126$ for $126 \le x < 140$,
$x \mapsto x$ for $x \ge 140$,
$0 \le j < 42$.
\ADFvfyParStart{(152, ((41, 42, ((126, 3), (14, 1), (12, 12)))), ((14, 9), (26, 1)))} 

\adfDgap
\noindent{\boldmath $ 14^{9} 29^{1} $}~
With the point set $Z_{155}$ partitioned into
 residue classes modulo $9$ for $\{0, 1, \dots, 125\}$, and
 $\{126, 127, \dots, 154\}$,
 the design is generated from

\adfLgap 
$(126, 58, 48, 19)$,
$(126, 47, 3, 86)$,
$(127, 11, 85, 14)$,
$(127, 111, 78, 40)$,\adfsplit
$(128, 49, 10, 116)$,
$(128, 42, 45, 53)$,
$(129, 123, 19, 102)$,
$(129, 32, 11, 58)$,\adfsplit
$(130, 96, 116, 34)$,
$(130, 63, 19, 89)$,
$(131, 104, 7, 47)$,
$(131, 108, 58, 123)$,\adfsplit
$(132, 35, 50, 16)$,
$(132, 66, 9, 49)$,
$(133, 44, 123, 112)$,
$(133, 103, 41, 54)$,\adfsplit
$(134, 99, 35, 104)$,
$(134, 73, 48, 22)$,
$(135, 17, 100, 96)$,
$(135, 115, 2, 105)$,\adfsplit
$(136, 58, 59, 20)$,
$(136, 78, 91, 27)$,
$(137, 89, 66, 56)$,
$(137, 87, 106, 49)$,\adfsplit
$(138, 9, 24, 118)$,
$(138, 110, 7, 35)$,
$(139, 113, 55, 88)$,
$(139, 51, 116, 84)$,\adfsplit
$(140, 94, 109, 77)$,
$(140, 0, 39, 8)$,
$(141, 67, 6, 101)$,
$(141, 74, 76, 45)$,\adfsplit
$(142, 28, 32, 120)$,
$(142, 93, 125, 121)$,
$(143, 85, 59, 44)$,
$(143, 52, 120, 51)$,\adfsplit
$(144, 60, 4, 65)$,
$(144, 15, 61, 86)$,
$(145, 99, 65, 43)$,
$(145, 14, 36, 28)$,\adfsplit
$(146, 24, 80, 4)$,
$(146, 99, 23, 7)$,
$(147, 17, 90, 50)$,
$(147, 9, 109, 34)$,\adfsplit
$(148, 92, 39, 42)$,
$(148, 17, 10, 7)$,
$(149, 1, 20, 24)$,
$(149, 59, 46, 3)$,\adfsplit
$(150, 23, 28, 72)$,
$(150, 80, 43, 93)$,
$(151, 90, 32, 69)$,
$(151, 10, 43, 125)$,\adfsplit
$(152, 100, 114, 103)$,
$(152, 27, 107, 86)$,
$(153, 97, 39, 50)$,
$(0, 16, 65, 153)$,\adfsplit
$(0, 1, 7, 31)$,
$(0, 6, 84, 98)$,
$(0, 12, 71, 85)$,
$(0, 30, 110, 121)$,\adfsplit
$(0, 29, 66, 107)$,
$(0, 2, 55, 102)$,
$(0, 35, 87, 101)$,
$(1, 33, 50, 85)$,\adfsplit
$(1, 13, 44, 61)$,
$(0, 19, 40, 154)$,
$(0, 17, 37, 52)$,
$(0, 22, 28, 125)$,\adfsplit
$(0, 46, 83, 124)$,
$(0, 44, 51, 119)$,
$(0, 62, 67, 74)$,
$(1, 47, 77, 89)$,\adfsplit
$(1, 2, 21, 62)$,
$(2, 3, 5, 154)$,
$(3, 70, 101, 125)$,
$(1, 70, 111, 123)$,\adfsplit
$(1, 3, 9, 106)$,
$(1, 15, 119, 125)$,
$(3, 16, 51, 89)$,
$(3, 23, 33, 82)$,\adfsplit
$(2, 27, 51, 93)$,
$(3, 10, 40, 124)$,
$(2, 34, 44, 100)$,
$(2, 10, 50, 112)$,\adfsplit
$(2, 8, 82, 104)$

\adfLgap \noindent by the mapping:
$x \mapsto x + 6 j \adfmod{126}$ for $x < 126$,
$x \mapsto x$ for $x \ge 126$,
$0 \le j < 21$.
\ADFvfyParStart{(155, ((85, 21, ((126, 6), (29, 29)))), ((14, 9), (29, 1)))} 

\adfDgap
\noindent{\boldmath $ 14^{9} 32^{1} $}~
With the point set $Z_{158}$ partitioned into
 residue classes modulo $9$ for $\{0, 1, \dots, 125\}$, and
 $\{126, 127, \dots, 157\}$,
 the design is generated from

\adfLgap 
$(140, 30, 83, 76)$,
$(141, 42, 76, 5)$,
$(142, 114, 113, 46)$,
$(143, 65, 75, 70)$,\adfsplit
$(144, 88, 59, 105)$,
$(145, 18, 91, 23)$,
$(146, 85, 101, 45)$,
$(147, 16, 89, 30)$,\adfsplit
$(148, 39, 118, 56)$,
$(149, 53, 21, 28)$,
$(150, 103, 113, 90)$,
$(151, 95, 87, 112)$,\adfsplit
$(152, 94, 9, 95)$,
$(153, 46, 86, 48)$,
$(154, 92, 79, 90)$,
$(155, 114, 83, 118)$,\adfsplit
$(156, 64, 56, 63)$,
$(157, 63, 77, 118)$,
$(126, 124, 27, 93)$,
$(126, 18, 95, 29)$,\adfsplit
$(126, 31, 66, 46)$,
$(126, 83, 22, 48)$,
$(126, 122, 39, 16)$,
$(126, 25, 85, 74)$,\adfsplit
$(0, 6, 28, 67)$,
$(0, 16, 19, 93)$,
$(0, 20, 43, 64)$,
$(0, 70, 82, 104)$,\adfsplit
$(1, 5, 29, 31)$,
$(0, 10, 88, 94)$,
$(1, 38, 52, 95)$,
$(1, 47, 77, 89)$,\adfsplit
$(1, 25, 94, 137)$,
$(1, 20, 71, 138)$,
$(0, 44, 118, 139)$,
$(1, 32, 80, 136)$,\adfsplit
$(0, 37, 87, 129)$,
$(0, 21, 78, 133)$,
$(2, 8, 41, 135)$,
$(0, 24, 74, 130)$,\adfsplit
$(0, 47, 62, 75)$,
$(0, 30, 101, 122)$,
$(0, 42, 107, 110)$,
$(0, 3, 29, 114)$

\adfLgap \noindent by the mapping:
$x \mapsto x + 3 j \adfmod{126}$ for $x < 126$,
$x \mapsto (x +  j \adfmod{14}) + 126$ for $126 \le x < 140$,
$x \mapsto x$ for $x \ge 140$,
$0 \le j < 42$.
\ADFvfyParStart{(158, ((44, 42, ((126, 3), (14, 1), (18, 18)))), ((14, 9), (32, 1)))} 

\adfDgap
\noindent{\boldmath $ 14^{9} 38^{1} $}~
With the point set $Z_{164}$ partitioned into
 residue classes modulo $9$ for $\{0, 1, \dots, 125\}$, and
 $\{126, 127, \dots, 163\}$,
 the design is generated from

\adfLgap 
$(147, 1, 107, 33)$,
$(148, 29, 87, 52)$,
$(149, 22, 110, 90)$,
$(150, 97, 92, 123)$,\adfsplit
$(151, 15, 76, 125)$,
$(152, 117, 122, 112)$,
$(153, 64, 87, 65)$,
$(154, 17, 36, 124)$,\adfsplit
$(155, 20, 31, 114)$,
$(156, 55, 15, 23)$,
$(157, 125, 73, 87)$,
$(158, 15, 43, 14)$,\adfsplit
$(159, 67, 30, 74)$,
$(160, 112, 104, 111)$,
$(161, 35, 81, 100)$,
$(162, 123, 32, 115)$,\adfsplit
$(163, 81, 22, 5)$,
$(126, 61, 63, 26)$,
$(126, 125, 114, 104)$,
$(126, 115, 103, 55)$,\adfsplit
$(126, 81, 47, 77)$,
$(126, 3, 74, 5)$,
$(126, 39, 92, 95)$,
$(126, 112, 65, 91)$,\adfsplit
$(126, 20, 64, 88)$,
$(126, 71, 57, 45)$,
$(126, 19, 69, 58)$,
$(0, 7, 29, 62)$,\adfsplit
$(0, 13, 47, 98)$,
$(1, 14, 38, 125)$,
$(1, 29, 41, 142)$,
$(2, 8, 50, 130)$,\adfsplit
$(0, 10, 16, 41)$,
$(0, 17, 83, 97)$,
$(0, 22, 86, 143)$,
$(1, 17, 70, 145)$,\adfsplit
$(0, 73, 77, 135)$,
$(0, 48, 113, 137)$,
$(0, 25, 42, 101)$,
$(0, 23, 79, 82)$,\adfsplit
$(0, 4, 34, 136)$,
$(1, 34, 76, 144)$,
$(0, 6, 55, 70)$,
$(0, 15, 46, 66)$,\adfsplit
$(0, 33, 85, 144)$,
$(0, 3, 24, 139)$,
$(0, 30, 69, 129)$

\adfLgap \noindent by the mapping:
$x \mapsto x + 3 j \adfmod{126}$ for $x < 126$,
$x \mapsto (x +  j \adfmod{21}) + 126$ for $126 \le x < 147$,
$x \mapsto x$ for $x \ge 147$,
$0 \le j < 42$.
\ADFvfyParStart{(164, ((47, 42, ((126, 3), (21, 1), (17, 17)))), ((14, 9), (38, 1)))} 

\adfDgap
\noindent{\boldmath $ 14^{9} 41^{1} $}~
With the point set $Z_{167}$ partitioned into
 residue classes modulo $9$ for $\{0, 1, \dots, 125\}$, and
 $\{126, 127, \dots, 166\}$,
 the design is generated from

\adfLgap 
$(165, 32, 34, 119)$,
$(165, 117, 24, 67)$,
$(166, 12, 92, 45)$,
$(166, 19, 106, 89)$,\adfsplit
$(126, 13, 66, 38)$,
$(126, 34, 108, 64)$,
$(126, 25, 33, 73)$,
$(126, 86, 112, 119)$,\adfsplit
$(126, 26, 77, 21)$,
$(126, 45, 35, 60)$,
$(127, 89, 22, 101)$,
$(127, 111, 69, 44)$,\adfsplit
$(127, 79, 96, 73)$,
$(127, 13, 28, 74)$,
$(127, 41, 27, 16)$,
$(127, 84, 54, 50)$,\adfsplit
$(128, 53, 110, 7)$,
$(128, 60, 40, 3)$,
$(128, 52, 116, 123)$,
$(128, 10, 12, 101)$,\adfsplit
$(128, 49, 0, 41)$,
$(128, 37, 27, 122)$,
$(129, 121, 52, 122)$,
$(129, 118, 41, 18)$,\adfsplit
$(129, 105, 101, 91)$,
$(129, 26, 74, 45)$,
$(129, 75, 58, 96)$,
$(129, 89, 30, 25)$,\adfsplit
$(130, 89, 110, 21)$,
$(130, 5, 51, 103)$,
$(130, 120, 0, 16)$,
$(130, 64, 98, 76)$,\adfsplit
$(130, 83, 63, 86)$,
$(130, 25, 24, 55)$,
$(131, 74, 94, 34)$,
$(131, 95, 18, 73)$,\adfsplit
$(131, 6, 66, 27)$,
$(131, 44, 33, 115)$,
$(131, 10, 13, 125)$,
$(131, 11, 122, 39)$,\adfsplit
$(132, 7, 120, 9)$,
$(132, 29, 124, 118)$,
$(132, 17, 6, 108)$,
$(132, 105, 76, 44)$,\adfsplit
$(132, 77, 110, 93)$,
$(132, 122, 109, 49)$,
$(133, 120, 33, 65)$,
$(133, 50, 97, 85)$,\adfsplit
$(133, 89, 63, 76)$,
$(133, 124, 26, 39)$,
$(133, 24, 72, 109)$,
$(133, 82, 92, 77)$,\adfsplit
$(134, 51, 112, 29)$,
$(134, 118, 99, 48)$,
$(134, 107, 78, 90)$,
$(134, 106, 62, 68)$,\adfsplit
$(134, 23, 115, 73)$,
$(134, 67, 93, 38)$,
$(135, 31, 47, 70)$,
$(135, 96, 64, 115)$,\adfsplit
$(135, 94, 117, 15)$,
$(135, 113, 0, 107)$,
$(135, 104, 74, 48)$,
$(135, 37, 93, 44)$,\adfsplit
$(0, 53, 95, 119)$,
$(1, 4, 5, 101)$,
$(0, 5, 25, 83)$,
$(2, 3, 41, 136)$,\adfsplit
$(0, 50, 64, 125)$,
$(0, 57, 65, 137)$,
$(0, 47, 91, 162)$,
$(1, 33, 93, 95)$,\adfsplit
$(3, 9, 34, 53)$,
$(1, 89, 123, 138)$,
$(3, 15, 89, 164)$,
$(1, 50, 53, 149)$,\adfsplit
$(1, 41, 112, 164)$,
$(2, 23, 76, 150)$,
$(1, 68, 125, 163)$,
$(1, 20, 25, 105)$,\adfsplit
$(1, 39, 69, 150)$,
$(1, 34, 38, 111)$,
$(1, 21, 94, 99)$,
$(3, 70, 94, 162)$,\adfsplit
$(0, 3, 10, 67)$,
$(0, 32, 123, 136)$,
$(0, 75, 76, 118)$,
$(0, 34, 112, 163)$,\adfsplit
$(0, 7, 28, 164)$,
$(0, 40, 61, 149)$,
$(0, 20, 104, 115)$,
$(0, 14, 97, 150)$,\adfsplit
$(0, 4, 62, 79)$,
$(0, 38, 46, 151)$,
$(0, 8, 58, 74)$,
$(0, 2, 116, 138)$,\adfsplit
$(0, 42, 86, 110)$

\adfLgap \noindent by the mapping:
$x \mapsto x + 6 j \adfmod{126}$ for $x < 126$,
$x \mapsto (x - 126 + 13 j \adfmod{39}) + 126$ for $126 \le x < 165$,
$x \mapsto x$ for $x \ge 165$,
$0 \le j < 21$.
\ADFvfyParStart{(167, ((97, 21, ((126, 6), (39, 13), (2, 2)))), ((14, 9), (41, 1)))} 

\adfDgap
\noindent{\boldmath $ 14^{9} 44^{1} $}~
With the point set $Z_{170}$ partitioned into
 residue classes modulo $9$ for $\{0, 1, \dots, 125\}$, and
 $\{126, 127, \dots, 169\}$,
 the design is generated from

\adfLgap 
$(168, 57, 43, 29)$,
$(169, 57, 67, 14)$,
$(126, 117, 22, 28)$,
$(126, 26, 72, 120)$,\adfsplit
$(126, 24, 124, 35)$,
$(126, 51, 30, 68)$,
$(126, 98, 122, 54)$,
$(126, 73, 83, 123)$,\adfsplit
$(126, 3, 44, 96)$,
$(126, 89, 106, 50)$,
$(126, 84, 90, 74)$,
$(126, 103, 107, 2)$,\adfsplit
$(126, 100, 116, 85)$,
$(126, 21, 110, 81)$,
$(126, 4, 43, 17)$,
$(126, 58, 63, 60)$,\adfsplit
$(126, 99, 38, 79)$,
$(126, 29, 41, 88)$,
$(126, 56, 59, 94)$,
$(126, 113, 0, 52)$,\adfsplit
$(126, 20, 12, 69)$,
$(126, 93, 16, 108)$,
$(126, 19, 11, 36)$,
$(126, 114, 119, 13)$,\adfsplit
$(0, 1, 35, 95)$,
$(0, 2, 50, 153)$,
$(0, 4, 23, 107)$,
$(0, 16, 62, 92)$,\adfsplit
$(0, 7, 59, 154)$,
$(0, 14, 125, 142)$,
$(1, 25, 83, 146)$,
$(0, 97, 119, 163)$,\adfsplit
$(0, 94, 122, 146)$,
$(2, 8, 77, 147)$,
$(0, 20, 53, 82)$,
$(1, 50, 61, 156)$,\adfsplit
$(1, 8, 79, 129)$,
$(0, 26, 28, 159)$,
$(0, 56, 61, 143)$,
$(0, 64, 104, 134)$,\adfsplit
$(0, 24, 71, 157)$,
$(1, 2, 85, 158)$,
$(0, 30, 73, 165)$,
$(0, 40, 115, 139)$,\adfsplit
$(0, 22, 79, 131)$,
$(0, 46, 67, 133)$,
$(0, 12, 51, 70)$,
$(0, 91, 103, 162)$,\adfsplit
$(0, 13, 84, 152)$,
$(0, 85, 88, 118)$

\adfLgap \noindent by the mapping:
$x \mapsto x + 3 j \adfmod{126}$ for $x < 126$,
$x \mapsto (x +  j \adfmod{42}) + 126$ for $126 \le x < 168$,
$x \mapsto x$ for $x \ge 168$,
$0 \le j < 42$.
\ADFvfyParStart{(170, ((50, 42, ((126, 3), (42, 1), (2, 2)))), ((14, 9), (44, 1)))} 

\adfDgap
\noindent{\boldmath $ 14^{9} 47^{1} $}~
With the point set $Z_{173}$ partitioned into
 residue classes modulo $9$ for $\{0, 1, \dots, 125\}$, and
 $\{126, 127, \dots, 172\}$,
 the design is generated from

\adfLgap 
$(171, 9, 102, 106)$,
$(171, 122, 67, 83)$,
$(172, 107, 73, 12)$,
$(172, 50, 81, 112)$,\adfsplit
$(126, 31, 116, 15)$,
$(126, 54, 39, 77)$,
$(126, 30, 4, 29)$,
$(126, 118, 115, 20)$,\adfsplit
$(126, 91, 70, 114)$,
$(126, 35, 27, 14)$,
$(127, 1, 79, 27)$,
$(127, 28, 87, 77)$,\adfsplit
$(127, 86, 110, 57)$,
$(127, 30, 89, 16)$,
$(127, 78, 22, 54)$,
$(127, 47, 67, 62)$,\adfsplit
$(128, 107, 87, 72)$,
$(128, 3, 94, 59)$,
$(128, 104, 47, 85)$,
$(128, 52, 118, 12)$,\adfsplit
$(128, 9, 26, 37)$,
$(128, 78, 20, 43)$,
$(129, 24, 26, 63)$,
$(129, 15, 41, 11)$,\adfsplit
$(129, 64, 68, 70)$,
$(129, 121, 48, 55)$,
$(129, 18, 40, 107)$,
$(129, 110, 3, 97)$,\adfsplit
$(130, 94, 74, 48)$,
$(130, 90, 82, 26)$,
$(130, 23, 83, 31)$,
$(130, 34, 37, 105)$,\adfsplit
$(130, 63, 114, 57)$,
$(130, 53, 25, 86)$,
$(131, 36, 53, 2)$,
$(131, 81, 33, 103)$,\adfsplit
$(131, 25, 120, 68)$,
$(131, 44, 4, 47)$,
$(131, 16, 46, 23)$,
$(131, 111, 24, 91)$,\adfsplit
$(132, 20, 67, 105)$,
$(132, 22, 5, 84)$,
$(132, 108, 111, 53)$,
$(132, 10, 65, 25)$,\adfsplit
$(132, 109, 60, 117)$,
$(132, 80, 32, 124)$,
$(133, 32, 7, 120)$,
$(133, 24, 44, 82)$,\adfsplit
$(133, 95, 18, 16)$,
$(133, 94, 45, 33)$,
$(133, 11, 17, 57)$,
$(133, 19, 31, 2)$,\adfsplit
$(134, 102, 74, 113)$,
$(134, 122, 62, 121)$,
$(134, 58, 29, 109)$,
$(134, 17, 99, 33)$,\adfsplit
$(134, 124, 108, 3)$,
$(134, 60, 97, 10)$,
$(135, 97, 1, 72)$,
$(135, 110, 17, 111)$,\adfsplit
$(135, 51, 4, 27)$,
$(135, 49, 14, 96)$,
$(135, 70, 82, 83)$,
$(135, 8, 30, 23)$,\adfsplit
$(136, 48, 43, 0)$,
$(136, 1, 123, 74)$,
$(136, 75, 32, 78)$,
$(136, 83, 95, 13)$,\adfsplit
$(136, 94, 62, 10)$,
$(3, 28, 65, 166)$,
$(0, 105, 110, 116)$,
$(0, 6, 14, 56)$,\adfsplit
$(0, 10, 29, 32)$,
$(0, 28, 86, 137)$,
$(0, 1, 51, 88)$,
$(0, 12, 42, 138)$,\adfsplit
$(0, 19, 34, 139)$,
$(1, 22, 63, 70)$,
$(0, 53, 66, 154)$,
$(0, 52, 85, 140)$,\adfsplit
$(0, 115, 122, 152)$,
$(0, 97, 101, 155)$,
$(0, 5, 107, 167)$,
$(0, 65, 93, 170)$,\adfsplit
$(0, 41, 83, 109)$,
$(1, 77, 125, 153)$,
$(1, 7, 87, 167)$,
$(2, 107, 112, 137)$,\adfsplit
$(3, 4, 87, 152)$,
$(1, 3, 85, 95)$,
$(1, 23, 34, 170)$,
$(1, 50, 58, 169)$,\adfsplit
$(2, 57, 59, 138)$,
$(1, 25, 117, 138)$,
$(1, 15, 65, 154)$,
$(1, 38, 52, 113)$,\adfsplit
$(3, 17, 112, 139)$,
$(3, 22, 33, 107)$,
$(2, 28, 52, 153)$,
$(2, 9, 82, 168)$,\adfsplit
$(2, 105, 118, 140)$,
$(2, 14, 81, 155)$,
$(2, 32, 93, 154)$

\adfLgap \noindent by the mapping:
$x \mapsto x + 6 j \adfmod{126}$ for $x < 126$,
$x \mapsto (x - 126 + 15 j \adfmod{45}) + 126$ for $126 \le x < 171$,
$x \mapsto x$ for $x \ge 171$,
$0 \le j < 21$.
\ADFvfyParStart{(173, ((103, 21, ((126, 6), (45, 15), (2, 2)))), ((14, 9), (47, 1)))} 

\adfDgap
\noindent{\boldmath $ 14^{9} 50^{1} $}~
With the point set $Z_{176}$ partitioned into
 residue classes modulo $9$ for $\{0, 1, \dots, 125\}$, and
 $\{126, 127, \dots, 175\}$,
 the design is generated from

\adfLgap 
$(168, 117, 107, 55)$,
$(169, 91, 90, 110)$,
$(170, 43, 54, 104)$,
$(171, 44, 64, 72)$,\adfsplit
$(172, 37, 104, 24)$,
$(173, 106, 75, 110)$,
$(174, 97, 111, 53)$,
$(175, 66, 118, 113)$,\adfsplit
$(126, 97, 12, 107)$,
$(126, 3, 47, 8)$,
$(126, 20, 67, 89)$,
$(126, 36, 84, 95)$,\adfsplit
$(126, 119, 23, 46)$,
$(126, 91, 34, 120)$,
$(126, 117, 123, 52)$,
$(126, 14, 43, 39)$,\adfsplit
$(126, 96, 37, 2)$,
$(126, 15, 108, 53)$,
$(126, 44, 61, 92)$,
$(126, 109, 111, 68)$,\adfsplit
$(126, 83, 7, 99)$,
$(126, 19, 115, 13)$,
$(126, 98, 4, 102)$,
$(126, 54, 17, 51)$,\adfsplit
$(126, 33, 110, 40)$,
$(0, 2, 8, 145)$,
$(0, 10, 26, 166)$,
$(0, 12, 41, 148)$,\adfsplit
$(0, 14, 17, 153)$,
$(0, 16, 23, 56)$,
$(0, 19, 65, 107)$,
$(0, 46, 74, 138)$,\adfsplit
$(0, 53, 79, 142)$,
$(1, 35, 56, 126)$,
$(1, 59, 119, 137)$,
$(0, 104, 119, 162)$,\adfsplit
$(0, 100, 125, 141)$,
$(1, 14, 38, 159)$,
$(1, 2, 116, 127)$,
$(0, 62, 113, 154)$,\adfsplit
$(1, 49, 113, 147)$,
$(0, 43, 86, 158)$,
$(1, 50, 52, 160)$,
$(0, 24, 84, 147)$,\adfsplit
$(0, 25, 75, 143)$,
$(0, 91, 103, 137)$,
$(0, 30, 87, 159)$,
$(0, 37, 58, 152)$,\adfsplit
$(0, 15, 88, 121)$,
$(0, 94, 109, 165)$,
$(1, 4, 43, 142)$,
$(0, 22, 82, 146)$,\adfsplit
$(0, 21, 70, 133)$

\adfLgap \noindent by the mapping:
$x \mapsto x + 3 j \adfmod{126}$ for $x < 126$,
$x \mapsto (x +  j \adfmod{42}) + 126$ for $126 \le x < 168$,
$x \mapsto x$ for $x \ge 168$,
$0 \le j < 42$.
\ADFvfyParStart{(176, ((53, 42, ((126, 3), (42, 1), (8, 8)))), ((14, 9), (50, 1)))} 

\adfDgap
\noindent{\boldmath $ 14^{9} 53^{1} $}~
With the point set $Z_{179}$ partitioned into
 residue classes modulo $9$ for $\{0, 1, \dots, 125\}$, and
 $\{126, 127, \dots, 178\}$,
 the design is generated from

\adfLgap 
$(177, 77, 94, 92)$,
$(177, 120, 117, 115)$,
$(178, 88, 65, 0)$,
$(178, 20, 3, 31)$,\adfsplit
$(126, 17, 6, 36)$,
$(126, 55, 41, 88)$,
$(126, 118, 3, 22)$,
$(126, 87, 74, 48)$,\adfsplit
$(126, 13, 86, 44)$,
$(126, 47, 79, 9)$,
$(127, 62, 86, 111)$,
$(127, 5, 38, 97)$,\adfsplit
$(127, 114, 47, 99)$,
$(127, 19, 72, 17)$,
$(127, 33, 12, 58)$,
$(127, 118, 67, 34)$,\adfsplit
$(128, 41, 26, 46)$,
$(128, 122, 56, 61)$,
$(128, 124, 22, 30)$,
$(128, 37, 108, 60)$,\adfsplit
$(128, 111, 87, 117)$,
$(128, 85, 101, 89)$,
$(129, 114, 8, 104)$,
$(129, 48, 92, 124)$,\adfsplit
$(129, 28, 101, 81)$,
$(129, 3, 67, 19)$,
$(129, 113, 25, 51)$,
$(129, 18, 58, 17)$,\adfsplit
$(130, 36, 89, 7)$,
$(130, 8, 102, 14)$,
$(130, 109, 110, 103)$,
$(130, 101, 46, 114)$,\adfsplit
$(130, 39, 27, 22)$,
$(130, 69, 113, 70)$,
$(131, 4, 82, 104)$,
$(131, 51, 102, 91)$,\adfsplit
$(131, 74, 124, 108)$,
$(131, 114, 3, 13)$,
$(131, 25, 80, 45)$,
$(131, 77, 47, 53)$,\adfsplit
$(132, 7, 107, 57)$,
$(132, 59, 36, 49)$,
$(132, 34, 94, 55)$,
$(132, 101, 123, 50)$,\adfsplit
$(132, 28, 120, 80)$,
$(132, 110, 63, 6)$,
$(133, 26, 45, 112)$,
$(133, 73, 95, 87)$,\adfsplit
$(133, 79, 64, 12)$,
$(133, 38, 67, 50)$,
$(133, 18, 107, 16)$,
$(133, 65, 93, 114)$,\adfsplit
$(134, 118, 62, 43)$,
$(134, 122, 60, 111)$,
$(134, 109, 48, 67)$,
$(134, 27, 92, 23)$,\adfsplit
$(134, 112, 124, 35)$,
$(134, 11, 54, 123)$,
$(135, 102, 60, 44)$,
$(135, 35, 38, 93)$,\adfsplit
$(135, 18, 67, 59)$,
$(135, 118, 105, 104)$,
$(135, 7, 47, 81)$,
$(135, 73, 76, 70)$,\adfsplit
$(136, 86, 6, 76)$,
$(136, 45, 71, 91)$,
$(136, 77, 93, 52)$,
$(136, 11, 51, 85)$,\adfsplit
$(136, 43, 36, 92)$,
$(136, 48, 28, 62)$,
$(137, 7, 53, 48)$,
$(137, 19, 63, 95)$,\adfsplit
$(137, 26, 65, 94)$,
$(137, 67, 124, 75)$,
$(137, 20, 100, 123)$,
$(0, 8, 12, 171)$,\adfsplit
$(1, 14, 103, 138)$,
$(0, 1, 31, 139)$,
$(0, 79, 91, 140)$,
$(1, 29, 71, 141)$,\adfsplit
$(1, 22, 26, 69)$,
$(0, 43, 109, 142)$,
$(1, 59, 105, 157)$,
$(1, 80, 88, 158)$,\adfsplit
$(1, 16, 104, 172)$,
$(1, 65, 86, 173)$,
$(1, 39, 44, 159)$,
$(1, 33, 70, 175)$,\adfsplit
$(0, 2, 33, 37)$,
$(0, 24, 74, 158)$,
$(2, 89, 99, 175)$,
$(0, 93, 100, 175)$,\adfsplit
$(2, 23, 80, 173)$,
$(2, 9, 87, 174)$,
$(2, 39, 100, 172)$,
$(2, 35, 46, 157)$,\adfsplit
$(2, 69, 112, 176)$,
$(3, 5, 100, 176)$,
$(3, 58, 59, 139)$,
$(0, 3, 87, 138)$,\adfsplit
$(3, 34, 69, 156)$,
$(0, 112, 119, 172)$,
$(0, 35, 101, 155)$,
$(0, 10, 29, 159)$,\adfsplit
$(0, 17, 82, 95)$,
$(2, 5, 64, 140)$,
$(0, 47, 98, 176)$,
$(0, 4, 66, 157)$,\adfsplit
$(0, 6, 28, 173)$

\adfLgap \noindent by the mapping:
$x \mapsto x + 6 j \adfmod{126}$ for $x < 126$,
$x \mapsto (x - 126 + 17 j \adfmod{51}) + 126$ for $126 \le x < 177$,
$x \mapsto x$ for $x \ge 177$,
$0 \le j < 21$.
\ADFvfyParStart{(179, ((109, 21, ((126, 6), (51, 17), (2, 2)))), ((14, 9), (53, 1)))} 

\section{4-GDDs for the proof of Lemma \ref{lem:4-GDD 20^u m^1}}
\label{app:4-GDD 20^u m^1}
\adfhide{
$ 20^9 11^1 $,
$ 20^9 17^1 $ and
$ 20^9 23^1 $.
}

\adfDgap
\noindent{\boldmath $ 20^{9} 11^{1} $}~
With the point set $Z_{191}$ partitioned into
 residue classes modulo $9$ for $\{0, 1, \dots, 179\}$, and
 $\{180, 181, \dots, 190\}$,
 the design is generated from

\adfLgap 
$(189, 37, 114, 99)$,
$(189, 82, 161, 50)$,
$(190, 118, 15, 0)$,
$(190, 5, 8, 103)$,\adfsplit
$(180, 94, 88, 109)$,
$(180, 71, 139, 158)$,
$(180, 136, 159, 147)$,
$(180, 5, 137, 62)$,\adfsplit
$(180, 18, 151, 30)$,
$(180, 164, 150, 153)$,
$(181, 54, 58, 74)$,
$(181, 104, 67, 69)$,\adfsplit
$(181, 78, 95, 134)$,
$(181, 46, 27, 101)$,
$(181, 147, 17, 16)$,
$(181, 91, 84, 79)$,\adfsplit
$(182, 32, 62, 82)$,
$(182, 55, 160, 15)$,
$(182, 117, 137, 121)$,
$(182, 110, 72, 59)$,\adfsplit
$(182, 39, 143, 150)$,
$(182, 115, 84, 76)$,
$(75, 103, 132, 62)$,
$(78, 111, 144, 125)$,\adfsplit
$(133, 64, 84, 130)$,
$(130, 123, 154, 32)$,
$(38, 142, 21, 126)$,
$(103, 144, 26, 57)$,\adfsplit
$(173, 175, 19, 140)$,
$(38, 77, 99, 170)$,
$(114, 112, 98, 64)$,
$(103, 154, 84, 20)$,\adfsplit
$(81, 172, 76, 105)$,
$(39, 87, 125, 145)$,
$(177, 122, 58, 125)$,
$(4, 154, 132, 101)$,\adfsplit
$(127, 80, 14, 57)$,
$(22, 176, 144, 21)$,
$(49, 27, 120, 128)$,
$(179, 149, 61, 177)$,\adfsplit
$(128, 135, 143, 121)$,
$(93, 177, 130, 89)$,
$(123, 152, 84, 126)$,
$(60, 142, 14, 4)$,\adfsplit
$(9, 48, 170, 50)$,
$(106, 94, 167, 132)$,
$(131, 173, 42, 178)$,
$(39, 146, 133, 163)$,\adfsplit
$(120, 4, 109, 36)$,
$(4, 150, 128, 75)$,
$(0, 1, 50, 74)$,
$(0, 13, 128, 140)$,\adfsplit
$(0, 10, 104, 146)$,
$(0, 11, 44, 129)$,
$(0, 6, 106, 176)$,
$(1, 4, 32, 110)$,\adfsplit
$(0, 21, 80, 86)$,
$(1, 9, 92, 94)$,
$(1, 7, 62, 158)$,
$(1, 2, 51, 139)$,\adfsplit
$(0, 23, 83, 98)$,
$(0, 29, 40, 152)$,
$(1, 33, 58, 176)$,
$(2, 23, 94, 177)$,\adfsplit
$(1, 26, 34, 77)$,
$(2, 105, 135, 178)$,
$(2, 3, 69, 82)$,
$(1, 59, 148, 170)$,\adfsplit
$(2, 40, 142, 161)$,
$(3, 130, 143, 167)$,
$(1, 21, 27, 137)$,
$(1, 47, 129, 171)$,\adfsplit
$(0, 51, 67, 111)$,
$(3, 58, 83, 118)$,
$(1, 40, 71, 81)$,
$(1, 22, 39, 115)$,\adfsplit
$(1, 53, 61, 130)$,
$(0, 24, 112, 119)$,
$(1, 16, 65, 131)$,
$(3, 35, 47, 81)$,\adfsplit
$(1, 35, 49, 147)$,
$(0, 91, 101, 159)$,
$(0, 5, 28, 132)$,
$(0, 43, 137, 166)$,\adfsplit
$(0, 77, 94, 179)$,
$(0, 61, 78, 157)$,
$(0, 107, 113, 145)$,
$(0, 25, 65, 120)$,\adfsplit
$(0, 55, 59, 155)$,
$(0, 30, 71, 127)$,
$(0, 37, 115, 143)$

\adfLgap \noindent by the mapping:
$x \mapsto x + 6 j \adfmod{180}$ for $x < 180$,
$x \mapsto (x + 3 j \adfmod{9}) + 180$ for $180 \le x < 189$,
$x \mapsto x$ for $x \ge 189$,
$0 \le j < 30$.
\ADFvfyParStart{(191, ((91, 30, ((180, 6), (9, 3), (2, 2)))), ((20, 9), (11, 1)))} 

\adfDgap
\noindent{\boldmath $ 20^{9} 17^{1} $}~
With the point set $Z_{197}$ partitioned into
 residue classes modulo $9$ for $\{0, 1, \dots, 179\}$, and
 $\{180, 181, \dots, 196\}$,
 the design is generated from

\adfLgap 
$(195, 177, 149, 157)$,
$(195, 108, 62, 112)$,
$(196, 108, 160, 32)$,
$(196, 3, 139, 119)$,\adfsplit
$(180, 147, 140, 90)$,
$(180, 40, 9, 23)$,
$(180, 87, 172, 164)$,
$(180, 11, 150, 134)$,\adfsplit
$(180, 31, 133, 17)$,
$(180, 178, 156, 145)$,
$(181, 162, 87, 95)$,
$(181, 104, 42, 109)$,\adfsplit
$(181, 112, 129, 133)$,
$(181, 16, 48, 62)$,
$(181, 83, 27, 161)$,
$(181, 110, 82, 13)$,\adfsplit
$(182, 175, 37, 61)$,
$(182, 44, 173, 76)$,
$(182, 32, 159, 99)$,
$(182, 172, 128, 17)$,\adfsplit
$(182, 149, 78, 30)$,
$(182, 57, 124, 90)$,
$(183, 168, 81, 102)$,
$(183, 31, 41, 91)$,\adfsplit
$(183, 108, 46, 79)$,
$(183, 83, 76, 33)$,
$(183, 178, 74, 8)$,
$(183, 165, 14, 107)$,\adfsplit
$(184, 123, 92, 148)$,
$(184, 21, 32, 28)$,
$(184, 127, 153, 5)$,
$(184, 78, 169, 134)$,\adfsplit
$(184, 71, 137, 67)$,
$(184, 72, 12, 52)$,
$(7, 60, 53, 56)$,
$(75, 116, 167, 114)$,\adfsplit
$(163, 124, 162, 60)$,
$(135, 52, 10, 92)$,
$(173, 25, 77, 3)$,
$(145, 160, 58, 59)$,\adfsplit
$(144, 104, 19, 25)$,
$(10, 26, 108, 25)$,
$(115, 21, 128, 149)$,
$(74, 145, 148, 122)$,\adfsplit
$(125, 148, 52, 36)$,
$(49, 110, 36, 51)$,
$(136, 6, 50, 31)$,
$(165, 67, 2, 5)$,\adfsplit
$(105, 63, 49, 86)$,
$(36, 147, 2, 35)$,
$(37, 78, 26, 75)$,
$(79, 84, 170, 0)$,\adfsplit
$(47, 80, 41, 144)$,
$(97, 41, 78, 81)$,
$(61, 148, 2, 167)$,
$(0, 7, 32, 38)$,\adfsplit
$(0, 5, 80, 110)$,
$(0, 20, 33, 158)$,
$(1, 9, 74, 152)$,
$(1, 40, 128, 140)$,\adfsplit
$(0, 8, 28, 85)$,
$(0, 68, 92, 152)$,
$(0, 10, 115, 122)$,
$(0, 21, 26, 73)$,\adfsplit
$(1, 13, 68, 93)$,
$(1, 22, 63, 164)$,
$(1, 23, 44, 58)$,
$(2, 4, 63, 125)$,\adfsplit
$(2, 17, 64, 93)$,
$(2, 3, 99, 112)$,
$(2, 41, 89, 118)$,
$(2, 87, 111, 160)$,\adfsplit
$(2, 113, 124, 135)$,
$(0, 31, 107, 131)$,
$(1, 49, 117, 178)$,
$(1, 29, 31, 171)$,\adfsplit
$(1, 75, 105, 153)$,
$(0, 49, 100, 137)$,
$(1, 17, 51, 85)$,
$(0, 6, 17, 43)$,\adfsplit
$(0, 39, 65, 109)$,
$(1, 33, 112, 143)$,
$(0, 12, 59, 157)$,
$(0, 75, 87, 133)$,\adfsplit
$(0, 23, 42, 163)$,
$(0, 24, 70, 125)$,
$(0, 58, 172, 178)$,
$(0, 30, 124, 136)$,\adfsplit
$(4, 17, 47, 136)$,
$(3, 83, 148, 178)$,
$(0, 51, 149, 154)$,
$(3, 4, 69, 160)$,\adfsplit
$(0, 69, 88, 155)$,
$(3, 40, 167, 179)$,
$(0, 35, 95, 166)$,
$(0, 29, 165, 167)$,\adfsplit
$(0, 76, 123, 129)$

\adfLgap \noindent by the mapping:
$x \mapsto x + 6 j \adfmod{180}$ for $x < 180$,
$x \mapsto (x + 5 j \adfmod{15}) + 180$ for $180 \le x < 195$,
$x \mapsto x$ for $x \ge 195$,
$0 \le j < 30$.
\ADFvfyParStart{(197, ((97, 30, ((180, 6), (15, 5), (2, 2)))), ((20, 9), (17, 1)))} 

\adfDgap
\noindent{\boldmath $ 20^{9} 23^{1} $}~
With the point set $Z_{203}$ partitioned into
 residue classes modulo $9$ for $\{0, 1, \dots, 179\}$, and
 $\{180, 181, \dots, 202\}$,
 the design is generated from

\adfLgap 
$(201, 108, 35, 14)$,
$(201, 82, 165, 97)$,
$(202, 169, 170, 48)$,
$(202, 77, 21, 172)$,\adfsplit
$(180, 100, 49, 24)$,
$(180, 75, 53, 140)$,
$(180, 106, 44, 37)$,
$(180, 54, 123, 56)$,\adfsplit
$(180, 61, 63, 155)$,
$(180, 59, 174, 130)$,
$(181, 61, 134, 14)$,
$(181, 3, 90, 51)$,\adfsplit
$(181, 78, 102, 28)$,
$(181, 81, 160, 73)$,
$(181, 121, 38, 71)$,
$(181, 112, 77, 65)$,\adfsplit
$(182, 162, 142, 113)$,
$(182, 39, 65, 109)$,
$(182, 42, 98, 74)$,
$(182, 45, 130, 133)$,\adfsplit
$(182, 10, 31, 17)$,
$(182, 84, 32, 87)$,
$(183, 148, 100, 164)$,
$(183, 13, 59, 124)$,\adfsplit
$(183, 78, 43, 68)$,
$(183, 162, 53, 120)$,
$(183, 159, 37, 75)$,
$(183, 26, 99, 29)$,\adfsplit
$(184, 162, 77, 174)$,
$(184, 161, 172, 148)$,
$(184, 74, 103, 99)$,
$(184, 52, 165, 119)$,\adfsplit
$(184, 122, 114, 80)$,
$(184, 61, 109, 123)$,
$(157, 17, 115, 11)$,
$(96, 151, 112, 118)$,\adfsplit
$(105, 125, 14, 124)$,
$(17, 140, 1, 128)$,
$(151, 116, 28, 75)$,
$(104, 6, 130, 91)$,\adfsplit
$(98, 15, 72, 199)$,
$(54, 134, 3, 105)$,
$(57, 64, 106, 74)$,
$(122, 64, 54, 94)$,\adfsplit
$(150, 139, 64, 44)$,
$(15, 83, 58, 179)$,
$(27, 29, 193, 14)$,
$(4, 200, 91, 75)$,\adfsplit
$(93, 62, 56, 73)$,
$(45, 125, 3, 173)$,
$(142, 165, 173, 141)$,
$(45, 83, 107, 55)$,\adfsplit
$(163, 67, 117, 132)$,
$(45, 109, 15, 50)$,
$(20, 105, 16, 138)$,
$(19, 51, 12, 62)$,\adfsplit
$(98, 141, 61, 192)$,
$(84, 118, 4, 159)$,
$(166, 164, 86, 125)$,
$(116, 121, 12, 7)$,\adfsplit
$(0, 14, 21, 33)$,
$(0, 37, 111, 177)$,
$(1, 4, 147, 153)$,
$(1, 23, 99, 159)$,\adfsplit
$(0, 47, 52, 87)$,
$(0, 23, 46, 105)$,
$(0, 29, 57, 118)$,
$(1, 27, 50, 80)$,\adfsplit
$(0, 15, 70, 92)$,
$(1, 57, 107, 158)$,
$(1, 45, 59, 92)$,
$(0, 148, 159, 185)$,\adfsplit
$(1, 68, 124, 129)$,
$(2, 46, 107, 111)$,
$(3, 16, 89, 119)$,
$(3, 106, 149, 166)$,\adfsplit
$(2, 16, 123, 186)$,
$(2, 3, 77, 130)$,
$(0, 38, 78, 142)$,
$(1, 13, 32, 116)$,\adfsplit
$(0, 82, 116, 166)$,
$(2, 10, 50, 88)$,
$(2, 100, 179, 199)$,
$(0, 4, 151, 172)$,\adfsplit
$(0, 28, 79, 134)$,
$(0, 17, 178, 192)$,
$(0, 110, 161, 176)$,
$(0, 44, 112, 167)$,\adfsplit
$(0, 59, 91, 152)$,
$(0, 20, 154, 193)$,
$(0, 6, 143, 164)$,
$(1, 16, 65, 121)$,\adfsplit
$(0, 5, 88, 125)$,
$(1, 11, 79, 199)$,
$(0, 19, 60, 179)$,
$(0, 11, 89, 114)$,\adfsplit
$(0, 35, 73, 173)$,
$(1, 5, 71, 193)$,
$(0, 53, 61, 84)$,
$(0, 1, 132, 186)$,\adfsplit
$(0, 41, 43, 67)$,
$(0, 13, 101, 163)$,
$(0, 30, 127, 133)$

\adfLgap \noindent by the mapping:
$x \mapsto x + 6 j \adfmod{180}$ for $x < 180$,
$x \mapsto (x - 180 + 7 j \adfmod{21}) + 180$ for $180 \le x < 201$,
$x \mapsto x$ for $x \ge 201$,
$0 \le j < 30$.
\ADFvfyParStart{(203, ((103, 30, ((180, 6), (21, 7), (2, 2)))), ((20, 9), (23, 1)))} 

\section{4-GDDs for the proof of Lemma \ref{lem:4-GDD 22^u m^1}}
\label{app:4-GDD 22^u m^1}
\adfhide{
$ 22^6 25^1 $,
$ 22^6 28^1 $,
$ 22^6 31^1 $,
$ 22^6 34^1 $,
$ 22^6 37^1 $,
$ 22^6 40^1 $,
$ 22^6 43^1 $,
$ 22^6 46^1 $,
$ 22^6 49^1 $,
$ 22^6 52^1 $,
$ 22^9 7^1 $,
$ 22^9 10^1 $,
$ 22^9 13^1 $,
$ 22^9 16^1 $,
$ 22^9 19^1 $,
$ 22^{12} 7^1 $,
$ 22^{12} 10^1 $,
$ 22^{12} 13^1 $,
$ 22^{12} 16^1 $,
$ 22^{12} 19^1 $,
$ 22^{15} 4^1 $,
$ 22^{15} 10^1 $,
$ 22^{15} 13^1 $,
$ 22^{15} 16^1 $,
$ 22^{15} 19^1 $,
$ 22^{18} 4^1 $,
$ 22^{18} 7^1 $,
$ 22^{18} 10^1 $,
$ 22^{18} 13^1 $,
$ 22^{18} 16^1 $ and
$ 22^{18} 19^1 $.
}

\adfDgap
\noindent{\boldmath $ 22^{6} 25^{1} $}~
With the point set $Z_{157}$ partitioned into
 residue classes modulo $6$ for $\{0, 1, \dots, 131\}$, and
 $\{132, 133, \dots, 156\}$,
 the design is generated from

\adfLgap 
$(132, 106, 48, 95)$,
$(132, 39, 128, 19)$,
$(133, 115, 90, 94)$,
$(133, 62, 63, 17)$,\adfsplit
$(134, 111, 52, 125)$,
$(134, 32, 18, 37)$,
$(135, 122, 29, 42)$,
$(135, 85, 21, 76)$,\adfsplit
$(136, 0, 101, 64)$,
$(136, 93, 8, 37)$,
$(137, 76, 55, 89)$,
$(137, 120, 117, 20)$,\adfsplit
$(138, 109, 60, 110)$,
$(138, 58, 9, 125)$,
$(0, 2, 55, 112)$,
$(0, 3, 7, 40)$,\adfsplit
$(0, 15, 107, 116)$,
$(0, 11, 26, 105)$,
$(0, 8, 41, 46)$,
$(0, 17, 34, 104)$,\adfsplit
$(0, 51, 61, 113)$,
$(0, 27, 35, 109)$,
$(0, 56, 125, 147)$,
$(0, 10, 81, 103)$,\adfsplit
$(0, 57, 89, 155)$,
$(0, 63, 65, 91)$,
$(156, 0, 44, 88)$,
$(156, 1, 45, 89)$

\adfLgap \noindent by the mapping:
$x \mapsto x + 2 j \adfmod{132}$ for $x < 132$,
$x \mapsto (x - 132 + 8 j \adfmod{24}) + 132$ for $132 \le x < 156$,
$156 \mapsto 156$,
$0 \le j < 66$
 for the first 26 blocks,
$0 \le j < 22$
 for the last two blocks.
\ADFvfyParStart{(157, ((26, 66, ((132, 2), (24, 8), (1, 1))), (2, 22, ((132, 2), (24, 8), (1, 1)))), ((22, 6), (25, 1)))} 

\adfDgap
\noindent{\boldmath $ 22^{6} 28^{1} $}~
With the point set $Z_{160}$ partitioned into
 residue classes modulo $6$ for $\{0, 1, \dots, 131\}$, and
 $\{132, 133, \dots, 159\}$,
 the design is generated from

\adfLgap 
$(132, 119, 86, 78)$,
$(132, 15, 10, 61)$,
$(133, 12, 80, 55)$,
$(133, 107, 111, 130)$,\adfsplit
$(134, 14, 34, 61)$,
$(134, 120, 129, 119)$,
$(135, 88, 60, 73)$,
$(135, 107, 57, 8)$,\adfsplit
$(136, 79, 122, 9)$,
$(136, 120, 5, 10)$,
$(137, 129, 115, 34)$,
$(137, 47, 110, 108)$,\adfsplit
$(138, 125, 2, 60)$,
$(138, 22, 109, 33)$,
$(0, 1, 14, 40)$,
$(0, 3, 67, 83)$,\adfsplit
$(0, 21, 50, 79)$,
$(1, 3, 29, 139)$,
$(0, 10, 56, 139)$,
$(0, 23, 31, 34)$,\adfsplit
$(0, 16, 61, 101)$,
$(0, 15, 37, 94)$,
$(0, 35, 55, 62)$,
$(0, 63, 97, 140)$,\adfsplit
$(0, 4, 77, 115)$,
$(0, 7, 39, 80)$,
$(0, 25, 100, 149)$,
$(159, 0, 44, 88)$,\adfsplit
$(159, 1, 45, 89)$

\adfLgap \noindent by the mapping:
$x \mapsto x + 2 j \adfmod{132}$ for $x < 132$,
$x \mapsto (x - 132 + 9 j \adfmod{27}) + 132$ for $132 \le x < 159$,
$159 \mapsto 159$,
$0 \le j < 66$
 for the first 27 blocks,
$0 \le j < 22$
 for the last two blocks.
\ADFvfyParStart{(160, ((27, 66, ((132, 2), (27, 9), (1, 1))), (2, 22, ((132, 2), (27, 9), (1, 1)))), ((22, 6), (28, 1)))} 

\adfDgap
\noindent{\boldmath $ 22^{6} 31^{1} $}~
With the point set $Z_{163}$ partitioned into
 residue classes modulo $6$ for $\{0, 1, \dots, 131\}$, and
 $\{132, 133, \dots, 162\}$,
 the design is generated from

\adfLgap 
$(132, 6, 10, 15)$,
$(132, 74, 7, 107)$,
$(133, 8, 77, 88)$,
$(133, 66, 115, 33)$,\adfsplit
$(134, 103, 65, 63)$,
$(134, 114, 10, 68)$,
$(135, 16, 35, 92)$,
$(135, 54, 111, 1)$,\adfsplit
$(136, 104, 101, 105)$,
$(136, 36, 61, 70)$,
$(137, 127, 120, 65)$,
$(137, 4, 128, 93)$,\adfsplit
$(138, 21, 101, 46)$,
$(138, 44, 127, 60)$,
$(139, 45, 118, 29)$,
$(0, 2, 13, 139)$,\adfsplit
$(0, 3, 26, 118)$,
$(0, 15, 32, 119)$,
$(0, 23, 31, 38)$,
$(0, 20, 47, 73)$,\adfsplit
$(0, 21, 35, 62)$,
$(0, 41, 51, 127)$,
$(0, 10, 95, 140)$,
$(0, 37, 101, 160)$,\adfsplit
$(0, 22, 93, 141)$,
$(0, 39, 68, 113)$,
$(0, 50, 111, 131)$,
$(0, 29, 63, 161)$,\adfsplit
$(162, 0, 44, 88)$,
$(162, 1, 45, 89)$

\adfLgap \noindent by the mapping:
$x \mapsto x + 2 j \adfmod{132}$ for $x < 132$,
$x \mapsto (x - 132 + 10 j \adfmod{30}) + 132$ for $132 \le x < 162$,
$162 \mapsto 162$,
$0 \le j < 66$
 for the first 28 blocks,
$0 \le j < 22$
 for the last two blocks.
\ADFvfyParStart{(163, ((28, 66, ((132, 2), (30, 10), (1, 1))), (2, 22, ((132, 2), (30, 10), (1, 1)))), ((22, 6), (31, 1)))} 

\adfDgap
\noindent{\boldmath $ 22^{6} 34^{1} $}~
With the point set $Z_{166}$ partitioned into
 residue classes modulo $6$ for $\{0, 1, \dots, 131\}$, and
 $\{132, 133, \dots, 165\}$,
 the design is generated from

\adfLgap 
$(132, 35, 63, 109)$,
$(132, 72, 127, 130)$,
$(132, 38, 67, 81)$,
$(132, 28, 45, 47)$,\adfsplit
$(132, 98, 30, 33)$,
$(132, 50, 18, 39)$,
$(132, 90, 73, 95)$,
$(132, 55, 125, 40)$,\adfsplit
$(132, 62, 46, 69)$,
$(132, 22, 8, 57)$,
$(132, 36, 9, 89)$,
$(132, 17, 52, 60)$,\adfsplit
$(132, 65, 103, 78)$,
$(132, 48, 86, 85)$,
$(132, 120, 70, 49)$,
$(132, 68, 34, 27)$,\adfsplit
$(0, 11, 92, 132)$,
$(1, 9, 41, 159)$,
$(0, 75, 95, 144)$,
$(0, 9, 28, 157)$,\adfsplit
$(0, 2, 41, 137)$,
$(0, 31, 70, 160)$,
$(0, 33, 56, 101)$,
$(0, 10, 69, 73)$,\adfsplit
$(0, 4, 26, 83)$,
$(0, 13, 47, 80)$,
$(0, 46, 117, 127)$,
$(0, 20, 107, 123)$,\adfsplit
$(0, 1, 27, 77)$,
$(165, 0, 44, 88)$,
$(165, 1, 45, 89)$

\adfLgap \noindent by the mapping:
$x \mapsto x + 2 j \adfmod{132}$ for $x < 132$,
$x \mapsto (x +  j \adfmod{33}) + 132$ for $132 \le x < 165$,
$165 \mapsto 165$,
$0 \le j < 66$
 for the first 29 blocks,
$0 \le j < 22$
 for the last two blocks.
\ADFvfyParStart{(166, ((29, 66, ((132, 2), (33, 1), (1, 1))), (2, 22, ((132, 2), (33, 1), (1, 1)))), ((22, 6), (34, 1)))} 

\adfDgap
\noindent{\boldmath $ 22^{6} 37^{1} $}~
With the point set $Z_{169}$ partitioned into
 residue classes modulo $6$ for $\{0, 1, \dots, 131\}$, and
 $\{132, 133, \dots, 168\}$,
 the design is generated from

\adfLgap 
$(165, 9, 23, 40)$,
$(132, 30, 10, 113)$,
$(132, 32, 58, 114)$,
$(132, 89, 67, 66)$,\adfsplit
$(132, 39, 41, 50)$,
$(132, 127, 68, 11)$,
$(132, 73, 27, 120)$,
$(0, 3, 8, 41)$,\adfsplit
$(0, 4, 32, 67)$,
$(0, 10, 81, 156)$,
$(0, 7, 34, 87)$,
$(0, 15, 107, 160)$,\adfsplit
$(0, 13, 77, 147)$,
$(0, 21, 95, 149)$,
$(0, 19, 89, 162)$,
$(168, 0, 44, 88)$

\adfLgap \noindent by the mapping:
$x \mapsto x +  j \adfmod{132}$ for $x < 132$,
$x \mapsto (x +  j \adfmod{33}) + 132$ for $132 \le x < 165$,
$x \mapsto (x +  j \adfmod{3}) + 165$ for $165 \le x < 168$,
$168 \mapsto 168$,
$0 \le j < 132$
 for the first 15 blocks,
$0 \le j < 44$
 for the last block.
\ADFvfyParStart{(169, ((15, 132, ((132, 1), (33, 1), (3, 1), (1, 1))), (1, 44, ((132, 1), (33, 1), (3, 1), (1, 1)))), ((22, 6), (37, 1)))} 

\adfDgap
\noindent{\boldmath $ 22^{6} 40^{1} $}~
With the point set $Z_{172}$ partitioned into
 residue classes modulo $6$ for $\{0, 1, \dots, 131\}$, and
 $\{132, 133, \dots, 171\}$,
 the design is generated from

\adfLgap 
$(165, 49, 65, 78)$,
$(165, 128, 10, 81)$,
$(166, 72, 19, 14)$,
$(166, 123, 47, 58)$,\adfsplit
$(132, 42, 63, 38)$,
$(132, 96, 39, 58)$,
$(132, 1, 118, 126)$,
$(132, 128, 66, 49)$,\adfsplit
$(132, 87, 77, 64)$,
$(132, 131, 79, 20)$,
$(132, 29, 94, 97)$,
$(132, 102, 103, 70)$,\adfsplit
$(132, 54, 82, 109)$,
$(132, 33, 56, 76)$,
$(132, 68, 17, 114)$,
$(132, 88, 98, 127)$,\adfsplit
$(132, 107, 116, 34)$,
$(132, 117, 85, 125)$,
$(0, 22, 129, 156)$,
$(0, 16, 56, 117)$,\adfsplit
$(0, 49, 87, 162)$,
$(0, 2, 11, 153)$,
$(0, 41, 45, 159)$,
$(1, 21, 71, 164)$,\adfsplit
$(0, 17, 51, 52)$,
$(0, 37, 63, 161)$,
$(0, 31, 68, 143)$,
$(0, 26, 69, 83)$,\adfsplit
$(0, 34, 125, 127)$,
$(0, 77, 99, 158)$,
$(0, 19, 47, 105)$,
$(171, 0, 44, 88)$,\adfsplit
$(171, 1, 45, 89)$

\adfLgap \noindent by the mapping:
$x \mapsto x + 2 j \adfmod{132}$ for $x < 132$,
$x \mapsto (x +  j \adfmod{33}) + 132$ for $132 \le x < 165$,
$x \mapsto (x - 165 + 2 j \adfmod{6}) + 165$ for $165 \le x < 171$,
$171 \mapsto 171$,
$0 \le j < 66$
 for the first 31 blocks,
$0 \le j < 22$
 for the last two blocks.
\ADFvfyParStart{(172, ((31, 66, ((132, 2), (33, 1), (6, 2), (1, 1))), (2, 22, ((132, 2), (33, 1), (6, 2), (1, 1)))), ((22, 6), (40, 1)))} 

\adfDgap
\noindent{\boldmath $ 22^{6} 43^{1} $}~
With the point set $Z_{175}$ partitioned into
 residue classes modulo $6$ for $\{0, 1, \dots, 131\}$, and
 $\{132, 133, \dots, 174\}$,
 the design is generated from

\adfLgap 
$(165, 7, 33, 71)$,
$(166, 105, 62, 28)$,
$(167, 14, 79, 21)$,
$(132, 105, 64, 113)$,\adfsplit
$(132, 62, 63, 49)$,
$(132, 81, 12, 104)$,
$(132, 103, 112, 93)$,
$(132, 127, 51, 66)$,\adfsplit
$(0, 2, 5, 33)$,
$(0, 4, 39, 50)$,
$(0, 17, 62, 139)$,
$(0, 22, 51, 163)$,\adfsplit
$(0, 21, 80, 154)$,
$(0, 25, 57, 148)$,
$(0, 20, 105, 164)$,
$(0, 16, 53, 162)$,\adfsplit
$(174, 0, 44, 88)$

\adfLgap \noindent by the mapping:
$x \mapsto x +  j \adfmod{132}$ for $x < 132$,
$x \mapsto (x +  j \adfmod{33}) + 132$ for $132 \le x < 165$,
$x \mapsto (x - 165 + 3 j \adfmod{9}) + 165$ for $165 \le x < 174$,
$174 \mapsto 174$,
$0 \le j < 132$
 for the first 16 blocks,
$0 \le j < 44$
 for the last block.
\ADFvfyParStart{(175, ((16, 132, ((132, 1), (33, 1), (9, 3), (1, 1))), (1, 44, ((132, 1), (33, 1), (9, 3), (1, 1)))), ((22, 6), (43, 1)))} 

\adfDgap
\noindent{\boldmath $ 22^{6} 46^{1} $}~
With the point set $Z_{178}$ partitioned into
 residue classes modulo $6$ for $\{0, 1, \dots, 131\}$, and
 $\{132, 133, \dots, 177\}$,
 the design is generated from

\adfLgap 
$(165, 102, 49, 44)$,
$(165, 47, 22, 75)$,
$(166, 4, 71, 56)$,
$(166, 7, 123, 54)$,\adfsplit
$(167, 116, 127, 71)$,
$(167, 118, 63, 96)$,
$(168, 86, 61, 66)$,
$(168, 52, 53, 123)$,\adfsplit
$(132, 26, 112, 30)$,
$(132, 110, 10, 107)$,
$(132, 63, 72, 29)$,
$(132, 88, 115, 119)$,\adfsplit
$(132, 122, 124, 37)$,
$(132, 86, 109, 45)$,
$(132, 108, 70, 31)$,
$(132, 105, 85, 66)$,\adfsplit
$(132, 117, 77, 36)$,
$(132, 99, 38, 106)$,
$(132, 100, 74, 84)$,
$(132, 67, 65, 82)$,\adfsplit
$(0, 21, 73, 164)$,
$(0, 7, 76, 140)$,
$(0, 75, 113, 159)$,
$(0, 8, 59, 109)$,\adfsplit
$(0, 3, 29, 139)$,
$(0, 121, 131, 151)$,
$(0, 9, 92, 133)$,
$(0, 37, 95, 149)$,\adfsplit
$(0, 13, 34, 62)$,
$(0, 14, 119, 141)$,
$(0, 17, 103, 138)$,
$(0, 33, 65, 158)$,\adfsplit
$(0, 35, 43, 57)$,
$(177, 0, 44, 88)$,
$(177, 1, 45, 89)$

\adfLgap \noindent by the mapping:
$x \mapsto x + 2 j \adfmod{132}$ for $x < 132$,
$x \mapsto (x +  j \adfmod{33}) + 132$ for $132 \le x < 165$,
$x \mapsto (x - 165 + 4 j \adfmod{12}) + 165$ for $165 \le x < 177$,
$177 \mapsto 177$,
$0 \le j < 66$
 for the first 33 blocks,
$0 \le j < 22$
 for the last two blocks.
\ADFvfyParStart{(178, ((33, 66, ((132, 2), (33, 1), (12, 4), (1, 1))), (2, 22, ((132, 2), (33, 1), (12, 4), (1, 1)))), ((22, 6), (46, 1)))} 

\adfDgap
\noindent{\boldmath $ 22^{6} 49^{1} $}~
With the point set $Z_{181}$ partitioned into
 residue classes modulo $6$ for $\{0, 1, \dots, 131\}$, and
 $\{132, 133, \dots, 180\}$,
 the design is generated from

\adfLgap 
$(165, 31, 54, 56)$,
$(166, 54, 82, 116)$,
$(167, 23, 126, 1)$,
$(168, 121, 65, 24)$,\adfsplit
$(169, 79, 78, 116)$,
$(132, 94, 49, 75)$,
$(132, 126, 79, 62)$,
$(132, 72, 45, 59)$,\adfsplit
$(132, 73, 77, 10)$,
$(0, 10, 43, 92)$,
$(0, 9, 121, 133)$,
$(0, 5, 80, 146)$,\adfsplit
$(0, 8, 81, 150)$,
$(0, 21, 53, 151)$,
$(0, 16, 93, 161)$,
$(0, 3, 61, 143)$,\adfsplit
$(0, 15, 101, 162)$,
$(180, 0, 44, 88)$

\adfLgap \noindent by the mapping:
$x \mapsto x +  j \adfmod{132}$ for $x < 132$,
$x \mapsto (x +  j \adfmod{33}) + 132$ for $132 \le x < 165$,
$x \mapsto (x + 5 j \adfmod{15}) + 165$ for $165 \le x < 180$,
$180 \mapsto 180$,
$0 \le j < 132$
 for the first 17 blocks,
$0 \le j < 44$
 for the last block.
\ADFvfyParStart{(181, ((17, 132, ((132, 1), (33, 1), (15, 5), (1, 1))), (1, 44, ((132, 1), (33, 1), (15, 5), (1, 1)))), ((22, 6), (49, 1)))} 

\adfDgap
\noindent{\boldmath $ 22^{6} 52^{1} $}~
With the point set $Z_{184}$ partitioned into
 residue classes modulo $6$ for $\{0, 1, \dots, 131\}$, and
 $\{132, 133, \dots, 183\}$,
 the design is generated from

\adfLgap 
$(165, 124, 57, 86)$,
$(165, 73, 5, 12)$,
$(166, 116, 114, 100)$,
$(166, 37, 29, 69)$,\adfsplit
$(167, 26, 107, 102)$,
$(167, 69, 124, 31)$,
$(168, 12, 1, 35)$,
$(168, 21, 40, 80)$,\adfsplit
$(169, 62, 97, 23)$,
$(169, 30, 100, 75)$,
$(170, 36, 125, 26)$,
$(170, 85, 100, 69)$,\adfsplit
$(132, 93, 80, 97)$,
$(132, 120, 127, 57)$,
$(132, 70, 17, 44)$,
$(132, 94, 20, 73)$,\adfsplit
$(132, 128, 100, 111)$,
$(132, 41, 106, 8)$,
$(132, 77, 114, 2)$,
$(132, 46, 129, 23)$,\adfsplit
$(0, 3, 80, 139)$,
$(0, 1, 4, 149)$,
$(0, 15, 25, 100)$,
$(0, 41, 127, 133)$,\adfsplit
$(0, 71, 85, 150)$,
$(0, 46, 97, 146)$,
$(0, 8, 29, 140)$,
$(0, 43, 119, 144)$,\adfsplit
$(0, 19, 131, 162)$,
$(0, 59, 68, 154)$,
$(0, 63, 91, 159)$,
$(0, 37, 82, 160)$,\adfsplit
$(1, 3, 83, 140)$,
$(0, 22, 49, 135)$,
$(0, 9, 31, 137)$,
$(183, 0, 44, 88)$,\adfsplit
$(183, 1, 45, 89)$

\adfLgap \noindent by the mapping:
$x \mapsto x + 2 j \adfmod{132}$ for $x < 132$,
$x \mapsto (x +  j \adfmod{33}) + 132$ for $132 \le x < 165$,
$x \mapsto (x - 165 + 6 j \adfmod{18}) + 165$ for $165 \le x < 183$,
$183 \mapsto 183$,
$0 \le j < 66$
 for the first 35 blocks,
$0 \le j < 22$
 for the last two blocks.
\ADFvfyParStart{(184, ((35, 66, ((132, 2), (33, 1), (18, 6), (1, 1))), (2, 22, ((132, 2), (33, 1), (18, 6), (1, 1)))), ((22, 6), (52, 1)))} 

\adfDgap
\noindent{\boldmath $ 22^{9} 7^{1} $}~
With the point set $Z_{205}$ partitioned into
 residue classes modulo $9$ for $\{0, 1, \dots, 197\}$, and
 $\{198, 199, \dots, 204\}$,
 the design is generated from

\adfLgap 
$(198, 77, 109, 112)$,
$(198, 81, 150, 98)$,
$(199, 143, 196, 3)$,
$(199, 73, 84, 92)$,\adfsplit
$(99, 132, 155, 95)$,
$(63, 157, 102, 43)$,
$(112, 143, 66, 24)$,
$(101, 183, 149, 58)$,\adfsplit
$(24, 144, 95, 10)$,
$(59, 38, 179, 150)$,
$(23, 53, 142, 141)$,
$(149, 53, 175, 2)$,\adfsplit
$(169, 154, 155, 130)$,
$(141, 54, 10, 20)$,
$(92, 60, 122, 72)$,
$(157, 163, 150, 74)$,\adfsplit
$(79, 76, 174, 134)$,
$(174, 176, 60, 161)$,
$(0, 6, 53, 176)$,
$(0, 37, 68, 138)$,\adfsplit
$(0, 69, 74, 169)$,
$(0, 4, 127, 137)$,
$(0, 26, 59, 106)$,
$(0, 16, 57, 191)$,\adfsplit
$(0, 11, 96, 157)$,
$(0, 17, 67, 150)$,
$(0, 35, 56, 160)$,
$(0, 19, 93, 105)$,\adfsplit
$(0, 49, 111, 155)$,
$(1, 9, 25, 47)$,
$(1, 3, 43, 71)$,
$(204, 0, 66, 132)$,\adfsplit
$(204, 1, 67, 133)$

\adfLgap \noindent by the mapping:
$x \mapsto x + 2 j \adfmod{198}$ for $x < 198$,
$x \mapsto (x + 2 j \adfmod{6}) + 198$ for $198 \le x < 204$,
$204 \mapsto 204$,
$0 \le j < 99$
 for the first 31 blocks,
$0 \le j < 33$
 for the last two blocks.
\ADFvfyParStart{(205, ((31, 99, ((198, 2), (6, 2), (1, 1))), (2, 33, ((198, 2), (6, 2), (1, 1)))), ((22, 9), (7, 1)))} 

\adfDgap
\noindent{\boldmath $ 22^{9} 10^{1} $}~
With the point set $Z_{208}$ partitioned into
 residue classes modulo $9$ for $\{0, 1, \dots, 197\}$, and
 $\{198, 199, \dots, 207\}$,
 the design is generated from

\adfLgap 
$(198, 67, 116, 108)$,
$(198, 181, 42, 88)$,
$(198, 5, 56, 174)$,
$(17, 114, 133, 23)$,\adfsplit
$(22, 155, 178, 53)$,
$(102, 123, 70, 90)$,
$(111, 180, 59, 137)$,
$(193, 178, 5, 63)$,\adfsplit
$(0, 1, 3, 194)$,
$(0, 11, 35, 48)$,
$(0, 16, 71, 138)$,
$(0, 38, 85, 141)$,\adfsplit
$(0, 30, 64, 104)$,
$(0, 14, 84, 112)$,
$(0, 22, 61, 111)$,
$(0, 17, 79, 123)$,\adfsplit
$(207, 0, 66, 132)$

\adfLgap \noindent by the mapping:
$x \mapsto x +  j \adfmod{198}$ for $x < 198$,
$x \mapsto (x +  j \adfmod{9}) + 198$ for $198 \le x < 207$,
$207 \mapsto 207$,
$0 \le j < 198$
 for the first 16 blocks,
$0 \le j < 66$
 for the last block.
\ADFvfyParStart{(208, ((16, 198, ((198, 1), (9, 1), (1, 1))), (1, 66, ((198, 1), (9, 1), (1, 1)))), ((22, 9), (10, 1)))} 

\adfDgap
\noindent{\boldmath $ 22^{9} 13^{1} $}~
With the point set $Z_{211}$ partitioned into
 residue classes modulo $9$ for $\{0, 1, \dots, 197\}$, and
 $\{198, 199, \dots, 210\}$,
 the design is generated from

\adfLgap 
$(207, 59, 85, 70)$,
$(207, 68, 147, 126)$,
$(198, 43, 31, 135)$,
$(198, 173, 62, 176)$,\adfsplit
$(198, 18, 100, 78)$,
$(198, 111, 148, 92)$,
$(198, 33, 41, 73)$,
$(198, 84, 125, 70)$,\adfsplit
$(43, 183, 14, 182)$,
$(155, 149, 62, 145)$,
$(174, 16, 68, 153)$,
$(39, 125, 172, 140)$,\adfsplit
$(155, 84, 153, 7)$,
$(129, 185, 54, 118)$,
$(169, 31, 53, 126)$,
$(194, 90, 98, 123)$,\adfsplit
$(135, 93, 190, 2)$,
$(65, 186, 93, 19)$,
$(176, 58, 134, 46)$,
$(29, 113, 60, 64)$,\adfsplit
$(0, 2, 39, 109)$,
$(0, 16, 73, 191)$,
$(0, 5, 35, 155)$,
$(0, 7, 23, 48)$,\adfsplit
$(0, 24, 147, 148)$,
$(0, 6, 44, 185)$,
$(0, 28, 89, 128)$,
$(0, 20, 46, 165)$,\adfsplit
$(0, 3, 62, 113)$,
$(0, 129, 149, 193)$,
$(0, 13, 47, 115)$,
$(1, 15, 39, 137)$,\adfsplit
$(0, 17, 34, 120)$,
$(210, 0, 66, 132)$,
$(210, 1, 67, 133)$

\adfLgap \noindent by the mapping:
$x \mapsto x + 2 j \adfmod{198}$ for $x < 198$,
$x \mapsto (x +  j \adfmod{9}) + 198$ for $198 \le x < 207$,
$x \mapsto (x +  j \adfmod{3}) + 207$ for $207 \le x < 210$,
$210 \mapsto 210$,
$0 \le j < 99$
 for the first 33 blocks,
$0 \le j < 33$
 for the last two blocks.
\ADFvfyParStart{(211, ((33, 99, ((198, 2), (9, 1), (3, 1), (1, 1))), (2, 33, ((198, 2), (9, 1), (3, 1), (1, 1)))), ((22, 9), (13, 1)))} 

\adfDgap
\noindent{\boldmath $ 22^{9} 16^{1} $}~
With the point set $Z_{214}$ partitioned into
 residue classes modulo $9$ for $\{0, 1, \dots, 197\}$, and
 $\{198, 199, \dots, 213\}$,
 the design is generated from

\adfLgap 
$(207, 127, 80, 6)$,
$(208, 15, 2, 100)$,
$(198, 30, 157, 99)$,
$(198, 41, 15, 190)$,\adfsplit
$(198, 137, 97, 62)$,
$(141, 174, 145, 40)$,
$(119, 5, 117, 112)$,
$(23, 173, 105, 93)$,\adfsplit
$(112, 191, 60, 17)$,
$(0, 1, 15, 111)$,
$(0, 37, 83, 143)$,
$(0, 16, 57, 78)$,\adfsplit
$(0, 6, 28, 38)$,
$(0, 39, 89, 133)$,
$(0, 11, 30, 167)$,
$(0, 8, 25, 59)$,\adfsplit
$(0, 3, 56, 76)$,
$(213, 0, 66, 132)$

\adfLgap \noindent by the mapping:
$x \mapsto x +  j \adfmod{198}$ for $x < 198$,
$x \mapsto (x +  j \adfmod{9}) + 198$ for $198 \le x < 207$,
$x \mapsto (x - 207 + 2 j \adfmod{6}) + 207$ for $207 \le x < 213$,
$213 \mapsto 213$,
$0 \le j < 198$
 for the first 17 blocks,
$0 \le j < 66$
 for the last block.
\ADFvfyParStart{(214, ((17, 198, ((198, 1), (9, 1), (6, 2), (1, 1))), (1, 66, ((198, 1), (9, 1), (6, 2), (1, 1)))), ((22, 9), (16, 1)))} 

\adfDgap
\noindent{\boldmath $ 22^{9} 19^{1} $}~
With the point set $Z_{217}$ partitioned into
 residue classes modulo $9$ for $\{0, 1, \dots, 197\}$, and
 $\{198, 199, \dots, 216\}$,
 the design is generated from

\adfLgap 
$(198, 118, 11, 134)$,
$(198, 110, 93, 131)$,
$(198, 186, 161, 157)$,
$(198, 112, 97, 102)$,\adfsplit
$(198, 195, 81, 55)$,
$(198, 86, 106, 54)$,
$(199, 14, 152, 12)$,
$(199, 100, 74, 99)$,\adfsplit
$(199, 7, 17, 103)$,
$(199, 36, 178, 129)$,
$(199, 150, 40, 113)$,
$(199, 83, 159, 181)$,\adfsplit
$(142, 107, 30, 67)$,
$(56, 100, 187, 159)$,
$(16, 195, 17, 164)$,
$(119, 22, 45, 105)$,\adfsplit
$(74, 17, 9, 87)$,
$(43, 41, 184, 189)$,
$(122, 89, 150, 20)$,
$(36, 107, 187, 28)$,\adfsplit
$(141, 107, 138, 59)$,
$(7, 144, 60, 53)$,
$(195, 148, 3, 152)$,
$(0, 29, 121, 185)$,\adfsplit
$(0, 17, 41, 85)$,
$(0, 7, 19, 62)$,
$(0, 15, 109, 125)$,
$(0, 95, 127, 157)$,\adfsplit
$(0, 35, 38, 177)$,
$(0, 6, 39, 124)$,
$(0, 14, 48, 115)$,
$(0, 11, 98, 120)$,\adfsplit
$(0, 24, 129, 152)$,
$(0, 40, 122, 187)$,
$(0, 12, 42, 106)$,
$(216, 0, 66, 132)$,\adfsplit
$(216, 1, 67, 133)$

\adfLgap \noindent by the mapping:
$x \mapsto x + 2 j \adfmod{198}$ for $x < 198$,
$x \mapsto (x + 2 j \adfmod{18}) + 198$ for $198 \le x < 216$,
$216 \mapsto 216$,
$0 \le j < 99$
 for the first 35 blocks,
$0 \le j < 33$
 for the last two blocks.
\ADFvfyParStart{(217, ((35, 99, ((198, 2), (18, 2), (1, 1))), (2, 33, ((198, 2), (18, 2), (1, 1)))), ((22, 9), (19, 1)))} 

\adfDgap
\noindent{\boldmath $ 22^{12} 7^{1} $}~
With the point set $Z_{271}$ partitioned into
 residue classes modulo $12$ for $\{0, 1, \dots, 263\}$, and
 $\{264, 265, \dots, 270\}$,
 the design is generated from

\adfLgap 
$(264, 73, 15, 233)$,
$(265, 194, 220, 105)$,
$(137, 81, 66, 200)$,
$(142, 193, 218, 260)$,\adfsplit
$(131, 53, 25, 122)$,
$(213, 60, 26, 139)$,
$(43, 93, 61, 255)$,
$(138, 121, 107, 230)$,\adfsplit
$(151, 60, 13, 56)$,
$(0, 1, 3, 136)$,
$(0, 5, 11, 112)$,
$(0, 7, 20, 105)$,\adfsplit
$(0, 8, 35, 122)$,
$(0, 29, 59, 207)$,
$(0, 40, 93, 161)$,
$(0, 39, 80, 220)$,\adfsplit
$(0, 16, 81, 198)$,
$(0, 10, 33, 55)$,
$(0, 49, 110, 174)$,
$(0, 19, 73, 94)$,\adfsplit
$(0, 37, 99, 137)$,
$(270, 0, 88, 176)$

\adfLgap \noindent by the mapping:
$x \mapsto x +  j \adfmod{264}$ for $x < 264$,
$x \mapsto (x + 2 j \adfmod{6}) + 264$ for $264 \le x < 270$,
$270 \mapsto 270$,
$0 \le j < 264$
 for the first 21 blocks,
$0 \le j < 88$
 for the last block.
\ADFvfyParStart{(271, ((21, 264, ((264, 1), (6, 2), (1, 1))), (1, 88, ((264, 1), (6, 2), (1, 1)))), ((22, 12), (7, 1)))} 

\adfDgap
\noindent{\boldmath $ 22^{12} 10^{1} $}~
With the point set $Z_{274}$ partitioned into
 residue classes modulo $12$ for $\{0, 1, \dots, 263\}$, and
 $\{264, 265, \dots, 273\}$,
 the design is generated from

\adfLgap 
$(264, 122, 28, 83)$,
$(264, 7, 204, 3)$,
$(265, 30, 148, 83)$,
$(265, 243, 224, 55)$,\adfsplit
$(266, 244, 116, 29)$,
$(266, 186, 175, 201)$,
$(167, 44, 136, 113)$,
$(126, 14, 17, 176)$,\adfsplit
$(167, 256, 218, 198)$,
$(90, 209, 170, 177)$,
$(89, 27, 104, 14)$,
$(58, 253, 117, 87)$,\adfsplit
$(135, 126, 257, 241)$,
$(49, 26, 243, 233)$,
$(191, 255, 120, 20)$,
$(67, 88, 198, 94)$,\adfsplit
$(7, 181, 233, 160)$,
$(48, 52, 116, 161)$,
$(72, 19, 206, 80)$,
$(223, 61, 107, 105)$,\adfsplit
$(224, 154, 182, 180)$,
$(21, 27, 164, 257)$,
$(182, 219, 136, 163)$,
$(30, 87, 191, 40)$,\adfsplit
$(117, 247, 252, 23)$,
$(263, 85, 19, 110)$,
$(0, 1, 106, 251)$,
$(0, 5, 14, 117)$,\adfsplit
$(0, 11, 165, 232)$,
$(0, 54, 143, 202)$,
$(0, 40, 101, 182)$,
$(0, 30, 179, 257)$,\adfsplit
$(0, 22, 191, 208)$,
$(0, 16, 201, 209)$,
$(0, 17, 167, 189)$,
$(0, 79, 219, 261)$,\adfsplit
$(0, 34, 85, 212)$,
$(0, 18, 99, 157)$,
$(0, 97, 98, 223)$,
$(0, 163, 181, 231)$,\adfsplit
$(0, 25, 74, 198)$,
$(0, 41, 76, 141)$,
$(0, 33, 107, 150)$,
$(273, 0, 88, 176)$,\adfsplit
$(273, 1, 89, 177)$

\adfLgap \noindent by the mapping:
$x \mapsto x + 2 j \adfmod{264}$ for $x < 264$,
$x \mapsto (x - 264 + 3 j \adfmod{9}) + 264$ for $264 \le x < 273$,
$273 \mapsto 273$,
$0 \le j < 132$
 for the first 43 blocks,
$0 \le j < 44$
 for the last two blocks.
\ADFvfyParStart{(274, ((43, 132, ((264, 2), (9, 3), (1, 1))), (2, 44, ((264, 2), (9, 3), (1, 1)))), ((22, 12), (10, 1)))} 

\adfDgap
\noindent{\boldmath $ 22^{12} 13^{1} $}~
With the point set $Z_{277}$ partitioned into
 residue classes modulo $12$ for $\{0, 1, \dots, 263\}$, and
 $\{264, 265, \dots, 276\}$,
 the design is generated from

\adfLgap 
$(264, 242, 151, 117)$,
$(265, 165, 92, 244)$,
$(266, 229, 26, 252)$,
$(267, 21, 88, 143)$,\adfsplit
$(232, 158, 121, 245)$,
$(112, 23, 14, 91)$,
$(188, 139, 113, 93)$,
$(68, 132, 183, 225)$,\adfsplit
$(260, 194, 81, 263)$,
$(84, 52, 114, 190)$,
$(0, 1, 5, 211)$,
$(0, 2, 8, 52)$,\adfsplit
$(0, 7, 17, 97)$,
$(0, 11, 29, 236)$,
$(0, 27, 92, 155)$,
$(0, 31, 71, 141)$,\adfsplit
$(0, 15, 114, 133)$,
$(0, 43, 102, 147)$,
$(0, 14, 100, 135)$,
$(0, 33, 116, 163)$,\adfsplit
$(0, 16, 41, 186)$,
$(0, 22, 103, 159)$,
$(276, 0, 88, 176)$

\adfLgap \noindent by the mapping:
$x \mapsto x +  j \adfmod{264}$ for $x < 264$,
$x \mapsto (x + 4 j \adfmod{12}) + 264$ for $264 \le x < 276$,
$276 \mapsto 276$,
$0 \le j < 264$
 for the first 22 blocks,
$0 \le j < 88$
 for the last block.
\ADFvfyParStart{(277, ((22, 264, ((264, 1), (12, 4), (1, 1))), (1, 88, ((264, 1), (12, 4), (1, 1)))), ((22, 12), (13, 1)))} 

\adfDgap
\noindent{\boldmath $ 22^{12} 16^{1} $}~
With the point set $Z_{280}$ partitioned into
 residue classes modulo $12$ for $\{0, 1, \dots, 263\}$, and
 $\{264, 265, \dots, 279\}$,
 the design is generated from

\adfLgap 
$(264, 29, 118, 68)$,
$(264, 163, 120, 21)$,
$(265, 83, 28, 63)$,
$(265, 182, 228, 97)$,\adfsplit
$(266, 53, 181, 66)$,
$(266, 129, 238, 122)$,
$(267, 90, 117, 161)$,
$(267, 253, 100, 8)$,\adfsplit
$(268, 254, 111, 178)$,
$(268, 143, 216, 109)$,
$(30, 46, 63, 124)$,
$(220, 260, 190, 93)$,\adfsplit
$(68, 249, 67, 54)$,
$(5, 216, 82, 14)$,
$(59, 150, 28, 159)$,
$(170, 35, 139, 111)$,\adfsplit
$(13, 198, 160, 59)$,
$(184, 27, 9, 80)$,
$(219, 230, 18, 67)$,
$(123, 133, 71, 226)$,\adfsplit
$(239, 55, 145, 185)$,
$(154, 54, 260, 195)$,
$(172, 226, 67, 254)$,
$(17, 175, 143, 222)$,\adfsplit
$(215, 53, 145, 151)$,
$(218, 65, 189, 216)$,
$(205, 148, 91, 204)$,
$(154, 160, 39, 80)$,\adfsplit
$(0, 3, 11, 189)$,
$(0, 8, 18, 83)$,
$(0, 19, 26, 138)$,
$(0, 25, 150, 231)$,\adfsplit
$(0, 15, 20, 110)$,
$(0, 37, 166, 227)$,
$(0, 213, 215, 229)$,
$(0, 22, 222, 261)$,\adfsplit
$(0, 21, 87, 102)$,
$(0, 44, 91, 147)$,
$(0, 23, 45, 66)$,
$(0, 29, 67, 177)$,\adfsplit
$(0, 4, 123, 140)$,
$(0, 34, 85, 127)$,
$(0, 5, 73, 219)$,
$(0, 69, 95, 99)$,\adfsplit
$(0, 32, 178, 241)$,
$(279, 0, 88, 176)$,
$(279, 1, 89, 177)$

\adfLgap \noindent by the mapping:
$x \mapsto x + 2 j \adfmod{264}$ for $x < 264$,
$x \mapsto (x - 264 + 5 j \adfmod{15}) + 264$ for $264 \le x < 279$,
$279 \mapsto 279$,
$0 \le j < 132$
 for the first 45 blocks,
$0 \le j < 44$
 for the last two blocks.
\ADFvfyParStart{(280, ((45, 132, ((264, 2), (15, 5), (1, 1))), (2, 44, ((264, 2), (15, 5), (1, 1)))), ((22, 12), (16, 1)))} 

\adfDgap
\noindent{\boldmath $ 22^{12} 19^{1} $}~
With the point set $Z_{283}$ partitioned into
 residue classes modulo $12$ for $\{0, 1, \dots, 263\}$, and
 $\{264, 265, \dots, 282\}$,
 the design is generated from

\adfLgap 
$(264, 243, 185, 172)$,
$(265, 259, 60, 137)$,
$(266, 106, 63, 86)$,
$(267, 66, 133, 155)$,\adfsplit
$(268, 262, 237, 134)$,
$(269, 227, 120, 88)$,
$(169, 155, 22, 4)$,
$(132, 91, 161, 98)$,\adfsplit
$(101, 193, 196, 19)$,
$(234, 50, 133, 29)$,
$(188, 179, 36, 64)$,
$(0, 1, 5, 191)$,\adfsplit
$(0, 2, 40, 66)$,
$(0, 6, 37, 52)$,
$(0, 8, 47, 57)$,
$(0, 45, 100, 150)$,\adfsplit
$(0, 30, 91, 188)$,
$(0, 11, 62, 141)$,
$(0, 42, 98, 179)$,
$(0, 19, 94, 148)$,\adfsplit
$(0, 33, 68, 178)$,
$(0, 17, 44, 155)$,
$(0, 16, 69, 162)$,
$(282, 0, 88, 176)$

\adfLgap \noindent by the mapping:
$x \mapsto x +  j \adfmod{264}$ for $x < 264$,
$x \mapsto (x - 264 + 6 j \adfmod{18}) + 264$ for $264 \le x < 282$,
$282 \mapsto 282$,
$0 \le j < 264$
 for the first 23 blocks,
$0 \le j < 88$
 for the last block.
\ADFvfyParStart{(283, ((23, 264, ((264, 1), (18, 6), (1, 1))), (1, 88, ((264, 1), (18, 6), (1, 1)))), ((22, 12), (19, 1)))} 

\adfDgap
\noindent{\boldmath $ 22^{15} 4^{1} $}~
With the point set $Z_{334}$ partitioned into
 residue classes modulo $15$ for $\{0, 1, \dots, 329\}$, and
 $\{330, 331, 332, 333\}$,
 the design is generated from

\adfLgap 
$(330, 98, 211, 138)$,
$(193, 86, 64, 143)$,
$(176, 2, 196, 297)$,
$(77, 267, 145, 19)$,\adfsplit
$(265, 232, 83, 318)$,
$(305, 26, 294, 318)$,
$(40, 120, 212, 213)$,
$(288, 47, 143, 329)$,\adfsplit
$(171, 29, 181, 188)$,
$(33, 52, 264, 38)$,
$(314, 72, 137, 90)$,
$(244, 195, 286, 192)$,\adfsplit
$(274, 278, 112, 14)$,
$(83, 199, 62, 329)$,
$(0, 2, 74, 83)$,
$(0, 6, 31, 302)$,\adfsplit
$(0, 8, 44, 211)$,
$(0, 12, 67, 191)$,
$(0, 16, 103, 215)$,
$(0, 37, 76, 184)$,\adfsplit
$(0, 61, 125, 259)$,
$(0, 32, 109, 187)$,
$(0, 26, 111, 228)$,
$(0, 43, 97, 197)$,\adfsplit
$(0, 23, 69, 192)$,
$(0, 27, 56, 216)$,
$(333, 0, 110, 220)$

\adfLgap \noindent by the mapping:
$x \mapsto x +  j \adfmod{330}$ for $x < 330$,
$x \mapsto (x +  j \adfmod{3}) + 330$ for $330 \le x < 333$,
$333 \mapsto 333$,
$0 \le j < 330$
 for the first 26 blocks,
$0 \le j < 110$
 for the last block.
\ADFvfyParStart{(334, ((26, 330, ((330, 1), (3, 1), (1, 1))), (1, 110, ((330, 1), (3, 1), (1, 1)))), ((22, 15), (4, 1)))} 

\adfDgap
\noindent{\boldmath $ 22^{15} 10^{1} $}~
With the point set $Z_{340}$ partitioned into
 residue classes modulo $15$ for $\{0, 1, \dots, 329\}$, and
 $\{330, 331, \dots, 339\}$,
 the design is generated from

\adfLgap 
$(330, 270, 257, 43)$,
$(331, 271, 65, 231)$,
$(332, 153, 152, 103)$,
$(26, 201, 47, 138)$,\adfsplit
$(55, 18, 174, 267)$,
$(148, 134, 226, 219)$,
$(238, 131, 327, 3)$,
$(226, 248, 253, 22)$,\adfsplit
$(59, 211, 268, 98)$,
$(211, 17, 52, 201)$,
$(253, 241, 179, 311)$,
$(265, 218, 134, 242)$,\adfsplit
$(154, 198, 150, 170)$,
$(229, 270, 91, 8)$,
$(239, 231, 34, 116)$,
$(0, 2, 19, 294)$,\adfsplit
$(0, 3, 34, 59)$,
$(0, 11, 122, 188)$,
$(0, 32, 101, 201)$,
$(0, 26, 72, 213)$,\adfsplit
$(0, 18, 97, 162)$,
$(0, 43, 130, 228)$,
$(0, 42, 94, 182)$,
$(0, 9, 76, 172)$,\adfsplit
$(0, 54, 127, 191)$,
$(0, 29, 80, 253)$,
$(0, 33, 86, 147)$,
$(339, 0, 110, 220)$

\adfLgap \noindent by the mapping:
$x \mapsto x +  j \adfmod{330}$ for $x < 330$,
$x \mapsto (x - 330 + 3 j \adfmod{9}) + 330$ for $330 \le x < 339$,
$339 \mapsto 339$,
$0 \le j < 330$
 for the first 27 blocks,
$0 \le j < 110$
 for the last block.
\ADFvfyParStart{(340, ((27, 330, ((330, 1), (9, 3), (1, 1))), (1, 110, ((330, 1), (9, 3), (1, 1)))), ((22, 15), (10, 1)))} 

\adfDgap
\noindent{\boldmath $ 22^{15} 13^{1} $}~
With the point set $Z_{343}$ partitioned into
 residue classes modulo $15$ for $\{0, 1, \dots, 329\}$, and
 $\{330, 331, \dots, 342\}$,
 the design is generated from

\adfLgap 
$(330, 107, 277, 16)$,
$(330, 74, 156, 273)$,
$(331, 274, 237, 80)$,
$(331, 275, 145, 324)$,\adfsplit
$(332, 17, 322, 246)$,
$(332, 213, 283, 206)$,
$(333, 152, 23, 132)$,
$(333, 262, 259, 171)$,\adfsplit
$(305, 6, 119, 307)$,
$(176, 263, 239, 48)$,
$(108, 60, 19, 86)$,
$(181, 102, 148, 169)$,\adfsplit
$(232, 288, 205, 124)$,
$(200, 46, 205, 41)$,
$(300, 242, 283, 86)$,
$(10, 150, 65, 136)$,\adfsplit
$(318, 131, 39, 31)$,
$(223, 321, 194, 45)$,
$(36, 296, 300, 63)$,
$(153, 275, 139, 171)$,\adfsplit
$(270, 304, 247, 13)$,
$(212, 114, 281, 321)$,
$(68, 179, 183, 221)$,
$(56, 273, 197, 299)$,\adfsplit
$(118, 39, 50, 86)$,
$(98, 10, 33, 34)$,
$(214, 67, 45, 162)$,
$(48, 286, 206, 244)$,\adfsplit
$(207, 124, 112, 61)$,
$(280, 206, 162, 134)$,
$(155, 16, 292, 311)$,
$(212, 252, 90, 193)$,\adfsplit
$(169, 120, 323, 205)$,
$(215, 47, 271, 214)$,
$(249, 255, 218, 124)$,
$(92, 139, 236, 320)$,\adfsplit
$(0, 6, 129, 230)$,
$(0, 3, 16, 268)$,
$(0, 10, 142, 192)$,
$(0, 2, 11, 61)$,\adfsplit
$(0, 18, 227, 234)$,
$(0, 53, 130, 275)$,
$(0, 8, 160, 177)$,
$(0, 71, 104, 211)$,\adfsplit
$(0, 35, 119, 218)$,
$(0, 121, 187, 269)$,
$(0, 155, 175, 249)$,
$(0, 173, 219, 287)$,\adfsplit
$(0, 277, 305, 321)$,
$(0, 179, 189, 237)$,
$(0, 89, 205, 291)$,
$(1, 63, 127, 259)$,\adfsplit
$(1, 35, 113, 193)$,
$(0, 13, 137, 244)$,
$(0, 124, 257, 309)$,
$(342, 0, 110, 220)$,\adfsplit
$(342, 1, 111, 221)$

\adfLgap \noindent by the mapping:
$x \mapsto x + 2 j \adfmod{330}$ for $x < 330$,
$x \mapsto (x - 330 + 4 j \adfmod{12}) + 330$ for $330 \le x < 342$,
$342 \mapsto 342$,
$0 \le j < 165$
 for the first 55 blocks,
$0 \le j < 55$
 for the last two blocks.
\ADFvfyParStart{(343, ((55, 165, ((330, 2), (12, 4), (1, 1))), (2, 55, ((330, 2), (12, 4), (1, 1)))), ((22, 15), (13, 1)))} 

\adfDgap
\noindent{\boldmath $ 22^{15} 16^{1} $}~
With the point set $Z_{346}$ partitioned into
 residue classes modulo $15$ for $\{0, 1, \dots, 329\}$, and
 $\{330, 331, \dots, 345\}$,
 the design is generated from

\adfLgap 
$(330, 236, 36, 268)$,
$(331, 158, 121, 204)$,
$(332, 193, 302, 57)$,
$(333, 257, 292, 18)$,\adfsplit
$(334, 180, 175, 251)$,
$(324, 72, 308, 3)$,
$(240, 65, 328, 197)$,
$(186, 190, 213, 305)$,\adfsplit
$(263, 19, 121, 47)$,
$(63, 150, 287, 290)$,
$(203, 54, 232, 225)$,
$(11, 310, 269, 124)$,\adfsplit
$(245, 302, 14, 26)$,
$(241, 28, 262, 288)$,
$(199, 41, 180, 303)$,
$(312, 75, 248, 298)$,\adfsplit
$(0, 1, 39, 63)$,
$(0, 2, 118, 176)$,
$(0, 8, 18, 52)$,
$(0, 17, 82, 250)$,\adfsplit
$(0, 33, 133, 184)$,
$(0, 40, 89, 166)$,
$(0, 6, 59, 183)$,
$(0, 11, 95, 203)$,\adfsplit
$(0, 13, 125, 161)$,
$(0, 20, 121, 187)$,
$(0, 55, 128, 189)$,
$(0, 48, 129, 208)$,\adfsplit
$(345, 0, 110, 220)$

\adfLgap \noindent by the mapping:
$x \mapsto x +  j \adfmod{330}$ for $x < 330$,
$x \mapsto (x + 5 j \adfmod{15}) + 330$ for $330 \le x < 345$,
$345 \mapsto 345$,
$0 \le j < 330$
 for the first 28 blocks,
$0 \le j < 110$
 for the last block.
\ADFvfyParStart{(346, ((28, 330, ((330, 1), (15, 5), (1, 1))), (1, 110, ((330, 1), (15, 5), (1, 1)))), ((22, 15), (16, 1)))} 

\adfDgap
\noindent{\boldmath $ 22^{15} 19^{1} $}~
With the point set $Z_{349}$ partitioned into
 residue classes modulo $15$ for $\{0, 1, \dots, 329\}$, and
 $\{330, 331, \dots, 348\}$,
 the design is generated from

\adfLgap 
$(330, 292, 44, 9)$,
$(330, 205, 252, 173)$,
$(331, 293, 267, 54)$,
$(331, 22, 86, 241)$,\adfsplit
$(332, 124, 0, 267)$,
$(332, 167, 7, 68)$,
$(333, 276, 74, 220)$,
$(333, 293, 63, 259)$,\adfsplit
$(334, 122, 24, 119)$,
$(334, 133, 315, 16)$,
$(335, 200, 229, 321)$,
$(335, 114, 88, 287)$,\adfsplit
$(17, 76, 176, 328)$,
$(275, 129, 28, 150)$,
$(57, 166, 50, 67)$,
$(225, 138, 68, 200)$,\adfsplit
$(11, 264, 30, 323)$,
$(195, 88, 37, 80)$,
$(161, 99, 217, 298)$,
$(235, 124, 226, 108)$,\adfsplit
$(82, 246, 258, 245)$,
$(27, 249, 180, 306)$,
$(38, 52, 43, 315)$,
$(31, 70, 258, 4)$,\adfsplit
$(131, 260, 237, 100)$,
$(313, 19, 152, 260)$,
$(79, 224, 151, 27)$,
$(314, 7, 93, 282)$,\adfsplit
$(289, 150, 163, 25)$,
$(109, 134, 138, 66)$,
$(198, 113, 209, 255)$,
$(298, 322, 143, 183)$,\adfsplit
$(48, 323, 94, 97)$,
$(279, 301, 65, 92)$,
$(245, 143, 154, 221)$,
$(187, 274, 284, 256)$,\adfsplit
$(104, 298, 276, 192)$,
$(178, 59, 107, 65)$,
$(0, 1, 104, 294)$,
$(0, 6, 44, 192)$,\adfsplit
$(0, 33, 250, 313)$,
$(0, 2, 50, 169)$,
$(0, 42, 123, 323)$,
$(0, 79, 156, 241)$,\adfsplit
$(0, 40, 92, 189)$,
$(0, 65, 130, 216)$,
$(0, 35, 218, 272)$,
$(0, 20, 61, 94)$,\adfsplit
$(0, 129, 145, 183)$,
$(0, 34, 73, 168)$,
$(0, 83, 153, 181)$,
$(0, 71, 207, 289)$,\adfsplit
$(0, 21, 223, 237)$,
$(0, 159, 203, 325)$,
$(1, 3, 77, 157)$,
$(1, 5, 69, 89)$,\adfsplit
$(1, 9, 141, 153)$,
$(348, 0, 110, 220)$,
$(348, 1, 111, 221)$

\adfLgap \noindent by the mapping:
$x \mapsto x + 2 j \adfmod{330}$ for $x < 330$,
$x \mapsto (x - 330 + 6 j \adfmod{18}) + 330$ for $330 \le x < 348$,
$348 \mapsto 348$,
$0 \le j < 165$
 for the first 57 blocks,
$0 \le j < 55$
 for the last two blocks.
\ADFvfyParStart{(349, ((57, 165, ((330, 2), (18, 6), (1, 1))), (2, 55, ((330, 2), (18, 6), (1, 1)))), ((22, 15), (19, 1)))} 

\adfDgap
\noindent{\boldmath $ 22^{18} 4^{1} $}~
With the point set $Z_{400}$ partitioned into
 residue classes modulo $18$ for $\{0, 1, \dots, 395\}$, and
 $\{396, 397, 398, 399\}$,
 the design is generated from

\adfLgap 
$(396, 202, 258, 139)$,
$(396, 297, 299, 56)$,
$(297, 106, 138, 296)$,
$(226, 167, 287, 250)$,\adfsplit
$(206, 187, 139, 36)$,
$(192, 33, 78, 124)$,
$(389, 159, 49, 221)$,
$(100, 63, 42, 52)$,\adfsplit
$(339, 81, 133, 57)$,
$(87, 350, 65, 319)$,
$(34, 291, 261, 187)$,
$(40, 306, 103, 209)$,\adfsplit
$(202, 161, 242, 319)$,
$(153, 251, 72, 157)$,
$(144, 84, 22, 223)$,
$(174, 267, 316, 301)$,\adfsplit
$(122, 357, 142, 370)$,
$(381, 169, 234, 299)$,
$(125, 361, 57, 301)$,
$(120, 162, 17, 50)$,\adfsplit
$(224, 103, 59, 86)$,
$(76, 247, 371, 374)$,
$(152, 109, 101, 308)$,
$(376, 168, 66, 289)$,\adfsplit
$(235, 112, 264, 186)$,
$(271, 342, 62, 349)$,
$(141, 70, 1, 114)$,
$(198, 271, 255, 121)$,\adfsplit
$(166, 20, 361, 12)$,
$(322, 217, 186, 373)$,
$(246, 261, 301, 62)$,
$(94, 98, 117, 330)$,\adfsplit
$(245, 266, 330, 156)$,
$(122, 172, 96, 197)$,
$(142, 279, 391, 136)$,
$(274, 96, 257, 141)$,\adfsplit
$(189, 148, 176, 41)$,
$(210, 212, 92, 106)$,
$(118, 106, 283, 325)$,
$(252, 50, 109, 66)$,\adfsplit
$(290, 237, 51, 202)$,
$(288, 23, 125, 122)$,
$(268, 373, 131, 188)$,
$(24, 13, 325, 296)$,\adfsplit
$(0, 5, 196, 265)$,
$(0, 82, 176, 285)$,
$(0, 22, 52, 167)$,
$(0, 150, 361, 371)$,\adfsplit
$(0, 66, 163, 389)$,
$(0, 83, 84, 391)$,
$(0, 9, 217, 256)$,
$(0, 95, 107, 224)$,\adfsplit
$(0, 33, 182, 229)$,
$(0, 38, 134, 173)$,
$(0, 34, 125, 175)$,
$(0, 129, 155, 225)$,\adfsplit
$(0, 67, 87, 281)$,
$(1, 7, 39, 317)$,
$(0, 213, 271, 335)$,
$(1, 29, 129, 175)$,\adfsplit
$(0, 99, 303, 317)$,
$(0, 53, 119, 128)$,
$(0, 92, 273, 296)$,
$(399, 0, 132, 264)$,\adfsplit
$(399, 1, 133, 265)$

\adfLgap \noindent by the mapping:
$x \mapsto x + 2 j \adfmod{396}$ for $x < 396$,
$x \mapsto (x +  j \adfmod{3}) + 396$ for $396 \le x < 399$,
$399 \mapsto 399$,
$0 \le j < 198$
 for the first 63 blocks,
$0 \le j < 66$
 for the last two blocks.
\ADFvfyParStart{(400, ((63, 198, ((396, 2), (3, 1), (1, 1))), (2, 66, ((396, 2), (3, 1), (1, 1)))), ((22, 18), (4, 1)))} 

\adfDgap
\noindent{\boldmath $ 22^{18} 7^{1} $}~
With the point set $Z_{403}$ partitioned into
 residue classes modulo $18$ for $\{0, 1, \dots, 395\}$, and
 $\{396, 397, \dots, 402\}$,
 the design is generated from

\adfLgap 
$(396, 243, 73, 311)$,
$(397, 304, 201, 248)$,
$(207, 317, 222, 320)$,
$(202, 290, 39, 80)$,\adfsplit
$(17, 193, 69, 101)$,
$(251, 32, 355, 165)$,
$(294, 373, 325, 259)$,
$(264, 65, 81, 325)$,\adfsplit
$(341, 30, 10, 343)$,
$(317, 327, 179, 257)$,
$(52, 314, 90, 367)$,
$(380, 139, 89, 379)$,\adfsplit
$(274, 373, 3, 246)$,
$(95, 116, 91, 348)$,
$(207, 229, 338, 196)$,
$(108, 286, 175, 101)$,\adfsplit
$(70, 271, 78, 301)$,
$(167, 143, 201, 244)$,
$(180, 347, 364, 267)$,
$(0, 5, 19, 194)$,\adfsplit
$(0, 6, 55, 302)$,
$(0, 12, 51, 181)$,
$(0, 44, 140, 261)$,
$(0, 40, 82, 236)$,\adfsplit
$(0, 23, 69, 228)$,
$(0, 75, 151, 279)$,
$(0, 9, 71, 116)$,
$(0, 57, 146, 303)$,\adfsplit
$(0, 29, 141, 278)$,
$(0, 27, 91, 214)$,
$(0, 37, 166, 225)$,
$(0, 13, 115, 235)$,\adfsplit
$(402, 0, 132, 264)$

\adfLgap \noindent by the mapping:
$x \mapsto x +  j \adfmod{396}$ for $x < 396$,
$x \mapsto (x + 2 j \adfmod{6}) + 396$ for $396 \le x < 402$,
$402 \mapsto 402$,
$0 \le j < 396$
 for the first 32 blocks,
$0 \le j < 132$
 for the last block.
\ADFvfyParStart{(403, ((32, 396, ((396, 1), (6, 2), (1, 1))), (1, 132, ((396, 1), (6, 2), (1, 1)))), ((22, 18), (7, 1)))} 

\adfDgap
\noindent{\boldmath $ 22^{18} 10^{1} $}~
With the point set $Z_{406}$ partitioned into
 residue classes modulo $18$ for $\{0, 1, \dots, 395\}$, and
 $\{396, 397, \dots, 405\}$,
 the design is generated from

\adfLgap 
$(396, 177, 78, 338)$,
$(396, 346, 379, 179)$,
$(397, 174, 40, 37)$,
$(397, 362, 161, 39)$,\adfsplit
$(398, 38, 303, 106)$,
$(398, 299, 240, 55)$,
$(59, 66, 158, 92)$,
$(319, 143, 218, 46)$,\adfsplit
$(211, 327, 318, 208)$,
$(260, 230, 355, 387)$,
$(146, 118, 287, 311)$,
$(173, 363, 228, 297)$,\adfsplit
$(25, 364, 46, 199)$,
$(89, 178, 115, 255)$,
$(181, 74, 248, 88)$,
$(321, 373, 369, 22)$,\adfsplit
$(337, 6, 71, 90)$,
$(39, 346, 315, 234)$,
$(385, 50, 23, 100)$,
$(136, 216, 215, 379)$,\adfsplit
$(204, 131, 121, 278)$,
$(344, 171, 385, 230)$,
$(40, 87, 32, 67)$,
$(284, 288, 34, 279)$,\adfsplit
$(327, 112, 282, 77)$,
$(355, 41, 150, 137)$,
$(226, 108, 247, 62)$,
$(89, 114, 317, 310)$,\adfsplit
$(298, 303, 143, 356)$,
$(42, 351, 236, 98)$,
$(89, 180, 87, 336)$,
$(310, 135, 276, 175)$,\adfsplit
$(86, 219, 42, 71)$,
$(66, 293, 285, 349)$,
$(140, 324, 153, 141)$,
$(5, 352, 202, 86)$,\adfsplit
$(269, 2, 211, 133)$,
$(140, 283, 116, 179)$,
$(156, 247, 60, 293)$,
$(22, 315, 273, 229)$,\adfsplit
$(158, 312, 264, 331)$,
$(236, 226, 311, 339)$,
$(107, 167, 57, 392)$,
$(298, 106, 374, 257)$,\adfsplit
$(246, 391, 329, 18)$,
$(65, 50, 102, 319)$,
$(0, 6, 22, 298)$,
$(0, 38, 178, 220)$,\adfsplit
$(0, 2, 124, 188)$,
$(0, 12, 32, 94)$,
$(0, 40, 100, 357)$,
$(0, 51, 275, 326)$,\adfsplit
$(0, 71, 230, 353)$,
$(0, 159, 206, 367)$,
$(0, 43, 158, 310)$,
$(0, 102, 291, 379)$,\adfsplit
$(0, 17, 31, 88)$,
$(0, 23, 53, 148)$,
$(0, 37, 105, 263)$,
$(0, 87, 181, 385)$,\adfsplit
$(0, 11, 117, 327)$,
$(0, 25, 109, 179)$,
$(1, 17, 93, 295)$,
$(1, 39, 113, 297)$,\adfsplit
$(1, 7, 135, 157)$,
$(405, 0, 132, 264)$,
$(405, 1, 133, 265)$

\adfLgap \noindent by the mapping:
$x \mapsto x + 2 j \adfmod{396}$ for $x < 396$,
$x \mapsto (x + 3 j \adfmod{9}) + 396$ for $396 \le x < 405$,
$405 \mapsto 405$,
$0 \le j < 198$
 for the first 65 blocks,
$0 \le j < 66$
 for the last two blocks.
\ADFvfyParStart{(406, ((65, 198, ((396, 2), (9, 3), (1, 1))), (2, 66, ((396, 2), (9, 3), (1, 1)))), ((22, 18), (10, 1)))} 

\adfDgap
\noindent{\boldmath $ 22^{18} 13^{1} $}~
With the point set $Z_{409}$ partitioned into
 residue classes modulo $18$ for $\{0, 1, \dots, 395\}$, and
 $\{396, 397, \dots, 408\}$,
 the design is generated from

\adfLgap 
$(396, 154, 189, 263)$,
$(397, 72, 290, 142)$,
$(398, 130, 87, 332)$,
$(399, 233, 21, 115)$,\adfsplit
$(385, 302, 365, 231)$,
$(276, 197, 379, 211)$,
$(55, 152, 335, 86)$,
$(158, 164, 126, 137)$,\adfsplit
$(285, 380, 31, 227)$,
$(191, 225, 139, 364)$,
$(121, 354, 93, 178)$,
$(72, 177, 64, 31)$,\adfsplit
$(367, 59, 327, 170)$,
$(64, 128, 364, 61)$,
$(106, 17, 299, 221)$,
$(235, 205, 254, 49)$,\adfsplit
$(20, 297, 161, 190)$,
$(83, 143, 310, 255)$,
$(115, 324, 265, 326)$,
$(278, 362, 41, 34)$,\adfsplit
$(106, 89, 162, 64)$,
$(0, 1, 77, 315)$,
$(0, 4, 16, 181)$,
$(0, 50, 149, 272)$,\adfsplit
$(0, 9, 100, 295)$,
$(0, 22, 45, 316)$,
$(0, 44, 92, 232)$,
$(0, 37, 106, 267)$,\adfsplit
$(0, 15, 137, 190)$,
$(0, 10, 127, 189)$,
$(0, 5, 51, 138)$,
$(0, 24, 145, 265)$,\adfsplit
$(0, 13, 39, 143)$,
$(408, 0, 132, 264)$

\adfLgap \noindent by the mapping:
$x \mapsto x +  j \adfmod{396}$ for $x < 396$,
$x \mapsto (x + 4 j \adfmod{12}) + 396$ for $396 \le x < 408$,
$408 \mapsto 408$,
$0 \le j < 396$
 for the first 33 blocks,
$0 \le j < 132$
 for the last block.
\ADFvfyParStart{(409, ((33, 396, ((396, 1), (12, 4), (1, 1))), (1, 132, ((396, 1), (12, 4), (1, 1)))), ((22, 18), (13, 1)))} 

\adfDgap
\noindent{\boldmath $ 22^{18} 16^{1} $}~
With the point set $Z_{412}$ partitioned into
 residue classes modulo $18$ for $\{0, 1, \dots, 395\}$, and
 $\{396, 397, \dots, 411\}$,
 the design is generated from

\adfLgap 
$(396, 70, 195, 236)$,
$(396, 330, 127, 155)$,
$(397, 176, 75, 41)$,
$(397, 133, 166, 282)$,\adfsplit
$(398, 260, 69, 289)$,
$(398, 16, 108, 359)$,
$(399, 304, 117, 299)$,
$(399, 390, 73, 350)$,\adfsplit
$(400, 269, 195, 308)$,
$(400, 85, 76, 390)$,
$(244, 37, 178, 44)$,
$(281, 343, 214, 74)$,\adfsplit
$(164, 101, 86, 217)$,
$(375, 385, 118, 24)$,
$(29, 239, 21, 354)$,
$(224, 229, 64, 281)$,\adfsplit
$(158, 84, 295, 5)$,
$(180, 89, 210, 328)$,
$(137, 111, 303, 352)$,
$(206, 349, 225, 64)$,\adfsplit
$(147, 171, 270, 190)$,
$(298, 37, 320, 139)$,
$(330, 223, 375, 274)$,
$(199, 282, 367, 350)$,\adfsplit
$(192, 150, 20, 302)$,
$(14, 282, 218, 232)$,
$(88, 177, 257, 164)$,
$(234, 37, 285, 221)$,\adfsplit
$(250, 353, 22, 45)$,
$(111, 324, 28, 177)$,
$(390, 254, 37, 303)$,
$(13, 63, 348, 314)$,\adfsplit
$(140, 186, 290, 16)$,
$(214, 97, 126, 247)$,
$(2, 301, 382, 321)$,
$(35, 67, 32, 8)$,\adfsplit
$(254, 275, 100, 202)$,
$(379, 303, 5, 238)$,
$(343, 354, 189, 190)$,
$(16, 384, 190, 14)$,\adfsplit
$(308, 356, 237, 137)$,
$(100, 297, 125, 15)$,
$(145, 315, 218, 345)$,
$(271, 274, 154, 185)$,\adfsplit
$(369, 332, 392, 322)$,
$(131, 342, 130, 333)$,
$(248, 216, 355, 227)$,
$(134, 76, 96, 243)$,\adfsplit
$(0, 115, 173, 240)$,
$(0, 12, 98, 339)$,
$(0, 8, 146, 190)$,
$(0, 84, 303, 349)$,\adfsplit
$(0, 4, 227, 239)$,
$(0, 77, 186, 337)$,
$(0, 62, 195, 300)$,
$(0, 105, 208, 379)$,\adfsplit
$(0, 87, 229, 233)$,
$(0, 69, 75, 187)$,
$(0, 6, 259, 290)$,
$(0, 213, 369, 371)$,\adfsplit
$(0, 7, 321, 359)$,
$(0, 99, 159, 263)$,
$(0, 39, 301, 341)$,
$(1, 17, 85, 223)$,\adfsplit
$(0, 41, 55, 381)$,
$(1, 49, 141, 237)$,
$(0, 81, 123, 201)$,
$(411, 0, 132, 264)$,\adfsplit
$(411, 1, 133, 265)$

\adfLgap \noindent by the mapping:
$x \mapsto x + 2 j \adfmod{396}$ for $x < 396$,
$x \mapsto (x - 396 + 5 j \adfmod{15}) + 396$ for $396 \le x < 411$,
$411 \mapsto 411$,
$0 \le j < 198$
 for the first 67 blocks,
$0 \le j < 66$
 for the last two blocks.
\ADFvfyParStart{(412, ((67, 198, ((396, 2), (15, 5), (1, 1))), (2, 66, ((396, 2), (15, 5), (1, 1)))), ((22, 18), (16, 1)))} 

\adfDgap
\noindent{\boldmath $ 22^{18} 19^{1} $}~
With the point set $Z_{415}$ partitioned into
 residue classes modulo $18$ for $\{0, 1, \dots, 395\}$, and
 $\{396, 397, \dots, 414\}$,
 the design is generated from

\adfLgap 
$(396, 278, 6, 136)$,
$(397, 199, 333, 296)$,
$(398, 209, 106, 291)$,
$(399, 53, 111, 259)$,\adfsplit
$(400, 391, 63, 62)$,
$(401, 248, 105, 322)$,
$(27, 41, 148, 12)$,
$(143, 44, 135, 236)$,\adfsplit
$(4, 61, 249, 174)$,
$(240, 192, 303, 272)$,
$(349, 249, 84, 327)$,
$(293, 12, 249, 19)$,\adfsplit
$(69, 251, 60, 90)$,
$(393, 80, 78, 172)$,
$(20, 390, 47, 216)$,
$(257, 390, 343, 141)$,\adfsplit
$(221, 202, 267, 150)$,
$(10, 234, 74, 199)$,
$(114, 259, 364, 2)$,
$(327, 343, 330, 386)$,\adfsplit
$(217, 40, 340, 352)$,
$(0, 4, 141, 158)$,
$(0, 5, 176, 186)$,
$(0, 20, 60, 354)$,\adfsplit
$(0, 50, 120, 317)$,
$(0, 11, 66, 298)$,
$(0, 39, 127, 212)$,
$(0, 25, 76, 228)$,\adfsplit
$(0, 6, 110, 155)$,
$(0, 24, 119, 233)$,
$(0, 28, 69, 156)$,
$(0, 49, 138, 278)$,\adfsplit
$(0, 33, 106, 183)$,
$(0, 23, 61, 218)$,
$(414, 0, 132, 264)$

\adfLgap \noindent by the mapping:
$x \mapsto x +  j \adfmod{396}$ for $x < 396$,
$x \mapsto (x + 6 j \adfmod{18}) + 396$ for $396 \le x < 414$,
$414 \mapsto 414$,
$0 \le j < 396$
 for the first 34 blocks,
$0 \le j < 132$
 for the last block.
\ADFvfyParStart{(415, ((34, 396, ((396, 1), (18, 6), (1, 1))), (1, 132, ((396, 1), (18, 6), (1, 1)))), ((22, 18), (19, 1)))} 

\section{4-GDDs for the proof of Lemma \ref{lem:4-GDD 26^u m^1}}
\label{app:4-GDD 26^u m^1}
\adfhide{
$ 26^6 8^1 $,
$ 26^6 11^1 $,
$ 26^6 14^1 $,
$ 26^6 17^1 $,
$ 26^6 20^1 $,
$ 26^6 23^1 $,
$ 26^6 29^1 $,
$ 26^6 32^1 $,
$ 26^6 35^1 $,
$ 26^6 38^1 $,
$ 26^6 41^1 $,
$ 26^6 44^1 $,
$ 26^6 47^1 $,
$ 26^6 50^1 $,
$ 26^6 53^1 $,
$ 26^6 56^1 $,
$ 26^6 59^1 $,
$ 26^6 62^1 $,
$ 26^9 11^1 $,
$ 26^9 14^1 $,
$ 26^9 17^1 $,
$ 26^9 20^1 $,
$ 26^9 23^1 $,
$ 26^{12} 14^1 $,
$ 26^{12} 17^1 $,
$ 26^{12} 20^1 $,
$ 26^{12} 23^1 $,
$ 26^{15} 17^1 $,
$ 26^{15} 20^1 $,
$ 26^{15} 23^1 $,
$ 26^{18} 20^1 $ and
$ 26^{18} 23^1 $.
}

\adfDgap
\noindent{\boldmath $ 26^{6} 8^{1} $}~
With the point set $Z_{164}$ partitioned into
 residue classes modulo $6$ for $\{0, 1, \dots, 155\}$, and
 $\{156, 157, \dots, 163\}$,
 the design is generated from

\adfLgap 
$(156, 0, 1, 2)$,
$(157, 0, 52, 50)$,
$(158, 0, 103, 155)$,
$(159, 0, 106, 53)$,\adfsplit
$(160, 0, 154, 104)$,
$(161, 0, 4, 11)$,
$(162, 0, 7, 152)$,
$(163, 0, 145, 149)$,\adfsplit
$(41, 115, 27, 22)$,
$(38, 59, 67, 102)$,
$(99, 138, 130, 61)$,
$(91, 146, 155, 132)$,\adfsplit
$(63, 18, 89, 73)$,
$(0, 3, 16, 97)$,
$(0, 10, 43, 134)$,
$(0, 13, 70, 128)$,\adfsplit
$(0, 25, 40, 123)$,
$(0, 44, 105, 133)$,
$(0, 26, 100, 135)$,
$(0, 38, 77, 111)$,\adfsplit
$(0, 20, 47, 69)$,
$(0, 67, 107, 153)$,
$(0, 34, 85, 129)$,
$(0, 15, 80, 139)$,\adfsplit
$(0, 17, 37, 46)$,
$(0, 57, 62, 137)$,
$(0, 31, 68, 131)$

\adfLgap \noindent by the mapping:
$x \mapsto x + 3 j \adfmod{156}$ for $x < 156$,
$x \mapsto x$ for $x \ge 156$,
$0 \le j < 52$
 for the first eight blocks;
$x \mapsto x + 2 j \adfmod{156}$ for $x < 156$,
$x \mapsto x$ for $x \ge 156$,
$0 \le j < 78$
 for the last 19 blocks.
\ADFvfyParStart{(164, ((8, 52, ((156, 3), (8, 8))), (19, 78, ((156, 2), (8, 8)))), ((26, 6), (8, 1)))} 

\adfDgap
\noindent{\boldmath $ 26^{6} 11^{1} $}~
With the point set $Z_{167}$ partitioned into
 residue classes modulo $6$ for $\{0, 1, \dots, 155\}$, and
 $\{156, 157, \dots, 166\}$,
 the design is generated from

\adfLgap 
$(156, 0, 1, 2)$,
$(157, 0, 52, 50)$,
$(158, 0, 103, 155)$,
$(159, 0, 106, 53)$,\adfsplit
$(160, 0, 154, 104)$,
$(161, 0, 4, 11)$,
$(162, 0, 7, 152)$,
$(163, 0, 145, 149)$,\adfsplit
$(164, 0, 10, 23)$,
$(165, 0, 13, 146)$,
$(166, 0, 133, 143)$,
$(0, 3, 8, 95)$,\adfsplit
$(0, 9, 25, 58)$,
$(0, 14, 34, 81)$,
$(0, 17, 44, 63)$,
$(0, 26, 65, 97)$,\adfsplit
$(0, 31, 74, 111)$,
$(0, 22, 51, 121)$,
$(0, 15, 55, 83)$,
$(0, 21, 62, 100)$

\adfLgap \noindent by the mapping:
$x \mapsto x + 3 j \adfmod{156}$ for $x < 156$,
$x \mapsto x$ for $x \ge 156$,
$0 \le j < 52$
 for the first 11 blocks;
$x \mapsto x +  j \adfmod{156}$ for $x < 156$,
$x \mapsto x$ for $x \ge 156$,
$0 \le j < 156$
 for the last nine blocks.
\ADFvfyParStart{(167, ((11, 52, ((156, 3), (11, 11))), (9, 156, ((156, 1), (11, 11)))), ((26, 6), (11, 1)))} 

\adfDgap
\noindent{\boldmath $ 26^{6} 14^{1} $}~
With the point set $Z_{170}$ partitioned into
 residue classes modulo $6$ for $\{0, 1, \dots, 155\}$, and
 $\{156, 157, \dots, 169\}$,
 the design is generated from

\adfLgap 
$(156, 0, 1, 2)$,
$(157, 0, 52, 50)$,
$(158, 0, 103, 155)$,
$(159, 0, 106, 53)$,\adfsplit
$(160, 0, 154, 104)$,
$(161, 0, 4, 11)$,
$(162, 0, 7, 152)$,
$(163, 0, 145, 149)$,\adfsplit
$(164, 0, 10, 23)$,
$(165, 0, 13, 146)$,
$(166, 0, 133, 143)$,
$(167, 0, 16, 35)$,\adfsplit
$(168, 0, 19, 140)$,
$(169, 0, 121, 137)$,
$(0, 3, 8, 83)$,
$(0, 5, 14, 123)$,\adfsplit
$(0, 9, 29, 73)$,
$(0, 15, 55, 130)$,
$(0, 17, 87, 136)$,
$(0, 21, 38, 112)$,\adfsplit
$(0, 27, 61, 124)$,
$(0, 57, 71, 98)$,
$(0, 31, 99, 131)$,
$(0, 28, 97, 105)$,\adfsplit
$(0, 25, 51, 80)$,
$(0, 47, 92, 141)$,
$(0, 43, 89, 117)$,
$(0, 33, 91, 94)$,\adfsplit
$(0, 45, 67, 88)$,
$(0, 22, 85, 122)$,
$(0, 39, 70, 116)$

\adfLgap \noindent by the mapping:
$x \mapsto x + 3 j \adfmod{156}$ for $x < 156$,
$x \mapsto x$ for $x \ge 156$,
$0 \le j < 52$
 for the first 14 blocks;
$x \mapsto x + 2 j \adfmod{156}$ for $x < 156$,
$x \mapsto x$ for $x \ge 156$,
$0 \le j < 78$
 for the last 17 blocks.
\ADFvfyParStart{(170, ((14, 52, ((156, 3), (14, 14))), (17, 78, ((156, 2), (14, 14)))), ((26, 6), (14, 1)))} 

\adfDgap
\noindent{\boldmath $ 26^{6} 17^{1} $}~
With the point set $Z_{173}$ partitioned into
 residue classes modulo $6$ for $\{0, 1, \dots, 155\}$, and
 $\{156, 157, \dots, 172\}$,
 the design is generated from

\adfLgap 
$(156, 0, 1, 2)$,
$(157, 0, 52, 50)$,
$(158, 0, 103, 155)$,
$(159, 0, 106, 53)$,\adfsplit
$(160, 0, 154, 104)$,
$(161, 0, 4, 11)$,
$(162, 0, 7, 152)$,
$(163, 0, 145, 149)$,\adfsplit
$(164, 0, 10, 23)$,
$(165, 0, 13, 146)$,
$(166, 0, 133, 143)$,
$(167, 0, 16, 35)$,\adfsplit
$(168, 0, 19, 140)$,
$(169, 0, 121, 137)$,
$(170, 0, 22, 5)$,
$(171, 0, 139, 134)$,\adfsplit
$(172, 0, 151, 17)$,
$(0, 3, 29, 110)$,
$(0, 8, 63, 88)$,
$(0, 9, 40, 109)$,\adfsplit
$(0, 21, 62, 95)$,
$(0, 28, 67, 111)$,
$(0, 27, 65, 97)$,
$(0, 15, 58, 92)$,\adfsplit
$(0, 14, 51, 71)$

\adfLgap \noindent by the mapping:
$x \mapsto x + 3 j \adfmod{156}$ for $x < 156$,
$x \mapsto x$ for $x \ge 156$,
$0 \le j < 52$
 for the first 17 blocks;
$x \mapsto x +  j \adfmod{156}$ for $x < 156$,
$x \mapsto x$ for $x \ge 156$,
$0 \le j < 156$
 for the last eight blocks.
\ADFvfyParStart{(173, ((17, 52, ((156, 3), (17, 17))), (8, 156, ((156, 1), (17, 17)))), ((26, 6), (17, 1)))} 

\adfDgap
\noindent{\boldmath $ 26^{6} 20^{1} $}~
With the point set $Z_{176}$ partitioned into
 residue classes modulo $6$ for $\{0, 1, \dots, 155\}$, and
 $\{156, 157, \dots, 175\}$,
 the design is generated from

\adfLgap 
$(168, 0, 1, 2)$,
$(169, 0, 52, 50)$,
$(170, 0, 103, 155)$,
$(171, 0, 106, 53)$,\adfsplit
$(172, 0, 154, 104)$,
$(173, 0, 4, 11)$,
$(174, 0, 7, 152)$,
$(175, 0, 145, 149)$,\adfsplit
$(156, 76, 21, 101)$,
$(156, 6, 140, 82)$,
$(156, 75, 7, 35)$,
$(156, 36, 13, 26)$,\adfsplit
$(157, 69, 136, 110)$,
$(157, 91, 17, 48)$,
$(157, 11, 118, 54)$,
$(157, 25, 44, 135)$,\adfsplit
$(58, 25, 74, 137)$,
$(147, 150, 133, 62)$,
$(0, 3, 8, 23)$,
$(0, 14, 33, 95)$,\adfsplit
$(0, 27, 61, 100)$,
$(0, 17, 38, 75)$,
$(0, 41, 57, 86)$,
$(0, 46, 97, 105)$,\adfsplit
$(0, 9, 44, 65)$,
$(0, 13, 28, 122)$,
$(0, 29, 116, 147)$,
$(0, 20, 93, 119)$,\adfsplit
$(0, 39, 109, 131)$,
$(0, 32, 67, 77)$,
$(0, 5, 74, 129)$

\adfLgap \noindent by the mapping:
$x \mapsto x + 3 j \adfmod{156}$ for $x < 156$,
$x \mapsto (x + 2 j \adfmod{12}) + 156$ for $156 \le x < 168$,
$x \mapsto x$ for $x \ge 168$,
$0 \le j < 52$
 for the first eight blocks;
$x \mapsto x + 2 j \adfmod{156}$ for $x < 156$,
$x \mapsto (x + 2 j \adfmod{12}) + 156$ for $156 \le x < 168$,
$x \mapsto x$ for $x \ge 168$,
$0 \le j < 78$
 for the last 23 blocks.
\ADFvfyParStart{(176, ((8, 52, ((156, 3), (12, 2), (8, 8))), (23, 78, ((156, 2), (12, 2), (8, 8)))), ((26, 6), (20, 1)))} 

\adfDgap
\noindent{\boldmath $ 26^{6} 23^{1} $}~
With the point set $Z_{179}$ partitioned into
 residue classes modulo $6$ for $\{0, 1, \dots, 155\}$, and
 $\{156, 157, \dots, 178\}$,
 the design is generated from

\adfLgap 
$(174, 0, 1, 2)$,
$(175, 0, 52, 50)$,
$(176, 0, 103, 155)$,
$(177, 0, 106, 53)$,\adfsplit
$(178, 0, 154, 104)$,
$(156, 83, 87, 106)$,
$(156, 140, 145, 18)$,
$(157, 14, 21, 114)$,\adfsplit
$(157, 151, 112, 53)$,
$(158, 28, 93, 137)$,
$(0, 8, 25, 158)$,
$(0, 3, 13, 82)$,\adfsplit
$(0, 15, 43, 83)$,
$(0, 9, 35, 80)$,
$(0, 11, 38, 105)$,
$(0, 16, 37, 57)$,\adfsplit
$(0, 22, 55, 86)$,
$(0, 14, 46, 95)$

\adfLgap \noindent by the mapping:
$x \mapsto x + 3 j \adfmod{156}$ for $x < 156$,
$x \mapsto (x - 156 + 3 j \adfmod{18}) + 156$ for $156 \le x < 174$,
$x \mapsto x$ for $x \ge 174$,
$0 \le j < 52$
 for the first five blocks;
$x \mapsto x +  j \adfmod{156}$ for $x < 156$,
$x \mapsto (x - 156 + 3 j \adfmod{18}) + 156$ for $156 \le x < 174$,
$x \mapsto x$ for $x \ge 174$,
$0 \le j < 156$
 for the last 13 blocks.
\ADFvfyParStart{(179, ((5, 52, ((156, 3), (18, 3), (5, 5))), (13, 156, ((156, 1), (18, 3), (5, 5)))), ((26, 6), (23, 1)))} 

\adfDgap
\noindent{\boldmath $ 26^{6} 29^{1} $}~
With the point set $Z_{185}$ partitioned into
 residue classes modulo $6$ for $\{0, 1, \dots, 155\}$, and
 $\{156, 157, \dots, 184\}$,
 the design is generated from

\adfLgap 
$(180, 0, 1, 2)$,
$(181, 0, 52, 50)$,
$(182, 0, 103, 155)$,
$(183, 0, 106, 53)$,\adfsplit
$(184, 0, 154, 104)$,
$(156, 154, 7, 12)$,
$(156, 81, 44, 5)$,
$(157, 90, 97, 124)$,\adfsplit
$(157, 104, 75, 137)$,
$(158, 55, 122, 11)$,
$(158, 106, 36, 123)$,
$(0, 3, 64, 68)$,\adfsplit
$(0, 8, 19, 143)$,
$(0, 22, 71, 159)$,
$(0, 10, 38, 93)$,
$(0, 23, 58, 105)$,\adfsplit
$(0, 20, 46, 77)$,
$(0, 16, 59, 171)$,
$(0, 15, 40, 115)$

\adfLgap \noindent by the mapping:
$x \mapsto x + 3 j \adfmod{156}$ for $x < 156$,
$x \mapsto (x - 156 + 4 j \adfmod{24}) + 156$ for $156 \le x < 180$,
$x \mapsto x$ for $x \ge 180$,
$0 \le j < 52$
 for the first five blocks;
$x \mapsto x +  j \adfmod{156}$ for $x < 156$,
$x \mapsto (x - 156 + 4 j \adfmod{24}) + 156$ for $156 \le x < 180$,
$x \mapsto x$ for $x \ge 180$,
$0 \le j < 156$
 for the last 14 blocks.
\ADFvfyParStart{(185, ((5, 52, ((156, 3), (24, 4), (5, 5))), (14, 156, ((156, 1), (24, 4), (5, 5)))), ((26, 6), (29, 1)))} 

\adfDgap
\noindent{\boldmath $ 26^{6} 32^{1} $}~
With the point set $Z_{188}$ partitioned into
 residue classes modulo $6$ for $\{0, 1, \dots, 155\}$, and
 $\{156, 157, \dots, 187\}$,
 the design is generated from

\adfLgap 
$(180, 0, 1, 2)$,
$(181, 0, 52, 50)$,
$(182, 0, 103, 155)$,
$(183, 0, 106, 53)$,\adfsplit
$(184, 0, 154, 104)$,
$(185, 0, 4, 11)$,
$(186, 0, 7, 152)$,
$(187, 0, 145, 149)$,\adfsplit
$(156, 85, 58, 101)$,
$(156, 83, 147, 96)$,
$(156, 92, 115, 153)$,
$(156, 42, 136, 122)$,\adfsplit
$(157, 72, 37, 34)$,
$(157, 124, 91, 26)$,
$(157, 113, 20, 87)$,
$(157, 90, 21, 119)$,\adfsplit
$(158, 83, 151, 38)$,
$(158, 68, 93, 142)$,
$(158, 100, 75, 84)$,
$(158, 17, 126, 61)$,\adfsplit
$(159, 115, 29, 60)$,
$(0, 5, 130, 151)$,
$(0, 8, 28, 41)$,
$(0, 10, 85, 122)$,\adfsplit
$(0, 22, 99, 133)$,
$(0, 19, 101, 129)$,
$(0, 46, 109, 141)$,
$(0, 69, 83, 86)$,\adfsplit
$(0, 9, 71, 92)$,
$(0, 15, 35, 115)$,
$(0, 17, 56, 167)$,
$(0, 39, 116, 171)$,\adfsplit
$(0, 32, 81, 89)$,
$(0, 97, 137, 163)$,
$(0, 37, 59, 88)$

\adfLgap \noindent by the mapping:
$x \mapsto x + 3 j \adfmod{156}$ for $x < 156$,
$x \mapsto (x - 156 + 4 j \adfmod{24}) + 156$ for $156 \le x < 180$,
$x \mapsto x$ for $x \ge 180$,
$0 \le j < 52$
 for the first eight blocks;
$x \mapsto x + 2 j \adfmod{156}$ for $x < 156$,
$x \mapsto (x - 156 + 4 j \adfmod{24}) + 156$ for $156 \le x < 180$,
$x \mapsto x$ for $x \ge 180$,
$0 \le j < 78$
 for the last 27 blocks.
\ADFvfyParStart{(188, ((8, 52, ((156, 3), (24, 4), (8, 8))), (27, 78, ((156, 2), (24, 4), (8, 8)))), ((26, 6), (32, 1)))} 

\adfDgap
\noindent{\boldmath $ 26^{6} 35^{1} $}~
With the point set $Z_{191}$ partitioned into
 residue classes modulo $6$ for $\{0, 1, \dots, 155\}$, and
 $\{156, 157, \dots, 190\}$,
 the design is generated from

\adfLgap 
$(186, 0, 1, 2)$,
$(187, 0, 52, 50)$,
$(188, 0, 103, 155)$,
$(189, 0, 106, 53)$,\adfsplit
$(190, 0, 154, 104)$,
$(156, 19, 132, 128)$,
$(156, 41, 123, 52)$,
$(157, 26, 99, 29)$,\adfsplit
$(157, 64, 115, 108)$,
$(158, 105, 85, 68)$,
$(158, 34, 132, 11)$,
$(159, 46, 147, 60)$,\adfsplit
$(0, 5, 15, 179)$,
$(0, 8, 34, 97)$,
$(0, 16, 41, 80)$,
$(0, 22, 68, 99)$,\adfsplit
$(0, 13, 45, 94)$,
$(0, 28, 61, 165)$,
$(0, 19, 135, 185)$,
$(0, 9, 38, 65)$

\adfLgap \noindent by the mapping:
$x \mapsto x + 3 j \adfmod{156}$ for $x < 156$,
$x \mapsto (x - 156 + 5 j \adfmod{30}) + 156$ for $156 \le x < 186$,
$x \mapsto x$ for $x \ge 186$,
$0 \le j < 52$
 for the first five blocks;
$x \mapsto x +  j \adfmod{156}$ for $x < 156$,
$x \mapsto (x - 156 + 5 j \adfmod{30}) + 156$ for $156 \le x < 186$,
$x \mapsto x$ for $x \ge 186$,
$0 \le j < 156$
 for the last 15 blocks.
\ADFvfyParStart{(191, ((5, 52, ((156, 3), (30, 5), (5, 5))), (15, 156, ((156, 1), (30, 5), (5, 5)))), ((26, 6), (35, 1)))} 

\adfDgap
\noindent{\boldmath $ 26^{6} 38^{1} $}~
With the point set $Z_{194}$ partitioned into
 residue classes modulo $6$ for $\{0, 1, \dots, 155\}$, and
 $\{156, 157, \dots, 193\}$,
 the design is generated from

\adfLgap 
$(186, 0, 1, 2)$,
$(187, 0, 52, 50)$,
$(188, 0, 103, 155)$,
$(189, 0, 106, 53)$,\adfsplit
$(190, 0, 154, 104)$,
$(191, 0, 4, 11)$,
$(192, 0, 7, 152)$,
$(193, 0, 145, 149)$,\adfsplit
$(156, 9, 50, 150)$,
$(156, 4, 17, 133)$,
$(156, 131, 75, 96)$,
$(156, 116, 55, 142)$,\adfsplit
$(157, 155, 82, 123)$,
$(157, 124, 14, 48)$,
$(157, 140, 67, 125)$,
$(157, 42, 117, 133)$,\adfsplit
$(158, 98, 131, 76)$,
$(158, 118, 49, 0)$,
$(158, 147, 139, 116)$,
$(158, 149, 18, 81)$,\adfsplit
$(159, 142, 71, 85)$,
$(159, 29, 38, 63)$,
$(159, 103, 52, 42)$,
$(0, 8, 45, 159)$,\adfsplit
$(0, 3, 92, 97)$,
$(0, 17, 20, 94)$,
$(0, 19, 116, 137)$,
$(0, 9, 88, 170)$,\adfsplit
$(0, 32, 143, 185)$,
$(0, 101, 123, 180)$,
$(0, 29, 39, 121)$,
$(0, 14, 57, 127)$,\adfsplit
$(0, 27, 107, 112)$,
$(0, 28, 86, 133)$,
$(0, 93, 119, 139)$,
$(0, 89, 117, 175)$,\adfsplit
$(0, 16, 81, 125)$

\adfLgap \noindent by the mapping:
$x \mapsto x + 3 j \adfmod{156}$ for $x < 156$,
$x \mapsto (x - 156 + 5 j \adfmod{30}) + 156$ for $156 \le x < 186$,
$x \mapsto x$ for $x \ge 186$,
$0 \le j < 52$
 for the first eight blocks;
$x \mapsto x + 2 j \adfmod{156}$ for $x < 156$,
$x \mapsto (x - 156 + 5 j \adfmod{30}) + 156$ for $156 \le x < 186$,
$x \mapsto x$ for $x \ge 186$,
$0 \le j < 78$
 for the last 29 blocks.
\ADFvfyParStart{(194, ((8, 52, ((156, 3), (30, 5), (8, 8))), (29, 78, ((156, 2), (30, 5), (8, 8)))), ((26, 6), (38, 1)))} 

\adfDgap
\noindent{\boldmath $ 26^{6} 41^{1} $}~
With the point set $Z_{197}$ partitioned into
 residue classes modulo $6$ for $\{0, 1, \dots, 155\}$, and
 $\{156, 157, \dots, 196\}$,
 the design is generated from

\adfLgap 
$(192, 0, 1, 2)$,
$(193, 0, 52, 50)$,
$(194, 0, 103, 155)$,
$(195, 0, 106, 53)$,\adfsplit
$(196, 0, 154, 104)$,
$(156, 109, 100, 65)$,
$(156, 51, 144, 116)$,
$(157, 61, 45, 0)$,\adfsplit
$(157, 86, 124, 65)$,
$(158, 139, 114, 82)$,
$(158, 74, 45, 125)$,
$(159, 27, 19, 108)$,\adfsplit
$(159, 11, 44, 112)$,
$(0, 3, 14, 34)$,
$(0, 10, 23, 92)$,
$(0, 5, 46, 160)$,\adfsplit
$(0, 4, 47, 178)$,
$(0, 7, 26, 173)$,
$(0, 15, 71, 98)$,
$(0, 37, 86, 191)$,\adfsplit
$(0, 17, 39, 79)$

\adfLgap \noindent by the mapping:
$x \mapsto x + 3 j \adfmod{156}$ for $x < 156$,
$x \mapsto (x - 156 + 6 j \adfmod{36}) + 156$ for $156 \le x < 192$,
$x \mapsto x$ for $x \ge 192$,
$0 \le j < 52$
 for the first five blocks;
$x \mapsto x +  j \adfmod{156}$ for $x < 156$,
$x \mapsto (x - 156 + 6 j \adfmod{36}) + 156$ for $156 \le x < 192$,
$x \mapsto x$ for $x \ge 192$,
$0 \le j < 156$
 for the last 16 blocks.
\ADFvfyParStart{(197, ((5, 52, ((156, 3), (36, 6), (5, 5))), (16, 156, ((156, 1), (36, 6), (5, 5)))), ((26, 6), (41, 1)))} 

\adfDgap
\noindent{\boldmath $ 26^{6} 44^{1} $}~
With the point set $Z_{200}$ partitioned into
 residue classes modulo $6$ for $\{0, 1, \dots, 155\}$, and
 $\{156, 157, \dots, 199\}$,
 the design is generated from

\adfLgap 
$(192, 0, 1, 2)$,
$(193, 0, 52, 50)$,
$(194, 0, 103, 155)$,
$(195, 0, 106, 53)$,\adfsplit
$(196, 0, 154, 104)$,
$(197, 0, 4, 11)$,
$(198, 0, 7, 152)$,
$(199, 0, 145, 149)$,\adfsplit
$(156, 142, 7, 96)$,
$(156, 73, 112, 153)$,
$(156, 11, 8, 39)$,
$(156, 138, 2, 89)$,\adfsplit
$(157, 83, 104, 141)$,
$(157, 42, 51, 64)$,
$(157, 121, 48, 146)$,
$(157, 125, 151, 154)$,\adfsplit
$(158, 143, 27, 43)$,
$(158, 8, 101, 88)$,
$(158, 50, 132, 109)$,
$(158, 126, 93, 70)$,\adfsplit
$(159, 18, 95, 26)$,
$(159, 65, 16, 21)$,
$(159, 44, 132, 91)$,
$(159, 109, 99, 154)$,\adfsplit
$(160, 44, 143, 1)$,
$(0, 10, 27, 95)$,
$(0, 14, 29, 151)$,
$(0, 34, 139, 161)$,\adfsplit
$(0, 62, 141, 179)$,
$(0, 35, 43, 116)$,
$(0, 25, 45, 119)$,
$(0, 57, 121, 185)$,\adfsplit
$(0, 33, 55, 191)$,
$(0, 26, 91, 118)$,
$(0, 32, 71, 166)$,
$(0, 16, 97, 128)$,\adfsplit
$(0, 63, 109, 184)$,
$(0, 19, 51, 89)$,
$(0, 61, 70, 190)$

\adfLgap \noindent by the mapping:
$x \mapsto x + 3 j \adfmod{156}$ for $x < 156$,
$x \mapsto (x - 156 + 6 j \adfmod{36}) + 156$ for $156 \le x < 192$,
$x \mapsto x$ for $x \ge 192$,
$0 \le j < 52$
 for the first eight blocks;
$x \mapsto x + 2 j \adfmod{156}$ for $x < 156$,
$x \mapsto (x - 156 + 6 j \adfmod{36}) + 156$ for $156 \le x < 192$,
$x \mapsto x$ for $x \ge 192$,
$0 \le j < 78$
 for the last 31 blocks.
\ADFvfyParStart{(200, ((8, 52, ((156, 3), (36, 6), (8, 8))), (31, 78, ((156, 2), (36, 6), (8, 8)))), ((26, 6), (44, 1)))} 

\adfDgap
\noindent{\boldmath $ 26^{6} 47^{1} $}~
With the point set $Z_{203}$ partitioned into
 residue classes modulo $6$ for $\{0, 1, \dots, 155\}$, and
 $\{156, 157, \dots, 202\}$,
 the design is generated from

\adfLgap 
$(198, 0, 1, 2)$,
$(199, 0, 52, 50)$,
$(200, 0, 103, 155)$,
$(201, 0, 106, 53)$,\adfsplit
$(202, 0, 154, 104)$,
$(156, 13, 80, 5)$,
$(156, 144, 75, 154)$,
$(157, 44, 154, 5)$,\adfsplit
$(157, 72, 117, 97)$,
$(158, 133, 96, 28)$,
$(158, 87, 101, 74)$,
$(159, 70, 101, 79)$,\adfsplit
$(159, 126, 110, 105)$,
$(160, 13, 53, 56)$,
$(0, 4, 63, 160)$,
$(0, 17, 64, 161)$,\adfsplit
$(0, 19, 74, 196)$,
$(0, 34, 83, 162)$,
$(0, 11, 44, 76)$,
$(0, 28, 57, 98)$,\adfsplit
$(0, 15, 38, 100)$,
$(0, 26, 61, 197)$

\adfLgap \noindent by the mapping:
$x \mapsto x + 3 j \adfmod{156}$ for $x < 156$,
$x \mapsto (x - 156 + 7 j \adfmod{42}) + 156$ for $156 \le x < 198$,
$x \mapsto x$ for $x \ge 198$,
$0 \le j < 52$
 for the first five blocks;
$x \mapsto x +  j \adfmod{156}$ for $x < 156$,
$x \mapsto (x - 156 + 7 j \adfmod{42}) + 156$ for $156 \le x < 198$,
$x \mapsto x$ for $x \ge 198$,
$0 \le j < 156$
 for the last 17 blocks.
\ADFvfyParStart{(203, ((5, 52, ((156, 3), (42, 7), (5, 5))), (17, 156, ((156, 1), (42, 7), (5, 5)))), ((26, 6), (47, 1)))} 

\adfDgap
\noindent{\boldmath $ 26^{6} 50^{1} $}~
With the point set $Z_{206}$ partitioned into
 residue classes modulo $6$ for $\{0, 1, \dots, 155\}$, and
 $\{156, 157, \dots, 205\}$,
 the design is generated from

\adfLgap 
$(198, 0, 1, 2)$,
$(199, 0, 52, 50)$,
$(200, 0, 103, 155)$,
$(201, 0, 106, 53)$,\adfsplit
$(202, 0, 154, 104)$,
$(203, 0, 4, 11)$,
$(204, 0, 7, 152)$,
$(205, 0, 145, 149)$,\adfsplit
$(156, 86, 51, 71)$,
$(156, 37, 6, 22)$,
$(156, 12, 129, 101)$,
$(156, 80, 136, 139)$,\adfsplit
$(157, 72, 135, 38)$,
$(157, 92, 109, 82)$,
$(157, 137, 79, 4)$,
$(157, 90, 57, 95)$,\adfsplit
$(158, 76, 2, 49)$,
$(158, 72, 140, 94)$,
$(158, 143, 103, 129)$,
$(158, 41, 102, 135)$,\adfsplit
$(159, 46, 75, 145)$,
$(159, 9, 119, 122)$,
$(159, 124, 132, 103)$,
$(159, 138, 89, 8)$,\adfsplit
$(160, 34, 127, 111)$,
$(160, 144, 68, 23)$,
$(160, 26, 138, 1)$,
$(0, 13, 101, 188)$,\adfsplit
$(0, 9, 28, 98)$,
$(0, 14, 39, 119)$,
$(0, 37, 45, 161)$,
$(0, 83, 139, 162)$,\adfsplit
$(0, 57, 79, 168)$,
$(0, 41, 73, 176)$,
$(0, 21, 55, 65)$,
$(0, 32, 147, 175)$,\adfsplit
$(0, 40, 109, 196)$,
$(0, 20, 91, 169)$,
$(0, 51, 64, 125)$,
$(0, 38, 87, 151)$,\adfsplit
$(0, 23, 94, 190)$

\adfLgap \noindent by the mapping:
$x \mapsto x + 3 j \adfmod{156}$ for $x < 156$,
$x \mapsto (x - 156 + 7 j \adfmod{42}) + 156$ for $156 \le x < 198$,
$x \mapsto x$ for $x \ge 198$,
$0 \le j < 52$
 for the first eight blocks;
$x \mapsto x + 2 j \adfmod{156}$ for $x < 156$,
$x \mapsto (x - 156 + 7 j \adfmod{42}) + 156$ for $156 \le x < 198$,
$x \mapsto x$ for $x \ge 198$,
$0 \le j < 78$
 for the last 33 blocks.
\ADFvfyParStart{(206, ((8, 52, ((156, 3), (42, 7), (8, 8))), (33, 78, ((156, 2), (42, 7), (8, 8)))), ((26, 6), (50, 1)))} 

\adfDgap
\noindent{\boldmath $ 26^{6} 53^{1} $}~
With the point set $Z_{209}$ partitioned into
 residue classes modulo $6$ for $\{0, 1, \dots, 155\}$, and
 $\{156, 157, \dots, 208\}$,
 the design is generated from

\adfLgap 
$(195, 0, 1, 2)$,
$(196, 0, 52, 50)$,
$(197, 0, 103, 155)$,
$(198, 0, 106, 53)$,\adfsplit
$(199, 0, 154, 104)$,
$(200, 0, 4, 11)$,
$(201, 0, 7, 152)$,
$(202, 0, 145, 149)$,\adfsplit
$(203, 0, 10, 23)$,
$(204, 0, 13, 146)$,
$(205, 0, 133, 143)$,
$(206, 0, 16, 35)$,\adfsplit
$(207, 0, 19, 140)$,
$(208, 0, 121, 137)$,
$(156, 127, 118, 18)$,
$(156, 9, 140, 83)$,\adfsplit
$(156, 28, 69, 36)$,
$(156, 104, 46, 19)$,
$(156, 32, 112, 144)$,
$(156, 91, 134, 53)$,\adfsplit
$(0, 3, 29, 68)$,
$(0, 20, 93, 156)$,
$(0, 22, 67, 192)$,
$(0, 5, 64, 179)$,\adfsplit
$(0, 14, 101, 171)$,
$(0, 17, 51, 79)$,
$(0, 15, 46, 173)$,
$(0, 37, 107, 189)$,\adfsplit
$(0, 21, 116, 183)$

\adfLgap \noindent by the mapping:
$x \mapsto x + 3 j \adfmod{156}$ for $x < 156$,
$x \mapsto (x +  j \adfmod{39}) + 156$ for $156 \le x < 195$,
$x \mapsto x$ for $x \ge 195$,
$0 \le j < 52$
 for the first 14 blocks;
$x \mapsto x +  j \adfmod{156}$ for $x < 156$,
$x \mapsto (x +  j \adfmod{39}) + 156$ for $156 \le x < 195$,
$x \mapsto x$ for $x \ge 195$,
$0 \le j < 156$
 for the last 15 blocks.
\ADFvfyParStart{(209, ((14, 52, ((156, 3), (39, 1), (14, 14))), (15, 156, ((156, 1), (39, 1), (14, 14)))), ((26, 6), (53, 1)))} 

\adfDgap
\noindent{\boldmath $ 26^{6} 56^{1} $}~
With the point set $Z_{212}$ partitioned into
 residue classes modulo $6$ for $\{0, 1, \dots, 155\}$, and
 $\{156, 157, \dots, 211\}$,
 the design is generated from

\adfLgap 
$(195, 0, 1, 2)$,
$(196, 0, 52, 50)$,
$(197, 0, 103, 155)$,
$(198, 0, 106, 53)$,\adfsplit
$(199, 0, 154, 104)$,
$(200, 0, 4, 11)$,
$(201, 0, 7, 152)$,
$(202, 0, 145, 149)$,\adfsplit
$(203, 0, 10, 23)$,
$(204, 0, 13, 146)$,
$(205, 0, 133, 143)$,
$(206, 0, 16, 35)$,\adfsplit
$(207, 0, 19, 140)$,
$(208, 0, 121, 137)$,
$(209, 0, 22, 5)$,
$(210, 0, 139, 134)$,\adfsplit
$(211, 0, 151, 17)$,
$(156, 94, 54, 140)$,
$(156, 113, 80, 37)$,
$(156, 21, 106, 151)$,\adfsplit
$(156, 59, 114, 34)$,
$(156, 15, 12, 49)$,
$(156, 89, 58, 43)$,
$(156, 101, 92, 133)$,\adfsplit
$(156, 63, 72, 1)$,
$(156, 22, 149, 85)$,
$(156, 124, 45, 152)$,
$(156, 25, 64, 120)$,\adfsplit
$(156, 135, 84, 19)$,
$(156, 29, 147, 66)$,
$(0, 14, 29, 111)$,
$(0, 20, 89, 185)$,\adfsplit
$(0, 43, 94, 191)$,
$(0, 39, 64, 169)$,
$(0, 38, 82, 173)$,
$(0, 109, 129, 190)$,\adfsplit
$(0, 55, 99, 182)$,
$(0, 21, 79, 179)$,
$(0, 26, 124, 180)$,
$(0, 34, 93, 107)$,\adfsplit
$(0, 57, 125, 186)$,
$(0, 65, 68, 135)$,
$(0, 27, 83, 156)$,
$(0, 47, 75, 160)$,\adfsplit
$(0, 8, 95, 165)$,
$(0, 115, 123, 170)$

\adfLgap \noindent by the mapping:
$x \mapsto x + 3 j \adfmod{156}$ for $x < 156$,
$x \mapsto (x +  j \adfmod{39}) + 156$ for $156 \le x < 195$,
$x \mapsto x$ for $x \ge 195$,
$0 \le j < 52$
 for the first 17 blocks;
$x \mapsto x + 2 j \adfmod{156}$ for $x < 156$,
$x \mapsto (x +  j \adfmod{39}) + 156$ for $156 \le x < 195$,
$x \mapsto x$ for $x \ge 195$,
$0 \le j < 78$
 for the last 29 blocks.
\ADFvfyParStart{(212, ((17, 52, ((156, 3), (39, 1), (17, 17))), (29, 78, ((156, 2), (39, 1), (17, 17)))), ((26, 6), (56, 1)))} 

\adfDgap
\noindent{\boldmath $ 26^{6} 59^{1} $}~
With the point set $Z_{215}$ partitioned into
 residue classes modulo $6$ for $\{0, 1, \dots, 155\}$, and
 $\{156, 157, \dots, 214\}$,
 the design is generated from

\adfLgap 
$(195, 0, 1, 2)$,
$(196, 0, 52, 50)$,
$(197, 0, 103, 155)$,
$(198, 0, 106, 53)$,\adfsplit
$(199, 0, 154, 104)$,
$(200, 0, 4, 11)$,
$(201, 0, 7, 152)$,
$(202, 0, 145, 149)$,\adfsplit
$(203, 0, 10, 23)$,
$(204, 0, 13, 146)$,
$(205, 0, 133, 143)$,
$(206, 0, 16, 35)$,\adfsplit
$(207, 0, 19, 140)$,
$(208, 0, 121, 137)$,
$(209, 0, 22, 5)$,
$(210, 0, 139, 134)$,\adfsplit
$(211, 0, 151, 17)$,
$(212, 0, 25, 56)$,
$(213, 0, 31, 131)$,
$(214, 0, 100, 125)$,\adfsplit
$(156, 66, 63, 106)$,
$(156, 126, 46, 80)$,
$(156, 32, 114, 17)$,
$(156, 93, 42, 10)$,\adfsplit
$(0, 8, 28, 117)$,
$(0, 26, 111, 156)$,
$(0, 37, 112, 194)$,
$(0, 33, 98, 170)$,\adfsplit
$(0, 21, 115, 179)$,
$(0, 27, 95, 165)$,
$(0, 29, 86, 172)$,
$(0, 14, 93, 166)$,\adfsplit
$(0, 38, 107, 173)$,
$(0, 9, 101, 184)$

\adfLgap \noindent by the mapping:
$x \mapsto x + 3 j \adfmod{156}$ for $x < 156$,
$x \mapsto (x +  j \adfmod{39}) + 156$ for $156 \le x < 195$,
$x \mapsto x$ for $x \ge 195$,
$0 \le j < 52$
 for the first 20 blocks;
$x \mapsto x +  j \adfmod{156}$ for $x < 156$,
$x \mapsto (x +  j \adfmod{39}) + 156$ for $156 \le x < 195$,
$x \mapsto x$ for $x \ge 195$,
$0 \le j < 156$
 for the last 14 blocks.
\ADFvfyParStart{(215, ((20, 52, ((156, 3), (39, 1), (20, 20))), (14, 156, ((156, 1), (39, 1), (20, 20)))), ((26, 6), (59, 1)))} 

\adfDgap
\noindent{\boldmath $ 26^{6} 62^{1} $}~
With the point set $Z_{218}$ partitioned into
 residue classes modulo $6$ for $\{0, 1, \dots, 155\}$, and
 $\{156, 157, \dots, 217\}$,
 the design is generated from

\adfLgap 
$(195, 0, 1, 2)$,
$(196, 0, 52, 50)$,
$(197, 0, 103, 155)$,
$(198, 0, 106, 53)$,\adfsplit
$(199, 0, 154, 104)$,
$(200, 0, 4, 11)$,
$(201, 0, 7, 152)$,
$(202, 0, 145, 149)$,\adfsplit
$(203, 0, 10, 23)$,
$(204, 0, 13, 146)$,
$(205, 0, 133, 143)$,
$(206, 0, 16, 35)$,\adfsplit
$(207, 0, 19, 140)$,
$(208, 0, 121, 137)$,
$(209, 0, 22, 5)$,
$(210, 0, 139, 134)$,\adfsplit
$(211, 0, 151, 17)$,
$(212, 0, 25, 56)$,
$(213, 0, 31, 131)$,
$(214, 0, 100, 125)$,\adfsplit
$(215, 0, 28, 8)$,
$(216, 0, 136, 128)$,
$(217, 0, 148, 20)$,
$(156, 144, 51, 22)$,\adfsplit
$(156, 80, 53, 150)$,
$(156, 41, 115, 26)$,
$(156, 92, 15, 77)$,
$(156, 136, 31, 75)$,\adfsplit
$(156, 45, 137, 90)$,
$(156, 27, 122, 65)$,
$(156, 69, 24, 116)$,
$(156, 91, 33, 132)$,\adfsplit
$(156, 85, 39, 124)$,
$(156, 29, 20, 108)$,
$(156, 64, 138, 101)$,
$(0, 62, 147, 178)$,\adfsplit
$(0, 39, 107, 181)$,
$(0, 93, 127, 192)$,
$(0, 38, 87, 175)$,
$(0, 55, 135, 179)$,\adfsplit
$(0, 44, 113, 191)$,
$(0, 75, 101, 170)$,
$(0, 21, 80, 158)$,
$(0, 3, 73, 161)$,\adfsplit
$(0, 14, 41, 81)$,
$(0, 43, 46, 177)$,
$(0, 33, 65, 171)$,
$(0, 26, 58, 164)$,\adfsplit
$(0, 40, 123, 187)$,
$(0, 91, 105, 160)$

\adfLgap \noindent by the mapping:
$x \mapsto x + 3 j \adfmod{156}$ for $x < 156$,
$x \mapsto (x +  j \adfmod{39}) + 156$ for $156 \le x < 195$,
$x \mapsto x$ for $x \ge 195$,
$0 \le j < 52$
 for the first 23 blocks;
$x \mapsto x + 2 j \adfmod{156}$ for $x < 156$,
$x \mapsto (x +  j \adfmod{39}) + 156$ for $156 \le x < 195$,
$x \mapsto x$ for $x \ge 195$,
$0 \le j < 78$
 for the last 27 blocks.
\ADFvfyParStart{(218, ((23, 52, ((156, 3), (39, 1), (23, 23))), (27, 78, ((156, 2), (39, 1), (23, 23)))), ((26, 6), (62, 1)))} 

\adfDgap
\noindent{\boldmath $ 26^{9} 11^{1} $}~
With the point set $Z_{245}$ partitioned into
 residue classes modulo $9$ for $\{0, 1, \dots, 233\}$, and
 $\{234, 235, \dots, 244\}$,
 the design is generated from

\adfLgap 
$(243, 138, 193, 123)$,
$(243, 164, 185, 4)$,
$(244, 75, 209, 82)$,
$(244, 158, 25, 96)$,\adfsplit
$(234, 225, 65, 197)$,
$(234, 185, 82, 210)$,
$(234, 3, 8, 13)$,
$(234, 105, 74, 32)$,\adfsplit
$(234, 232, 184, 199)$,
$(234, 36, 168, 61)$,
$(235, 4, 210, 60)$,
$(235, 62, 193, 97)$,\adfsplit
$(235, 146, 216, 51)$,
$(235, 59, 147, 16)$,
$(235, 217, 101, 32)$,
$(235, 45, 143, 118)$,\adfsplit
$(236, 45, 187, 40)$,
$(236, 178, 113, 28)$,
$(236, 65, 18, 183)$,
$(236, 87, 74, 48)$,\adfsplit
$(236, 107, 194, 31)$,
$(236, 98, 114, 37)$,
$(40, 56, 223, 19)$,
$(43, 125, 5, 57)$,\adfsplit
$(0, 42, 227, 214)$,
$(122, 37, 134, 79)$,
$(97, 125, 117, 154)$,
$(105, 79, 233, 118)$,\adfsplit
$(214, 89, 41, 164)$,
$(231, 8, 1, 185)$,
$(171, 14, 111, 47)$,
$(11, 180, 16, 181)$,\adfsplit
$(152, 104, 20, 142)$,
$(62, 155, 225, 169)$,
$(139, 5, 78, 21)$,
$(152, 100, 155, 32)$,\adfsplit
$(83, 27, 210, 193)$,
$(180, 182, 152, 61)$,
$(155, 186, 215, 138)$,
$(123, 94, 152, 18)$,\adfsplit
$(39, 59, 9, 196)$,
$(157, 162, 104, 129)$,
$(15, 140, 134, 27)$,
$(231, 82, 169, 221)$,\adfsplit
$(184, 217, 129, 180)$,
$(30, 41, 53, 7)$,
$(211, 86, 171, 219)$,
$(146, 113, 222, 46)$,\adfsplit
$(111, 143, 160, 68)$,
$(147, 125, 76, 106)$,
$(117, 231, 35, 50)$,
$(136, 204, 43, 113)$,\adfsplit
$(63, 65, 193, 109)$,
$(67, 120, 0, 233)$,
$(233, 141, 176, 7)$,
$(53, 43, 46, 166)$,\adfsplit
$(65, 151, 68, 69)$,
$(95, 90, 79, 4)$,
$(191, 31, 35, 0)$,
$(185, 22, 170, 196)$,\adfsplit
$(145, 160, 22, 0)$,
$(139, 170, 91, 126)$,
$(73, 133, 30, 9)$,
$(118, 36, 119, 211)$,\adfsplit
$(147, 92, 115, 190)$,
$(0, 6, 95, 179)$,
$(0, 7, 215, 221)$,
$(0, 19, 131, 155)$,\adfsplit
$(0, 14, 53, 119)$,
$(0, 12, 149, 151)$,
$(0, 80, 109, 167)$,
$(0, 41, 71, 92)$,\adfsplit
$(0, 15, 59, 121)$,
$(1, 7, 81, 95)$,
$(1, 41, 99, 130)$,
$(0, 21, 79, 101)$,\adfsplit
$(1, 3, 25, 203)$,
$(3, 28, 65, 161)$,
$(4, 10, 161, 203)$,
$(1, 45, 131, 178)$,\adfsplit
$(0, 32, 93, 197)$,
$(2, 21, 112, 143)$,
$(2, 53, 147, 208)$,
$(2, 68, 171, 197)$,\adfsplit
$(3, 22, 124, 197)$,
$(1, 2, 13, 79)$,
$(0, 91, 193, 231)$,
$(0, 49, 74, 169)$,\adfsplit
$(1, 87, 166, 219)$,
$(0, 33, 85, 230)$,
$(1, 51, 194, 196)$,
$(1, 159, 201, 218)$,\adfsplit
$(0, 112, 133, 152)$,
$(0, 146, 202, 205)$,
$(0, 68, 75, 187)$,
$(0, 97, 110, 170)$,\adfsplit
$(1, 106, 184, 188)$,
$(0, 129, 130, 195)$,
$(0, 3, 141, 188)$,
$(3, 9, 87, 226)$,\adfsplit
$(0, 111, 194, 226)$,
$(2, 106, 130, 213)$,
$(0, 88, 123, 147)$,
$(0, 136, 159, 200)$,\adfsplit
$(0, 30, 87, 184)$,
$(2, 81, 148, 190)$,
$(0, 24, 122, 142)$,
$(0, 46, 116, 140)$,\adfsplit
$(0, 10, 86, 182)$,
$(0, 20, 34, 100)$,
$(0, 16, 38, 60)$,
$(0, 50, 94, 128)$,\adfsplit
$(0, 8, 124, 138)$,
$(0, 56, 64, 66)$,
$(0, 40, 52, 78)$

\adfLgap \noindent by the mapping:
$x \mapsto x + 6 j \adfmod{234}$ for $x < 234$,
$x \mapsto (x + 3 j \adfmod{9}) + 234$ for $234 \le x < 243$,
$x \mapsto x$ for $x \ge 243$,
$0 \le j < 39$.
\ADFvfyParStart{(245, ((115, 39, ((234, 6), (9, 3), (2, 2)))), ((26, 9), (11, 1)))} 

\adfDgap
\noindent{\boldmath $ 26^{9} 14^{1} $}~
With the point set $Z_{248}$ partitioned into
 residue classes modulo $9$ for $\{0, 1, \dots, 233\}$, and
 $\{234, 235, \dots, 247\}$,
 the design is generated from

\adfLgap 
$(246, 171, 107, 73)$,
$(247, 121, 98, 99)$,
$(234, 7, 75, 149)$,
$(234, 134, 130, 74)$,\adfsplit
$(234, 81, 114, 181)$,
$(235, 102, 188, 193)$,
$(235, 38, 90, 51)$,
$(235, 34, 59, 100)$,\adfsplit
$(236, 39, 142, 164)$,
$(236, 161, 140, 208)$,
$(236, 198, 33, 193)$,
$(237, 80, 214, 60)$,\adfsplit
$(237, 175, 174, 217)$,
$(237, 72, 164, 212)$,
$(162, 220, 30, 55)$,
$(152, 183, 67, 34)$,\adfsplit
$(142, 53, 77, 102)$,
$(215, 92, 225, 85)$,
$(118, 7, 144, 215)$,
$(5, 84, 37, 76)$,\adfsplit
$(40, 192, 160, 99)$,
$(129, 209, 126, 76)$,
$(104, 232, 233, 36)$,
$(61, 154, 76, 200)$,\adfsplit
$(112, 24, 101, 140)$,
$(12, 64, 62, 200)$,
$(182, 61, 189, 1)$,
$(100, 135, 49, 51)$,\adfsplit
$(84, 203, 137, 121)$,
$(143, 149, 193, 223)$,
$(41, 31, 79, 98)$,
$(202, 56, 183, 196)$,\adfsplit
$(206, 1, 214, 41)$,
$(0, 2, 64, 223)$,
$(0, 5, 112, 169)$,
$(0, 4, 151, 228)$,\adfsplit
$(0, 8, 193, 205)$,
$(0, 11, 130, 214)$,
$(0, 7, 31, 186)$,
$(0, 17, 76, 178)$,\adfsplit
$(0, 12, 28, 133)$,
$(0, 73, 176, 211)$,
$(0, 23, 111, 145)$,
$(1, 4, 80, 221)$,\adfsplit
$(0, 109, 167, 200)$,
$(0, 14, 139, 210)$,
$(0, 104, 124, 179)$,
$(0, 29, 143, 220)$,\adfsplit
$(0, 32, 115, 177)$,
$(0, 56, 142, 206)$,
$(0, 35, 147, 217)$,
$(0, 42, 191, 194)$,\adfsplit
$(0, 59, 137, 159)$,
$(0, 131, 161, 173)$,
$(0, 21, 65, 204)$,
$(0, 60, 156, 218)$,\adfsplit
$(0, 15, 120, 230)$,
$(0, 47, 98, 168)$,
$(0, 26, 41, 128)$

\adfLgap \noindent by the mapping:
$x \mapsto x + 3 j \adfmod{234}$ for $x < 234$,
$x \mapsto (x - 234 + 4 j \adfmod{12}) + 234$ for $234 \le x < 246$,
$x \mapsto x$ for $x \ge 246$,
$0 \le j < 78$.
\ADFvfyParStart{(248, ((59, 78, ((234, 3), (12, 4), (2, 2)))), ((26, 9), (14, 1)))} 

\adfDgap
\noindent{\boldmath $ 26^{9} 17^{1} $}~
With the point set $Z_{251}$ partitioned into
 residue classes modulo $9$ for $\{0, 1, \dots, 233\}$, and
 $\{234, 235, \dots, 250\}$,
 the design is generated from

\adfLgap 
$(249, 228, 91, 142)$,
$(249, 189, 197, 68)$,
$(250, 131, 195, 97)$,
$(250, 190, 8, 30)$,\adfsplit
$(234, 91, 122, 159)$,
$(234, 146, 219, 153)$,
$(234, 125, 54, 149)$,
$(234, 170, 34, 175)$,\adfsplit
$(234, 228, 83, 76)$,
$(234, 154, 7, 60)$,
$(235, 2, 117, 163)$,
$(235, 159, 219, 88)$,\adfsplit
$(235, 224, 205, 226)$,
$(235, 89, 211, 95)$,
$(235, 101, 198, 168)$,
$(235, 228, 50, 184)$,\adfsplit
$(236, 179, 225, 150)$,
$(236, 91, 131, 106)$,
$(236, 116, 82, 164)$,
$(236, 22, 155, 48)$,\adfsplit
$(236, 25, 21, 139)$,
$(236, 195, 140, 90)$,
$(237, 146, 122, 172)$,
$(237, 37, 42, 94)$,\adfsplit
$(237, 225, 223, 95)$,
$(237, 52, 30, 157)$,
$(237, 129, 33, 62)$,
$(237, 108, 125, 47)$,\adfsplit
$(238, 4, 194, 120)$,
$(238, 28, 157, 79)$,
$(238, 81, 37, 60)$,
$(238, 123, 16, 180)$,\adfsplit
$(238, 71, 110, 3)$,
$(238, 131, 47, 116)$,
$(148, 92, 125, 86)$,
$(19, 104, 168, 81)$,\adfsplit
$(155, 150, 143, 85)$,
$(174, 155, 90, 1)$,
$(151, 138, 154, 107)$,
$(39, 77, 46, 29)$,\adfsplit
$(155, 170, 133, 67)$,
$(108, 174, 172, 166)$,
$(3, 98, 118, 155)$,
$(20, 72, 194, 208)$,\adfsplit
$(0, 4, 116, 20)$,
$(233, 19, 232, 31)$,
$(124, 87, 57, 17)$,
$(133, 217, 105, 215)$,\adfsplit
$(132, 166, 56, 98)$,
$(193, 7, 176, 118)$,
$(231, 188, 49, 183)$,
$(207, 182, 168, 201)$,\adfsplit
$(146, 25, 81, 53)$,
$(109, 75, 59, 49)$,
$(73, 11, 123, 224)$,
$(139, 0, 32, 110)$,\adfsplit
$(193, 59, 162, 8)$,
$(140, 129, 53, 6)$,
$(85, 233, 117, 160)$,
$(191, 72, 16, 71)$,\adfsplit
$(7, 111, 152, 231)$,
$(156, 9, 115, 62)$,
$(0, 103, 141, 133)$,
$(78, 16, 104, 28)$,\adfsplit
$(85, 71, 12, 225)$,
$(218, 217, 201, 7)$,
$(128, 37, 203, 57)$,
$(85, 36, 60, 46)$,\adfsplit
$(133, 220, 114, 225)$,
$(31, 100, 92, 35)$,
$(59, 6, 200, 192)$,
$(127, 229, 221, 12)$,\adfsplit
$(0, 2, 12, 121)$,
$(0, 41, 67, 114)$,
$(0, 3, 146, 157)$,
$(0, 1, 7, 40)$,\adfsplit
$(0, 15, 37, 79)$,
$(0, 55, 69, 151)$,
$(0, 38, 51, 175)$,
$(0, 43, 68, 152)$,\adfsplit
$(0, 86, 163, 206)$,
$(0, 11, 44, 199)$,
$(0, 113, 164, 223)$,
$(1, 56, 188, 220)$,\adfsplit
$(0, 88, 91, 98)$,
$(1, 65, 194, 197)$,
$(1, 165, 178, 200)$,
$(1, 68, 137, 179)$,\adfsplit
$(1, 14, 75, 159)$,
$(1, 87, 131, 134)$,
$(0, 112, 205, 221)$,
$(1, 147, 161, 166)$,\adfsplit
$(1, 53, 124, 164)$,
$(1, 77, 143, 177)$,
$(0, 62, 92, 93)$,
$(0, 42, 218, 230)$,\adfsplit
$(0, 128, 214, 231)$,
$(0, 28, 60, 138)$,
$(0, 23, 102, 159)$,
$(0, 6, 197, 201)$,\adfsplit
$(0, 83, 166, 196)$,
$(0, 76, 104, 142)$,
$(2, 68, 153, 177)$,
$(2, 51, 82, 106)$,\adfsplit
$(2, 94, 142, 161)$,
$(2, 76, 99, 172)$,
$(2, 89, 118, 232)$,
$(2, 130, 159, 215)$,\adfsplit
$(0, 100, 149, 179)$,
$(0, 125, 183, 185)$,
$(0, 130, 143, 165)$,
$(0, 35, 46, 129)$,\adfsplit
$(0, 77, 154, 219)$,
$(3, 35, 88, 131)$,
$(3, 53, 148, 190)$,
$(3, 29, 81, 160)$,\adfsplit
$(3, 28, 89, 178)$,
$(3, 15, 77, 112)$,
$(3, 23, 64, 142)$,
$(3, 4, 45, 136)$,\adfsplit
$(3, 52, 105, 226)$

\adfLgap \noindent by the mapping:
$x \mapsto x + 6 j \adfmod{234}$ for $x < 234$,
$x \mapsto (x - 234 + 5 j \adfmod{15}) + 234$ for $234 \le x < 249$,
$x \mapsto x$ for $x \ge 249$,
$0 \le j < 39$.
\ADFvfyParStart{(251, ((121, 39, ((234, 6), (15, 5), (2, 2)))), ((26, 9), (17, 1)))} 

\adfDgap
\noindent{\boldmath $ 26^{9} 20^{1} $}~
With the point set $Z_{254}$ partitioned into
 residue classes modulo $9$ for $\{0, 1, \dots, 233\}$, and
 $\{234, 235, \dots, 253\}$,
 the design is generated from

\adfLgap 
$(252, 215, 186, 217)$,
$(253, 232, 126, 53)$,
$(234, 162, 172, 123)$,
$(234, 92, 221, 70)$,\adfsplit
$(234, 76, 80, 21)$,
$(235, 112, 209, 87)$,
$(235, 32, 117, 16)$,
$(235, 201, 224, 109)$,\adfsplit
$(236, 184, 106, 32)$,
$(236, 144, 69, 164)$,
$(236, 170, 163, 102)$,
$(237, 174, 34, 91)$,\adfsplit
$(237, 14, 80, 139)$,
$(237, 33, 146, 162)$,
$(238, 44, 81, 102)$,
$(238, 13, 32, 199)$,\adfsplit
$(238, 200, 204, 169)$,
$(239, 103, 23, 187)$,
$(239, 143, 141, 9)$,
$(239, 12, 155, 28)$,\adfsplit
$(74, 202, 39, 115)$,
$(195, 73, 176, 13)$,
$(107, 220, 208, 216)$,
$(192, 154, 69, 146)$,\adfsplit
$(211, 14, 26, 46)$,
$(99, 58, 59, 114)$,
$(29, 163, 176, 18)$,
$(19, 115, 87, 176)$,\adfsplit
$(150, 81, 184, 3)$,
$(12, 24, 38, 86)$,
$(171, 83, 51, 221)$,
$(53, 147, 211, 155)$,\adfsplit
$(52, 38, 197, 23)$,
$(5, 15, 170, 205)$,
$(60, 161, 212, 156)$,
$(199, 107, 79, 213)$,\adfsplit
$(185, 206, 102, 160)$,
$(89, 86, 119, 120)$,
$(0, 1, 41, 65)$,
$(0, 3, 110, 167)$,\adfsplit
$(0, 6, 53, 137)$,
$(0, 7, 80, 119)$,
$(0, 13, 71, 182)$,
$(0, 38, 116, 154)$,\adfsplit
$(0, 98, 124, 212)$,
$(1, 4, 95, 188)$,
$(1, 80, 86, 103)$,
$(0, 206, 211, 217)$,\adfsplit
$(0, 48, 169, 221)$,
$(0, 82, 92, 187)$,
$(0, 17, 60, 127)$,
$(0, 44, 86, 148)$,\adfsplit
$(0, 46, 79, 227)$,
$(1, 31, 124, 167)$,
$(0, 19, 24, 209)$,
$(0, 66, 136, 157)$,\adfsplit
$(0, 30, 118, 160)$,
$(0, 37, 57, 232)$,
$(0, 43, 93, 202)$,
$(0, 40, 208, 223)$,\adfsplit
$(0, 22, 51, 84)$,
$(0, 42, 115, 139)$

\adfLgap \noindent by the mapping:
$x \mapsto x + 3 j \adfmod{234}$ for $x < 234$,
$x \mapsto (x + 6 j \adfmod{18}) + 234$ for $234 \le x < 252$,
$x \mapsto x$ for $x \ge 252$,
$0 \le j < 78$.
\ADFvfyParStart{(254, ((62, 78, ((234, 3), (18, 6), (2, 2)))), ((26, 9), (20, 1)))} 

\adfDgap
\noindent{\boldmath $ 26^{9} 23^{1} $}~
With the point set $Z_{257}$ partitioned into
 residue classes modulo $9$ for $\{0, 1, \dots, 233\}$, and
 $\{234, 235, \dots, 256\}$,
 the design is generated from

\adfLgap 
$(255, 183, 161, 184)$,
$(255, 37, 98, 24)$,
$(256, 194, 132, 29)$,
$(256, 88, 175, 171)$,\adfsplit
$(234, 186, 35, 173)$,
$(234, 109, 32, 202)$,
$(234, 15, 39, 80)$,
$(234, 117, 120, 131)$,\adfsplit
$(234, 146, 151, 229)$,
$(234, 136, 162, 196)$,
$(235, 20, 21, 64)$,
$(235, 99, 150, 106)$,\adfsplit
$(235, 206, 33, 91)$,
$(235, 131, 130, 18)$,
$(235, 233, 173, 30)$,
$(235, 104, 139, 169)$,\adfsplit
$(236, 224, 93, 199)$,
$(236, 208, 49, 20)$,
$(236, 72, 196, 213)$,
$(236, 223, 59, 71)$,\adfsplit
$(236, 174, 58, 45)$,
$(236, 194, 150, 65)$,
$(237, 98, 0, 96)$,
$(237, 64, 56, 7)$,\adfsplit
$(237, 111, 91, 214)$,
$(237, 89, 229, 30)$,
$(237, 104, 101, 45)$,
$(237, 95, 130, 195)$,\adfsplit
$(238, 66, 52, 73)$,
$(238, 79, 143, 184)$,
$(238, 9, 131, 159)$,
$(238, 186, 74, 201)$,\adfsplit
$(238, 36, 13, 140)$,
$(238, 154, 206, 65)$,
$(239, 2, 143, 220)$,
$(239, 212, 45, 219)$,\adfsplit
$(239, 118, 87, 174)$,
$(239, 227, 60, 97)$,
$(239, 145, 0, 77)$,
$(239, 142, 193, 44)$,\adfsplit
$(240, 137, 112, 120)$,
$(240, 164, 96, 127)$,
$(240, 223, 149, 3)$,
$(240, 215, 32, 180)$,\adfsplit
$(240, 170, 232, 195)$,
$(240, 193, 190, 63)$,
$(44, 225, 222, 145)$,
$(228, 2, 176, 196)$,\adfsplit
$(130, 24, 221, 74)$,
$(50, 139, 36, 3)$,
$(103, 178, 77, 24)$,
$(158, 110, 37, 9)$,\adfsplit
$(20, 96, 58, 100)$,
$(58, 214, 192, 164)$,
$(54, 183, 25, 47)$,
$(3, 151, 117, 119)$,\adfsplit
$(186, 147, 104, 136)$,
$(71, 121, 128, 132)$,
$(107, 45, 12, 137)$,
$(80, 144, 219, 46)$,\adfsplit
$(146, 69, 133, 212)$,
$(152, 175, 28, 41)$,
$(80, 200, 5, 54)$,
$(179, 7, 108, 50)$,\adfsplit
$(222, 148, 54, 136)$,
$(195, 113, 103, 61)$,
$(232, 112, 17, 33)$,
$(69, 28, 107, 229)$,\adfsplit
$(146, 64, 158, 18)$,
$(213, 126, 136, 166)$,
$(62, 99, 148, 51)$,
$(166, 32, 150, 27)$,\adfsplit
$(184, 159, 110, 125)$,
$(102, 186, 108, 40)$,
$(226, 40, 142, 89)$,
$(52, 50, 67, 213)$,\adfsplit
$(125, 169, 172, 195)$,
$(0, 12, 70, 200)$,
$(1, 5, 16, 22)$,
$(0, 57, 142, 218)$,\adfsplit
$(0, 47, 76, 194)$,
$(0, 5, 20, 232)$,
$(1, 2, 70, 139)$,
$(1, 40, 45, 178)$,\adfsplit
$(0, 28, 38, 52)$,
$(0, 64, 92, 123)$,
$(1, 35, 56, 202)$,
$(0, 21, 110, 134)$,\adfsplit
$(1, 7, 110, 188)$,
$(2, 8, 21, 140)$,
$(0, 1, 32, 224)$,
$(2, 28, 86, 197)$,\adfsplit
$(1, 39, 164, 194)$,
$(2, 53, 59, 232)$,
$(0, 73, 97, 164)$,
$(0, 80, 101, 102)$,\adfsplit
$(1, 153, 176, 179)$,
$(2, 57, 124, 153)$,
$(2, 35, 166, 203)$,
$(2, 77, 201, 207)$,\adfsplit
$(2, 94, 167, 219)$,
$(1, 3, 121, 140)$,
$(3, 118, 184, 227)$,
$(0, 23, 130, 215)$,\adfsplit
$(3, 23, 94, 207)$,
$(1, 29, 123, 142)$,
$(3, 83, 166, 197)$,
$(1, 77, 141, 196)$,\adfsplit
$(0, 69, 85, 214)$,
$(1, 47, 52, 175)$,
$(1, 57, 69, 225)$,
$(1, 59, 63, 137)$,\adfsplit
$(1, 33, 119, 165)$,
$(0, 48, 137, 213)$,
$(1, 189, 197, 221)$,
$(0, 43, 51, 93)$,\adfsplit
$(0, 39, 91, 107)$,
$(3, 47, 131, 179)$,
$(0, 49, 159, 181)$,
$(0, 25, 42, 219)$,\adfsplit
$(0, 30, 191, 193)$,
$(0, 29, 109, 121)$,
$(0, 19, 114, 169)$,
$(0, 60, 179, 187)$,\adfsplit
$(0, 24, 65, 175)$,
$(0, 61, 209, 229)$,
$(0, 67, 115, 155)$

\adfLgap \noindent by the mapping:
$x \mapsto x + 6 j \adfmod{234}$ for $x < 234$,
$x \mapsto (x - 234 + 7 j \adfmod{21}) + 234$ for $234 \le x < 255$,
$x \mapsto x$ for $x \ge 255$,
$0 \le j < 39$.
\ADFvfyParStart{(257, ((127, 39, ((234, 6), (21, 7), (2, 2)))), ((26, 9), (23, 1)))} 

\adfDgap
\noindent{\boldmath $ 26^{12} 14^{1} $}~
With the point set $Z_{326}$ partitioned into
 residue classes modulo $12$ for $\{0, 1, \dots, 311\}$, and
 $\{312, 313, \dots, 325\}$,
 the design is generated from

\adfLgap 
$(312, 0, 1, 2)$,
$(313, 0, 103, 209)$,
$(314, 0, 106, 104)$,
$(315, 0, 208, 311)$,\adfsplit
$(316, 0, 310, 206)$,
$(317, 0, 4, 11)$,
$(318, 0, 7, 308)$,
$(319, 0, 301, 305)$,\adfsplit
$(320, 0, 10, 23)$,
$(321, 0, 13, 302)$,
$(322, 0, 289, 299)$,
$(323, 0, 16, 35)$,\adfsplit
$(324, 0, 19, 296)$,
$(325, 0, 277, 293)$,
$(79, 62, 248, 11)$,
$(73, 6, 123, 213)$,\adfsplit
$(311, 202, 51, 115)$,
$(271, 145, 243, 52)$,
$(239, 86, 284, 270)$,
$(300, 142, 78, 83)$,\adfsplit
$(89, 235, 272, 122)$,
$(63, 289, 206, 107)$,
$(27, 245, 156, 203)$,
$(102, 212, 134, 241)$,\adfsplit
$(49, 188, 107, 264)$,
$(237, 106, 303, 294)$,
$(110, 104, 57, 291)$,
$(202, 252, 230, 301)$,\adfsplit
$(264, 54, 224, 57)$,
$(301, 6, 248, 100)$,
$(184, 48, 10, 205)$,
$(272, 181, 127, 149)$,\adfsplit
$(298, 176, 120, 255)$,
$(203, 111, 138, 129)$,
$(307, 184, 282, 107)$,
$(294, 273, 74, 253)$,\adfsplit
$(97, 56, 89, 12)$,
$(137, 183, 222, 122)$,
$(205, 150, 84, 280)$,
$(0, 18, 81, 86)$,\adfsplit
$(0, 8, 42, 62)$,
$(0, 43, 73, 152)$,
$(0, 26, 172, 211)$,
$(0, 51, 112, 263)$,\adfsplit
$(0, 80, 243, 257)$,
$(0, 52, 111, 245)$,
$(0, 27, 30, 224)$,
$(0, 69, 223, 229)$,\adfsplit
$(0, 74, 215, 249)$,
$(0, 133, 159, 283)$,
$(0, 46, 147, 217)$,
$(0, 58, 149, 297)$,\adfsplit
$(0, 127, 165, 241)$,
$(0, 45, 125, 247)$,
$(0, 37, 230, 287)$,
$(0, 31, 87, 205)$,\adfsplit
$(1, 41, 129, 211)$

\adfLgap \noindent by the mapping:
$x \mapsto x + 3 j \adfmod{312}$ for $x < 312$,
$x \mapsto x$ for $x \ge 312$,
$0 \le j < 104$
 for the first 14 blocks;
$x \mapsto x + 2 j \adfmod{312}$ for $x < 312$,
$x \mapsto x$ for $x \ge 312$,
$0 \le j < 156$
 for the last 43 blocks.
\ADFvfyParStart{(326, ((14, 104, ((312, 3), (14, 14))), (43, 156, ((312, 2), (14, 14)))), ((26, 12), (14, 1)))} 

\adfDgap
\noindent{\boldmath $ 26^{12} 17^{1} $}~
With the point set $Z_{329}$ partitioned into
 residue classes modulo $12$ for $\{0, 1, \dots, 311\}$, and
 $\{312, 313, \dots, 328\}$,
 the design is generated from

\adfLgap 
$(312, 0, 1, 2)$,
$(313, 0, 103, 209)$,
$(314, 0, 106, 104)$,
$(315, 0, 208, 311)$,\adfsplit
$(316, 0, 310, 206)$,
$(317, 0, 4, 11)$,
$(318, 0, 7, 308)$,
$(319, 0, 301, 305)$,\adfsplit
$(320, 0, 10, 23)$,
$(321, 0, 13, 302)$,
$(322, 0, 289, 299)$,
$(323, 0, 16, 35)$,\adfsplit
$(324, 0, 19, 296)$,
$(325, 0, 277, 293)$,
$(326, 0, 22, 5)$,
$(327, 0, 295, 290)$,\adfsplit
$(328, 0, 307, 17)$,
$(286, 122, 195, 81)$,
$(7, 145, 254, 34)$,
$(236, 181, 274, 131)$,\adfsplit
$(20, 151, 159, 286)$,
$(94, 246, 32, 183)$,
$(188, 146, 1, 4)$,
$(286, 306, 187, 212)$,\adfsplit
$(164, 173, 261, 78)$,
$(149, 115, 196, 70)$,
$(226, 102, 156, 257)$,
$(18, 141, 167, 98)$,\adfsplit
$(0, 37, 77, 236)$,
$(0, 21, 85, 137)$,
$(0, 15, 71, 225)$,
$(0, 59, 134, 200)$,\adfsplit
$(0, 30, 83, 140)$,
$(0, 28, 118, 150)$,
$(0, 14, 58, 179)$,
$(0, 6, 67, 182)$,\adfsplit
$(0, 18, 51, 100)$,
$(0, 29, 68, 146)$

\adfLgap \noindent by the mapping:
$x \mapsto x + 3 j \adfmod{312}$ for $x < 312$,
$x \mapsto x$ for $x \ge 312$,
$0 \le j < 104$
 for the first 17 blocks;
$x \mapsto x +  j \adfmod{312}$ for $x < 312$,
$x \mapsto x$ for $x \ge 312$,
$0 \le j < 312$
 for the last 21 blocks.
\ADFvfyParStart{(329, ((17, 104, ((312, 3), (17, 17))), (21, 312, ((312, 1), (17, 17)))), ((26, 12), (17, 1)))} 

\adfDgap
\noindent{\boldmath $ 26^{12} 20^{1} $}~
With the point set $Z_{332}$ partitioned into
 residue classes modulo $12$ for $\{0, 1, \dots, 311\}$, and
 $\{312, 313, \dots, 331\}$,
 the design is generated from

\adfLgap 
$(312, 0, 1, 2)$,
$(313, 0, 103, 209)$,
$(314, 0, 106, 104)$,
$(315, 0, 208, 311)$,\adfsplit
$(316, 0, 310, 206)$,
$(317, 0, 4, 11)$,
$(318, 0, 7, 308)$,
$(319, 0, 301, 305)$,\adfsplit
$(320, 0, 10, 23)$,
$(321, 0, 13, 302)$,
$(322, 0, 289, 299)$,
$(323, 0, 16, 35)$,\adfsplit
$(324, 0, 19, 296)$,
$(325, 0, 277, 293)$,
$(326, 0, 22, 5)$,
$(327, 0, 295, 290)$,\adfsplit
$(328, 0, 307, 17)$,
$(329, 0, 25, 53)$,
$(330, 0, 28, 287)$,
$(331, 0, 259, 284)$,\adfsplit
$(73, 7, 268, 182)$,
$(162, 12, 148, 241)$,
$(162, 147, 278, 60)$,
$(188, 36, 77, 273)$,\adfsplit
$(223, 64, 269, 263)$,
$(232, 106, 86, 1)$,
$(74, 4, 289, 245)$,
$(212, 102, 169, 58)$,\adfsplit
$(75, 175, 293, 301)$,
$(168, 310, 301, 138)$,
$(107, 158, 44, 89)$,
$(76, 233, 264, 270)$,\adfsplit
$(259, 121, 256, 168)$,
$(214, 138, 253, 165)$,
$(63, 265, 16, 113)$,
$(162, 122, 267, 309)$,\adfsplit
$(193, 3, 283, 122)$,
$(103, 294, 11, 41)$,
$(105, 27, 250, 72)$,
$(140, 11, 49, 88)$,\adfsplit
$(1, 118, 186, 259)$,
$(279, 46, 259, 210)$,
$(285, 172, 143, 270)$,
$(0, 8, 139, 197)$,\adfsplit
$(0, 9, 42, 149)$,
$(0, 18, 50, 250)$,
$(0, 37, 58, 248)$,
$(0, 61, 135, 222)$,\adfsplit
$(0, 46, 123, 184)$,
$(0, 43, 182, 256)$,
$(0, 82, 191, 247)$,
$(0, 34, 153, 246)$,\adfsplit
$(0, 38, 217, 258)$,
$(0, 21, 78, 253)$,
$(0, 75, 223, 257)$,
$(0, 55, 125, 239)$,\adfsplit
$(1, 53, 135, 211)$,
$(0, 65, 211, 309)$,
$(0, 57, 83, 245)$,
$(0, 26, 155, 169)$,\adfsplit
$(0, 31, 95, 231)$

\adfLgap \noindent by the mapping:
$x \mapsto x + 3 j \adfmod{312}$ for $x < 312$,
$x \mapsto x$ for $x \ge 312$,
$0 \le j < 104$
 for the first 20 blocks;
$x \mapsto x + 2 j \adfmod{312}$ for $x < 312$,
$x \mapsto x$ for $x \ge 312$,
$0 \le j < 156$
 for the last 41 blocks.
\ADFvfyParStart{(332, ((20, 104, ((312, 3), (20, 20))), (41, 156, ((312, 2), (20, 20)))), ((26, 12), (20, 1)))} 

\adfDgap
\noindent{\boldmath $ 26^{12} 23^{1} $}~
With the point set $Z_{335}$ partitioned into
 residue classes modulo $12$ for $\{0, 1, \dots, 311\}$, and
 $\{312, 313, \dots, 334\}$,
 the design is generated from

\adfLgap 
$(312, 0, 1, 2)$,
$(313, 0, 103, 209)$,
$(314, 0, 106, 104)$,
$(315, 0, 208, 311)$,\adfsplit
$(316, 0, 310, 206)$,
$(317, 0, 4, 11)$,
$(318, 0, 7, 308)$,
$(319, 0, 301, 305)$,\adfsplit
$(320, 0, 10, 23)$,
$(321, 0, 13, 302)$,
$(322, 0, 289, 299)$,
$(323, 0, 16, 35)$,\adfsplit
$(324, 0, 19, 296)$,
$(325, 0, 277, 293)$,
$(326, 0, 22, 5)$,
$(327, 0, 295, 290)$,\adfsplit
$(328, 0, 307, 17)$,
$(329, 0, 25, 53)$,
$(330, 0, 28, 287)$,
$(331, 0, 259, 284)$,\adfsplit
$(332, 0, 31, 65)$,
$(333, 0, 34, 281)$,
$(334, 0, 247, 278)$,
$(230, 275, 281, 238)$,\adfsplit
$(220, 158, 104, 217)$,
$(50, 165, 244, 133)$,
$(41, 303, 246, 176)$,
$(202, 109, 138, 38)$,\adfsplit
$(244, 63, 83, 132)$,
$(183, 156, 64, 17)$,
$(290, 20, 153, 276)$,
$(187, 33, 121, 106)$,\adfsplit
$(141, 16, 290, 223)$,
$(0, 18, 140, 170)$,
$(0, 40, 114, 223)$,
$(0, 63, 143, 218)$,\adfsplit
$(0, 39, 134, 210)$,
$(0, 41, 126, 225)$,
$(0, 21, 98, 195)$,
$(0, 33, 91, 244)$,\adfsplit
$(0, 44, 130, 191)$,
$(0, 9, 55, 145)$,
$(0, 26, 78, 188)$

\adfLgap \noindent by the mapping:
$x \mapsto x + 3 j \adfmod{312}$ for $x < 312$,
$x \mapsto x$ for $x \ge 312$,
$0 \le j < 104$
 for the first 23 blocks;
$x \mapsto x +  j \adfmod{312}$ for $x < 312$,
$x \mapsto x$ for $x \ge 312$,
$0 \le j < 312$
 for the last 20 blocks.
\ADFvfyParStart{(335, ((23, 104, ((312, 3), (23, 23))), (20, 312, ((312, 1), (23, 23)))), ((26, 12), (23, 1)))} 

\adfDgap
\noindent{\boldmath $ 26^{15} 17^{1} $}~
With the point set $Z_{407}$ partitioned into
 residue classes modulo $15$ for $\{0, 1, \dots, 389\}$, and
 $\{390, 391, \dots, 406\}$,
 the design is generated from

\adfLgap 
$(390, 0, 1, 2)$,
$(391, 0, 130, 128)$,
$(392, 0, 259, 389)$,
$(393, 0, 262, 131)$,\adfsplit
$(394, 0, 388, 260)$,
$(395, 0, 4, 11)$,
$(396, 0, 7, 386)$,
$(397, 0, 379, 383)$,\adfsplit
$(398, 0, 10, 23)$,
$(399, 0, 13, 380)$,
$(400, 0, 367, 377)$,
$(401, 0, 16, 35)$,\adfsplit
$(402, 0, 19, 374)$,
$(403, 0, 355, 371)$,
$(404, 0, 22, 5)$,
$(405, 0, 373, 368)$,\adfsplit
$(406, 0, 385, 17)$,
$(228, 370, 192, 311)$,
$(24, 124, 193, 0)$,
$(225, 259, 103, 378)$,\adfsplit
$(371, 267, 151, 219)$,
$(132, 243, 331, 225)$,
$(263, 6, 150, 112)$,
$(78, 389, 260, 167)$,\adfsplit
$(157, 85, 343, 314)$,
$(388, 151, 84, 362)$,
$(358, 141, 86, 280)$,
$(128, 99, 74, 276)$,\adfsplit
$(119, 87, 205, 24)$,
$(108, 178, 158, 349)$,
$(103, 216, 161, 53)$,
$(241, 134, 386, 45)$,\adfsplit
$(132, 159, 79, 190)$,
$(51, 15, 129, 307)$,
$(25, 248, 239, 90)$,
$(361, 81, 108, 224)$,\adfsplit
$(324, 154, 97, 57)$,
$(19, 255, 362, 52)$,
$(387, 245, 70, 346)$,
$(233, 136, 94, 359)$,\adfsplit
$(127, 362, 226, 201)$,
$(291, 342, 74, 166)$,
$(308, 231, 137, 95)$,
$(117, 364, 58, 241)$,\adfsplit
$(246, 199, 378, 311)$,
$(156, 377, 159, 213)$,
$(49, 132, 368, 144)$,
$(253, 14, 319, 220)$,\adfsplit
$(389, 350, 147, 86)$,
$(356, 32, 248, 262)$,
$(289, 373, 227, 83)$,
$(81, 202, 236, 312)$,\adfsplit
$(0, 6, 46, 334)$,
$(0, 9, 156, 207)$,
$(0, 8, 109, 200)$,
$(0, 28, 72, 322)$,\adfsplit
$(0, 32, 341, 353)$,
$(0, 64, 152, 286)$,
$(0, 77, 186, 327)$,
$(0, 52, 133, 209)$,\adfsplit
$(0, 37, 98, 342)$,
$(0, 121, 164, 349)$,
$(0, 43, 201, 316)$,
$(0, 18, 145, 215)$,\adfsplit
$(0, 91, 147, 218)$,
$(0, 87, 189, 228)$,
$(0, 21, 308, 387)$,
$(0, 53, 73, 261)$,\adfsplit
$(0, 273, 287, 369)$,
$(0, 31, 205, 231)$,
$(1, 7, 45, 53)$,
$(1, 29, 93, 253)$

\adfLgap \noindent by the mapping:
$x \mapsto x + 3 j \adfmod{390}$ for $x < 390$,
$x \mapsto x$ for $x \ge 390$,
$0 \le j < 130$
 for the first 17 blocks;
$x \mapsto x + 2 j \adfmod{390}$ for $x < 390$,
$x \mapsto x$ for $x \ge 390$,
$0 \le j < 195$
 for the last 55 blocks.
\ADFvfyParStart{(407, ((17, 130, ((390, 3), (17, 17))), (55, 195, ((390, 2), (17, 17)))), ((26, 15), (17, 1)))} 

\adfDgap
\noindent{\boldmath $ 26^{15} 20^{1} $}~
With the point set $Z_{410}$ partitioned into
 residue classes modulo $15$ for $\{0, 1, \dots, 389\}$, and
 $\{390, 391, \dots, 409\}$,
 the design is generated from

\adfLgap 
$(390, 0, 1, 2)$,
$(391, 0, 130, 128)$,
$(392, 0, 259, 389)$,
$(393, 0, 262, 131)$,\adfsplit
$(394, 0, 388, 260)$,
$(395, 0, 4, 11)$,
$(396, 0, 7, 386)$,
$(397, 0, 379, 383)$,\adfsplit
$(398, 0, 10, 23)$,
$(399, 0, 13, 380)$,
$(400, 0, 367, 377)$,
$(401, 0, 16, 35)$,\adfsplit
$(402, 0, 19, 374)$,
$(403, 0, 355, 371)$,
$(404, 0, 22, 5)$,
$(405, 0, 373, 368)$,\adfsplit
$(406, 0, 385, 17)$,
$(407, 0, 25, 53)$,
$(408, 0, 28, 365)$,
$(409, 0, 337, 362)$,\adfsplit
$(383, 322, 31, 170)$,
$(16, 114, 364, 163)$,
$(140, 183, 286, 358)$,
$(56, 256, 114, 70)$,\adfsplit
$(299, 376, 43, 130)$,
$(26, 112, 138, 190)$,
$(77, 341, 159, 111)$,
$(27, 256, 332, 210)$,\adfsplit
$(264, 284, 255, 23)$,
$(44, 20, 340, 47)$,
$(102, 43, 227, 382)$,
$(292, 331, 25, 144)$,\adfsplit
$(386, 187, 155, 293)$,
$(323, 78, 166, 290)$,
$(0, 6, 18, 327)$,
$(0, 8, 62, 151)$,\adfsplit
$(0, 36, 109, 217)$,
$(0, 67, 171, 272)$,
$(0, 37, 102, 316)$,
$(0, 47, 115, 179)$,\adfsplit
$(0, 41, 107, 277)$,
$(0, 31, 168, 223)$,
$(0, 56, 127, 290)$,
$(0, 40, 136, 228)$,\adfsplit
$(0, 79, 193, 273)$,
$(0, 21, 116, 237)$,
$(0, 50, 133, 299)$

\adfLgap \noindent by the mapping:
$x \mapsto x + 3 j \adfmod{390}$ for $x < 390$,
$x \mapsto x$ for $x \ge 390$,
$0 \le j < 130$
 for the first 20 blocks;
$x \mapsto x +  j \adfmod{390}$ for $x < 390$,
$x \mapsto x$ for $x \ge 390$,
$0 \le j < 390$
 for the last 27 blocks.
\ADFvfyParStart{(410, ((20, 130, ((390, 3), (20, 20))), (27, 390, ((390, 1), (20, 20)))), ((26, 15), (20, 1)))} 

\adfDgap
\noindent{\boldmath $ 26^{15} 23^{1} $}~
With the point set $Z_{413}$ partitioned into
 residue classes modulo $15$ for $\{0, 1, \dots, 389\}$, and
 $\{390, 391, \dots, 412\}$,
 the design is generated from

\adfLgap 
$(390, 0, 1, 2)$,
$(391, 0, 130, 128)$,
$(392, 0, 259, 389)$,
$(393, 0, 262, 131)$,\adfsplit
$(394, 0, 388, 260)$,
$(395, 0, 4, 11)$,
$(396, 0, 7, 386)$,
$(397, 0, 379, 383)$,\adfsplit
$(398, 0, 10, 23)$,
$(399, 0, 13, 380)$,
$(400, 0, 367, 377)$,
$(401, 0, 16, 35)$,\adfsplit
$(402, 0, 19, 374)$,
$(403, 0, 355, 371)$,
$(404, 0, 22, 5)$,
$(405, 0, 373, 368)$,\adfsplit
$(406, 0, 385, 17)$,
$(407, 0, 25, 53)$,
$(408, 0, 28, 365)$,
$(409, 0, 337, 362)$,\adfsplit
$(410, 0, 31, 65)$,
$(411, 0, 34, 359)$,
$(412, 0, 325, 356)$,
$(102, 254, 196, 338)$,\adfsplit
$(175, 252, 374, 91)$,
$(45, 367, 201, 347)$,
$(322, 230, 221, 45)$,
$(87, 348, 180, 265)$,\adfsplit
$(349, 323, 88, 85)$,
$(56, 360, 239, 144)$,
$(135, 231, 193, 18)$,
$(387, 168, 358, 308)$,\adfsplit
$(206, 347, 73, 320)$,
$(43, 0, 51, 319)$,
$(23, 148, 256, 337)$,
$(157, 23, 295, 221)$,\adfsplit
$(244, 175, 378, 312)$,
$(115, 185, 342, 26)$,
$(159, 88, 389, 327)$,
$(279, 173, 229, 320)$,\adfsplit
$(216, 44, 56, 103)$,
$(386, 102, 172, 305)$,
$(13, 224, 334, 278)$,
$(184, 86, 141, 14)$,\adfsplit
$(107, 192, 324, 25)$,
$(363, 118, 280, 342)$,
$(367, 89, 116, 193)$,
$(169, 75, 350, 301)$,\adfsplit
$(30, 69, 251, 36)$,
$(335, 212, 133, 151)$,
$(91, 142, 19, 43)$,
$(349, 163, 121, 54)$,\adfsplit
$(191, 228, 306, 254)$,
$(281, 188, 15, 142)$,
$(299, 17, 157, 346)$,
$(40, 237, 323, 225)$,\adfsplit
$(0, 3, 32, 151)$,
$(0, 8, 177, 279)$,
$(0, 36, 99, 317)$,
$(0, 96, 198, 369)$,\adfsplit
$(0, 101, 145, 181)$,
$(0, 37, 103, 335)$,
$(0, 121, 153, 231)$,
$(0, 193, 207, 247)$,\adfsplit
$(0, 49, 149, 303)$,
$(1, 7, 53, 197)$,
$(0, 38, 246, 287)$,
$(0, 57, 227, 286)$,\adfsplit
$(0, 18, 155, 308)$,
$(0, 14, 111, 178)$,
$(0, 42, 241, 274)$,
$(0, 40, 186, 333)$,\adfsplit
$(0, 48, 124, 242)$,
$(0, 20, 44, 156)$,
$(0, 9, 188, 326)$,
$(0, 80, 167, 206)$

\adfLgap \noindent by the mapping:
$x \mapsto x + 3 j \adfmod{390}$ for $x < 390$,
$x \mapsto x$ for $x \ge 390$,
$0 \le j < 130$
 for the first 23 blocks;
$x \mapsto x + 2 j \adfmod{390}$ for $x < 390$,
$x \mapsto x$ for $x \ge 390$,
$0 \le j < 195$
 for the last 53 blocks.
\ADFvfyParStart{(413, ((23, 130, ((390, 3), (23, 23))), (53, 195, ((390, 2), (23, 23)))), ((26, 15), (23, 1)))} 

\adfDgap
\noindent{\boldmath $ 26^{18} 20^{1} $}~
With the point set $Z_{488}$ partitioned into
 residue classes modulo $18$ for $\{0, 1, \dots, 467\}$, and
 $\{468, 469, \dots, 487\}$,
 the design is generated from

\adfLgap 
$(486, 0, 1, 2)$,
$(487, 0, 157, 314)$,
$(468, 196, 133, 365)$,
$(468, 117, 354, 134)$,\adfsplit
$(469, 202, 143, 234)$,
$(469, 278, 343, 231)$,
$(470, 212, 245, 357)$,
$(470, 378, 205, 406)$,\adfsplit
$(471, 341, 178, 198)$,
$(471, 235, 122, 51)$,
$(472, 208, 235, 24)$,
$(472, 177, 281, 236)$,\adfsplit
$(473, 301, 148, 461)$,
$(473, 54, 326, 237)$,
$(364, 347, 119, 44)$,
$(211, 70, 350, 189)$,\adfsplit
$(223, 255, 354, 142)$,
$(272, 420, 129, 248)$,
$(210, 181, 275, 15)$,
$(206, 56, 149, 267)$,\adfsplit
$(232, 81, 15, 367)$,
$(399, 139, 86, 1)$,
$(429, 443, 150, 13)$,
$(366, 92, 224, 52)$,\adfsplit
$(247, 47, 4, 230)$,
$(293, 105, 236, 192)$,
$(212, 183, 42, 208)$,
$(232, 192, 370, 33)$,\adfsplit
$(4, 153, 133, 233)$,
$(95, 353, 10, 368)$,
$(417, 151, 139, 174)$,
$(329, 465, 194, 337)$,\adfsplit
$(31, 408, 433, 424)$,
$(297, 78, 93, 339)$,
$(302, 223, 62, 309)$,
$(409, 127, 188, 238)$,\adfsplit
$(74, 240, 283, 166)$,
$(401, 454, 350, 451)$,
$(168, 333, 38, 433)$,
$(121, 126, 15, 197)$,\adfsplit
$(318, 55, 288, 167)$,
$(119, 5, 276, 112)$,
$(364, 16, 446, 162)$,
$(100, 30, 305, 365)$,\adfsplit
$(408, 463, 383, 345)$,
$(265, 181, 142, 220)$,
$(91, 45, 439, 78)$,
$(166, 101, 463, 444)$,\adfsplit
$(404, 187, 185, 435)$,
$(152, 248, 179, 369)$,
$(283, 12, 194, 79)$,
$(273, 441, 82, 76)$,\adfsplit
$(355, 145, 42, 279)$,
$(418, 192, 87, 131)$,
$(281, 108, 420, 457)$,
$(44, 87, 349, 405)$,\adfsplit
$(414, 300, 304, 436)$,
$(416, 126, 318, 71)$,
$(0, 5, 193, 346)$,
$(0, 11, 14, 188)$,\adfsplit
$(0, 7, 78, 246)$,
$(0, 17, 102, 343)$,
$(0, 87, 97, 350)$,
$(0, 9, 105, 319)$,\adfsplit
$(0, 6, 244, 426)$,
$(0, 69, 88, 344)$,
$(0, 50, 109, 167)$,
$(0, 86, 341, 417)$,\adfsplit
$(0, 49, 129, 262)$,
$(0, 155, 189, 376)$,
$(0, 47, 203, 370)$,
$(0, 147, 233, 321)$,\adfsplit
$(0, 172, 227, 239)$,
$(0, 21, 207, 408)$,
$(0, 118, 149, 179)$,
$(0, 31, 130, 437)$,\adfsplit
$(0, 29, 52, 447)$,
$(0, 80, 159, 455)$,
$(3, 11, 35, 244)$,
$(3, 16, 87, 382)$,\adfsplit
$(3, 51, 100, 376)$

\adfLgap \noindent by the mapping:
$x \mapsto x \oplus (3 j)$ for $x < 468$,
$x \mapsto (x + 6 j \adfmod{18}) + 468$ for $468 \le x < 486$,
$x \mapsto x$ for $x \ge 486$,
$0 \le j < 156$
 for the first two blocks;
$x \mapsto x \oplus j \oplus j$ for $x < 468$,
$x \mapsto (x + 6 j \adfmod{18}) + 468$ for $468 \le x < 486$,
$x \mapsto x$ for $x \ge 486$,
$0 \le j < 234$
 for the last 79 blocks.
\ADFvfyParStart{(488, ((2, 156, ((468, 3, (156, 3)), (18, 6), (2, 2))), (79, 234, ((468, 2, (156, 3)), (18, 6), (2, 2)))), ((26, 18), (20, 1)))} 

\adfDgap
\noindent{\boldmath $ 26^{18} 23^{1} $}~
With the point set $Z_{491}$ partitioned into
 residue classes modulo $18$ for $\{0, 1, \dots, 467\}$, and
 $\{468, 469, \dots, 490\}$,
 the design is generated from

\adfLgap 
$(489, 130, 191, 42)$,
$(490, 221, 112, 249)$,
$(468, 375, 355, 9)$,
$(468, 131, 100, 430)$,\adfsplit
$(468, 335, 98, 156)$,
$(469, 74, 288, 312)$,
$(469, 138, 241, 91)$,
$(469, 140, 350, 409)$,\adfsplit
$(470, 343, 96, 405)$,
$(470, 128, 138, 4)$,
$(470, 95, 152, 115)$,
$(471, 155, 217, 274)$,\adfsplit
$(471, 331, 417, 9)$,
$(471, 294, 248, 107)$,
$(472, 291, 177, 436)$,
$(472, 59, 47, 423)$,\adfsplit
$(472, 430, 386, 343)$,
$(473, 73, 401, 161)$,
$(473, 160, 459, 132)$,
$(473, 31, 389, 246)$,\adfsplit
$(474, 68, 394, 373)$,
$(474, 281, 314, 444)$,
$(474, 307, 252, 195)$,
$(370, 158, 430, 198)$,\adfsplit
$(48, 399, 420, 397)$,
$(256, 355, 165, 16)$,
$(261, 302, 109, 301)$,
$(416, 197, 262, 126)$,\adfsplit
$(82, 261, 77, 276)$,
$(43, 123, 232, 437)$,
$(13, 28, 332, 329)$,
$(149, 446, 163, 160)$,\adfsplit
$(236, 383, 284, 373)$,
$(39, 466, 434, 68)$,
$(420, 209, 92, 5)$,
$(362, 153, 386, 226)$,\adfsplit
$(361, 269, 427, 414)$,
$(171, 293, 137, 0)$,
$(383, 29, 13, 20)$,
$(331, 89, 198, 67)$,\adfsplit
$(70, 359, 442, 104)$,
$(301, 454, 326, 421)$,
$(172, 9, 287, 25)$,
$(159, 289, 279, 283)$,\adfsplit
$(351, 332, 33, 162)$,
$(319, 95, 290, 328)$,
$(184, 443, 249, 381)$,
$(197, 151, 103, 210)$,\adfsplit
$(38, 33, 40, 306)$,
$(55, 346, 293, 212)$,
$(6, 334, 25, 467)$,
$(158, 339, 301, 108)$,\adfsplit
$(45, 90, 251, 286)$,
$(117, 282, 314, 381)$,
$(255, 223, 446, 447)$,
$(392, 47, 217, 117)$,\adfsplit
$(351, 187, 429, 373)$,
$(329, 379, 34, 264)$,
$(189, 257, 154, 328)$,
$(346, 417, 97, 279)$,\adfsplit
$(249, 459, 37, 282)$,
$(254, 228, 191, 52)$,
$(312, 61, 209, 154)$,
$(0, 462, 326, 240)$,\adfsplit
$(119, 115, 40, 456)$,
$(288, 167, 413, 368)$,
$(313, 142, 99, 255)$,
$(27, 235, 52, 121)$,\adfsplit
$(255, 360, 397, 140)$,
$(337, 286, 453, 437)$,
$(67, 106, 251, 79)$,
$(116, 192, 200, 222)$,\adfsplit
$(136, 384, 426, 453)$,
$(268, 80, 429, 226)$,
$(186, 267, 206, 262)$,
$(89, 439, 375, 376)$,\adfsplit
$(142, 277, 161, 41)$,
$(214, 169, 335, 387)$,
$(193, 82, 52, 237)$,
$(121, 289, 90, 16)$,\adfsplit
$(145, 452, 156, 173)$,
$(182, 167, 29, 421)$,
$(394, 213, 308, 17)$,
$(111, 19, 265, 36)$,\adfsplit
$(436, 270, 193, 185)$,
$(82, 160, 437, 210)$,
$(367, 265, 435, 5)$,
$(0, 2, 49, 385)$,\adfsplit
$(0, 14, 211, 442)$,
$(0, 34, 118, 235)$,
$(0, 61, 267, 373)$,
$(0, 35, 160, 325)$,\adfsplit
$(0, 74, 205, 460)$,
$(0, 97, 121, 371)$,
$(0, 70, 280, 381)$,
$(1, 14, 262, 314)$,\adfsplit
$(1, 59, 65, 82)$,
$(0, 277, 374, 404)$,
$(1, 113, 191, 242)$,
$(0, 79, 128, 183)$,\adfsplit
$(0, 221, 263, 343)$,
$(0, 85, 174, 377)$,
$(0, 110, 187, 260)$,
$(0, 47, 186, 268)$,\adfsplit
$(0, 224, 335, 439)$,
$(0, 84, 173, 454)$,
$(1, 266, 362, 401)$,
$(0, 9, 184, 251)$,\adfsplit
$(0, 158, 317, 358)$,
$(0, 98, 164, 321)$,
$(0, 23, 66, 315)$,
$(0, 3, 188, 215)$,\adfsplit
$(0, 11, 93, 261)$,
$(0, 77, 146, 320)$,
$(0, 134, 155, 437)$,
$(0, 123, 389, 464)$,\adfsplit
$(0, 59, 152, 344)$,
$(0, 51, 107, 167)$,
$(0, 48, 111, 350)$,
$(0, 39, 83, 419)$,\adfsplit
$(0, 12, 255, 368)$,
$(0, 101, 308, 333)$

\adfLgap \noindent by the mapping:
$x \mapsto x + 3 j \adfmod{468}$ for $x < 468$,
$x \mapsto (x - 468 + 7 j \adfmod{21}) + 468$ for $468 \le x < 489$,
$x \mapsto x$ for $x \ge 489$,
$0 \le j < 156$.
\ADFvfyParStart{(491, ((122, 156, ((468, 3), (21, 7), (2, 2)))), ((26, 18), (23, 1)))} 

\section{4-GDDs for the proof of Lemma \ref{lem:4-GDD 28^u m^1}}
\label{app:4-GDD 28^u m^1}
\adfhide{
$ 28^9 19^1 $,
$ 28^9 25^1 $ and
$ 28^9 31^1 $.
}

\adfDgap
\noindent{\boldmath $ 28^{9} 19^{1} $}~
With the point set $Z_{271}$ partitioned into
 residue classes modulo $9$ for $\{0, 1, \dots, 251\}$, and
 $\{252, 253, \dots, 270\}$,
 the design is generated from

\adfLgap 
$(252, 230, 231, 66)$,
$(252, 189, 193, 38)$,
$(252, 8, 118, 132)$,
$(252, 73, 201, 155)$,\adfsplit
$(252, 52, 112, 0)$,
$(252, 205, 215, 239)$,
$(253, 108, 194, 251)$,
$(253, 147, 88, 247)$,\adfsplit
$(253, 189, 80, 195)$,
$(253, 238, 172, 168)$,
$(253, 151, 245, 217)$,
$(253, 221, 74, 120)$,\adfsplit
$(44, 91, 171, 220)$,
$(190, 193, 83, 215)$,
$(6, 73, 120, 88)$,
$(219, 146, 186, 143)$,\adfsplit
$(116, 207, 132, 219)$,
$(63, 74, 196, 80)$,
$(176, 98, 119, 250)$,
$(28, 7, 78, 146)$,\adfsplit
$(55, 74, 203, 39)$,
$(237, 172, 81, 17)$,
$(9, 241, 211, 212)$,
$(237, 224, 25, 77)$,\adfsplit
$(197, 13, 61, 195)$,
$(35, 248, 4, 16)$,
$(47, 210, 248, 220)$,
$(203, 217, 159, 18)$,\adfsplit
$(0, 2, 179, 230)$,
$(0, 5, 146, 172)$,
$(0, 83, 139, 213)$,
$(0, 30, 71, 150)$,\adfsplit
$(0, 7, 48, 92)$,
$(0, 43, 104, 197)$,
$(0, 15, 98, 229)$,
$(0, 58, 183, 245)$,\adfsplit
$(0, 64, 133, 239)$,
$(0, 137, 215, 223)$,
$(0, 34, 95, 190)$,
$(0, 35, 94, 149)$,\adfsplit
$(0, 23, 56, 163)$,
$(0, 17, 77, 227)$,
$(0, 11, 37, 42)$,
$(270, 0, 84, 168)$,\adfsplit
$(270, 1, 85, 169)$

\adfLgap \noindent by the mapping:
$x \mapsto x + 2 j \adfmod{252}$ for $x < 252$,
$x \mapsto (x + 2 j \adfmod{18}) + 252$ for $252 \le x < 270$,
$270 \mapsto 270$,
$0 \le j < 126$
 for the first 43 blocks,
$0 \le j < 42$
 for the last two blocks.
\ADFvfyParStart{(271, ((43, 126, ((252, 2), (18, 2), (1, 1))), (2, 42, ((252, 2), (18, 2), (1, 1)))), ((28, 9), (19, 1)))} 

\adfDgap
\noindent{\boldmath $ 28^{9} 25^{1} $}~
With the point set $Z_{277}$ partitioned into
 residue classes modulo $9$ for $\{0, 1, \dots, 251\}$, and
 $\{252, 253, \dots, 276\}$,
 the design is generated from

\adfLgap 
$(270, 168, 147, 77)$,
$(270, 212, 157, 232)$,
$(271, 232, 153, 175)$,
$(271, 155, 170, 72)$,\adfsplit
$(252, 133, 48, 193)$,
$(252, 19, 237, 137)$,
$(252, 64, 113, 180)$,
$(252, 4, 26, 123)$,\adfsplit
$(252, 63, 20, 17)$,
$(252, 34, 176, 204)$,
$(253, 196, 73, 75)$,
$(253, 83, 107, 99)$,\adfsplit
$(253, 94, 20, 87)$,
$(253, 42, 248, 121)$,
$(253, 151, 126, 188)$,
$(253, 48, 208, 203)$,\adfsplit
$(75, 22, 25, 74)$,
$(199, 231, 57, 162)$,
$(161, 81, 33, 110)$,
$(62, 176, 31, 1)$,\adfsplit
$(157, 168, 215, 158)$,
$(128, 140, 161, 25)$,
$(161, 230, 174, 232)$,
$(82, 66, 175, 18)$,\adfsplit
$(167, 73, 171, 8)$,
$(197, 33, 184, 59)$,
$(59, 48, 135, 222)$,
$(177, 8, 103, 154)$,\adfsplit
$(118, 150, 238, 77)$,
$(79, 161, 132, 126)$,
$(0, 59, 73, 111)$,
$(1, 13, 41, 197)$,\adfsplit
$(0, 55, 61, 167)$,
$(0, 31, 123, 133)$,
$(0, 39, 187, 229)$,
$(1, 21, 65, 187)$,\adfsplit
$(0, 7, 139, 156)$,
$(0, 4, 70, 213)$,
$(0, 19, 44, 130)$,
$(0, 5, 24, 38)$,\adfsplit
$(0, 17, 76, 118)$,
$(0, 30, 71, 202)$,
$(0, 15, 152, 192)$,
$(0, 8, 113, 148)$,\adfsplit
$(0, 26, 94, 128)$,
$(276, 0, 84, 168)$,
$(276, 1, 85, 169)$

\adfLgap \noindent by the mapping:
$x \mapsto x + 2 j \adfmod{252}$ for $x < 252$,
$x \mapsto (x + 2 j \adfmod{18}) + 252$ for $252 \le x < 270$,
$x \mapsto (x + 2 j \adfmod{6}) + 270$ for $270 \le x < 276$,
$276 \mapsto 276$,
$0 \le j < 126$
 for the first 45 blocks,
$0 \le j < 42$
 for the last two blocks.
\ADFvfyParStart{(277, ((45, 126, ((252, 2), (18, 2), (6, 2), (1, 1))), (2, 42, ((252, 2), (18, 2), (6, 2), (1, 1)))), ((28, 9), (25, 1)))} 

\adfDgap
\noindent{\boldmath $ 28^{9} 31^{1} $}~
With the point set $Z_{283}$ partitioned into
 residue classes modulo $9$ for $\{0, 1, \dots, 251\}$, and
 $\{252, 253, \dots, 282\}$,
 the design is generated from

\adfLgap 
$(279, 105, 71, 130)$,
$(279, 62, 216, 139)$,
$(252, 55, 90, 241)$,
$(252, 215, 103, 245)$,\adfsplit
$(252, 24, 250, 147)$,
$(252, 38, 33, 242)$,
$(252, 208, 167, 189)$,
$(252, 104, 156, 112)$,\adfsplit
$(253, 1, 103, 77)$,
$(253, 237, 83, 79)$,
$(253, 22, 32, 16)$,
$(253, 118, 80, 81)$,\adfsplit
$(253, 54, 132, 12)$,
$(253, 159, 92, 161)$,
$(254, 124, 139, 177)$,
$(254, 194, 181, 76)$,\adfsplit
$(254, 179, 150, 135)$,
$(254, 80, 192, 209)$,
$(254, 131, 25, 180)$,
$(254, 218, 190, 129)$,\adfsplit
$(111, 176, 76, 52)$,
$(145, 77, 250, 156)$,
$(69, 81, 127, 251)$,
$(9, 19, 130, 3)$,\adfsplit
$(94, 154, 203, 227)$,
$(247, 207, 182, 240)$,
$(133, 130, 48, 50)$,
$(45, 77, 163, 3)$,\adfsplit
$(226, 158, 65, 197)$,
$(199, 44, 85, 147)$,
$(58, 24, 244, 119)$,
$(6, 26, 128, 57)$,\adfsplit
$(0, 5, 13, 161)$,
$(0, 71, 119, 235)$,
$(0, 22, 55, 115)$,
$(0, 57, 104, 179)$,\adfsplit
$(0, 21, 101, 121)$,
$(0, 12, 177, 178)$,
$(0, 30, 143, 199)$,
$(0, 46, 182, 229)$,\adfsplit
$(0, 40, 96, 185)$,
$(0, 14, 190, 209)$,
$(0, 167, 231, 245)$,
$(0, 110, 221, 249)$,\adfsplit
$(0, 11, 106, 213)$,
$(0, 4, 164, 201)$,
$(0, 23, 138, 188)$,
$(282, 0, 84, 168)$,\adfsplit
$(282, 1, 85, 169)$

\adfLgap \noindent by the mapping:
$x \mapsto x + 2 j \adfmod{252}$ for $x < 252$,
$x \mapsto (x - 252 + 3 j \adfmod{27}) + 252$ for $252 \le x < 279$,
$x \mapsto (x +  j \adfmod{3}) + 279$ for $279 \le x < 282$,
$282 \mapsto 282$,
$0 \le j < 126$
 for the first 47 blocks,
$0 \le j < 42$
 for the last two blocks.
\ADFvfyParStart{(283, ((47, 126, ((252, 2), (27, 3), (3, 1), (1, 1))), (2, 42, ((252, 2), (27, 3), (3, 1), (1, 1)))), ((28, 9), (31, 1)))} 

\section{4-GDDs for the proof of Lemma \ref{lem:4-GDD 34^u m^1}}
\label{app:4-GDD 34^u m^1}
\adfhide{
$ 34^6 7^1 $,
$ 34^6 10^1 $,
$ 34^6 13^1 $,
$ 34^6 16^1 $,
$ 34^6 19^1 $,
$ 34^6 22^1 $,
$ 34^6 73^1 $,
$ 34^6 76^1 $,
$ 34^6 79^1 $,
$ 34^6 82^1 $,
$ 34^9 19^1 $,
$ 34^{12} 25^1 $,
$ 34^{12} 28^1 $,
$ 34^{12} 31^1 $ and
$ 34^{15} 31^1 $.
}

\adfDgap
\noindent{\boldmath $ 34^{6} 7^{1} $}~
With the point set $Z_{211}$ partitioned into
 residue classes modulo $6$ for $\{0, 1, \dots, 203\}$, and
 $\{204, 205, \dots, 210\}$,
 the design is generated from

\adfLgap 
$(204, 121, 135, 58)$,
$(204, 2, 95, 132)$,
$(180, 136, 23, 81)$,
$(135, 191, 194, 97)$,\adfsplit
$(172, 71, 150, 157)$,
$(76, 161, 192, 159)$,
$(0, 1, 5, 196)$,
$(0, 10, 21, 164)$,\adfsplit
$(0, 32, 73, 155)$,
$(0, 23, 62, 129)$,
$(0, 28, 71, 128)$,
$(0, 25, 52, 117)$,\adfsplit
$(0, 20, 46, 115)$,
$(0, 17, 51, 70)$,
$(0, 16, 45, 80)$,
$(210, 0, 68, 136)$

\adfLgap \noindent by the mapping:
$x \mapsto x +  j \adfmod{204}$ for $x < 204$,
$x \mapsto (x +  j \adfmod{6}) + 204$ for $204 \le x < 210$,
$210 \mapsto 210$,
$0 \le j < 204$
 for the first 15 blocks,
$0 \le j < 68$
 for the last block.
\ADFvfyParStart{(211, ((15, 204, ((204, 1), (6, 1), (1, 1))), (1, 68, ((204, 1), (6, 1), (1, 1)))), ((34, 6), (7, 1)))} 

\adfDgap
\noindent{\boldmath $ 34^{6} 10^{1} $}~
With the point set $Z_{214}$ partitioned into
 residue classes modulo $5$ for $\{0, 1, \dots, 169\}$,
 $\{170, 171, \dots, 203\}$, and
 $\{204, 205, \dots, 213\}$,
 the design is generated from

\adfLgap 
$(206, 197, 45, 102)$,
$(206, 172, 130, 78)$,
$(181, 166, 14, 60)$,
$(175, 87, 119, 163)$,\adfsplit
$(192, 151, 53, 167)$,
$(198, 87, 39, 128)$,
$(197, 33, 91, 37)$,
$(174, 140, 164, 127)$,\adfsplit
$(174, 135, 61, 134)$,
$(203, 15, 106, 53)$,
$(0, 2, 151, 202)$,
$(0, 22, 88, 191)$,\adfsplit
$(0, 23, 49, 142)$,
$(0, 7, 34, 43)$,
$(0, 3, 17, 86)$,
$(0, 12, 123, 205)$,\adfsplit
$(0, 6, 39, 68)$,
$(0, 8, 107, 207)$,
$(0, 11, 42, 103)$

\adfLgap \noindent by the mapping:
$x \mapsto x +  j \adfmod{170}$ for $x < 170$,
$x \mapsto (x +  j \adfmod{34}) + 170$ for $170 \le x < 204$,
$x \mapsto (x - 204 +  j \adfmod{10}) + 204$ for $x \ge 204$,
$0 \le j < 170$.
\ADFvfyParStart{(214, ((19, 170, ((170, 1), (34, 1), (10, 1)))), ((34, 5), (34, 1), (10, 1)))} 

\adfDgap
\noindent{\boldmath $ 34^{6} 13^{1} $}~
With the point set $Z_{217}$ partitioned into
 residue classes modulo $6$ for $\{0, 1, \dots, 203\}$, and
 $\{204, 205, \dots, 216\}$,
 the design is generated from

\adfLgap 
$(204, 168, 5, 85)$,
$(204, 40, 201, 50)$,
$(205, 95, 187, 112)$,
$(205, 174, 129, 134)$,\adfsplit
$(133, 164, 34, 132)$,
$(68, 89, 66, 117)$,
$(90, 28, 25, 104)$,
$(0, 4, 13, 117)$,\adfsplit
$(0, 7, 77, 93)$,
$(0, 19, 57, 101)$,
$(0, 26, 81, 133)$,
$(0, 20, 59, 109)$,\adfsplit
$(0, 11, 33, 67)$,
$(0, 15, 61, 146)$,
$(0, 8, 35, 175)$,
$(0, 25, 88, 135)$,\adfsplit
$(216, 0, 68, 136)$

\adfLgap \noindent by the mapping:
$x \mapsto x +  j \adfmod{204}$ for $x < 204$,
$x \mapsto (x + 2 j \adfmod{12}) + 204$ for $204 \le x < 216$,
$216 \mapsto 216$,
$0 \le j < 204$
 for the first 16 blocks,
$0 \le j < 68$
 for the last block.
\ADFvfyParStart{(217, ((16, 204, ((204, 1), (12, 2), (1, 1))), (1, 68, ((204, 1), (12, 2), (1, 1)))), ((34, 6), (13, 1)))} 

\adfDgap
\noindent{\boldmath $ 34^{6} 16^{1} $}~
With the point set $Z_{220}$ partitioned into
 residue classes modulo $6$ for $\{0, 1, \dots, 203\}$, and
 $\{204, 205, \dots, 219\}$,
 the design is generated from

\adfLgap 
$(204, 32, 109, 153)$,
$(204, 138, 196, 179)$,
$(205, 33, 108, 121)$,
$(205, 28, 35, 128)$,\adfsplit
$(206, 110, 83, 162)$,
$(206, 97, 34, 57)$,
$(207, 156, 142, 39)$,
$(207, 158, 11, 43)$,\adfsplit
$(208, 98, 28, 12)$,
$(208, 117, 191, 25)$,
$(201, 115, 6, 191)$,
$(140, 42, 81, 76)$,\adfsplit
$(147, 77, 176, 10)$,
$(141, 187, 68, 94)$,
$(164, 133, 132, 131)$,
$(80, 34, 129, 151)$,\adfsplit
$(174, 10, 170, 123)$,
$(124, 119, 12, 67)$,
$(0, 2, 17, 124)$,
$(0, 3, 142, 151)$,\adfsplit
$(0, 11, 31, 76)$,
$(0, 20, 81, 181)$,
$(0, 8, 131, 135)$,
$(0, 10, 37, 179)$,\adfsplit
$(0, 29, 43, 50)$,
$(0, 75, 103, 116)$,
$(0, 45, 56, 155)$,
$(0, 25, 51, 176)$,\adfsplit
$(0, 21, 110, 143)$,
$(0, 19, 22, 189)$,
$(0, 35, 74, 133)$,
$(0, 69, 85, 149)$,\adfsplit
$(0, 83, 91, 141)$,
$(219, 0, 68, 136)$,
$(219, 1, 69, 137)$

\adfLgap \noindent by the mapping:
$x \mapsto x + 2 j \adfmod{204}$ for $x < 204$,
$x \mapsto (x - 204 + 5 j \adfmod{15}) + 204$ for $204 \le x < 219$,
$219 \mapsto 219$,
$0 \le j < 102$
 for the first 33 blocks,
$0 \le j < 34$
 for the last two blocks.
\ADFvfyParStart{(220, ((33, 102, ((204, 2), (15, 5), (1, 1))), (2, 34, ((204, 2), (15, 5), (1, 1)))), ((34, 6), (16, 1)))} 

\adfDgap
\noindent{\boldmath $ 34^{6} 19^{1} $}~
With the point set $Z_{223}$ partitioned into
 residue classes modulo $6$ for $\{0, 1, \dots, 203\}$, and
 $\{204, 205, \dots, 222\}$,
 the design is generated from

\adfLgap 
$(204, 180, 79, 194)$,
$(204, 76, 119, 147)$,
$(205, 180, 57, 119)$,
$(205, 134, 52, 37)$,\adfsplit
$(206, 200, 87, 83)$,
$(206, 174, 79, 88)$,
$(112, 143, 133, 110)$,
$(39, 28, 144, 179)$,\adfsplit
$(141, 191, 190, 146)$,
$(0, 3, 32, 40)$,
$(0, 34, 75, 134)$,
$(0, 25, 77, 146)$,\adfsplit
$(0, 22, 73, 147)$,
$(0, 16, 63, 128)$,
$(0, 7, 26, 184)$,
$(0, 13, 80, 119)$,\adfsplit
$(0, 17, 55, 111)$,
$(222, 0, 68, 136)$

\adfLgap \noindent by the mapping:
$x \mapsto x +  j \adfmod{204}$ for $x < 204$,
$x \mapsto (x - 204 + 3 j \adfmod{18}) + 204$ for $204 \le x < 222$,
$222 \mapsto 222$,
$0 \le j < 204$
 for the first 17 blocks,
$0 \le j < 68$
 for the last block.
\ADFvfyParStart{(223, ((17, 204, ((204, 1), (18, 3), (1, 1))), (1, 68, ((204, 1), (18, 3), (1, 1)))), ((34, 6), (19, 1)))} 

\adfDgap
\noindent{\boldmath $ 34^{6} 22^{1} $}~
With the point set $Z_{226}$ partitioned into
 residue classes modulo $6$ for $\{0, 1, \dots, 203\}$, and
 $\{204, 205, \dots, 225\}$,
 the design is generated from

\adfLgap 
$(204, 122, 93, 113)$,
$(204, 132, 67, 160)$,
$(205, 132, 44, 197)$,
$(205, 22, 3, 67)$,\adfsplit
$(206, 139, 202, 90)$,
$(206, 147, 200, 53)$,
$(207, 36, 137, 94)$,
$(207, 45, 169, 38)$,\adfsplit
$(208, 60, 170, 19)$,
$(208, 136, 99, 131)$,
$(209, 146, 97, 179)$,
$(209, 172, 75, 120)$,\adfsplit
$(210, 86, 168, 130)$,
$(210, 107, 117, 151)$,
$(49, 8, 71, 88)$,
$(102, 95, 141, 164)$,\adfsplit
$(135, 10, 8, 29)$,
$(110, 25, 40, 113)$,
$(0, 1, 5, 100)$,
$(0, 4, 15, 149)$,\adfsplit
$(0, 9, 17, 130)$,
$(0, 8, 178, 191)$,
$(0, 16, 56, 85)$,
$(0, 35, 75, 158)$,\adfsplit
$(0, 20, 47, 133)$,
$(0, 77, 79, 129)$,
$(0, 14, 64, 161)$,
$(0, 103, 177, 203)$,\adfsplit
$(0, 10, 61, 99)$,
$(0, 71, 87, 118)$,
$(0, 25, 106, 137)$,
$(0, 22, 59, 115)$,\adfsplit
$(0, 67, 95, 128)$,
$(0, 23, 169, 172)$,
$(0, 117, 179, 193)$,
$(225, 0, 68, 136)$,\adfsplit
$(225, 1, 69, 137)$

\adfLgap \noindent by the mapping:
$x \mapsto x + 2 j \adfmod{204}$ for $x < 204$,
$x \mapsto (x - 204 + 7 j \adfmod{21}) + 204$ for $204 \le x < 225$,
$225 \mapsto 225$,
$0 \le j < 102$
 for the first 35 blocks,
$0 \le j < 34$
 for the last two blocks.
\ADFvfyParStart{(226, ((35, 102, ((204, 2), (21, 7), (1, 1))), (2, 34, ((204, 2), (21, 7), (1, 1)))), ((34, 6), (22, 1)))} 

\adfDgap
\noindent{\boldmath $ 34^{6} 73^{1} $}~
With the point set $Z_{277}$ partitioned into
 residue classes modulo $6$ for $\{0, 1, \dots, 203\}$, and
 $\{204, 205, \dots, 276\}$,
 the design is generated from

\adfLgap 
$(255, 22, 69, 119)$,
$(256, 38, 177, 22)$,
$(257, 93, 31, 59)$,
$(258, 163, 167, 111)$,\adfsplit
$(259, 80, 57, 151)$,
$(260, 194, 124, 165)$,
$(261, 133, 131, 144)$,
$(204, 196, 171, 116)$,\adfsplit
$(204, 121, 188, 131)$,
$(204, 37, 126, 95)$,
$(204, 198, 135, 179)$,
$(204, 5, 90, 130)$,\adfsplit
$(204, 66, 103, 21)$,
$(204, 104, 40, 191)$,
$(204, 124, 81, 55)$,
$(204, 201, 125, 98)$,\adfsplit
$(0, 3, 35, 86)$,
$(0, 5, 38, 243)$,
$(0, 20, 95, 221)$,
$(0, 8, 17, 247)$,\adfsplit
$(0, 7, 111, 213)$,
$(0, 1, 106, 246)$,
$(0, 22, 61, 219)$,
$(0, 21, 113, 235)$,\adfsplit
$(0, 46, 123, 244)$,
$(0, 14, 73, 88)$,
$(276, 0, 68, 136)$

\adfLgap \noindent by the mapping:
$x \mapsto x +  j \adfmod{204}$ for $x < 204$,
$x \mapsto (x +  j \adfmod{51}) + 204$ for $204 \le x < 255$,
$x \mapsto (x - 255 + 7 j \adfmod{21}) + 255$ for $255 \le x < 276$,
$276 \mapsto 276$,
$0 \le j < 204$
 for the first 26 blocks,
$0 \le j < 68$
 for the last block.
\ADFvfyParStart{(277, ((26, 204, ((204, 1), (51, 1), (21, 7), (1, 1))), (1, 68, ((204, 1), (51, 1), (21, 7), (1, 1)))), ((34, 6), (73, 1)))} 

\adfDgap
\noindent{\boldmath $ 34^{6} 76^{1} $}~
With the point set $Z_{280}$ partitioned into
 residue classes modulo $6$ for $\{0, 1, \dots, 203\}$, and
 $\{204, 205, \dots, 279\}$,
 the design is generated from

\adfLgap 
$(255, 112, 5, 81)$,
$(255, 139, 104, 166)$,
$(255, 133, 87, 24)$,
$(255, 95, 86, 102)$,\adfsplit
$(256, 183, 170, 79)$,
$(256, 29, 96, 25)$,
$(256, 202, 44, 18)$,
$(256, 93, 143, 88)$,\adfsplit
$(257, 23, 88, 86)$,
$(257, 5, 175, 6)$,
$(257, 105, 94, 24)$,
$(257, 97, 8, 39)$,\adfsplit
$(258, 83, 136, 55)$,
$(258, 60, 45, 94)$,
$(258, 14, 181, 174)$,
$(258, 176, 147, 65)$,\adfsplit
$(204, 156, 55, 184)$,
$(204, 175, 149, 84)$,
$(204, 163, 161, 182)$,
$(204, 167, 33, 46)$,\adfsplit
$(204, 132, 170, 115)$,
$(204, 193, 190, 0)$,
$(204, 140, 67, 165)$,
$(204, 85, 23, 28)$,\adfsplit
$(204, 43, 4, 183)$,
$(204, 164, 1, 160)$,
$(204, 26, 133, 180)$,
$(204, 5, 142, 21)$,\adfsplit
$(204, 36, 83, 27)$,
$(204, 143, 124, 151)$,
$(204, 153, 104, 6)$,
$(0, 1, 100, 117)$,\adfsplit
$(0, 8, 61, 233)$,
$(0, 10, 32, 225)$,
$(0, 15, 124, 145)$,
$(0, 29, 73, 223)$,\adfsplit
$(0, 79, 165, 217)$,
$(0, 82, 201, 248)$,
$(0, 43, 88, 238)$,
$(0, 85, 118, 239)$,\adfsplit
$(0, 143, 163, 212)$,
$(0, 33, 71, 243)$,
$(0, 101, 111, 251)$,
$(0, 51, 140, 246)$,\adfsplit
$(1, 15, 95, 212)$,
$(0, 23, 74, 247)$,
$(0, 59, 128, 219)$,
$(0, 56, 181, 250)$,\adfsplit
$(0, 40, 127, 218)$,
$(1, 41, 153, 220)$,
$(0, 37, 69, 112)$,
$(0, 77, 99, 205)$,\adfsplit
$(0, 52, 110, 234)$,
$(279, 0, 68, 136)$,
$(279, 1, 69, 137)$

\adfLgap \noindent by the mapping:
$x \mapsto x + 2 j \adfmod{204}$ for $x < 204$,
$x \mapsto (x +  j \adfmod{51}) + 204$ for $204 \le x < 255$,
$x \mapsto (x - 255 + 4 j \adfmod{24}) + 255$ for $255 \le x < 279$,
$279 \mapsto 279$,
$0 \le j < 102$
 for the first 53 blocks,
$0 \le j < 34$
 for the last two blocks.
\ADFvfyParStart{(280, ((53, 102, ((204, 2), (51, 1), (24, 4), (1, 1))), (2, 34, ((204, 2), (51, 1), (24, 4), (1, 1)))), ((34, 6), (76, 1)))} 

\adfDgap
\noindent{\boldmath $ 34^{6} 79^{1} $}~
With the point set $Z_{283}$ partitioned into
 residue classes modulo $6$ for $\{0, 1, \dots, 203\}$, and
 $\{204, 205, \dots, 282\}$,
 the design is generated from

\adfLgap 
$(255, 146, 37, 42)$,
$(256, 8, 4, 75)$,
$(257, 112, 68, 159)$,
$(258, 195, 97, 131)$,\adfsplit
$(259, 133, 173, 90)$,
$(260, 64, 140, 18)$,
$(261, 97, 60, 68)$,
$(262, 110, 109, 126)$,\adfsplit
$(263, 197, 46, 96)$,
$(204, 64, 33, 96)$,
$(204, 85, 70, 195)$,
$(204, 122, 41, 148)$,\adfsplit
$(204, 48, 89, 177)$,
$(204, 95, 56, 108)$,
$(204, 4, 15, 60)$,
$(204, 87, 114, 2)$,\adfsplit
$(204, 193, 76, 78)$,
$(0, 7, 58, 80)$,
$(0, 25, 130, 233)$,
$(0, 3, 62, 223)$,\adfsplit
$(0, 19, 166, 237)$,
$(0, 23, 134, 226)$,
$(0, 28, 61, 244)$,
$(0, 10, 149, 248)$,\adfsplit
$(0, 9, 86, 241)$,
$(0, 14, 35, 239)$,
$(0, 20, 69, 247)$,
$(282, 0, 68, 136)$

\adfLgap \noindent by the mapping:
$x \mapsto x +  j \adfmod{204}$ for $x < 204$,
$x \mapsto (x +  j \adfmod{51}) + 204$ for $204 \le x < 255$,
$x \mapsto (x - 255 + 9 j \adfmod{27}) + 255$ for $255 \le x < 282$,
$282 \mapsto 282$,
$0 \le j < 204$
 for the first 27 blocks,
$0 \le j < 68$
 for the last block.
\ADFvfyParStart{(283, ((27, 204, ((204, 1), (51, 1), (27, 9), (1, 1))), (1, 68, ((204, 1), (51, 1), (27, 9), (1, 1)))), ((34, 6), (79, 1)))} 

\adfDgap
\noindent{\boldmath $ 34^{6} 82^{1} $}~
With the point set $Z_{286}$ partitioned into
 residue classes modulo $6$ for $\{0, 1, \dots, 203\}$, and
 $\{204, 205, \dots, 285\}$,
 the design is generated from

\adfLgap 
$(255, 38, 41, 22)$,
$(255, 128, 90, 33)$,
$(255, 103, 71, 120)$,
$(255, 87, 172, 145)$,\adfsplit
$(256, 15, 36, 106)$,
$(256, 76, 150, 119)$,
$(256, 92, 125, 139)$,
$(256, 158, 85, 57)$,\adfsplit
$(257, 159, 84, 133)$,
$(257, 199, 8, 153)$,
$(257, 40, 18, 161)$,
$(257, 34, 158, 23)$,\adfsplit
$(258, 110, 156, 201)$,
$(258, 188, 100, 67)$,
$(258, 41, 133, 94)$,
$(258, 63, 59, 42)$,\adfsplit
$(259, 200, 37, 161)$,
$(259, 148, 134, 31)$,
$(259, 183, 126, 22)$,
$(259, 107, 0, 9)$,\adfsplit
$(204, 138, 170, 161)$,
$(204, 165, 168, 35)$,
$(204, 33, 80, 144)$,
$(204, 129, 173, 91)$,\adfsplit
$(204, 12, 123, 188)$,
$(204, 118, 200, 133)$,
$(204, 166, 57, 162)$,
$(204, 157, 95, 90)$,\adfsplit
$(204, 191, 94, 74)$,
$(204, 149, 8, 172)$,
$(204, 201, 107, 19)$,
$(204, 10, 23, 87)$,\adfsplit
$(204, 52, 79, 108)$,
$(204, 65, 81, 30)$,
$(204, 2, 3, 155)$,
$(204, 145, 20, 18)$,\adfsplit
$(0, 10, 199, 216)$,
$(0, 25, 65, 106)$,
$(0, 59, 115, 236)$,
$(0, 7, 8, 239)$,\adfsplit
$(1, 21, 71, 231)$,
$(0, 79, 197, 231)$,
$(0, 29, 105, 219)$,
$(0, 55, 89, 244)$,\adfsplit
$(0, 31, 86, 248)$,
$(0, 73, 81, 226)$,
$(0, 85, 185, 204)$,
$(0, 53, 128, 227)$,\adfsplit
$(0, 159, 169, 232)$,
$(0, 34, 167, 241)$,
$(0, 61, 63, 228)$,
$(0, 44, 179, 209)$,\adfsplit
$(0, 52, 110, 217)$,
$(0, 50, 112, 207)$,
$(0, 11, 178, 230)$,
$(285, 0, 68, 136)$,\adfsplit
$(285, 1, 69, 137)$

\adfLgap \noindent by the mapping:
$x \mapsto x + 2 j \adfmod{204}$ for $x < 204$,
$x \mapsto (x +  j \adfmod{51}) + 204$ for $204 \le x < 255$,
$x \mapsto (x - 255 + 5 j \adfmod{30}) + 255$ for $255 \le x < 285$,
$285 \mapsto 285$,
$0 \le j < 102$
 for the first 55 blocks,
$0 \le j < 34$
 for the last two blocks.
\ADFvfyParStart{(286, ((55, 102, ((204, 2), (51, 1), (30, 5), (1, 1))), (2, 34, ((204, 2), (51, 1), (30, 5), (1, 1)))), ((34, 6), (82, 1)))} 

\adfDgap
\noindent{\boldmath $ 34^{9} 19^{1} $}~
With the point set $Z_{325}$ partitioned into
 residue classes modulo $9$ for $\{0, 1, \dots, 305\}$, and
 $\{306, 307, \dots, 324\}$,
 the design is generated from

\adfLgap 
$(306, 260, 49, 33)$,
$(306, 100, 234, 251)$,
$(307, 35, 86, 81)$,
$(307, 294, 220, 79)$,\adfsplit
$(308, 159, 56, 229)$,
$(308, 66, 106, 263)$,
$(309, 248, 187, 180)$,
$(309, 227, 286, 99)$,\adfsplit
$(310, 106, 50, 71)$,
$(310, 51, 90, 283)$,
$(311, 182, 244, 42)$,
$(311, 235, 29, 279)$,\adfsplit
$(54, 215, 102, 114)$,
$(121, 284, 286, 138)$,
$(57, 217, 31, 297)$,
$(254, 292, 69, 75)$,\adfsplit
$(8, 14, 106, 153)$,
$(27, 196, 141, 175)$,
$(217, 26, 193, 75)$,
$(193, 101, 296, 105)$,\adfsplit
$(221, 270, 240, 271)$,
$(73, 168, 273, 44)$,
$(69, 10, 206, 282)$,
$(210, 139, 29, 296)$,\adfsplit
$(199, 231, 23, 65)$,
$(202, 178, 201, 137)$,
$(145, 207, 85, 182)$,
$(169, 154, 176, 227)$,\adfsplit
$(170, 103, 117, 276)$,
$(151, 123, 208, 103)$,
$(280, 196, 228, 265)$,
$(61, 289, 281, 204)$,\adfsplit
$(64, 20, 31, 61)$,
$(0, 4, 61, 199)$,
$(0, 5, 94, 129)$,
$(0, 65, 177, 215)$,\adfsplit
$(0, 33, 55, 187)$,
$(0, 19, 87, 170)$,
$(0, 159, 263, 275)$,
$(0, 97, 181, 233)$,\adfsplit
$(0, 71, 281, 283)$,
$(0, 107, 183, 194)$,
$(0, 43, 53, 174)$,
$(0, 3, 227, 296)$,\adfsplit
$(0, 14, 123, 152)$,
$(0, 46, 128, 259)$,
$(0, 26, 205, 218)$,
$(0, 8, 75, 150)$,\adfsplit
$(0, 16, 96, 116)$,
$(0, 58, 122, 188)$,
$(0, 28, 70, 256)$,
$(324, 0, 102, 204)$,\adfsplit
$(324, 1, 103, 205)$

\adfLgap \noindent by the mapping:
$x \mapsto x + 2 j \adfmod{306}$ for $x < 306$,
$x \mapsto (x + 6 j \adfmod{18}) + 306$ for $306 \le x < 324$,
$324 \mapsto 324$,
$0 \le j < 153$
 for the first 51 blocks,
$0 \le j < 51$
 for the last two blocks.
\ADFvfyParStart{(325, ((51, 153, ((306, 2), (18, 6), (1, 1))), (2, 51, ((306, 2), (18, 6), (1, 1)))), ((34, 9), (19, 1)))} 

\adfDgap
\noindent{\boldmath $ 34^{12} 25^{1} $}~
With the point set $Z_{433}$ partitioned into
 residue classes modulo $12$ for $\{0, 1, \dots, 407\}$, and
 $\{408, 409, \dots, 432\}$,
 the design is generated from

\adfLgap 
$(408, 224, 75, 241)$,
$(408, 299, 298, 24)$,
$(409, 115, 154, 192)$,
$(409, 8, 65, 303)$,\adfsplit
$(410, 89, 52, 206)$,
$(410, 165, 84, 277)$,
$(411, 166, 161, 361)$,
$(411, 242, 33, 144)$,\adfsplit
$(82, 171, 223, 107)$,
$(20, 107, 51, 91)$,
$(14, 392, 401, 204)$,
$(22, 33, 136, 267)$,\adfsplit
$(184, 31, 250, 252)$,
$(176, 133, 58, 203)$,
$(83, 145, 90, 189)$,
$(404, 321, 202, 354)$,\adfsplit
$(239, 293, 66, 242)$,
$(127, 392, 109, 168)$,
$(102, 267, 144, 7)$,
$(386, 321, 229, 401)$,\adfsplit
$(202, 366, 293, 189)$,
$(0, 4, 142, 150)$,
$(0, 6, 175, 185)$,
$(0, 20, 102, 130)$,\adfsplit
$(0, 14, 63, 121)$,
$(0, 19, 88, 210)$,
$(0, 23, 76, 162)$,
$(0, 22, 101, 127)$,\adfsplit
$(0, 74, 159, 268)$,
$(0, 29, 196, 230)$,
$(0, 46, 93, 330)$,
$(0, 61, 155, 222)$,\adfsplit
$(0, 35, 135, 261)$,
$(0, 45, 160, 250)$,
$(0, 32, 129, 257)$,
$(432, 0, 136, 272)$

\adfLgap \noindent by the mapping:
$x \mapsto x +  j \adfmod{408}$ for $x < 408$,
$x \mapsto (x + 4 j \adfmod{24}) + 408$ for $408 \le x < 432$,
$432 \mapsto 432$,
$0 \le j < 408$
 for the first 35 blocks,
$0 \le j < 136$
 for the last block.
\ADFvfyParStart{(433, ((35, 408, ((408, 1), (24, 4), (1, 1))), (1, 136, ((408, 1), (24, 4), (1, 1)))), ((34, 12), (25, 1)))} 

\adfDgap
\noindent{\boldmath $ 34^{12} 28^{1} $}~
With the point set $Z_{436}$ partitioned into
 residue classes modulo $12$ for $\{0, 1, \dots, 407\}$, and
 $\{408, 409, \dots, 435\}$,
 the design is generated from

\adfLgap 
$(408, 46, 235, 72)$,
$(408, 128, 21, 119)$,
$(409, 366, 33, 253)$,
$(409, 238, 41, 134)$,\adfsplit
$(410, 311, 265, 112)$,
$(410, 285, 348, 308)$,
$(411, 282, 56, 292)$,
$(411, 75, 259, 287)$,\adfsplit
$(412, 288, 386, 161)$,
$(412, 334, 277, 255)$,
$(413, 322, 90, 365)$,
$(413, 290, 325, 375)$,\adfsplit
$(414, 145, 10, 129)$,
$(414, 305, 204, 56)$,
$(415, 252, 67, 68)$,
$(415, 339, 275, 370)$,\adfsplit
$(416, 400, 141, 14)$,
$(416, 313, 198, 29)$,
$(213, 164, 301, 352)$,
$(263, 345, 363, 158)$,\adfsplit
$(110, 152, 173, 244)$,
$(105, 54, 124, 111)$,
$(140, 174, 85, 391)$,
$(29, 371, 366, 292)$,\adfsplit
$(148, 5, 302, 187)$,
$(205, 355, 50, 39)$,
$(344, 192, 385, 242)$,
$(256, 355, 279, 128)$,\adfsplit
$(176, 15, 385, 204)$,
$(389, 362, 251, 22)$,
$(292, 214, 63, 169)$,
$(117, 334, 64, 83)$,\adfsplit
$(44, 41, 15, 150)$,
$(378, 80, 211, 288)$,
$(184, 19, 298, 398)$,
$(153, 221, 358, 188)$,\adfsplit
$(186, 398, 67, 111)$,
$(12, 340, 20, 347)$,
$(12, 218, 219, 382)$,
$(25, 72, 233, 218)$,\adfsplit
$(193, 339, 388, 406)$,
$(43, 328, 117, 373)$,
$(387, 78, 218, 331)$,
$(100, 195, 225, 85)$,\adfsplit
$(78, 125, 16, 360)$,
$(205, 71, 302, 201)$,
$(66, 308, 41, 369)$,
$(293, 358, 151, 33)$,\adfsplit
$(0, 3, 30, 356)$,
$(0, 2, 6, 394)$,
$(0, 9, 46, 158)$,
$(0, 11, 13, 354)$,\adfsplit
$(0, 25, 142, 364)$,
$(0, 32, 118, 278)$,
$(0, 37, 210, 375)$,
$(0, 17, 76, 261)$,\adfsplit
$(0, 56, 159, 342)$,
$(0, 73, 94, 218)$,
$(0, 58, 174, 233)$,
$(0, 91, 325, 339)$,\adfsplit
$(0, 87, 139, 178)$,
$(0, 93, 133, 255)$,
$(0, 31, 89, 267)$,
$(0, 81, 171, 258)$,\adfsplit
$(0, 97, 347, 401)$,
$(0, 277, 287, 403)$,
$(0, 117, 323, 355)$,
$(0, 107, 221, 283)$,\adfsplit
$(1, 9, 95, 223)$,
$(0, 237, 279, 391)$,
$(0, 147, 167, 365)$,
$(435, 0, 136, 272)$,\adfsplit
$(435, 1, 137, 273)$

\adfLgap \noindent by the mapping:
$x \mapsto x + 2 j \adfmod{408}$ for $x < 408$,
$x \mapsto (x - 408 + 9 j \adfmod{27}) + 408$ for $408 \le x < 435$,
$435 \mapsto 435$,
$0 \le j < 204$
 for the first 71 blocks,
$0 \le j < 68$
 for the last two blocks.
\ADFvfyParStart{(436, ((71, 204, ((408, 2), (27, 9), (1, 1))), (2, 68, ((408, 2), (27, 9), (1, 1)))), ((34, 12), (28, 1)))} 

\adfDgap
\noindent{\boldmath $ 34^{12} 31^{1} $}~
With the point set $Z_{439}$ partitioned into
 residue classes modulo $12$ for $\{0, 1, \dots, 407\}$, and
 $\{408, 409, \dots, 438\}$,
 the design is generated from

\adfLgap 
$(408, 248, 179, 322)$,
$(408, 252, 103, 75)$,
$(409, 140, 394, 143)$,
$(409, 127, 132, 267)$,\adfsplit
$(410, 237, 59, 160)$,
$(410, 170, 229, 276)$,
$(411, 358, 351, 294)$,
$(411, 350, 305, 337)$,\adfsplit
$(412, 237, 215, 103)$,
$(412, 236, 396, 346)$,
$(48, 353, 169, 104)$,
$(4, 301, 177, 119)$,\adfsplit
$(378, 50, 139, 116)$,
$(242, 119, 334, 145)$,
$(144, 244, 37, 297)$,
$(226, 24, 265, 314)$,\adfsplit
$(335, 302, 165, 378)$,
$(278, 196, 58, 102)$,
$(188, 137, 172, 158)$,
$(192, 64, 65, 181)$,\adfsplit
$(22, 283, 152, 245)$,
$(0, 2, 10, 196)$,
$(0, 4, 19, 286)$,
$(0, 6, 40, 250)$,\adfsplit
$(0, 17, 54, 275)$,
$(0, 31, 102, 183)$,
$(0, 9, 61, 266)$,
$(0, 29, 174, 229)$,\adfsplit
$(0, 75, 162, 247)$,
$(0, 18, 91, 181)$,
$(0, 20, 99, 289)$,
$(0, 46, 109, 171)$,\adfsplit
$(0, 41, 83, 294)$,
$(0, 27, 105, 322)$,
$(0, 67, 165, 233)$,
$(0, 25, 95, 304)$,\adfsplit
$(438, 0, 136, 272)$

\adfLgap \noindent by the mapping:
$x \mapsto x +  j \adfmod{408}$ for $x < 408$,
$x \mapsto (x - 408 + 5 j \adfmod{30}) + 408$ for $408 \le x < 438$,
$438 \mapsto 438$,
$0 \le j < 408$
 for the first 36 blocks,
$0 \le j < 136$
 for the last block.
\ADFvfyParStart{(439, ((36, 408, ((408, 1), (30, 5), (1, 1))), (1, 136, ((408, 1), (30, 5), (1, 1)))), ((34, 12), (31, 1)))} 

\adfDgap
\noindent{\boldmath $ 34^{15} 31^{1} $}~
With the point set $Z_{541}$ partitioned into
 residue classes modulo $15$ for $\{0, 1, \dots, 509\}$, and
 $\{510, 511, \dots, 540\}$,
 the design is generated from

\adfLgap 
$(510, 246, 332, 147)$,
$(510, 319, 89, 88)$,
$(511, 123, 169, 311)$,
$(511, 298, 446, 102)$,\adfsplit
$(512, 356, 162, 175)$,
$(512, 195, 334, 263)$,
$(513, 273, 94, 222)$,
$(513, 295, 11, 410)$,\adfsplit
$(514, 408, 141, 296)$,
$(514, 499, 503, 436)$,
$(515, 473, 7, 464)$,
$(515, 417, 4, 456)$,\adfsplit
$(516, 42, 64, 434)$,
$(516, 5, 205, 459)$,
$(517, 80, 358, 402)$,
$(517, 489, 395, 61)$,\adfsplit
$(518, 44, 413, 190)$,
$(518, 180, 387, 145)$,
$(519, 163, 378, 29)$,
$(519, 302, 195, 172)$,\adfsplit
$(94, 356, 119, 223)$,
$(266, 304, 106, 380)$,
$(70, 63, 256, 124)$,
$(287, 465, 53, 416)$,\adfsplit
$(406, 54, 338, 228)$,
$(443, 102, 319, 269)$,
$(37, 276, 509, 338)$,
$(506, 421, 276, 43)$,\adfsplit
$(89, 111, 482, 502)$,
$(53, 358, 454, 417)$,
$(61, 7, 326, 154)$,
$(219, 112, 178, 35)$,\adfsplit
$(162, 106, 39, 36)$,
$(507, 146, 349, 224)$,
$(147, 190, 298, 341)$,
$(191, 129, 246, 462)$,\adfsplit
$(80, 280, 453, 52)$,
$(394, 296, 83, 429)$,
$(251, 353, 346, 88)$,
$(106, 497, 208, 414)$,\adfsplit
$(290, 18, 119, 98)$,
$(205, 102, 404, 378)$,
$(286, 25, 323, 399)$,
$(153, 76, 372, 298)$,\adfsplit
$(90, 124, 43, 67)$,
$(230, 431, 244, 242)$,
$(16, 100, 116, 297)$,
$(2, 223, 281, 464)$,\adfsplit
$(493, 46, 483, 101)$,
$(318, 259, 329, 27)$,
$(216, 368, 289, 211)$,
$(383, 455, 302, 384)$,\adfsplit
$(412, 256, 469, 205)$,
$(467, 329, 324, 160)$,
$(365, 492, 272, 141)$,
$(57, 258, 352, 349)$,\adfsplit
$(368, 20, 257, 447)$,
$(190, 27, 489, 317)$,
$(244, 428, 420, 231)$,
$(135, 393, 231, 341)$,\adfsplit
$(398, 29, 416, 292)$,
$(438, 251, 233, 502)$,
$(334, 190, 285, 200)$,
$(0, 19, 242, 266)$,\adfsplit
$(0, 6, 33, 394)$,
$(0, 4, 40, 92)$,
$(0, 17, 29, 46)$,
$(0, 32, 131, 292)$,\adfsplit
$(0, 39, 183, 374)$,
$(0, 89, 264, 406)$,
$(0, 50, 115, 204)$,
$(0, 137, 251, 372)$,\adfsplit
$(0, 168, 469, 501)$,
$(0, 42, 328, 483)$,
$(0, 185, 211, 298)$,
$(0, 87, 254, 489)$,\adfsplit
$(0, 31, 157, 190)$,
$(0, 61, 72, 479)$,
$(0, 229, 385, 491)$,
$(0, 69, 261, 457)$,\adfsplit
$(0, 253, 267, 419)$,
$(0, 257, 331, 485)$,
$(0, 227, 293, 409)$,
$(0, 121, 357, 445)$,\adfsplit
$(1, 3, 87, 217)$,
$(1, 17, 37, 447)$,
$(1, 9, 43, 351)$,
$(1, 7, 205, 245)$,\adfsplit
$(1, 29, 141, 289)$,
$(540, 0, 170, 340)$,
$(540, 1, 171, 341)$

\adfLgap \noindent by the mapping:
$x \mapsto x + 2 j \adfmod{510}$ for $x < 510$,
$x \mapsto (x + 10 j \adfmod{30}) + 510$ for $510 \le x < 540$,
$540 \mapsto 540$,
$0 \le j < 255$
 for the first 89 blocks,
$0 \le j < 85$
 for the last two blocks.
\ADFvfyParStart{(541, ((89, 255, ((510, 2), (30, 10), (1, 1))), (2, 85, ((510, 2), (30, 10), (1, 1)))), ((34, 15), (31, 1)))} 

\section{4-GDDs for the proof of Lemma \ref{lem:4-GDD 38^u m^1}}
\label{app:4-GDD 38^u m^1}
\adfhide{
$ 38^6 5^1 $,
$ 38^6 8^1 $,
$ 38^6 11^1 $,
$ 38^6 14^1 $,
$ 38^6 17^1 $,
$ 38^6 20^1 $,
$ 38^6 23^1 $,
$ 38^6 26^1 $,
$ 38^6 83^1 $,
$ 38^6 86^1 $,
$ 38^6 89^1 $,
$ 38^6 92^1 $,
$ 38^9 11^1 $,
$ 38^9 14^1 $,
$ 38^9 17^1 $,
$ 38^9 20^1 $,
$ 38^9 23^1 $,
$ 38^{12} 14^1 $,
$ 38^{12} 17^1 $,
$ 38^{12} 20^1 $,
$ 38^{12} 23^1 $,
$ 38^{15} 17^1 $,
$ 38^{15} 20^1 $,
$ 38^{15} 23^1 $,
$ 38^{15} 26^1 $,
$ 38^{15} 29^1 $,
$ 38^{18} 20^1 $,
$ 38^{18} 23^1 $,
$ 38^{18} 26^1 $,
$ 38^{18} 29^1 $,
$ 38^{18} 32^1 $ and
$ 38^{18} 35^1 $.
}

\adfDgap
\noindent{\boldmath $ 38^{6} 5^{1} $}~
With the point set $Z_{233}$ partitioned into
 residue classes modulo $6$ for $\{0, 1, \dots, 227\}$, and
 $\{228, 229, 230, 231, 232\}$,
 the design is generated from

\adfLgap 
$(228, 0, 1, 2)$,
$(229, 0, 76, 74)$,
$(230, 0, 151, 227)$,
$(231, 0, 154, 77)$,\adfsplit
$(232, 0, 226, 152)$,
$(61, 112, 147, 23)$,
$(35, 207, 8, 166)$,
$(113, 26, 21, 66)$,\adfsplit
$(173, 61, 68, 78)$,
$(177, 42, 157, 29)$,
$(0, 3, 14, 46)$,
$(0, 4, 63, 71)$,\adfsplit
$(0, 9, 31, 64)$,
$(0, 21, 79, 167)$,
$(0, 15, 34, 145)$,
$(0, 39, 91, 166)$,\adfsplit
$(0, 50, 103, 171)$,
$(0, 23, 49, 122)$,
$(0, 16, 85, 110)$,
$(0, 28, 65, 109)$

\adfLgap \noindent by the mapping:
$x \mapsto x + 3 j \adfmod{228}$ for $x < 228$,
$x \mapsto x$ for $x \ge 228$,
$0 \le j < 76$
 for the first five blocks;
$x \mapsto x +  j \adfmod{228}$ for $x < 228$,
$x \mapsto x$ for $x \ge 228$,
$0 \le j < 228$
 for the last 15 blocks.
\ADFvfyParStart{(233, ((5, 76, ((228, 3), (5, 5))), (15, 228, ((228, 1), (5, 5)))), ((38, 6), (5, 1)))} 

\adfDgap
\noindent{\boldmath $ 38^{6} 8^{1} $}~
With the point set $Z_{236}$ partitioned into
 residue classes modulo $6$ for $\{0, 1, \dots, 227\}$, and
 $\{228, 229, \dots, 235\}$,
 the design is generated from

\adfLgap 
$(228, 0, 1, 2)$,
$(229, 0, 76, 74)$,
$(230, 0, 151, 227)$,
$(231, 0, 154, 77)$,\adfsplit
$(232, 0, 226, 152)$,
$(233, 0, 4, 11)$,
$(234, 0, 7, 224)$,
$(235, 0, 217, 221)$,\adfsplit
$(128, 43, 138, 190)$,
$(93, 161, 142, 102)$,
$(122, 76, 99, 19)$,
$(3, 97, 94, 198)$,\adfsplit
$(38, 106, 155, 19)$,
$(104, 29, 70, 99)$,
$(134, 217, 189, 172)$,
$(163, 100, 53, 116)$,\adfsplit
$(76, 49, 33, 18)$,
$(162, 39, 49, 170)$,
$(160, 24, 127, 182)$,
$(0, 5, 20, 148)$,\adfsplit
$(0, 9, 164, 191)$,
$(0, 13, 21, 130)$,
$(0, 25, 50, 184)$,
$(0, 26, 61, 142)$,\adfsplit
$(0, 32, 123, 167)$,
$(0, 41, 75, 139)$,
$(0, 43, 101, 140)$,
$(0, 28, 157, 177)$,\adfsplit
$(0, 67, 99, 155)$,
$(0, 113, 163, 225)$,
$(0, 39, 65, 172)$,
$(0, 51, 89, 110)$,\adfsplit
$(0, 57, 71, 199)$,
$(0, 14, 159, 211)$,
$(0, 87, 109, 122)$,
$(0, 37, 82, 161)$,\adfsplit
$(0, 53, 93, 175)$

\adfLgap \noindent by the mapping:
$x \mapsto x + 3 j \adfmod{228}$ for $x < 228$,
$x \mapsto x$ for $x \ge 228$,
$0 \le j < 76$
 for the first eight blocks;
$x \mapsto x + 2 j \adfmod{228}$ for $x < 228$,
$x \mapsto x$ for $x \ge 228$,
$0 \le j < 114$
 for the last 29 blocks.
\ADFvfyParStart{(236, ((8, 76, ((228, 3), (8, 8))), (29, 114, ((228, 2), (8, 8)))), ((38, 6), (8, 1)))} 

\adfDgap
\noindent{\boldmath $ 38^{6} 11^{1} $}~
With the point set $Z_{239}$ partitioned into
 residue classes modulo $6$ for $\{0, 1, \dots, 227\}$, and
 $\{228, 229, \dots, 238\}$,
 the design is generated from

\adfLgap 
$(234, 0, 1, 2)$,
$(235, 0, 76, 74)$,
$(236, 0, 151, 227)$,
$(237, 0, 154, 77)$,\adfsplit
$(238, 0, 226, 152)$,
$(228, 25, 6, 77)$,
$(228, 166, 45, 128)$,
$(161, 48, 70, 20)$,\adfsplit
$(205, 99, 88, 65)$,
$(44, 11, 24, 27)$,
$(47, 38, 87, 156)$,
$(0, 4, 85, 93)$,\adfsplit
$(0, 5, 37, 64)$,
$(0, 21, 62, 101)$,
$(0, 35, 82, 133)$,
$(0, 10, 63, 163)$,\adfsplit
$(0, 25, 56, 161)$,
$(0, 26, 94, 155)$,
$(0, 46, 103, 158)$,
$(0, 7, 86, 131)$,\adfsplit
$(0, 14, 43, 58)$

\adfLgap \noindent by the mapping:
$x \mapsto x + 3 j \adfmod{228}$ for $x < 228$,
$x \mapsto (x +  j \adfmod{6}) + 228$ for $228 \le x < 234$,
$x \mapsto x$ for $x \ge 234$,
$0 \le j < 76$
 for the first five blocks;
$x \mapsto x +  j \adfmod{228}$ for $x < 228$,
$x \mapsto (x +  j \adfmod{6}) + 228$ for $228 \le x < 234$,
$x \mapsto x$ for $x \ge 234$,
$0 \le j < 228$
 for the last 16 blocks.
\ADFvfyParStart{(239, ((5, 76, ((228, 3), (6, 1), (5, 5))), (16, 228, ((228, 1), (6, 1), (5, 5)))), ((38, 6), (11, 1)))} 

\adfDgap
\noindent{\boldmath $ 38^{6} 14^{1} $}~
With the point set $Z_{242}$ partitioned into
 residue classes modulo $6$ for $\{0, 1, \dots, 227\}$, and
 $\{228, 229, \dots, 241\}$,
 the design is generated from

\adfLgap 
$(234, 0, 1, 2)$,
$(235, 0, 76, 74)$,
$(236, 0, 151, 227)$,
$(237, 0, 154, 77)$,\adfsplit
$(238, 0, 226, 152)$,
$(239, 0, 4, 11)$,
$(240, 0, 7, 224)$,
$(241, 0, 217, 221)$,\adfsplit
$(228, 54, 227, 2)$,
$(228, 202, 115, 15)$,
$(228, 113, 0, 157)$,
$(228, 4, 176, 213)$,\adfsplit
$(89, 180, 148, 223)$,
$(214, 25, 200, 105)$,
$(125, 116, 94, 174)$,
$(226, 23, 96, 79)$,\adfsplit
$(195, 217, 86, 226)$,
$(108, 219, 193, 203)$,
$(77, 182, 64, 69)$,
$(147, 210, 128, 107)$,\adfsplit
$(155, 194, 168, 52)$,
$(0, 3, 17, 68)$,
$(0, 8, 23, 201)$,
$(0, 16, 62, 83)$,\adfsplit
$(0, 27, 94, 194)$,
$(0, 33, 145, 190)$,
$(0, 50, 101, 129)$,
$(0, 59, 91, 195)$,\adfsplit
$(0, 44, 149, 187)$,
$(0, 73, 135, 203)$,
$(0, 47, 93, 122)$,
$(0, 69, 89, 158)$,\adfsplit
$(0, 40, 104, 167)$,
$(0, 28, 57, 175)$,
$(0, 35, 87, 121)$,
$(0, 10, 65, 223)$,\adfsplit
$(0, 45, 92, 191)$,
$(0, 49, 107, 171)$,
$(0, 20, 117, 205)$

\adfLgap \noindent by the mapping:
$x \mapsto x + 3 j \adfmod{228}$ for $x < 228$,
$x \mapsto (x +  j \adfmod{6}) + 228$ for $228 \le x < 234$,
$x \mapsto x$ for $x \ge 234$,
$0 \le j < 76$
 for the first eight blocks;
$x \mapsto x + 2 j \adfmod{228}$ for $x < 228$,
$x \mapsto (x +  j \adfmod{6}) + 228$ for $228 \le x < 234$,
$x \mapsto x$ for $x \ge 234$,
$0 \le j < 114$
 for the last 31 blocks.
\ADFvfyParStart{(242, ((8, 76, ((228, 3), (6, 1), (8, 8))), (31, 114, ((228, 2), (6, 1), (8, 8)))), ((38, 6), (14, 1)))} 

\adfDgap
\noindent{\boldmath $ 38^{6} 17^{1} $}~
With the point set $Z_{245}$ partitioned into
 residue classes modulo $6$ for $\{0, 1, \dots, 227\}$, and
 $\{228, 229, \dots, 244\}$,
 the design is generated from

\adfLgap 
$(240, 0, 1, 2)$,
$(241, 0, 76, 74)$,
$(242, 0, 151, 227)$,
$(243, 0, 154, 77)$,\adfsplit
$(244, 0, 226, 152)$,
$(228, 220, 7, 135)$,
$(228, 179, 62, 96)$,
$(229, 168, 215, 148)$,\adfsplit
$(229, 98, 159, 1)$,
$(209, 128, 201, 88)$,
$(123, 58, 17, 79)$,
$(143, 50, 196, 192)$,\adfsplit
$(0, 3, 10, 139)$,
$(0, 5, 14, 27)$,
$(0, 16, 33, 71)$,
$(0, 31, 63, 164)$,\adfsplit
$(0, 28, 57, 116)$,
$(0, 23, 75, 110)$,
$(0, 45, 91, 170)$,
$(0, 19, 56, 124)$,\adfsplit
$(0, 25, 51, 94)$,
$(0, 11, 50, 130)$

\adfLgap \noindent by the mapping:
$x \mapsto x + 3 j \adfmod{228}$ for $x < 228$,
$x \mapsto (x + 2 j \adfmod{12}) + 228$ for $228 \le x < 240$,
$x \mapsto x$ for $x \ge 240$,
$0 \le j < 76$
 for the first five blocks;
$x \mapsto x +  j \adfmod{228}$ for $x < 228$,
$x \mapsto (x + 2 j \adfmod{12}) + 228$ for $228 \le x < 240$,
$x \mapsto x$ for $x \ge 240$,
$0 \le j < 228$
 for the last 17 blocks.
\ADFvfyParStart{(245, ((5, 76, ((228, 3), (12, 2), (5, 5))), (17, 228, ((228, 1), (12, 2), (5, 5)))), ((38, 6), (17, 1)))} 

\adfDgap
\noindent{\boldmath $ 38^{6} 20^{1} $}~
With the point set $Z_{248}$ partitioned into
 residue classes modulo $5$ for $\{0, 1, \dots, 189\}$,
 $\{190, 191, \dots, 227\}$, and
 $\{228, 229, \dots, 247\}$,
 the design is generated from

\adfLgap 
$(235, 215, 160, 158)$,
$(243, 220, 33, 100)$,
$(244, 210, 91, 73)$,
$(234, 214, 125, 61)$,\adfsplit
$(226, 44, 65, 37)$,
$(199, 108, 54, 171)$,
$(204, 137, 179, 18)$,
$(204, 11, 14, 95)$,\adfsplit
$(0, 1, 178, 236)$,
$(0, 8, 19, 246)$,
$(0, 4, 76, 82)$,
$(0, 9, 23, 47)$,\adfsplit
$(0, 16, 62, 93)$,
$(0, 17, 43, 146)$,
$(0, 49, 137, 211)$,
$(0, 48, 99, 202)$,\adfsplit
$(0, 32, 121, 204)$,
$(0, 33, 131, 201)$,
$(0, 34, 86, 208)$,
$(0, 22, 154, 216)$,\adfsplit
$(0, 37, 116, 241)$,
$(0, 39, 133, 231)$,
$(0, 27, 68, 134)$

\adfLgap \noindent by the mapping:
$x \mapsto x +  j \adfmod{190}$ for $x < 190$,
$x \mapsto (x +  j \adfmod{38}) + 190$ for $190 \le x < 228$,
$x \mapsto (x - 228 + 2 j \adfmod{20}) + 228$ for $x \ge 228$,
$0 \le j < 190$.
\ADFvfyParStart{(248, ((23, 190, ((190, 1), (38, 1), (20, 2)))), ((38, 5), (38, 1), (20, 1)))} 

\adfDgap
\noindent{\boldmath $ 38^{6} 23^{1} $}~
With the point set $Z_{251}$ partitioned into
 residue classes modulo $6$ for $\{0, 1, \dots, 227\}$, and
 $\{228, 229, \dots, 250\}$,
 the design is generated from

\adfLgap 
$(246, 0, 1, 2)$,
$(247, 0, 76, 74)$,
$(248, 0, 151, 227)$,
$(249, 0, 154, 77)$,\adfsplit
$(250, 0, 226, 152)$,
$(228, 163, 23, 196)$,
$(228, 0, 116, 177)$,
$(229, 6, 17, 81)$,\adfsplit
$(229, 170, 202, 61)$,
$(230, 87, 77, 139)$,
$(230, 6, 92, 46)$,
$(52, 48, 39, 145)$,\adfsplit
$(111, 41, 126, 104)$,
$(0, 3, 8, 211)$,
$(0, 14, 45, 193)$,
$(0, 16, 57, 185)$,\adfsplit
$(0, 28, 101, 139)$,
$(0, 27, 98, 145)$,
$(0, 26, 95, 129)$,
$(0, 29, 82, 121)$,\adfsplit
$(0, 21, 65, 115)$,
$(0, 19, 56, 123)$,
$(0, 23, 81, 160)$

\adfLgap \noindent by the mapping:
$x \mapsto x + 3 j \adfmod{228}$ for $x < 228$,
$x \mapsto (x - 228 + 3 j \adfmod{18}) + 228$ for $228 \le x < 246$,
$x \mapsto x$ for $x \ge 246$,
$0 \le j < 76$
 for the first five blocks;
$x \mapsto x +  j \adfmod{228}$ for $x < 228$,
$x \mapsto (x - 228 + 3 j \adfmod{18}) + 228$ for $228 \le x < 246$,
$x \mapsto x$ for $x \ge 246$,
$0 \le j < 228$
 for the last 18 blocks.
\ADFvfyParStart{(251, ((5, 76, ((228, 3), (18, 3), (5, 5))), (18, 228, ((228, 1), (18, 3), (5, 5)))), ((38, 6), (23, 1)))} 

\adfDgap
\noindent{\boldmath $ 38^{6} 26^{1} $}~
With the point set $Z_{254}$ partitioned into
 residue classes modulo $6$ for $\{0, 1, \dots, 227\}$, and
 $\{228, 229, \dots, 253\}$,
 the design is generated from

\adfLgap 
$(246, 0, 1, 2)$,
$(247, 0, 76, 74)$,
$(248, 0, 151, 227)$,
$(249, 0, 154, 77)$,\adfsplit
$(250, 0, 226, 152)$,
$(251, 0, 4, 11)$,
$(252, 0, 7, 224)$,
$(253, 0, 217, 221)$,\adfsplit
$(228, 49, 156, 179)$,
$(228, 151, 219, 197)$,
$(228, 202, 188, 165)$,
$(228, 138, 182, 112)$,\adfsplit
$(229, 175, 84, 119)$,
$(229, 51, 164, 222)$,
$(229, 10, 73, 117)$,
$(229, 76, 17, 122)$,\adfsplit
$(230, 89, 216, 49)$,
$(230, 224, 95, 112)$,
$(230, 225, 62, 222)$,
$(230, 202, 127, 15)$,\adfsplit
$(204, 117, 14, 4)$,
$(13, 27, 174, 221)$,
$(104, 193, 226, 138)$,
$(62, 156, 70, 13)$,\adfsplit
$(192, 93, 109, 128)$,
$(0, 5, 183, 208)$,
$(0, 9, 17, 98)$,
$(0, 13, 117, 212)$,\adfsplit
$(0, 15, 146, 173)$,
$(0, 21, 149, 176)$,
$(0, 37, 40, 177)$,
$(0, 31, 59, 165)$,\adfsplit
$(0, 53, 100, 219)$,
$(0, 39, 80, 125)$,
$(0, 43, 135, 178)$,
$(0, 69, 136, 215)$,\adfsplit
$(0, 51, 124, 199)$,
$(0, 19, 71, 189)$,
$(0, 33, 65, 166)$,
$(0, 32, 143, 207)$,\adfsplit
$(0, 56, 213, 223)$,
$(0, 22, 105, 131)$,
$(0, 49, 87, 118)$

\adfLgap \noindent by the mapping:
$x \mapsto x + 3 j \adfmod{228}$ for $x < 228$,
$x \mapsto (x - 228 + 3 j \adfmod{18}) + 228$ for $228 \le x < 246$,
$x \mapsto x$ for $x \ge 246$,
$0 \le j < 76$
 for the first eight blocks;
$x \mapsto x + 2 j \adfmod{228}$ for $x < 228$,
$x \mapsto (x - 228 + 3 j \adfmod{18}) + 228$ for $228 \le x < 246$,
$x \mapsto x$ for $x \ge 246$,
$0 \le j < 114$
 for the last 35 blocks.
\ADFvfyParStart{(254, ((8, 76, ((228, 3), (18, 3), (8, 8))), (35, 114, ((228, 2), (18, 3), (8, 8)))), ((38, 6), (26, 1)))} 

\adfDgap
\noindent{\boldmath $ 38^{6} 83^{1} $}~
With the point set $Z_{311}$ partitioned into
 residue classes modulo $6$ for $\{0, 1, \dots, 227\}$, and
 $\{228, 229, \dots, 310\}$,
 the design is generated from

\adfLgap 
$(306, 0, 1, 2)$,
$(307, 0, 76, 74)$,
$(308, 0, 151, 227)$,
$(309, 0, 154, 77)$,\adfsplit
$(310, 0, 226, 152)$,
$(228, 182, 52, 175)$,
$(228, 141, 167, 156)$,
$(229, 212, 123, 120)$,\adfsplit
$(229, 37, 161, 202)$,
$(230, 144, 203, 163)$,
$(230, 64, 39, 56)$,
$(231, 120, 86, 221)$,\adfsplit
$(231, 15, 10, 67)$,
$(232, 99, 136, 126)$,
$(232, 79, 179, 200)$,
$(233, 206, 159, 197)$,\adfsplit
$(233, 36, 4, 91)$,
$(234, 200, 196, 51)$,
$(234, 95, 175, 210)$,
$(235, 210, 141, 226)$,\adfsplit
$(235, 97, 110, 11)$,
$(236, 58, 44, 125)$,
$(236, 201, 151, 18)$,
$(0, 23, 56, 117)$,\adfsplit
$(0, 29, 91, 250)$,
$(0, 22, 71, 276)$,
$(0, 58, 122, 238)$,
$(0, 20, 88, 251)$,\adfsplit
$(0, 44, 97, 265)$,
$(0, 46, 119, 278)$,
$(0, 31, 70, 253)$,
$(0, 51, 163, 305)$,\adfsplit
$(0, 28, 103, 146)$

\adfLgap \noindent by the mapping:
$x \mapsto x + 3 j \adfmod{228}$ for $x < 228$,
$x \mapsto (x - 228 + 13 j \adfmod{78}) + 228$ for $228 \le x < 306$,
$x \mapsto x$ for $x \ge 306$,
$0 \le j < 76$
 for the first five blocks;
$x \mapsto x +  j \adfmod{228}$ for $x < 228$,
$x \mapsto (x - 228 + 13 j \adfmod{78}) + 228$ for $228 \le x < 306$,
$x \mapsto x$ for $x \ge 306$,
$0 \le j < 228$
 for the last 28 blocks.
\ADFvfyParStart{(311, ((5, 76, ((228, 3), (78, 13), (5, 5))), (28, 228, ((228, 1), (78, 13), (5, 5)))), ((38, 6), (83, 1)))} 

\adfDgap
\noindent{\boldmath $ 38^{6} 86^{1} $}~
With the point set $Z_{314}$ partitioned into
 residue classes modulo $6$ for $\{0, 1, \dots, 227\}$, and
 $\{228, 229, \dots, 313\}$,
 the design is generated from

\adfLgap 
$(306, 0, 1, 2)$,
$(307, 0, 76, 74)$,
$(308, 0, 151, 227)$,
$(309, 0, 154, 77)$,\adfsplit
$(310, 0, 226, 152)$,
$(311, 0, 4, 11)$,
$(312, 0, 7, 224)$,
$(313, 0, 217, 221)$,\adfsplit
$(228, 18, 40, 35)$,
$(228, 62, 127, 34)$,
$(228, 68, 3, 101)$,
$(228, 129, 120, 37)$,\adfsplit
$(229, 202, 87, 0)$,
$(229, 179, 91, 44)$,
$(229, 194, 161, 181)$,
$(229, 117, 162, 16)$,\adfsplit
$(230, 181, 32, 153)$,
$(230, 202, 114, 139)$,
$(230, 192, 160, 197)$,
$(230, 74, 171, 59)$,\adfsplit
$(231, 101, 171, 38)$,
$(231, 78, 119, 205)$,
$(231, 201, 64, 0)$,
$(231, 128, 94, 151)$,\adfsplit
$(232, 44, 97, 153)$,
$(232, 192, 23, 39)$,
$(232, 136, 19, 66)$,
$(232, 182, 202, 221)$,\adfsplit
$(233, 113, 212, 69)$,
$(233, 121, 30, 70)$,
$(233, 144, 175, 135)$,
$(233, 88, 74, 227)$,\adfsplit
$(234, 17, 213, 187)$,
$(234, 227, 108, 181)$,
$(234, 188, 27, 64)$,
$(234, 226, 50, 42)$,\adfsplit
$(235, 171, 5, 220)$,
$(235, 7, 120, 117)$,
$(235, 217, 32, 203)$,
$(235, 106, 222, 206)$,\adfsplit
$(236, 134, 54, 190)$,
$(236, 128, 105, 41)$,
$(236, 4, 127, 47)$,
$(236, 25, 219, 204)$,\adfsplit
$(237, 57, 217, 30)$,
$(0, 3, 38, 302)$,
$(0, 10, 71, 237)$,
$(0, 35, 211, 263)$,\adfsplit
$(0, 21, 106, 238)$,
$(0, 46, 125, 207)$,
$(0, 45, 118, 173)$,
$(0, 29, 58, 239)$,\adfsplit
$(0, 86, 189, 278)$,
$(0, 89, 99, 291)$,
$(0, 83, 177, 304)$,
$(0, 105, 209, 240)$,\adfsplit
$(0, 159, 167, 266)$,
$(0, 147, 197, 253)$,
$(0, 62, 130, 292)$,
$(0, 131, 169, 290)$,\adfsplit
$(0, 50, 157, 277)$,
$(0, 69, 94, 175)$,
$(0, 95, 117, 303)$

\adfLgap \noindent by the mapping:
$x \mapsto x + 3 j \adfmod{228}$ for $x < 228$,
$x \mapsto (x - 228 + 13 j \adfmod{78}) + 228$ for $228 \le x < 306$,
$x \mapsto x$ for $x \ge 306$,
$0 \le j < 76$
 for the first eight blocks;
$x \mapsto x + 2 j \adfmod{228}$ for $x < 228$,
$x \mapsto (x - 228 + 13 j \adfmod{78}) + 228$ for $228 \le x < 306$,
$x \mapsto x$ for $x \ge 306$,
$0 \le j < 114$
 for the last 55 blocks.
\ADFvfyParStart{(314, ((8, 76, ((228, 3), (78, 13), (8, 8))), (55, 114, ((228, 2), (78, 13), (8, 8)))), ((38, 6), (86, 1)))} 

\adfDgap
\noindent{\boldmath $ 38^{6} 89^{1} $}~
With the point set $Z_{317}$ partitioned into
 residue classes modulo $6$ for $\{0, 1, \dots, 227\}$, and
 $\{228, 229, \dots, 316\}$,
 the design is generated from

\adfLgap 
$(312, 0, 1, 2)$,
$(313, 0, 76, 74)$,
$(314, 0, 151, 227)$,
$(315, 0, 154, 77)$,\adfsplit
$(316, 0, 226, 152)$,
$(228, 111, 64, 86)$,
$(228, 54, 152, 121)$,
$(228, 84, 119, 178)$,\adfsplit
$(228, 69, 185, 103)$,
$(229, 34, 18, 175)$,
$(229, 170, 27, 53)$,
$(229, 36, 49, 117)$,\adfsplit
$(229, 83, 76, 80)$,
$(230, 13, 59, 96)$,
$(230, 10, 74, 21)$,
$(230, 99, 138, 200)$,\adfsplit
$(230, 163, 4, 53)$,
$(231, 21, 35, 174)$,
$(231, 68, 154, 51)$,
$(231, 0, 43, 5)$,\adfsplit
$(231, 4, 133, 110)$,
$(232, 65, 93, 188)$,
$(232, 88, 145, 38)$,
$(232, 202, 66, 139)$,\adfsplit
$(0, 9, 176, 295)$,
$(0, 20, 135, 254)$,
$(0, 21, 109, 241)$,
$(0, 44, 100, 255)$,\adfsplit
$(0, 33, 91, 269)$,
$(0, 8, 209, 248)$,
$(0, 32, 163, 261)$,
$(0, 10, 51, 80)$,\adfsplit
$(0, 15, 188, 247)$,
$(0, 45, 149, 303)$

\adfLgap \noindent by the mapping:
$x \mapsto x + 3 j \adfmod{228}$ for $x < 228$,
$x \mapsto (x - 228 + 7 j \adfmod{84}) + 228$ for $228 \le x < 312$,
$x \mapsto x$ for $x \ge 312$,
$0 \le j < 76$
 for the first five blocks;
$x \mapsto x +  j \adfmod{228}$ for $x < 228$,
$x \mapsto (x - 228 + 7 j \adfmod{84}) + 228$ for $228 \le x < 312$,
$x \mapsto x$ for $x \ge 312$,
$0 \le j < 228$
 for the last 29 blocks.
\ADFvfyParStart{(317, ((5, 76, ((228, 3), (84, 7), (5, 5))), (29, 228, ((228, 1), (84, 7), (5, 5)))), ((38, 6), (89, 1)))} 

\adfDgap
\noindent{\boldmath $ 38^{6} 92^{1} $}~
With the point set $Z_{320}$ partitioned into
 residue classes modulo $6$ for $\{0, 1, \dots, 227\}$, and
 $\{228, 229, \dots, 319\}$,
 the design is generated from

\adfLgap 
$(312, 0, 1, 2)$,
$(313, 0, 76, 74)$,
$(314, 0, 151, 227)$,
$(315, 0, 154, 77)$,\adfsplit
$(316, 0, 226, 152)$,
$(317, 0, 4, 11)$,
$(318, 0, 7, 224)$,
$(319, 0, 217, 221)$,\adfsplit
$(228, 89, 222, 86)$,
$(228, 31, 129, 190)$,
$(228, 171, 208, 120)$,
$(228, 176, 167, 97)$,\adfsplit
$(229, 81, 191, 186)$,
$(229, 109, 158, 46)$,
$(229, 187, 144, 5)$,
$(229, 128, 100, 3)$,\adfsplit
$(230, 208, 221, 128)$,
$(230, 71, 205, 12)$,
$(230, 57, 194, 70)$,
$(230, 103, 18, 123)$,\adfsplit
$(231, 65, 158, 117)$,
$(231, 160, 92, 207)$,
$(231, 114, 61, 167)$,
$(231, 115, 190, 204)$,\adfsplit
$(232, 112, 50, 169)$,
$(232, 92, 71, 175)$,
$(232, 183, 168, 5)$,
$(232, 226, 18, 105)$,\adfsplit
$(233, 216, 200, 43)$,
$(233, 197, 118, 3)$,
$(233, 194, 155, 93)$,
$(233, 162, 172, 193)$,\adfsplit
$(234, 193, 76, 221)$,
$(234, 153, 131, 186)$,
$(234, 3, 62, 199)$,
$(234, 68, 226, 132)$,\adfsplit
$(235, 89, 14, 177)$,
$(235, 188, 226, 144)$,
$(235, 28, 37, 47)$,
$(235, 111, 103, 126)$,\adfsplit
$(236, 45, 61, 88)$,
$(236, 154, 128, 59)$,
$(236, 185, 146, 204)$,
$(236, 223, 39, 42)$,\adfsplit
$(237, 40, 9, 197)$,
$(0, 8, 25, 237)$,
$(0, 29, 196, 279)$,
$(0, 37, 101, 307)$,\adfsplit
$(0, 23, 110, 183)$,
$(0, 49, 129, 238)$,
$(0, 27, 109, 239)$,
$(0, 161, 199, 240)$,\adfsplit
$(1, 15, 101, 268)$,
$(0, 34, 86, 268)$,
$(0, 46, 143, 296)$,
$(0, 22, 67, 253)$,\adfsplit
$(0, 100, 211, 294)$,
$(0, 35, 147, 241)$,
$(0, 33, 125, 309)$,
$(0, 41, 98, 295)$,\adfsplit
$(0, 165, 223, 252)$,
$(0, 40, 121, 308)$,
$(0, 177, 203, 255)$,
$(0, 99, 155, 311)$,\adfsplit
$(0, 50, 106, 283)$

\adfLgap \noindent by the mapping:
$x \mapsto x + 3 j \adfmod{228}$ for $x < 228$,
$x \mapsto (x - 228 + 14 j \adfmod{84}) + 228$ for $228 \le x < 312$,
$x \mapsto x$ for $x \ge 312$,
$0 \le j < 76$
 for the first eight blocks;
$x \mapsto x + 2 j \adfmod{228}$ for $x < 228$,
$x \mapsto (x - 228 + 14 j \adfmod{84}) + 228$ for $228 \le x < 312$,
$x \mapsto x$ for $x \ge 312$,
$0 \le j < 114$
 for the last 57 blocks.
\ADFvfyParStart{(320, ((8, 76, ((228, 3), (84, 14), (8, 8))), (57, 114, ((228, 2), (84, 14), (8, 8)))), ((38, 6), (92, 1)))} 

\adfDgap
\noindent{\boldmath $ 38^{9} 11^{1} $}~
With the point set $Z_{353}$ partitioned into
 residue classes modulo $9$ for $\{0, 1, \dots, 341\}$, and
 $\{342, 343, \dots, 352\}$,
 the design is generated from

\adfLgap 
$(351, 0, 1, 2)$,
$(352, 0, 115, 230)$,
$(342, 24, 292, 57)$,
$(342, 38, 139, 287)$,\adfsplit
$(343, 140, 263, 222)$,
$(343, 34, 285, 283)$,
$(344, 47, 188, 297)$,
$(344, 55, 274, 30)$,\adfsplit
$(311, 166, 134, 25)$,
$(157, 92, 268, 257)$,
$(35, 329, 289, 135)$,
$(261, 325, 239, 30)$,\adfsplit
$(324, 181, 12, 176)$,
$(115, 38, 9, 80)$,
$(180, 113, 327, 256)$,
$(296, 92, 205, 18)$,\adfsplit
$(7, 280, 114, 76)$,
$(218, 152, 203, 276)$,
$(222, 128, 147, 215)$,
$(203, 1, 123, 290)$,\adfsplit
$(38, 14, 130, 88)$,
$(121, 63, 125, 254)$,
$(258, 277, 244, 71)$,
$(279, 87, 98, 262)$,\adfsplit
$(228, 79, 58, 182)$,
$(253, 221, 148, 281)$,
$(327, 178, 266, 272)$,
$(91, 74, 231, 9)$,\adfsplit
$(319, 187, 224, 3)$,
$(108, 57, 316, 104)$,
$(137, 62, 247, 30)$,
$(204, 134, 230, 112)$,\adfsplit
$(293, 211, 209, 3)$,
$(237, 136, 143, 124)$,
$(0, 10, 44, 122)$,
$(0, 11, 211, 332)$,\adfsplit
$(0, 3, 34, 150)$,
$(0, 15, 174, 314)$,
$(0, 43, 129, 145)$,
$(0, 23, 120, 259)$,\adfsplit
$(0, 17, 57, 60)$,
$(0, 14, 186, 317)$,
$(0, 12, 83, 197)$,
$(0, 55, 114, 173)$,\adfsplit
$(0, 79, 163, 293)$,
$(0, 47, 48, 113)$,
$(0, 70, 149, 253)$,
$(0, 105, 190, 203)$,\adfsplit
$(0, 52, 262, 292)$,
$(0, 16, 94, 305)$,
$(0, 101, 285, 322)$,
$(0, 39, 196, 220)$,\adfsplit
$(3, 52, 99, 208)$,
$(3, 11, 130, 171)$,
$(3, 9, 35, 70)$

\adfLgap \noindent by the mapping:
$x \mapsto x \oplus (3 j)$ for $x < 342$,
$x \mapsto (x + 3 j \adfmod{9}) + 342$ for $342 \le x < 351$,
$x \mapsto x$ for $x \ge 351$,
$0 \le j < 114$
 for the first two blocks;
$x \mapsto x \oplus j \oplus j$ for $x < 342$,
$x \mapsto (x + 3 j \adfmod{9}) + 342$ for $342 \le x < 351$,
$x \mapsto x$ for $x \ge 351$,
$0 \le j < 171$
 for the last 53 blocks.
\ADFvfyParStart{(353, ((2, 114, ((342, 3, (114, 3)), (9, 3), (2, 2))), (53, 171, ((342, 2, (114, 3)), (9, 3), (2, 2)))), ((38, 9), (11, 1)))} 

\adfDgap
\noindent{\boldmath $ 38^{9} 14^{1} $}~
With the point set $Z_{356}$ partitioned into
 residue classes modulo $9$ for $\{0, 1, \dots, 341\}$, and
 $\{342, 343, \dots, 355\}$,
 the design is generated from

\adfLgap 
$(354, 0, 1, 2)$,
$(355, 0, 115, 230)$,
$(342, 274, 73, 105)$,
$(342, 221, 108, 146)$,\adfsplit
$(343, 67, 75, 47)$,
$(343, 226, 150, 56)$,
$(19, 60, 301, 158)$,
$(152, 101, 31, 96)$,\adfsplit
$(2, 301, 296, 330)$,
$(99, 213, 292, 304)$,
$(58, 260, 42, 126)$,
$(158, 273, 9, 78)$,\adfsplit
$(87, 64, 49, 137)$,
$(49, 341, 147, 88)$,
$(23, 102, 227, 260)$,
$(7, 339, 167, 139)$,\adfsplit
$(33, 55, 200, 80)$,
$(130, 278, 156, 320)$,
$(0, 3, 77, 168)$,
$(0, 17, 109, 128)$,\adfsplit
$(0, 6, 89, 93)$,
$(0, 20, 57, 296)$,
$(0, 59, 156, 242)$,
$(0, 4, 116, 212)$,\adfsplit
$(0, 24, 58, 215)$,
$(0, 13, 65, 129)$,
$(0, 21, 123, 184)$,
$(0, 30, 85, 209)$,\adfsplit
$(0, 46, 110, 238)$

\adfLgap \noindent by the mapping:
$x \mapsto x \oplus (3 j)$ for $x < 342$,
$x \mapsto (x - 342 + 2 j \adfmod{12}) + 342$ for $342 \le x < 354$,
$x \mapsto x$ for $x \ge 354$,
$0 \le j < 114$
 for the first two blocks;
$x \mapsto x \oplus j$ for $x < 342$,
$x \mapsto (x - 342 + 2 j \adfmod{12}) + 342$ for $342 \le x < 354$,
$x \mapsto x$ for $x \ge 354$,
$0 \le j < 342$
 for the last 27 blocks.
\ADFvfyParStart{(356, ((2, 114, ((342, 3, (114, 3)), (12, 2), (2, 2))), (27, 342, ((342, 1, (114, 3)), (12, 2), (2, 2)))), ((38, 9), (14, 1)))} 

\adfDgap
\noindent{\boldmath $ 38^{9} 17^{1} $}~
With the point set $Z_{359}$ partitioned into
 residue classes modulo $9$ for $\{0, 1, \dots, 341\}$, and
 $\{342, 343, \dots, 358\}$,
 the design is generated from

\adfLgap 
$(357, 0, 1, 2)$,
$(358, 0, 115, 230)$,
$(342, 40, 171, 134)$,
$(342, 313, 78, 293)$,\adfsplit
$(343, 77, 117, 92)$,
$(343, 70, 114, 31)$,
$(344, 329, 192, 40)$,
$(344, 229, 63, 134)$,\adfsplit
$(345, 289, 284, 341)$,
$(345, 3, 52, 6)$,
$(346, 88, 65, 67)$,
$(346, 336, 213, 332)$,\adfsplit
$(321, 257, 61, 57)$,
$(257, 298, 119, 148)$,
$(339, 200, 43, 230)$,
$(33, 319, 302, 260)$,\adfsplit
$(335, 42, 196, 275)$,
$(203, 287, 6, 79)$,
$(153, 160, 73, 57)$,
$(86, 110, 268, 179)$,\adfsplit
$(246, 306, 249, 218)$,
$(153, 70, 282, 251)$,
$(257, 189, 12, 73)$,
$(172, 200, 141, 297)$,\adfsplit
$(132, 176, 243, 58)$,
$(21, 216, 313, 60)$,
$(241, 266, 198, 253)$,
$(29, 274, 23, 133)$,\adfsplit
$(332, 138, 112, 5)$,
$(238, 218, 212, 273)$,
$(100, 96, 75, 176)$,
$(188, 100, 45, 166)$,\adfsplit
$(45, 197, 255, 267)$,
$(304, 206, 337, 131)$,
$(262, 279, 42, 194)$,
$(211, 182, 201, 61)$,\adfsplit
$(0, 5, 67, 173)$,
$(0, 10, 242, 250)$,
$(0, 16, 287, 334)$,
$(0, 17, 51, 165)$,\adfsplit
$(0, 89, 166, 236)$,
$(0, 33, 203, 262)$,
$(0, 47, 95, 304)$,
$(0, 94, 145, 249)$,\adfsplit
$(0, 19, 148, 178)$,
$(0, 167, 172, 183)$,
$(0, 77, 201, 214)$,
$(0, 35, 231, 273)$,\adfsplit
$(0, 65, 102, 276)$,
$(0, 105, 140, 328)$,
$(0, 41, 127, 228)$,
$(0, 37, 176, 255)$,\adfsplit
$(0, 48, 123, 229)$,
$(0, 14, 160, 210)$,
$(0, 26, 79, 246)$,
$(0, 38, 85, 163)$,\adfsplit
$(0, 49, 120, 258)$

\adfLgap \noindent by the mapping:
$x \mapsto x \oplus (3 j)$ for $x < 342$,
$x \mapsto (x - 342 + 5 j \adfmod{15}) + 342$ for $342 \le x < 357$,
$x \mapsto x$ for $x \ge 357$,
$0 \le j < 114$
 for the first two blocks;
$x \mapsto x \oplus j \oplus j$ for $x < 342$,
$x \mapsto (x - 342 + 5 j \adfmod{15}) + 342$ for $342 \le x < 357$,
$x \mapsto x$ for $x \ge 357$,
$0 \le j < 171$
 for the last 55 blocks.
\ADFvfyParStart{(359, ((2, 114, ((342, 3, (114, 3)), (15, 5), (2, 2))), (55, 171, ((342, 2, (114, 3)), (15, 5), (2, 2)))), ((38, 9), (17, 1)))} 

\adfDgap
\noindent{\boldmath $ 38^{9} 20^{1} $}~
With the point set $Z_{362}$ partitioned into
 residue classes modulo $9$ for $\{0, 1, \dots, 341\}$, and
 $\{342, 343, \dots, 361\}$,
 the design is generated from

\adfLgap 
$(360, 0, 1, 2)$,
$(361, 0, 115, 230)$,
$(342, 165, 216, 80)$,
$(342, 29, 22, 181)$,\adfsplit
$(343, 227, 297, 14)$,
$(343, 55, 36, 142)$,
$(344, 35, 235, 196)$,
$(344, 48, 134, 123)$,\adfsplit
$(122, 106, 39, 213)$,
$(89, 117, 11, 57)$,
$(7, 129, 321, 91)$,
$(69, 107, 59, 300)$,\adfsplit
$(184, 164, 214, 71)$,
$(95, 272, 258, 305)$,
$(171, 268, 46, 177)$,
$(118, 50, 42, 0)$,\adfsplit
$(180, 31, 131, 245)$,
$(31, 16, 161, 235)$,
$(194, 109, 130, 191)$,
$(0, 4, 141, 158)$,\adfsplit
$(0, 5, 22, 124)$,
$(0, 29, 66, 169)$,
$(0, 44, 121, 178)$,
$(0, 24, 98, 157)$,\adfsplit
$(0, 34, 89, 185)$,
$(0, 26, 136, 182)$,
$(0, 43, 148, 218)$,
$(0, 10, 35, 273)$,\adfsplit
$(0, 37, 95, 184)$,
$(0, 12, 68, 238)$

\adfLgap \noindent by the mapping:
$x \mapsto x \oplus (3 j)$ for $x < 342$,
$x \mapsto (x + 3 j \adfmod{18}) + 342$ for $342 \le x < 360$,
$x \mapsto x$ for $x \ge 360$,
$0 \le j < 114$
 for the first two blocks;
$x \mapsto x \oplus j$ for $x < 342$,
$x \mapsto (x + 3 j \adfmod{18}) + 342$ for $342 \le x < 360$,
$x \mapsto x$ for $x \ge 360$,
$0 \le j < 342$
 for the last 28 blocks.
\ADFvfyParStart{(362, ((2, 114, ((342, 3, (114, 3)), (18, 3), (2, 2))), (28, 342, ((342, 1, (114, 3)), (18, 3), (2, 2)))), ((38, 9), (20, 1)))} 

\adfDgap
\noindent{\boldmath $ 38^{9} 23^{1} $}~
With the point set $Z_{365}$ partitioned into
 residue classes modulo $9$ for $\{0, 1, \dots, 341\}$, and
 $\{342, 343, \dots, 364\}$,
 the design is generated from

\adfLgap 
$(363, 0, 1, 2)$,
$(364, 0, 115, 230)$,
$(342, 95, 270, 118)$,
$(342, 242, 21, 151)$,\adfsplit
$(343, 333, 5, 289)$,
$(343, 180, 310, 152)$,
$(344, 98, 234, 169)$,
$(344, 107, 189, 232)$,\adfsplit
$(345, 158, 252, 160)$,
$(345, 197, 165, 211)$,
$(346, 252, 95, 331)$,
$(346, 86, 297, 190)$,\adfsplit
$(347, 89, 284, 252)$,
$(347, 235, 310, 165)$,
$(348, 50, 87, 292)$,
$(348, 240, 263, 217)$,\adfsplit
$(132, 142, 126, 215)$,
$(136, 301, 41, 282)$,
$(136, 186, 193, 117)$,
$(277, 62, 216, 230)$,\adfsplit
$(314, 120, 279, 331)$,
$(338, 6, 146, 73)$,
$(186, 208, 16, 130)$,
$(14, 15, 106, 166)$,\adfsplit
$(218, 243, 98, 258)$,
$(228, 86, 135, 177)$,
$(144, 55, 287, 221)$,
$(175, 212, 147, 109)$,\adfsplit
$(111, 198, 253, 249)$,
$(138, 40, 325, 32)$,
$(30, 258, 205, 288)$,
$(79, 198, 316, 94)$,\adfsplit
$(182, 215, 50, 166)$,
$(104, 250, 44, 56)$,
$(167, 267, 206, 334)$,
$(162, 131, 179, 15)$,\adfsplit
$(0, 3, 283, 304)$,
$(0, 8, 104, 186)$,
$(0, 25, 59, 93)$,
$(0, 35, 169, 264)$,\adfsplit
$(0, 39, 44, 240)$,
$(0, 29, 229, 316)$,
$(0, 53, 204, 274)$,
$(0, 42, 235, 323)$,\adfsplit
$(0, 111, 121, 233)$,
$(0, 123, 212, 334)$,
$(0, 105, 128, 298)$,
$(0, 38, 221, 305)$,\adfsplit
$(0, 24, 202, 263)$,
$(0, 69, 251, 268)$,
$(0, 100, 227, 257)$,
$(0, 160, 184, 219)$,\adfsplit
$(0, 137, 220, 321)$,
$(0, 201, 213, 226)$,
$(0, 154, 231, 287)$,
$(0, 129, 136, 339)$,\adfsplit
$(0, 76, 262, 309)$,
$(3, 9, 77, 105)$,
$(3, 11, 100, 179)$

\adfLgap \noindent by the mapping:
$x \mapsto x \oplus (3 j)$ for $x < 342$,
$x \mapsto (x - 342 + 7 j \adfmod{21}) + 342$ for $342 \le x < 363$,
$x \mapsto x$ for $x \ge 363$,
$0 \le j < 114$
 for the first two blocks;
$x \mapsto x \oplus j \oplus j$ for $x < 342$,
$x \mapsto (x - 342 + 7 j \adfmod{21}) + 342$ for $342 \le x < 363$,
$x \mapsto x$ for $x \ge 363$,
$0 \le j < 171$
 for the last 57 blocks.
\ADFvfyParStart{(365, ((2, 114, ((342, 3, (114, 3)), (21, 7), (2, 2))), (57, 171, ((342, 2, (114, 3)), (21, 7), (2, 2)))), ((38, 9), (23, 1)))} 

\adfDgap
\noindent{\boldmath $ 38^{12} 14^{1} $}~
With the point set $Z_{470}$ partitioned into
 residue classes modulo $12$ for $\{0, 1, \dots, 455\}$, and
 $\{456, 457, \dots, 469\}$,
 the design is generated from

\adfLgap 
$(456, 0, 1, 2)$,
$(457, 0, 151, 305)$,
$(458, 0, 154, 152)$,
$(459, 0, 304, 455)$,\adfsplit
$(460, 0, 454, 302)$,
$(461, 0, 4, 11)$,
$(462, 0, 7, 452)$,
$(463, 0, 445, 449)$,\adfsplit
$(464, 0, 10, 23)$,
$(465, 0, 13, 446)$,
$(466, 0, 433, 443)$,
$(467, 0, 16, 35)$,\adfsplit
$(468, 0, 19, 440)$,
$(469, 0, 421, 437)$,
$(354, 254, 413, 232)$,
$(208, 187, 167, 453)$,\adfsplit
$(439, 270, 357, 131)$,
$(85, 423, 53, 136)$,
$(288, 258, 49, 329)$,
$(296, 387, 143, 381)$,\adfsplit
$(106, 363, 112, 45)$,
$(365, 452, 331, 88)$,
$(71, 312, 237, 380)$,
$(8, 399, 329, 145)$,\adfsplit
$(165, 329, 276, 414)$,
$(399, 168, 59, 265)$,
$(26, 92, 53, 267)$,
$(382, 161, 319, 205)$,\adfsplit
$(325, 304, 130, 399)$,
$(7, 364, 63, 318)$,
$(216, 405, 114, 374)$,
$(275, 377, 10, 228)$,\adfsplit
$(110, 150, 349, 236)$,
$(14, 205, 129, 22)$,
$(87, 7, 412, 454)$,
$(231, 414, 417, 257)$,\adfsplit
$(244, 18, 361, 442)$,
$(336, 266, 364, 305)$,
$(47, 208, 115, 398)$,
$(427, 70, 316, 113)$,\adfsplit
$(413, 81, 76, 391)$,
$(68, 25, 273, 451)$,
$(390, 13, 195, 55)$,
$(202, 69, 77, 122)$,\adfsplit
$(24, 433, 256, 404)$,
$(396, 308, 282, 148)$,
$(73, 364, 254, 87)$,
$(88, 51, 413, 313)$,\adfsplit
$(18, 83, 121, 420)$,
$(90, 371, 171, 334)$,
$(45, 133, 204, 365)$,
$(0, 311, 365, 82)$,\adfsplit
$(357, 186, 397, 406)$,
$(313, 207, 442, 150)$,
$(0, 14, 75, 333)$,
$(0, 18, 63, 153)$,\adfsplit
$(0, 17, 139, 167)$,
$(0, 15, 77, 143)$,
$(0, 25, 123, 187)$,
$(0, 129, 301, 351)$,\adfsplit
$(0, 109, 136, 355)$,
$(0, 133, 166, 307)$,
$(0, 55, 101, 453)$,
$(0, 203, 329, 439)$,\adfsplit
$(0, 67, 263, 268)$,
$(0, 83, 271, 326)$,
$(0, 44, 353, 431)$,
$(0, 73, 125, 344)$,\adfsplit
$(0, 58, 176, 399)$,
$(0, 50, 259, 366)$,
$(0, 74, 267, 359)$,
$(0, 249, 361, 379)$,\adfsplit
$(0, 34, 127, 206)$,
$(0, 49, 64, 278)$,
$(0, 33, 62, 256)$,
$(0, 69, 186, 328)$,\adfsplit
$(0, 52, 146, 202)$,
$(0, 20, 124, 294)$,
$(0, 38, 116, 222)$

\adfLgap \noindent by the mapping:
$x \mapsto x + 3 j \adfmod{456}$ for $x < 456$,
$x \mapsto x$ for $x \ge 456$,
$0 \le j < 152$
 for the first 14 blocks;
$x \mapsto x + 2 j \adfmod{456}$ for $x < 456$,
$x \mapsto x$ for $x \ge 456$,
$0 \le j < 228$
 for the last 65 blocks.
\ADFvfyParStart{(470, ((14, 152, ((456, 3), (14, 14))), (65, 228, ((456, 2), (14, 14)))), ((38, 12), (14, 1)))} 

\adfDgap
\noindent{\boldmath $ 38^{12} 17^{1} $}~
With the point set $Z_{473}$ partitioned into
 residue classes modulo $12$ for $\{0, 1, \dots, 455\}$, and
 $\{456, 457, \dots, 472\}$,
 the design is generated from

\adfLgap 
$(468, 0, 1, 2)$,
$(469, 0, 151, 305)$,
$(470, 0, 154, 152)$,
$(471, 0, 304, 455)$,\adfsplit
$(472, 0, 454, 302)$,
$(456, 249, 407, 243)$,
$(456, 26, 181, 246)$,
$(456, 67, 144, 5)$,\adfsplit
$(456, 82, 64, 140)$,
$(396, 205, 449, 323)$,
$(229, 238, 116, 323)$,
$(200, 434, 211, 311)$,\adfsplit
$(384, 255, 227, 145)$,
$(107, 88, 435, 174)$,
$(46, 260, 319, 375)$,
$(444, 137, 417, 25)$,\adfsplit
$(256, 18, 190, 452)$,
$(119, 157, 370, 427)$,
$(437, 114, 44, 397)$,
$(33, 122, 215, 54)$,\adfsplit
$(58, 326, 72, 423)$,
$(445, 252, 77, 213)$,
$(0, 3, 7, 297)$,
$(0, 5, 13, 446)$,\adfsplit
$(0, 16, 33, 411)$,
$(0, 20, 52, 255)$,
$(0, 25, 51, 165)$,
$(0, 29, 163, 198)$,\adfsplit
$(0, 30, 71, 117)$,
$(0, 49, 124, 219)$,
$(0, 55, 138, 245)$,
$(0, 34, 181, 282)$,\adfsplit
$(0, 42, 187, 289)$,
$(0, 47, 197, 277)$,
$(0, 44, 142, 357)$,
$(0, 69, 185, 315)$,\adfsplit
$(0, 43, 135, 272)$,
$(0, 31, 121, 310)$,
$(0, 54, 160, 285)$,
$(0, 74, 153, 352)$,\adfsplit
$(0, 50, 131, 250)$

\adfLgap \noindent by the mapping:
$x \mapsto x + 3 j \adfmod{456}$ for $x < 456$,
$x \mapsto (x +  j \adfmod{12}) + 456$ for $456 \le x < 468$,
$x \mapsto x$ for $x \ge 468$,
$0 \le j < 152$
 for the first five blocks;
$x \mapsto x +  j \adfmod{456}$ for $x < 456$,
$x \mapsto (x +  j \adfmod{12}) + 456$ for $456 \le x < 468$,
$x \mapsto x$ for $x \ge 468$,
$0 \le j < 456$
 for the last 36 blocks.
\ADFvfyParStart{(473, ((5, 152, ((456, 3), (12, 1), (5, 5))), (36, 456, ((456, 1), (12, 1), (5, 5)))), ((38, 12), (17, 1)))} 

\adfDgap
\noindent{\boldmath $ 38^{12} 20^{1} $}~
With the point set $Z_{476}$ partitioned into
 residue classes modulo $12$ for $\{0, 1, \dots, 455\}$, and
 $\{456, 457, \dots, 475\}$,
 the design is generated from

\adfLgap 
$(468, 0, 1, 2)$,
$(469, 0, 151, 305)$,
$(470, 0, 154, 152)$,
$(471, 0, 304, 455)$,\adfsplit
$(472, 0, 454, 302)$,
$(473, 0, 4, 11)$,
$(474, 0, 7, 452)$,
$(475, 0, 445, 449)$,\adfsplit
$(456, 343, 57, 228)$,
$(456, 6, 428, 119)$,
$(456, 339, 401, 265)$,
$(456, 82, 371, 135)$,\adfsplit
$(456, 280, 293, 181)$,
$(456, 76, 32, 358)$,
$(456, 258, 374, 165)$,
$(456, 283, 192, 26)$,\adfsplit
$(307, 75, 287, 450)$,
$(451, 129, 210, 101)$,
$(208, 176, 366, 446)$,
$(125, 31, 344, 262)$,\adfsplit
$(132, 44, 328, 403)$,
$(189, 387, 10, 174)$,
$(222, 213, 168, 455)$,
$(57, 316, 324, 99)$,\adfsplit
$(416, 71, 366, 33)$,
$(75, 386, 25, 8)$,
$(417, 209, 150, 444)$,
$(415, 242, 138, 13)$,\adfsplit
$(37, 34, 141, 242)$,
$(389, 90, 412, 202)$,
$(422, 164, 445, 283)$,
$(301, 236, 412, 372)$,\adfsplit
$(178, 203, 319, 219)$,
$(364, 257, 378, 47)$,
$(7, 419, 256, 54)$,
$(41, 434, 294, 55)$,\adfsplit
$(21, 79, 272, 442)$,
$(43, 306, 217, 86)$,
$(437, 332, 408, 454)$,
$(298, 435, 43, 349)$,\adfsplit
$(167, 414, 29, 64)$,
$(277, 189, 449, 23)$,
$(294, 67, 112, 226)$,
$(159, 196, 72, 227)$,\adfsplit
$(28, 219, 97, 432)$,
$(318, 201, 435, 276)$,
$(284, 231, 354, 142)$,
$(186, 213, 439, 166)$,\adfsplit
$(333, 206, 428, 425)$,
$(384, 256, 231, 318)$,
$(188, 115, 147, 137)$,
$(207, 309, 276, 1)$,\adfsplit
$(0, 5, 418, 440)$,
$(0, 6, 265, 321)$,
$(0, 9, 85, 167)$,
$(0, 10, 102, 160)$,\adfsplit
$(0, 19, 148, 426)$,
$(0, 26, 135, 392)$,
$(0, 18, 55, 74)$,
$(0, 39, 250, 307)$,\adfsplit
$(0, 133, 188, 443)$,
$(0, 61, 79, 194)$,
$(0, 99, 337, 427)$,
$(0, 86, 184, 395)$,\adfsplit
$(0, 100, 226, 379)$,
$(0, 157, 343, 351)$,
$(0, 94, 269, 318)$,
$(0, 110, 295, 310)$,\adfsplit
$(0, 118, 283, 353)$,
$(0, 215, 393, 399)$,
$(0, 129, 169, 329)$,
$(0, 101, 291, 417)$,\adfsplit
$(0, 311, 391, 425)$,
$(0, 97, 389, 435)$,
$(0, 83, 214, 397)$,
$(0, 28, 177, 451)$,\adfsplit
$(0, 195, 261, 371)$,
$(0, 249, 275, 373)$,
$(0, 171, 223, 301)$

\adfLgap \noindent by the mapping:
$x \mapsto x + 3 j \adfmod{456}$ for $x < 456$,
$x \mapsto (x +  j \adfmod{12}) + 456$ for $456 \le x < 468$,
$x \mapsto x$ for $x \ge 468$,
$0 \le j < 152$
 for the first eight blocks;
$x \mapsto x + 2 j \adfmod{456}$ for $x < 456$,
$x \mapsto (x +  j \adfmod{12}) + 456$ for $456 \le x < 468$,
$x \mapsto x$ for $x \ge 468$,
$0 \le j < 228$
 for the last 71 blocks.
\ADFvfyParStart{(476, ((8, 152, ((456, 3), (12, 1), (8, 8))), (71, 228, ((456, 2), (12, 1), (8, 8)))), ((38, 12), (20, 1)))} 

\adfDgap
\noindent{\boldmath $ 38^{12} 23^{1} $}~
With the point set $Z_{479}$ partitioned into
 residue classes modulo $12$ for $\{0, 1, \dots, 455\}$, and
 $\{456, 457, \dots, 478\}$,
 the design is generated from

\adfLgap 
$(468, 0, 1, 2)$,
$(469, 0, 151, 305)$,
$(470, 0, 154, 152)$,
$(471, 0, 304, 455)$,\adfsplit
$(472, 0, 454, 302)$,
$(473, 0, 4, 11)$,
$(474, 0, 7, 452)$,
$(475, 0, 445, 449)$,\adfsplit
$(476, 0, 10, 23)$,
$(477, 0, 13, 446)$,
$(478, 0, 433, 443)$,
$(456, 348, 174, 413)$,\adfsplit
$(456, 37, 104, 95)$,
$(456, 331, 22, 388)$,
$(456, 315, 62, 165)$,
$(148, 130, 49, 164)$,\adfsplit
$(60, 208, 343, 13)$,
$(394, 338, 120, 305)$,
$(448, 19, 226, 229)$,
$(187, 58, 444, 281)$,\adfsplit
$(441, 255, 178, 413)$,
$(455, 142, 21, 276)$,
$(421, 335, 208, 137)$,
$(278, 353, 94, 390)$,\adfsplit
$(29, 316, 355, 273)$,
$(69, 114, 95, 192)$,
$(0, 5, 100, 106)$,
$(0, 8, 62, 76)$,\adfsplit
$(0, 15, 32, 61)$,
$(0, 20, 41, 289)$,
$(0, 25, 63, 364)$,
$(0, 30, 191, 226)$,\adfsplit
$(0, 31, 119, 159)$,
$(0, 42, 122, 232)$,
$(0, 44, 153, 317)$,
$(0, 69, 200, 314)$,\adfsplit
$(0, 55, 188, 275)$,
$(0, 98, 205, 345)$,
$(0, 49, 165, 311)$,
$(0, 79, 170, 294)$,\adfsplit
$(0, 59, 177, 290)$,
$(0, 50, 102, 331)$,
$(0, 74, 157, 250)$,
$(0, 66, 171, 320)$,\adfsplit
$(0, 64, 137, 278)$,
$(0, 51, 104, 189)$

\adfLgap \noindent by the mapping:
$x \mapsto x + 3 j \adfmod{456}$ for $x < 456$,
$x \mapsto (x +  j \adfmod{12}) + 456$ for $456 \le x < 468$,
$x \mapsto x$ for $x \ge 468$,
$0 \le j < 152$
 for the first 11 blocks;
$x \mapsto x +  j \adfmod{456}$ for $x < 456$,
$x \mapsto (x +  j \adfmod{12}) + 456$ for $456 \le x < 468$,
$x \mapsto x$ for $x \ge 468$,
$0 \le j < 456$
 for the last 35 blocks.
\ADFvfyParStart{(479, ((11, 152, ((456, 3), (12, 1), (11, 11))), (35, 456, ((456, 1), (12, 1), (11, 11)))), ((38, 12), (23, 1)))} 

\adfDgap
\noindent{\boldmath $ 38^{15} 17^{1} $}~
With the point set $Z_{587}$ partitioned into
 residue classes modulo $15$ for $\{0, 1, \dots, 569\}$, and
 $\{570, 571, \dots, 586\}$,
 the design is generated from

\adfLgap 
$(582, 0, 1, 2)$,
$(583, 0, 190, 188)$,
$(584, 0, 379, 569)$,
$(585, 0, 382, 191)$,\adfsplit
$(586, 0, 568, 380)$,
$(570, 285, 170, 562)$,
$(570, 450, 535, 557)$,
$(571, 189, 58, 145)$,\adfsplit
$(571, 449, 120, 32)$,
$(572, 420, 374, 475)$,
$(572, 358, 441, 383)$,
$(573, 243, 511, 274)$,\adfsplit
$(573, 419, 386, 462)$,
$(179, 353, 567, 415)$,
$(46, 65, 494, 206)$,
$(294, 120, 268, 87)$,\adfsplit
$(288, 169, 96, 16)$,
$(252, 146, 403, 361)$,
$(485, 413, 206, 187)$,
$(406, 485, 317, 324)$,\adfsplit
$(335, 216, 360, 131)$,
$(345, 242, 129, 38)$,
$(125, 236, 552, 188)$,
$(90, 182, 446, 144)$,\adfsplit
$(176, 332, 481, 323)$,
$(508, 479, 429, 307)$,
$(550, 219, 395, 44)$,
$(30, 406, 8, 25)$,\adfsplit
$(171, 31, 490, 315)$,
$(353, 435, 469, 501)$,
$(490, 472, 327, 314)$,
$(76, 81, 407, 514)$,\adfsplit
$(436, 201, 104, 562)$,
$(353, 245, 176, 468)$,
$(313, 352, 362, 519)$,
$(145, 549, 282, 278)$,\adfsplit
$(323, 40, 441, 209)$,
$(396, 462, 420, 91)$,
$(469, 9, 440, 222)$,
$(390, 289, 560, 191)$,\adfsplit
$(149, 137, 394, 266)$,
$(464, 24, 183, 94)$,
$(306, 373, 468, 121)$,
$(446, 138, 336, 457)$,\adfsplit
$(291, 365, 522, 248)$,
$(164, 400, 260, 106)$,
$(198, 158, 181, 516)$,
$(146, 197, 232, 43)$,\adfsplit
$(275, 539, 236, 557)$,
$(509, 405, 317, 205)$,
$(329, 490, 288, 534)$,
$(15, 366, 422, 545)$,\adfsplit
$(141, 188, 384, 70)$,
$(442, 514, 387, 189)$,
$(435, 73, 144, 487)$,
$(380, 37, 255, 171)$,\adfsplit
$(311, 211, 51, 248)$,
$(422, 472, 299, 118)$,
$(294, 110, 537, 262)$,
$(307, 20, 330, 69)$,\adfsplit
$(461, 290, 368, 514)$,
$(218, 342, 134, 3)$,
$(28, 14, 517, 482)$,
$(389, 380, 475, 159)$,\adfsplit
$(366, 57, 144, 374)$,
$(391, 488, 429, 550)$,
$(386, 515, 128, 511)$,
$(0, 12, 28, 59)$,\adfsplit
$(0, 3, 6, 344)$,
$(0, 21, 57, 406)$,
$(0, 7, 34, 61)$,
$(0, 65, 432, 509)$,\adfsplit
$(0, 68, 142, 470)$,
$(0, 36, 134, 320)$,
$(0, 52, 197, 456)$,
$(0, 99, 263, 336)$,\adfsplit
$(0, 104, 217, 297)$,
$(0, 108, 367, 388)$,
$(0, 20, 301, 342)$,
$(0, 233, 244, 513)$,\adfsplit
$(0, 136, 323, 557)$,
$(0, 94, 431, 487)$,
$(0, 53, 220, 493)$,
$(0, 173, 249, 419)$,\adfsplit
$(0, 133, 211, 353)$,
$(0, 213, 391, 399)$,
$(0, 81, 205, 519)$,
$(0, 127, 155, 501)$,\adfsplit
$(0, 183, 461, 477)$,
$(0, 377, 479, 505)$,
$(0, 37, 295, 359)$,
$(0, 277, 423, 471)$,\adfsplit
$(1, 11, 107, 369)$,
$(1, 7, 21, 349)$,
$(1, 47, 185, 321)$,
$(1, 25, 93, 501)$

\adfLgap \noindent by the mapping:
$x \mapsto x + 3 j \adfmod{570}$ for $x < 570$,
$x \mapsto (x - 570 + 4 j \adfmod{12}) + 570$ for $570 \le x < 582$,
$x \mapsto x$ for $x \ge 582$,
$0 \le j < 190$
 for the first five blocks;
$x \mapsto x + 2 j \adfmod{570}$ for $x < 570$,
$x \mapsto (x - 570 + 4 j \adfmod{12}) + 570$ for $570 \le x < 582$,
$x \mapsto x$ for $x \ge 582$,
$0 \le j < 285$
 for the last 91 blocks.
\ADFvfyParStart{(587, ((5, 190, ((570, 3), (12, 4), (5, 5))), (91, 285, ((570, 2), (12, 4), (5, 5)))), ((38, 15), (17, 1)))} 

\adfDgap
\noindent{\boldmath $ 38^{15} 20^{1} $}~
With the point set $Z_{590}$ partitioned into
 residue classes modulo $15$ for $\{0, 1, \dots, 569\}$, and
 $\{570, 571, \dots, 589\}$,
 the design is generated from

\adfLgap 
$(585, 0, 1, 2)$,
$(586, 0, 190, 188)$,
$(587, 0, 379, 569)$,
$(588, 0, 382, 191)$,\adfsplit
$(589, 0, 568, 380)$,
$(570, 38, 183, 169)$,
$(570, 400, 232, 376)$,
$(570, 569, 283, 279)$,\adfsplit
$(570, 225, 177, 272)$,
$(570, 426, 401, 80)$,
$(434, 103, 219, 421)$,
$(457, 391, 423, 537)$,\adfsplit
$(489, 386, 93, 235)$,
$(138, 553, 234, 23)$,
$(328, 436, 344, 4)$,
$(551, 37, 375, 568)$,\adfsplit
$(21, 151, 277, 105)$,
$(260, 478, 473, 555)$,
$(268, 460, 434, 215)$,
$(244, 330, 36, 371)$,\adfsplit
$(508, 410, 205, 438)$,
$(530, 160, 247, 321)$,
$(395, 191, 154, 134)$,
$(333, 150, 534, 86)$,\adfsplit
$(138, 509, 266, 451)$,
$(105, 527, 128, 491)$,
$(185, 495, 261, 167)$,
$(254, 351, 145, 245)$,\adfsplit
$(41, 413, 91, 190)$,
$(121, 347, 14, 189)$,
$(13, 466, 414, 554)$,
$(27, 319, 46, 15)$,\adfsplit
$(0, 3, 121, 262)$,
$(0, 6, 49, 437)$,
$(0, 7, 78, 167)$,
$(0, 8, 81, 91)$,\adfsplit
$(0, 11, 33, 302)$,
$(0, 21, 59, 317)$,
$(0, 44, 123, 457)$,
$(0, 61, 154, 381)$,\adfsplit
$(0, 65, 137, 418)$,
$(0, 51, 194, 350)$,
$(0, 27, 159, 214)$,
$(0, 104, 228, 392)$,\adfsplit
$(0, 62, 163, 265)$,
$(0, 42, 212, 348)$,
$(0, 85, 196, 417)$,
$(0, 35, 147, 216)$,\adfsplit
$(0, 63, 197, 326)$,
$(0, 54, 173, 445)$,
$(0, 67, 177, 408)$

\adfLgap \noindent by the mapping:
$x \mapsto x + 3 j \adfmod{570}$ for $x < 570$,
$x \mapsto (x +  j \adfmod{15}) + 570$ for $570 \le x < 585$,
$x \mapsto x$ for $x \ge 585$,
$0 \le j < 190$
 for the first five blocks;
$x \mapsto x +  j \adfmod{570}$ for $x < 570$,
$x \mapsto (x +  j \adfmod{15}) + 570$ for $570 \le x < 585$,
$x \mapsto x$ for $x \ge 585$,
$0 \le j < 570$
 for the last 46 blocks.
\ADFvfyParStart{(590, ((5, 190, ((570, 3), (15, 1), (5, 5))), (46, 570, ((570, 1), (15, 1), (5, 5)))), ((38, 15), (20, 1)))} 

\adfDgap
\noindent{\boldmath $ 38^{15} 23^{1} $}~
With the point set $Z_{593}$ partitioned into
 residue classes modulo $15$ for $\{0, 1, \dots, 569\}$, and
 $\{570, 571, \dots, 592\}$,
 the design is generated from

\adfLgap 
$(585, 0, 1, 2)$,
$(586, 0, 190, 188)$,
$(587, 0, 379, 569)$,
$(588, 0, 382, 191)$,\adfsplit
$(589, 0, 568, 380)$,
$(590, 0, 4, 11)$,
$(591, 0, 7, 566)$,
$(592, 0, 559, 563)$,\adfsplit
$(570, 20, 278, 348)$,
$(570, 15, 246, 54)$,
$(570, 70, 146, 557)$,
$(570, 543, 510, 242)$,\adfsplit
$(570, 471, 252, 493)$,
$(570, 404, 207, 39)$,
$(570, 148, 359, 71)$,
$(570, 307, 481, 563)$,\adfsplit
$(570, 139, 406, 95)$,
$(570, 412, 4, 55)$,
$(115, 457, 199, 477)$,
$(388, 161, 222, 347)$,\adfsplit
$(498, 57, 191, 239)$,
$(260, 40, 19, 386)$,
$(138, 411, 296, 43)$,
$(215, 498, 107, 75)$,\adfsplit
$(235, 326, 380, 486)$,
$(148, 539, 504, 242)$,
$(4, 465, 141, 128)$,
$(315, 36, 157, 429)$,\adfsplit
$(226, 119, 289, 216)$,
$(334, 501, 163, 286)$,
$(98, 117, 32, 463)$,
$(527, 530, 463, 267)$,\adfsplit
$(560, 160, 387, 256)$,
$(85, 426, 337, 299)$,
$(33, 50, 486, 25)$,
$(311, 228, 280, 233)$,\adfsplit
$(50, 14, 8, 516)$,
$(381, 181, 475, 124)$,
$(323, 439, 532, 62)$,
$(43, 7, 521, 420)$,\adfsplit
$(270, 162, 5, 206)$,
$(473, 52, 72, 261)$,
$(396, 380, 248, 445)$,
$(213, 146, 72, 113)$,\adfsplit
$(332, 509, 54, 493)$,
$(435, 63, 398, 325)$,
$(363, 360, 80, 506)$,
$(137, 331, 254, 156)$,\adfsplit
$(15, 72, 276, 294)$,
$(235, 50, 304, 15)$,
$(275, 451, 539, 396)$,
$(456, 354, 235, 473)$,\adfsplit
$(392, 378, 193, 24)$,
$(483, 156, 147, 217)$,
$(558, 238, 10, 416)$,
$(18, 100, 68, 432)$,\adfsplit
$(126, 408, 39, 45)$,
$(50, 371, 249, 317)$,
$(478, 364, 487, 20)$,
$(215, 146, 186, 61)$,\adfsplit
$(280, 272, 200, 24)$,
$(204, 353, 120, 455)$,
$(257, 356, 243, 394)$,
$(83, 238, 120, 55)$,\adfsplit
$(446, 140, 492, 174)$,
$(416, 527, 141, 367)$,
$(188, 145, 66, 362)$,
$(129, 95, 225, 276)$,\adfsplit
$(203, 22, 332, 216)$,
$(0, 23, 372, 460)$,
$(0, 21, 208, 492)$,
$(0, 27, 58, 86)$,\adfsplit
$(0, 43, 154, 336)$,
$(0, 12, 59, 442)$,
$(0, 24, 95, 196)$,
$(0, 39, 340, 437)$,\adfsplit
$(0, 25, 173, 236)$,
$(0, 138, 370, 541)$,
$(0, 93, 246, 409)$,
$(0, 81, 152, 243)$,\adfsplit
$(0, 151, 326, 565)$,
$(0, 145, 155, 402)$,
$(0, 92, 245, 386)$,
$(0, 26, 357, 384)$,\adfsplit
$(0, 53, 217, 440)$,
$(0, 56, 293, 355)$,
$(0, 103, 247, 333)$,
$(0, 127, 253, 265)$,\adfsplit
$(0, 363, 491, 515)$,
$(0, 89, 325, 547)$,
$(0, 107, 249, 381)$,
$(0, 136, 367, 407)$,\adfsplit
$(0, 87, 113, 159)$,
$(0, 139, 275, 517)$,
$(1, 43, 147, 365)$,
$(0, 223, 289, 341)$,\adfsplit
$(1, 51, 125, 405)$,
$(1, 19, 77, 321)$,
$(1, 81, 179, 285)$

\adfLgap \noindent by the mapping:
$x \mapsto x + 3 j \adfmod{570}$ for $x < 570$,
$x \mapsto (x +  j \adfmod{15}) + 570$ for $570 \le x < 585$,
$x \mapsto x$ for $x \ge 585$,
$0 \le j < 190$
 for the first eight blocks;
$x \mapsto x + 2 j \adfmod{570}$ for $x < 570$,
$x \mapsto (x +  j \adfmod{15}) + 570$ for $570 \le x < 585$,
$x \mapsto x$ for $x \ge 585$,
$0 \le j < 285$
 for the last 91 blocks.
\ADFvfyParStart{(593, ((8, 190, ((570, 3), (15, 1), (8, 8))), (91, 285, ((570, 2), (15, 1), (8, 8)))), ((38, 15), (23, 1)))} 

\adfDgap
\noindent{\boldmath $ 38^{15} 26^{1} $}~
With the point set $Z_{596}$ partitioned into
 residue classes modulo $15$ for $\{0, 1, \dots, 569\}$, and
 $\{570, 571, \dots, 595\}$,
 the design is generated from

\adfLgap 
$(585, 0, 1, 2)$,
$(586, 0, 190, 188)$,
$(587, 0, 379, 569)$,
$(588, 0, 382, 191)$,\adfsplit
$(589, 0, 568, 380)$,
$(590, 0, 4, 11)$,
$(591, 0, 7, 566)$,
$(592, 0, 559, 563)$,\adfsplit
$(593, 0, 10, 23)$,
$(594, 0, 13, 560)$,
$(595, 0, 547, 557)$,
$(570, 452, 313, 284)$,\adfsplit
$(570, 271, 188, 372)$,
$(570, 381, 469, 265)$,
$(570, 322, 56, 545)$,
$(570, 198, 495, 159)$,\adfsplit
$(263, 236, 550, 421)$,
$(201, 478, 175, 443)$,
$(310, 521, 186, 497)$,
$(52, 231, 132, 428)$,\adfsplit
$(280, 517, 53, 228)$,
$(296, 380, 82, 238)$,
$(541, 560, 253, 108)$,
$(320, 134, 154, 498)$,\adfsplit
$(177, 165, 269, 100)$,
$(462, 90, 355, 49)$,
$(141, 302, 204, 50)$,
$(52, 566, 513, 113)$,\adfsplit
$(3, 20, 345, 263)$,
$(72, 429, 363, 319)$,
$(362, 308, 170, 232)$,
$(441, 479, 352, 128)$,\adfsplit
$(553, 503, 47, 254)$,
$(146, 49, 216, 265)$,
$(325, 319, 466, 277)$,
$(562, 23, 95, 460)$,\adfsplit
$(448, 402, 41, 533)$,
$(0, 3, 8, 239)$,
$(0, 9, 25, 484)$,
$(0, 18, 144, 436)$,\adfsplit
$(0, 14, 69, 136)$,
$(0, 21, 43, 410)$,
$(0, 28, 115, 447)$,
$(0, 47, 202, 295)$,\adfsplit
$(0, 74, 174, 350)$,
$(0, 71, 244, 387)$,
$(0, 33, 146, 286)$,
$(0, 37, 162, 377)$,\adfsplit
$(0, 34, 128, 324)$,
$(0, 73, 222, 393)$,
$(0, 96, 208, 449)$,
$(0, 32, 164, 232)$,\adfsplit
$(0, 40, 148, 349)$,
$(0, 57, 229, 308)$,
$(0, 36, 218, 269)$,
$(0, 59, 212, 371)$

\adfLgap \noindent by the mapping:
$x \mapsto x + 3 j \adfmod{570}$ for $x < 570$,
$x \mapsto (x +  j \adfmod{15}) + 570$ for $570 \le x < 585$,
$x \mapsto x$ for $x \ge 585$,
$0 \le j < 190$
 for the first 11 blocks;
$x \mapsto x +  j \adfmod{570}$ for $x < 570$,
$x \mapsto (x +  j \adfmod{15}) + 570$ for $570 \le x < 585$,
$x \mapsto x$ for $x \ge 585$,
$0 \le j < 570$
 for the last 45 blocks.
\ADFvfyParStart{(596, ((11, 190, ((570, 3), (15, 1), (11, 11))), (45, 570, ((570, 1), (15, 1), (11, 11)))), ((38, 15), (26, 1)))} 

\adfDgap
\noindent{\boldmath $ 38^{15} 29^{1} $}~
With the point set $Z_{599}$ partitioned into
 residue classes modulo $15$ for $\{0, 1, \dots, 569\}$, and
 $\{570, 571, \dots, 598\}$,
 the design is generated from

\adfLgap 
$(585, 0, 1, 2)$,
$(586, 0, 190, 188)$,
$(587, 0, 379, 569)$,
$(588, 0, 382, 191)$,\adfsplit
$(589, 0, 568, 380)$,
$(590, 0, 4, 11)$,
$(591, 0, 7, 566)$,
$(592, 0, 559, 563)$,\adfsplit
$(593, 0, 10, 23)$,
$(594, 0, 13, 560)$,
$(595, 0, 547, 557)$,
$(596, 0, 16, 35)$,\adfsplit
$(597, 0, 19, 554)$,
$(598, 0, 535, 551)$,
$(570, 549, 479, 488)$,
$(570, 484, 522, 280)$,\adfsplit
$(570, 352, 411, 242)$,
$(570, 127, 326, 208)$,
$(570, 503, 153, 527)$,
$(570, 20, 237, 343)$,\adfsplit
$(570, 546, 331, 349)$,
$(570, 105, 286, 408)$,
$(570, 251, 420, 295)$,
$(570, 444, 164, 365)$,\adfsplit
$(84, 225, 559, 348)$,
$(458, 189, 201, 322)$,
$(350, 442, 479, 484)$,
$(252, 184, 403, 67)$,\adfsplit
$(37, 279, 327, 347)$,
$(385, 461, 549, 240)$,
$(260, 419, 6, 167)$,
$(248, 39, 559, 151)$,\adfsplit
$(334, 11, 403, 466)$,
$(485, 430, 511, 477)$,
$(518, 527, 97, 252)$,
$(475, 29, 87, 350)$,\adfsplit
$(231, 562, 329, 550)$,
$(201, 196, 159, 268)$,
$(201, 8, 450, 521)$,
$(111, 446, 419, 263)$,\adfsplit
$(98, 9, 274, 310)$,
$(147, 334, 443, 276)$,
$(56, 535, 421, 158)$,
$(272, 455, 43, 510)$,\adfsplit
$(380, 564, 336, 150)$,
$(487, 324, 423, 19)$,
$(568, 63, 171, 25)$,
$(123, 522, 351, 22)$,\adfsplit
$(146, 392, 435, 441)$,
$(307, 3, 85, 64)$,
$(559, 399, 196, 131)$,
$(567, 430, 483, 191)$,\adfsplit
$(462, 188, 171, 406)$,
$(41, 204, 472, 555)$,
$(481, 435, 138, 343)$,
$(391, 558, 64, 219)$,\adfsplit
$(286, 113, 217, 264)$,
$(131, 301, 384, 205)$,
$(290, 164, 523, 379)$,
$(234, 72, 568, 365)$,\adfsplit
$(23, 561, 209, 320)$,
$(75, 558, 264, 524)$,
$(55, 498, 281, 539)$,
$(543, 261, 515, 436)$,\adfsplit
$(417, 23, 269, 95)$,
$(91, 164, 482, 568)$,
$(177, 466, 379, 415)$,
$(27, 157, 443, 140)$,\adfsplit
$(317, 370, 200, 56)$,
$(415, 240, 412, 76)$,
$(497, 63, 267, 511)$,
$(346, 287, 60, 34)$,\adfsplit
$(523, 413, 391, 204)$,
$(91, 131, 238, 464)$,
$(0, 14, 77, 439)$,
$(0, 6, 153, 417)$,\adfsplit
$(0, 24, 283, 467)$,
$(0, 8, 33, 549)$,
$(0, 28, 421, 499)$,
$(0, 29, 32, 245)$,\adfsplit
$(0, 31, 97, 539)$,
$(0, 51, 177, 257)$,
$(0, 95, 287, 388)$,
$(0, 199, 333, 545)$,\adfsplit
$(0, 157, 357, 451)$,
$(0, 229, 427, 527)$,
$(0, 91, 206, 429)$,
$(0, 88, 207, 521)$,\adfsplit
$(0, 369, 485, 537)$,
$(0, 67, 185, 346)$,
$(0, 48, 133, 347)$,
$(0, 40, 113, 152)$,\adfsplit
$(0, 78, 189, 488)$,
$(0, 123, 200, 368)$,
$(0, 39, 178, 466)$,
$(0, 57, 98, 446)$,\adfsplit
$(0, 54, 138, 430)$,
$(0, 18, 192, 338)$,
$(0, 64, 130, 272)$,
$(0, 52, 114, 464)$,\adfsplit
$(0, 94, 248, 356)$,
$(0, 46, 96, 244)$,
$(0, 20, 100, 216)$

\adfLgap \noindent by the mapping:
$x \mapsto x + 3 j \adfmod{570}$ for $x < 570$,
$x \mapsto (x +  j \adfmod{15}) + 570$ for $570 \le x < 585$,
$x \mapsto x$ for $x \ge 585$,
$0 \le j < 190$
 for the first 14 blocks;
$x \mapsto x + 2 j \adfmod{570}$ for $x < 570$,
$x \mapsto (x +  j \adfmod{15}) + 570$ for $570 \le x < 585$,
$x \mapsto x$ for $x \ge 585$,
$0 \le j < 285$
 for the last 89 blocks.
\ADFvfyParStart{(599, ((14, 190, ((570, 3), (15, 1), (14, 14))), (89, 285, ((570, 2), (15, 1), (14, 14)))), ((38, 15), (29, 1)))} 

\adfDgap
\noindent{\boldmath $ 38^{18} 20^{1} $}~
With the point set $Z_{704}$ partitioned into
 residue classes modulo $18$ for $\{0, 1, \dots, 683\}$, and
 $\{684, 685, \dots, 703\}$,
 the design is generated from

\adfLgap 
$(702, 0, 1, 2)$,
$(703, 0, 229, 458)$,
$(684, 367, 484, 490)$,
$(684, 138, 506, 430)$,\adfsplit
$(684, 247, 501, 3)$,
$(684, 658, 549, 23)$,
$(684, 288, 568, 89)$,
$(684, 428, 112, 636)$,\adfsplit
$(684, 410, 663, 353)$,
$(684, 186, 487, 365)$,
$(684, 193, 467, 660)$,
$(684, 450, 443, 61)$,\adfsplit
$(684, 524, 116, 165)$,
$(684, 423, 397, 422)$,
$(653, 244, 462, 88)$,
$(309, 481, 171, 446)$,\adfsplit
$(59, 278, 645, 615)$,
$(419, 184, 396, 235)$,
$(70, 1, 312, 67)$,
$(161, 250, 67, 6)$,\adfsplit
$(620, 265, 311, 394)$,
$(486, 595, 395, 15)$,
$(49, 369, 493, 304)$,
$(394, 544, 463, 347)$,\adfsplit
$(66, 108, 191, 320)$,
$(67, 245, 469, 178)$,
$(476, 653, 97, 517)$,
$(356, 20, 15, 301)$,\adfsplit
$(216, 244, 439, 572)$,
$(43, 184, 182, 27)$,
$(56, 377, 474, 553)$,
$(674, 517, 228, 136)$,\adfsplit
$(390, 55, 470, 303)$,
$(545, 158, 665, 480)$,
$(482, 672, 147, 632)$,
$(120, 155, 96, 334)$,\adfsplit
$(533, 394, 307, 222)$,
$(45, 122, 12, 257)$,
$(350, 203, 363, 258)$,
$(135, 267, 403, 95)$,\adfsplit
$(215, 227, 202, 678)$,
$(368, 353, 361, 670)$,
$(366, 281, 567, 107)$,
$(227, 372, 329, 32)$,\adfsplit
$(647, 216, 679, 639)$,
$(303, 647, 261, 632)$,
$(534, 640, 544, 376)$,
$(659, 140, 293, 25)$,\adfsplit
$(371, 296, 629, 96)$,
$(660, 511, 212, 263)$,
$(551, 631, 121, 485)$,
$(88, 656, 353, 51)$,\adfsplit
$(226, 463, 530, 173)$,
$(612, 91, 118, 209)$,
$(123, 444, 679, 400)$,
$(474, 200, 40, 68)$,\adfsplit
$(300, 657, 565, 232)$,
$(363, 569, 585, 618)$,
$(84, 214, 311, 666)$,
$(55, 17, 465, 535)$,\adfsplit
$(564, 117, 100, 518)$,
$(224, 485, 149, 138)$,
$(309, 331, 624, 468)$,
$(553, 433, 526, 314)$,\adfsplit
$(640, 96, 530, 261)$,
$(541, 171, 85, 180)$,
$(234, 153, 233, 86)$,
$(31, 84, 278, 70)$,\adfsplit
$(18, 117, 156, 623)$,
$(461, 482, 361, 648)$,
$(177, 133, 395, 634)$,
$(407, 246, 73, 427)$,\adfsplit
$(418, 344, 81, 330)$,
$(121, 118, 178, 613)$,
$(213, 539, 463, 62)$,
$(161, 397, 457, 172)$,\adfsplit
$(54, 272, 158, 676)$,
$(90, 272, 381, 568)$,
$(245, 51, 75, 346)$,
$(472, 80, 558, 579)$,\adfsplit
$(467, 224, 571, 426)$,
$(299, 537, 638, 244)$,
$(626, 42, 453, 313)$,
$(109, 331, 200, 368)$,\adfsplit
$(154, 284, 674, 307)$,
$(0, 9, 127, 625)$,
$(0, 20, 50, 474)$,
$(0, 12, 438, 483)$,\adfsplit
$(0, 22, 239, 656)$,
$(0, 128, 298, 512)$,
$(0, 17, 188, 585)$,
$(0, 74, 170, 389)$,\adfsplit
$(0, 124, 406, 529)$,
$(0, 19, 48, 116)$,
$(0, 6, 496, 606)$,
$(0, 13, 207, 454)$,\adfsplit
$(0, 8, 312, 575)$,
$(0, 64, 199, 503)$,
$(0, 73, 171, 579)$,
$(0, 63, 196, 541)$,\adfsplit
$(0, 285, 545, 623)$,
$(0, 58, 543, 658)$,
$(0, 225, 273, 318)$,
$(0, 203, 231, 621)$,\adfsplit
$(0, 40, 455, 573)$,
$(0, 34, 137, 365)$,
$(0, 322, 449, 453)$,
$(0, 142, 353, 653)$,\adfsplit
$(0, 148, 341, 628)$,
$(0, 335, 417, 581)$,
$(3, 64, 256, 315)$,
$(3, 47, 287, 358)$,\adfsplit
$(3, 41, 226, 371)$,
$(3, 76, 155, 477)$,
$(3, 29, 143, 172)$

\adfLgap \noindent by the mapping:
$x \mapsto x \oplus (3 j)$ for $x < 684$,
$x \mapsto (x +  j \adfmod{18}) + 684$ for $684 \le x < 702$,
$x \mapsto x$ for $x \ge 702$,
$0 \le j < 228$
 for the first two blocks;
$x \mapsto x \oplus j \oplus j$ for $x < 684$,
$x \mapsto (x +  j \adfmod{18}) + 684$ for $684 \le x < 702$,
$x \mapsto x$ for $x \ge 702$,
$0 \le j < 342$
 for the last 113 blocks.
\ADFvfyParStart{(704, ((2, 228, ((684, 3, (228, 3)), (18, 1), (2, 2))), (113, 342, ((684, 2, (228, 3)), (18, 1), (2, 2)))), ((38, 18), (20, 1)))} 

\adfDgap
\noindent{\boldmath $ 38^{18} 23^{1} $}~
With the point set $Z_{707}$ partitioned into
 residue classes modulo $18$ for $\{0, 1, \dots, 683\}$, and
 $\{684, 685, \dots, 706\}$,
 the design is generated from

\adfLgap 
$(705, 0, 1, 2)$,
$(706, 0, 229, 458)$,
$(684, 220, 224, 153)$,
$(685, 3, 590, 142)$,\adfsplit
$(686, 609, 272, 127)$,
$(687, 394, 387, 347)$,
$(688, 621, 115, 335)$,
$(689, 171, 584, 40)$,\adfsplit
$(690, 137, 522, 361)$,
$(301, 141, 287, 403)$,
$(394, 260, 419, 472)$,
$(16, 625, 364, 570)$,\adfsplit
$(512, 597, 12, 194)$,
$(337, 353, 584, 325)$,
$(485, 301, 412, 677)$,
$(143, 316, 375, 77)$,\adfsplit
$(113, 603, 520, 245)$,
$(98, 380, 89, 417)$,
$(82, 552, 337, 105)$,
$(14, 330, 647, 624)$,\adfsplit
$(258, 102, 371, 340)$,
$(652, 430, 583, 389)$,
$(230, 5, 624, 483)$,
$(584, 202, 145, 299)$,\adfsplit
$(312, 345, 15, 364)$,
$(634, 168, 10, 319)$,
$(132, 311, 466, 318)$,
$(297, 380, 97, 600)$,\adfsplit
$(84, 126, 231, 675)$,
$(193, 606, 29, 294)$,
$(658, 238, 12, 581)$,
$(454, 325, 128, 291)$,\adfsplit
$(369, 427, 68, 610)$,
$(669, 340, 624, 277)$,
$(370, 193, 43, 414)$,
$(27, 160, 65, 637)$,\adfsplit
$(206, 319, 223, 408)$,
$(529, 15, 623, 422)$,
$(564, 213, 88, 4)$,
$(0, 3, 11, 253)$,\adfsplit
$(0, 5, 32, 53)$,
$(0, 6, 196, 257)$,
$(0, 14, 101, 120)$,
$(0, 13, 127, 230)$,\adfsplit
$(0, 24, 260, 352)$,
$(0, 68, 272, 389)$,
$(0, 39, 175, 513)$,
$(0, 49, 244, 300)$,\adfsplit
$(0, 30, 249, 438)$,
$(0, 81, 233, 401)$,
$(0, 115, 243, 518)$,
$(0, 29, 211, 296)$,\adfsplit
$(0, 46, 213, 385)$,
$(0, 149, 314, 510)$,
$(0, 37, 123, 316)$,
$(0, 70, 258, 380)$,\adfsplit
$(0, 15, 121, 394)$,
$(0, 35, 157, 295)$,
$(0, 64, 140, 596)$

\adfLgap \noindent by the mapping:
$x \mapsto x \oplus (3 j)$ for $x < 684$,
$x \mapsto (x - 684 + 7 j \adfmod{21}) + 684$ for $684 \le x < 705$,
$x \mapsto x$ for $x \ge 705$,
$0 \le j < 228$
 for the first two blocks;
$x \mapsto x \oplus j$ for $x < 684$,
$x \mapsto (x - 684 + 7 j \adfmod{21}) + 684$ for $684 \le x < 705$,
$x \mapsto x$ for $x \ge 705$,
$0 \le j < 684$
 for the last 57 blocks.
\ADFvfyParStart{(707, ((2, 228, ((684, 3, (228, 3)), (21, 7), (2, 2))), (57, 684, ((684, 1, (228, 3)), (21, 7), (2, 2)))), ((38, 18), (23, 1)))} 

\adfDgap
\noindent{\boldmath $ 38^{18} 26^{1} $}~
With the point set $Z_{710}$ partitioned into
 residue classes modulo $18$ for $\{0, 1, \dots, 683\}$, and
 $\{684, 685, \dots, 709\}$,
 the design is generated from

\adfLgap 
$(708, 0, 1, 2)$,
$(709, 0, 229, 458)$,
$(684, 622, 615, 470)$,
$(684, 54, 172, 120)$,\adfsplit
$(684, 525, 587, 487)$,
$(684, 341, 320, 325)$,
$(685, 122, 419, 583)$,
$(685, 609, 64, 510)$,\adfsplit
$(685, 289, 420, 579)$,
$(685, 32, 286, 581)$,
$(686, 135, 609, 264)$,
$(686, 53, 535, 208)$,\adfsplit
$(686, 37, 478, 407)$,
$(686, 416, 426, 422)$,
$(687, 358, 313, 668)$,
$(687, 14, 53, 567)$,\adfsplit
$(687, 546, 573, 355)$,
$(687, 664, 276, 371)$,
$(70, 4, 128, 519)$,
$(456, 507, 23, 358)$,\adfsplit
$(322, 410, 444, 407)$,
$(20, 515, 90, 93)$,
$(356, 431, 64, 205)$,
$(4, 489, 522, 282)$,\adfsplit
$(187, 398, 671, 663)$,
$(296, 218, 23, 298)$,
$(353, 665, 440, 636)$,
$(87, 272, 72, 489)$,\adfsplit
$(397, 133, 147, 598)$,
$(359, 423, 486, 379)$,
$(244, 497, 384, 235)$,
$(114, 327, 383, 591)$,\adfsplit
$(174, 32, 261, 363)$,
$(588, 322, 658, 418)$,
$(4, 135, 235, 576)$,
$(543, 413, 105, 97)$,\adfsplit
$(441, 618, 482, 559)$,
$(392, 563, 38, 18)$,
$(468, 590, 169, 659)$,
$(542, 251, 291, 511)$,\adfsplit
$(150, 479, 148, 20)$,
$(225, 514, 482, 472)$,
$(396, 408, 675, 543)$,
$(241, 609, 20, 541)$,\adfsplit
$(86, 525, 357, 253)$,
$(642, 96, 674, 562)$,
$(298, 250, 511, 120)$,
$(582, 532, 157, 213)$,\adfsplit
$(615, 636, 210, 263)$,
$(133, 480, 377, 122)$,
$(303, 291, 474, 61)$,
$(184, 233, 291, 478)$,\adfsplit
$(555, 274, 430, 498)$,
$(188, 510, 532, 218)$,
$(665, 423, 138, 347)$,
$(226, 207, 462, 272)$,\adfsplit
$(576, 464, 545, 67)$,
$(426, 276, 90, 251)$,
$(215, 205, 201, 558)$,
$(294, 408, 245, 91)$,\adfsplit
$(558, 519, 435, 464)$,
$(179, 192, 605, 382)$,
$(30, 487, 536, 335)$,
$(510, 329, 629, 266)$,\adfsplit
$(218, 413, 1, 90)$,
$(586, 637, 106, 553)$,
$(325, 581, 600, 40)$,
$(421, 605, 58, 72)$,\adfsplit
$(540, 379, 318, 492)$,
$(474, 358, 51, 275)$,
$(192, 501, 479, 590)$,
$(12, 217, 303, 272)$,\adfsplit
$(407, 61, 476, 660)$,
$(251, 316, 502, 131)$,
$(578, 676, 417, 554)$,
$(601, 364, 560, 170)$,\adfsplit
$(242, 151, 262, 427)$,
$(285, 25, 16, 390)$,
$(56, 118, 221, 619)$,
$(480, 231, 136, 348)$,\adfsplit
$(80, 23, 466, 683)$,
$(404, 301, 75, 640)$,
$(2, 507, 582, 233)$,
$(172, 353, 180, 236)$,\adfsplit
$(222, 540, 371, 157)$,
$(386, 229, 457, 562)$,
$(0, 25, 582, 610)$,
$(0, 4, 284, 307)$,\adfsplit
$(0, 38, 98, 474)$,
$(0, 40, 373, 516)$,
$(0, 46, 588, 652)$,
$(0, 82, 192, 312)$,\adfsplit
$(0, 50, 480, 640)$,
$(0, 176, 422, 669)$,
$(0, 19, 196, 219)$,
$(0, 42, 182, 364)$,\adfsplit
$(0, 67, 141, 469)$,
$(0, 136, 404, 628)$,
$(0, 86, 382, 529)$,
$(0, 94, 242, 475)$,\adfsplit
$(0, 79, 455, 528)$,
$(0, 43, 292, 591)$,
$(0, 129, 304, 518)$,
$(0, 88, 123, 215)$,\adfsplit
$(0, 277, 605, 639)$,
$(0, 353, 383, 400)$,
$(0, 281, 387, 538)$,
$(0, 167, 580, 657)$,\adfsplit
$(0, 143, 172, 497)$,
$(0, 93, 106, 369)$,
$(0, 280, 443, 617)$,
$(0, 137, 352, 531)$,\adfsplit
$(3, 47, 117, 574)$,
$(3, 9, 100, 575)$,
$(3, 63, 167, 389)$,
$(3, 58, 196, 537)$,\adfsplit
$(3, 34, 179, 407)$

\adfLgap \noindent by the mapping:
$x \mapsto x \oplus (3 j)$ for $x < 684$,
$x \mapsto (x - 684 + 4 j \adfmod{24}) + 684$ for $684 \le x < 708$,
$x \mapsto x$ for $x \ge 708$,
$0 \le j < 228$
 for the first two blocks;
$x \mapsto x \oplus j \oplus j$ for $x < 684$,
$x \mapsto (x - 684 + 4 j \adfmod{24}) + 684$ for $684 \le x < 708$,
$x \mapsto x$ for $x \ge 708$,
$0 \le j < 342$
 for the last 115 blocks.
\ADFvfyParStart{(710, ((2, 228, ((684, 3, (228, 3)), (24, 4), (2, 2))), (115, 342, ((684, 2, (228, 3)), (24, 4), (2, 2)))), ((38, 18), (26, 1)))} 

\adfDgap
\noindent{\boldmath $ 38^{18} 29^{1} $}~
With the point set $Z_{713}$ partitioned into
 residue classes modulo $18$ for $\{0, 1, \dots, 683\}$, and
 $\{684, 685, \dots, 712\}$,
 the design is generated from

\adfLgap 
$(711, 0, 1, 2)$,
$(712, 0, 229, 458)$,
$(684, 603, 464, 413)$,
$(684, 625, 457, 471)$,\adfsplit
$(684, 352, 69, 182)$,
$(685, 612, 122, 640)$,
$(685, 159, 664, 472)$,
$(685, 155, 260, 354)$,\adfsplit
$(686, 678, 528, 340)$,
$(686, 9, 677, 23)$,
$(686, 334, 76, 251)$,
$(581, 623, 224, 161)$,\adfsplit
$(621, 328, 345, 332)$,
$(637, 584, 110, 522)$,
$(302, 54, 473, 555)$,
$(265, 489, 126, 28)$,\adfsplit
$(649, 599, 285, 230)$,
$(170, 365, 462, 164)$,
$(514, 161, 436, 491)$,
$(358, 605, 525, 419)$,\adfsplit
$(490, 583, 339, 558)$,
$(524, 250, 46, 294)$,
$(490, 14, 147, 229)$,
$(398, 653, 594, 30)$,\adfsplit
$(630, 466, 272, 245)$,
$(383, 597, 510, 559)$,
$(470, 237, 29, 91)$,
$(38, 566, 611, 31)$,\adfsplit
$(501, 571, 268, 107)$,
$(143, 620, 529, 389)$,
$(49, 484, 382, 596)$,
$(462, 331, 230, 184)$,\adfsplit
$(438, 567, 281, 560)$,
$(376, 31, 189, 158)$,
$(529, 107, 553, 486)$,
$(380, 569, 400, 168)$,\adfsplit
$(387, 461, 287, 284)$,
$(584, 42, 466, 419)$,
$(303, 563, 9, 548)$,
$(558, 525, 206, 403)$,\adfsplit
$(0, 8, 83, 92)$,
$(0, 5, 13, 461)$,
$(0, 12, 316, 337)$,
$(0, 29, 267, 302)$,\adfsplit
$(0, 22, 95, 181)$,
$(0, 32, 146, 185)$,
$(0, 38, 320, 393)$,
$(0, 119, 242, 370)$,\adfsplit
$(0, 96, 213, 311)$,
$(0, 99, 271, 434)$,
$(0, 48, 203, 384)$,
$(0, 79, 275, 407)$,\adfsplit
$(0, 44, 110, 394)$,
$(0, 52, 277, 418)$,
$(0, 37, 94, 481)$,
$(0, 81, 205, 340)$,\adfsplit
$(0, 31, 312, 380)$,
$(0, 69, 190, 499)$,
$(0, 60, 137, 571)$,
$(0, 40, 128, 546)$

\adfLgap \noindent by the mapping:
$x \mapsto x \oplus (3 j)$ for $x < 684$,
$x \mapsto (x - 684 + 3 j \adfmod{27}) + 684$ for $684 \le x < 711$,
$x \mapsto x$ for $x \ge 711$,
$0 \le j < 228$
 for the first two blocks;
$x \mapsto x \oplus j$ for $x < 684$,
$x \mapsto (x - 684 + 3 j \adfmod{27}) + 684$ for $684 \le x < 711$,
$x \mapsto x$ for $x \ge 711$,
$0 \le j < 684$
 for the last 58 blocks.
\ADFvfyParStart{(713, ((2, 228, ((684, 3, (228, 3)), (27, 3), (2, 2))), (58, 684, ((684, 1, (228, 3)), (27, 3), (2, 2)))), ((38, 18), (29, 1)))} 

\adfDgap
\noindent{\boldmath $ 38^{18} 32^{1} $}~
With the point set $Z_{716}$ partitioned into
 residue classes modulo $18$ for $\{0, 1, \dots, 683\}$, and
 $\{684, 685, \dots, 715\}$,
 the design is generated from

\adfLgap 
$(714, 0, 1, 2)$,
$(715, 0, 229, 458)$,
$(684, 214, 432, 33)$,
$(684, 642, 379, 628)$,\adfsplit
$(684, 349, 563, 183)$,
$(684, 677, 50, 308)$,
$(685, 303, 348, 566)$,
$(685, 163, 357, 203)$,\adfsplit
$(685, 298, 112, 553)$,
$(685, 582, 368, 533)$,
$(686, 309, 479, 448)$,
$(686, 389, 266, 483)$,\adfsplit
$(686, 444, 439, 342)$,
$(686, 577, 32, 538)$,
$(687, 254, 404, 635)$,
$(687, 607, 529, 178)$,\adfsplit
$(687, 507, 81, 408)$,
$(687, 138, 220, 485)$,
$(688, 304, 270, 243)$,
$(688, 655, 81, 298)$,\adfsplit
$(688, 347, 360, 1)$,
$(688, 368, 245, 554)$,
$(88, 349, 75, 379)$,
$(129, 265, 74, 161)$,\adfsplit
$(52, 385, 182, 159)$,
$(563, 129, 280, 29)$,
$(410, 125, 569, 228)$,
$(610, 10, 603, 651)$,\adfsplit
$(139, 220, 396, 20)$,
$(679, 34, 434, 222)$,
$(6, 681, 146, 629)$,
$(139, 604, 669, 68)$,\adfsplit
$(194, 647, 9, 519)$,
$(211, 415, 663, 182)$,
$(167, 621, 273, 97)$,
$(287, 212, 262, 78)$,\adfsplit
$(435, 344, 567, 10)$,
$(358, 555, 117, 242)$,
$(635, 400, 467, 48)$,
$(51, 152, 86, 639)$,\adfsplit
$(287, 663, 223, 42)$,
$(179, 362, 138, 477)$,
$(238, 15, 86, 527)$,
$(339, 62, 219, 658)$,\adfsplit
$(362, 50, 459, 669)$,
$(170, 163, 122, 537)$,
$(618, 423, 679, 133)$,
$(129, 484, 386, 153)$,\adfsplit
$(516, 599, 55, 658)$,
$(523, 368, 378, 424)$,
$(402, 167, 554, 587)$,
$(358, 572, 18, 408)$,\adfsplit
$(461, 147, 449, 312)$,
$(565, 555, 353, 147)$,
$(208, 332, 155, 234)$,
$(556, 331, 282, 577)$,\adfsplit
$(31, 681, 462, 491)$,
$(10, 505, 88, 557)$,
$(677, 215, 662, 411)$,
$(312, 437, 349, 177)$,\adfsplit
$(472, 200, 85, 225)$,
$(154, 506, 518, 148)$,
$(645, 211, 155, 514)$,
$(294, 371, 401, 10)$,\adfsplit
$(357, 148, 291, 564)$,
$(176, 290, 363, 658)$,
$(421, 513, 634, 312)$,
$(66, 448, 33, 113)$,\adfsplit
$(668, 172, 246, 671)$,
$(162, 370, 222, 483)$,
$(206, 93, 668, 475)$,
$(218, 347, 47, 55)$,\adfsplit
$(274, 560, 155, 582)$,
$(415, 499, 314, 209)$,
$(360, 595, 81, 604)$,
$(63, 383, 600, 205)$,\adfsplit
$(523, 69, 120, 367)$,
$(482, 509, 200, 327)$,
$(568, 255, 466, 204)$,
$(592, 506, 627, 388)$,\adfsplit
$(411, 556, 369, 335)$,
$(675, 120, 483, 290)$,
$(66, 370, 637, 641)$,
$(673, 161, 221, 44)$,\adfsplit
$(119, 650, 667, 588)$,
$(137, 307, 102, 654)$,
$(0, 5, 154, 526)$,
$(0, 10, 21, 637)$,\adfsplit
$(0, 11, 80, 93)$,
$(0, 6, 497, 598)$,
$(0, 17, 42, 682)$,
$(0, 14, 345, 459)$,\adfsplit
$(0, 22, 537, 611)$,
$(0, 43, 281, 681)$,
$(0, 45, 65, 395)$,
$(0, 23, 256, 531)$,\adfsplit
$(3, 29, 76, 159)$,
$(0, 141, 435, 510)$,
$(0, 69, 207, 525)$,
$(0, 105, 263, 318)$,\adfsplit
$(0, 73, 359, 646)$,
$(0, 161, 166, 568)$,
$(0, 85, 358, 447)$,
$(0, 86, 444, 683)$,\adfsplit
$(0, 183, 298, 391)$,
$(0, 135, 157, 511)$,
$(0, 63, 406, 500)$,
$(0, 24, 343, 567)$,\adfsplit
$(0, 171, 241, 416)$,
$(0, 131, 146, 516)$,
$(0, 96, 397, 503)$,
$(0, 26, 221, 410)$,\adfsplit
$(0, 92, 291, 428)$,
$(0, 111, 151, 278)$,
$(0, 56, 120, 433)$,
$(0, 31, 200, 476)$,\adfsplit
$(0, 19, 110, 302)$,
$(0, 38, 103, 367)$,
$(0, 121, 331, 456)$

\adfLgap \noindent by the mapping:
$x \mapsto x \oplus (3 j)$ for $x < 684$,
$x \mapsto (x - 684 + 5 j \adfmod{30}) + 684$ for $684 \le x < 714$,
$x \mapsto x$ for $x \ge 714$,
$0 \le j < 228$
 for the first two blocks;
$x \mapsto x \oplus j \oplus j$ for $x < 684$,
$x \mapsto (x - 684 + 5 j \adfmod{30}) + 684$ for $684 \le x < 714$,
$x \mapsto x$ for $x \ge 714$,
$0 \le j < 342$
 for the last 117 blocks.
\ADFvfyParStart{(716, ((2, 228, ((684, 3, (228, 3)), (30, 5), (2, 2))), (117, 342, ((684, 2, (228, 3)), (30, 5), (2, 2)))), ((38, 18), (32, 1)))} 

\adfDgap
\noindent{\boldmath $ 38^{18} 35^{1} $}~
With the point set $Z_{719}$ partitioned into
 residue classes modulo $18$ for $\{0, 1, \dots, 683\}$, and
 $\{684, 685, \dots, 718\}$,
 the design is generated from

\adfLgap 
$(717, 0, 1, 2)$,
$(718, 0, 229, 458)$,
$(684, 223, 573, 617)$,
$(685, 223, 623, 183)$,\adfsplit
$(686, 222, 115, 113)$,
$(687, 250, 314, 240)$,
$(688, 662, 273, 52)$,
$(689, 568, 360, 473)$,\adfsplit
$(690, 184, 476, 486)$,
$(691, 437, 39, 220)$,
$(692, 518, 153, 676)$,
$(693, 12, 430, 467)$,\adfsplit
$(694, 589, 638, 252)$,
$(305, 561, 274, 362)$,
$(616, 86, 311, 277)$,
$(302, 269, 509, 672)$,\adfsplit
$(667, 533, 471, 355)$,
$(259, 379, 184, 383)$,
$(515, 534, 292, 422)$,
$(385, 554, 70, 26)$,\adfsplit
$(529, 196, 328, 631)$,
$(71, 55, 562, 210)$,
$(195, 191, 633, 448)$,
$(546, 169, 567, 667)$,\adfsplit
$(17, 128, 245, 369)$,
$(549, 376, 20, 619)$,
$(140, 59, 560, 413)$,
$(29, 462, 77, 80)$,\adfsplit
$(411, 334, 672, 544)$,
$(475, 239, 566, 219)$,
$(632, 481, 671, 101)$,
$(149, 668, 207, 633)$,\adfsplit
$(494, 418, 13, 110)$,
$(105, 268, 210, 399)$,
$(366, 121, 444, 340)$,
$(625, 307, 531, 598)$,\adfsplit
$(105, 482, 646, 60)$,
$(439, 495, 568, 481)$,
$(596, 27, 383, 664)$,
$(236, 5, 407, 529)$,\adfsplit
$(618, 515, 333, 86)$,
$(0, 6, 149, 174)$,
$(0, 8, 68, 238)$,
$(0, 9, 272, 280)$,\adfsplit
$(0, 12, 26, 55)$,
$(0, 23, 148, 363)$,
$(0, 15, 185, 277)$,
$(0, 89, 222, 409)$,\adfsplit
$(0, 30, 209, 476)$,
$(0, 52, 197, 266)$,
$(0, 63, 235, 358)$,
$(0, 25, 184, 230)$,\adfsplit
$(0, 53, 257, 370)$,
$(0, 82, 191, 274)$,
$(0, 92, 309, 408)$,
$(0, 66, 150, 404)$,\adfsplit
$(0, 56, 175, 472)$,
$(0, 24, 65, 608)$,
$(0, 96, 236, 371)$,
$(0, 38, 139, 348)$,\adfsplit
$(0, 19, 157, 237)$

\adfLgap \noindent by the mapping:
$x \mapsto x \oplus (3 j)$ for $x < 684$,
$x \mapsto (x - 684 + 11 j \adfmod{33}) + 684$ for $684 \le x < 717$,
$x \mapsto x$ for $x \ge 717$,
$0 \le j < 228$
 for the first two blocks;
$x \mapsto x \oplus j$ for $x < 684$,
$x \mapsto (x - 684 + 11 j \adfmod{33}) + 684$ for $684 \le x < 717$,
$x \mapsto x$ for $x \ge 717$,
$0 \le j < 684$
 for the last 59 blocks.
\ADFvfyParStart{(719, ((2, 228, ((684, 3, (228, 3)), (33, 11), (2, 2))), (59, 684, ((684, 1, (228, 3)), (33, 11), (2, 2)))), ((38, 18), (35, 1)))} 

\section{4-GDDs for the proof of Lemma \ref{lem:4-GDD 40^u m^1}}
\label{app:4-GDD 40^u m^1}
\adfhide{
$ 40^6 73^1 $,
$ 40^6 79^1 $,
$ 40^6 82^1 $,
$ 40^6 91^1 $,
$ 40^6 94^1 $,
$ 40^6 97^1 $,
$ 40^9 139^1 $,
$ 40^9 142^1 $,
$ 40^9 151^1 $,
$ 40^9 154^1 $ and
$ 40^9 157^1 $.
}

\adfDgap
\noindent{\boldmath $ 40^{6} 73^{1} $}~
With the point set $Z_{313}$ partitioned into
 residue classes modulo $6$ for $\{0, 1, \dots, 239\}$, and
 $\{240, 241, \dots, 312\}$,
 the design is generated from

\adfLgap 
$(300, 202, 126, 183)$,
$(300, 204, 172, 50)$,
$(300, 139, 141, 215)$,
$(300, 140, 75, 203)$,\adfsplit
$(300, 38, 73, 129)$,
$(300, 18, 137, 16)$,
$(300, 175, 77, 224)$,
$(300, 48, 229, 142)$,\adfsplit
$(240, 233, 48, 136)$,
$(240, 183, 42, 98)$,
$(240, 126, 52, 207)$,
$(240, 137, 225, 58)$,\adfsplit
$(240, 4, 227, 216)$,
$(240, 150, 155, 223)$,
$(240, 101, 145, 132)$,
$(240, 76, 26, 109)$,\adfsplit
$(240, 158, 201, 67)$,
$(240, 143, 19, 116)$,
$(240, 161, 151, 24)$,
$(240, 65, 122, 97)$,\adfsplit
$(240, 43, 110, 83)$,
$(240, 198, 70, 29)$,
$(240, 237, 186, 179)$,
$(240, 112, 157, 180)$,\adfsplit
$(240, 202, 86, 141)$,
$(240, 40, 135, 20)$,
$(240, 234, 34, 1)$,
$(240, 212, 220, 49)$,\adfsplit
$(240, 199, 214, 215)$,
$(240, 71, 93, 188)$,
$(240, 167, 171, 130)$,
$(240, 111, 173, 181)$,\adfsplit
$(240, 69, 192, 115)$,
$(240, 28, 222, 44)$,
$(240, 238, 123, 104)$,
$(240, 208, 73, 99)$,\adfsplit
$(240, 194, 197, 64)$,
$(240, 125, 153, 18)$,
$(0, 21, 92, 282)$,
$(0, 23, 177, 269)$,\adfsplit
$(0, 4, 201, 275)$,
$(0, 14, 58, 284)$,
$(0, 64, 229, 272)$,
$(0, 189, 227, 289)$,\adfsplit
$(0, 49, 195, 258)$,
$(0, 26, 129, 260)$,
$(0, 61, 113, 293)$,
$(0, 39, 187, 246)$,\adfsplit
$(0, 29, 38, 147)$,
$(0, 205, 219, 239)$,
$(0, 9, 139, 218)$,
$(0, 10, 25, 75)$,\adfsplit
$(0, 17, 98, 157)$,
$(0, 53, 100, 211)$,
$(0, 82, 171, 235)$,
$(0, 47, 151, 188)$,\adfsplit
$(0, 34, 101, 104)$,
$(312, 0, 80, 160)$,
$(312, 1, 81, 161)$

\adfLgap \noindent by the mapping:
$x \mapsto x + 2 j \adfmod{240}$ for $x < 240$,
$x \mapsto (x +  j \adfmod{60}) + 240$ for $240 \le x < 300$,
$x \mapsto (x +  j \adfmod{12}) + 300$ for $300 \le x < 312$,
$312 \mapsto 312$,
$0 \le j < 120$
 for the first 57 blocks,
$0 \le j < 40$
 for the last two blocks.
\ADFvfyParStart{(313, ((57, 120, ((240, 2), (60, 1), (12, 1), (1, 1))), (2, 40, ((240, 2), (60, 1), (12, 1), (1, 1)))), ((40, 6), (73, 1)))} 

\adfDgap
\noindent{\boldmath $ 40^{6} 79^{1} $}~
With the point set $Z_{319}$ partitioned into
 residue classes modulo $6$ for $\{0, 1, \dots, 239\}$, and
 $\{240, 241, \dots, 318\}$,
 the design is generated from

\adfLgap 
$(315, 125, 236, 195)$,
$(315, 103, 82, 30)$,
$(240, 129, 67, 29)$,
$(240, 91, 206, 78)$,\adfsplit
$(240, 163, 32, 130)$,
$(240, 46, 194, 11)$,
$(240, 34, 135, 175)$,
$(240, 177, 169, 114)$,\adfsplit
$(240, 227, 112, 110)$,
$(240, 113, 98, 192)$,
$(240, 63, 185, 210)$,
$(240, 156, 51, 148)$,\adfsplit
$(241, 185, 27, 206)$,
$(241, 21, 98, 7)$,
$(241, 221, 204, 64)$,
$(241, 138, 142, 89)$,\adfsplit
$(241, 103, 195, 50)$,
$(241, 199, 28, 92)$,
$(241, 216, 61, 159)$,
$(241, 23, 175, 136)$,\adfsplit
$(241, 70, 42, 3)$,
$(241, 74, 60, 227)$,
$(242, 206, 25, 168)$,
$(242, 70, 43, 59)$,\adfsplit
$(242, 30, 113, 14)$,
$(242, 225, 49, 215)$,
$(242, 158, 124, 6)$,
$(242, 159, 12, 178)$,\adfsplit
$(242, 71, 84, 91)$,
$(242, 97, 207, 182)$,
$(242, 107, 63, 16)$,
$(242, 201, 80, 22)$,\adfsplit
$(243, 178, 57, 181)$,
$(243, 8, 201, 197)$,
$(243, 138, 76, 32)$,
$(243, 193, 161, 159)$,\adfsplit
$(243, 85, 130, 20)$,
$(243, 116, 66, 135)$,
$(243, 74, 213, 67)$,
$(243, 54, 143, 229)$,\adfsplit
$(0, 10, 77, 288)$,
$(0, 1, 56, 308)$,
$(0, 9, 37, 154)$,
$(0, 20, 46, 197)$,\adfsplit
$(0, 22, 49, 231)$,
$(0, 11, 225, 244)$,
$(0, 23, 79, 269)$,
$(0, 109, 159, 249)$,\adfsplit
$(0, 41, 70, 235)$,
$(0, 43, 207, 304)$,
$(0, 103, 239, 264)$,
$(0, 75, 181, 203)$,\adfsplit
$(0, 31, 200, 259)$,
$(0, 32, 137, 294)$,
$(0, 35, 164, 274)$,
$(0, 29, 136, 217)$,\adfsplit
$(0, 5, 158, 254)$,
$(0, 45, 113, 116)$,
$(0, 51, 68, 299)$,
$(318, 0, 80, 160)$,\adfsplit
$(318, 1, 81, 161)$

\adfLgap \noindent by the mapping:
$x \mapsto x + 2 j \adfmod{240}$ for $x < 240$,
$x \mapsto (x - 240 + 5 j \adfmod{75}) + 240$ for $240 \le x < 315$,
$x \mapsto (x +  j \adfmod{3}) + 315$ for $315 \le x < 318$,
$318 \mapsto 318$,
$0 \le j < 120$
 for the first 59 blocks,
$0 \le j < 40$
 for the last two blocks.
\ADFvfyParStart{(319, ((59, 120, ((240, 2), (75, 5), (3, 1), (1, 1))), (2, 40, ((240, 2), (75, 5), (3, 1), (1, 1)))), ((40, 6), (79, 1)))} 

\adfDgap
\noindent{\boldmath $ 40^{6} 82^{1} $}~
With the point set $Z_{322}$ partitioned into
 residue classes modulo $6$ for $\{0, 1, \dots, 239\}$, and
 $\{240, 241, \dots, 321\}$,
 the design is generated from

\adfLgap 
$(300, 107, 79, 117)$,
$(301, 161, 136, 96)$,
$(302, 85, 29, 234)$,
$(303, 114, 209, 97)$,\adfsplit
$(304, 23, 210, 133)$,
$(305, 206, 118, 18)$,
$(306, 76, 213, 143)$,
$(240, 65, 201, 192)$,\adfsplit
$(240, 60, 83, 82)$,
$(240, 78, 129, 71)$,
$(240, 197, 73, 100)$,
$(240, 89, 26, 207)$,\adfsplit
$(240, 33, 139, 188)$,
$(240, 104, 159, 223)$,
$(240, 134, 55, 63)$,
$(240, 109, 238, 113)$,\adfsplit
$(240, 96, 7, 76)$,
$(240, 15, 234, 172)$,
$(240, 190, 37, 236)$,
$(0, 2, 39, 197)$,\adfsplit
$(0, 5, 98, 258)$,
$(0, 13, 81, 154)$,
$(0, 16, 50, 266)$,
$(0, 14, 33, 179)$,\adfsplit
$(0, 11, 225, 254)$,
$(0, 57, 166, 298)$,
$(0, 29, 164, 299)$,
$(0, 3, 196, 275)$,\adfsplit
$(0, 31, 123, 265)$,
$(0, 32, 139, 268)$,
$(321, 0, 80, 160)$

\adfLgap \noindent by the mapping:
$x \mapsto x +  j \adfmod{240}$ for $x < 240$,
$x \mapsto (x +  j \adfmod{60}) + 240$ for $240 \le x < 300$,
$x \mapsto (x - 300 + 7 j \adfmod{21}) + 300$ for $300 \le x < 321$,
$321 \mapsto 321$,
$0 \le j < 240$
 for the first 30 blocks,
$0 \le j < 80$
 for the last block.
\ADFvfyParStart{(322, ((30, 240, ((240, 1), (60, 1), (21, 7), (1, 1))), (1, 80, ((240, 1), (60, 1), (21, 7), (1, 1)))), ((40, 6), (82, 1)))} 

\adfDgap
\noindent{\boldmath $ 40^{6} 91^{1} $}~
With the point set $Z_{331}$ partitioned into
 residue classes modulo $6$ for $\{0, 1, \dots, 239\}$, and
 $\{240, 241, \dots, 330\}$,
 the design is generated from

\adfLgap 
$(240, 207, 223, 102)$,
$(240, 45, 80, 173)$,
$(240, 6, 209, 208)$,
$(240, 151, 107, 81)$,\adfsplit
$(240, 79, 182, 0)$,
$(240, 26, 166, 123)$,
$(240, 204, 221, 111)$,
$(240, 234, 93, 215)$,\adfsplit
$(240, 143, 109, 216)$,
$(240, 230, 15, 210)$,
$(240, 177, 25, 232)$,
$(240, 157, 12, 118)$,\adfsplit
$(240, 239, 98, 130)$,
$(240, 56, 67, 77)$,
$(240, 13, 191, 138)$,
$(240, 224, 159, 82)$,\adfsplit
$(240, 181, 214, 212)$,
$(240, 184, 68, 125)$,
$(240, 168, 14, 40)$,
$(240, 235, 9, 76)$,\adfsplit
$(241, 56, 42, 47)$,
$(241, 83, 0, 196)$,
$(241, 201, 161, 94)$,
$(241, 172, 207, 131)$,\adfsplit
$(241, 62, 179, 225)$,
$(241, 183, 97, 124)$,
$(241, 166, 115, 212)$,
$(241, 194, 43, 125)$,\adfsplit
$(241, 178, 26, 174)$,
$(241, 108, 79, 219)$,
$(241, 90, 73, 137)$,
$(241, 144, 231, 80)$,\adfsplit
$(241, 96, 88, 151)$,
$(241, 8, 18, 70)$,
$(241, 195, 145, 173)$,
$(241, 170, 169, 189)$,\adfsplit
$(241, 237, 95, 98)$,
$(241, 149, 164, 142)$,
$(241, 33, 192, 1)$,
$(241, 127, 100, 66)$,\adfsplit
$(242, 169, 128, 45)$,
$(242, 173, 194, 100)$,
$(242, 210, 71, 123)$,
$(242, 218, 91, 147)$,\adfsplit
$(242, 219, 144, 215)$,
$(242, 138, 208, 203)$,
$(0, 9, 155, 212)$,
$(0, 13, 122, 242)$,\adfsplit
$(0, 16, 119, 177)$,
$(0, 51, 217, 269)$,
$(0, 15, 149, 299)$,
$(0, 40, 135, 227)$,\adfsplit
$(0, 23, 31, 263)$,
$(0, 193, 195, 254)$,
$(0, 29, 201, 293)$,
$(0, 68, 191, 257)$,\adfsplit
$(0, 91, 129, 314)$,
$(0, 50, 215, 251)$,
$(0, 3, 158, 248)$,
$(0, 43, 104, 284)$,\adfsplit
$(0, 125, 229, 278)$,
$(0, 69, 76, 260)$,
$(0, 56, 130, 317)$,
$(330, 0, 80, 160)$,\adfsplit
$(330, 1, 81, 161)$

\adfLgap \noindent by the mapping:
$x \mapsto x + 2 j \adfmod{240}$ for $x < 240$,
$x \mapsto (x - 240 + 3 j \adfmod{90}) + 240$ for $240 \le x < 330$,
$330 \mapsto 330$,
$0 \le j < 120$
 for the first 63 blocks,
$0 \le j < 40$
 for the last two blocks.
\ADFvfyParStart{(331, ((63, 120, ((240, 2), (90, 3), (1, 1))), (2, 40, ((240, 2), (90, 3), (1, 1)))), ((40, 6), (91, 1)))} 

\adfDgap
\noindent{\boldmath $ 40^{6} 94^{1} $}~
With the point set $Z_{334}$ partitioned into
 residue classes modulo $6$ for $\{0, 1, \dots, 239\}$, and
 $\{240, 241, \dots, 333\}$,
 the design is generated from

\adfLgap 
$(330, 119, 73, 57)$,
$(240, 68, 149, 217)$,
$(240, 17, 234, 207)$,
$(240, 64, 143, 122)$,\adfsplit
$(240, 18, 118, 133)$,
$(240, 101, 236, 129)$,
$(240, 200, 192, 70)$,
$(240, 76, 175, 21)$,\adfsplit
$(240, 31, 95, 104)$,
$(240, 195, 49, 126)$,
$(240, 112, 0, 3)$,
$(241, 1, 135, 167)$,\adfsplit
$(241, 104, 193, 149)$,
$(241, 88, 83, 12)$,
$(241, 3, 56, 226)$,
$(241, 41, 184, 128)$,\adfsplit
$(241, 147, 140, 160)$,
$(241, 212, 127, 186)$,
$(241, 99, 48, 95)$,
$(241, 120, 169, 81)$,\adfsplit
$(241, 112, 235, 204)$,
$(0, 10, 129, 272)$,
$(0, 11, 226, 296)$,
$(0, 19, 207, 251)$,\adfsplit
$(0, 34, 147, 266)$,
$(0, 22, 197, 281)$,
$(0, 2, 37, 329)$,
$(0, 40, 103, 317)$,\adfsplit
$(0, 61, 165, 302)$,
$(0, 41, 139, 326)$,
$(0, 29, 67, 124)$,
$(0, 1, 158, 314)$,\adfsplit
$(333, 0, 80, 160)$

\adfLgap \noindent by the mapping:
$x \mapsto x +  j \adfmod{240}$ for $x < 240$,
$x \mapsto (x - 240 + 3 j \adfmod{90}) + 240$ for $240 \le x < 330$,
$x \mapsto (x +  j \adfmod{3}) + 330$ for $330 \le x < 333$,
$333 \mapsto 333$,
$0 \le j < 240$
 for the first 32 blocks,
$0 \le j < 80$
 for the last block.
\ADFvfyParStart{(334, ((32, 240, ((240, 1), (90, 3), (3, 1), (1, 1))), (1, 80, ((240, 1), (90, 3), (3, 1), (1, 1)))), ((40, 6), (94, 1)))} 

\adfDgap
\noindent{\boldmath $ 40^{6} 97^{1} $}~
With the point set $Z_{337}$ partitioned into
 residue classes modulo $6$ for $\{0, 1, \dots, 239\}$, and
 $\{240, 241, \dots, 336\}$,
 the design is generated from

\adfLgap 
$(330, 185, 133, 110)$,
$(330, 71, 45, 154)$,
$(330, 123, 31, 12)$,
$(330, 224, 114, 220)$,\adfsplit
$(240, 222, 103, 136)$,
$(240, 69, 211, 140)$,
$(240, 228, 112, 55)$,
$(240, 95, 144, 88)$,\adfsplit
$(240, 153, 36, 13)$,
$(240, 86, 154, 1)$,
$(240, 50, 82, 83)$,
$(240, 11, 220, 127)$,\adfsplit
$(240, 219, 182, 184)$,
$(240, 75, 72, 119)$,
$(240, 74, 201, 139)$,
$(240, 114, 183, 190)$,\adfsplit
$(240, 210, 226, 225)$,
$(240, 161, 18, 238)$,
$(240, 56, 53, 66)$,
$(240, 97, 8, 231)$,\adfsplit
$(240, 65, 180, 85)$,
$(240, 89, 152, 87)$,
$(240, 169, 17, 164)$,
$(240, 158, 167, 57)$,\adfsplit
$(241, 50, 22, 11)$,
$(241, 70, 168, 147)$,
$(241, 76, 122, 139)$,
$(241, 159, 146, 24)$,\adfsplit
$(241, 158, 216, 25)$,
$(241, 201, 233, 163)$,
$(241, 180, 101, 28)$,
$(241, 229, 143, 171)$,\adfsplit
$(241, 97, 63, 184)$,
$(241, 178, 137, 133)$,
$(241, 154, 107, 8)$,
$(241, 129, 198, 55)$,\adfsplit
$(241, 116, 102, 35)$,
$(241, 150, 195, 211)$,
$(241, 172, 117, 61)$,
$(241, 164, 234, 89)$,\adfsplit
$(241, 67, 166, 126)$,
$(241, 119, 105, 134)$,
$(241, 213, 80, 132)$,
$(241, 212, 100, 185)$,\adfsplit
$(242, 126, 76, 181)$,
$(242, 159, 10, 212)$,
$(242, 34, 185, 174)$,
$(242, 184, 141, 120)$,\adfsplit
$(0, 41, 92, 205)$,
$(0, 22, 104, 245)$,
$(0, 29, 166, 290)$,
$(0, 8, 34, 260)$,\adfsplit
$(0, 53, 62, 278)$,
$(0, 27, 221, 293)$,
$(0, 109, 203, 254)$,
$(0, 51, 209, 320)$,\adfsplit
$(0, 115, 137, 305)$,
$(0, 25, 215, 266)$,
$(0, 167, 235, 287)$,
$(0, 59, 123, 296)$,\adfsplit
$(0, 39, 79, 314)$,
$(0, 44, 101, 329)$,
$(1, 11, 123, 266)$,
$(0, 83, 91, 272)$,\adfsplit
$(0, 43, 179, 257)$,
$(336, 0, 80, 160)$,
$(336, 1, 81, 161)$

\adfLgap \noindent by the mapping:
$x \mapsto x + 2 j \adfmod{240}$ for $x < 240$,
$x \mapsto (x - 240 + 3 j \adfmod{90}) + 240$ for $240 \le x < 330$,
$x \mapsto (x +  j \adfmod{6}) + 330$ for $330 \le x < 336$,
$336 \mapsto 336$,
$0 \le j < 120$
 for the first 65 blocks,
$0 \le j < 40$
 for the last two blocks.
\ADFvfyParStart{(337, ((65, 120, ((240, 2), (90, 3), (6, 1), (1, 1))), (2, 40, ((240, 2), (90, 3), (6, 1), (1, 1)))), ((40, 6), (97, 1)))} 

\adfDgap
\noindent{\boldmath $ 40^{9} 139^{1} $}~
With the point set $Z_{499}$ partitioned into
 residue classes modulo $9$ for $\{0, 1, \dots, 359\}$, and
 $\{360, 361, \dots, 498\}$,
 the design is generated from

\adfLgap 
$(480, 87, 352, 143)$,
$(480, 65, 22, 345)$,
$(480, 213, 77, 253)$,
$(480, 120, 162, 218)$,\adfsplit
$(480, 296, 16, 86)$,
$(480, 128, 45, 241)$,
$(480, 72, 103, 187)$,
$(480, 10, 23, 142)$,\adfsplit
$(480, 134, 276, 199)$,
$(480, 112, 135, 174)$,
$(480, 248, 89, 330)$,
$(480, 47, 111, 301)$,\adfsplit
$(360, 333, 138, 177)$,
$(360, 196, 249, 117)$,
$(360, 254, 163, 342)$,
$(360, 107, 223, 126)$,\adfsplit
$(360, 49, 81, 61)$,
$(360, 274, 69, 269)$,
$(360, 251, 106, 252)$,
$(360, 133, 296, 200)$,\adfsplit
$(360, 185, 87, 52)$,
$(360, 150, 101, 135)$,
$(360, 268, 195, 17)$,
$(360, 146, 230, 99)$,\adfsplit
$(360, 204, 79, 74)$,
$(360, 308, 16, 109)$,
$(360, 108, 323, 244)$,
$(360, 232, 118, 225)$,\adfsplit
$(360, 23, 220, 127)$,
$(360, 326, 310, 187)$,
$(360, 111, 94, 140)$,
$(360, 119, 77, 313)$,\adfsplit
$(360, 277, 218, 285)$,
$(360, 267, 248, 47)$,
$(360, 39, 328, 353)$,
$(360, 0, 44, 275)$,\adfsplit
$(360, 33, 151, 210)$,
$(360, 2, 78, 337)$,
$(360, 59, 291, 82)$,
$(360, 24, 236, 62)$,\adfsplit
$(360, 265, 60, 250)$,
$(360, 19, 63, 89)$,
$(360, 282, 332, 38)$,
$(360, 85, 142, 71)$,\adfsplit
$(360, 354, 1, 186)$,
$(360, 32, 54, 123)$,
$(360, 261, 96, 156)$,
$(360, 355, 224, 55)$,\adfsplit
$(360, 53, 160, 168)$,
$(360, 215, 5, 58)$,
$(360, 41, 91, 312)$,
$(360, 166, 290, 64)$,\adfsplit
$(361, 136, 42, 359)$,
$(361, 96, 236, 284)$,
$(361, 169, 82, 185)$,
$(361, 95, 160, 170)$,\adfsplit
$(361, 37, 268, 346)$,
$(361, 251, 264, 240)$,
$(361, 343, 347, 319)$,
$(361, 57, 342, 149)$,\adfsplit
$(361, 100, 355, 97)$,
$(361, 123, 86, 22)$,
$(361, 104, 190, 93)$,
$(361, 234, 255, 47)$,\adfsplit
$(361, 295, 316, 317)$,
$(361, 33, 121, 88)$,
$(361, 323, 45, 275)$,
$(361, 151, 13, 10)$,\adfsplit
$(361, 207, 318, 64)$,
$(361, 43, 351, 138)$,
$(0, 2, 127, 273)$,
$(0, 26, 55, 149)$,\adfsplit
$(0, 6, 20, 199)$,
$(0, 40, 251, 257)$,
$(0, 4, 61, 122)$,
$(0, 32, 73, 169)$,\adfsplit
$(0, 34, 194, 397)$,
$(0, 12, 164, 473)$,
$(0, 110, 248, 409)$,
$(0, 52, 156, 383)$,\adfsplit
$(0, 74, 202, 389)$,
$(0, 28, 58, 415)$,
$(0, 49, 176, 469)$,
$(0, 147, 178, 371)$,\adfsplit
$(0, 83, 100, 425)$,
$(0, 227, 327, 479)$,
$(0, 85, 206, 421)$,
$(0, 92, 203, 447)$,\adfsplit
$(0, 185, 319, 373)$,
$(0, 177, 291, 391)$,
$(0, 247, 325, 407)$,
$(0, 187, 309, 449)$,\adfsplit
$(0, 47, 121, 423)$,
$(0, 77, 219, 369)$,
$(1, 31, 203, 369)$,
$(1, 39, 149, 207)$,\adfsplit
$(1, 11, 285, 441)$,
$(1, 63, 175, 473)$,
$(1, 3, 69, 477)$,
$(498, 0, 120, 240)$,\adfsplit
$(498, 1, 121, 241)$

\adfLgap \noindent by the mapping:
$x \mapsto x + 2 j \adfmod{360}$ for $x < 360$,
$x \mapsto (x + 2 j \adfmod{120}) + 360$ for $360 \le x < 480$,
$x \mapsto (x - 480 +  j \adfmod{18}) + 480$ for $480 \le x < 498$,
$498 \mapsto 498$,
$0 \le j < 180$
 for the first 99 blocks,
$0 \le j < 60$
 for the last two blocks.
\ADFvfyParStart{(499, ((99, 180, ((360, 2), (120, 2), (18, 1), (1, 1))), (2, 60, ((360, 2), (120, 2), (18, 1), (1, 1)))), ((40, 9), (139, 1)))} 

\adfDgap
\noindent{\boldmath $ 40^{9} 142^{1} $}~
With the point set $Z_{502}$ partitioned into
 residue classes modulo $9$ for $\{0, 1, \dots, 359\}$, and
 $\{360, 361, \dots, 501\}$,
 the design is generated from

\adfLgap 
$(480, 285, 263, 184)$,
$(481, 348, 266, 325)$,
$(482, 316, 119, 303)$,
$(483, 215, 208, 132)$,\adfsplit
$(484, 270, 53, 1)$,
$(485, 313, 260, 69)$,
$(486, 62, 103, 342)$,
$(360, 15, 289, 134)$,\adfsplit
$(360, 332, 294, 128)$,
$(360, 336, 227, 315)$,
$(360, 259, 78, 147)$,
$(360, 316, 220, 61)$,\adfsplit
$(360, 207, 114, 292)$,
$(360, 258, 342, 268)$,
$(360, 209, 166, 68)$,
$(360, 355, 39, 277)$,\adfsplit
$(360, 22, 24, 7)$,
$(360, 249, 218, 215)$,
$(360, 133, 3, 23)$,
$(360, 63, 62, 2)$,\adfsplit
$(360, 156, 288, 307)$,
$(360, 339, 41, 244)$,
$(360, 271, 42, 194)$,
$(360, 343, 210, 85)$,\adfsplit
$(360, 305, 79, 256)$,
$(360, 125, 131, 238)$,
$(360, 81, 175, 10)$,
$(360, 203, 153, 97)$,\adfsplit
$(360, 260, 228, 160)$,
$(360, 270, 146, 197)$,
$(0, 4, 29, 334)$,
$(0, 5, 145, 211)$,\adfsplit
$(0, 8, 344, 429)$,
$(0, 11, 39, 374)$,
$(0, 14, 47, 161)$,
$(0, 12, 200, 360)$,\adfsplit
$(0, 42, 256, 410)$,
$(0, 75, 167, 375)$,
$(0, 70, 218, 386)$,
$(0, 67, 237, 361)$,\adfsplit
$(0, 89, 263, 392)$,
$(0, 87, 224, 430)$,
$(0, 64, 202, 371)$,
$(0, 40, 232, 451)$,\adfsplit
$(0, 46, 295, 422)$,
$(0, 35, 231, 424)$,
$(0, 103, 221, 379)$,
$(0, 57, 302, 416)$,\adfsplit
$(0, 37, 210, 407)$,
$(0, 48, 233, 435)$,
$(501, 0, 120, 240)$

\adfLgap \noindent by the mapping:
$x \mapsto x +  j \adfmod{360}$ for $x < 360$,
$x \mapsto (x +  j \adfmod{120}) + 360$ for $360 \le x < 480$,
$x \mapsto (x - 480 + 7 j \adfmod{21}) + 480$ for $480 \le x < 501$,
$501 \mapsto 501$,
$0 \le j < 360$
 for the first 50 blocks,
$0 \le j < 120$
 for the last block.
\ADFvfyParStart{(502, ((50, 360, ((360, 1), (120, 1), (21, 7), (1, 1))), (1, 120, ((360, 1), (120, 1), (21, 7), (1, 1)))), ((40, 9), (142, 1)))} 

\adfDgap
\noindent{\boldmath $ 40^{9} 151^{1} $}~
With the point set $Z_{511}$ partitioned into
 residue classes modulo $9$ for $\{0, 1, \dots, 359\}$, and
 $\{360, 361, \dots, 510\}$,
 the design is generated from

\adfLgap 
$(480, 186, 299, 230)$,
$(480, 92, 72, 163)$,
$(480, 91, 93, 141)$,
$(480, 240, 19, 238)$,\adfsplit
$(480, 154, 99, 65)$,
$(480, 195, 256, 74)$,
$(480, 41, 90, 51)$,
$(480, 355, 278, 209)$,\adfsplit
$(480, 267, 297, 102)$,
$(480, 112, 224, 275)$,
$(480, 176, 197, 193)$,
$(480, 88, 202, 183)$,\adfsplit
$(480, 184, 160, 20)$,
$(480, 205, 36, 348)$,
$(480, 130, 173, 105)$,
$(480, 83, 68, 204)$,\adfsplit
$(480, 169, 294, 206)$,
$(480, 286, 11, 181)$,
$(480, 62, 78, 307)$,
$(480, 287, 129, 217)$,\adfsplit
$(360, 15, 148, 122)$,
$(360, 40, 37, 169)$,
$(360, 237, 324, 123)$,
$(360, 225, 251, 146)$,\adfsplit
$(360, 116, 272, 300)$,
$(360, 314, 85, 195)$,
$(360, 211, 349, 158)$,
$(360, 142, 41, 48)$,\adfsplit
$(360, 279, 238, 186)$,
$(360, 275, 214, 63)$,
$(360, 210, 221, 18)$,
$(360, 347, 295, 307)$,\adfsplit
$(360, 223, 29, 179)$,
$(360, 153, 218, 293)$,
$(360, 100, 197, 178)$,
$(360, 333, 73, 309)$,\adfsplit
$(360, 355, 92, 88)$,
$(360, 71, 150, 145)$,
$(360, 81, 245, 222)$,
$(360, 274, 257, 323)$,\adfsplit
$(360, 287, 23, 228)$,
$(360, 264, 267, 151)$,
$(360, 50, 19, 110)$,
$(360, 283, 252, 318)$,\adfsplit
$(360, 89, 156, 166)$,
$(360, 136, 326, 72)$,
$(360, 6, 140, 134)$,
$(360, 199, 297, 250)$,\adfsplit
$(360, 172, 45, 87)$,
$(360, 336, 240, 182)$,
$(360, 239, 184, 351)$,
$(360, 106, 248, 253)$,\adfsplit
$(360, 219, 301, 322)$,
$(360, 68, 305, 247)$,
$(360, 354, 171, 104)$,
$(360, 80, 21, 113)$,\adfsplit
$(360, 241, 316, 42)$,
$(360, 294, 232, 95)$,
$(360, 217, 249, 44)$,
$(360, 296, 124, 310)$,\adfsplit
$(361, 213, 331, 9)$,
$(361, 61, 279, 145)$,
$(361, 120, 229, 16)$,
$(361, 20, 123, 138)$,\adfsplit
$(361, 184, 192, 266)$,
$(361, 234, 11, 12)$,
$(361, 79, 44, 255)$,
$(361, 122, 99, 353)$,\adfsplit
$(361, 332, 96, 289)$,
$(361, 264, 231, 217)$,
$(361, 80, 10, 191)$,
$(361, 160, 335, 167)$,\adfsplit
$(361, 183, 189, 161)$,
$(361, 246, 288, 193)$,
$(361, 338, 297, 76)$,
$(361, 28, 78, 185)$,\adfsplit
$(0, 1, 65, 210)$,
$(0, 29, 40, 123)$,
$(0, 12, 25, 197)$,
$(0, 37, 57, 363)$,\adfsplit
$(0, 39, 143, 371)$,
$(0, 251, 331, 365)$,
$(0, 119, 127, 435)$,
$(0, 68, 183, 393)$,\adfsplit
$(0, 101, 287, 381)$,
$(0, 145, 201, 473)$,
$(0, 125, 208, 379)$,
$(1, 79, 207, 365)$,\adfsplit
$(0, 73, 303, 439)$,
$(1, 87, 209, 431)$,
$(0, 87, 149, 425)$,
$(1, 47, 225, 461)$,\adfsplit
$(0, 187, 289, 429)$,
$(1, 17, 301, 439)$,
$(0, 89, 249, 413)$,
$(0, 133, 146, 373)$,\adfsplit
$(0, 32, 116, 447)$,
$(0, 80, 228, 387)$,
$(0, 38, 130, 407)$,
$(0, 30, 314, 449)$,\adfsplit
$(0, 34, 200, 391)$,
$(0, 22, 260, 443)$,
$(0, 56, 258, 395)$,
$(510, 0, 120, 240)$,\adfsplit
$(510, 1, 121, 241)$

\adfLgap \noindent by the mapping:
$x \mapsto x + 2 j \adfmod{360}$ for $x < 360$,
$x \mapsto (x + 2 j \adfmod{120}) + 360$ for $360 \le x < 480$,
$x \mapsto (x +  j \adfmod{30}) + 480$ for $480 \le x < 510$,
$510 \mapsto 510$,
$0 \le j < 180$
 for the first 103 blocks,
$0 \le j < 60$
 for the last two blocks.
\ADFvfyParStart{(511, ((103, 180, ((360, 2), (120, 2), (30, 1), (1, 1))), (2, 60, ((360, 2), (120, 2), (30, 1), (1, 1)))), ((40, 9), (151, 1)))} 

\adfDgap
\noindent{\boldmath $ 40^{9} 154^{1} $}~
With the point set $Z_{514}$ partitioned into
 residue classes modulo $9$ for $\{0, 1, \dots, 359\}$, and
 $\{360, 361, \dots, 513\}$,
 the design is generated from

\adfLgap 
$(480, 228, 106, 5)$,
$(481, 76, 252, 137)$,
$(482, 37, 249, 335)$,
$(483, 234, 49, 104)$,\adfsplit
$(484, 281, 192, 28)$,
$(485, 253, 306, 149)$,
$(486, 205, 209, 348)$,
$(487, 303, 70, 350)$,\adfsplit
$(488, 78, 112, 98)$,
$(489, 1, 156, 155)$,
$(490, 197, 289, 0)$,
$(360, 87, 200, 104)$,\adfsplit
$(360, 298, 306, 275)$,
$(360, 238, 287, 357)$,
$(360, 208, 258, 70)$,
$(360, 261, 92, 194)$,\adfsplit
$(360, 304, 42, 100)$,
$(360, 48, 195, 305)$,
$(360, 342, 77, 103)$,
$(360, 34, 31, 46)$,\adfsplit
$(360, 251, 191, 230)$,
$(360, 354, 270, 329)$,
$(360, 299, 301, 133)$,
$(360, 68, 245, 359)$,\adfsplit
$(360, 290, 149, 215)$,
$(360, 248, 97, 171)$,
$(360, 324, 273, 268)$,
$(360, 285, 126, 53)$,\adfsplit
$(360, 76, 228, 187)$,
$(360, 249, 164, 216)$,
$(360, 62, 127, 27)$,
$(360, 157, 189, 356)$,\adfsplit
$(360, 63, 169, 292)$,
$(360, 294, 300, 252)$,
$(360, 256, 213, 226)$,
$(360, 210, 123, 86)$,\adfsplit
$(360, 83, 355, 262)$,
$(0, 7, 136, 227)$,
$(0, 64, 173, 398)$,
$(0, 19, 303, 442)$,\adfsplit
$(0, 46, 281, 381)$,
$(0, 11, 160, 448)$,
$(0, 10, 228, 476)$,
$(0, 40, 145, 407)$,\adfsplit
$(0, 38, 116, 440)$,
$(0, 22, 248, 478)$,
$(0, 29, 292, 389)$,
$(0, 82, 165, 386)$,\adfsplit
$(0, 44, 190, 439)$,
$(0, 16, 202, 441)$,
$(0, 28, 210, 371)$,
$(0, 24, 266, 399)$,\adfsplit
$(513, 0, 120, 240)$

\adfLgap \noindent by the mapping:
$x \mapsto x +  j \adfmod{360}$ for $x < 360$,
$x \mapsto (x +  j \adfmod{120}) + 360$ for $360 \le x < 480$,
$x \mapsto (x - 480 + 11 j \adfmod{33}) + 480$ for $480 \le x < 513$,
$513 \mapsto 513$,
$0 \le j < 360$
 for the first 52 blocks,
$0 \le j < 120$
 for the last block.
\ADFvfyParStart{(514, ((52, 360, ((360, 1), (120, 1), (33, 11), (1, 1))), (1, 120, ((360, 1), (120, 1), (33, 11), (1, 1)))), ((40, 9), (154, 1)))} 

\adfDgap
\noindent{\boldmath $ 40^{9} 157^{1} $}~
With the point set $Z_{517}$ partitioned into
 residue classes modulo $9$ for $\{0, 1, \dots, 359\}$, and
 $\{360, 361, \dots, 516\}$,
 the design is generated from

\adfLgap 
$(480, 191, 3, 302)$,
$(480, 355, 172, 224)$,
$(480, 268, 321, 68)$,
$(480, 307, 245, 31)$,\adfsplit
$(480, 161, 276, 142)$,
$(480, 200, 79, 197)$,
$(480, 294, 84, 40)$,
$(480, 193, 251, 298)$,\adfsplit
$(480, 64, 332, 157)$,
$(480, 76, 345, 215)$,
$(480, 238, 297, 137)$,
$(480, 182, 0, 102)$,\adfsplit
$(480, 77, 331, 290)$,
$(480, 54, 194, 170)$,
$(480, 168, 279, 176)$,
$(480, 164, 253, 277)$,\adfsplit
$(480, 42, 185, 111)$,
$(480, 271, 261, 346)$,
$(480, 66, 15, 46)$,
$(480, 252, 239, 195)$,\adfsplit
$(480, 99, 357, 165)$,
$(480, 18, 11, 106)$,
$(480, 131, 97, 278)$,
$(480, 160, 336, 1)$,\adfsplit
$(360, 359, 342, 155)$,
$(360, 305, 124, 228)$,
$(360, 162, 241, 210)$,
$(360, 101, 21, 154)$,\adfsplit
$(360, 189, 231, 19)$,
$(360, 345, 187, 5)$,
$(360, 253, 252, 6)$,
$(360, 217, 211, 84)$,\adfsplit
$(360, 336, 15, 47)$,
$(360, 355, 120, 81)$,
$(360, 257, 39, 117)$,
$(360, 118, 20, 340)$,\adfsplit
$(360, 264, 191, 268)$,
$(360, 44, 180, 174)$,
$(360, 22, 179, 3)$,
$(360, 9, 23, 254)$,\adfsplit
$(360, 326, 62, 136)$,
$(360, 106, 171, 280)$,
$(360, 46, 178, 296)$,
$(360, 151, 184, 138)$,\adfsplit
$(360, 250, 87, 224)$,
$(360, 78, 338, 53)$,
$(360, 70, 107, 219)$,
$(360, 277, 273, 11)$,\adfsplit
$(360, 92, 283, 267)$,
$(360, 88, 172, 229)$,
$(360, 110, 276, 72)$,
$(360, 7, 158, 94)$,\adfsplit
$(360, 197, 335, 234)$,
$(360, 213, 8, 103)$,
$(360, 170, 236, 41)$,
$(360, 30, 73, 329)$,\adfsplit
$(360, 289, 122, 29)$,
$(360, 322, 200, 203)$,
$(360, 165, 288, 205)$,
$(360, 266, 295, 57)$,\adfsplit
$(360, 152, 74, 63)$,
$(360, 265, 353, 306)$,
$(360, 68, 196, 301)$,
$(360, 319, 232, 75)$,\adfsplit
$(361, 96, 119, 325)$,
$(361, 87, 156, 133)$,
$(361, 317, 262, 9)$,
$(361, 203, 271, 94)$,\adfsplit
$(361, 136, 131, 39)$,
$(361, 299, 73, 337)$,
$(361, 316, 240, 79)$,
$(361, 254, 63, 268)$,\adfsplit
$(361, 124, 269, 145)$,
$(361, 175, 38, 201)$,
$(361, 306, 65, 328)$,
$(361, 301, 107, 250)$,\adfsplit
$(361, 161, 105, 70)$,
$(361, 352, 351, 187)$,
$(361, 322, 47, 324)$,
$(361, 354, 200, 52)$,\adfsplit
$(0, 5, 7, 469)$,
$(0, 15, 50, 345)$,
$(0, 193, 257, 363)$,
$(0, 33, 317, 389)$,\adfsplit
$(1, 9, 137, 381)$,
$(1, 61, 175, 447)$,
$(0, 73, 121, 407)$,
$(1, 23, 333, 361)$,\adfsplit
$(0, 11, 221, 455)$,
$(0, 211, 343, 431)$,
$(1, 71, 279, 429)$,
$(0, 281, 293, 419)$,\adfsplit
$(0, 67, 266, 449)$,
$(0, 49, 70, 473)$,
$(0, 86, 301, 379)$,
$(0, 124, 311, 381)$,\adfsplit
$(0, 16, 305, 433)$,
$(0, 34, 159, 475)$,
$(0, 112, 331, 383)$,
$(0, 60, 278, 375)$,\adfsplit
$(0, 28, 192, 437)$,
$(0, 62, 214, 421)$,
$(0, 12, 68, 439)$,
$(0, 10, 328, 461)$,\adfsplit
$(0, 30, 202, 425)$,
$(516, 0, 120, 240)$,
$(516, 1, 121, 241)$

\adfLgap \noindent by the mapping:
$x \mapsto x + 2 j \adfmod{360}$ for $x < 360$,
$x \mapsto (x + 2 j \adfmod{120}) + 360$ for $360 \le x < 480$,
$x \mapsto (x - 480 +  j \adfmod{36}) + 480$ for $480 \le x < 516$,
$516 \mapsto 516$,
$0 \le j < 180$
 for the first 105 blocks,
$0 \le j < 60$
 for the last two blocks.
\ADFvfyParStart{(517, ((105, 180, ((360, 2), (120, 2), (36, 1), (1, 1))), (2, 60, ((360, 2), (120, 2), (36, 1), (1, 1)))), ((40, 9), (157, 1)))} 

\section{4-GDDs for the proof of Lemma \ref{lem:4-GDD 46^u m^1}}
\label{app:4-GDD 46^u m^1}
\adfhide{
$ 46^6 13^1 $,
$ 46^6 16^1 $,
$ 46^6 19^1 $,
$ 46^6 22^1 $,
$ 46^6 73^1 $,
$ 46^6 76^1 $,
$ 46^6 79^1 $,
$ 46^6 82^1 $,
$ 46^6 85^1 $,
$ 46^6 88^1 $,
$ 46^6 91^1 $,
$ 46^6 94^1 $,
$ 46^6 97^1 $,
$ 46^6 100^1 $,
$ 46^6 103^1 $,
$ 46^6 106^1 $,
$ 46^6 109^1 $,
$ 46^6 112^1 $,
$ 46^9 19^1 $,
$ 46^{12} 25^1 $,
$ 46^{12} 28^1 $,
$ 46^{12} 31^1 $,
$ 46^{12} 34^1 $ and
$ 46^{15} 31^1 $.
}

\adfDgap
\noindent{\boldmath $ 46^{6} 13^{1} $}~
With the point set $Z_{289}$ partitioned into
 residue classes modulo $6$ for $\{0, 1, \dots, 275\}$, and
 $\{276, 277, \dots, 288\}$,
 the design is generated from

\adfLgap 
$(276, 273, 34, 221)$,
$(276, 63, 248, 187)$,
$(276, 229, 136, 246)$,
$(276, 107, 266, 60)$,\adfsplit
$(86, 243, 30, 211)$,
$(177, 162, 197, 68)$,
$(186, 267, 112, 212)$,
$(139, 54, 185, 141)$,\adfsplit
$(102, 59, 80, 177)$,
$(162, 165, 193, 28)$,
$(108, 172, 158, 159)$,
$(0, 4, 9, 80)$,\adfsplit
$(0, 29, 62, 223)$,
$(0, 39, 98, 199)$,
$(0, 10, 79, 146)$,
$(0, 16, 57, 143)$,\adfsplit
$(0, 7, 34, 45)$,
$(0, 19, 68, 173)$,
$(0, 23, 106, 164)$,
$(0, 8, 73, 177)$,\adfsplit
$(0, 25, 113, 153)$,
$(288, 0, 92, 184)$

\adfLgap \noindent by the mapping:
$x \mapsto x +  j \adfmod{276}$ for $x < 276$,
$x \mapsto (x +  j \adfmod{12}) + 276$ for $276 \le x < 288$,
$288 \mapsto 288$,
$0 \le j < 276$
 for the first 21 blocks,
$0 \le j < 92$
 for the last block.
\ADFvfyParStart{(289, ((21, 276, ((276, 1), (12, 1), (1, 1))), (1, 92, ((276, 1), (12, 1), (1, 1)))), ((46, 6), (13, 1)))} 

\adfDgap
\noindent{\boldmath $ 46^{6} 16^{1} $}~
With the point set $Z_{292}$ partitioned into
 residue classes modulo $6$ for $\{0, 1, \dots, 275\}$, and
 $\{276, 277, \dots, 291\}$,
 the design is generated from

\adfLgap 
$(288, 172, 74, 197)$,
$(288, 7, 186, 183)$,
$(276, 257, 268, 164)$,
$(276, 103, 105, 250)$,\adfsplit
$(276, 25, 35, 108)$,
$(276, 218, 66, 219)$,
$(277, 163, 182, 178)$,
$(277, 87, 258, 5)$,\adfsplit
$(277, 57, 155, 184)$,
$(277, 152, 240, 85)$,
$(236, 114, 75, 215)$,
$(45, 259, 178, 108)$,\adfsplit
$(80, 138, 154, 3)$,
$(123, 222, 172, 170)$,
$(4, 186, 80, 23)$,
$(26, 65, 174, 79)$,\adfsplit
$(11, 64, 132, 97)$,
$(11, 52, 199, 90)$,
$(194, 30, 241, 81)$,
$(160, 235, 192, 195)$,\adfsplit
$(258, 178, 235, 93)$,
$(53, 31, 159, 128)$,
$(133, 202, 3, 185)$,
$(206, 232, 201, 66)$,\adfsplit
$(123, 95, 127, 24)$,
$(22, 255, 101, 84)$,
$(79, 260, 268, 21)$,
$(0, 5, 236, 243)$,\adfsplit
$(0, 9, 22, 89)$,
$(0, 46, 251, 267)$,
$(0, 11, 55, 100)$,
$(0, 34, 107, 194)$,\adfsplit
$(0, 20, 64, 85)$,
$(0, 15, 158, 217)$,
$(0, 27, 145, 146)$,
$(0, 28, 119, 165)$,\adfsplit
$(0, 10, 169, 195)$,
$(0, 14, 175, 239)$,
$(0, 61, 69, 173)$,
$(0, 49, 117, 269)$,\adfsplit
$(0, 41, 134, 190)$,
$(0, 13, 63, 83)$,
$(0, 139, 215, 249)$,
$(291, 0, 92, 184)$,\adfsplit
$(291, 1, 93, 185)$

\adfLgap \noindent by the mapping:
$x \mapsto x + 2 j \adfmod{276}$ for $x < 276$,
$x \mapsto (x + 2 j \adfmod{12}) + 276$ for $276 \le x < 288$,
$x \mapsto (x +  j \adfmod{3}) + 288$ for $288 \le x < 291$,
$291 \mapsto 291$,
$0 \le j < 138$
 for the first 43 blocks;
$x \mapsto x + 2 j \adfmod{276}$ for $x < 276$,
$x \mapsto (x +  j \adfmod{12}) + 276$ for $276 \le x < 288$,
$x \mapsto (x +  j \adfmod{3}) + 288$ for $288 \le x < 291$,
$291 \mapsto 291$,
$0 \le j < 46$
 for the last two blocks.
\ADFvfyParStart{(292, ((43, 138, ((276, 2), (12, 2), (3, 1), (1, 1))), (2, 46, ((276, 2), (12, 1), (3, 1), (1, 1)))), ((46, 6), (16, 1)))} 

\adfDgap
\noindent{\boldmath $ 46^{6} 19^{1} $}~
With the point set $Z_{295}$ partitioned into
 residue classes modulo $6$ for $\{0, 1, \dots, 275\}$, and
 $\{276, 277, \dots, 294\}$,
 the design is generated from

\adfLgap 
$(288, 59, 27, 4)$,
$(288, 19, 30, 266)$,
$(276, 235, 120, 221)$,
$(276, 135, 114, 182)$,\adfsplit
$(276, 70, 152, 261)$,
$(276, 241, 119, 268)$,
$(51, 114, 52, 251)$,
$(128, 135, 185, 118)$,\adfsplit
$(96, 243, 38, 137)$,
$(94, 187, 213, 174)$,
$(191, 182, 1, 36)$,
$(196, 193, 168, 68)$,\adfsplit
$(0, 2, 113, 118)$,
$(0, 15, 31, 202)$,
$(0, 38, 91, 135)$,
$(0, 19, 83, 117)$,\adfsplit
$(0, 4, 56, 107)$,
$(0, 46, 133, 182)$,
$(0, 33, 70, 145)$,
$(0, 8, 69, 203)$,\adfsplit
$(0, 20, 79, 124)$,
$(0, 22, 65, 188)$,
$(294, 0, 92, 184)$

\adfLgap \noindent by the mapping:
$x \mapsto x +  j \adfmod{276}$ for $x < 276$,
$x \mapsto (x +  j \adfmod{12}) + 276$ for $276 \le x < 288$,
$x \mapsto (x +  j \adfmod{6}) + 288$ for $288 \le x < 294$,
$294 \mapsto 294$,
$0 \le j < 276$
 for the first 22 blocks,
$0 \le j < 92$
 for the last block.
\ADFvfyParStart{(295, ((22, 276, ((276, 1), (12, 1), (6, 1), (1, 1))), (1, 92, ((276, 1), (12, 1), (6, 1), (1, 1)))), ((46, 6), (19, 1)))} 

\adfDgap
\noindent{\boldmath $ 46^{6} 22^{1} $}~
With the point set $Z_{298}$ partitioned into
 residue classes modulo $6$ for $\{0, 1, \dots, 275\}$, and
 $\{276, 277, \dots, 297\}$,
 the design is generated from

\adfLgap 
$(294, 70, 37, 161)$,
$(294, 48, 56, 189)$,
$(276, 236, 154, 49)$,
$(276, 57, 112, 96)$,\adfsplit
$(276, 242, 185, 30)$,
$(276, 31, 23, 171)$,
$(277, 93, 184, 18)$,
$(277, 41, 25, 27)$,\adfsplit
$(277, 168, 248, 23)$,
$(277, 82, 79, 194)$,
$(278, 178, 249, 230)$,
$(278, 65, 164, 3)$,\adfsplit
$(278, 6, 229, 100)$,
$(278, 23, 0, 43)$,
$(50, 197, 123, 187)$,
$(195, 241, 8, 119)$,\adfsplit
$(170, 30, 87, 143)$,
$(189, 20, 124, 155)$,
$(212, 241, 114, 135)$,
$(84, 3, 22, 230)$,\adfsplit
$(146, 132, 35, 76)$,
$(18, 35, 56, 271)$,
$(108, 185, 68, 135)$,
$(160, 8, 114, 161)$,\adfsplit
$(77, 189, 256, 156)$,
$(259, 42, 129, 260)$,
$(254, 81, 169, 264)$,
$(79, 154, 161, 123)$,\adfsplit
$(254, 52, 39, 17)$,
$(268, 29, 87, 205)$,
$(0, 50, 207, 239)$,
$(0, 11, 63, 128)$,\adfsplit
$(0, 25, 53, 231)$,
$(0, 26, 81, 85)$,
$(0, 76, 183, 251)$,
$(0, 35, 225, 232)$,\adfsplit
$(0, 3, 83, 154)$,
$(0, 4, 167, 271)$,
$(0, 5, 34, 139)$,
$(0, 45, 118, 227)$,\adfsplit
$(0, 15, 142, 265)$,
$(0, 69, 86, 229)$,
$(0, 13, 28, 116)$,
$(0, 9, 119, 244)$,\adfsplit
$(0, 2, 101, 256)$,
$(297, 0, 92, 184)$,
$(297, 1, 93, 185)$

\adfLgap \noindent by the mapping:
$x \mapsto x + 2 j \adfmod{276}$ for $x < 276$,
$x \mapsto (x - 276 + 3 j \adfmod{18}) + 276$ for $276 \le x < 294$,
$x \mapsto (x +  j \adfmod{3}) + 294$ for $294 \le x < 297$,
$297 \mapsto 297$,
$0 \le j < 138$
 for the first 45 blocks,
$0 \le j < 46$
 for the last two blocks.
\ADFvfyParStart{(298, ((45, 138, ((276, 2), (18, 3), (3, 1), (1, 1))), (2, 46, ((276, 2), (18, 3), (3, 1), (1, 1)))), ((46, 6), (22, 1)))} 

\adfDgap
\noindent{\boldmath $ 46^{6} 73^{1} $}~
With the point set $Z_{349}$ partitioned into
 residue classes modulo $6$ for $\{0, 1, \dots, 275\}$, and
 $\{276, 277, \dots, 348\}$,
 the design is generated from

\adfLgap 
$(276, 187, 83, 207)$,
$(276, 225, 218, 12)$,
$(276, 25, 222, 208)$,
$(276, 29, 56, 190)$,\adfsplit
$(277, 117, 16, 263)$,
$(277, 188, 12, 211)$,
$(277, 193, 2, 75)$,
$(277, 222, 154, 89)$,\adfsplit
$(278, 126, 265, 255)$,
$(278, 212, 264, 10)$,
$(278, 52, 89, 187)$,
$(278, 225, 134, 131)$,\adfsplit
$(279, 16, 59, 33)$,
$(279, 2, 89, 195)$,
$(279, 150, 190, 247)$,
$(279, 73, 180, 224)$,\adfsplit
$(280, 120, 69, 242)$,
$(280, 76, 97, 164)$,
$(280, 111, 227, 70)$,
$(280, 210, 5, 43)$,\adfsplit
$(281, 86, 139, 105)$,
$(0, 1, 59, 128)$,
$(0, 4, 32, 163)$,
$(0, 11, 110, 123)$,\adfsplit
$(0, 49, 214, 281)$,
$(0, 2, 47, 82)$,
$(0, 5, 81, 287)$,
$(0, 15, 61, 136)$,\adfsplit
$(0, 9, 95, 317)$,
$(0, 16, 55, 105)$,
$(0, 8, 33, 64)$,
$(348, 0, 92, 184)$

\adfLgap \noindent by the mapping:
$x \mapsto x +  j \adfmod{276}$ for $x < 276$,
$x \mapsto (x - 276 + 6 j \adfmod{72}) + 276$ for $276 \le x < 348$,
$348 \mapsto 348$,
$0 \le j < 276$
 for the first 31 blocks,
$0 \le j < 92$
 for the last block.
\ADFvfyParStart{(349, ((31, 276, ((276, 1), (72, 6), (1, 1))), (1, 92, ((276, 1), (72, 6), (1, 1)))), ((46, 6), (73, 1)))} 

\adfDgap
\noindent{\boldmath $ 46^{6} 76^{1} $}~
With the point set $Z_{352}$ partitioned into
 residue classes modulo $6$ for $\{0, 1, \dots, 275\}$, and
 $\{276, 277, \dots, 351\}$,
 the design is generated from

\adfLgap 
$(348, 117, 13, 156)$,
$(348, 166, 167, 116)$,
$(276, 60, 235, 269)$,
$(276, 153, 275, 42)$,\adfsplit
$(276, 92, 123, 94)$,
$(276, 158, 265, 64)$,
$(277, 199, 125, 160)$,
$(277, 274, 158, 246)$,\adfsplit
$(277, 272, 252, 63)$,
$(277, 45, 121, 203)$,
$(278, 162, 171, 196)$,
$(278, 233, 260, 273)$,\adfsplit
$(278, 121, 98, 202)$,
$(278, 115, 275, 228)$,
$(279, 63, 52, 175)$,
$(279, 144, 239, 237)$,\adfsplit
$(279, 102, 41, 128)$,
$(279, 98, 94, 157)$,
$(280, 120, 43, 200)$,
$(280, 167, 57, 98)$,\adfsplit
$(280, 3, 16, 101)$,
$(280, 181, 270, 202)$,
$(281, 153, 98, 47)$,
$(281, 205, 156, 209)$,\adfsplit
$(281, 91, 184, 128)$,
$(281, 99, 238, 270)$,
$(282, 38, 109, 83)$,
$(282, 236, 189, 221)$,\adfsplit
$(282, 130, 204, 135)$,
$(282, 112, 90, 7)$,
$(283, 151, 23, 208)$,
$(283, 21, 50, 156)$,\adfsplit
$(283, 20, 205, 78)$,
$(283, 77, 94, 63)$,
$(284, 95, 274, 140)$,
$(284, 65, 159, 182)$,\adfsplit
$(284, 145, 252, 112)$,
$(284, 175, 153, 138)$,
$(285, 16, 131, 138)$,
$(285, 89, 109, 117)$,\adfsplit
$(285, 182, 252, 199)$,
$(285, 248, 118, 243)$,
$(286, 139, 183, 26)$,
$(286, 58, 209, 12)$,\adfsplit
$(286, 81, 54, 268)$,
$(286, 212, 97, 11)$,
$(287, 62, 3, 22)$,
$(0, 7, 16, 57)$,\adfsplit
$(0, 10, 200, 221)$,
$(0, 8, 155, 166)$,
$(0, 25, 77, 148)$,
$(0, 79, 149, 152)$,\adfsplit
$(0, 81, 139, 227)$,
$(1, 11, 135, 299)$,
$(0, 35, 98, 243)$,
$(0, 3, 100, 143)$,\adfsplit
$(0, 61, 99, 275)$,
$(0, 44, 117, 173)$,
$(0, 101, 165, 181)$,
$(0, 64, 167, 311)$,\adfsplit
$(0, 14, 52, 135)$,
$(0, 19, 65, 164)$,
$(0, 82, 191, 323)$,
$(351, 0, 92, 184)$,\adfsplit
$(351, 1, 93, 185)$

\adfLgap \noindent by the mapping:
$x \mapsto x + 2 j \adfmod{276}$ for $x < 276$,
$x \mapsto (x - 276 + 12 j \adfmod{72}) + 276$ for $276 \le x < 348$,
$x \mapsto (x +  j \adfmod{3}) + 348$ for $348 \le x < 351$,
$351 \mapsto 351$,
$0 \le j < 138$
 for the first 63 blocks,
$0 \le j < 46$
 for the last two blocks.
\ADFvfyParStart{(352, ((63, 138, ((276, 2), (72, 12), (3, 1), (1, 1))), (2, 46, ((276, 2), (72, 12), (3, 1), (1, 1)))), ((46, 6), (76, 1)))} 

\adfDgap
\noindent{\boldmath $ 46^{6} 79^{1} $}~
With the point set $Z_{355}$ partitioned into
 residue classes modulo $6$ for $\{0, 1, \dots, 275\}$, and
 $\{276, 277, \dots, 354\}$,
 the design is generated from

\adfLgap 
$(348, 107, 128, 21)$,
$(348, 91, 28, 90)$,
$(276, 210, 125, 201)$,
$(276, 178, 27, 248)$,\adfsplit
$(276, 193, 215, 230)$,
$(276, 7, 124, 60)$,
$(277, 48, 55, 10)$,
$(277, 6, 185, 242)$,\adfsplit
$(277, 76, 92, 273)$,
$(277, 217, 71, 99)$,
$(278, 28, 188, 223)$,
$(278, 106, 129, 217)$,\adfsplit
$(278, 95, 228, 26)$,
$(278, 138, 257, 243)$,
$(279, 137, 2, 69)$,
$(279, 171, 4, 181)$,\adfsplit
$(279, 48, 151, 95)$,
$(279, 258, 82, 80)$,
$(280, 154, 266, 79)$,
$(280, 47, 236, 73)$,\adfsplit
$(280, 129, 252, 257)$,
$(280, 4, 159, 198)$,
$(0, 4, 17, 249)$,
$(0, 25, 59, 287)$,\adfsplit
$(0, 3, 46, 152)$,
$(0, 8, 139, 299)$,
$(0, 29, 94, 305)$,
$(0, 11, 91, 140)$,\adfsplit
$(0, 19, 52, 218)$,
$(0, 32, 115, 311)$,
$(0, 50, 101, 172)$,
$(0, 20, 61, 203)$,\adfsplit
$(354, 0, 92, 184)$

\adfLgap \noindent by the mapping:
$x \mapsto x +  j \adfmod{276}$ for $x < 276$,
$x \mapsto (x - 276 + 6 j \adfmod{72}) + 276$ for $276 \le x < 348$,
$x \mapsto (x +  j \adfmod{6}) + 348$ for $348 \le x < 354$,
$354 \mapsto 354$,
$0 \le j < 276$
 for the first 32 blocks,
$0 \le j < 92$
 for the last block.
\ADFvfyParStart{(355, ((32, 276, ((276, 1), (72, 6), (6, 1), (1, 1))), (1, 92, ((276, 1), (72, 6), (6, 1), (1, 1)))), ((46, 6), (79, 1)))} 

\adfDgap
\noindent{\boldmath $ 46^{6} 82^{1} $}~
With the point set $Z_{358}$ partitioned into
 residue classes modulo $6$ for $\{0, 1, \dots, 275\}$, and
 $\{276, 277, \dots, 357\}$,
 the design is generated from

\adfLgap 
$(354, 141, 131, 210)$,
$(354, 116, 115, 196)$,
$(276, 76, 225, 108)$,
$(276, 97, 53, 170)$,\adfsplit
$(276, 234, 47, 58)$,
$(276, 103, 80, 171)$,
$(277, 85, 105, 35)$,
$(277, 139, 74, 243)$,\adfsplit
$(277, 137, 24, 196)$,
$(277, 20, 214, 258)$,
$(278, 190, 53, 133)$,
$(278, 189, 143, 184)$,\adfsplit
$(278, 188, 252, 123)$,
$(278, 2, 103, 54)$,
$(279, 262, 107, 200)$,
$(279, 121, 62, 195)$,\adfsplit
$(279, 90, 148, 127)$,
$(279, 216, 185, 261)$,
$(280, 257, 24, 103)$,
$(280, 126, 153, 268)$,\adfsplit
$(280, 164, 145, 142)$,
$(280, 11, 206, 27)$,
$(281, 117, 83, 62)$,
$(281, 67, 202, 90)$,\adfsplit
$(281, 112, 29, 140)$,
$(281, 217, 27, 72)$,
$(282, 176, 159, 36)$,
$(282, 28, 265, 29)$,\adfsplit
$(282, 66, 142, 237)$,
$(282, 59, 50, 175)$,
$(283, 129, 239, 1)$,
$(283, 38, 244, 123)$,\adfsplit
$(283, 120, 139, 106)$,
$(283, 164, 161, 210)$,
$(284, 209, 126, 75)$,
$(284, 32, 67, 154)$,\adfsplit
$(284, 23, 1, 141)$,
$(284, 110, 216, 256)$,
$(285, 33, 62, 172)$,
$(285, 71, 7, 272)$,\adfsplit
$(285, 0, 221, 229)$,
$(285, 54, 111, 58)$,
$(286, 98, 0, 249)$,
$(286, 219, 215, 104)$,\adfsplit
$(286, 100, 29, 199)$,
$(286, 246, 37, 262)$,
$(287, 125, 78, 25)$,
$(287, 251, 226, 69)$,\adfsplit
$(287, 206, 7, 156)$,
$(287, 172, 159, 260)$,
$(288, 55, 132, 57)$,
$(0, 74, 209, 261)$,\adfsplit
$(0, 73, 158, 267)$,
$(0, 20, 107, 163)$,
$(0, 7, 68, 131)$,
$(0, 13, 116, 177)$,\adfsplit
$(0, 15, 29, 124)$,
$(0, 2, 10, 41)$,
$(0, 56, 241, 269)$,
$(0, 86, 243, 314)$,\adfsplit
$(0, 17, 43, 105)$,
$(0, 71, 129, 327)$,
$(0, 93, 239, 271)$,
$(0, 34, 103, 128)$,\adfsplit
$(0, 26, 179, 353)$,
$(357, 0, 92, 184)$,
$(357, 1, 93, 185)$

\adfLgap \noindent by the mapping:
$x \mapsto x + 2 j \adfmod{276}$ for $x < 276$,
$x \mapsto (x - 276 + 13 j \adfmod{78}) + 276$ for $276 \le x < 354$,
$x \mapsto (x +  j \adfmod{3}) + 354$ for $354 \le x < 357$,
$357 \mapsto 357$,
$0 \le j < 138$
 for the first 65 blocks,
$0 \le j < 46$
 for the last two blocks.
\ADFvfyParStart{(358, ((65, 138, ((276, 2), (78, 13), (3, 1), (1, 1))), (2, 46, ((276, 2), (78, 13), (3, 1), (1, 1)))), ((46, 6), (82, 1)))} 

\adfDgap
\noindent{\boldmath $ 46^{6} 85^{1} $}~
With the point set $Z_{361}$ partitioned into
 residue classes modulo $6$ for $\{0, 1, \dots, 275\}$, and
 $\{276, 277, \dots, 360\}$,
 the design is generated from

\adfLgap 
$(276, 236, 95, 126)$,
$(276, 226, 189, 86)$,
$(276, 15, 139, 5)$,
$(276, 133, 252, 208)$,\adfsplit
$(277, 227, 255, 157)$,
$(277, 262, 69, 62)$,
$(277, 270, 221, 188)$,
$(277, 199, 172, 60)$,\adfsplit
$(278, 97, 272, 123)$,
$(278, 112, 55, 242)$,
$(278, 210, 239, 94)$,
$(278, 273, 29, 228)$,\adfsplit
$(279, 123, 98, 136)$,
$(279, 31, 46, 201)$,
$(279, 72, 179, 157)$,
$(279, 68, 162, 221)$,\adfsplit
$(280, 157, 275, 116)$,
$(280, 81, 209, 232)$,
$(280, 270, 219, 254)$,
$(280, 7, 24, 178)$,\adfsplit
$(281, 166, 152, 261)$,
$(281, 137, 98, 171)$,
$(281, 84, 227, 19)$,
$(0, 2, 273, 337)$,\adfsplit
$(0, 19, 40, 93)$,
$(0, 46, 115, 282)$,
$(0, 81, 185, 310)$,
$(0, 11, 61, 224)$,\adfsplit
$(0, 4, 62, 71)$,
$(0, 1, 189, 352)$,
$(0, 8, 55, 220)$,
$(0, 43, 129, 303)$,\adfsplit
$(0, 20, 99, 196)$,
$(360, 0, 92, 184)$

\adfLgap \noindent by the mapping:
$x \mapsto x +  j \adfmod{276}$ for $x < 276$,
$x \mapsto (x - 276 + 7 j \adfmod{84}) + 276$ for $276 \le x < 360$,
$360 \mapsto 360$,
$0 \le j < 276$
 for the first 33 blocks,
$0 \le j < 92$
 for the last block.
\ADFvfyParStart{(361, ((33, 276, ((276, 1), (84, 7), (1, 1))), (1, 92, ((276, 1), (84, 7), (1, 1)))), ((46, 6), (85, 1)))} 

\adfDgap
\noindent{\boldmath $ 46^{6} 88^{1} $}~
With the point set $Z_{364}$ partitioned into
 residue classes modulo $6$ for $\{0, 1, \dots, 275\}$, and
 $\{276, 277, \dots, 363\}$,
 the design is generated from

\adfLgap 
$(360, 42, 244, 206)$,
$(360, 195, 199, 23)$,
$(276, 265, 40, 218)$,
$(276, 138, 45, 77)$,\adfsplit
$(276, 155, 231, 236)$,
$(276, 211, 250, 24)$,
$(277, 170, 241, 168)$,
$(277, 267, 94, 270)$,\adfsplit
$(277, 153, 197, 235)$,
$(277, 167, 4, 92)$,
$(278, 60, 31, 118)$,
$(278, 245, 195, 126)$,\adfsplit
$(278, 237, 119, 188)$,
$(278, 74, 73, 196)$,
$(279, 157, 170, 94)$,
$(279, 138, 81, 268)$,\adfsplit
$(279, 11, 31, 200)$,
$(279, 192, 87, 269)$,
$(280, 21, 235, 89)$,
$(280, 150, 196, 59)$,\adfsplit
$(280, 182, 39, 156)$,
$(280, 85, 32, 46)$,
$(281, 126, 98, 166)$,
$(281, 171, 181, 40)$,\adfsplit
$(281, 180, 31, 17)$,
$(281, 44, 191, 249)$,
$(282, 180, 19, 155)$,
$(282, 51, 205, 94)$,\adfsplit
$(282, 42, 20, 5)$,
$(282, 2, 172, 21)$,
$(283, 42, 158, 199)$,
$(283, 251, 145, 28)$,\adfsplit
$(283, 250, 39, 260)$,
$(283, 245, 21, 264)$,
$(284, 264, 169, 8)$,
$(284, 189, 266, 222)$,\adfsplit
$(284, 226, 143, 115)$,
$(284, 148, 207, 173)$,
$(285, 59, 211, 66)$,
$(285, 243, 200, 34)$,\adfsplit
$(285, 208, 84, 61)$,
$(285, 33, 245, 230)$,
$(286, 72, 254, 181)$,
$(286, 139, 78, 213)$,\adfsplit
$(286, 99, 94, 101)$,
$(286, 8, 275, 40)$,
$(287, 14, 25, 95)$,
$(287, 3, 104, 88)$,\adfsplit
$(287, 30, 43, 34)$,
$(287, 185, 177, 156)$,
$(288, 217, 190, 69)$,
$(0, 3, 34, 91)$,\adfsplit
$(0, 17, 82, 177)$,
$(0, 86, 229, 255)$,
$(0, 97, 137, 172)$,
$(0, 8, 93, 148)$,\adfsplit
$(0, 31, 80, 289)$,
$(0, 56, 118, 344)$,
$(0, 35, 51, 206)$,
$(0, 83, 142, 345)$,\adfsplit
$(0, 52, 151, 249)$,
$(0, 149, 259, 331)$,
$(0, 1, 64, 231)$,
$(0, 23, 103, 303)$,\adfsplit
$(0, 37, 123, 179)$,
$(0, 101, 265, 330)$,
$(0, 45, 67, 316)$,
$(363, 0, 92, 184)$,\adfsplit
$(363, 1, 93, 185)$

\adfLgap \noindent by the mapping:
$x \mapsto x + 2 j \adfmod{276}$ for $x < 276$,
$x \mapsto (x - 276 + 14 j \adfmod{84}) + 276$ for $276 \le x < 360$,
$x \mapsto (x +  j \adfmod{3}) + 360$ for $360 \le x < 363$,
$363 \mapsto 363$,
$0 \le j < 138$
 for the first 67 blocks,
$0 \le j < 46$
 for the last two blocks.
\ADFvfyParStart{(364, ((67, 138, ((276, 2), (84, 14), (3, 1), (1, 1))), (2, 46, ((276, 2), (84, 14), (3, 1), (1, 1)))), ((46, 6), (88, 1)))} 

\adfDgap
\noindent{\boldmath $ 46^{6} 91^{1} $}~
With the point set $Z_{367}$ partitioned into
 residue classes modulo $6$ for $\{0, 1, \dots, 275\}$, and
 $\{276, 277, \dots, 366\}$,
 the design is generated from

\adfLgap 
$(360, 170, 17, 120)$,
$(360, 15, 187, 136)$,
$(276, 162, 136, 165)$,
$(276, 238, 84, 197)$,\adfsplit
$(276, 181, 47, 218)$,
$(276, 44, 243, 103)$,
$(277, 165, 34, 38)$,
$(277, 59, 159, 144)$,\adfsplit
$(277, 16, 91, 90)$,
$(277, 161, 128, 37)$,
$(278, 129, 212, 73)$,
$(278, 106, 99, 115)$,\adfsplit
$(278, 24, 35, 256)$,
$(278, 182, 246, 5)$,
$(279, 91, 53, 110)$,
$(279, 274, 123, 92)$,\adfsplit
$(279, 256, 246, 59)$,
$(279, 216, 261, 193)$,
$(280, 81, 222, 64)$,
$(280, 29, 235, 75)$,\adfsplit
$(280, 169, 74, 96)$,
$(280, 275, 262, 68)$,
$(281, 29, 198, 135)$,
$(281, 169, 8, 141)$,\adfsplit
$(0, 14, 81, 351)$,
$(0, 43, 190, 281)$,
$(0, 20, 52, 117)$,
$(0, 5, 39, 205)$,\adfsplit
$(0, 2, 49, 111)$,
$(0, 25, 189, 289)$,
$(0, 27, 175, 345)$,
$(0, 8, 88, 303)$,\adfsplit
$(0, 53, 146, 352)$,
$(0, 21, 61, 119)$,
$(366, 0, 92, 184)$

\adfLgap \noindent by the mapping:
$x \mapsto x +  j \adfmod{276}$ for $x < 276$,
$x \mapsto (x - 276 + 7 j \adfmod{84}) + 276$ for $276 \le x < 360$,
$x \mapsto (x +  j \adfmod{6}) + 360$ for $360 \le x < 366$,
$366 \mapsto 366$,
$0 \le j < 276$
 for the first 34 blocks,
$0 \le j < 92$
 for the last block.
\ADFvfyParStart{(367, ((34, 276, ((276, 1), (84, 7), (6, 1), (1, 1))), (1, 92, ((276, 1), (84, 7), (6, 1), (1, 1)))), ((46, 6), (91, 1)))} 

\adfDgap
\noindent{\boldmath $ 46^{6} 94^{1} $}~
With the point set $Z_{370}$ partitioned into
 residue classes modulo $6$ for $\{0, 1, \dots, 275\}$, and
 $\{276, 277, \dots, 369\}$,
 the design is generated from

\adfLgap 
$(366, 264, 196, 152)$,
$(366, 111, 85, 107)$,
$(276, 113, 268, 114)$,
$(276, 155, 117, 13)$,\adfsplit
$(276, 24, 104, 19)$,
$(276, 94, 242, 195)$,
$(277, 156, 208, 111)$,
$(277, 154, 210, 188)$,\adfsplit
$(277, 213, 151, 257)$,
$(277, 110, 107, 61)$,
$(278, 171, 55, 248)$,
$(278, 41, 146, 40)$,\adfsplit
$(278, 0, 169, 237)$,
$(278, 58, 90, 107)$,
$(279, 215, 172, 93)$,
$(279, 39, 176, 19)$,\adfsplit
$(279, 168, 205, 154)$,
$(279, 2, 101, 162)$,
$(280, 210, 85, 8)$,
$(280, 259, 257, 243)$,\adfsplit
$(280, 264, 256, 98)$,
$(280, 274, 105, 11)$,
$(281, 121, 24, 128)$,
$(281, 235, 159, 268)$,\adfsplit
$(281, 215, 50, 238)$,
$(281, 45, 162, 17)$,
$(282, 23, 144, 124)$,
$(282, 74, 138, 223)$,\adfsplit
$(282, 105, 32, 97)$,
$(282, 94, 135, 65)$,
$(283, 226, 44, 79)$,
$(283, 254, 189, 53)$,\adfsplit
$(283, 270, 87, 47)$,
$(283, 121, 4, 156)$,
$(284, 168, 7, 244)$,
$(284, 53, 226, 128)$,\adfsplit
$(284, 111, 174, 229)$,
$(284, 9, 50, 203)$,
$(285, 176, 139, 226)$,
$(285, 205, 42, 232)$,\adfsplit
$(285, 9, 65, 62)$,
$(285, 35, 60, 87)$,
$(286, 140, 94, 83)$,
$(286, 243, 48, 184)$,\adfsplit
$(286, 13, 86, 77)$,
$(286, 66, 213, 223)$,
$(287, 133, 59, 264)$,
$(287, 76, 259, 209)$,\adfsplit
$(287, 162, 177, 152)$,
$(287, 230, 15, 202)$,
$(288, 76, 217, 80)$,
$(288, 225, 120, 139)$,\adfsplit
$(288, 10, 123, 203)$,
$(288, 258, 29, 218)$,
$(289, 227, 258, 20)$,
$(0, 16, 233, 320)$,\adfsplit
$(0, 23, 123, 218)$,
$(0, 21, 206, 263)$,
$(0, 7, 185, 304)$,
$(0, 63, 187, 365)$,\adfsplit
$(0, 33, 143, 214)$,
$(0, 5, 135, 250)$,
$(0, 29, 177, 305)$,
$(0, 45, 100, 209)$,\adfsplit
$(0, 79, 111, 130)$,
$(0, 67, 82, 255)$,
$(0, 2, 11, 364)$,
$(0, 69, 127, 289)$,\adfsplit
$(0, 125, 142, 335)$,
$(369, 0, 92, 184)$,
$(369, 1, 93, 185)$

\adfLgap \noindent by the mapping:
$x \mapsto x + 2 j \adfmod{276}$ for $x < 276$,
$x \mapsto (x - 276 + 15 j \adfmod{90}) + 276$ for $276 \le x < 366$,
$x \mapsto (x +  j \adfmod{3}) + 366$ for $366 \le x < 369$,
$369 \mapsto 369$,
$0 \le j < 138$
 for the first 69 blocks,
$0 \le j < 46$
 for the last two blocks.
\ADFvfyParStart{(370, ((69, 138, ((276, 2), (90, 15), (3, 1), (1, 1))), (2, 46, ((276, 2), (90, 15), (3, 1), (1, 1)))), ((46, 6), (94, 1)))} 

\adfDgap
\noindent{\boldmath $ 46^{6} 97^{1} $}~
With the point set $Z_{373}$ partitioned into
 residue classes modulo $6$ for $\{0, 1, \dots, 275\}$, and
 $\{276, 277, \dots, 372\}$,
 the design is generated from

\adfLgap 
$(276, 220, 117, 248)$,
$(276, 154, 77, 147)$,
$(276, 2, 7, 234)$,
$(276, 0, 241, 83)$,\adfsplit
$(277, 160, 270, 47)$,
$(277, 46, 24, 45)$,
$(277, 110, 217, 77)$,
$(277, 127, 243, 188)$,\adfsplit
$(278, 7, 240, 142)$,
$(278, 121, 56, 167)$,
$(278, 237, 6, 136)$,
$(278, 242, 87, 89)$,\adfsplit
$(279, 170, 99, 118)$,
$(279, 72, 129, 197)$,
$(279, 203, 104, 91)$,
$(279, 241, 258, 172)$,\adfsplit
$(280, 31, 116, 203)$,
$(280, 27, 265, 4)$,
$(280, 114, 273, 26)$,
$(280, 12, 166, 197)$,\adfsplit
$(281, 124, 44, 83)$,
$(281, 165, 101, 175)$,
$(281, 61, 240, 243)$,
$(281, 218, 258, 274)$,\adfsplit
$(282, 113, 176, 139)$,
$(0, 20, 47, 79)$,
$(0, 8, 58, 362)$,
$(0, 34, 149, 290)$,\adfsplit
$(0, 81, 187, 330)$,
$(0, 11, 93, 283)$,
$(0, 51, 124, 347)$,
$(0, 25, 100, 167)$,\adfsplit
$(0, 14, 76, 171)$,
$(0, 4, 133, 363)$,
$(0, 9, 148, 371)$,
$(372, 0, 92, 184)$

\adfLgap \noindent by the mapping:
$x \mapsto x +  j \adfmod{276}$ for $x < 276$,
$x \mapsto (x - 276 + 8 j \adfmod{96}) + 276$ for $276 \le x < 372$,
$372 \mapsto 372$,
$0 \le j < 276$
 for the first 35 blocks,
$0 \le j < 92$
 for the last block.
\ADFvfyParStart{(373, ((35, 276, ((276, 1), (96, 8), (1, 1))), (1, 92, ((276, 1), (96, 8), (1, 1)))), ((46, 6), (97, 1)))} 

\adfDgap
\noindent{\boldmath $ 46^{6} 100^{1} $}~
With the point set $Z_{376}$ partitioned into
 residue classes modulo $5$ for $\{0, 1, \dots, 229\}$,
 $\{230, 231, \dots, 275\}$, and
 $\{276, 277, \dots, 375\}$,
 the design is generated from

\adfLgap 
$(326, 250, 218, 46)$,
$(362, 109, 88, 176)$,
$(335, 256, 181, 124)$,
$(366, 25, 189, 3)$,\adfsplit
$(292, 268, 137, 113)$,
$(281, 95, 166, 134)$,
$(288, 244, 18, 51)$,
$(333, 85, 153, 164)$,\adfsplit
$(297, 252, 178, 80)$,
$(359, 84, 85, 178)$,
$(300, 257, 75, 157)$,
$(336, 171, 93, 224)$,\adfsplit
$(306, 253, 32, 45)$,
$(335, 17, 60, 86)$,
$(348, 247, 15, 183)$,
$(292, 218, 134, 220)$,\adfsplit
$(293, 264, 188, 75)$,
$(314, 23, 4, 192)$,
$(374, 243, 56, 147)$,
$(319, 45, 49, 146)$,\adfsplit
$(329, 251, 134, 98)$,
$(320, 83, 146, 155)$,
$(339, 239, 74, 80)$,
$(371, 112, 10, 39)$,\adfsplit
$(353, 263, 103, 22)$,
$(305, 12, 66, 190)$,
$(322, 256, 19, 138)$,
$(277, 205, 191, 164)$,\adfsplit
$(0, 8, 59, 353)$,
$(0, 16, 231, 315)$,
$(0, 12, 242, 311)$,
$(0, 56, 246, 361)$,\adfsplit
$(0, 38, 154, 265)$,
$(0, 7, 103, 277)$,
$(0, 23, 240, 287)$,
$(0, 47, 244, 280)$,\adfsplit
$(0, 83, 266, 354)$,
$(0, 48, 122, 320)$,
$(0, 28, 77, 253)$,
$(0, 18, 64, 328)$,\adfsplit
$(0, 3, 37, 126)$,
$(0, 17, 109, 334)$,
$(0, 31, 143, 368)$

\adfLgap \noindent by the mapping:
$x \mapsto x +  j \adfmod{230}$ for $x < 230$,
$x \mapsto (x +  j \adfmod{46}) + 230$ for $230 \le x < 276$,
$x \mapsto (x - 276 + 10 j \adfmod{100}) + 276$ for $x \ge 276$,
$0 \le j < 230$.
\ADFvfyParStart{(376, ((43, 230, ((230, 1), (46, 1), (100, 10)))), ((46, 5), (46, 1), (100, 1)))} 

\adfDgap
\noindent{\boldmath $ 46^{6} 103^{1} $}~
With the point set $Z_{379}$ partitioned into
 residue classes modulo $6$ for $\{0, 1, \dots, 275\}$, and
 $\{276, 277, \dots, 378\}$,
 the design is generated from

\adfLgap 
$(372, 255, 146, 274)$,
$(372, 221, 138, 193)$,
$(276, 37, 226, 201)$,
$(276, 200, 0, 256)$,\adfsplit
$(276, 182, 151, 275)$,
$(276, 3, 269, 210)$,
$(277, 23, 52, 63)$,
$(277, 32, 217, 238)$,\adfsplit
$(277, 98, 132, 173)$,
$(277, 42, 115, 33)$,
$(278, 83, 184, 228)$,
$(278, 262, 205, 254)$,\adfsplit
$(278, 257, 140, 69)$,
$(278, 7, 3, 6)$,
$(279, 175, 160, 122)$,
$(279, 101, 234, 169)$,\adfsplit
$(279, 261, 200, 22)$,
$(279, 35, 120, 231)$,
$(280, 203, 219, 85)$,
$(280, 134, 4, 273)$,\adfsplit
$(280, 89, 43, 216)$,
$(280, 186, 188, 238)$,
$(281, 48, 25, 266)$,
$(281, 202, 235, 138)$,\adfsplit
$(281, 111, 128, 173)$,
$(281, 131, 136, 21)$,
$(0, 14, 171, 306)$,
$(0, 13, 113, 135)$,\adfsplit
$(0, 32, 99, 307)$,
$(0, 74, 169, 330)$,
$(0, 47, 170, 347)$,
$(0, 39, 125, 315)$,\adfsplit
$(0, 79, 160, 298)$,
$(0, 27, 182, 370)$,
$(0, 43, 147, 323)$,
$(0, 26, 77, 213)$,\adfsplit
$(378, 0, 92, 184)$

\adfLgap \noindent by the mapping:
$x \mapsto x +  j \adfmod{276}$ for $x < 276$,
$x \mapsto (x - 276 + 8 j \adfmod{96}) + 276$ for $276 \le x < 372$,
$x \mapsto (x +  j \adfmod{6}) + 372$ for $372 \le x < 378$,
$378 \mapsto 378$,
$0 \le j < 276$
 for the first 36 blocks,
$0 \le j < 92$
 for the last block.
\ADFvfyParStart{(379, ((36, 276, ((276, 1), (96, 8), (6, 1), (1, 1))), (1, 92, ((276, 1), (96, 8), (6, 1), (1, 1)))), ((46, 6), (103, 1)))} 

\adfDgap
\noindent{\boldmath $ 46^{6} 106^{1} $}~
With the point set $Z_{382}$ partitioned into
 residue classes modulo $6$ for $\{0, 1, \dots, 275\}$, and
 $\{276, 277, \dots, 381\}$,
 the design is generated from

\adfLgap 
$(378, 257, 187, 38)$,
$(378, 60, 28, 123)$,
$(276, 23, 104, 171)$,
$(276, 139, 182, 17)$,\adfsplit
$(276, 141, 16, 24)$,
$(276, 181, 246, 10)$,
$(277, 197, 39, 73)$,
$(277, 247, 274, 26)$,\adfsplit
$(277, 186, 260, 45)$,
$(277, 184, 83, 84)$,
$(278, 27, 58, 2)$,
$(278, 224, 60, 25)$,\adfsplit
$(278, 259, 100, 245)$,
$(278, 107, 45, 222)$,
$(279, 89, 231, 110)$,
$(279, 177, 176, 203)$,\adfsplit
$(279, 265, 196, 114)$,
$(279, 118, 216, 139)$,
$(280, 135, 247, 214)$,
$(280, 268, 273, 257)$,\adfsplit
$(280, 50, 72, 169)$,
$(280, 191, 174, 140)$,
$(281, 105, 142, 181)$,
$(281, 63, 162, 136)$,\adfsplit
$(281, 257, 164, 19)$,
$(281, 98, 144, 23)$,
$(282, 111, 96, 131)$,
$(282, 116, 169, 257)$,\adfsplit
$(282, 198, 16, 146)$,
$(282, 211, 189, 250)$,
$(283, 79, 218, 102)$,
$(283, 25, 256, 203)$,\adfsplit
$(283, 201, 245, 72)$,
$(283, 106, 87, 272)$,
$(284, 48, 239, 31)$,
$(284, 30, 77, 85)$,\adfsplit
$(284, 129, 142, 212)$,
$(284, 76, 134, 27)$,
$(285, 246, 74, 241)$,
$(285, 63, 137, 31)$,\adfsplit
$(285, 104, 179, 100)$,
$(285, 24, 33, 130)$,
$(286, 32, 63, 22)$,
$(286, 21, 217, 30)$,\adfsplit
$(286, 108, 119, 232)$,
$(286, 245, 38, 247)$,
$(287, 229, 189, 275)$,
$(287, 80, 168, 123)$,\adfsplit
$(287, 170, 173, 34)$,
$(287, 246, 208, 175)$,
$(288, 102, 215, 211)$,
$(288, 133, 17, 189)$,\adfsplit
$(288, 254, 256, 240)$,
$(288, 226, 51, 92)$,
$(289, 252, 83, 98)$,
$(289, 52, 81, 210)$,\adfsplit
$(289, 235, 154, 125)$,
$(289, 8, 73, 15)$,
$(0, 44, 225, 290)$,
$(0, 59, 208, 325)$,\adfsplit
$(0, 49, 185, 292)$,
$(0, 73, 76, 341)$,
$(0, 80, 183, 308)$,
$(0, 105, 115, 376)$,\adfsplit
$(0, 83, 133, 359)$,
$(0, 37, 226, 326)$,
$(0, 13, 20, 143)$,
$(0, 23, 86, 214)$,\adfsplit
$(0, 19, 71, 360)$,
$(0, 64, 153, 309)$,
$(0, 165, 229, 358)$,
$(0, 189, 217, 375)$,\adfsplit
$(0, 57, 157, 251)$,
$(381, 0, 92, 184)$,
$(381, 1, 93, 185)$

\adfLgap \noindent by the mapping:
$x \mapsto x + 2 j \adfmod{276}$ for $x < 276$,
$x \mapsto (x - 276 + 17 j \adfmod{102}) + 276$ for $276 \le x < 378$,
$x \mapsto (x +  j \adfmod{3}) + 378$ for $378 \le x < 381$,
$381 \mapsto 381$,
$0 \le j < 138$
 for the first 73 blocks,
$0 \le j < 46$
 for the last two blocks.
\ADFvfyParStart{(382, ((73, 138, ((276, 2), (102, 17), (3, 1), (1, 1))), (2, 46, ((276, 2), (102, 17), (3, 1), (1, 1)))), ((46, 6), (106, 1)))} 

\adfDgap
\noindent{\boldmath $ 46^{6} 109^{1} $}~
With the point set $Z_{385}$ partitioned into
 residue classes modulo $6$ for $\{0, 1, \dots, 275\}$, and
 $\{276, 277, \dots, 384\}$,
 the design is generated from

\adfLgap 
$(276, 174, 161, 200)$,
$(276, 211, 240, 215)$,
$(276, 265, 274, 111)$,
$(276, 64, 146, 165)$,\adfsplit
$(277, 108, 259, 227)$,
$(277, 207, 270, 173)$,
$(277, 268, 92, 177)$,
$(277, 61, 230, 94)$,\adfsplit
$(278, 198, 100, 266)$,
$(278, 267, 224, 19)$,
$(278, 245, 237, 10)$,
$(278, 155, 25, 108)$,\adfsplit
$(279, 44, 168, 261)$,
$(279, 169, 119, 99)$,
$(279, 114, 190, 125)$,
$(279, 232, 103, 218)$,\adfsplit
$(280, 198, 94, 147)$,
$(280, 140, 85, 213)$,
$(280, 84, 148, 179)$,
$(280, 115, 254, 137)$,\adfsplit
$(281, 71, 94, 132)$,
$(281, 76, 9, 229)$,
$(281, 79, 173, 158)$,
$(281, 186, 260, 87)$,\adfsplit
$(282, 29, 64, 74)$,
$(282, 103, 82, 191)$,
$(282, 45, 252, 200)$,
$(282, 61, 147, 66)$,\adfsplit
$(0, 3, 80, 310)$,
$(0, 1, 58, 293)$,
$(0, 37, 201, 301)$,
$(0, 40, 145, 373)$,\adfsplit
$(0, 16, 165, 355)$,
$(0, 7, 142, 365)$,
$(0, 2, 118, 356)$,
$(0, 46, 133, 329)$,\adfsplit
$(0, 17, 44, 106)$,
$(384, 0, 92, 184)$

\adfLgap \noindent by the mapping:
$x \mapsto x +  j \adfmod{276}$ for $x < 276$,
$x \mapsto (x - 276 + 9 j \adfmod{108}) + 276$ for $276 \le x < 384$,
$384 \mapsto 384$,
$0 \le j < 276$
 for the first 37 blocks,
$0 \le j < 92$
 for the last block.
\ADFvfyParStart{(385, ((37, 276, ((276, 1), (108, 9), (1, 1))), (1, 92, ((276, 1), (108, 9), (1, 1)))), ((46, 6), (109, 1)))} 

\adfDgap
\noindent{\boldmath $ 46^{6} 112^{1} $}~
With the point set $Z_{388}$ partitioned into
 residue classes modulo $6$ for $\{0, 1, \dots, 275\}$, and
 $\{276, 277, \dots, 387\}$,
 the design is generated from

\adfLgap 
$(384, 190, 258, 221)$,
$(384, 127, 75, 182)$,
$(276, 30, 211, 182)$,
$(276, 117, 52, 92)$,\adfsplit
$(276, 267, 17, 265)$,
$(276, 82, 216, 71)$,
$(277, 57, 116, 148)$,
$(277, 109, 113, 267)$,\adfsplit
$(277, 204, 239, 230)$,
$(277, 138, 82, 151)$,
$(278, 168, 235, 5)$,
$(278, 232, 191, 188)$,\adfsplit
$(278, 181, 81, 210)$,
$(278, 207, 226, 2)$,
$(279, 149, 62, 234)$,
$(279, 43, 178, 179)$,\adfsplit
$(279, 236, 88, 9)$,
$(279, 51, 133, 156)$,
$(280, 154, 120, 231)$,
$(280, 157, 5, 74)$,\adfsplit
$(280, 103, 138, 52)$,
$(280, 263, 189, 32)$,
$(281, 254, 269, 36)$,
$(281, 31, 111, 4)$,\adfsplit
$(281, 93, 157, 140)$,
$(281, 22, 59, 210)$,
$(282, 55, 93, 240)$,
$(282, 113, 152, 27)$,\adfsplit
$(282, 52, 71, 230)$,
$(282, 262, 150, 85)$,
$(283, 275, 261, 204)$,
$(283, 77, 250, 170)$,\adfsplit
$(283, 15, 121, 88)$,
$(283, 223, 80, 150)$,
$(284, 120, 9, 4)$,
$(284, 235, 63, 119)$,\adfsplit
$(284, 10, 49, 164)$,
$(284, 62, 126, 137)$,
$(285, 179, 124, 49)$,
$(285, 177, 156, 199)$,\adfsplit
$(285, 267, 92, 186)$,
$(285, 14, 34, 233)$,
$(286, 92, 59, 207)$,
$(286, 258, 268, 137)$,\adfsplit
$(286, 25, 26, 166)$,
$(286, 228, 235, 57)$,
$(287, 49, 132, 29)$,
$(287, 19, 162, 251)$,\adfsplit
$(287, 184, 122, 117)$,
$(287, 248, 118, 111)$,
$(288, 79, 132, 232)$,
$(288, 255, 42, 221)$,\adfsplit
$(288, 167, 266, 157)$,
$(288, 92, 142, 189)$,
$(289, 258, 67, 83)$,
$(289, 238, 68, 72)$,\adfsplit
$(289, 141, 53, 217)$,
$(289, 147, 26, 28)$,
$(290, 209, 20, 271)$,
$(290, 97, 105, 186)$,\adfsplit
$(0, 95, 127, 362)$,
$(0, 16, 38, 308)$,
$(0, 93, 227, 291)$,
$(0, 109, 159, 292)$,\adfsplit
$(0, 41, 249, 293)$,
$(0, 28, 74, 347)$,
$(0, 23, 63, 365)$,
$(0, 153, 263, 383)$,\adfsplit
$(0, 79, 261, 310)$,
$(0, 215, 273, 327)$,
$(0, 76, 225, 345)$,
$(0, 14, 59, 381)$,\adfsplit
$(0, 118, 255, 328)$,
$(0, 8, 61, 267)$,
$(0, 82, 245, 364)$,
$(387, 0, 92, 184)$,\adfsplit
$(387, 1, 93, 185)$

\adfLgap \noindent by the mapping:
$x \mapsto x + 2 j \adfmod{276}$ for $x < 276$,
$x \mapsto (x - 276 + 18 j \adfmod{108}) + 276$ for $276 \le x < 384$,
$x \mapsto (x +  j \adfmod{3}) + 384$ for $384 \le x < 387$,
$387 \mapsto 387$,
$0 \le j < 138$
 for the first 75 blocks,
$0 \le j < 46$
 for the last two blocks.
\ADFvfyParStart{(388, ((75, 138, ((276, 2), (108, 18), (3, 1), (1, 1))), (2, 46, ((276, 2), (108, 18), (3, 1), (1, 1)))), ((46, 6), (112, 1)))} 

\adfDgap
\noindent{\boldmath $ 46^{9} 19^{1} $}~
With the point set $Z_{433}$ partitioned into
 residue classes modulo $9$ for $\{0, 1, \dots, 413\}$, and
 $\{414, 415, \dots, 432\}$,
 the design is generated from

\adfLgap 
$(414, 15, 64, 408)$,
$(414, 3, 127, 86)$,
$(414, 241, 58, 329)$,
$(414, 42, 38, 224)$,\adfsplit
$(414, 178, 335, 180)$,
$(414, 99, 67, 53)$,
$(415, 82, 128, 303)$,
$(415, 271, 369, 399)$,\adfsplit
$(415, 112, 96, 396)$,
$(415, 239, 355, 335)$,
$(415, 250, 35, 156)$,
$(415, 169, 314, 230)$,\adfsplit
$(116, 121, 9, 3)$,
$(82, 101, 90, 34)$,
$(172, 250, 390, 8)$,
$(302, 375, 250, 333)$,\adfsplit
$(84, 60, 271, 40)$,
$(226, 183, 400, 378)$,
$(54, 278, 249, 146)$,
$(96, 20, 153, 75)$,\adfsplit
$(335, 346, 145, 336)$,
$(133, 95, 2, 186)$,
$(89, 129, 160, 24)$,
$(195, 81, 178, 280)$,\adfsplit
$(110, 116, 270, 391)$,
$(74, 9, 203, 253)$,
$(210, 133, 123, 297)$,
$(266, 225, 352, 385)$,\adfsplit
$(47, 318, 181, 409)$,
$(60, 266, 400, 345)$,
$(259, 264, 200, 91)$,
$(302, 156, 90, 141)$,\adfsplit
$(33, 148, 325, 387)$,
$(108, 193, 372, 230)$,
$(142, 83, 239, 1)$,
$(344, 181, 373, 147)$,\adfsplit
$(128, 190, 102, 197)$,
$(280, 308, 35, 268)$,
$(254, 28, 41, 88)$,
$(333, 298, 218, 345)$,\adfsplit
$(144, 202, 102, 393)$,
$(168, 269, 118, 338)$,
$(43, 312, 35, 373)$,
$(267, 83, 292, 295)$,\adfsplit
$(382, 11, 177, 225)$,
$(6, 390, 271, 315)$,
$(0, 1, 14, 308)$,
$(0, 15, 256, 380)$,\adfsplit
$(0, 68, 229, 246)$,
$(0, 77, 156, 332)$,
$(0, 37, 204, 302)$,
$(0, 53, 104, 222)$,\adfsplit
$(0, 38, 148, 341)$,
$(0, 39, 96, 379)$,
$(0, 89, 214, 325)$,
$(0, 75, 227, 347)$,\adfsplit
$(0, 25, 128, 267)$,
$(0, 71, 375, 391)$,
$(0, 185, 255, 321)$,
$(0, 237, 263, 395)$,\adfsplit
$(0, 23, 163, 381)$,
$(0, 109, 167, 259)$,
$(0, 113, 116, 319)$,
$(1, 3, 103, 107)$,\adfsplit
$(0, 165, 233, 313)$,
$(0, 123, 147, 407)$,
$(1, 23, 87, 233)$,
$(432, 0, 138, 276)$,\adfsplit
$(432, 1, 139, 277)$

\adfLgap \noindent by the mapping:
$x \mapsto x + 2 j \adfmod{414}$ for $x < 414$,
$x \mapsto (x + 2 j \adfmod{18}) + 414$ for $414 \le x < 432$,
$432 \mapsto 432$,
$0 \le j < 207$
 for the first 67 blocks,
$0 \le j < 69$
 for the last two blocks.
\ADFvfyParStart{(433, ((67, 207, ((414, 2), (18, 2), (1, 1))), (2, 69, ((414, 2), (18, 2), (1, 1)))), ((46, 9), (19, 1)))} 

\adfDgap
\noindent{\boldmath $ 46^{12} 25^{1} $}~
With the point set $Z_{577}$ partitioned into
 residue classes modulo $12$ for $\{0, 1, \dots, 551\}$, and
 $\{552, 553, \dots, 576\}$,
 the design is generated from

\adfLgap 
$(552, 401, 371, 93)$,
$(552, 19, 282, 36)$,
$(552, 95, 32, 367)$,
$(552, 225, 157, 266)$,\adfsplit
$(552, 166, 413, 507)$,
$(552, 6, 160, 313)$,
$(552, 543, 154, 4)$,
$(552, 168, 284, 518)$,\adfsplit
$(334, 344, 403, 541)$,
$(424, 461, 295, 480)$,
$(459, 499, 233, 501)$,
$(466, 481, 117, 254)$,\adfsplit
$(441, 466, 450, 179)$,
$(282, 383, 44, 73)$,
$(526, 273, 433, 291)$,
$(72, 263, 198, 481)$,\adfsplit
$(47, 338, 255, 444)$,
$(291, 211, 370, 188)$,
$(234, 176, 298, 137)$,
$(526, 104, 180, 436)$,\adfsplit
$(293, 304, 128, 421)$,
$(142, 291, 116, 186)$,
$(396, 369, 404, 75)$,
$(357, 75, 269, 157)$,\adfsplit
$(300, 266, 541, 102)$,
$(222, 297, 317, 368)$,
$(65, 61, 398, 8)$,
$(267, 300, 193, 245)$,\adfsplit
$(149, 34, 288, 148)$,
$(485, 351, 306, 84)$,
$(503, 13, 144, 322)$,
$(287, 373, 483, 100)$,\adfsplit
$(0, 3, 121, 232)$,
$(0, 46, 100, 236)$,
$(0, 21, 87, 148)$,
$(0, 85, 174, 315)$,\adfsplit
$(0, 6, 104, 239)$,
$(0, 5, 78, 255)$,
$(0, 67, 224, 394)$,
$(0, 14, 91, 433)$,\adfsplit
$(0, 43, 145, 292)$,
$(0, 47, 172, 295)$,
$(0, 7, 38, 221)$,
$(0, 28, 195, 381)$,\adfsplit
$(0, 49, 99, 400)$,
$(0, 32, 124, 205)$,
$(576, 0, 184, 368)$

\adfLgap \noindent by the mapping:
$x \mapsto x +  j \adfmod{552}$ for $x < 552$,
$x \mapsto (x +  j \adfmod{24}) + 552$ for $552 \le x < 576$,
$576 \mapsto 576$,
$0 \le j < 552$
 for the first 46 blocks,
$0 \le j < 184$
 for the last block.
\ADFvfyParStart{(577, ((46, 552, ((552, 1), (24, 1), (1, 1))), (1, 184, ((552, 1), (24, 1), (1, 1)))), ((46, 12), (25, 1)))} 

\adfDgap
\noindent{\boldmath $ 46^{12} 28^{1} $}~
With the point set $Z_{580}$ partitioned into
 residue classes modulo $12$ for $\{0, 1, \dots, 551\}$, and
 $\{552, 553, \dots, 579\}$,
 the design is generated from

\adfLgap 
$(576, 89, 36, 303)$,
$(576, 517, 362, 394)$,
$(552, 497, 114, 436)$,
$(552, 275, 39, 300)$,\adfsplit
$(552, 199, 489, 51)$,
$(552, 173, 397, 432)$,
$(552, 140, 112, 409)$,
$(552, 14, 310, 115)$,\adfsplit
$(552, 514, 93, 50)$,
$(552, 479, 414, 272)$,
$(553, 97, 242, 235)$,
$(553, 214, 245, 80)$,\adfsplit
$(553, 107, 391, 350)$,
$(553, 150, 420, 345)$,
$(553, 213, 47, 88)$,
$(553, 205, 220, 432)$,\adfsplit
$(553, 90, 363, 89)$,
$(553, 284, 82, 495)$,
$(388, 361, 194, 379)$,
$(44, 408, 79, 47)$,\adfsplit
$(297, 12, 342, 259)$,
$(367, 255, 66, 47)$,
$(266, 260, 281, 371)$,
$(277, 14, 247, 375)$,\adfsplit
$(140, 202, 506, 285)$,
$(413, 106, 287, 433)$,
$(356, 219, 522, 401)$,
$(459, 525, 395, 320)$,\adfsplit
$(517, 315, 246, 343)$,
$(51, 223, 323, 510)$,
$(211, 42, 389, 203)$,
$(293, 231, 73, 319)$,\adfsplit
$(266, 459, 208, 404)$,
$(485, 123, 441, 2)$,
$(108, 104, 484, 281)$,
$(398, 511, 438, 188)$,\adfsplit
$(498, 514, 373, 523)$,
$(367, 250, 8, 141)$,
$(220, 279, 157, 130)$,
$(178, 336, 19, 76)$,\adfsplit
$(450, 219, 95, 25)$,
$(252, 542, 547, 249)$,
$(72, 122, 79, 513)$,
$(190, 234, 448, 533)$,\adfsplit
$(26, 125, 226, 313)$,
$(61, 237, 332, 268)$,
$(530, 119, 130, 413)$,
$(369, 206, 295, 449)$,\adfsplit
$(385, 29, 45, 39)$,
$(307, 47, 89, 442)$,
$(353, 324, 275, 450)$,
$(310, 32, 438, 457)$,\adfsplit
$(530, 425, 379, 46)$,
$(446, 324, 537, 375)$,
$(414, 481, 36, 424)$,
$(291, 378, 364, 325)$,\adfsplit
$(358, 122, 159, 293)$,
$(21, 466, 515, 372)$,
$(327, 341, 395, 538)$,
$(490, 452, 121, 398)$,\adfsplit
$(481, 171, 371, 486)$,
$(313, 432, 398, 126)$,
$(257, 308, 0, 217)$,
$(178, 294, 180, 409)$,\adfsplit
$(137, 133, 189, 370)$,
$(340, 411, 155, 241)$,
$(515, 52, 10, 492)$,
$(94, 373, 175, 198)$,\adfsplit
$(0, 1, 18, 153)$,
$(0, 8, 159, 397)$,
$(0, 11, 13, 459)$,
$(0, 17, 339, 506)$,\adfsplit
$(0, 227, 449, 531)$,
$(0, 137, 277, 354)$,
$(0, 119, 261, 363)$,
$(0, 110, 289, 539)$,\adfsplit
$(0, 22, 351, 515)$,
$(0, 103, 239, 268)$,
$(0, 33, 100, 303)$,
$(0, 79, 232, 422)$,\adfsplit
$(0, 403, 497, 519)$,
$(0, 52, 421, 471)$,
$(0, 47, 346, 470)$,
$(0, 118, 266, 491)$,\adfsplit
$(0, 39, 305, 381)$,
$(0, 26, 245, 337)$,
$(0, 95, 255, 334)$,
$(0, 106, 215, 344)$,\adfsplit
$(0, 30, 328, 405)$,
$(0, 76, 317, 474)$,
$(0, 20, 170, 390)$,
$(0, 56, 136, 234)$,\adfsplit
$(0, 66, 140, 226)$,
$(579, 0, 184, 368)$,
$(579, 1, 185, 369)$

\adfLgap \noindent by the mapping:
$x \mapsto x + 2 j \adfmod{552}$ for $x < 552$,
$x \mapsto (x + 2 j \adfmod{24}) + 552$ for $552 \le x < 576$,
$x \mapsto (x +  j \adfmod{3}) + 576$ for $576 \le x < 579$,
$579 \mapsto 579$,
$0 \le j < 276$
 for the first 93 blocks,
$0 \le j < 92$
 for the last two blocks.
\ADFvfyParStart{(580, ((93, 276, ((552, 2), (24, 2), (3, 1), (1, 1))), (2, 92, ((552, 2), (24, 2), (3, 1), (1, 1)))), ((46, 12), (28, 1)))} 

\adfDgap
\noindent{\boldmath $ 46^{12} 31^{1} $}~
With the point set $Z_{583}$ partitioned into
 residue classes modulo $12$ for $\{0, 1, \dots, 551\}$, and
 $\{552, 553, \dots, 582\}$,
 the design is generated from

\adfLgap 
$(576, 349, 32, 105)$,
$(576, 420, 520, 89)$,
$(552, 367, 274, 257)$,
$(552, 261, 164, 287)$,\adfsplit
$(552, 324, 399, 18)$,
$(552, 174, 364, 125)$,
$(552, 440, 163, 338)$,
$(552, 0, 181, 395)$,\adfsplit
$(552, 537, 147, 358)$,
$(552, 302, 169, 184)$,
$(98, 419, 252, 400)$,
$(103, 481, 59, 34)$,\adfsplit
$(398, 543, 133, 444)$,
$(375, 36, 235, 238)$,
$(243, 148, 438, 494)$,
$(429, 232, 447, 395)$,\adfsplit
$(123, 353, 376, 418)$,
$(161, 320, 455, 120)$,
$(380, 157, 345, 358)$,
$(178, 303, 458, 189)$,\adfsplit
$(351, 70, 410, 318)$,
$(288, 350, 479, 190)$,
$(530, 64, 66, 260)$,
$(451, 545, 541, 394)$,\adfsplit
$(3, 23, 169, 282)$,
$(377, 511, 550, 191)$,
$(113, 493, 439, 362)$,
$(7, 370, 431, 116)$,\adfsplit
$(244, 181, 320, 101)$,
$(145, 114, 113, 423)$,
$(380, 447, 343, 36)$,
$(545, 171, 467, 298)$,\adfsplit
$(188, 307, 408, 413)$,
$(0, 6, 14, 224)$,
$(0, 45, 203, 267)$,
$(0, 16, 74, 461)$,\adfsplit
$(0, 10, 40, 509)$,
$(0, 55, 126, 207)$,
$(0, 9, 112, 236)$,
$(0, 29, 79, 470)$,\adfsplit
$(0, 27, 209, 260)$,
$(0, 28, 66, 465)$,
$(0, 7, 183, 268)$,
$(0, 116, 238, 402)$,\adfsplit
$(0, 21, 170, 375)$,
$(0, 47, 136, 365)$,
$(0, 68, 138, 435)$,
$(582, 0, 184, 368)$

\adfLgap \noindent by the mapping:
$x \mapsto x +  j \adfmod{552}$ for $x < 552$,
$x \mapsto (x +  j \adfmod{24}) + 552$ for $552 \le x < 576$,
$x \mapsto (x +  j \adfmod{6}) + 576$ for $576 \le x < 582$,
$582 \mapsto 582$,
$0 \le j < 552$
 for the first 47 blocks,
$0 \le j < 184$
 for the last block.
\ADFvfyParStart{(583, ((47, 552, ((552, 1), (24, 1), (6, 1), (1, 1))), (1, 184, ((552, 1), (24, 1), (6, 1), (1, 1)))), ((46, 12), (31, 1)))} 

\adfDgap
\noindent{\boldmath $ 46^{12} 34^{1} $}~
With the point set $Z_{586}$ partitioned into
 residue classes modulo $12$ for $\{0, 1, \dots, 551\}$, and
 $\{552, 553, \dots, 585\}$,
 the design is generated from

\adfLgap 
$(576, 39, 232, 251)$,
$(579, 40, 233, 252)$,
$(577, 155, 418, 186)$,
$(580, 156, 419, 187)$,\adfsplit
$(578, 383, 339, 139)$,
$(581, 384, 340, 140)$,
$(552, 225, 100, 467)$,
$(553, 226, 101, 468)$,\adfsplit
$(552, 459, 265, 536)$,
$(553, 460, 266, 537)$,
$(552, 497, 282, 283)$,
$(553, 498, 283, 284)$,\adfsplit
$(552, 538, 237, 101)$,
$(553, 539, 238, 102)$,
$(552, 408, 64, 102)$,
$(553, 409, 65, 103)$,\adfsplit
$(552, 218, 116, 79)$,
$(553, 219, 117, 80)$,
$(552, 156, 382, 95)$,
$(553, 157, 383, 96)$,\adfsplit
$(552, 253, 110, 327)$,
$(553, 254, 111, 328)$,
$(549, 100, 134, 299)$,
$(550, 101, 135, 300)$,\adfsplit
$(21, 493, 418, 47)$,
$(22, 494, 419, 48)$,
$(444, 269, 523, 11)$,
$(445, 270, 524, 12)$,\adfsplit
$(173, 27, 21, 79)$,
$(174, 28, 22, 80)$,
$(243, 486, 310, 139)$,
$(244, 487, 311, 140)$,\adfsplit
$(341, 42, 356, 271)$,
$(342, 43, 357, 272)$,
$(47, 146, 54, 113)$,
$(48, 147, 55, 114)$,\adfsplit
$(192, 246, 105, 292)$,
$(193, 247, 106, 293)$,
$(332, 134, 367, 137)$,
$(333, 135, 368, 138)$,\adfsplit
$(241, 525, 94, 371)$,
$(242, 526, 95, 372)$,
$(334, 103, 2, 306)$,
$(335, 104, 3, 307)$,\adfsplit
$(178, 296, 174, 47)$,
$(179, 297, 175, 48)$,
$(406, 374, 13, 5)$,
$(407, 375, 14, 6)$,\adfsplit
$(128, 307, 217, 16)$,
$(129, 308, 218, 17)$,
$(459, 41, 537, 510)$,
$(460, 42, 538, 511)$,\adfsplit
$(316, 440, 26, 298)$,
$(317, 441, 27, 299)$,
$(60, 234, 541, 320)$,
$(61, 235, 542, 321)$,\adfsplit
$(463, 544, 277, 275)$,
$(464, 545, 278, 276)$,
$(243, 81, 480, 221)$,
$(244, 82, 481, 222)$,\adfsplit
$(428, 364, 311, 475)$,
$(429, 365, 312, 476)$,
$(317, 260, 4, 99)$,
$(318, 261, 5, 100)$,\adfsplit
$(155, 92, 305, 410)$,
$(156, 93, 306, 411)$,
$(351, 341, 508, 181)$,
$(352, 342, 509, 182)$,\adfsplit
$(131, 161, 120, 170)$,
$(132, 162, 121, 171)$,
$(506, 151, 83, 156)$,
$(0, 13, 342, 356)$,\adfsplit
$(0, 16, 98, 123)$,
$(0, 29, 43, 206)$,
$(0, 113, 222, 429)$,
$(0, 62, 211, 407)$,\adfsplit
$(0, 49, 163, 374)$,
$(0, 23, 158, 487)$,
$(0, 65, 88, 370)$,
$(0, 73, 419, 436)$,\adfsplit
$(0, 69, 343, 386)$,
$(0, 91, 257, 527)$,
$(0, 45, 114, 286)$,
$(0, 20, 325, 423)$,\adfsplit
$(1, 21, 179, 351)$,
$(0, 145, 273, 355)$,
$(0, 128, 445, 507)$,
$(0, 173, 363, 439)$,\adfsplit
$(0, 133, 224, 503)$,
$(0, 5, 316, 333)$,
$(0, 97, 523, 539)$,
$(0, 42, 110, 219)$,\adfsplit
$(0, 21, 278, 404)$,
$(0, 93, 135, 190)$,
$(0, 55, 76, 459)$,
$(585, 0, 184, 368)$,\adfsplit
$(585, 1, 185, 369)$

\adfLgap \noindent by the mapping:
$x \mapsto x + 2 j \adfmod{552}$ for $x < 552$,
$x \mapsto (x + 2 j \adfmod{24}) + 552$ for $552 \le x < 576$,
$x \mapsto (x + 6 j \adfmod{9}) + 576$ for $576 \le x < 585$,
$585 \mapsto 585$,
$0 \le j < 276$
 for the first 95 blocks,
$0 \le j < 92$
 for the last two blocks.
\ADFvfyParStart{(586, ((95, 276, ((552, 2), (24, 2), (9, 6), (1, 1))), (2, 92, ((552, 2), (24, 2), (9, 6), (1, 1)))), ((46, 12), (34, 1)))} 

\adfDgap
\noindent{\boldmath $ 46^{15} 31^{1} $}~
With the point set $Z_{721}$ partitioned into
 residue classes modulo $15$ for $\{0, 1, \dots, 689\}$, and
 $\{690, 691, \dots, 720\}$,
 the design is generated from

\adfLgap 
$(690, 463, 460, 229)$,
$(691, 464, 461, 230)$,
$(690, 414, 388, 147)$,
$(691, 415, 389, 148)$,\adfsplit
$(690, 247, 475, 91)$,
$(691, 248, 476, 92)$,
$(690, 292, 623, 291)$,
$(691, 293, 624, 292)$,\adfsplit
$(690, 326, 518, 35)$,
$(691, 327, 519, 36)$,
$(690, 434, 311, 258)$,
$(691, 435, 312, 259)$,\adfsplit
$(690, 60, 539, 426)$,
$(691, 61, 540, 427)$,
$(690, 422, 500, 72)$,
$(691, 423, 501, 73)$,\adfsplit
$(690, 399, 136, 617)$,
$(691, 400, 137, 618)$,
$(690, 484, 645, 303)$,
$(691, 485, 646, 304)$,\adfsplit
$(467, 139, 151, 33)$,
$(468, 140, 152, 34)$,
$(81, 143, 453, 282)$,
$(82, 144, 454, 283)$,\adfsplit
$(128, 16, 176, 557)$,
$(129, 17, 177, 558)$,
$(7, 193, 527, 184)$,
$(8, 194, 528, 185)$,\adfsplit
$(289, 186, 40, 317)$,
$(290, 187, 41, 318)$,
$(430, 320, 107, 247)$,
$(431, 321, 108, 248)$,\adfsplit
$(670, 344, 395, 68)$,
$(671, 345, 396, 69)$,
$(13, 459, 676, 2)$,
$(14, 460, 677, 3)$,\adfsplit
$(261, 358, 430, 434)$,
$(262, 359, 431, 435)$,
$(400, 301, 366, 319)$,
$(401, 302, 367, 320)$,\adfsplit
$(502, 201, 116, 343)$,
$(503, 202, 117, 344)$,
$(11, 633, 442, 572)$,
$(12, 634, 443, 573)$,\adfsplit
$(398, 19, 671, 279)$,
$(399, 20, 672, 280)$,
$(337, 83, 54, 150)$,
$(338, 84, 55, 151)$,\adfsplit
$(25, 682, 645, 227)$,
$(26, 683, 646, 228)$,
$(423, 676, 659, 337)$,
$(424, 677, 660, 338)$,\adfsplit
$(406, 18, 332, 120)$,
$(407, 19, 333, 121)$,
$(236, 288, 307, 293)$,
$(237, 289, 308, 294)$,\adfsplit
$(234, 142, 553, 257)$,
$(235, 143, 554, 258)$,
$(153, 109, 627, 26)$,
$(154, 110, 628, 27)$,\adfsplit
$(141, 13, 524, 221)$,
$(142, 14, 525, 222)$,
$(680, 490, 672, 343)$,
$(681, 491, 673, 344)$,\adfsplit
$(211, 235, 352, 11)$,
$(212, 236, 353, 12)$,
$(635, 511, 674, 266)$,
$(636, 512, 675, 267)$,\adfsplit
$(613, 592, 18, 414)$,
$(614, 593, 19, 415)$,
$(158, 26, 224, 193)$,
$(159, 27, 225, 194)$,\adfsplit
$(595, 329, 398, 155)$,
$(596, 330, 399, 156)$,
$(549, 559, 374, 415)$,
$(550, 560, 375, 416)$,\adfsplit
$(351, 394, 242, 39)$,
$(352, 395, 243, 40)$,
$(154, 612, 320, 557)$,
$(155, 613, 321, 558)$,\adfsplit
$(671, 366, 457, 503)$,
$(672, 367, 458, 504)$,
$(275, 490, 652, 544)$,
$(276, 491, 653, 545)$,\adfsplit
$(227, 314, 469, 221)$,
$(228, 315, 470, 222)$,
$(102, 138, 503, 419)$,
$(103, 139, 504, 420)$,\adfsplit
$(505, 247, 341, 405)$,
$(0, 7, 133, 471)$,
$(0, 25, 395, 579)$,
$(0, 13, 151, 425)$,\adfsplit
$(0, 295, 393, 677)$,
$(0, 49, 288, 545)$,
$(0, 63, 501, 533)$,
$(1, 57, 205, 411)$,\adfsplit
$(1, 3, 23, 291)$,
$(0, 251, 439, 497)$,
$(0, 107, 189, 229)$,
$(0, 143, 247, 297)$,\adfsplit
$(0, 223, 419, 461)$,
$(0, 121, 274, 467)$,
$(0, 122, 412, 565)$,
$(0, 111, 333, 410)$,\adfsplit
$(0, 221, 246, 485)$,
$(0, 136, 569, 683)$,
$(0, 59, 138, 484)$,
$(0, 158, 422, 627)$,\adfsplit
$(0, 56, 157, 308)$,
$(0, 2, 79, 470)$,
$(0, 50, 148, 502)$,
$(0, 40, 82, 576)$,\adfsplit
$(0, 58, 284, 641)$,
$(0, 32, 370, 589)$,
$(0, 64, 164, 258)$,
$(0, 22, 293, 586)$,\adfsplit
$(0, 20, 145, 204)$,
$(720, 0, 230, 460)$,
$(720, 1, 231, 461)$

\adfLgap \noindent by the mapping:
$x \mapsto x + 2 j \adfmod{690}$ for $x < 690$,
$x \mapsto (x + 2 j \adfmod{30}) + 690$ for $690 \le x < 720$,
$720 \mapsto 720$,
$0 \le j < 345$
 for the first 117 blocks,
$0 \le j < 115$
 for the last two blocks.
\ADFvfyParStart{(721, ((117, 345, ((690, 2), (30, 2), (1, 1))), (2, 115, ((690, 2), (30, 2), (1, 1)))), ((46, 15), (31, 1)))} 

\section{4-GDDs for the proof of Lemma \ref{lem:4-GDD 50^u m^1}}
\label{app:4-GDD 50^u m^1}
\adfhide{
$ 50^6 8^1 $,
$ 50^6 11^1 $,
$ 50^6 14^1 $,
$ 50^6 17^1 $,
$ 50^6 23^1 $,
$ 50^6 26^1 $,
$ 50^6 29^1 $,
$ 50^6 32^1 $,
$ 50^6 38^1 $,
$ 50^6 113^1 $,
$ 50^6 116^1 $,
$ 50^6 119^1 $,
$ 50^6 122^1 $,
$ 50^9 11^1 $,
$ 50^9 17^1 $ and
$ 50^9 23^1 $.
}

\adfDgap
\noindent{\boldmath $ 50^{6} 8^{1} $}~
With the point set $Z_{308}$ partitioned into
 residue classes modulo $6$ for $\{0, 1, \dots, 299\}$, and
 $\{300, 301, \dots, 307\}$,
 the design is generated from

\adfLgap 
$(300, 0, 1, 2)$,
$(301, 0, 100, 98)$,
$(302, 0, 199, 299)$,
$(303, 0, 202, 101)$,\adfsplit
$(304, 0, 298, 200)$,
$(305, 0, 4, 11)$,
$(306, 0, 7, 296)$,
$(307, 0, 289, 293)$,\adfsplit
$(136, 59, 145, 56)$,
$(91, 68, 274, 119)$,
$(268, 108, 71, 189)$,
$(273, 98, 36, 197)$,\adfsplit
$(62, 221, 243, 12)$,
$(230, 102, 160, 295)$,
$(59, 111, 284, 175)$,
$(189, 234, 70, 158)$,\adfsplit
$(31, 281, 10, 87)$,
$(155, 214, 171, 102)$,
$(206, 166, 179, 21)$,
$(201, 191, 38, 283)$,\adfsplit
$(54, 103, 177, 64)$,
$(76, 215, 48, 67)$,
$(272, 219, 6, 175)$,
$(274, 188, 139, 29)$,\adfsplit
$(215, 253, 93, 26)$,
$(81, 112, 56, 73)$,
$(144, 164, 179, 15)$,
$(0, 5, 8, 45)$,\adfsplit
$(0, 14, 253, 268)$,
$(0, 22, 176, 205)$,
$(0, 16, 194, 235)$,
$(0, 26, 64, 157)$,\adfsplit
$(0, 87, 133, 208)$,
$(0, 61, 74, 195)$,
$(0, 71, 129, 283)$,
$(0, 33, 104, 148)$,\adfsplit
$(0, 27, 116, 259)$,
$(0, 79, 149, 182)$,
$(0, 63, 97, 190)$,
$(0, 85, 105, 281)$,\adfsplit
$(0, 43, 68, 215)$,
$(0, 83, 109, 166)$,
$(0, 47, 82, 177)$,
$(0, 52, 111, 125)$,\adfsplit
$(0, 201, 233, 295)$,
$(0, 91, 142, 279)$,
$(0, 57, 170, 277)$

\adfLgap \noindent by the mapping:
$x \mapsto x + 3 j \adfmod{300}$ for $x < 300$,
$x \mapsto x$ for $x \ge 300$,
$0 \le j < 100$
 for the first eight blocks;
$x \mapsto x + 2 j \adfmod{300}$ for $x < 300$,
$x \mapsto x$ for $x \ge 300$,
$0 \le j < 150$
 for the last 39 blocks.
\ADFvfyParStart{(308, ((8, 100, ((300, 3), (8, 8))), (39, 150, ((300, 2), (8, 8)))), ((50, 6), (8, 1)))} 

\adfDgap
\noindent{\boldmath $ 50^{6} 11^{1} $}~
With the point set $Z_{311}$ partitioned into
 residue classes modulo $6$ for $\{0, 1, \dots, 299\}$, and
 $\{300, 301, \dots, 310\}$,
 the design is generated from

\adfLgap 
$(306, 0, 1, 2)$,
$(307, 0, 100, 98)$,
$(308, 0, 199, 299)$,
$(309, 0, 202, 101)$,\adfsplit
$(310, 0, 298, 200)$,
$(300, 274, 225, 102)$,
$(300, 109, 215, 56)$,
$(144, 88, 5, 157)$,\adfsplit
$(22, 295, 107, 228)$,
$(151, 125, 232, 278)$,
$(113, 246, 75, 205)$,
$(288, 249, 220, 260)$,\adfsplit
$(46, 32, 123, 65)$,
$(45, 28, 107, 20)$,
$(181, 290, 166, 35)$,
$(0, 3, 23, 268)$,\adfsplit
$(0, 4, 65, 122)$,
$(0, 5, 64, 230)$,
$(0, 7, 50, 143)$,
$(0, 22, 73, 117)$,\adfsplit
$(0, 9, 89, 160)$,
$(0, 47, 110, 184)$,
$(0, 34, 86, 189)$,
$(0, 10, 31, 135)$,\adfsplit
$(0, 16, 113, 201)$,
$(0, 37, 119, 195)$

\adfLgap \noindent by the mapping:
$x \mapsto x + 3 j \adfmod{300}$ for $x < 300$,
$x \mapsto (x +  j \adfmod{6}) + 300$ for $300 \le x < 306$,
$x \mapsto x$ for $x \ge 306$,
$0 \le j < 100$
 for the first five blocks;
$x \mapsto x +  j \adfmod{300}$ for $x < 300$,
$x \mapsto (x +  j \adfmod{6}) + 300$ for $300 \le x < 306$,
$x \mapsto x$ for $x \ge 306$,
$0 \le j < 300$
 for the last 21 blocks.
\ADFvfyParStart{(311, ((5, 100, ((300, 3), (6, 1), (5, 5))), (21, 300, ((300, 1), (6, 1), (5, 5)))), ((50, 6), (11, 1)))} 

\adfDgap
\noindent{\boldmath $ 50^{6} 14^{1} $}~
With the point set $Z_{314}$ partitioned into
 residue classes modulo $6$ for $\{0, 1, \dots, 299\}$, and
 $\{300, 301, \dots, 313\}$,
 the design is generated from

\adfLgap 
$(306, 0, 1, 2)$,
$(307, 0, 100, 98)$,
$(308, 0, 199, 299)$,
$(309, 0, 202, 101)$,\adfsplit
$(310, 0, 298, 200)$,
$(311, 0, 4, 11)$,
$(312, 0, 7, 296)$,
$(313, 0, 289, 293)$,\adfsplit
$(300, 202, 229, 275)$,
$(300, 266, 162, 100)$,
$(300, 273, 132, 17)$,
$(300, 272, 163, 75)$,\adfsplit
$(131, 208, 193, 158)$,
$(173, 289, 117, 18)$,
$(159, 190, 38, 299)$,
$(48, 179, 93, 73)$,\adfsplit
$(75, 174, 113, 260)$,
$(157, 5, 216, 286)$,
$(65, 284, 228, 112)$,
$(3, 282, 142, 31)$,\adfsplit
$(204, 256, 123, 149)$,
$(105, 178, 90, 145)$,
$(185, 153, 122, 106)$,
$(118, 98, 31, 81)$,\adfsplit
$(239, 56, 96, 3)$,
$(237, 288, 79, 245)$,
$(127, 258, 190, 69)$,
$(92, 241, 286, 159)$,\adfsplit
$(0, 3, 122, 127)$,
$(0, 8, 17, 85)$,
$(0, 10, 32, 69)$,
$(0, 13, 26, 291)$,\adfsplit
$(0, 14, 58, 209)$,
$(0, 23, 92, 211)$,
$(0, 33, 113, 236)$,
$(0, 83, 157, 279)$,\adfsplit
$(0, 28, 71, 163)$,
$(0, 38, 112, 235)$,
$(0, 29, 170, 295)$,
$(0, 41, 51, 176)$,\adfsplit
$(0, 65, 117, 136)$,
$(0, 105, 146, 275)$,
$(0, 75, 145, 221)$,
$(0, 39, 53, 206)$,\adfsplit
$(0, 61, 118, 225)$,
$(0, 87, 110, 205)$,
$(0, 34, 251, 254)$,
$(0, 93, 187, 203)$,\adfsplit
$(0, 57, 142, 218)$

\adfLgap \noindent by the mapping:
$x \mapsto x + 3 j \adfmod{300}$ for $x < 300$,
$x \mapsto (x +  j \adfmod{6}) + 300$ for $300 \le x < 306$,
$x \mapsto x$ for $x \ge 306$,
$0 \le j < 100$
 for the first eight blocks;
$x \mapsto x + 2 j \adfmod{300}$ for $x < 300$,
$x \mapsto (x +  j \adfmod{6}) + 300$ for $300 \le x < 306$,
$x \mapsto x$ for $x \ge 306$,
$0 \le j < 150$
 for the last 41 blocks.
\ADFvfyParStart{(314, ((8, 100, ((300, 3), (6, 1), (8, 8))), (41, 150, ((300, 2), (6, 1), (8, 8)))), ((50, 6), (14, 1)))} 

\adfDgap
\noindent{\boldmath $ 50^{6} 17^{1} $}~
With the point set $Z_{317}$ partitioned into
 residue classes modulo $6$ for $\{0, 1, \dots, 299\}$, and
 $\{300, 301, \dots, 316\}$,
 the design is generated from

\adfLgap 
$(312, 0, 1, 2)$,
$(313, 0, 100, 98)$,
$(314, 0, 199, 299)$,
$(315, 0, 202, 101)$,\adfsplit
$(316, 0, 298, 200)$,
$(300, 260, 93, 157)$,
$(300, 6, 209, 160)$,
$(301, 0, 191, 169)$,\adfsplit
$(301, 50, 112, 255)$,
$(200, 131, 253, 66)$,
$(110, 6, 155, 291)$,
$(148, 113, 164, 187)$,\adfsplit
$(17, 202, 169, 206)$,
$(167, 238, 14, 43)$,
$(143, 200, 6, 189)$,
$(125, 42, 290, 295)$,\adfsplit
$(0, 3, 17, 110)$,
$(0, 7, 28, 68)$,
$(0, 20, 70, 245)$,
$(0, 13, 38, 94)$,\adfsplit
$(0, 19, 63, 142)$,
$(0, 10, 77, 218)$,
$(0, 27, 85, 116)$,
$(0, 32, 73, 160)$,\adfsplit
$(0, 9, 43, 182)$,
$(0, 26, 112, 171)$,
$(0, 8, 88, 209)$

\adfLgap \noindent by the mapping:
$x \mapsto x + 3 j \adfmod{300}$ for $x < 300$,
$x \mapsto (x + 2 j \adfmod{12}) + 300$ for $300 \le x < 312$,
$x \mapsto x$ for $x \ge 312$,
$0 \le j < 100$
 for the first five blocks;
$x \mapsto x +  j \adfmod{300}$ for $x < 300$,
$x \mapsto (x + 2 j \adfmod{12}) + 300$ for $300 \le x < 312$,
$x \mapsto x$ for $x \ge 312$,
$0 \le j < 300$
 for the last 22 blocks.
\ADFvfyParStart{(317, ((5, 100, ((300, 3), (12, 2), (5, 5))), (22, 300, ((300, 1), (12, 2), (5, 5)))), ((50, 6), (17, 1)))} 

\adfDgap
\noindent{\boldmath $ 50^{6} 23^{1} $}~
With the point set $Z_{323}$ partitioned into
 residue classes modulo $6$ for $\{0, 1, \dots, 299\}$, and
 $\{300, 301, \dots, 322\}$,
 the design is generated from

\adfLgap 
$(318, 0, 1, 2)$,
$(319, 0, 100, 98)$,
$(320, 0, 199, 299)$,
$(321, 0, 202, 101)$,\adfsplit
$(322, 0, 298, 200)$,
$(300, 286, 83, 207)$,
$(300, 84, 206, 265)$,
$(301, 3, 226, 265)$,\adfsplit
$(301, 83, 138, 20)$,
$(302, 288, 58, 227)$,
$(302, 243, 116, 133)$,
$(197, 128, 216, 82)$,\adfsplit
$(228, 93, 256, 253)$,
$(198, 208, 235, 123)$,
$(19, 53, 106, 174)$,
$(244, 103, 72, 255)$,\adfsplit
$(270, 176, 292, 219)$,
$(0, 4, 33, 259)$,
$(0, 5, 13, 57)$,
$(0, 14, 106, 205)$,\adfsplit
$(0, 9, 32, 274)$,
$(0, 15, 71, 154)$,
$(0, 16, 81, 211)$,
$(0, 21, 107, 157)$,\adfsplit
$(0, 20, 123, 187)$,
$(0, 49, 125, 207)$,
$(0, 7, 47, 196)$,
$(0, 62, 129, 209)$

\adfLgap \noindent by the mapping:
$x \mapsto x + 3 j \adfmod{300}$ for $x < 300$,
$x \mapsto (x - 300 + 3 j \adfmod{18}) + 300$ for $300 \le x < 318$,
$x \mapsto x$ for $x \ge 318$,
$0 \le j < 100$
 for the first five blocks;
$x \mapsto x +  j \adfmod{300}$ for $x < 300$,
$x \mapsto (x - 300 + 3 j \adfmod{18}) + 300$ for $300 \le x < 318$,
$x \mapsto x$ for $x \ge 318$,
$0 \le j < 300$
 for the last 23 blocks.
\ADFvfyParStart{(323, ((5, 100, ((300, 3), (18, 3), (5, 5))), (23, 300, ((300, 1), (18, 3), (5, 5)))), ((50, 6), (23, 1)))} 

\adfDgap
\noindent{\boldmath $ 50^{6} 26^{1} $}~
With the point set $Z_{326}$ partitioned into
 residue classes modulo $6$ for $\{0, 1, \dots, 299\}$, and
 $\{300, 301, \dots, 325\}$,
 the design is generated from

\adfLgap 
$(318, 0, 1, 2)$,
$(319, 0, 100, 98)$,
$(320, 0, 199, 299)$,
$(321, 0, 202, 101)$,\adfsplit
$(322, 0, 298, 200)$,
$(323, 0, 4, 11)$,
$(324, 0, 7, 296)$,
$(325, 0, 289, 293)$,\adfsplit
$(300, 203, 266, 69)$,
$(300, 126, 58, 139)$,
$(300, 220, 27, 36)$,
$(300, 5, 253, 296)$,\adfsplit
$(301, 119, 86, 36)$,
$(301, 67, 39, 164)$,
$(301, 273, 142, 121)$,
$(301, 210, 53, 256)$,\adfsplit
$(302, 249, 110, 174)$,
$(302, 68, 217, 12)$,
$(302, 106, 17, 235)$,
$(302, 244, 207, 11)$,\adfsplit
$(223, 261, 44, 287)$,
$(177, 107, 242, 282)$,
$(109, 46, 176, 299)$,
$(7, 153, 233, 134)$,\adfsplit
$(271, 4, 290, 150)$,
$(297, 36, 140, 53)$,
$(35, 262, 174, 211)$,
$(161, 44, 210, 237)$,\adfsplit
$(0, 22, 47, 32)$,
$(296, 45, 203, 58)$,
$(38, 219, 259, 82)$,
$(198, 274, 251, 26)$,\adfsplit
$(0, 3, 194, 274)$,
$(0, 5, 8, 93)$,
$(0, 16, 71, 230)$,
$(0, 21, 65, 136)$,\adfsplit
$(0, 23, 159, 188)$,
$(0, 34, 69, 77)$,
$(0, 61, 95, 111)$,
$(0, 155, 201, 223)$,\adfsplit
$(0, 39, 124, 169)$,
$(0, 167, 187, 273)$,
$(0, 113, 127, 182)$,
$(0, 51, 92, 239)$,\adfsplit
$(0, 79, 110, 197)$,
$(0, 59, 91, 219)$,
$(0, 57, 119, 241)$,
$(0, 82, 171, 265)$,\adfsplit
$(0, 20, 142, 295)$,
$(0, 31, 41, 226)$,
$(0, 28, 191, 283)$,
$(0, 38, 189, 247)$,\adfsplit
$(0, 58, 152, 285)$

\adfLgap \noindent by the mapping:
$x \mapsto x + 3 j \adfmod{300}$ for $x < 300$,
$x \mapsto (x - 300 + 3 j \adfmod{18}) + 300$ for $300 \le x < 318$,
$x \mapsto x$ for $x \ge 318$,
$0 \le j < 100$
 for the first eight blocks;
$x \mapsto x + 2 j \adfmod{300}$ for $x < 300$,
$x \mapsto (x - 300 + 3 j \adfmod{18}) + 300$ for $300 \le x < 318$,
$x \mapsto x$ for $x \ge 318$,
$0 \le j < 150$
 for the last 45 blocks.
\ADFvfyParStart{(326, ((8, 100, ((300, 3), (18, 3), (8, 8))), (45, 150, ((300, 2), (18, 3), (8, 8)))), ((50, 6), (26, 1)))} 

\adfDgap
\noindent{\boldmath $ 50^{6} 29^{1} $}~
With the point set $Z_{329}$ partitioned into
 residue classes modulo $6$ for $\{0, 1, \dots, 299\}$, and
 $\{300, 301, \dots, 328\}$,
 the design is generated from

\adfLgap 
$(324, 0, 1, 2)$,
$(325, 0, 100, 98)$,
$(326, 0, 199, 299)$,
$(327, 0, 202, 101)$,\adfsplit
$(328, 0, 298, 200)$,
$(300, 275, 87, 103)$,
$(300, 138, 194, 280)$,
$(301, 94, 234, 63)$,\adfsplit
$(301, 187, 176, 245)$,
$(302, 292, 21, 125)$,
$(302, 54, 199, 272)$,
$(303, 75, 125, 26)$,\adfsplit
$(303, 18, 175, 292)$,
$(20, 227, 295, 117)$,
$(107, 222, 194, 63)$,
$(269, 254, 163, 166)$,\adfsplit
$(180, 142, 33, 109)$,
$(144, 197, 189, 274)$,
$(0, 4, 67, 74)$,
$(0, 14, 65, 219)$,\adfsplit
$(0, 35, 92, 211)$,
$(0, 17, 111, 166)$,
$(0, 13, 47, 165)$,
$(0, 9, 32, 193)$,\adfsplit
$(0, 5, 46, 225)$,
$(0, 19, 39, 260)$,
$(0, 21, 43, 238)$,
$(0, 10, 37, 173)$,\adfsplit
$(0, 52, 113, 177)$

\adfLgap \noindent by the mapping:
$x \mapsto x + 3 j \adfmod{300}$ for $x < 300$,
$x \mapsto (x - 300 + 4 j \adfmod{24}) + 300$ for $300 \le x < 324$,
$x \mapsto x$ for $x \ge 324$,
$0 \le j < 100$
 for the first five blocks;
$x \mapsto x +  j \adfmod{300}$ for $x < 300$,
$x \mapsto (x - 300 + 4 j \adfmod{24}) + 300$ for $300 \le x < 324$,
$x \mapsto x$ for $x \ge 324$,
$0 \le j < 300$
 for the last 24 blocks.
\ADFvfyParStart{(329, ((5, 100, ((300, 3), (24, 4), (5, 5))), (24, 300, ((300, 1), (24, 4), (5, 5)))), ((50, 6), (29, 1)))} 

\adfDgap
\noindent{\boldmath $ 50^{6} 32^{1} $}~
With the point set $Z_{332}$ partitioned into
 residue classes modulo $6$ for $\{0, 1, \dots, 299\}$, and
 $\{300, 301, \dots, 331\}$,
 the design is generated from

\adfLgap 
$(324, 0, 1, 2)$,
$(325, 0, 100, 98)$,
$(326, 0, 199, 299)$,
$(327, 0, 202, 101)$,\adfsplit
$(328, 0, 298, 200)$,
$(329, 0, 4, 11)$,
$(330, 0, 7, 296)$,
$(331, 0, 289, 293)$,\adfsplit
$(300, 234, 183, 2)$,
$(300, 216, 248, 151)$,
$(300, 52, 299, 141)$,
$(300, 101, 298, 37)$,\adfsplit
$(301, 167, 271, 76)$,
$(301, 85, 266, 48)$,
$(301, 210, 281, 237)$,
$(301, 188, 255, 10)$,\adfsplit
$(302, 0, 134, 113)$,
$(302, 70, 23, 246)$,
$(302, 200, 243, 151)$,
$(302, 4, 213, 133)$,\adfsplit
$(303, 170, 153, 131)$,
$(303, 66, 257, 280)$,
$(303, 58, 55, 20)$,
$(303, 240, 13, 231)$,\adfsplit
$(103, 160, 93, 252)$,
$(109, 244, 161, 0)$,
$(127, 216, 53, 32)$,
$(131, 254, 264, 249)$,\adfsplit
$(285, 155, 216, 175)$,
$(91, 63, 242, 131)$,
$(41, 154, 282, 195)$,
$(199, 184, 224, 269)$,\adfsplit
$(140, 195, 88, 59)$,
$(257, 284, 204, 130)$,
$(0, 3, 16, 241)$,
$(0, 5, 139, 255)$,\adfsplit
$(0, 8, 57, 133)$,
$(0, 13, 29, 182)$,
$(0, 9, 17, 164)$,
$(0, 14, 34, 197)$,\adfsplit
$(0, 70, 175, 212)$,
$(0, 26, 76, 123)$,
$(0, 28, 79, 207)$,
$(0, 22, 115, 237)$,\adfsplit
$(0, 155, 169, 201)$,
$(0, 117, 143, 229)$,
$(0, 61, 104, 185)$,
$(0, 58, 157, 188)$,\adfsplit
$(0, 44, 106, 265)$,
$(0, 23, 236, 267)$,
$(0, 33, 193, 227)$,
$(0, 25, 63, 190)$,\adfsplit
$(0, 75, 223, 281)$,
$(0, 46, 111, 206)$,
$(0, 19, 148, 231)$

\adfLgap \noindent by the mapping:
$x \mapsto x + 3 j \adfmod{300}$ for $x < 300$,
$x \mapsto (x - 300 + 4 j \adfmod{24}) + 300$ for $300 \le x < 324$,
$x \mapsto x$ for $x \ge 324$,
$0 \le j < 100$
 for the first eight blocks;
$x \mapsto x + 2 j \adfmod{300}$ for $x < 300$,
$x \mapsto (x - 300 + 4 j \adfmod{24}) + 300$ for $300 \le x < 324$,
$x \mapsto x$ for $x \ge 324$,
$0 \le j < 150$
 for the last 47 blocks.
\ADFvfyParStart{(332, ((8, 100, ((300, 3), (24, 4), (8, 8))), (47, 150, ((300, 2), (24, 4), (8, 8)))), ((50, 6), (32, 1)))} 

\adfDgap
\noindent{\boldmath $ 50^{6} 38^{1} $}~
With the point set $Z_{338}$ partitioned into
 residue classes modulo $6$ for $\{0, 1, \dots, 299\}$, and
 $\{300, 301, \dots, 337\}$,
 the design is generated from

\adfLgap 
$(330, 0, 1, 2)$,
$(331, 0, 100, 98)$,
$(332, 0, 199, 299)$,
$(333, 0, 202, 101)$,\adfsplit
$(334, 0, 298, 200)$,
$(335, 0, 4, 11)$,
$(336, 0, 7, 296)$,
$(337, 0, 289, 293)$,\adfsplit
$(300, 206, 1, 107)$,
$(300, 261, 202, 290)$,
$(300, 44, 138, 189)$,
$(300, 228, 225, 53)$,\adfsplit
$(300, 210, 100, 51)$,
$(300, 76, 109, 183)$,
$(300, 263, 294, 130)$,
$(300, 218, 46, 240)$,\adfsplit
$(300, 156, 237, 296)$,
$(300, 192, 7, 15)$,
$(300, 238, 246, 65)$,
$(300, 223, 244, 299)$,\adfsplit
$(300, 11, 265, 200)$,
$(300, 277, 197, 92)$,
$(300, 172, 271, 41)$,
$(300, 8, 162, 274)$,\adfsplit
$(300, 84, 219, 254)$,
$(300, 235, 269, 153)$,
$(300, 28, 215, 19)$,
$(300, 193, 207, 122)$,\adfsplit
$(90, 207, 194, 77)$,
$(241, 185, 46, 294)$,
$(222, 254, 249, 125)$,
$(109, 78, 135, 167)$,\adfsplit
$(276, 193, 11, 44)$,
$(29, 70, 186, 123)$,
$(185, 117, 160, 102)$,
$(245, 40, 30, 133)$,\adfsplit
$(174, 157, 93, 113)$,
$(0, 3, 148, 157)$,
$(0, 5, 20, 43)$,
$(0, 14, 63, 250)$,\adfsplit
$(0, 17, 67, 238)$,
$(0, 44, 178, 275)$,
$(0, 21, 211, 233)$,
$(0, 77, 225, 235)$,\adfsplit
$(0, 37, 159, 262)$,
$(0, 39, 175, 230)$,
$(0, 47, 86, 207)$,
$(0, 40, 177, 263)$,\adfsplit
$(0, 61, 80, 153)$,
$(0, 29, 56, 277)$,
$(0, 19, 76, 158)$,
$(0, 91, 131, 182)$,\adfsplit
$(0, 45, 92, 179)$,
$(0, 16, 167, 229)$,
$(0, 41, 69, 176)$,
$(0, 28, 74, 255)$,\adfsplit
$(0, 26, 173, 189)$

\adfLgap \noindent by the mapping:
$x \mapsto x + 3 j \adfmod{300}$ for $x < 300$,
$x \mapsto (x +  j \adfmod{30}) + 300$ for $300 \le x < 330$,
$x \mapsto x$ for $x \ge 330$,
$0 \le j < 100$
 for the first eight blocks;
$x \mapsto x + 2 j \adfmod{300}$ for $x < 300$,
$x \mapsto (x +  j \adfmod{30}) + 300$ for $300 \le x < 330$,
$x \mapsto x$ for $x \ge 330$,
$0 \le j < 150$
 for the last 49 blocks.
\ADFvfyParStart{(338, ((8, 100, ((300, 3), (30, 1), (8, 8))), (49, 150, ((300, 2), (30, 1), (8, 8)))), ((50, 6), (38, 1)))} 

\adfDgap
\noindent{\boldmath $ 50^{6} 113^{1} $}~
With the point set $Z_{413}$ partitioned into
 residue classes modulo $6$ for $\{0, 1, \dots, 299\}$, and
 $\{300, 301, \dots, 412\}$,
 the design is generated from

\adfLgap 
$(408, 0, 1, 2)$,
$(409, 0, 100, 98)$,
$(410, 0, 199, 299)$,
$(411, 0, 202, 101)$,\adfsplit
$(412, 0, 298, 200)$,
$(390, 39, 90, 73)$,
$(390, 113, 200, 64)$,
$(391, 166, 153, 276)$,\adfsplit
$(391, 173, 176, 187)$,
$(392, 91, 146, 124)$,
$(392, 0, 59, 201)$,
$(300, 152, 58, 33)$,\adfsplit
$(300, 275, 159, 168)$,
$(300, 247, 294, 287)$,
$(300, 96, 106, 241)$,
$(300, 87, 50, 119)$,\adfsplit
$(300, 109, 191, 124)$,
$(300, 141, 224, 162)$,
$(300, 270, 23, 193)$,
$(300, 26, 165, 52)$,\adfsplit
$(300, 55, 218, 100)$,
$(301, 222, 142, 183)$,
$(301, 133, 214, 81)$,
$(301, 272, 96, 7)$,\adfsplit
$(301, 290, 198, 271)$,
$(301, 116, 145, 256)$,
$(301, 23, 0, 129)$,
$(301, 8, 251, 139)$,\adfsplit
$(301, 209, 87, 178)$,
$(301, 250, 54, 125)$,
$(0, 27, 88, 340)$,
$(0, 20, 63, 166)$,\adfsplit
$(0, 5, 157, 323)$,
$(0, 58, 230, 305)$,
$(0, 16, 44, 362)$,
$(0, 76, 173, 386)$,\adfsplit
$(0, 4, 50, 359)$,
$(0, 74, 149, 329)$,
$(0, 86, 195, 371)$,
$(0, 64, 159, 380)$,\adfsplit
$(0, 56, 121, 344)$,
$(0, 8, 93, 374)$,
$(0, 38, 117, 185)$

\adfLgap \noindent by the mapping:
$x \mapsto x + 3 j \adfmod{300}$ for $x < 300$,
$x \mapsto (x - 300 + 3 j \adfmod{90}) + 300$ for $300 \le x < 390$,
$x \mapsto (x - 390 + 3 j \adfmod{18}) + 390$ for $390 \le x < 408$,
$x \mapsto x$ for $x \ge 408$,
$0 \le j < 100$
 for the first five blocks;
$x \mapsto x +  j \adfmod{300}$ for $x < 300$,
$x \mapsto (x - 300 + 3 j \adfmod{90}) + 300$ for $300 \le x < 390$,
$x \mapsto (x - 390 + 3 j \adfmod{18}) + 390$ for $390 \le x < 408$,
$x \mapsto x$ for $x \ge 408$,
$0 \le j < 300$
 for the last 38 blocks.
\ADFvfyParStart{(413, ((5, 100, ((300, 3), (90, 3), (18, 3), (5, 5))), (38, 300, ((300, 1), (90, 3), (18, 3), (5, 5)))), ((50, 6), (113, 1)))} 

\adfDgap
\noindent{\boldmath $ 50^{6} 116^{1} $}~
With the point set $Z_{416}$ partitioned into
 residue classes modulo $6$ for $\{0, 1, \dots, 299\}$, and
 $\{300, 301, \dots, 415\}$,
 the design is generated from

\adfLgap 
$(405, 0, 1, 2)$,
$(406, 0, 100, 98)$,
$(407, 0, 199, 299)$,
$(408, 0, 202, 101)$,\adfsplit
$(409, 0, 298, 200)$,
$(410, 0, 4, 11)$,
$(411, 0, 7, 296)$,
$(412, 0, 289, 293)$,\adfsplit
$(413, 0, 10, 23)$,
$(414, 0, 13, 290)$,
$(415, 0, 277, 287)$,
$(300, 287, 46, 282)$,\adfsplit
$(300, 138, 31, 51)$,
$(300, 32, 54, 255)$,
$(300, 53, 13, 123)$,
$(300, 8, 237, 155)$,\adfsplit
$(300, 299, 88, 175)$,
$(300, 224, 232, 169)$,
$(300, 10, 127, 266)$,
$(300, 0, 80, 244)$,\adfsplit
$(300, 276, 219, 71)$,
$(301, 55, 240, 249)$,
$(301, 217, 132, 225)$,
$(301, 50, 231, 23)$,\adfsplit
$(301, 2, 16, 123)$,
$(301, 134, 64, 17)$,
$(301, 218, 41, 280)$,
$(301, 169, 56, 264)$,\adfsplit
$(301, 118, 66, 209)$,
$(301, 82, 288, 61)$,
$(301, 27, 215, 193)$,
$(302, 19, 210, 95)$,\adfsplit
$(302, 213, 256, 239)$,
$(302, 157, 173, 154)$,
$(302, 87, 43, 48)$,
$(302, 162, 287, 39)$,\adfsplit
$(302, 202, 176, 281)$,
$(302, 260, 156, 285)$,
$(302, 104, 121, 84)$,
$(302, 160, 145, 248)$,\adfsplit
$(302, 92, 261, 208)$,
$(303, 104, 25, 256)$,
$(303, 239, 20, 180)$,
$(303, 225, 54, 8)$,\adfsplit
$(303, 103, 32, 281)$,
$(303, 298, 116, 143)$,
$(303, 211, 250, 108)$,
$(303, 37, 6, 201)$,\adfsplit
$(303, 249, 274, 17)$,
$(303, 139, 153, 65)$,
$(303, 292, 207, 162)$,
$(304, 69, 220, 98)$,\adfsplit
$(304, 90, 176, 52)$,
$(304, 65, 214, 33)$,
$(304, 140, 109, 71)$,
$(304, 28, 105, 246)$,\adfsplit
$(304, 42, 217, 89)$,
$(304, 272, 84, 223)$,
$(304, 284, 16, 51)$,
$(304, 241, 287, 108)$,\adfsplit
$(304, 113, 147, 175)$,
$(305, 74, 220, 277)$,
$(305, 51, 211, 196)$,
$(0, 21, 190, 265)$,\adfsplit
$(0, 16, 127, 207)$,
$(0, 28, 260, 305)$,
$(0, 33, 74, 239)$,
$(0, 55, 173, 306)$,\adfsplit
$(0, 65, 224, 362)$,
$(0, 58, 247, 319)$,
$(0, 51, 235, 333)$,
$(0, 89, 153, 389)$,\adfsplit
$(0, 135, 281, 368)$,
$(0, 97, 255, 396)$,
$(0, 134, 297, 347)$,
$(0, 81, 167, 382)$,\adfsplit
$(0, 29, 225, 390)$,
$(0, 63, 233, 376)$,
$(0, 99, 172, 404)$,
$(0, 34, 83, 383)$,\adfsplit
$(0, 137, 187, 334)$,
$(0, 106, 263, 397)$,
$(0, 209, 267, 348)$,
$(0, 41, 50, 320)$

\adfLgap \noindent by the mapping:
$x \mapsto x + 3 j \adfmod{300}$ for $x < 300$,
$x \mapsto (x - 300 + 7 j \adfmod{105}) + 300$ for $300 \le x < 405$,
$x \mapsto x$ for $x \ge 405$,
$0 \le j < 100$
 for the first 11 blocks;
$x \mapsto x + 2 j \adfmod{300}$ for $x < 300$,
$x \mapsto (x - 300 + 7 j \adfmod{105}) + 300$ for $300 \le x < 405$,
$x \mapsto x$ for $x \ge 405$,
$0 \le j < 150$
 for the last 73 blocks.
\ADFvfyParStart{(416, ((11, 100, ((300, 3), (105, 7), (11, 11))), (73, 150, ((300, 2), (105, 7), (11, 11)))), ((50, 6), (116, 1)))} 

\adfDgap
\noindent{\boldmath $ 50^{6} 119^{1} $}~
With the point set $Z_{419}$ partitioned into
 residue classes modulo $6$ for $\{0, 1, \dots, 299\}$, and
 $\{300, 301, \dots, 418\}$,
 the design is generated from

\adfLgap 
$(414, 0, 1, 2)$,
$(415, 0, 100, 98)$,
$(416, 0, 199, 299)$,
$(417, 0, 202, 101)$,\adfsplit
$(418, 0, 298, 200)$,
$(390, 269, 264, 46)$,
$(390, 183, 26, 223)$,
$(391, 87, 52, 104)$,\adfsplit
$(391, 271, 120, 59)$,
$(392, 21, 167, 44)$,
$(392, 168, 202, 43)$,
$(393, 166, 44, 235)$,\adfsplit
$(393, 276, 161, 27)$,
$(300, 267, 252, 293)$,
$(300, 202, 221, 258)$,
$(300, 156, 166, 128)$,\adfsplit
$(300, 15, 59, 2)$,
$(300, 197, 85, 290)$,
$(300, 178, 157, 254)$,
$(300, 160, 129, 13)$,\adfsplit
$(300, 271, 206, 34)$,
$(300, 139, 210, 201)$,
$(300, 123, 84, 215)$,
$(301, 46, 296, 157)$,\adfsplit
$(301, 131, 255, 271)$,
$(301, 245, 253, 82)$,
$(301, 68, 293, 282)$,
$(301, 64, 290, 258)$,\adfsplit
$(301, 126, 209, 130)$,
$(0, 3, 58, 328)$,
$(0, 7, 53, 385)$,
$(0, 14, 267, 382)$,\adfsplit
$(0, 29, 196, 307)$,
$(0, 27, 209, 302)$,
$(0, 73, 158, 347)$,
$(0, 25, 135, 314)$,\adfsplit
$(0, 64, 219, 344)$,
$(0, 20, 87, 326)$,
$(0, 45, 113, 193)$,
$(0, 49, 230, 338)$,\adfsplit
$(0, 43, 173, 389)$,
$(0, 22, 201, 353)$,
$(0, 59, 164, 383)$,
$(0, 89, 206, 362)$

\adfLgap \noindent by the mapping:
$x \mapsto x + 3 j \adfmod{300}$ for $x < 300$,
$x \mapsto (x - 300 + 3 j \adfmod{90}) + 300$ for $300 \le x < 390$,
$x \mapsto (x - 390 + 4 j \adfmod{24}) + 390$ for $390 \le x < 414$,
$x \mapsto x$ for $x \ge 414$,
$0 \le j < 100$
 for the first five blocks;
$x \mapsto x +  j \adfmod{300}$ for $x < 300$,
$x \mapsto (x - 300 + 3 j \adfmod{90}) + 300$ for $300 \le x < 390$,
$x \mapsto (x - 390 + 4 j \adfmod{24}) + 390$ for $390 \le x < 414$,
$x \mapsto x$ for $x \ge 414$,
$0 \le j < 300$
 for the last 39 blocks.
\ADFvfyParStart{(419, ((5, 100, ((300, 3), (90, 3), (24, 4), (5, 5))), (39, 300, ((300, 1), (90, 3), (24, 4), (5, 5)))), ((50, 6), (119, 1)))} 

\adfDgap
\noindent{\boldmath $ 50^{6} 122^{1} $}~
With the point set $Z_{422}$ partitioned into
 residue classes modulo $6$ for $\{0, 1, \dots, 299\}$, and
 $\{300, 301, \dots, 421\}$,
 the design is generated from

\adfLgap 
$(405, 0, 1, 2)$,
$(406, 0, 100, 98)$,
$(407, 0, 199, 299)$,
$(408, 0, 202, 101)$,\adfsplit
$(409, 0, 298, 200)$,
$(410, 0, 4, 11)$,
$(411, 0, 7, 296)$,
$(412, 0, 289, 293)$,\adfsplit
$(413, 0, 10, 23)$,
$(414, 0, 13, 290)$,
$(415, 0, 277, 287)$,
$(416, 0, 16, 35)$,\adfsplit
$(417, 0, 19, 284)$,
$(418, 0, 265, 281)$,
$(419, 0, 22, 5)$,
$(420, 0, 283, 278)$,\adfsplit
$(421, 0, 295, 17)$,
$(300, 117, 34, 295)$,
$(300, 11, 208, 253)$,
$(300, 120, 68, 40)$,\adfsplit
$(300, 21, 61, 140)$,
$(300, 139, 166, 243)$,
$(300, 233, 42, 194)$,
$(300, 86, 129, 232)$,\adfsplit
$(300, 17, 246, 187)$,
$(300, 24, 209, 62)$,
$(300, 108, 65, 15)$,
$(301, 85, 266, 16)$,\adfsplit
$(301, 186, 22, 43)$,
$(301, 230, 131, 3)$,
$(301, 218, 137, 67)$,
$(301, 120, 51, 14)$,\adfsplit
$(301, 147, 95, 162)$,
$(301, 204, 269, 31)$,
$(301, 178, 62, 138)$,
$(301, 9, 94, 229)$,\adfsplit
$(301, 135, 53, 160)$,
$(302, 237, 76, 79)$,
$(302, 201, 41, 242)$,
$(302, 215, 235, 288)$,\adfsplit
$(302, 60, 4, 183)$,
$(302, 107, 66, 163)$,
$(302, 252, 219, 31)$,
$(302, 15, 234, 44)$,\adfsplit
$(302, 67, 142, 230)$,
$(302, 233, 206, 28)$,
$(302, 239, 100, 68)$,
$(303, 191, 36, 103)$,\adfsplit
$(303, 284, 190, 264)$,
$(303, 47, 279, 158)$,
$(303, 93, 211, 282)$,
$(303, 135, 258, 196)$,\adfsplit
$(303, 90, 145, 176)$,
$(303, 34, 247, 80)$,
$(303, 199, 233, 58)$,
$(303, 5, 171, 142)$,\adfsplit
$(303, 269, 237, 212)$,
$(304, 148, 179, 84)$,
$(304, 253, 77, 44)$,
$(304, 288, 112, 75)$,\adfsplit
$(304, 277, 94, 201)$,
$(304, 155, 286, 61)$,
$(304, 186, 278, 220)$,
$(304, 183, 192, 122)$,\adfsplit
$(304, 236, 131, 295)$,
$(304, 147, 109, 293)$,
$(304, 50, 210, 159)$,
$(305, 266, 100, 17)$,\adfsplit
$(0, 9, 158, 333)$,
$(0, 8, 125, 211)$,
$(0, 47, 182, 313)$,
$(0, 15, 172, 312)$,\adfsplit
$(0, 26, 130, 389)$,
$(0, 187, 251, 396)$,
$(0, 93, 245, 347)$,
$(0, 105, 131, 361)$,\adfsplit
$(0, 145, 237, 368)$,
$(1, 29, 75, 396)$,
$(0, 14, 63, 326)$,
$(0, 113, 223, 362)$,\adfsplit
$(0, 89, 133, 376)$,
$(0, 235, 243, 348)$,
$(0, 68, 153, 320)$,
$(0, 82, 255, 341)$,\adfsplit
$(0, 44, 297, 334)$,
$(0, 53, 159, 383)$,
$(0, 91, 112, 404)$,
$(0, 115, 129, 355)$

\adfLgap \noindent by the mapping:
$x \mapsto x + 3 j \adfmod{300}$ for $x < 300$,
$x \mapsto (x - 300 + 7 j \adfmod{105}) + 300$ for $300 \le x < 405$,
$x \mapsto x$ for $x \ge 405$,
$0 \le j < 100$
 for the first 17 blocks;
$x \mapsto x + 2 j \adfmod{300}$ for $x < 300$,
$x \mapsto (x - 300 + 7 j \adfmod{105}) + 300$ for $300 \le x < 405$,
$x \mapsto x$ for $x \ge 405$,
$0 \le j < 150$
 for the last 71 blocks.
\ADFvfyParStart{(422, ((17, 100, ((300, 3), (105, 7), (17, 17))), (71, 150, ((300, 2), (105, 7), (17, 17)))), ((50, 6), (122, 1)))} 

\adfDgap
\noindent{\boldmath $ 50^{9} 11^{1} $}~
With the point set $Z_{461}$ partitioned into
 residue classes modulo $9$ for $\{0, 1, \dots, 449\}$, and
 $\{450, 451, \dots, 460\}$,
 the design is generated from

\adfLgap 
$(459, 0, 1, 2)$,
$(460, 0, 151, 302)$,
$(450, 378, 129, 356)$,
$(450, 224, 408, 293)$,\adfsplit
$(450, 164, 297, 114)$,
$(450, 69, 157, 119)$,
$(450, 52, 19, 382)$,
$(450, 233, 154, 439)$,\adfsplit
$(220, 92, 70, 195)$,
$(375, 102, 22, 29)$,
$(220, 323, 252, 385)$,
$(271, 432, 13, 70)$,\adfsplit
$(21, 173, 176, 448)$,
$(263, 12, 446, 125)$,
$(261, 427, 28, 165)$,
$(408, 388, 421, 119)$,\adfsplit
$(270, 213, 436, 281)$,
$(110, 251, 166, 381)$,
$(47, 302, 379, 152)$,
$(5, 98, 46, 159)$,\adfsplit
$(445, 109, 302, 314)$,
$(283, 95, 298, 92)$,
$(410, 373, 250, 309)$,
$(107, 2, 157, 73)$,\adfsplit
$(52, 261, 191, 251)$,
$(7, 191, 21, 235)$,
$(147, 355, 105, 395)$,
$(235, 196, 212, 252)$,\adfsplit
$(240, 410, 436, 352)$,
$(370, 255, 365, 423)$,
$(294, 22, 153, 350)$,
$(216, 52, 4, 422)$,\adfsplit
$(216, 237, 373, 8)$,
$(418, 92, 433, 43)$,
$(428, 386, 148, 150)$,
$(67, 161, 168, 72)$,\adfsplit
$(130, 185, 201, 312)$,
$(191, 205, 279, 361)$,
$(236, 247, 363, 144)$,
$(238, 151, 226, 41)$,\adfsplit
$(8, 367, 157, 15)$,
$(164, 342, 346, 69)$,
$(215, 394, 283, 163)$,
$(27, 404, 264, 325)$,\adfsplit
$(39, 184, 343, 278)$,
$(73, 327, 303, 263)$,
$(0, 5, 102, 124)$,
$(0, 6, 316, 347)$,\adfsplit
$(0, 16, 137, 143)$,
$(0, 7, 65, 154)$,
$(0, 23, 179, 281)$,
$(0, 39, 263, 346)$,\adfsplit
$(0, 53, 159, 237)$,
$(0, 95, 172, 246)$,
$(3, 35, 69, 243)$,
$(0, 100, 364, 381)$,\adfsplit
$(0, 28, 160, 276)$,
$(0, 47, 208, 399)$,
$(0, 93, 127, 429)$,
$(0, 255, 269, 304)$,\adfsplit
$(0, 220, 224, 442)$,
$(0, 44, 167, 424)$,
$(0, 98, 292, 322)$,
$(0, 46, 168, 250)$,\adfsplit
$(0, 86, 265, 407)$,
$(0, 107, 132, 254)$,
$(0, 32, 70, 110)$,
$(0, 73, 186, 314)$,\adfsplit
$(0, 30, 97, 145)$,
$(0, 55, 121, 290)$,
$(0, 24, 62, 200)$

\adfLgap \noindent by the mapping:
$x \mapsto x \oplus (3 j)$ for $x < 450$,
$x \mapsto (x +  j \adfmod{9}) + 450$ for $450 \le x < 459$,
$x \mapsto x$ for $x \ge 459$,
$0 \le j < 150$
 for the first two blocks;
$x \mapsto x \oplus j \oplus j$ for $x < 450$,
$x \mapsto (x +  j \adfmod{9}) + 450$ for $450 \le x < 459$,
$x \mapsto x$ for $x \ge 459$,
$0 \le j < 225$
 for the last 69 blocks.
\ADFvfyParStart{(461, ((2, 150, ((450, 3, (150, 3)), (9, 1), (2, 2))), (69, 225, ((450, 2, (150, 3)), (9, 1), (2, 2)))), ((50, 9), (11, 1)))} 

\adfDgap
\noindent{\boldmath $ 50^{9} 17^{1} $}~
With the point set $Z_{467}$ partitioned into
 residue classes modulo $9$ for $\{0, 1, \dots, 449\}$, and
 $\{450, 451, \dots, 466\}$,
 the design is generated from

\adfLgap 
$(465, 0, 1, 2)$,
$(466, 0, 151, 302)$,
$(450, 199, 249, 308)$,
$(450, 421, 404, 148)$,\adfsplit
$(450, 291, 26, 365)$,
$(450, 262, 247, 117)$,
$(450, 403, 419, 84)$,
$(450, 371, 192, 203)$,\adfsplit
$(450, 46, 15, 440)$,
$(450, 310, 33, 126)$,
$(450, 107, 150, 318)$,
$(450, 122, 145, 304)$,\adfsplit
$(359, 255, 13, 105)$,
$(20, 446, 118, 439)$,
$(167, 116, 75, 316)$,
$(210, 412, 359, 256)$,\adfsplit
$(315, 26, 57, 312)$,
$(279, 410, 375, 138)$,
$(51, 171, 341, 299)$,
$(442, 256, 26, 303)$,\adfsplit
$(139, 86, 171, 395)$,
$(218, 251, 343, 338)$,
$(198, 335, 410, 352)$,
$(132, 157, 331, 84)$,\adfsplit
$(81, 420, 12, 145)$,
$(310, 375, 44, 349)$,
$(226, 131, 65, 368)$,
$(328, 349, 192, 406)$,\adfsplit
$(427, 17, 174, 429)$,
$(20, 386, 141, 54)$,
$(122, 241, 227, 53)$,
$(293, 263, 256, 330)$,\adfsplit
$(184, 415, 313, 108)$,
$(264, 157, 353, 5)$,
$(258, 103, 0, 431)$,
$(182, 35, 43, 369)$,\adfsplit
$(196, 229, 361, 128)$,
$(305, 320, 414, 364)$,
$(8, 295, 105, 172)$,
$(200, 260, 337, 384)$,\adfsplit
$(91, 320, 126, 292)$,
$(269, 423, 196, 192)$,
$(151, 428, 206, 246)$,
$(333, 94, 341, 84)$,\adfsplit
$(424, 372, 60, 437)$,
$(265, 263, 303, 354)$,
$(318, 253, 155, 331)$,
$(0, 6, 23, 91)$,\adfsplit
$(0, 14, 75, 434)$,
$(0, 21, 204, 392)$,
$(0, 12, 44, 379)$,
$(0, 39, 211, 260)$,\adfsplit
$(0, 58, 106, 284)$,
$(0, 66, 301, 371)$,
$(0, 62, 155, 310)$,
$(0, 130, 150, 449)$,\adfsplit
$(0, 113, 223, 337)$,
$(0, 107, 183, 240)$,
$(0, 79, 221, 244)$,
$(0, 67, 353, 437)$,\adfsplit
$(0, 22, 119, 354)$,
$(0, 123, 170, 443)$,
$(0, 161, 238, 264)$,
$(0, 83, 220, 401)$,\adfsplit
$(0, 172, 178, 365)$,
$(0, 148, 291, 294)$,
$(0, 131, 143, 250)$,
$(3, 27, 112, 365)$,\adfsplit
$(0, 160, 209, 388)$,
$(0, 112, 316, 345)$,
$(3, 17, 58, 141)$,
$(3, 22, 191, 213)$,\adfsplit
$(3, 63, 124, 238)$

\adfLgap \noindent by the mapping:
$x \mapsto x \oplus (3 j)$ for $x < 450$,
$x \mapsto (x +  j \adfmod{15}) + 450$ for $450 \le x < 465$,
$x \mapsto x$ for $x \ge 465$,
$0 \le j < 150$
 for the first two blocks;
$x \mapsto x \oplus j \oplus j$ for $x < 450$,
$x \mapsto (x +  j \adfmod{15}) + 450$ for $450 \le x < 465$,
$x \mapsto x$ for $x \ge 465$,
$0 \le j < 225$
 for the last 71 blocks.
\ADFvfyParStart{(467, ((2, 150, ((450, 3, (150, 3)), (15, 1), (2, 2))), (71, 225, ((450, 2, (150, 3)), (15, 1), (2, 2)))), ((50, 9), (17, 1)))} 

\adfDgap
\noindent{\boldmath $ 50^{9} 23^{1} $}~
With the point set $Z_{473}$ partitioned into
 residue classes modulo $9$ for $\{0, 1, \dots, 449\}$, and
 $\{450, 451, \dots, 472\}$,
 the design is generated from

\adfLgap 
$(471, 0, 1, 2)$,
$(472, 0, 151, 302)$,
$(450, 208, 404, 299)$,
$(450, 33, 390, 157)$,\adfsplit
$(451, 392, 168, 161)$,
$(451, 442, 51, 43)$,
$(452, 299, 284, 33)$,
$(452, 256, 198, 181)$,\adfsplit
$(453, 303, 311, 378)$,
$(453, 43, 110, 130)$,
$(454, 170, 66, 196)$,
$(454, 411, 271, 137)$,\adfsplit
$(455, 169, 292, 108)$,
$(455, 165, 365, 254)$,
$(456, 0, 303, 121)$,
$(456, 365, 314, 196)$,\adfsplit
$(215, 229, 131, 425)$,
$(72, 61, 236, 275)$,
$(306, 287, 83, 313)$,
$(270, 164, 96, 349)$,\adfsplit
$(341, 435, 358, 128)$,
$(165, 8, 293, 281)$,
$(221, 402, 256, 359)$,
$(71, 432, 240, 273)$,\adfsplit
$(387, 169, 245, 172)$,
$(270, 14, 197, 254)$,
$(140, 225, 145, 236)$,
$(72, 201, 177, 439)$,\adfsplit
$(183, 396, 112, 430)$,
$(358, 318, 425, 284)$,
$(61, 366, 200, 117)$,
$(363, 23, 1, 245)$,\adfsplit
$(250, 86, 24, 351)$,
$(31, 28, 155, 369)$,
$(94, 354, 370, 402)$,
$(435, 17, 100, 54)$,\adfsplit
$(233, 325, 409, 243)$,
$(288, 137, 69, 388)$,
$(183, 394, 23, 208)$,
$(9, 217, 309, 62)$,\adfsplit
$(31, 66, 41, 169)$,
$(67, 294, 199, 203)$,
$(93, 245, 204, 323)$,
$(348, 93, 254, 229)$,\adfsplit
$(220, 100, 422, 102)$,
$(252, 191, 125, 136)$,
$(263, 312, 244, 162)$,
$(15, 185, 227, 63)$,\adfsplit
$(48, 395, 270, 133)$,
$(0, 5, 35, 438)$,
$(0, 21, 28, 71)$,
$(0, 6, 112, 161)$,\adfsplit
$(0, 13, 363, 419)$,
$(0, 42, 280, 437)$,
$(0, 49, 113, 305)$,
$(3, 64, 149, 339)$,\adfsplit
$(0, 30, 208, 376)$,
$(0, 172, 267, 311)$,
$(0, 52, 148, 154)$,
$(0, 94, 246, 411)$,\adfsplit
$(0, 249, 275, 309)$,
$(0, 137, 193, 429)$,
$(0, 159, 163, 340)$,
$(0, 50, 127, 197)$,\adfsplit
$(0, 69, 293, 430)$,
$(0, 65, 321, 336)$,
$(0, 57, 328, 361)$,
$(0, 26, 255, 286)$,\adfsplit
$(0, 24, 206, 412)$,
$(0, 53, 146, 330)$,
$(0, 44, 110, 283)$,
$(0, 56, 134, 205)$,\adfsplit
$(0, 55, 115, 296)$,
$(0, 31, 156, 199)$,
$(0, 32, 102, 218)$

\adfLgap \noindent by the mapping:
$x \mapsto x \oplus (3 j)$ for $x < 450$,
$x \mapsto (x - 450 + 7 j \adfmod{21}) + 450$ for $450 \le x < 471$,
$x \mapsto x$ for $x \ge 471$,
$0 \le j < 150$
 for the first two blocks;
$x \mapsto x \oplus j \oplus j$ for $x < 450$,
$x \mapsto (x - 450 + 7 j \adfmod{21}) + 450$ for $450 \le x < 471$,
$x \mapsto x$ for $x \ge 471$,
$0 \le j < 225$
 for the last 73 blocks.
\ADFvfyParStart{(473, ((2, 150, ((450, 3, (150, 3)), (21, 7), (2, 2))), (73, 225, ((450, 2, (150, 3)), (21, 7), (2, 2)))), ((50, 9), (23, 1)))} 

\section{4-GDDs for the proof of Lemma \ref{lem:4-GDD 52^u m^1}}
\label{app:4-GDD 52^u m^1}
\adfhide{
$ 52^{15} 34^1 $.
}

\adfDgap
\noindent{\boldmath $ 52^{15} 34^{1} $}~
With the point set $Z_{814}$ partitioned into
 residue classes modulo $15$ for $\{0, 1, \dots, 779\}$, and
 $\{780, 781, \dots, 813\}$,
 the design is generated from

\adfLgap 
$(780, 459, 337, 668)$,
$(781, 226, 656, 132)$,
$(782, 644, 202, 300)$,
$(783, 635, 361, 363)$,\adfsplit
$(784, 21, 449, 379)$,
$(785, 24, 238, 179)$,
$(786, 391, 455, 360)$,
$(787, 125, 168, 766)$,\adfsplit
$(788, 196, 170, 588)$,
$(789, 291, 742, 263)$,
$(790, 761, 774, 550)$,
$(353, 301, 245, 321)$,\adfsplit
$(626, 549, 469, 423)$,
$(644, 656, 490, 376)$,
$(107, 90, 517, 665)$,
$(606, 695, 317, 568)$,\adfsplit
$(308, 467, 409, 749)$,
$(725, 760, 444, 361)$,
$(220, 282, 479, 695)$,
$(732, 479, 478, 300)$,\adfsplit
$(564, 764, 682, 331)$,
$(79, 367, 358, 280)$,
$(196, 357, 495, 129)$,
$(139, 158, 292, 582)$,\adfsplit
$(375, 618, 467, 217)$,
$(24, 127, 398, 447)$,
$(555, 518, 439, 551)$,
$(86, 384, 554, 50)$,\adfsplit
$(12, 183, 486, 34)$,
$(44, 50, 533, 731)$,
$(208, 258, 634, 83)$,
$(235, 42, 559, 666)$,\adfsplit
$(200, 418, 580, 607)$,
$(705, 267, 448, 380)$,
$(626, 436, 53, 439)$,
$(312, 531, 708, 422)$,\adfsplit
$(78, 325, 178, 15)$,
$(685, 350, 619, 296)$,
$(423, 681, 394, 227)$,
$(66, 11, 278, 343)$,\adfsplit
$(748, 737, 199, 730)$,
$(153, 322, 426, 234)$,
$(196, 523, 402, 230)$,
$(754, 383, 105, 697)$,\adfsplit
$(80, 505, 577, 696)$,
$(0, 5, 96, 612)$,
$(0, 8, 140, 184)$,
$(0, 10, 246, 387)$,\adfsplit
$(0, 16, 308, 359)$,
$(0, 14, 53, 709)$,
$(0, 23, 317, 462)$,
$(0, 128, 265, 572)$,\adfsplit
$(0, 25, 73, 484)$,
$(0, 61, 295, 532)$,
$(0, 69, 186, 647)$,
$(0, 42, 194, 241)$,\adfsplit
$(0, 24, 245, 458)$,
$(0, 41, 215, 519)$,
$(0, 21, 204, 244)$,
$(0, 115, 238, 620)$,\adfsplit
$(0, 102, 232, 467)$,
$(0, 86, 230, 447)$,
$(0, 74, 185, 469)$,
$(0, 84, 220, 637)$,\adfsplit
$(0, 97, 226, 638)$,
$(0, 106, 252, 514)$,
$(813, 0, 260, 520)$

\adfLgap \noindent by the mapping:
$x \mapsto x +  j \adfmod{780}$ for $x < 780$,
$x \mapsto (x - 780 + 11 j \adfmod{33}) + 780$ for $780 \le x < 813$,
$813 \mapsto 813$,
$0 \le j < 780$
 for the first 66 blocks,
$0 \le j < 260$
 for the last block.
\ADFvfyParStart{(814, ((66, 780, ((780, 1), (33, 11), (1, 1))), (1, 260, ((780, 1), (33, 11), (1, 1)))), ((52, 15), (34, 1)))} 

\section{4-GDDs for the proof of Lemma \ref{lem:4-GDD 58^u m^1}}
\label{app:4-GDD 58^u m^1}
\adfhide{
$ 58^6 133^1 $,
$ 58^6 136^1 $,
$ 58^6 139^1 $,
$ 58^6 142^1 $ and
$ 58^9 19^1 $.
}

\adfDgap
\noindent{\boldmath $ 58^{6} 133^{1} $}~
With the point set $Z_{481}$ partitioned into
 residue classes modulo $6$ for $\{0, 1, \dots, 347\}$, and
 $\{348, 349, \dots, 480\}$,
 the design is generated from

\adfLgap 
$(348, 256, 301, 251)$,
$(348, 305, 296, 195)$,
$(348, 314, 94, 18)$,
$(348, 105, 120, 223)$,\adfsplit
$(349, 200, 43, 297)$,
$(349, 220, 53, 290)$,
$(349, 51, 346, 30)$,
$(349, 36, 143, 1)$,\adfsplit
$(350, 62, 238, 179)$,
$(350, 329, 291, 55)$,
$(350, 104, 337, 288)$,
$(350, 321, 30, 208)$,\adfsplit
$(351, 81, 347, 204)$,
$(351, 6, 40, 313)$,
$(351, 89, 176, 111)$,
$(351, 142, 91, 2)$,\adfsplit
$(352, 76, 320, 120)$,
$(352, 299, 295, 315)$,
$(352, 34, 273, 246)$,
$(352, 302, 217, 341)$,\adfsplit
$(353, 7, 282, 71)$,
$(353, 148, 161, 147)$,
$(353, 106, 261, 308)$,
$(353, 2, 168, 313)$,\adfsplit
$(354, 294, 248, 117)$,
$(354, 175, 77, 274)$,
$(354, 100, 263, 183)$,
$(354, 182, 48, 277)$,\adfsplit
$(355, 21, 314, 126)$,
$(355, 163, 188, 5)$,
$(355, 274, 251, 337)$,
$(355, 39, 316, 120)$,\adfsplit
$(356, 210, 181, 104)$,
$(356, 293, 154, 295)$,
$(0, 3, 11, 356)$,
$(0, 7, 17, 389)$,\adfsplit
$(0, 19, 62, 88)$,
$(0, 91, 218, 379)$,
$(0, 31, 226, 358)$,
$(0, 58, 125, 369)$,\adfsplit
$(0, 92, 213, 368)$,
$(0, 61, 280, 467)$,
$(0, 79, 248, 424)$,
$(0, 33, 187, 402)$,\adfsplit
$(0, 40, 133, 189)$,
$(0, 28, 175, 390)$,
$(480, 0, 116, 232)$

\adfLgap \noindent by the mapping:
$x \mapsto x +  j \adfmod{348}$ for $x < 348$,
$x \mapsto (x - 348 + 11 j \adfmod{132}) + 348$ for $348 \le x < 480$,
$480 \mapsto 480$,
$0 \le j < 348$
 for the first 46 blocks,
$0 \le j < 116$
 for the last block.
\ADFvfyParStart{(481, ((46, 348, ((348, 1), (132, 11), (1, 1))), (1, 116, ((348, 1), (132, 11), (1, 1)))), ((58, 6), (133, 1)))} 

\adfDgap
\noindent{\boldmath $ 58^{6} 136^{1} $}~
With the point set $Z_{484}$ partitioned into
 residue classes modulo $6$ for $\{0, 1, \dots, 347\}$, and
 $\{348, 349, \dots, 483\}$,
 the design is generated from

\adfLgap 
$(480, 331, 207, 86)$,
$(480, 324, 113, 334)$,
$(348, 140, 42, 220)$,
$(348, 151, 321, 48)$,\adfsplit
$(348, 170, 347, 337)$,
$(348, 58, 243, 221)$,
$(349, 78, 172, 133)$,
$(349, 11, 297, 56)$,\adfsplit
$(349, 192, 87, 118)$,
$(349, 283, 218, 101)$,
$(350, 58, 338, 187)$,
$(350, 77, 189, 138)$,\adfsplit
$(350, 107, 253, 228)$,
$(350, 116, 172, 159)$,
$(351, 150, 191, 110)$,
$(351, 184, 291, 103)$,\adfsplit
$(351, 168, 190, 121)$,
$(351, 104, 101, 297)$,
$(352, 286, 236, 281)$,
$(352, 192, 307, 184)$,\adfsplit
$(352, 114, 213, 146)$,
$(352, 35, 109, 207)$,
$(353, 329, 256, 92)$,
$(353, 103, 216, 21)$,\adfsplit
$(353, 47, 250, 230)$,
$(353, 147, 61, 198)$,
$(354, 155, 87, 4)$,
$(354, 258, 265, 297)$,\adfsplit
$(354, 92, 118, 156)$,
$(354, 182, 17, 127)$,
$(355, 71, 48, 326)$,
$(355, 97, 114, 340)$,\adfsplit
$(355, 166, 201, 43)$,
$(355, 293, 291, 320)$,
$(356, 313, 112, 284)$,
$(356, 201, 163, 186)$,\adfsplit
$(356, 298, 207, 119)$,
$(356, 192, 245, 50)$,
$(357, 306, 145, 287)$,
$(357, 48, 260, 317)$,\adfsplit
$(357, 136, 279, 74)$,
$(357, 130, 31, 285)$,
$(358, 236, 40, 288)$,
$(358, 117, 142, 13)$,\adfsplit
$(358, 299, 135, 283)$,
$(358, 281, 278, 186)$,
$(359, 167, 225, 220)$,
$(359, 318, 217, 58)$,\adfsplit
$(359, 187, 51, 14)$,
$(359, 60, 44, 53)$,
$(360, 284, 269, 333)$,
$(360, 205, 87, 226)$,\adfsplit
$(360, 318, 184, 170)$,
$(360, 95, 235, 312)$,
$(361, 336, 15, 145)$,
$(361, 310, 53, 79)$,\adfsplit
$(361, 266, 270, 203)$,
$(361, 200, 165, 40)$,
$(362, 225, 4, 116)$,
$(362, 193, 87, 74)$,\adfsplit
$(362, 94, 191, 336)$,
$(362, 149, 115, 282)$,
$(363, 329, 309, 259)$,
$(363, 128, 256, 207)$,\adfsplit
$(363, 106, 330, 95)$,
$(363, 38, 240, 301)$,
$(0, 1, 238, 272)$,
$(0, 17, 21, 289)$,\adfsplit
$(0, 19, 213, 364)$,
$(0, 2, 33, 408)$,
$(0, 11, 290, 474)$,
$(0, 89, 141, 430)$,\adfsplit
$(0, 46, 305, 365)$,
$(0, 28, 189, 453)$,
$(0, 44, 261, 431)$,
$(1, 93, 215, 387)$,\adfsplit
$(0, 179, 223, 366)$,
$(0, 47, 262, 454)$,
$(0, 75, 190, 388)$,
$(0, 239, 339, 432)$,\adfsplit
$(0, 139, 140, 367)$,
$(0, 135, 149, 433)$,
$(0, 59, 130, 455)$,
$(0, 147, 175, 389)$,\adfsplit
$(0, 85, 118, 368)$,
$(0, 251, 291, 369)$,
$(0, 71, 154, 391)$,
$(0, 255, 311, 479)$,\adfsplit
$(0, 82, 253, 435)$,
$(0, 229, 275, 434)$,
$(0, 101, 166, 456)$,
$(0, 199, 207, 412)$,\adfsplit
$(0, 63, 104, 191)$,
$(483, 0, 116, 232)$,
$(483, 1, 117, 233)$

\adfLgap \noindent by the mapping:
$x \mapsto x + 2 j \adfmod{348}$ for $x < 348$,
$x \mapsto (x - 348 + 22 j \adfmod{132}) + 348$ for $348 \le x < 480$,
$x \mapsto (x +  j \adfmod{3}) + 480$ for $480 \le x < 483$,
$483 \mapsto 483$,
$0 \le j < 174$
 for the first 93 blocks,
$0 \le j < 58$
 for the last two blocks.
\ADFvfyParStart{(484, ((93, 174, ((348, 2), (132, 22), (3, 1), (1, 1))), (2, 58, ((348, 2), (132, 22), (3, 1), (1, 1)))), ((58, 6), (136, 1)))} 

\adfDgap
\noindent{\boldmath $ 58^{6} 139^{1} $}~
With the point set $Z_{487}$ partitioned into
 residue classes modulo $6$ for $\{0, 1, \dots, 347\}$, and
 $\{348, 349, \dots, 486\}$,
 the design is generated from

\adfLgap 
$(480, 87, 260, 58)$,
$(480, 59, 229, 90)$,
$(348, 203, 12, 220)$,
$(348, 314, 113, 177)$,\adfsplit
$(348, 1, 80, 298)$,
$(348, 171, 18, 199)$,
$(349, 1, 16, 315)$,
$(349, 19, 257, 8)$,\adfsplit
$(349, 114, 314, 189)$,
$(349, 312, 227, 322)$,
$(350, 285, 199, 239)$,
$(350, 92, 318, 147)$,\adfsplit
$(350, 169, 82, 5)$,
$(350, 52, 50, 300)$,
$(351, 244, 155, 108)$,
$(351, 277, 56, 17)$,\adfsplit
$(351, 90, 27, 199)$,
$(351, 146, 213, 334)$,
$(352, 284, 193, 150)$,
$(352, 319, 274, 191)$,\adfsplit
$(352, 324, 201, 89)$,
$(352, 134, 15, 292)$,
$(353, 260, 157, 136)$,
$(353, 123, 180, 254)$,\adfsplit
$(353, 257, 234, 91)$,
$(353, 71, 141, 34)$,
$(354, 103, 144, 112)$,
$(354, 275, 261, 126)$,\adfsplit
$(354, 346, 3, 290)$,
$(354, 41, 68, 49)$,
$(355, 135, 102, 170)$,
$(355, 236, 121, 77)$,\adfsplit
$(355, 167, 160, 271)$,
$(355, 310, 12, 177)$,
$(356, 63, 218, 101)$,
$(0, 1, 4, 26)$,\adfsplit
$(0, 58, 152, 356)$,
$(0, 13, 219, 422)$,
$(0, 69, 187, 444)$,
$(0, 73, 179, 368)$,\adfsplit
$(0, 16, 267, 402)$,
$(0, 20, 82, 390)$,
$(0, 52, 145, 467)$,
$(0, 65, 272, 401)$,\adfsplit
$(0, 53, 154, 435)$,
$(0, 59, 151, 413)$,
$(0, 80, 243, 479)$,
$(486, 0, 116, 232)$

\adfLgap \noindent by the mapping:
$x \mapsto x +  j \adfmod{348}$ for $x < 348$,
$x \mapsto (x - 348 + 11 j \adfmod{132}) + 348$ for $348 \le x < 480$,
$x \mapsto (x +  j \adfmod{6}) + 480$ for $480 \le x < 486$,
$486 \mapsto 486$,
$0 \le j < 348$
 for the first 47 blocks,
$0 \le j < 116$
 for the last block.
\ADFvfyParStart{(487, ((47, 348, ((348, 1), (132, 11), (6, 1), (1, 1))), (1, 116, ((348, 1), (132, 11), (6, 1), (1, 1)))), ((58, 6), (139, 1)))} 

\adfDgap
\noindent{\boldmath $ 58^{6} 142^{1} $}~
With the point set $Z_{490}$ partitioned into
 residue classes modulo $6$ for $\{0, 1, \dots, 347\}$, and
 $\{348, 349, \dots, 489\}$,
 the design is generated from

\adfLgap 
$(486, 82, 180, 99)$,
$(486, 311, 14, 73)$,
$(348, 16, 99, 104)$,
$(348, 59, 168, 127)$,\adfsplit
$(348, 34, 158, 318)$,
$(348, 333, 341, 97)$,
$(349, 260, 24, 25)$,
$(349, 90, 266, 211)$,\adfsplit
$(349, 75, 286, 323)$,
$(349, 137, 9, 160)$,
$(350, 62, 211, 117)$,
$(350, 185, 340, 212)$,\adfsplit
$(350, 214, 72, 11)$,
$(350, 325, 246, 3)$,
$(351, 243, 172, 157)$,
$(351, 78, 273, 155)$,\adfsplit
$(351, 8, 118, 312)$,
$(351, 247, 329, 182)$,
$(352, 5, 36, 284)$,
$(352, 299, 206, 247)$,\adfsplit
$(352, 42, 268, 219)$,
$(352, 129, 310, 313)$,
$(353, 300, 244, 19)$,
$(353, 188, 179, 87)$,\adfsplit
$(353, 297, 5, 210)$,
$(353, 242, 34, 85)$,
$(354, 202, 69, 215)$,
$(354, 206, 65, 132)$,\adfsplit
$(354, 207, 100, 13)$,
$(354, 200, 198, 163)$,
$(355, 281, 112, 60)$,
$(355, 287, 9, 152)$,\adfsplit
$(355, 166, 73, 291)$,
$(355, 126, 31, 146)$,
$(356, 317, 290, 318)$,
$(356, 20, 123, 151)$,\adfsplit
$(356, 239, 345, 10)$,
$(356, 169, 60, 328)$,
$(357, 86, 167, 82)$,
$(357, 73, 42, 172)$,\adfsplit
$(357, 84, 269, 9)$,
$(357, 211, 284, 75)$,
$(358, 309, 127, 84)$,
$(358, 88, 26, 174)$,\adfsplit
$(358, 207, 250, 191)$,
$(358, 284, 149, 193)$,
$(359, 52, 203, 282)$,
$(359, 149, 130, 145)$,\adfsplit
$(359, 195, 242, 96)$,
$(359, 116, 105, 139)$,
$(360, 137, 325, 304)$,
$(360, 141, 186, 254)$,\adfsplit
$(360, 204, 279, 226)$,
$(360, 163, 284, 227)$,
$(361, 282, 86, 40)$,
$(361, 166, 195, 216)$,\adfsplit
$(361, 305, 295, 260)$,
$(361, 95, 157, 237)$,
$(362, 199, 36, 3)$,
$(362, 138, 309, 85)$,\adfsplit
$(362, 326, 190, 287)$,
$(362, 248, 101, 232)$,
$(363, 76, 0, 266)$,
$(363, 202, 13, 234)$,\adfsplit
$(363, 5, 195, 19)$,
$(363, 321, 143, 248)$,
$(364, 293, 16, 21)$,
$(364, 319, 326, 111)$,\adfsplit
$(364, 138, 272, 311)$,
$(364, 97, 262, 288)$,
$(365, 185, 178, 282)$,
$(0, 8, 33, 365)$,\adfsplit
$(0, 9, 209, 480)$,
$(0, 11, 57, 457)$,
$(0, 40, 129, 263)$,
$(0, 115, 237, 366)$,\adfsplit
$(0, 10, 165, 458)$,
$(0, 14, 63, 412)$,
$(0, 47, 345, 435)$,
$(0, 91, 111, 367)$,\adfsplit
$(0, 38, 241, 391)$,
$(0, 70, 164, 369)$,
$(0, 161, 259, 370)$,
$(0, 329, 331, 438)$,\adfsplit
$(0, 179, 219, 461)$,
$(0, 199, 231, 484)$,
$(0, 101, 178, 393)$,
$(0, 245, 283, 483)$,\adfsplit
$(0, 92, 279, 439)$,
$(0, 117, 175, 462)$,
$(0, 119, 141, 482)$,
$(0, 58, 323, 390)$,\adfsplit
$(0, 34, 95, 436)$,
$(0, 153, 182, 460)$,
$(0, 211, 285, 414)$,
$(489, 0, 116, 232)$,\adfsplit
$(489, 1, 117, 233)$

\adfLgap \noindent by the mapping:
$x \mapsto x + 2 j \adfmod{348}$ for $x < 348$,
$x \mapsto (x - 348 + 23 j \adfmod{138}) + 348$ for $348 \le x < 486$,
$x \mapsto (x +  j \adfmod{3}) + 486$ for $486 \le x < 489$,
$489 \mapsto 489$,
$0 \le j < 174$
 for the first 95 blocks,
$0 \le j < 58$
 for the last two blocks.
\ADFvfyParStart{(490, ((95, 174, ((348, 2), (138, 23), (3, 1), (1, 1))), (2, 58, ((348, 2), (138, 23), (3, 1), (1, 1)))), ((58, 6), (142, 1)))} 

\adfDgap
\noindent{\boldmath $ 58^{9} 19^{1} $}~
With the point set $Z_{541}$ partitioned into
 residue classes modulo $9$ for $\{0, 1, \dots, 521\}$, and
 $\{522, 523, \dots, 540\}$,
 the design is generated from

\adfLgap 
$(522, 32, 382, 186)$,
$(522, 120, 141, 333)$,
$(522, 239, 219, 163)$,
$(522, 236, 17, 486)$,\adfsplit
$(522, 151, 65, 337)$,
$(522, 80, 124, 100)$,
$(523, 100, 326, 288)$,
$(523, 152, 219, 59)$,\adfsplit
$(523, 177, 493, 193)$,
$(523, 372, 243, 16)$,
$(523, 176, 181, 35)$,
$(523, 389, 168, 436)$,\adfsplit
$(396, 281, 41, 444)$,
$(417, 316, 123, 31)$,
$(132, 9, 344, 211)$,
$(187, 164, 77, 382)$,\adfsplit
$(9, 84, 15, 22)$,
$(392, 47, 366, 296)$,
$(30, 182, 88, 105)$,
$(301, 205, 14, 173)$,\adfsplit
$(272, 324, 332, 55)$,
$(62, 201, 292, 63)$,
$(307, 420, 243, 53)$,
$(303, 70, 127, 504)$,\adfsplit
$(133, 450, 91, 384)$,
$(508, 37, 179, 88)$,
$(399, 306, 133, 416)$,
$(393, 227, 511, 262)$,\adfsplit
$(362, 358, 291, 180)$,
$(133, 498, 236, 514)$,
$(88, 168, 369, 143)$,
$(195, 349, 247, 345)$,\adfsplit
$(322, 266, 152, 166)$,
$(495, 454, 132, 503)$,
$(178, 175, 212, 51)$,
$(414, 208, 427, 167)$,\adfsplit
$(274, 366, 176, 399)$,
$(36, 300, 451, 13)$,
$(452, 145, 266, 45)$,
$(363, 236, 431, 365)$,\adfsplit
$(35, 272, 69, 183)$,
$(478, 219, 63, 344)$,
$(468, 390, 33, 232)$,
$(472, 336, 56, 216)$,\adfsplit
$(186, 434, 21, 499)$,
$(298, 256, 213, 203)$,
$(117, 196, 128, 442)$,
$(239, 422, 447, 407)$,\adfsplit
$(328, 8, 425, 97)$,
$(119, 181, 230, 399)$,
$(3, 266, 352, 81)$,
$(24, 22, 95, 485)$,\adfsplit
$(432, 302, 516, 283)$,
$(441, 479, 304, 480)$,
$(52, 23, 336, 346)$,
$(0, 3, 104, 494)$,\adfsplit
$(0, 6, 124, 146)$,
$(0, 11, 32, 472)$,
$(0, 37, 168, 290)$,
$(0, 7, 12, 150)$,\adfsplit
$(0, 15, 29, 220)$,
$(0, 39, 184, 439)$,
$(0, 69, 355, 410)$,
$(0, 47, 59, 448)$,\adfsplit
$(0, 30, 194, 240)$,
$(0, 147, 343, 374)$,
$(0, 64, 192, 449)$,
$(0, 19, 231, 311)$,\adfsplit
$(0, 103, 185, 367)$,
$(0, 40, 307, 446)$,
$(0, 43, 300, 413)$,
$(0, 277, 335, 365)$,\adfsplit
$(0, 161, 325, 375)$,
$(0, 83, 465, 489)$,
$(0, 183, 211, 341)$,
$(0, 123, 149, 425)$,\adfsplit
$(0, 77, 142, 247)$,
$(0, 85, 179, 298)$,
$(1, 47, 121, 225)$,
$(1, 61, 233, 339)$,\adfsplit
$(0, 95, 143, 318)$,
$(0, 35, 283, 353)$,
$(1, 23, 135, 335)$,
$(540, 0, 174, 348)$,\adfsplit
$(540, 1, 175, 349)$

\adfLgap \noindent by the mapping:
$x \mapsto x + 2 j \adfmod{522}$ for $x < 522$,
$x \mapsto (x + 2 j \adfmod{18}) + 522$ for $522 \le x < 540$,
$540 \mapsto 540$,
$0 \le j < 261$
 for the first 83 blocks,
$0 \le j < 87$
 for the last two blocks.
\ADFvfyParStart{(541, ((83, 261, ((522, 2), (18, 2), (1, 1))), (2, 87, ((522, 2), (18, 2), (1, 1)))), ((58, 9), (19, 1)))} 

\section{4-GDDs for the proof of Lemma \ref{lem:4-GDD 62^u m^1}}
\label{app:4-GDD 62^u m^1}
\adfhide{
$ 62^6 23^1 $,
$ 62^6 26^1 $,
$ 62^6 29^1 $,
$ 62^6 32^1 $,
$ 62^6 35^1 $,
$ 62^6 38^1 $,
$ 62^6 41^1 $,
$ 62^6 44^1 $,
$ 62^6 47^1 $,
$ 62^6 143^1 $,
$ 62^6 146^1 $,
$ 62^6 149^1 $ and
$ 62^6 152^1 $.
}

\adfDgap
\noindent{\boldmath $ 62^{6} 23^{1} $}~
With the point set $Z_{395}$ partitioned into
 residue classes modulo $6$ for $\{0, 1, \dots, 371\}$, and
 $\{372, 373, \dots, 394\}$,
 the design is generated from

\adfLgap 
$(384, 0, 1, 2)$,
$(385, 0, 124, 122)$,
$(386, 0, 247, 371)$,
$(387, 0, 250, 125)$,\adfsplit
$(388, 0, 370, 248)$,
$(389, 0, 4, 11)$,
$(390, 0, 7, 368)$,
$(391, 0, 361, 365)$,\adfsplit
$(392, 0, 10, 23)$,
$(393, 0, 13, 362)$,
$(394, 0, 349, 359)$,
$(372, 359, 145, 242)$,\adfsplit
$(372, 281, 363, 144)$,
$(372, 298, 127, 320)$,
$(372, 57, 222, 256)$,
$(205, 178, 224, 93)$,\adfsplit
$(83, 144, 199, 88)$,
$(72, 261, 31, 80)$,
$(207, 328, 35, 97)$,
$(41, 61, 344, 4)$,\adfsplit
$(160, 321, 203, 228)$,
$(125, 368, 60, 219)$,
$(248, 311, 172, 198)$,
$(166, 89, 321, 68)$,\adfsplit
$(0, 3, 38, 47)$,
$(0, 14, 29, 226)$,
$(0, 16, 86, 281)$,
$(0, 80, 163, 244)$,\adfsplit
$(0, 52, 147, 221)$,
$(0, 45, 133, 269)$,
$(0, 71, 176, 249)$,
$(0, 31, 135, 188)$,\adfsplit
$(0, 33, 134, 185)$,
$(0, 39, 145, 245)$,
$(0, 17, 75, 280)$,
$(0, 40, 99, 242)$,\adfsplit
$(0, 28, 115, 182)$

\adfLgap \noindent by the mapping:
$x \mapsto x + 3 j \adfmod{372}$ for $x < 372$,
$x \mapsto (x +  j \adfmod{12}) + 372$ for $372 \le x < 384$,
$x \mapsto x$ for $x \ge 384$,
$0 \le j < 124$
 for the first 11 blocks;
$x \mapsto x +  j \adfmod{372}$ for $x < 372$,
$x \mapsto (x +  j \adfmod{12}) + 372$ for $372 \le x < 384$,
$x \mapsto x$ for $x \ge 384$,
$0 \le j < 372$
 for the last 26 blocks.
\ADFvfyParStart{(395, ((11, 124, ((372, 3), (12, 1), (11, 11))), (26, 372, ((372, 1), (12, 1), (11, 11)))), ((62, 6), (23, 1)))} 

\adfDgap
\noindent{\boldmath $ 62^{6} 26^{1} $}~
With the point set $Z_{398}$ partitioned into
 residue classes modulo $6$ for $\{0, 1, \dots, 371\}$, and
 $\{372, 373, \dots, 397\}$,
 the design is generated from

\adfLgap 
$(393, 0, 1, 2)$,
$(394, 0, 124, 122)$,
$(395, 0, 247, 371)$,
$(396, 0, 250, 125)$,\adfsplit
$(397, 0, 370, 248)$,
$(372, 79, 366, 236)$,
$(372, 287, 106, 333)$,
$(373, 87, 94, 97)$,\adfsplit
$(373, 83, 318, 98)$,
$(374, 114, 249, 209)$,
$(374, 100, 25, 146)$,
$(375, 86, 370, 139)$,\adfsplit
$(375, 126, 57, 101)$,
$(376, 321, 13, 86)$,
$(376, 142, 347, 162)$,
$(377, 4, 170, 161)$,\adfsplit
$(377, 25, 102, 345)$,
$(378, 223, 184, 110)$,
$(378, 216, 5, 183)$,
$(4, 115, 158, 47)$,\adfsplit
$(349, 27, 371, 16)$,
$(27, 212, 220, 18)$,
$(191, 256, 342, 356)$,
$(369, 53, 346, 157)$,\adfsplit
$(267, 169, 35, 18)$,
$(170, 165, 53, 259)$,
$(70, 215, 99, 133)$,
$(44, 112, 303, 365)$,\adfsplit
$(171, 100, 197, 324)$,
$(288, 335, 81, 172)$,
$(92, 292, 21, 102)$,
$(247, 186, 159, 182)$,\adfsplit
$(215, 16, 348, 145)$,
$(315, 362, 166, 35)$,
$(54, 332, 10, 337)$,
$(12, 127, 281, 302)$,\adfsplit
$(192, 273, 37, 302)$,
$(210, 231, 125, 49)$,
$(99, 154, 42, 74)$,
$(352, 259, 210, 15)$,\adfsplit
$(347, 144, 332, 339)$,
$(358, 194, 133, 59)$,
$(0, 13, 76, 128)$,
$(0, 19, 33, 314)$,\adfsplit
$(0, 27, 31, 344)$,
$(0, 83, 231, 280)$,
$(0, 133, 146, 285)$,
$(0, 75, 106, 277)$,\adfsplit
$(0, 117, 275, 313)$,
$(0, 51, 214, 337)$,
$(0, 26, 93, 193)$,
$(0, 153, 173, 361)$,\adfsplit
$(0, 70, 225, 335)$,
$(0, 103, 119, 160)$,
$(0, 104, 213, 293)$,
$(0, 56, 118, 319)$,\adfsplit
$(0, 37, 69, 136)$,
$(0, 34, 121, 267)$,
$(0, 16, 38, 143)$,
$(0, 41, 99, 271)$,\adfsplit
$(0, 64, 223, 353)$,
$(0, 98, 229, 232)$

\adfLgap \noindent by the mapping:
$x \mapsto x + 3 j \adfmod{372}$ for $x < 372$,
$x \mapsto (x - 372 + 7 j \adfmod{21}) + 372$ for $372 \le x < 393$,
$x \mapsto x$ for $x \ge 393$,
$0 \le j < 124$
 for the first five blocks;
$x \mapsto x + 2 j \adfmod{372}$ for $x < 372$,
$x \mapsto (x - 372 + 7 j \adfmod{21}) + 372$ for $372 \le x < 393$,
$x \mapsto x$ for $x \ge 393$,
$0 \le j < 186$
 for the last 57 blocks.
\ADFvfyParStart{(398, ((5, 124, ((372, 3), (21, 7), (5, 5))), (57, 186, ((372, 2), (21, 7), (5, 5)))), ((62, 6), (26, 1)))} 

\adfDgap
\noindent{\boldmath $ 62^{6} 29^{1} $}~
With the point set $Z_{401}$ partitioned into
 residue classes modulo $6$ for $\{0, 1, \dots, 371\}$, and
 $\{372, 373, \dots, 400\}$,
 the design is generated from

\adfLgap 
$(396, 0, 1, 2)$,
$(397, 0, 124, 122)$,
$(398, 0, 247, 371)$,
$(399, 0, 250, 125)$,\adfsplit
$(400, 0, 370, 248)$,
$(372, 223, 134, 24)$,
$(372, 298, 133, 365)$,
$(372, 270, 308, 333)$,\adfsplit
$(372, 15, 215, 268)$,
$(373, 6, 279, 236)$,
$(373, 307, 202, 285)$,
$(373, 100, 120, 221)$,\adfsplit
$(373, 218, 145, 71)$,
$(159, 98, 265, 322)$,
$(314, 97, 51, 328)$,
$(145, 83, 74, 316)$,\adfsplit
$(83, 122, 117, 90)$,
$(35, 231, 172, 2)$,
$(303, 266, 322, 139)$,
$(172, 252, 13, 203)$,\adfsplit
$(258, 268, 14, 365)$,
$(0, 3, 157, 161)$,
$(0, 11, 26, 285)$,
$(0, 13, 64, 92)$,\adfsplit
$(0, 55, 166, 243)$,
$(0, 23, 68, 291)$,
$(0, 8, 58, 249)$,
$(0, 40, 115, 185)$,\adfsplit
$(0, 29, 94, 146)$,
$(0, 44, 91, 179)$,
$(0, 82, 175, 260)$,
$(0, 16, 116, 219)$,\adfsplit
$(0, 17, 86, 238)$,
$(0, 35, 76, 212)$

\adfLgap \noindent by the mapping:
$x \mapsto x + 3 j \adfmod{372}$ for $x < 372$,
$x \mapsto (x - 372 + 2 j \adfmod{24}) + 372$ for $372 \le x < 396$,
$x \mapsto x$ for $x \ge 396$,
$0 \le j < 124$
 for the first five blocks;
$x \mapsto x +  j \adfmod{372}$ for $x < 372$,
$x \mapsto (x - 372 + 2 j \adfmod{24}) + 372$ for $372 \le x < 396$,
$x \mapsto x$ for $x \ge 396$,
$0 \le j < 372$
 for the last 29 blocks.
\ADFvfyParStart{(401, ((5, 124, ((372, 3), (24, 2), (5, 5))), (29, 372, ((372, 1), (24, 2), (5, 5)))), ((62, 6), (29, 1)))} 

\adfDgap
\noindent{\boldmath $ 62^{6} 32^{1} $}~
With the point set $Z_{404}$ partitioned into
 residue classes modulo $6$ for $\{0, 1, \dots, 371\}$, and
 $\{372, 373, \dots, 403\}$,
 the design is generated from

\adfLgap 
$(399, 0, 1, 2)$,
$(400, 0, 124, 122)$,
$(401, 0, 247, 371)$,
$(402, 0, 250, 125)$,\adfsplit
$(403, 0, 370, 248)$,
$(372, 5, 198, 260)$,
$(372, 117, 46, 211)$,
$(373, 216, 85, 316)$,\adfsplit
$(373, 302, 141, 257)$,
$(374, 319, 30, 15)$,
$(374, 194, 251, 136)$,
$(375, 191, 360, 86)$,\adfsplit
$(375, 343, 40, 57)$,
$(376, 185, 39, 157)$,
$(376, 262, 288, 278)$,
$(377, 264, 31, 39)$,\adfsplit
$(377, 125, 118, 122)$,
$(378, 91, 135, 66)$,
$(378, 346, 83, 8)$,
$(379, 351, 76, 32)$,\adfsplit
$(379, 223, 161, 96)$,
$(380, 115, 60, 11)$,
$(380, 334, 135, 254)$,
$(184, 360, 149, 21)$,\adfsplit
$(55, 164, 126, 292)$,
$(62, 369, 40, 85)$,
$(109, 290, 99, 29)$,
$(94, 295, 2, 324)$,\adfsplit
$(363, 18, 196, 349)$,
$(350, 135, 246, 299)$,
$(118, 255, 194, 295)$,
$(279, 317, 20, 48)$,\adfsplit
$(54, 340, 203, 188)$,
$(295, 12, 75, 125)$,
$(284, 126, 107, 94)$,
$(292, 339, 120, 205)$,\adfsplit
$(365, 316, 129, 108)$,
$(134, 244, 163, 323)$,
$(165, 316, 264, 131)$,
$(352, 26, 77, 259)$,\adfsplit
$(172, 36, 95, 207)$,
$(160, 200, 283, 267)$,
$(4, 229, 236, 24)$,
$(305, 186, 268, 259)$,\adfsplit
$(0, 8, 183, 349)$,
$(0, 11, 14, 103)$,
$(0, 5, 154, 284)$,
$(0, 9, 94, 227)$,\adfsplit
$(0, 91, 146, 267)$,
$(0, 74, 281, 339)$,
$(0, 19, 41, 213)$,
$(0, 77, 116, 309)$,\adfsplit
$(0, 245, 277, 351)$,
$(0, 129, 155, 313)$,
$(0, 70, 215, 325)$,
$(0, 64, 159, 170)$,\adfsplit
$(0, 87, 151, 188)$,
$(0, 229, 249, 305)$,
$(0, 99, 197, 254)$,
$(0, 27, 148, 347)$,\adfsplit
$(0, 31, 131, 260)$,
$(0, 39, 43, 56)$,
$(0, 81, 163, 311)$,
$(0, 68, 237, 367)$

\adfLgap \noindent by the mapping:
$x \mapsto x + 3 j \adfmod{372}$ for $x < 372$,
$x \mapsto (x - 372 + 9 j \adfmod{27}) + 372$ for $372 \le x < 399$,
$x \mapsto x$ for $x \ge 399$,
$0 \le j < 124$
 for the first five blocks;
$x \mapsto x + 2 j \adfmod{372}$ for $x < 372$,
$x \mapsto (x - 372 + 9 j \adfmod{27}) + 372$ for $372 \le x < 399$,
$x \mapsto x$ for $x \ge 399$,
$0 \le j < 186$
 for the last 59 blocks.
\ADFvfyParStart{(404, ((5, 124, ((372, 3), (27, 9), (5, 5))), (59, 186, ((372, 2), (27, 9), (5, 5)))), ((62, 6), (32, 1)))} 

\adfDgap
\noindent{\boldmath $ 62^{6} 35^{1} $}~
With the point set $Z_{407}$ partitioned into
 residue classes modulo $6$ for $\{0, 1, \dots, 371\}$, and
 $\{372, 373, \dots, 406\}$,
 the design is generated from

\adfLgap 
$(396, 0, 1, 2)$,
$(397, 0, 124, 122)$,
$(398, 0, 247, 371)$,
$(399, 0, 250, 125)$,\adfsplit
$(400, 0, 370, 248)$,
$(401, 0, 4, 11)$,
$(402, 0, 7, 368)$,
$(403, 0, 361, 365)$,\adfsplit
$(404, 0, 10, 23)$,
$(405, 0, 13, 362)$,
$(406, 0, 349, 359)$,
$(372, 260, 240, 151)$,\adfsplit
$(372, 181, 17, 147)$,
$(372, 138, 237, 34)$,
$(372, 251, 302, 232)$,
$(373, 199, 70, 359)$,\adfsplit
$(373, 159, 0, 277)$,
$(373, 116, 113, 256)$,
$(373, 242, 309, 114)$,
$(158, 89, 121, 244)$,\adfsplit
$(268, 186, 15, 299)$,
$(200, 333, 48, 317)$,
$(36, 41, 86, 69)$,
$(329, 13, 285, 338)$,\adfsplit
$(0, 245, 320, 147)$,
$(285, 175, 214, 29)$,
$(0, 8, 22, 81)$,
$(0, 15, 41, 238)$,\adfsplit
$(0, 29, 91, 280)$,
$(0, 35, 142, 200)$,
$(0, 55, 166, 260)$,
$(0, 25, 182, 209)$,\adfsplit
$(0, 49, 106, 267)$,
$(0, 21, 64, 257)$,
$(0, 40, 101, 181)$,
$(0, 38, 135, 214)$,\adfsplit
$(0, 63, 137, 202)$,
$(0, 68, 145, 221)$,
$(0, 46, 131, 224)$

\adfLgap \noindent by the mapping:
$x \mapsto x + 3 j \adfmod{372}$ for $x < 372$,
$x \mapsto (x - 372 + 2 j \adfmod{24}) + 372$ for $372 \le x < 396$,
$x \mapsto x$ for $x \ge 396$,
$0 \le j < 124$
 for the first 11 blocks;
$x \mapsto x +  j \adfmod{372}$ for $x < 372$,
$x \mapsto (x - 372 + 2 j \adfmod{24}) + 372$ for $372 \le x < 396$,
$x \mapsto x$ for $x \ge 396$,
$0 \le j < 372$
 for the last 28 blocks.
\ADFvfyParStart{(407, ((11, 124, ((372, 3), (24, 2), (11, 11))), (28, 372, ((372, 1), (24, 2), (11, 11)))), ((62, 6), (35, 1)))} 

\adfDgap
\noindent{\boldmath $ 62^{6} 38^{1} $}~
With the point set $Z_{410}$ partitioned into
 residue classes modulo $6$ for $\{0, 1, \dots, 371\}$, and
 $\{372, 373, \dots, 409\}$,
 the design is generated from

\adfLgap 
$(402, 0, 1, 2)$,
$(403, 0, 124, 122)$,
$(404, 0, 247, 371)$,
$(405, 0, 250, 125)$,\adfsplit
$(406, 0, 370, 248)$,
$(407, 0, 4, 11)$,
$(408, 0, 7, 368)$,
$(409, 0, 361, 365)$,\adfsplit
$(372, 127, 87, 257)$,
$(372, 311, 357, 50)$,
$(372, 132, 260, 184)$,
$(372, 106, 366, 193)$,\adfsplit
$(373, 317, 103, 260)$,
$(373, 135, 323, 208)$,
$(373, 181, 274, 318)$,
$(373, 45, 300, 266)$,\adfsplit
$(374, 360, 328, 51)$,
$(374, 356, 42, 273)$,
$(374, 137, 349, 98)$,
$(374, 307, 94, 287)$,\adfsplit
$(375, 11, 2, 226)$,
$(375, 43, 152, 125)$,
$(375, 24, 273, 301)$,
$(375, 88, 147, 186)$,\adfsplit
$(376, 209, 226, 54)$,
$(376, 320, 271, 129)$,
$(376, 144, 275, 291)$,
$(376, 4, 109, 242)$,\adfsplit
$(61, 341, 117, 354)$,
$(9, 161, 34, 80)$,
$(207, 62, 335, 132)$,
$(228, 75, 166, 361)$,\adfsplit
$(327, 227, 126, 62)$,
$(204, 233, 184, 303)$,
$(136, 167, 31, 162)$,
$(12, 25, 28, 323)$,\adfsplit
$(36, 289, 352, 5)$,
$(26, 298, 109, 312)$,
$(346, 11, 176, 123)$,
$(81, 100, 259, 122)$,\adfsplit
$(107, 97, 56, 273)$,
$(21, 78, 218, 286)$,
$(143, 111, 180, 283)$,
$(36, 3, 202, 221)$,\adfsplit
$(326, 48, 245, 237)$,
$(178, 127, 104, 243)$,
$(0, 3, 254, 271)$,
$(0, 8, 55, 77)$,\adfsplit
$(0, 21, 89, 332)$,
$(0, 10, 221, 255)$,
$(0, 53, 163, 327)$,
$(0, 91, 236, 329)$,\adfsplit
$(0, 27, 50, 142)$,
$(0, 45, 109, 194)$,
$(0, 104, 275, 313)$,
$(0, 106, 285, 343)$,\adfsplit
$(0, 33, 88, 218)$,
$(0, 143, 158, 325)$,
$(0, 141, 176, 367)$,
$(0, 229, 243, 305)$,\adfsplit
$(0, 35, 182, 334)$,
$(0, 15, 269, 290)$,
$(0, 28, 113, 188)$,
$(0, 43, 292, 363)$,\adfsplit
$(0, 67, 146, 262)$,
$(0, 111, 161, 205)$,
$(0, 153, 233, 259)$

\adfLgap \noindent by the mapping:
$x \mapsto x + 3 j \adfmod{372}$ for $x < 372$,
$x \mapsto (x - 372 + 5 j \adfmod{30}) + 372$ for $372 \le x < 402$,
$x \mapsto x$ for $x \ge 402$,
$0 \le j < 124$
 for the first eight blocks;
$x \mapsto x + 2 j \adfmod{372}$ for $x < 372$,
$x \mapsto (x - 372 + 5 j \adfmod{30}) + 372$ for $372 \le x < 402$,
$x \mapsto x$ for $x \ge 402$,
$0 \le j < 186$
 for the last 59 blocks.
\ADFvfyParStart{(410, ((8, 124, ((372, 3), (30, 5), (8, 8))), (59, 186, ((372, 2), (30, 5), (8, 8)))), ((62, 6), (38, 1)))} 

\adfDgap
\noindent{\boldmath $ 62^{6} 41^{1} $}~
With the point set $Z_{413}$ partitioned into
 residue classes modulo $6$ for $\{0, 1, \dots, 371\}$, and
 $\{372, 373, \dots, 412\}$,
 the design is generated from

\adfLgap 
$(408, 0, 1, 2)$,
$(409, 0, 124, 122)$,
$(410, 0, 247, 371)$,
$(411, 0, 250, 125)$,\adfsplit
$(412, 0, 370, 248)$,
$(372, 166, 201, 35)$,
$(372, 195, 41, 26)$,
$(372, 124, 354, 55)$,\adfsplit
$(372, 120, 92, 13)$,
$(373, 87, 199, 274)$,
$(373, 105, 282, 53)$,
$(373, 112, 332, 144)$,\adfsplit
$(373, 182, 47, 25)$,
$(374, 315, 289, 40)$,
$(374, 2, 90, 309)$,
$(374, 227, 368, 286)$,\adfsplit
$(374, 353, 295, 204)$,
$(274, 303, 156, 113)$,
$(19, 58, 111, 50)$,
$(253, 309, 70, 173)$,\adfsplit
$(250, 206, 299, 253)$,
$(100, 200, 209, 96)$,
$(296, 337, 353, 162)$,
$(0, 5, 76, 86)$,\adfsplit
$(0, 7, 40, 242)$,
$(0, 13, 111, 178)$,
$(0, 14, 37, 99)$,
$(0, 45, 115, 209)$,\adfsplit
$(0, 25, 89, 176)$,
$(0, 20, 83, 121)$,
$(0, 51, 119, 224)$,
$(0, 19, 158, 232)$,\adfsplit
$(0, 21, 55, 200)$,
$(0, 17, 127, 243)$,
$(0, 11, 117, 212)$,
$(0, 27, 77, 244)$

\adfLgap \noindent by the mapping:
$x \mapsto x + 3 j \adfmod{372}$ for $x < 372$,
$x \mapsto (x - 372 + 3 j \adfmod{36}) + 372$ for $372 \le x < 408$,
$x \mapsto x$ for $x \ge 408$,
$0 \le j < 124$
 for the first five blocks;
$x \mapsto x +  j \adfmod{372}$ for $x < 372$,
$x \mapsto (x - 372 + 3 j \adfmod{36}) + 372$ for $372 \le x < 408$,
$x \mapsto x$ for $x \ge 408$,
$0 \le j < 372$
 for the last 31 blocks.
\ADFvfyParStart{(413, ((5, 124, ((372, 3), (36, 3), (5, 5))), (31, 372, ((372, 1), (36, 3), (5, 5)))), ((62, 6), (41, 1)))} 

\adfDgap
\noindent{\boldmath $ 62^{6} 44^{1} $}~
With the point set $Z_{416}$ partitioned into
 residue classes modulo $6$ for $\{0, 1, \dots, 371\}$, and
 $\{372, 373, \dots, 415\}$,
 the design is generated from

\adfLgap 
$(408, 0, 1, 2)$,
$(409, 0, 124, 122)$,
$(410, 0, 247, 371)$,
$(411, 0, 250, 125)$,\adfsplit
$(412, 0, 370, 248)$,
$(413, 0, 4, 11)$,
$(414, 0, 7, 368)$,
$(415, 0, 361, 365)$,\adfsplit
$(372, 341, 277, 369)$,
$(372, 72, 344, 256)$,
$(372, 254, 35, 319)$,
$(372, 190, 294, 39)$,\adfsplit
$(373, 367, 171, 5)$,
$(373, 73, 340, 138)$,
$(373, 118, 36, 128)$,
$(373, 371, 194, 333)$,\adfsplit
$(374, 277, 350, 243)$,
$(374, 96, 321, 40)$,
$(374, 53, 346, 331)$,
$(374, 150, 332, 323)$,\adfsplit
$(375, 303, 136, 98)$,
$(375, 225, 126, 109)$,
$(375, 47, 368, 334)$,
$(375, 300, 31, 137)$,\adfsplit
$(376, 266, 301, 293)$,
$(376, 252, 238, 143)$,
$(376, 183, 340, 20)$,
$(376, 67, 213, 42)$,\adfsplit
$(377, 67, 102, 371)$,
$(377, 74, 111, 106)$,
$(377, 213, 340, 257)$,
$(377, 32, 360, 109)$,\adfsplit
$(292, 123, 67, 212)$,
$(223, 342, 244, 177)$,
$(336, 83, 266, 105)$,
$(176, 335, 22, 205)$,\adfsplit
$(29, 291, 193, 56)$,
$(57, 368, 306, 77)$,
$(162, 98, 347, 238)$,
$(23, 105, 259, 166)$,\adfsplit
$(192, 146, 167, 327)$,
$(173, 218, 141, 100)$,
$(93, 2, 293, 18)$,
$(46, 307, 105, 192)$,\adfsplit
$(139, 136, 35, 210)$,
$(223, 183, 4, 140)$,
$(175, 224, 30, 10)$,
$(109, 150, 20, 29)$,\adfsplit
$(0, 13, 257, 332)$,
$(0, 8, 31, 208)$,
$(0, 17, 68, 87)$,
$(0, 15, 28, 67)$,\adfsplit
$(0, 47, 94, 152)$,
$(0, 26, 127, 160)$,
$(0, 111, 161, 343)$,
$(0, 112, 287, 349)$,\adfsplit
$(0, 151, 303, 329)$,
$(0, 22, 155, 213)$,
$(0, 45, 259, 322)$,
$(0, 63, 266, 335)$,\adfsplit
$(0, 33, 217, 256)$,
$(0, 110, 223, 239)$,
$(0, 71, 128, 369)$,
$(0, 49, 142, 273)$,\adfsplit
$(0, 107, 193, 224)$,
$(0, 86, 243, 319)$,
$(0, 57, 199, 317)$,
$(0, 166, 353, 367)$,\adfsplit
$(0, 55, 196, 293)$

\adfLgap \noindent by the mapping:
$x \mapsto x + 3 j \adfmod{372}$ for $x < 372$,
$x \mapsto (x - 372 + 6 j \adfmod{36}) + 372$ for $372 \le x < 408$,
$x \mapsto x$ for $x \ge 408$,
$0 \le j < 124$
 for the first eight blocks;
$x \mapsto x + 2 j \adfmod{372}$ for $x < 372$,
$x \mapsto (x - 372 + 6 j \adfmod{36}) + 372$ for $372 \le x < 408$,
$x \mapsto x$ for $x \ge 408$,
$0 \le j < 186$
 for the last 61 blocks.
\ADFvfyParStart{(416, ((8, 124, ((372, 3), (36, 6), (8, 8))), (61, 186, ((372, 2), (36, 6), (8, 8)))), ((62, 6), (44, 1)))} 

\adfDgap
\noindent{\boldmath $ 62^{6} 47^{1} $}~
With the point set $Z_{419}$ partitioned into
 residue classes modulo $6$ for $\{0, 1, \dots, 371\}$, and
 $\{372, 373, \dots, 418\}$,
 the design is generated from

\adfLgap 
$(408, 0, 1, 2)$,
$(409, 0, 124, 122)$,
$(410, 0, 247, 371)$,
$(411, 0, 250, 125)$,\adfsplit
$(412, 0, 370, 248)$,
$(413, 0, 4, 11)$,
$(414, 0, 7, 368)$,
$(415, 0, 361, 365)$,\adfsplit
$(416, 0, 10, 23)$,
$(417, 0, 13, 362)$,
$(418, 0, 349, 359)$,
$(372, 211, 48, 94)$,\adfsplit
$(372, 113, 50, 309)$,
$(372, 198, 87, 47)$,
$(372, 296, 88, 193)$,
$(373, 14, 159, 331)$,\adfsplit
$(373, 148, 308, 216)$,
$(373, 318, 185, 141)$,
$(373, 263, 370, 73)$,
$(374, 335, 236, 148)$,\adfsplit
$(374, 108, 142, 279)$,
$(374, 305, 253, 225)$,
$(374, 150, 278, 367)$,
$(314, 178, 331, 117)$,\adfsplit
$(14, 259, 208, 107)$,
$(283, 240, 17, 352)$,
$(342, 173, 109, 244)$,
$(313, 179, 132, 140)$,\adfsplit
$(0, 3, 22, 290)$,
$(0, 5, 31, 345)$,
$(0, 9, 109, 251)$,
$(0, 15, 77, 110)$,\adfsplit
$(0, 14, 154, 179)$,
$(0, 16, 57, 241)$,
$(0, 20, 73, 143)$,
$(0, 29, 79, 170)$,\adfsplit
$(0, 45, 119, 206)$,
$(0, 76, 157, 243)$,
$(0, 35, 118, 183)$,
$(0, 49, 116, 275)$,\adfsplit
$(0, 21, 59, 115)$

\adfLgap \noindent by the mapping:
$x \mapsto x + 3 j \adfmod{372}$ for $x < 372$,
$x \mapsto (x - 372 + 3 j \adfmod{36}) + 372$ for $372 \le x < 408$,
$x \mapsto x$ for $x \ge 408$,
$0 \le j < 124$
 for the first 11 blocks;
$x \mapsto x +  j \adfmod{372}$ for $x < 372$,
$x \mapsto (x - 372 + 3 j \adfmod{36}) + 372$ for $372 \le x < 408$,
$x \mapsto x$ for $x \ge 408$,
$0 \le j < 372$
 for the last 30 blocks.
\ADFvfyParStart{(419, ((11, 124, ((372, 3), (36, 3), (11, 11))), (30, 372, ((372, 1), (36, 3), (11, 11)))), ((62, 6), (47, 1)))} 

\adfDgap
\noindent{\boldmath $ 62^{6} 143^{1} $}~
With the point set $Z_{515}$ partitioned into
 residue classes modulo $6$ for $\{0, 1, \dots, 371\}$, and
 $\{372, 373, \dots, 514\}$,
 the design is generated from

\adfLgap 
$(504, 0, 1, 2)$,
$(505, 0, 124, 122)$,
$(506, 0, 247, 371)$,
$(507, 0, 250, 125)$,\adfsplit
$(508, 0, 370, 248)$,
$(509, 0, 4, 11)$,
$(510, 0, 7, 368)$,
$(511, 0, 361, 365)$,\adfsplit
$(512, 0, 10, 23)$,
$(513, 0, 13, 362)$,
$(514, 0, 349, 359)$,
$(372, 67, 160, 326)$,\adfsplit
$(372, 59, 264, 51)$,
$(372, 82, 277, 140)$,
$(372, 65, 366, 297)$,
$(373, 213, 14, 124)$,\adfsplit
$(373, 354, 205, 334)$,
$(373, 96, 231, 53)$,
$(373, 67, 287, 140)$,
$(374, 71, 332, 169)$,\adfsplit
$(374, 235, 209, 28)$,
$(374, 322, 183, 228)$,
$(374, 66, 122, 141)$,
$(375, 204, 368, 241)$,\adfsplit
$(375, 230, 87, 269)$,
$(375, 94, 129, 43)$,
$(375, 143, 246, 184)$,
$(376, 28, 336, 283)$,\adfsplit
$(376, 73, 33, 248)$,
$(376, 107, 123, 46)$,
$(376, 29, 182, 342)$,
$(377, 305, 218, 171)$,\adfsplit
$(377, 85, 52, 348)$,
$(377, 330, 226, 212)$,
$(377, 7, 345, 191)$,
$(378, 260, 306, 178)$,\adfsplit
$(378, 31, 98, 123)$,
$(378, 371, 316, 129)$,
$(378, 204, 209, 97)$,
$(379, 314, 33, 305)$,\adfsplit
$(379, 63, 92, 107)$,
$(379, 100, 258, 127)$,
$(379, 48, 118, 217)$,
$(0, 17, 292, 435)$,\adfsplit
$(0, 3, 88, 145)$,
$(0, 22, 101, 224)$,
$(0, 21, 344, 446)$,
$(0, 31, 83, 402)$,\adfsplit
$(0, 32, 221, 424)$,
$(0, 95, 200, 381)$,
$(0, 65, 211, 459)$,
$(0, 38, 106, 415)$,\adfsplit
$(0, 115, 231, 503)$,
$(0, 63, 239, 481)$,
$(0, 74, 193, 480)$,
$(0, 81, 217, 502)$,\adfsplit
$(0, 50, 251, 436)$

\adfLgap \noindent by the mapping:
$x \mapsto x + 3 j \adfmod{372}$ for $x < 372$,
$x \mapsto (x - 372 + 11 j \adfmod{132}) + 372$ for $372 \le x < 504$,
$x \mapsto x$ for $x \ge 504$,
$0 \le j < 124$
 for the first 11 blocks;
$x \mapsto x +  j \adfmod{372}$ for $x < 372$,
$x \mapsto (x - 372 + 11 j \adfmod{132}) + 372$ for $372 \le x < 504$,
$x \mapsto x$ for $x \ge 504$,
$0 \le j < 372$
 for the last 46 blocks.
\ADFvfyParStart{(515, ((11, 124, ((372, 3), (132, 11), (11, 11))), (46, 372, ((372, 1), (132, 11), (11, 11)))), ((62, 6), (143, 1)))} 

\adfDgap
\noindent{\boldmath $ 62^{6} 146^{1} $}~
With the point set $Z_{518}$ partitioned into
 residue classes modulo $6$ for $\{0, 1, \dots, 371\}$, and
 $\{372, 373, \dots, 517\}$,
 the design is generated from

\adfLgap 
$(510, 0, 1, 2)$,
$(511, 0, 124, 122)$,
$(512, 0, 247, 371)$,
$(513, 0, 250, 125)$,\adfsplit
$(514, 0, 370, 248)$,
$(515, 0, 4, 11)$,
$(516, 0, 7, 368)$,
$(517, 0, 361, 365)$,\adfsplit
$(372, 251, 87, 366)$,
$(372, 302, 264, 305)$,
$(372, 130, 9, 121)$,
$(372, 200, 340, 163)$,\adfsplit
$(373, 23, 141, 226)$,
$(373, 324, 61, 136)$,
$(373, 137, 350, 78)$,
$(373, 39, 224, 91)$,\adfsplit
$(374, 341, 123, 8)$,
$(374, 158, 175, 328)$,
$(374, 263, 60, 142)$,
$(374, 333, 282, 265)$,\adfsplit
$(375, 326, 151, 159)$,
$(375, 82, 227, 13)$,
$(375, 348, 328, 357)$,
$(375, 116, 186, 305)$,\adfsplit
$(376, 242, 135, 22)$,
$(376, 241, 48, 100)$,
$(376, 167, 93, 128)$,
$(376, 126, 283, 173)$,\adfsplit
$(377, 224, 13, 370)$,
$(377, 156, 319, 287)$,
$(377, 162, 122, 33)$,
$(377, 15, 160, 209)$,\adfsplit
$(378, 278, 234, 85)$,
$(378, 39, 55, 371)$,
$(378, 225, 252, 244)$,
$(378, 365, 248, 94)$,\adfsplit
$(379, 360, 188, 124)$,
$(379, 334, 7, 359)$,
$(379, 318, 147, 241)$,
$(379, 101, 369, 350)$,\adfsplit
$(380, 71, 61, 327)$,
$(380, 21, 44, 317)$,
$(380, 258, 242, 352)$,
$(380, 250, 319, 84)$,\adfsplit
$(381, 34, 302, 199)$,
$(381, 203, 13, 364)$,
$(381, 258, 147, 161)$,
$(381, 309, 92, 60)$,\adfsplit
$(382, 95, 181, 324)$,
$(382, 356, 306, 261)$,
$(382, 10, 7, 159)$,
$(382, 293, 218, 196)$,\adfsplit
$(383, 118, 333, 42)$,
$(383, 284, 255, 0)$,
$(383, 235, 340, 47)$,
$(383, 289, 362, 161)$,\adfsplit
$(384, 162, 193, 208)$,
$(384, 71, 358, 120)$,
$(384, 305, 248, 243)$,
$(384, 369, 170, 343)$,\adfsplit
$(385, 148, 161, 296)$,
$(385, 115, 198, 38)$,
$(385, 250, 36, 171)$,
$(385, 145, 345, 287)$,\adfsplit
$(386, 60, 215, 361)$,
$(386, 333, 74, 31)$,
$(386, 142, 185, 207)$,
$(386, 366, 8, 172)$,\adfsplit
$(387, 103, 65, 196)$,
$(387, 145, 206, 243)$,
$(387, 131, 312, 20)$,
$(387, 186, 21, 130)$,\adfsplit
$(388, 243, 307, 362)$,
$(388, 359, 296, 40)$,
$(388, 261, 29, 234)$,
$(388, 216, 334, 277)$,\adfsplit
$(389, 322, 350, 31)$,
$(389, 281, 180, 363)$,
$(389, 28, 270, 249)$,
$(389, 169, 332, 299)$,\adfsplit
$(0, 23, 74, 307)$,
$(0, 33, 266, 413)$,
$(0, 10, 68, 391)$,
$(0, 5, 286, 392)$,\adfsplit
$(0, 35, 260, 393)$,
$(0, 55, 182, 274)$,
$(0, 83, 127, 394)$,
$(0, 67, 347, 417)$,\adfsplit
$(0, 62, 133, 509)$,
$(0, 34, 137, 440)$,
$(0, 99, 305, 459)$,
$(0, 26, 177, 461)$,\adfsplit
$(0, 95, 142, 482)$,
$(0, 89, 176, 508)$,
$(0, 87, 128, 460)$,
$(1, 29, 177, 414)$,\adfsplit
$(0, 231, 281, 484)$,
$(0, 181, 341, 485)$,
$(0, 225, 359, 507)$,
$(0, 129, 175, 462)$,\adfsplit
$(0, 73, 107, 309)$,
$(0, 213, 313, 483)$,
$(0, 105, 185, 436)$

\adfLgap \noindent by the mapping:
$x \mapsto x + 3 j \adfmod{372}$ for $x < 372$,
$x \mapsto (x - 372 + 23 j \adfmod{138}) + 372$ for $372 \le x < 510$,
$x \mapsto x$ for $x \ge 510$,
$0 \le j < 124$
 for the first eight blocks;
$x \mapsto x + 2 j \adfmod{372}$ for $x < 372$,
$x \mapsto (x - 372 + 23 j \adfmod{138}) + 372$ for $372 \le x < 510$,
$x \mapsto x$ for $x \ge 510$,
$0 \le j < 186$
 for the last 95 blocks.
\ADFvfyParStart{(518, ((8, 124, ((372, 3), (138, 23), (8, 8))), (95, 186, ((372, 2), (138, 23), (8, 8)))), ((62, 6), (146, 1)))} 

\adfDgap
\noindent{\boldmath $ 62^{6} 149^{1} $}~
With the point set $Z_{521}$ partitioned into
 residue classes modulo $6$ for $\{0, 1, \dots, 371\}$, and
 $\{372, 373, \dots, 520\}$,
 the design is generated from

\adfLgap 
$(516, 0, 1, 2)$,
$(517, 0, 124, 122)$,
$(518, 0, 247, 371)$,
$(519, 0, 250, 125)$,\adfsplit
$(520, 0, 370, 248)$,
$(372, 343, 128, 71)$,
$(372, 313, 53, 86)$,
$(372, 346, 231, 24)$,\adfsplit
$(372, 153, 126, 16)$,
$(373, 260, 358, 139)$,
$(373, 158, 11, 57)$,
$(373, 1, 276, 171)$,\adfsplit
$(373, 293, 52, 114)$,
$(374, 233, 91, 78)$,
$(374, 236, 157, 311)$,
$(374, 262, 96, 51)$,\adfsplit
$(374, 88, 165, 110)$,
$(375, 298, 114, 357)$,
$(375, 304, 19, 251)$,
$(375, 325, 219, 308)$,\adfsplit
$(375, 29, 98, 72)$,
$(376, 67, 160, 150)$,
$(376, 353, 46, 12)$,
$(376, 57, 152, 311)$,\adfsplit
$(376, 315, 277, 86)$,
$(377, 262, 176, 3)$,
$(377, 352, 206, 157)$,
$(377, 173, 102, 213)$,\adfsplit
$(377, 251, 331, 360)$,
$(378, 33, 172, 161)$,
$(378, 19, 327, 167)$,
$(378, 362, 82, 246)$,\adfsplit
$(378, 108, 152, 193)$,
$(379, 162, 134, 185)$,
$(379, 192, 40, 59)$,
$(379, 214, 115, 27)$,\adfsplit
$(379, 33, 25, 152)$,
$(380, 115, 346, 42)$,
$(380, 131, 145, 320)$,
$(380, 291, 101, 338)$,\adfsplit
$(0, 3, 7, 416)$,
$(0, 9, 278, 405)$,
$(0, 16, 74, 430)$,
$(0, 37, 107, 395)$,\adfsplit
$(0, 15, 209, 453)$,
$(0, 25, 196, 513)$,
$(0, 67, 236, 429)$,
$(0, 52, 134, 442)$,\adfsplit
$(0, 39, 281, 514)$,
$(0, 35, 158, 490)$,
$(0, 32, 255, 455)$,
$(0, 21, 221, 431)$,\adfsplit
$(0, 63, 167, 479)$,
$(0, 5, 61, 81)$

\adfLgap \noindent by the mapping:
$x \mapsto x + 3 j \adfmod{372}$ for $x < 372$,
$x \mapsto (x - 372 + 12 j \adfmod{144}) + 372$ for $372 \le x < 516$,
$x \mapsto x$ for $x \ge 516$,
$0 \le j < 124$
 for the first five blocks;
$x \mapsto x +  j \adfmod{372}$ for $x < 372$,
$x \mapsto (x - 372 + 12 j \adfmod{144}) + 372$ for $372 \le x < 516$,
$x \mapsto x$ for $x \ge 516$,
$0 \le j < 372$
 for the last 49 blocks.
\ADFvfyParStart{(521, ((5, 124, ((372, 3), (144, 12), (5, 5))), (49, 372, ((372, 1), (144, 12), (5, 5)))), ((62, 6), (149, 1)))} 

\adfDgap
\noindent{\boldmath $ 62^{6} 152^{1} $}~
With the point set $Z_{524}$ partitioned into
 residue classes modulo $6$ for $\{0, 1, \dots, 371\}$, and
 $\{372, 373, \dots, 523\}$,
 the design is generated from

\adfLgap 
$(516, 0, 1, 2)$,
$(517, 0, 124, 122)$,
$(518, 0, 247, 371)$,
$(519, 0, 250, 125)$,\adfsplit
$(520, 0, 370, 248)$,
$(521, 0, 4, 11)$,
$(522, 0, 7, 368)$,
$(523, 0, 361, 365)$,\adfsplit
$(372, 23, 102, 225)$,
$(372, 62, 118, 15)$,
$(372, 193, 164, 77)$,
$(372, 76, 247, 108)$,\adfsplit
$(373, 127, 299, 122)$,
$(373, 176, 85, 231)$,
$(373, 24, 64, 81)$,
$(373, 298, 54, 245)$,\adfsplit
$(374, 93, 116, 276)$,
$(374, 302, 95, 193)$,
$(374, 280, 354, 199)$,
$(374, 257, 231, 130)$,\adfsplit
$(375, 78, 283, 63)$,
$(375, 74, 155, 93)$,
$(375, 356, 328, 325)$,
$(375, 336, 161, 70)$,\adfsplit
$(376, 42, 28, 299)$,
$(376, 204, 7, 315)$,
$(376, 302, 118, 305)$,
$(376, 237, 80, 13)$,\adfsplit
$(377, 124, 350, 39)$,
$(377, 248, 130, 71)$,
$(377, 365, 103, 66)$,
$(377, 205, 180, 249)$,\adfsplit
$(378, 193, 215, 254)$,
$(378, 204, 357, 124)$,
$(378, 7, 138, 245)$,
$(378, 80, 183, 310)$,\adfsplit
$(379, 88, 278, 69)$,
$(379, 145, 252, 274)$,
$(379, 66, 215, 31)$,
$(379, 327, 176, 5)$,\adfsplit
$(380, 253, 42, 267)$,
$(380, 242, 83, 24)$,
$(380, 46, 79, 320)$,
$(380, 329, 249, 28)$,\adfsplit
$(381, 299, 68, 268)$,
$(381, 265, 183, 156)$,
$(381, 317, 213, 31)$,
$(381, 234, 10, 2)$,\adfsplit
$(382, 175, 272, 293)$,
$(382, 310, 85, 228)$,
$(382, 167, 75, 304)$,
$(382, 321, 114, 278)$,\adfsplit
$(383, 40, 110, 156)$,
$(383, 83, 272, 334)$,
$(383, 271, 117, 281)$,
$(383, 145, 162, 111)$,\adfsplit
$(384, 219, 227, 154)$,
$(384, 306, 257, 369)$,
$(384, 368, 199, 132)$,
$(384, 230, 364, 253)$,\adfsplit
$(385, 57, 84, 133)$,
$(385, 166, 90, 271)$,
$(385, 347, 27, 110)$,
$(385, 232, 53, 128)$,\adfsplit
$(386, 341, 321, 127)$,
$(386, 193, 98, 48)$,
$(386, 280, 236, 66)$,
$(386, 159, 322, 119)$,\adfsplit
$(387, 359, 163, 93)$,
$(387, 230, 10, 288)$,
$(387, 64, 291, 337)$,
$(387, 366, 29, 20)$,\adfsplit
$(388, 328, 53, 369)$,
$(388, 102, 296, 7)$,
$(388, 157, 299, 144)$,
$(388, 51, 338, 238)$,\adfsplit
$(389, 80, 132, 299)$,
$(389, 259, 221, 6)$,
$(389, 63, 262, 229)$,
$(389, 26, 136, 249)$,\adfsplit
$(390, 141, 67, 52)$,
$(390, 132, 166, 123)$,
$(0, 10, 87, 438)$,
$(0, 39, 71, 510)$,\adfsplit
$(0, 47, 88, 141)$,
$(0, 51, 286, 391)$,
$(0, 20, 303, 392)$,
$(0, 16, 133, 464)$,\adfsplit
$(0, 129, 242, 512)$,
$(1, 29, 129, 392)$,
$(0, 45, 166, 393)$,
$(0, 159, 196, 442)$,\adfsplit
$(0, 75, 317, 490)$,
$(0, 68, 161, 418)$,
$(0, 38, 347, 511)$,
$(0, 119, 351, 514)$,\adfsplit
$(0, 179, 315, 463)$,
$(0, 239, 327, 439)$,
$(0, 79, 92, 465)$,
$(0, 249, 307, 441)$,\adfsplit
$(0, 99, 115, 489)$,
$(0, 199, 267, 395)$,
$(0, 64, 343, 491)$,
$(0, 135, 295, 419)$,\adfsplit
$(0, 112, 367, 515)$

\adfLgap \noindent by the mapping:
$x \mapsto x + 3 j \adfmod{372}$ for $x < 372$,
$x \mapsto (x - 372 + 24 j \adfmod{144}) + 372$ for $372 \le x < 516$,
$x \mapsto x$ for $x \ge 516$,
$0 \le j < 124$
 for the first eight blocks;
$x \mapsto x + 2 j \adfmod{372}$ for $x < 372$,
$x \mapsto (x - 372 + 24 j \adfmod{144}) + 372$ for $372 \le x < 516$,
$x \mapsto x$ for $x \ge 516$,
$0 \le j < 186$
 for the last 97 blocks.
\ADFvfyParStart{(524, ((8, 124, ((372, 3), (144, 24), (8, 8))), (97, 186, ((372, 2), (144, 24), (8, 8)))), ((62, 6), (152, 1)))} 

\section{4-GDDs for the proof of Lemma \ref{lem:4-GDD 76^u m^1}}
\label{app:4-GDD 76^u m^1}
\adfhide{
$ 76^6 13^1 $ and
$ 76^6 16^1 $.
}

\adfDgap
\noindent{\boldmath $ 76^{6} 13^{1} $}~
With the point set $Z_{469}$ partitioned into
 residue classes modulo $6$ for $\{0, 1, \dots, 455\}$, and
 $\{456, 457, \dots, 468\}$,
 the design is generated from

\adfLgap 
$(405, 84, 52, 456)$,
$(233, 156, 19, 457)$,
$(364, 386, 423, 457)$,
$(301, 154, 450, 458)$,\adfsplit
$(249, 338, 395, 458)$,
$(252, 109, 203, 459)$,
$(326, 100, 27, 459)$,
$(202, 207, 365, 115)$,\adfsplit
$(4, 2, 277, 209)$,
$(176, 388, 121, 81)$,
$(416, 405, 409, 306)$,
$(155, 147, 326, 421)$,\adfsplit
$(176, 291, 217, 288)$,
$(134, 177, 64, 199)$,
$(201, 56, 421, 328)$,
$(335, 453, 286, 210)$,\adfsplit
$(247, 201, 377, 238)$,
$(440, 203, 424, 192)$,
$(400, 234, 296, 279)$,
$(344, 198, 269, 183)$,\adfsplit
$(398, 298, 264, 405)$,
$(401, 172, 303, 386)$,
$(110, 433, 105, 348)$,
$(256, 227, 140, 420)$,\adfsplit
$(314, 222, 17, 183)$,
$(164, 351, 245, 322)$,
$(422, 247, 419, 340)$,
$(303, 128, 148, 179)$,\adfsplit
$(388, 320, 33, 61)$,
$(116, 341, 63, 72)$,
$(265, 88, 299, 324)$,
$(251, 272, 130, 90)$,\adfsplit
$(239, 156, 439, 206)$,
$(88, 19, 260, 287)$,
$(38, 189, 240, 388)$,
$(412, 36, 109, 411)$,\adfsplit
$(210, 65, 4, 175)$,
$(437, 421, 225, 402)$,
$(121, 144, 339, 284)$,
$(135, 269, 252, 319)$,\adfsplit
$(380, 227, 48, 229)$,
$(0, 1, 196, 209)$,
$(0, 21, 91, 430)$,
$(0, 10, 63, 266)$,\adfsplit
$(0, 4, 227, 429)$,
$(0, 8, 328, 351)$,
$(0, 14, 111, 382)$,
$(0, 58, 326, 415)$,\adfsplit
$(0, 75, 286, 391)$,
$(0, 118, 237, 349)$,
$(0, 64, 247, 411)$,
$(0, 52, 317, 393)$,\adfsplit
$(0, 86, 289, 377)$,
$(0, 46, 345, 355)$,
$(0, 337, 399, 437)$,
$(0, 241, 267, 389)$,\adfsplit
$(0, 137, 193, 297)$,
$(0, 123, 143, 259)$,
$(0, 38, 239, 409)$,
$(0, 56, 154, 327)$,\adfsplit
$(0, 39, 278, 443)$,
$(0, 333, 413, 464)$,
$(0, 109, 191, 423)$,
$(0, 69, 133, 194)$,\adfsplit
$(0, 301, 315, 359)$,
$(0, 25, 255, 287)$,
$(0, 28, 122, 249)$,
$(468, 0, 152, 304)$,\adfsplit
$(468, 1, 153, 305)$

\adfLgap \noindent by the mapping:
$x \mapsto x + 2 j \adfmod{456}$ for $x < 456$,
$x \mapsto (x + 4 j \adfmod{12}) + 456$ for $456 \le x < 468$,
$468 \mapsto 468$,
$0 \le j < 228$
 for the first 67 blocks,
$0 \le j < 76$
 for the last two blocks.
\ADFvfyParStart{(469, ((67, 228, ((456, 2), (12, 4), (1, 1))), (2, 76, ((456, 2), (12, 4), (1, 1)))), ((76, 6), (13, 1)))} 

\adfDgap
\noindent{\boldmath $ 76^{6} 16^{1} $}~
With the point set $Z_{472}$ partitioned into
 residue classes modulo $6$ for $\{0, 1, \dots, 455\}$, and
 $\{456, 457, \dots, 471\}$,
 the design is generated from

\adfLgap 
$(456, 438, 419, 136)$,
$(457, 409, 162, 2)$,
$(458, 165, 385, 332)$,
$(459, 167, 223, 72)$,\adfsplit
$(460, 413, 249, 421)$,
$(424, 57, 427, 320)$,
$(38, 81, 83, 277)$,
$(229, 329, 404, 126)$,\adfsplit
$(381, 350, 330, 167)$,
$(364, 66, 183, 194)$,
$(211, 45, 216, 332)$,
$(320, 304, 205, 389)$,\adfsplit
$(237, 376, 294, 215)$,
$(293, 164, 163, 418)$,
$(12, 447, 424, 409)$,
$(45, 277, 302, 400)$,\adfsplit
$(222, 129, 181, 430)$,
$(299, 66, 292, 332)$,
$(24, 374, 291, 53)$,
$(0, 4, 13, 398)$,\adfsplit
$(0, 10, 27, 347)$,
$(0, 14, 77, 147)$,
$(0, 26, 81, 131)$,
$(0, 32, 142, 269)$,\adfsplit
$(0, 61, 148, 261)$,
$(0, 47, 111, 229)$,
$(0, 68, 211, 303)$,
$(0, 67, 191, 279)$,\adfsplit
$(0, 34, 179, 259)$,
$(0, 91, 188, 334)$,
$(0, 28, 165, 299)$,
$(0, 37, 76, 149)$,\adfsplit
$(0, 46, 140, 205)$,
$(0, 35, 176, 250)$,
$(471, 0, 152, 304)$

\adfLgap \noindent by the mapping:
$x \mapsto x +  j \adfmod{456}$ for $x < 456$,
$x \mapsto (x - 456 + 5 j \adfmod{15}) + 456$ for $456 \le x < 471$,
$471 \mapsto 471$,
$0 \le j < 456$
 for the first 34 blocks,
$0 \le j < 152$
 for the last block.
\ADFvfyParStart{(472, ((34, 456, ((456, 1), (15, 5), (1, 1))), (1, 152, ((456, 1), (15, 5), (1, 1)))), ((76, 6), (16, 1)))} 


\begin{thebibliography}{99}

\bibitem{AbelColbournDinitz2007} R. J. R. Abel, C. J. Colbourn and J. H. Dinitz,
    Mutually Orthogonal Latin Squares (MOLS),
    \textit{Handbook of Combinatorial Designs}, second edition (C. J. Colbourn and J. H. Dinitz, eds), Chapman \& Hall/CRC Press, London, 2007, 160--193.

\bibitem{BrouwerSchrijverHanani1977} A. E. Brouwer, A. Schrijver and H. Hanani, Group divisible designs with block size four,
    \textit{Discrete Math.} \textbf{20} (1977), 1--10.

\bibitem{CaoWangWei2009} H. Cao, L. Wang and R. Wei,
    The existence of HGDDs with block size four and its application to double frames,
    \textit{Discrete Mathematics} \textbf{2009} (2002), 945--949.

\bibitem{ColbournHoffmanRees1992} C. J. Colbourn, D. G. Hoffman and R. S. Rees,
    A new class of group divisible designs with block size three,
    \textit{J. Combin. Theory, Series A} \textbf{59} (1992), 73--89.

\bibitem{Forbes2019} A. D. Forbes,
    Group divisible designs with block size four and type $g^u m^1$ -- II,
    \textit{J. Combin. Des.} \textbf{27} (2019), 311--349.

\bibitem{ForbesForbes2018} A. D. Forbes and K. A. Forbes,
    Group divisible designs with block size 4 and type $g^u m^1$,
    \textit{J. Combin. Des.} \textbf{26} (2018), 519--539.

\bibitem{Ge2007} G. Ge, Group Divisible Designs,
    \textit{Handbook of Combinatorial Designs}, second edition (ed. C. J. Colbourn and J. H. Dinitz), Chapman \& Hall/CRC Press, 2007, 255--260.

\bibitem{GeLing2004} G. Ge and A. C. H. Ling,
    Group divisible designs with block size four and group type $g^u m^1$ for small $g$,
    \textit{Discrete Math.} \textbf{285} (2004), 97--120.

\bibitem{GeLing2005} G. Ge and A. C. H. Ling,
    Group divisible designs with block size four and group type $g^u m^1$ with minimum $m$,
    \textit{Designs, Codes and Cryptography} \textbf{34} (2005), 117--126.

\bibitem{GeRees2002} G. Ge and R. S. Rees,
    On group-divisible designs with block size four and group-type $g^u m^1$,
    \textit{Designs, Codes and Cryptography} \textbf{27} (2002), 5--24.

\bibitem{GeRees2004} G. Ge and R. S. Rees,
    On group-divisible designs with block size four and group type $6^u m^1$,
    \textit{Discrete Math.} \textbf{279} (2004), 247--265.

\bibitem{GeReesZhu2002} G. Ge, R. S. Rees and L. Zhu,
    Group-divisible designs with block size four and group-type $g^u m^1$ with $m$ as large or as small as possible,
    \textit{J. Combin. Theory, Series A} \textbf{98} (2002), 357--376.

\bibitem{GeWei2004} G. Ge, R. Wei,
    HGDDs with block size four,
    \textit{Discrete Math.} \textbf{279} (2004), 267--276.

\bibitem{Hanini1975} H. Hanani,
   Balanced incomplete block designs and related designs,
   \textit{Discrete Math.} \textbf{11} (1975), 255--369.

\bibitem{Schuster2010} E. Schuster,
    Group divisible designs with block size four and group type $g^u m^1$ where $g$ is a multiple of $8$,
    \textit{Discrete Math.} \textbf{310} (2010), 2258--2270.

\bibitem{Schuster2014} E. Schuster,
    New classes of group divisible designs with block size 4 and group type $g^u m^1$,
    \textit{J. Combin. Math. Combin. Comput.} \textbf{91} (2014), 65--105.

\bibitem{WeiGe2013} H. Wei and G. Ge,
    Group divisible designs with block size four and group type $g^u m^1$ for more small $g$,
    \textit{Discrete Math.} \textbf{313} (2013), 2065--2083.

\bibitem{WeiGe2014} H. Wei and G. Ge,
    Group Divisible designs with block size four and group type $g^u m^1$ for $g \equiv 0 \adfmod{6}$,
    \textit{J. Combin. Des.} \textbf{22} (2014), 26--52.

\bibitem{WeiGe2015} H. Wei and G. Ge,
    Group divisible designs with block size four and group type $g^u m^1$,
    \textit{Designs, Codes and Cryptography} \textbf{74} (2015), 243--282.

\end{thebibliography}
\end{document}